\documentclass[12pt,a4paper,twoside,openright]{amsbook}

\usepackage[latin1]{inputenc}
\usepackage[french,english]{babel}
\usepackage{indentfirst}
\usepackage{paralist}
\usepackage{amssymb}
\usepackage{amsthm}
\usepackage{amsmath}
\usepackage{mathrsfs}
\usepackage{amscd}
\usepackage{graphicx}
\usepackage{hyperref}

\usepackage{wrapfig}

\usepackage[margin=1in]{geometry}
\usepackage{cite}
\usepackage{lineno} 
\usepackage{enumitem} 
\usepackage{accents} 

\usepackage{minitoc}
\hypersetup{
	colorlinks=true,
	linkcolor=blue,
	citecolor=red,
	urlcolor=black,
	linktoc=all
}

\theoremstyle{plain} 
\newtheorem{theorem}{Theorem}[section] 
\newtheorem{corollary}[theorem]{Corollary} 
\newtheorem{lemma}[theorem]{Lemma}
\newtheorem{prop}[theorem]{Proposition}
\newtheorem{defn}[theorem]{Definition}

\theoremstyle{definition}
\newtheorem{example}{Example}[section]
\newtheorem{remark}[theorem]{Remark}

\theoremstyle{plain} 
\newtheorem{frtheorem}{Th\'eor\`eme}[section] 
\newtheorem{frcorollary}[frtheorem]{Corollaire} 

\newtheorem{frprop}[frtheorem]{Proposition}
\newtheorem{frdefn}[frtheorem]{D\'efinition}

\theoremstyle{definition}

\newtheorem{taggedtheoremx}{Remark}
\newenvironment{taggedtheorem}[1]
{\renewcommand\thetaggedtheoremx{#1}\taggedtheoremx}
{\endtaggedtheoremx}

\numberwithin{section}{chapter}
\numberwithin{subsection}{section}
\numberwithin{subsubsection}{subsection}

\numberwithin{equation}{chapter}

\renewcommand{\emptyset}{\varnothing}

\newcommand{\pan}{\mathfrak{P}}
\newcommand{\vel}{\mathcal{V}}
\newcommand{\wgotica}{\mathfrak{w}}
\newcommand{\sgotica}{\mathfrak{s}}

\renewcommand{\epsilon}{\varepsilon}
\newcommand{\e}{\varepsilon}

\newcommand{\R}{\ensuremath{\mathbb{R}}}
\newcommand{\Rn}{\ensuremath{\mathbb{R}^n}}
\newcommand{\N}{\ensuremath{\mathbb{N}}}
\newcommand{\eps}{\ensuremath{\varepsilon}}

\newcommand{\Ha}{\mathcal H}
\newcommand{\Ll}{\mathcal L}
\newcommand{\Co}{\mathcal C}
\DeclareMathOperator{\Dim}{Dim}
\newcommand{\Fc}{\mathcal E}
\newcommand{\F}{\mathcal F}
\newcommand{\s}{\mathbb S}
\newcommand{\Op}{\mathcal O}
\newcommand{\Sg}{\mathcal{S}g}
\newcommand{\I}{H}

\newcommand{\Ue}{\overline{u}}
\newcommand{\Uc}{\mathcal U}
\newcommand{\lf}{\mathfrak F}
\newcommand{\Vf}{\mathcal V}
\newcommand{\Gs}{\mathfrak G}
\newcommand{\BaLL}{\mathcal B}
\newcommand{\Cf}{\mathfrak C}
\newcommand{\Qc}{\mathcal Q}
\newcommand{\We}{\overline{w}}
\newcommand{\Pe}{\mathcal P}
\newcommand{\Xe}{\mathcal X}

\newcommand{\C}{\mathscr{C}}
\renewcommand{\L}{\mathscr{L}}
\newcommand{\loc}{{\rm loc}}
\newcommand{\PV}{\mbox{\normalfont P.V.}}
\newcommand{\Haus}{\mathcal{H}}
\newcommand{\red}{{\partial^*}}
\DeclareMathOperator{\Per}{Per}
\newcommand{\cu}{\mathscr H}

\newcommand{\G}{\mathcal G}
\newcommand{\A}{\mathcal A}
\newcommand{\Nl}{\mathcal N}
\newcommand{\K}{\mathcal K}
\DeclareMathOperator{\diam}{diam}
\DeclareMathOperator{\Tail}{Tail}
\newcommand{\h}{\mathscr{H}}
\newcommand{\U}{\mathcal U}
\newcommand{\lra}{\longrightarrow}
\newcommand{\B}{\mathfrak B}
\newcommand{\W}{\mathcal W}
\newcommand{\kers}{|x-y|^{n-1+s}}

\newcommand{\dkers}{\frac{dx\,dy}{|x-y|^{n-1+s}}}
\newcommand{\dKers}{\frac{dX\,dY}{|X-Y|^{n+1+s}}}
\newcommand{\ubar}[1]{\underaccent{\bar}{#1}}

\def\Xint#1{\mathchoice
{\XXint\displaystyle\textstyle{#1}}%
{\XXint\textstyle\scriptstyle{#1}}%
{\XXint\scriptstyle\scriptscriptstyle{#1}}%
{\XXint\scriptscriptstyle\scriptscriptstyle{#1}}%
\!\int}
\def\XXint#1#2#3{{\setbox0=\hbox{$#1{#2#3}{\int}$ }
\vcenter{\hbox{$#2#3$ }}\kern-.6\wd0}}

\def\dashint{\Xint-}

\newlength{\dhatheight}

\newcommand{\lr}[1]{\left( #1 \right)}

\newcommand{\eqlab}[1]{\begin{equation}  \begin{aligned}#1 \end{aligned}\end{equation}} 
\newcommand{\bgs}[1]{\begin{equation*} \begin{aligned}#1\end{aligned}\end{equation*}} 
  \newcommand{\sys}[2][]{\begin{equation*}#1  \left\{\begin{aligned}#2\end{aligned}\right.\end{equation*}}
\newcommand{\alig}[1] {\left\{\begin{aligned}#1 \end{aligned}\right.}

\begin{document}

\dominitoc[n]

\frontmatter

\begin{titlepage}
	
\begin{center}
	
{
	\renewcommand{\arraystretch}{1.6}
	\begin{tabular}{c@{\hskip 10pt}c}
		\includegraphics[width=.15\textwidth]{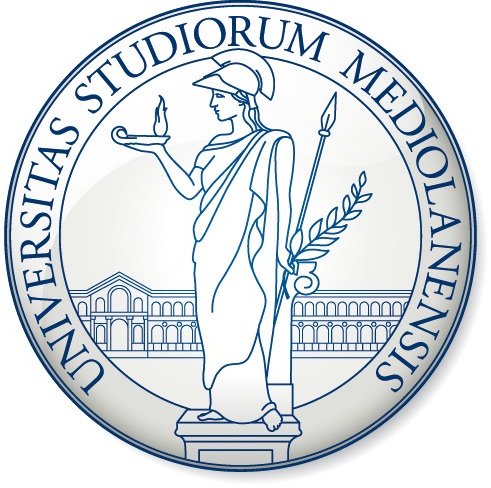} &
		\includegraphics[width=.15\textwidth]{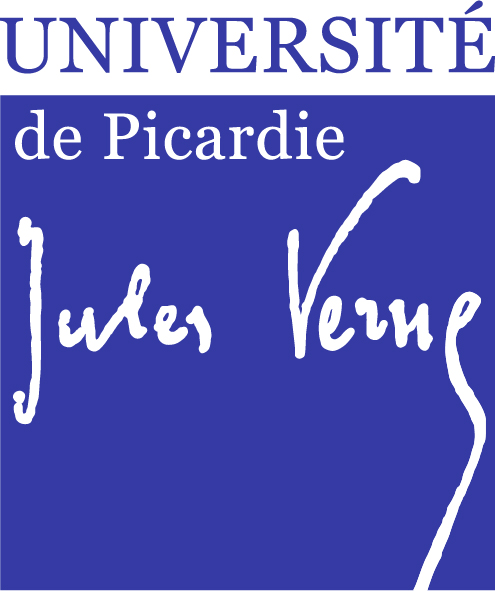} \\
		\renewcommand{\arraystretch}{0.8}
		\begin{tabular}{c}
			\textsc{\fontsize{8}{10}\selectfont Universit\`a degli Studi di Milano} \\
			\fontsize{7}{8}\selectfont Scuola di Dottorato di Ricerca in Matematica
		\end{tabular}
		&
		\renewcommand{\arraystretch}{0.8}
		\begin{tabular}{c}
			\textsc{\fontsize{8}{10}\selectfont Universit\'e de Picardie Jules Verne} \\
			\fontsize{7}{8}\selectfont Ecole Doctorale en Sciences, Technologie et Sant\'e (ED 585)
		\end{tabular}
	\end{tabular}
}

\vspace{24pt}

\textbf{\fontsize{20}{24}\selectfont Tesi di Dottorato in Cotutela}

\vspace{20pt}

DIPARTIMENTO DI MATEMATICA \emph{``FEDERIGO ENRIQUES"}

\vspace{6pt}	

CORSO DI DOTTORATO IN MATEMATICA, CICLO XXXI

\vspace{18pt}

LABORATOIRE AMI\'ENOIS DE MATH\'EMATIQUE FONDAMENTALE ET APPLIQU\'EE (LAMFA)

\vspace{20pt}

\fbox{
	\parbox{0.96\textwidth}{
		\textbf{
			\begin{center}
					\textit{\fontsize{18}{22}\selectfont
						Minimization Problems Involving Nonlocal Functionals: Nonlocal Minimal Surfaces and a Free Boundary Problem}
				\end{center}
			}
	}
}

\vspace{12pt}

	MAT/05

\end{center}

\vspace{18pt}
\begin{flushright}
	DOTTORANDO
	
	\vspace{6pt}
	
	Luca LOMBARDINI
	
\end{flushright}
	
	\vspace{30pt}
	
TUTOR

\vspace{6pt}

Enrico VALDINOCI (Milano)

\vspace{12pt}

Alberto FARINA (Amiens)

\vspace{60pt}

COORDINATORE DEL DOTTORATO (Milano)

\vspace{6pt}

Vieri MASTROPIETRO
	
\begin{center}

\vfill

{\small	(ANNO ACCADEMICO 2017-2018) }


\end{center}

\end{titlepage}


\title{Minimization Problems Involving Nonlocal Functionals: Nonlocal Minimal Surfaces and a Free Boundary Problem}
\author{Luca Lombardini}
\email{luca.lombardini@unimi.it}
\address{Universit\`a degli Studi di Milano, Dipartimento di Matematica, Via Cesare Saldini 50, 20133 Milano, Italy\\
 and Universit\'e de Picardie Jules Verne, 
Facult\'e des Sciences, 33 Rue Saint Leu, 80039 Amiens CEDEX 1, France}
\maketitle

\begin{abstract}
	
This doctoral thesis is devoted to the analysis of some minimization problems that involve nonlocal functionals. We are mainly concerned with the $s$-fractional perimeter and its minimizers, the $s$-minimal sets. We investigate the behavior of sets having (locally) finite fractional perimeter and we establish existence and compactness results for (locally) $s$-minimal sets.
We study the $s$-minimal sets in highly nonlocal regimes, that correspond to small values of the fractional parameter $s$. We introduce a functional framework for studying those $s$-minimal sets that can be globally written as subgraphs. In particular, we prove existence and uniqueness results for minimizers of a fractional version of the classical area functional and we show the equivalence between minimizers and various notions of solution of the fractional mean curvature equation. We also prove a flatness result for entire nonlocal minimal graphs having some partial derivatives bounded from either above or below.\\
Moreover, we consider a free boundary problem, which consists in the minimization of a functional defined as the sum of a nonlocal energy, plus the classical perimeter. Concerning this problem, we prove uniform energy estimates and we study the blow-up sequence of a minimizer---in particular establishing a Weiss-type monotonicity formula.
\end{abstract}

\tableofcontents

\begin{chapter}*{Introduction}
	\adjustmtc

\section{Summary}

This doctoral thesis is devoted to the analysis of some minimization problems that involve nonlocal functionals. Nonlocal operators have attracted an increasing attention in the latest years, both because of their mathematical interest and for their applications---e.g., in modelling anomalous diffusion processes or long-range phase transitions.
We refer the interested reader to \cite{bucval} for an introduction to nonlocal problems.

\smallskip

In this thesis, we are mainly concerned with the $s$-fractional perimeter---which can be considered as a fractional and nonlocal version of the classical perimeter introduced by De Giorgi and Caccioppoli---and its minimizers, the $s$-minimal sets, that were first considered in \cite{CRS10}. The boundaries of these $s$-minimal sets are usually referred to as nonlocal minimal surfaces.
In particular:
\begin{itemize}
	\item we investigate the behavior of sets having (locally) finite fractional perimeter, proving the density of smooth open sets, an optimal asymptotic result for $s\to1^-$, and studying the connection existing between the fractional perimeter and sets having fractal boundaries.
	\item We establish existence and compactness results for minimizers of the fractional perimeter, that extend those proved in \cite{CRS10}.
	\item We study the $s$-minimal sets in highly nonlocal regimes, that correspond to small values of the fractional parameter $s$. We show that, in this case, the minimizers exhibit a behavior completely different from that of their local counterparts---the (classical) minimal surfaces.
	\item We introduce a functional framework for studying those $s$-minimal sets that can be globally written as subgraphs. In particular, we prove existence and uniqueness results for minimizers of a fractional version of the classical area functional and we prove a rearrangement inequality that implies that the subgraphs of these minimizers are minimizing for the fractional perimeter. We refer to the boundaries of such minimizers as nonlocal minimal graphs. We also show the equivalence between minimizers and various notions of solution---namely, weak solutions, viscosity solutions and smooth pointwise solutions---of the fractional mean curvature equation.
	\item We prove a flatness result for entire nonlocal minimal graphs having some partial derivatives bounded from either above or below---thus, in particular, extending to the fractional framework classical theorems due to Bernstein and Moser.
\end{itemize}

\smallskip

We also consider a free boundary problem, which consists in the minimization of a functional defined as the sum of a nonlocal energy, plus the classical perimeter of the interface of separation between the two phases. 
Concerning this problem:
\begin{itemize}
	\item we prove the existence of minimizers and we introduce an equivalent minimization problem which has a ``local nature"---through the extension technique of \cite{CS07}.
	\item We prove uniform energy estimates and we study the blow-up sequence of a minimizer. In particular, we establish a monotonicity formula that implies that blow-up limits are homogeneous.
	\item We investigate the regularity of the free boundary in the case in which the perimeter has a dominant role over the nonlocal energy.
\end{itemize}

\smallskip

We also mention that the last chapter of the thesis consists in a paper that provides a mathematical model which describes the formation of groups of penguins on the shore at sunset. During the occasion of a research trip at the University of Melbourne, we observed the Phillip Island penguin parade and we were so fascinated by the peculiar behavior of the little penguins that we decided to try and describe it mathematically.

\medskip

The thesis is divided into seven chapters, each of which is based on one of the following research articles, that I have written---together with collaborators---during my PhD:
\begin{enumerate}
	\item \emph{Fractional perimeters from a fractal perspective}, published in Advanced Nonlinear Studies---see \cite{Myfractal}. 
	\item \emph{Approximation of sets of finite fractional perimeter by smooth sets and comparison of local and global $s$-minimal surfaces}, published in Interfaces and Free Boundaries---see \cite{mine_cyl_stuff}.
	\item \emph{Complete stickiness of nonlocal minimal surfaces for small values of the fractional parameter}, joint work with C. Bucur and E. Valdinoci, published in Annales de l'Institut Henri Poincar\'e Analyse Non Lin\'eaire---see \cite{BLV16}.
	\item \emph{On nonlocal minimal graphs}, joint work with M. Cozzi, currently in preparation.
	\item \emph{Bernstein-Moser-type results for nonlocal minimal graphs}, joint work with M. Cozzi and A. Farina, currently under submission---see \cite{CFL18}.
	\item A partial, preliminary, version of the article \emph{A free boundary problem: superposition of nonlocal energy plus classical perimeter}, joint work with S. Dipierro and E. Valdinoci, currently in preparation.
	\item \emph{The Phillip Island penguin parade (a mathematical treatment)}, joint work with S. Dipierro, P. Miraglio and E. Valdinoci, published in ANZIAM Journal---see \cite{DMLV16}.
\end{enumerate}

The Appendix contains some auxiliary results that have been exploited throughout the thesis.

\section{A more detailed introduction}

We now proceed to give a detailed description of the contents and main results of this thesis. We observe that each topic has its own, more in-depth, presentation, at the beginning of the corresponding chapter. Moreover, each chapter has its own table of contents, to help the reader navigate through the sections.

\subsection{Sets of (locally) finite fractional perimeter}

The $s$-fractional perimeter and its minimizers, the $s$-minimal sets,
were introduced in \cite{CRS10}, in 2010, mainly motivated by applications to phase transition problems in the presence of long-range interactions.
In the subsequent years, they have attracted a lot of interest,
especially concerning the regularity theory and the qualitative behavior of the boundaries of the $s$-minimal sets, which are the so-called nonlocal minimal surfaces. We refer the interested reader to
\cite{V13} and~\cite[Chapter~6]{bucval}
for an introduction, and to the survey \cite{DV18} for some recent developments.

In particular, we mention that, even if finding the optimal regularity of nonlocal minimal surfaces is still an engaging open problem,
it is known that nonlocal minimal surfaces are $(n-1)$-rectifiable.
More precisely, they are smooth,
except possibly for a singular set of Hausdorff dimension at most equal to $n-3$ (see \cite{CRS10}, \cite{SV13} and \cite{FV17}).
As a consequence, an $s$-minimal set
has (locally) finite perimeter (in the sense of De Giorgi and Caccioppoli)---and actually some uniform estimates for the (classical) perimeter of $s$-minimal sets are available (see \cite{CSV16}).


\smallskip

On the other hand, the boundary of a generic set $E$ having finite $s$-perimeter can be very irregular and indeed it can be ``nowhere rectifiable'', like in the case of the von Koch snowflake.

Actually, the $s$-perimeter can be used (following the seminal paper \cite{Visintin}) to define a ``fractal dimension'' for the measure theoretic boundary
\begin{equation*}
\partial^-E:=\{x\in\R^n\,|\,0<|E\cap B_r(x)|<\omega_nr^n\textrm{ for every }r>0\},
\end{equation*}
of a set $E\subseteq\R^n$.

Before going on, we recall the definition of the $s$-perimeter.
Given a fractional parameter $s\in(0,1)$, we define the interaction
\begin{equation*}
\mathcal L_s(A,B):=\int_A\int_B\frac{1}{|x-y|^{n+s}}\,dx\,dy,
\end{equation*}
for every couple of disjoint sets $A,\,B\subseteq\mathbb R^n$.
Then the \emph{$s$-perimeter} of a set $E\subseteq\mathbb R^n$ in an open set $\Omega\subseteq\R^n$ is defined as
\begin{equation*}
\Per_s(E,\Omega):=\mathcal L_s(E\cap\Omega,\Co E\cap\Omega)+
\mathcal L_s(E\cap\Omega,\Co E\setminus\Omega)+
\mathcal L_s(E\setminus\Omega,\Co E\cap\Omega).
\end{equation*}
We simply write $\Per_s(E):=\Per_s(E,\R^n)$.

We say that a set $E\subseteq\R^n$ has \emph{locally finite $s$-perimeter} in an open set $\Omega\subseteq\R^n$ if
\begin{equation*}
\Per_s(E,\Omega')<\infty\qquad\textrm{for every open set }\Omega'\Subset\Omega.
\end{equation*}

We observe that we can rewrite the $s$-perimeter as
\begin{equation}\label{CH:1:func_per_form}
\Per_s(E,\Omega)=\frac{1}{2}\iint_{\R^{2n}\setminus(\Co\Omega)^2}\frac{|\chi_E(x)-\chi_E(y)|}{|x-y|^{n+s}}\,dx\,dy.
\end{equation}
Formula \eqref{CH:1:func_per_form} shows that the fractional perimeter is, roughly speaking, the $\Omega$-contribution
to the $W^{s,1}$-seminorm of the characteristic function $\chi_E$.

This functional is nonlocal, in that we need to know the set $E$
in the whole of $\R^n$ even to compute its $s$-perimeter in a small bounded domain $\Omega$
(contrary to what happens with the classical perimeter or the $\Ha^{n-1}$ measure, which are local functionals).
Moreover, the $s$-perimeter is ``fractional'', in the sense that the $W^{s,1}$-seminorm measures a fractional order of regularity.

We also observe that we can split the $s$-perimeter as
\bgs{
\Per_s(E,\Omega)=\Per_s^L(E,\Omega)+\Per_s^{NL}(E,\Omega),
}
where
\bgs{
\Per_s^L(E,\Omega):=\Ll_s(E\cap\Omega,\Co E\cap\Omega)=\frac{1}{2}[\chi_E]_{W^{s,1}(\Omega)}
}
can be thought of as the ``local part" of the fractional perimeter, and
\bgs{
\Per_s^{NL}(E,\Omega)&:=\Ll_s(E\cap\Omega,\Co E\setminus\Omega)+
\Ll_s(E\setminus\Omega,\Co E\cap\Omega)\\
&
=\int_\Omega\int_{\Co\Omega}\frac{|\chi_E(x)-\chi_E(y)|}{|x-y|^{n+s}}\,dx\,dy,
}
can be thought of as the ``nonlocal part".

\subsubsection{\textbf{Fractal boundaries}}\label{CH:0:Fractal_sec_Intro}



In 1991, in the paper \cite{Visintin} the author suggested using the index $s$ of the fractional seminorm
$[\chi_E]_{W^{s,1}(\Omega)}$ (and more general continuous families of functionals satisfying appropriate
generalized coarea formulas) as a way to measure the codimension of the measure theoretic boundary
$\partial^-E$ of the set $E$ in $\Omega$.
He proved that the fractal dimension obtained in this way,
\begin{equation*}
\Dim_F(\partial^-E,\Omega):=n-\sup\{s\in(0,1)\,|\,[\chi_E]_{W^{s,1}(\Omega)}<\infty\},
\end{equation*}
is less than or equal to the (upper) Minkowski dimension.

The relationship between the Minkowski dimension of the boundary of $E$
and the fractional regularity (in the sense of Besov spaces)
of the characteristic function $\chi_E$ was investigated also in \cite{Sickel}, in 1999.
In particular, in \cite[Remark 3.10]{Sickel}, the author proved that
the dimension $\Dim_F$ of the von Koch snowflake $S$ coincides with its Minkowski dimension, exploiting the fact that
$S$ is a John domain.

The Sobolev regularity of a characteristic function $\chi_E$ was further studied in \cite{FarRog}, in 2013,
where the authors consider the case in which the set $E$ is a quasiball. Since the von Koch snowflake $S$ is a typical
example of quasiball, the authors were able to prove that the dimension $\Dim_F$ of $S$ coincides with its Minkowski dimension.

\medskip

In Chapter \ref{CH_Fractals} we compute the dimension $\Dim_F$ of the von Koch snowflake $S$ in an elementary way, using only the
roto-translation invariance and the scaling property of the $s$-perimeter and the ``self-similarity'' of $S$.
More precisely, we show that
\bgs{
\Per_s(S)<\infty,\qquad\forall\,s\in\left(
0,2-\frac{\log 4}{\log 3}\right),
}
and
\bgs{
\Per_s(S)=\infty,\qquad\forall\,s\in\left[2-\frac{\log 4}{\log 3},1\right).
}
The proof can be extended in a natural way to all sets which can be defined in a recursive way similar to that of the von Koch snowflake.
As a consequence, we compute the dimension $\Dim_F$ of all such sets,
without having to require them to be John domains or quasiballs.

Furthermore, we show that we can easily obtain a lot of sets of this kind by appropriately modifying well known self-similar fractals like e.g. the von Koch snowflake, the Sierpinski triangle and the Menger sponge.
An example is depicted in Figure \ref{CH:1:buffo_triang}.

\begin{figure}[htbp]
	\begin{center}
		\includegraphics[width=60mm]{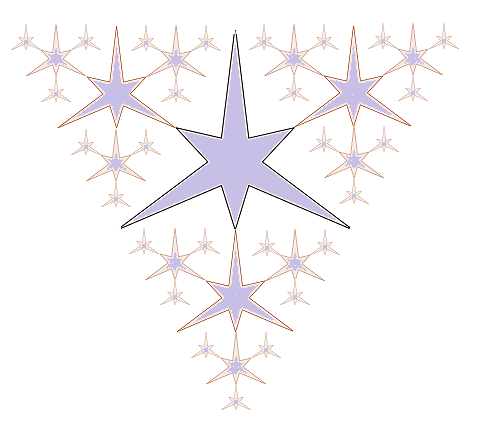}
		\caption{{\it Example of a ``fractal'' set constructed exploiting the structure of the Sierpinski triangle (seen at the fourth iterative step).}}
		\label{CH:1:buffo_triang}
	\end{center}
\end{figure}

\subsubsection{\textbf{Asymptotics $s\to1^-$}}\label{CH:0:asympto1_subsec}

The previous discussion shows that the $s$-perimeter of a set $E$ with an irregular, eventually fractal, boundary
can be finite for $s$ below some threshold, $s<\sigma$, and infinite for $s\in(\sigma,1)$.
On the other hand, it is well known that sets with a regular boundary have finite $s$-perimeter for every $s$ and actually
their $s$-perimeter converges, as $s$ tends to 1, to the classical perimeter,
both in the classical sense (see, e.g.,  \cite{uniform}) and in the $\Gamma$-convergence sense (see, e.g.,  \cite{Gamma}
and also \cite{Ponce} for related results).

In Chapter \ref{CH_Fractals}
we exploit \cite[Theorem 1]{Davila} to prove an optimal version of this asymptotic property for a set $E$ having finite classical
perimeter in a bounded open set with Lipschitz boundary.
More precisely, we prove that if $E$ has finite classical perimeter in a neighborhood of $\Omega$, then
\begin{equation*}
\lim_{s\to1}(1-s)\Per_s(E,\Omega)=\omega_{n-1}\Per(E,\overline\Omega).
\end{equation*}

We observe that we lower the regularity
requested in \cite{uniform}, where the authors required the boundary $\partial E$ to be $C^{1,\alpha}$,
to the optimal regularity (asking $E$ to have only finite perimeter). Moreover,
we do not have to ask $E$ to intersect $\partial\Omega$ ``transversally'', i.e. we do not require
\begin{equation*}
\Ha^{n-1}(\partial^*E\cap\partial\Omega)=0,
\end{equation*}
with $\partial^*E$ denoting the reduced boundary of $E$.

Indeed, we prove that the nonlocal part of the $s$-perimeter converges to
the perimeter on the boundary of $\Omega$, i.e. we prove that
\begin{equation*}
\lim_{s\to1}(1-s)\Per_s^{NL}(E,\Omega)=\omega_{n-1}\Ha^{n-1}(\partial^*E\cap\partial\Omega),
\end{equation*}
which is, to the best of the author's knowledge, a new result.

\subsubsection{\textbf{Approximation by smooth open sets}}\label{CH:0:Appro_subsection}

As we have observed in Section \ref{CH:0:Fractal_sec_Intro}, sets having finite fractional perimeter can have a very rough boundary, which may indeed be a nowhere rectifiable fractal (like the von Koch snowflake).

This represents a dramatic difference between the fractional and the classical perimeter, since Caccioppoli sets have a ``big'' portion of the boundary, the so-called reduced boundary, which is $(n-1)$-rectifiable (by De Giorgi's structure Theorem).

Still, we prove in the first part of Chapter \ref{CH_Appro_Min} that a set has (locally) finite fractional perimeter if and only if it can be approximated (in an appropriate way) by smooth open sets.
More precisely, we prove the following:
\begin{theorem}\label{CH:0:density_smooth_thm}
	Let $\Omega\subseteq\R^n$ be an open set. A set $E\subseteq\R^n$ has locally finite $s$-perimeter in $\Omega$
	if and only if there exists a sequence $E_h\subseteq\R^n$ of open sets with smooth boundary and $\eps_h\to0^+$ such that
	\begin{equation*}\begin{split}
	& (i)\quad E_h\xrightarrow{loc}E,\qquad\sup_{h\in\mathbb N}\Per_s(E_h,\Omega')<\infty\quad\textrm{for every }\Omega'\Subset\Omega,\\
	& (ii)\quad\lim_{h\to\infty}\Per_s(E_h,\Omega')=\Per_s(E,\Omega')\quad\textrm{for every }\Omega'\Subset\Omega,\\
	& (iii)\quad\partial E_h\subseteq N_{\eps_h}(\partial E).
	\end{split}
	\end{equation*}
	Moreover, if $\Omega=\R^n$ and the set $E$ is such that $|E|<\infty$ and $\Per_s(E)<\infty$, then
	\begin{equation*}
	|E_h\Delta E|\to0,\qquad\quad\qquad \lim_{h\to\infty}\Per_s(E_h)=\Per_s(E),
	\end{equation*}
	and we can require each set $E_h$ to be bounded (instead of asking $(iii)$).
\end{theorem}

Here above, $N_\delta(\partial E)$ denotes the tubular $\delta$-neighborhood of $\partial E$.

Such a result is well known for Caccioppoli sets (see, e.g.,  \cite{Maggi}) and indeed this density property can be used to define the (classical) perimeter functional as the relaxation---with respect to $L^1_{\loc}$ convergence---of the $\Ha^{n-1}$ measure of boundaries of smooth open sets, that is
\begin{equation}\label{CH:0:liminfclassical}\begin{split}\Per(E,\Omega)=\inf\Big\{\liminf_{k\to\infty}\Ha^{n-1}(\partial E_h&\cap\Omega)\,\big|\,
E_h\subseteq\R^n\textrm{ open with smooth}\\
&
\textrm{boundary, s.t. }E_h\xrightarrow{loc}E\Big\}.
\end{split}
\end{equation}

It is interesting to observe that in \cite{DV18} the authors have proved, by exploiting the divergence Theorem, that
if $E\subseteq\R^n$ is a bounded open set with smooth boundary, then
\begin{equation}\label{CH:0:liminf_fract_formula}
\Per_s(E)=c_{n,s}\int_{\partial E}\int_{\partial E}
\frac{2-|\nu_E(x)-\nu_E(y)|^2}{|x-y|^{n+s-2}}d\Ha^{n-1}_xd\Ha^{n-1}_y,
\end{equation}
where $\nu_E$ denotes the external normal of $E$ and
\begin{equation*}
c_{n,s}:=\frac{1}{2s(n+s-2)}.
\end{equation*}
By exploiting equality \eqref{CH:0:liminf_fract_formula}, the lower semicontinuity of the $s$-perimeter
and Theorem \ref{CH:0:density_smooth_thm}, we find that, if
$E\subseteq\R^n$ is such that $|E|<\infty$, then
\begin{equation*}\begin{split}
\Per_s(E)=\inf\bigg\{&\liminf_{h\to\infty}c_{n,s}\int_{\partial E_h}\int_{\partial E_h}
\frac{2-|\nu_{E_h}(x)-\nu_{E_h}(y)|^2}{|x-y|^{n+s-2}}d\Ha^{n-1}_xd\Ha^{n-1}_y\,\big|\\
&\quad
E_h\subseteq\R^n\textrm{ bounded open set with smooth
	boundary, s.t. }E_h\xrightarrow{loc}E\bigg\}.
\end{split}
\end{equation*}
This can be thought of
as an analogue of \eqref{CH:0:liminfclassical} in the fractional setting.

We also mention that in Section \ref{CH:4:Appro_Section} we will prove that a subgraph having locally finite $s$-perimeter in a cylinder $\Omega\times\R$ can be approximated by the subgraphs of smooth functions---and not just by arbitrary smooth open sets.

\subsection{Nonlocal minimal surfaces}\label{CH:0:NMS_Sec}

The second part of Chapter \ref{CH_Appro_Min} is concerned with sets minimizing the fractional perimeter.
The boundaries of these minimizers are often referred to as nonlocal minimal surfaces and naturally arise as
limit interfaces of long-range interaction phase transition models. In particular, in regimes where
the long-range interaction is dominant, the nonlocal Allen-Cahn energy functional $\Gamma$-converges
to the fractional perimeter (see, e.g., \cite{SV12})
and the minimal interfaces of the corresponding
Allen-Cahn equation approach locally uniformly
the nonlocal minimal surfaces (see, e.g., \cite{SV14}).

We now recall the definition of minimizing sets introduced in \cite{CRS10}.
\begin{defn}\label{CH:0:sMini_def}
Let $\Omega\subseteq\R^n$ be an open set and let $s\in(0,1)$. We say that a set $E\subseteq\R^n$ is \emph{$s$-minimal} in $\Omega$ if $\Per_s(E,\Omega)<\infty$ and
\[
\Per_s(E,\Omega)\leq\Per_s(F,\Omega)\quad\mbox{ for every }F\subseteq\R^n\mbox{ s.t. }F\setminus\Omega=E\setminus\Omega.
\]
\end{defn}
Among the many results, in \cite{CRS10} the authors have proved that,
if $\Omega\subseteq\R^n$ is a bounded open set with Lipschitz boundary, then for every fixed set $E_0\subseteq\Co\Omega$ there exists a set
$E\subseteq\R^n$ which is $s$-minimal in $\Omega$ and such that $E\setminus\Omega=E_0$.
The set $E_0$ is sometimes referred to as \emph{exterior data} and the set $E$ is said to be $s$-minimal in $\Omega$ with respect to the exterior data $E_0$.


We extend the aforementioned existence result, by proving that, in a generic open set $\Omega$, there exists an $s$-minimal set with respect to some fixed exterior data $E_0\subseteq\Co\Omega$ if and only if there exists a competitor having finite $s$-perimeter in $\Omega$. More precisely:

\begin{theorem}\label{CH:0:Global_Existence}
	Let $s\in(0,1)$, let $\Omega\subseteq\R^n$ be an open set and let $E_0\subseteq\Co\Omega$. Then, there exists a set $E\subseteq\R^n$ which is $s$-minimal in $\Omega$, with $E\setminus\Omega=E_0$, if and only if there exists a set $F\subseteq\R^n$ such that $F\setminus\Omega=E_0$ and $\Per_s(F,\Omega)<\infty$.
\end{theorem}

As a consequence, we observe that if $\Per_s(\Omega)<\infty$, then there always exists an $s$-minimal set with respect to the exterior data $E_0$, for every set $E_0\subseteq\Co\Omega$.

\medskip

Let us now turn the attention to the case in which the domain of minimization is not bounded. In this situation, it is convenient to introduce the notion of local minimizer.
\begin{defn}
Let $\Omega\subseteq\R^n$ be an open set and let $s\in(0,1)$. We say that a set $E\subseteq\R^n$ is \emph{locally $s$-minimal} in $\Omega$ if $E$ is $s$-minimal in every open set $\Omega'\Subset\Omega$.
\end{defn}

Notice in particular that we are only requiring $E$ to be of locally finite $s$-perimeter in $\Omega$ and not to have finite $s$-perimeter in the whole domain.
Indeed, the main reason
for the introduction of locally $s$-minimal sets is given by the fact that, in general, the $s$-perimeter of a set is not finite in unbounded domains.

We have seen in Theorem \ref{CH:0:Global_Existence} that the only obstacle to the existence of an $s$-minimal set, with respect to some fixed exterior data $E_0\subseteq\Co\Omega$, is the existence of a competitor having finite $s$-perimeter. On the other hand, we prove that a locally $s$-minimal set always exists, no matter what the domain $\Omega$ and the exterior data are.

\begin{theorem}\label{CH:0:Local_Existence}
	Let $s\in(0,1)$, let $\Omega\subseteq\R^n$ be an open set and let $E_0\subseteq\Co\Omega$. Then, there exists a set $E\subseteq\R^n$ which is locally $s$-minimal in $\Omega$, with $E\setminus\Omega=E_0$.
\end{theorem}
When $\Omega$ is a bounded open set with Lipschitz boundary, we show that the two notions of minimizer coincide. That is, if $\Omega\subseteq\R^n$ is a bounded open set with Lipschitz boundary and $E\subseteq\R^n$, then
\[
E\mbox{ is $s$-minimal in }\Omega\quad\Longleftrightarrow\quad E\mbox{ is locally $s$-minimal in }\Omega.
\]

However, we observe that this is not true in an arbitrary open set $\Omega$, since an $s$-minimal set---in the sense of Definition \ref{CH:0:sMini_def}---may not exist.

As an example, we consider the situation in which the domain of minimization is
the cylinder
\[
\Omega^\infty:=\Omega\times\R\subseteq\R^{n+1},
\]
with $\Omega\subseteq\R^n$ a bounded open set with regular boundary. We are interested in exterior data given by the subgraph of some measurable function $\varphi:\R^n\to\R$. That is, we consider
the subgraph
\[
\Sg(\varphi):=\left\{(x,x_{n+1})\in\R^{n+1}\,|\,x_{n+1}<\varphi(x)\right\},
\]
and we want to find a set $E\subseteq\R^{n+1}$ that minimizes---in some sense---the $s$-perimeter in the cylinder $\Omega^\infty$, with respect to the exterior data $E\setminus\Omega^\infty=\Sg(\varphi)\setminus\Omega^\infty$.

A motivation for considering such a minimization problem is given by the recent article \cite{graph}, where the authors have proved that if such a minimizing set $E$ exists---and if $\varphi$ is a continuous function---then $E$ is actually a global subgraph. More precisely, there exists a function $u:\R^n\to\R$, with $u=\varphi$ in $\R^n\setminus\overline{\Omega}$ and $u\in C(\overline{\Omega})$ such that
\[
E=\Sg(u).
\]

It is readily seen that if a function $u:\R^n\to\R$ is well behaved in $\Omega$, e.g., if $u\in BV(\Omega)\cap L^\infty(\Omega)$, then the local part of the $s$-perimeter of the subgraph of $u$ is finite,
\[
\Per_s^L(\Sg(u),\Omega^\infty)<\infty.
\]
On the other hand, the nonlocal part of the $s$-perimeter, in general, is infinite, even for very regular functions $u$. Indeed, we prove that if $u\in L^\infty(\R^n)$, then
\[
\Per_s^{NL}(\Sg(u),\Omega^\infty)=\infty.
\]

A first consequence of this observation---and of the apriori bound on the ``vertical variation" of a minimizing set provided by \cite[Lemma 3.3]{graph}---is the fact that, if $\varphi\in C(\R^n)\cap L^\infty(\R^n)$, then there can not exist a set $E$ which is $s$-minimal in $\Omega^\infty$---in the sense of Definition \ref{CH:0:sMini_def}---with respect to the exterior data $\Sg(\varphi)\setminus\Omega^\infty$.

Nevertheless, Theorem \ref{CH:0:Local_Existence} guarantees the existence of a set $E\subseteq\R^{n+1}$ that is locally $s$-minimal in $\Omega^\infty$ and such that $E\setminus\Omega^\infty=\Sg(\varphi)\setminus\Omega^\infty$. Therefore, Theorem \ref{CH:0:Local_Existence} and \cite[Theorem 1.1]{graph} together imply the existence of subgraphs (locally) minimizing the $s$-perimeter, that is, namely, nonparametric nonlocal minimal surfaces.

A second consequence consists in the fact that we can not define a naive fractional version of the classical area functional as
\[
\mathscr{A}_s(u,\Omega):=\Per_s(\Sg(u),\Omega^\infty),
\]
since this would be infinite even for a function $u\in C^\infty_c(\R^n)$. In Chapter \ref{CH_Nonparametric} we will get around this issue by introducing an appropriate functional setting for working with subgraphs.

\subsection{Stickiness effects for small values of $s$}

Chapter \ref{Asympto0_CH_label} is devoted to the study of $s$-minimal sets in highly nonlocal regimes, i.e. in the case in which the fractional parameter $s\in(0,1)$ is very small. We prove that in this situation the behavior of $s$-minimal sets, in some sense, degenerates.

\medskip

Let us first recall some known results concerning the asymptotics as $s\to1^-$.\\
We have already observed in Section \ref{CH:0:asympto1_subsec} that the $s$-perimeter converges to the classical perimeter as $s\to1^-$.
Moreover, as $s\to1^-$, $s$-minimal sets converge to minimizers of the classical perimeter, both in a ``uniform sense'' (see \cite{uniform, regularity}) and in the
$\Gamma$-convergence sense (see \cite{Gamma}).
As a consequence, one is able to prove (see \cite{regularity}) that for $s$ sufficiently close to 1, nonlocal minimal surfaces have the same regularity of classical minimal surfaces.
See also \cite{DV18} for a recent and quite comprehensive survey of the properties of $s$-minimal sets when $s$ is close to 1. 

Furthermore, we observe that also the fractional mean curvature converges, as $s\to1^-$, to its classical counterpart. To be more precise, let us first recall that the $s$-fractional mean curvature of a set $E$ at a point $q\in\partial E$ is defined as the principal value integral
\[\I_s[E](q):=\PV\int_{\R^n}\frac{\chi_{\Co E}(y)-\chi_E(y)}{|y-q|^{n+s}}\,dy,\]
that is
\[\I_s[E](q):=\lim_{\varrho\to0^+}\I_s^\varrho[E](q),\qquad\textrm{where}\qquad
\I_s^\varrho[E](q):=\int_{
	\Co B_\varrho(q)}\frac{\chi_{\Co E}(y)-\chi_E(y)}{|y-q|^{n+s}}\,dy.\] 
Let us remark that it is indeed necessary to interpret the above integral in the principal value sense, since the integrand is singular and not integrable in a neighborhood of $q$. On the other hand, if there is enough cancellation between $E$ and $\Co E$ in a neighborhood of $q$---e.g., if $\partial E$ is of class $C^2$ around $q$---then the integral is well defined in the principal value sense.

The fractional mean curvature was introduced in \cite{CRS10}, where the authors proved that it is the Euler-Lagrange operator appearing in the minimization of the $s$-perimeter. Indeed,
if $E\subseteq\R^n$ is $s$-minimal in an open set $\Omega$, then 
\[
\I_s[E]=0\quad\mbox{on } \partial E,
\]
in an appropriate viscosity sense---for more details see, e.g., Appendix \ref{CH:3:brr2}.

It is known (see, e.g., \cite[Theorem 12]{Abaty} and \cite{regularity}) that if $E\subseteq \Rn$ is a set with $C^2$ boundary, and $n\ge2$, then for any $x\in \partial E$ one has that
\[ \lim_{s \to 1} (1-s)\I_s[E] (x) = \varpi_{n-1}H[E](x).\]
Here above $H$ denotes the classical mean curvature of $E$ at the point $x$---with the convention that we take $H$ such that the curvature of the ball is a positive quantity---and
\[
\varpi_k:=\Ha^{k-1}(\{x\in\R^k\,|\,|x|=1\}),
\]
for every $k\ge1$. Let us also define $\varpi_0:=0$. We observe that for $n=1$, we have that
\[  \lim_{s \to 1} (1-s)\I_s[E] (x) = 0,\]
which is consistent with the notation $\varpi_0=0$---see also Remark \ref{CH:3:nuno}.

\medskip

As $s\to 0^+$, the asymptotics are more involved and present some surprising behavior.  This is due to the fact that as $s$ gets smaller, the nonlocal contribution to the $s$-perimeter becomes more and more important, while the local contribution loses influence. Some precise
results in this sense were achieved in \cite{DFPV13}. There, in order to encode the behavior at infinity of a set, the authors have introduced the quantity
\[
\alpha(E)=\lim_{s\to 0^+} s\int_{\Co B_1} \frac{\chi_E(y)}{|y|^{n+s}}\, dy,
\]
which appears naturally when looking at the
asymptotics as $s\to0^+$ of the fractional perimeter. Indeed, in \cite{DFPV13} the authors proved that, if $\Omega$ is a bounded open set with $C^{1,\gamma}$ boundary, for some $\gamma\in (0,1]$, $E \subseteq \Rn$ has finite $s_0$-perimeter in $\Omega$, for some $s_0\in (0,1)$, and $\alpha(E)$ exists, then
\bgs{
\lim_{s\to 0^+} s\Per_s(E,\Omega)=  \alpha(\Co E) |E\cap \Omega| + \alpha(E) |\Co E \cap \Omega|.
}

On the other hand, the asymptotic behavior for $s\to 0^+$ of the fractional mean curvature is studied in Chapter \ref{Asympto0_CH_label}
(see also \cite{DV18} for the particular case in which the set $E$ is bounded).
First of all, since the quantity $\alpha(E)$ may not exist---see \cite[Example 2.8 and 2.9]{DFPV13}---we define
\[
\overline \alpha (E):= \limsup_{s\to 0^+} s\int_{\Co B_1} \frac{\chi_E(y)}{|y|^{n+s}}\, dy
\quad\mbox{and} \quad
\underline \alpha(E) := \liminf_{s\to 0^+} s\int_{\Co B_1} \frac{\chi_E(y)}{|y|^{n+s}}\, dy.
\]
We prove that, when $s\to0^+$, the $s$-fractional mean curvature  becomes completely indifferent to the local geometry of the boundary $\partial E$, and indeed the limit value only depends on the behavior at infinity of the set $E$. More precisely, if $E\subseteq\Rn$ and $p\in\partial E$ is such that $\partial E$ is $C^{1,\gamma}$ near $p$,
for some $\gamma\in(0,1]$, then
\eqlab{\label{CH:0:liminf_curvat}
\liminf_{s\to0^+} s\,\I_s[E](p) =\varpi_n -2 \overline \alpha(E),
}
and
\[
	\limsup_{s\to0^+} s\,\I_s[E](p) =\varpi_n-2 \underline\alpha(E).
\]
We remark in particular that if $E$ is bounded, then $\alpha(E)$ exists and $\alpha(E)=0$. Hence, if $E\subseteq\R^n$ is a bounded open set with $C^{1,\gamma}$ boundary, the asymptotics is simply
\[
\lim_{s\to0^+}s\,\I_s[E](p)=\varpi_n,
\]
for every $p\in\partial E$---see also \cite[Appendix~B]{DV18}.

In Section \ref{CH:3:sectexamples} we compute the contribution from infinity $\alpha(E)$ of some sets. To have a few examples in mind, we mention here the following cases:
\begin{itemize}
	\item let $S\subseteq\s^{n-1}$ and consider the cone
	\[
	C:=\{t\sigma\in\R^n\,|\,t\geq0,\,\sigma\in S\}.
	\]
	Then, $\alpha(C)=\Ha^{n-1}(S)$.
	\item If $u\in L^\infty(\R^n)$, then
	$\alpha(\Sg(u))=\varpi_{n+1}/2$. More in general, if $u:\R^n\to\R$ is such that
	\[
	\lim_{|x|\to\infty}\frac{|u(x)|}{|x|}=0,
	\]
	then $\alpha(\Sg(u))=\varpi_{n+1}/2$.
	\item Let $u:\R^n\to\R$ be such that $u(x)\leq-|x|^2$, for every $x\in\R^n\setminus B_R$, for some $R>0$. Then $\alpha(\Sg(u))=0$.
\end{itemize}
Roughly speaking, from the above examples we see that $\alpha(E)$ does not depend on the local geometry or regularity of $E$, but only on its behavior at infinity.

\smallskip

Now we observe that, as $s \to 0^+$, $s$-minimal sets exhibit a rather unexpected behavior.

For instance, in \cite[Theorem 1.3]{boundary} it is proved that if we fix the first quadrant of the plane
as exterior data, then, quite surprisingly, when $s$ is small enough the $s$-minimal set in $B_1\subseteq\R^2$
is empty in $B_1$.
The main results of Chapter \ref{Asympto0_CH_label} take their inspiration from this result.


Heuristically, in order to generalize
\cite[Theorem 1.3]{boundary} we want to prove that, if $\Omega\subseteq\R^n$ is a bounded and connected open set with smooth boundary and if we fix as exterior data
a set $E_0\subseteq\Co\Omega$ such that $\overline{\alpha}(E_0)<\varpi_n/2$, then there is a contradiction between the Euler-Lagrange equation of an $s$-minimal set and the asymptotics
of the $s$-fractional mean curvature as $s\to0^+$.

To motivate why we expect such a contradiction, we observe that
the asymptotics \eqref{CH:0:liminf_curvat} seems to suggest that, if $s$ is small enough, then an $s$-minimal set $E$ having exterior data $E_0$ and such that $\partial E\cap\Omega\not=\emptyset$ should have some point $p\in\partial E\cap\Omega$ such that $\I_s[E](p)>0$---which would contradict the Euler-Lagrange equation. To avoid such a contradiction, we would then conclude that $\partial E=\emptyset$ in $\Omega$, meaning that either $E\cap\Omega=\Omega$ or $E\cap\Omega=\emptyset$. 

In order to turn this idea into a rigorous argument, we first prove that we can estimate the fractional mean curvature from below uniformly with respect to the radius of an exterior
tangent ball to $E$.
More precisely:
\begin{theorem}\label{CH:0:positive_curvature}
	Let $\Omega\subseteq\Rn$ be a bounded open set. 
	Let $E_0\subseteq\Co\Omega$ be such that
	\bgs{
		\overline \alpha(E_0)<\frac{\varpi_n}2,
	}
	and let
	\[
	\beta=\beta(E_0):=\frac{\varpi_n-2\overline \alpha(E_0)}4.
	\] We define 
	\bgs{
		\delta_s=\delta_s(E_0):=e^{-\frac{1}{s}\log \frac{\varpi_n+2\beta}{\varpi_n+\beta}},
	}
	for every $s\in(0,1)$.
	Then, there exists $s_0=s_0(E_0,\Omega)\in(0,\frac{1}{2}]$ such that, if $E\subseteq\Rn$ is such that $E\setminus\Omega=E_0$
	and $E$ has an exterior tangent ball
	of radius (at least) $\delta_\sigma$, for some $\sigma\in(0,s_0)$, at some point $q\in\partial E\cap\overline{\Omega}$, then
	\bgs{
		\liminf_{\varrho\to0^+}\I_s^\varrho[E](q)\geq\frac{\beta}{s}>0,\qquad\forall\,s\in(0,\sigma].
	}
\end{theorem}

Let us now introduce the following definition.
\begin{defn}\label{CH:0:deltadense}
	Let $\Omega\subseteq \Rn$ be a bounded open set.
	We say that a set $E$ is \emph{$\delta$-dense} in $\Omega$, for some fixed $\delta>0$, if $|B_\delta(x)\cap E|>0$ for any $x\in \Omega$ for which $B_\delta(x)\Subset\Omega$.
\end{defn}

By exploiting a careful geometric argument and Theorem \ref{CH:0:positive_curvature},
we can then pursue the heuristic idea outlined above and prove the following classification result:

\begin{theorem}\label{CH:0:Main_THM}
	Let $\Omega\subseteq\R^n$ be a bounded and connected open set with $C^2$ boundary. Let $E_0\subseteq \Co \Omega$ such that
	\[
	\overline{\alpha}(E_0)<\frac{\varpi_n}{2}.
	\]  
	Then, the following two results hold true.\\
	A)  Let $s_0$ and $\delta_s$ be as in Theorem \ref{CH:0:positive_curvature}. There exists
	$s_1=s_1(E_0,\Omega)\in (0,s_0]$ such that if $s<s_1$ and $E$ is an $s$-minimal set in $\Omega$ with exterior data $E_0$, then either
	\bgs{
		(A.1) \;  E\cap \Omega=\emptyset \quad  \mbox{ or} \quad\; (A.2)\;  E \mbox{ is } \delta_s-\mbox{dense in }\Omega.
	}
\noindent
	B) Either \\
	(B.1) there exists
	$\tilde s=\tilde s(E_0,\Omega)\in (0,1)$ such that if $E$ is an $s$-minimal set in $\Omega$ with exterior data $E_0$ and $s\in(0,\tilde s)$, then
	\bgs{
		E\cap \Omega=\emptyset,
	}
	or \\
	(B.2)    there exist  $\delta_k \searrow 0$, $s_k \searrow 0$ and a sequence of sets  $E_k$ such that each $E_k$ is $s_k$-minimal in $\Omega$ with exterior data $E_0$ and for every $k$
	\bgs{
		\partial E_k \cap B_{\delta_k}(x) \neq \emptyset \quad \mbox{for every } B_{\delta_k}(x)\Subset \Omega.
	}
\end{theorem}

Roughly speaking, either the $s$-minimal sets are empty in $\Omega$ when $s$ is small enough, or we can find a sequence $E_k$ of $s_k$-minimal sets, with $s_k\searrow0$, whose boundaries tend to (topologically) fill the domain $\Omega$ in the limit $k\to\infty$.

We point out that the typical behavior consists in being empty. Indeed, if the exterior data $E_0\subseteq\Co\Omega$ does not completely surround the domain $\Omega$, we have the following result:

 \begin{theorem}\label{CH:0:Not_Surrounded}
	Let $\Omega$ be a  bounded and  connected open set with $C^2$ boundary. Let $E_0\subseteq \Co \Omega$ such that 
	\[
	\overline{\alpha}(E_0)<\frac{\varpi_n}{2},
	\]
	and let $s_1$ be as in Theorem \ref{CH:0:Main_THM}. Suppose that there exists $R>0$ and $x_0\in \partial \Omega$ such that
	\[
	B_R(x_0)\setminus \Omega \subseteq \Co E_0.
	\]
	Then, there exists $s_3=s_3(E_0,\Omega)\in(0,s_1]$ such that if $s<s_3$ and $E$ is an $s$-minimal set in $\Omega$ with exterior data $E_0$, then 
	\[  E\cap \Omega=\emptyset .\]
\end{theorem}

We observe that the condition $\overline{\alpha}(E_0)<\varpi_n/2$
is somehow optimal. Indeed,
when $\alpha(E_0)$ exists and 
\[ \alpha(E_0)=\frac{\varpi_n}2,\]
several configurations may occur, depending on the position of $\Omega$ with respect to the exterior data $E_0\setminus \Omega$---we provide various examples in
Chapter \ref{Asympto0_CH_label}.

Moreover, notice that when $E$ is $s$-minimal in $\Omega$ with respect to $E_0$, then $\Co E$ is $s$-minimal in $\Omega$ with respect to $\Co E_0$. Also,
\[
\underline \alpha(E_0) >\frac{\varpi_n}{2} \qquad \implies \qquad \overline \alpha (\Co E_0)< \frac{\varpi_n}{2}.
\]
Thus, in this case we can apply Theorems \ref{CH:0:positive_curvature}, \ref{CH:0:Main_THM} and \ref{CH:0:Not_Surrounded} to $\Co E$ with respect to the exterior data $\Co E_0$. For instance, if
$E$ is $s$-minimal in $\Omega$ with exterior data $E_0$ with
\[ \underline \alpha(E_0) >\frac{\varpi_n}{2}, \]
and $s<s_1(\Co E_0, \Omega)$,
then either
\[ E\cap \Omega=\Omega \qquad \mbox{ or }  \qquad  \Co E \; \mbox{ is } \; \delta_s(\Co E_0)-\mbox{dense}.\]
The analogues of the just mentioned Theorems can be obtained similarly.

Therefore, from our main results and the above observations, we have a complete classification of nonlocal minimal surfaces when $s$ is small, whenever
\[ \alpha(E_0)\neq  \frac{\varpi_n}{2} .\]

\medskip

We point out that the stickiness phenomena described in \cite{boundary}
and in Chapter \ref{Asympto0_CH_label} are specific for nonlocal minimal surfaces---since
classical minimal surfaces cross transversally the boundary of a 
convex domain. 

Interestingly, these stickiness phenomena
are not present in the case of the fractional Laplacian,
where the boundary datum of the Dirichlet problem is attained continuously under rather
general assumptions, see \cite{MR3168912},
though solutions of $s$-Laplace
equations are in general no better than $C^s$ at the boundary, hence
the uniform continuity degenerates as~$s\to0^+$.

On the other hand, in case of fractional harmonic functions, a partial counterpart of
the stickiness phenomenon is, in a sense, given by the boundary explosive
solutions constructed in \cite{MR3393247,MR2985500}
(namely, in this case, the boundary of the subgraph of the fractional harmonic function
contains vertical walls).

We also mention that stickiness phenomena for nonlocal minimal graphs---eventually in the presence of obstacles---will be studied
in the forthcoming article \cite{LuCla}.

\medskip

In the final part of Chapter \ref{Asympto0_CH_label} we prove that the fractional mean curvature is continuous with respect to all variables.

To simplify a little the situation, suppose that $E_k,\,E\subseteq\R^n$ are sets with $C^{1,\gamma}$ boundaries, for some $\gamma\in(0,1]$, such that the boundaries $\partial E_k$ locally converge in the $C^{1,\gamma}$ sense to the boundary of $E$, as $k\to\infty$.
Then we prove that, if we have a sequence of points $x_k\in\partial E_k$ such that $x_k\to x\in\partial E$ and a sequence of indexes $s_k,\,s\in(0,\gamma)$ such that $s_k\to s$, it holds
\[
\lim_{k\to\infty}\I_{s_k}[E_k](x_k)=\I_s[E](x).
\]
Furthermore, we appropriately extend this convergence result in order to cover also the cases in which $s_k\to1$ or $s_k\to0$.

In particular, let us consider a set $E\subseteq\R^n$ such that $\alpha(E)$ exists and $\partial E$ is of class $C^2$. Then, if we define
\sys[
\tilde \I_s  {[}E{]} (x):=]{ &s(1-s)\I_s[E](x),  & \mbox{ for } &s\in (0,1) \\
&{\varpi_{n-1}} H[E](x), &\mbox{ for } &s=1\\
&\varpi_n-2\alpha(E), &\mbox{ for }
&s=0,}
the function
\[
\tilde\I_{(\,\cdot\,)}[E](\,\cdot\,)
:[0,1]\times\partial E\longrightarrow\R,
\qquad(s,x)\longmapsto \tilde \I_s[E](x),
\]
is continuous.
It is interesting to observe that the fractional mean curvature at a fixed point $q\in\partial E$ may change sign as $s$ varies from 0 to 1. Also---as a consequence of the continuity in the fractional parameter $s$---in such a case there exists an index $\sigma\in(0,1)$ such that $\I_\sigma[E](q)=0$.

\subsection{Nonparametric setting}

In Chapter \ref{CH_Nonparametric} we introduce a functional framework to study minimizers of the fractional perimeter which can be globally written as the subgraph
\[
\Sg(u)=\left\{(x,x_{n+1})\in\R^{n+1}\,|\,x_{n+1}<u(x)\right\},
\]
of some measurable function $u:\R^n\to\R$. We refer to the boundaries of such minimizers as \emph{nonlocal minimal graphs}.

We define a fractional version of the classical area functional and we study its functional and geometric properties.
Then we focus on minimizers and we prove existence and uniqueness results with respect to a large class of exterior data,
which includes locally bounded functions.

Furthermore, one of the main contributions
of Chapter \ref{CH_Nonparametric} consists in proving the equivalence of:
\begin{itemize}
	\item minimizers of the fractional area functional,
	\item minimizers of the fractional perimeter,
	\item weak solutions of the fractional mean curvature equation,
	\item viscosity solutions of the fractional mean curvature equation,
	\item smooth functions solving pointwise the fractional mean curvature equation.
\end{itemize}

\medskip

Before giving a detailed overview of the main results,
let us recall the definition of the classical area functional.
Given a bounded open set $\Omega\subseteq\R^n$ with Lipschitz boundary, the area functional is defined as
\[
\mathscr A(u,\Omega):=\int_\Omega\sqrt{1+|\nabla u|^2}\,dx=\Ha^n\left(\left\{(x,u(x))\in\R^{n+1}\,|\,x\in\Omega\right\}\right),
\]
for every Lipschitz function $u:\overline{\Omega}\to\R$. One then extends this functional, by defining the relaxed area functional of a function $u\in L^1(\Omega)$ as
\[
\mathscr A(u,\Omega):=\inf\left\{\liminf_{k\to\infty}\mathscr A(u_k,\Omega)\,|\,u_k\in C^1(\overline{\Omega}),\,\|u-u_k\|_{L^1(\Omega)}\to0\right\}.
\]
It is readily seen that, if $u\in L^1(\Omega)$, then
\eqlab{\label{CH:0:obs1_locarea}
\mathscr A(u,\Omega)<\infty\quad\Longleftrightarrow\quad u\in BV(\Omega),
}
in which case
\eqlab{\label{CH:0:obs2_locarea}
\mathscr A(u,\Omega)=\Per\left(\Sg(u),\Omega\times\R\right).
}

Roughly speaking, the functions of bounded variation are precisely those integrable functions whose subgraphs have finite perimeter---for the details see, e.g., \cite{Giusti,GiaMart12}.

We could thus be tempted to try and define a fractional version of the area functional, by considering the $s$-perimeter in place of the classical perimeter, setting, for a measurable function $u:\R^n\to\R$,
\[
\mathscr A_s(u,\Omega):=\Per_s(\Sg(u),\Omega\times\R).
\]
However, as we observed in the end of Section \ref{CH:0:NMS_Sec}, such a definition can not work, because
\[
\Per_s^{NL}(\Sg(u),\Omega\times\R)=\infty,
\]
even if $u\in C^\infty_c(\R^n)$.

Before going on, a couple of observations are in order.
Even if the nonlocal part of the fractional perimeter in the cylinder $\Omega^\infty:=\Omega\times\R$ is infinite, we recall that we know---see the end of Section \ref{CH:0:NMS_Sec}---that the local part is finite, provided the function $u$ is regular enough in $\Omega$.

If the function $u$ is bounded in $\Omega$, then we can consider the fractional perimeter in the ``truncated cylinder" $\Omega^M:=\Omega\times(-M,M)$, with $M\geq\|u\|_{L^\infty(\Omega)}$, instead of in the whole cylinder $\Omega^\infty$. As we will see below, by pursuing this idea we obtain a family of fractional area functionals $\F^M_s(\,\cdot\,,\Omega)$.

\smallskip

On the other hand, there is another possibility to come up with a definition of a fractional area functional.
In \cite{regularity}, the authors have observed that when~$E\subseteq\R^{n+1}$ is the subgraph of a function~$u$,
its fractional mean curvature can be written as an integrodifferential operator acting on~$u$. More precisely, letting~$u: \R^n \to \R$ be a function of, say, class~$C^{1, 1}$ in a neighborhood of a point~$x \in \R^n$, we have that
\[
H_s[\Sg(u)](x,u(x)) = \h_s u(x),
\]
with
\[
\h_s u(x) := 2 \, \PV \int_{\R^n} G_s \left( \frac{u(x)-u(y)}{|x - y|} \right) \frac{dy}{|x - y|^{n+s}},
\]
and
\[
G_s(t):=\int_0^t g_s(\tau)\,d\tau,\qquad g_s(t):=\frac{1}{(1+t^2)^\frac{n+1+s}{2}} \quad \mbox{for } t\in\R.
\]
We now show that $\h_s$ is the Euler-Lagrange operator associated to a (convex) functional $\F_s(\,\cdot\,,\Omega)$, which we will then consider as the $s$-fractional area functional.

Let us begin by remarking that, when~$u$ is not regular enough around~$x$, the quantity~$\h_s u(x)$ is in general not well-defined, due to the lack of cancellation required for the principal value to converge. Nevertheless,
we can understand the operator~$\h_s$ as defined in the following weak (distributional) sense. Given a measurable function~$u:\R^n \to \R$, we set
\[
	\langle \h_s u, v\rangle:=\int_{\Rn}\int_{\Rn} G_s \left( \frac{u(x)-u(y)}{|x-y|} \right) \big( v(x)-v(y) \big) \, \frac{dx\, dy}{|x-y|^{n+s}}
\]
for every $v\in C^\infty_c(\R^n)$. More generally, it is immediate to see--- by taking advantage of the fact that~$G_s$ is bounded---that this definition is well-posed for every~$v\in W^{s,1}(\Rn)$. Indeed, one has that
\[
|\langle \h_s u ,v \rangle|\le\frac{\Lambda_{n,s}}{2} \, [v]_{W^{s,1}(\Rn)},
\]
where
\[
\Lambda_{n,s}:=\int_\R g_s(t)\,dt<\infty.
\]
Hence,~$\h_s u$ can be interpreted as a continuous linear functional $\langle \h_s u,\,\cdot\, \rangle \in ( W^{s, 1}(\R^n))^*$.
Remarkably, this holds for every measurable function~$u: \R^n \to \R$, regardless of its regularity.

We now set
\[
\G_s(t):=\int_0^t G_s(\tau)\,d\tau\quad\mbox{for }t\in\R,
\]
and, given a measurable function $u:\R^n\to\R$ and an open set $\Omega\subseteq\R^n$, we define the \emph{$s$-fractional area functional}
\[
\F_s(u,\Omega):=\iint_{\R^{2n}\setminus(\Co\Omega)^2}\G_s\left(\frac{u(x)-u(y)}{|x-y|}\right)\frac{dx\,dy}{|x-y|^{n-1+s}}.
\]
Then, at least formally, we have
\[
\frac{d}{d\eps}\Big|_{\eps=0}\F_s(u+\eps v,\Omega)=\langle\h_s u,v\rangle\quad\mbox{for every }v\in C^\infty_c(\Omega).
\]

\smallskip

We remark that in Chapter \ref{CH_Nonparametric} we will actually consider more general functionals of fractional area-type---by taking in the above definitions a continuous and even function $g:\R\to(0,1]$ satisfying an appropriate integrability condition, and the corresponding functions $G$ and $\G$, in place of $g_s,\,G_s$ and $\G_s$, respectively. However, for simplicity in this introduction we stick to the ``geometric case" corresponding to the choice $g=g_s$.

\smallskip

Let us now get to the functional properties of $\F_s(\,\cdot\,,\Omega)$ and to its relationship with the fractional perimeter.

From now on, we fix $n\ge1$, $s\in(0,1)$ and a bounded open set $\Omega\subseteq\R^n$ with Lipschitz boundary.

It is convenient to split the fractional area functional as the sum of its local and nonlocal parts, that is
\[
\F_s(u,\Omega)=\A_s(u,\Omega)+\Nl_s(u,\Omega),
\]
with
\[
\A_s(u,\Omega):=\int_\Omega\int_\Omega
\G_s\left(\frac{u(x)-u(y)}{|x-y|}\right)\frac{dx\,dy}{|x-y|^{n-1+s}}
\]
and
\[
\Nl_s(u,\Omega):=2\int_\Omega\int_{\Co\Omega}
\G_s\left(\frac{u(x)-u(y)}{|x-y|}\right)\frac{dx\,dy}{|x-y|^{n-1+s}}.
\]
Let us first mention the following interesting observation---see, e.g., Lemma \ref{CH:A:usef_ineq_hit}.
If $u:\Omega\to\R$ is a measurable function, then
\[
[u]_{W^{s,1}(\Omega)}<\infty\quad\implies\quad\|u\|_{L^1(\Omega)}<\infty.
\]

Concerning the local part of the fractional area functional, we prove that, if $u:\Omega\to\R$ is a measurable function, then
\bgs{
\A_s(u,\Omega)<\infty\quad&\Longleftrightarrow\quad u\in W^{s,1}(\Omega)\\
&
\Longleftrightarrow\quad
\Per_s^L(\Sg(u),\Omega\times\R)<\infty.
}
Moreover, if $u\in W^{s,1}(\Omega)$, then
\[
\Per_s^L(\Sg(u),\Omega\times\R)
=\A_s(u,\Omega)+c,
\]
for some constant $c=c(n,s,\Omega)\geq0$.
These results can be thought of as the fractional counterparts of \eqref{CH:0:obs1_locarea} and \eqref{CH:0:obs2_locarea}.

On the other hand, in order for the nonlocal part to be finite, we have to impose some integrability condition on $u$ at infinity, namely
\eqlab{\label{CH:0:Infty_int_cond}
\int_\Omega\left(\int_{\Co\Omega}
\frac{|u(y)|}{|x-y|^{n+s}}\,dy\right)dx<\infty.
}
Such a condition is satisfied, e.g., if $u$ is globally bounded in $\R^n$ and, in general, it implies that the function $u$ must grow strictly sublinearly at infinity. It is thus a very restrictive condition.

Indeed, we remark that the operator $\h_s u$ is well-defined at a point $x$---provided $u$ is regular enough in a neighborhood of $x$---without having to impose any condition on $u$ at infinity. Moreover, as we have observed in Section
\ref{CH:0:NMS_Sec}, by Theorem \ref{CH:0:Local_Existence} and \cite[Theorem 1.1]{graph} we know that, fixed any continuous function $\varphi:\R^n\to\R$, there exists a function $u:\R^n\to\R$ such that $u=\varphi$ in $\R^n\setminus\overline{\Omega}$, $u\in C(\overline{\Omega})$ and $\Sg(u)$ is locally $s$-minimal in $\Omega^\infty$. Let us stress that no condition on $\varphi$ at infinity is required.

For these reasons, condition \eqref{CH:0:Infty_int_cond} seems to be unnaturally restrictive in our framework---even if at first glance it looks necessary, since it is needed to guarantee that $\F_s$ is well-defined.

In order to avoid imposing condition \eqref{CH:0:Infty_int_cond}, we define---see \eqref{CH:4:NMldef}---for every $M\geq0$, the 
``truncated" nonlocal part $\Nl_s^M(u,\Omega)$ and the truncated area functional
\[
\F_s^M(u,\Omega):=\A_s(u,\Omega)+\Nl_s^M(u,\Omega).
\]
Roughly speaking, the idea consists in adding, inside the double integral defining the nonlocal part, a term which balances the contribution coming from outside $\Omega$. For example, in the simplest case $M=0$, we have
\[
\Nl_s^0(u,\Omega)=2\int_\Omega\left\{\int_{\Co\Omega}\left[\G_s\left(\frac{u(x)-u(y)}{|x-y|}\right)-\G_s\left(\frac{u(y)}{|x-y|}\right)\right]\frac{dy}{|x-y|^{n-1+s}}\right\}dx.
\]
Remarkably, given a measurable function $u:\R^n\to\R$, we have
\[
|\Nl_s^M(u,\Omega)|<\infty\quad\mbox{if }u|_\Omega\in W^{s,1}(\Omega),
\]
regardless of the behavior of $u$ in $\Co\Omega$. However, we remark that, in general, the truncated nonlocal part can be negative, unless we require $u$ to be bounded in $\Omega$ and we take $M\ge\|u\|_{L^\infty(\Omega)}$. From a geometric point of view, the truncated area functionals correspond to considering the fractional perimeter in the truncated cylinder $\Omega^M$.

Indeed, if $u:\R^n\to\R$ is a measurable function such that $u|_\Omega\in W^{s,1}(\Omega)\cap L^\infty(\Omega)$, and $M\ge\|u\|_{L^\infty(\Omega)}$, we have
\[
\F_s^M(u,\Omega)=\Per_s	\big(\Sg(u),\Omega\times(-M,M)\big)+c_M,
\]
for some constant $c_M=c_M(n,s,\Omega)\geq0$.

\medskip

We now proceed to study the minimizers of the fractional area functional.

Given a measurable function $\varphi:\Co\Omega\to\R$, we define the space
\[
\W_\varphi^s(\Omega):=\left\{
u:\R^n\to\R\,|\,u|_\Omega\in W^{s,1}(\Omega)\mbox{ and }u=\varphi\mbox{ a.e. in }\Co\Omega\right\},
\]
and we say that $u\in\W^s_\varphi(\Omega)$ is a \emph{minimizer} of $\F_s$ in $\W^s_\varphi(\Omega)$ if
\[
\iint_{Q(\Omega)} \left\{ \G_s \left( \frac{u(x) - u(y)}{|x - y|} \right) - \G_s \left( \frac{v(x) - v(y)}{|x - y|} \right) \right\} \frac{dx\,dy}{|x - y|^{n - 1 + s}} \le 0
\]
for every $v\in\W^s_\varphi(\Omega)$. Here above, we have used the notation $Q(\Omega):=\R^{2n}\setminus(\Co\Omega)^2$.
Let us stress that such a definition is well-posed without having to impose conditions on the \emph{exterior data} $\varphi$, as indeed---thanks to the fractional Hardy-type inequality of Theorem \ref{CH:4:FHI}---we have
\[
\iint_{Q(\Omega)}\left|\G_s\left(\frac{u(x)-u(y)}{|x-y|}\right)
-\G_s\left(\frac{v(x)-v(y)}{|x-y|}\right)\right|\frac{dx\,dy}{|x-y|^{n-1+s}}
\le C\,\Lambda_{n,s}\|u-v\|_{W^{s,1}(\Omega)},
\]
for every~$u,\,v\in\W^s_\varphi(\Omega)$,
for some constant $C=C(n,s,\Omega)>0$.

We prove the existence of minimizers with respect to exterior data satisfying an appropriate integrability condition in a neighborhood of the domain $\Omega$. More precisely, given an open set~$\Op\subseteq\R^n$
such that~$\Omega\Subset\Op$, we define
the \emph{truncated tail} of~$\varphi:\Co\Omega\to\R$ at a point~$x\in\Omega$ as
\bgs{
	\Tail_s(\varphi, \Op\setminus\Omega;x):=\int_{\Op\setminus\Omega}\frac{|\varphi(y)|}{|x-y|^{n+s}}\,dy.
}
We also use the notation
\[
\Omega_\varrho:=\{x\in\R^n\,|\,d(x,\Omega)<\varrho\},
\]
for $\varrho>0$, to denote the $\varrho$-neighborhood of $\Omega$.
Then, we prove the following:
\begin{theorem} 
	There is a constant~$\Theta > 1$, depending only on~$n$ and~$s$, such that,
	given any function~$\varphi: \Co\Omega\to \R$
	with~$\Tail_s(\varphi, \Omega_{\Theta \diam(\Omega)} \setminus \Omega;\,\cdot\,) \in L^1(\Omega)$, there exists a unique minimizer~$u$ of~$\F_s$ within~$\W^s_\varphi(\Omega)$. Moreover,~$u$ satisfies
	\[
	\| u \|_{W^{s, 1}(\Omega)} \le C \left( \left\| \Tail_s(\varphi,\Omega_{\Theta \diam(\Omega)}\setminus \Omega;\,\cdot\,) \right\|_{L^1(\Omega)} + 1 \right),
	\]
	for some constant~$C=C(n,s,\Omega)>0$.
\end{theorem}
We observe that the condition on the integrability of the tail is much weaker than \eqref{CH:0:Infty_int_cond}, since we are not requiring anything on the behavior of $\varphi$ outside $\Omega_{\Theta \diam(\Omega)}$. 

We also mention that, roughly speaking, the integrability of the tail amounts to the integrability of $\varphi$ plus some regularity condition near the boundary of $\partial\Omega$. For example, if $\varphi\in L^1(\Omega_{\Theta \diam(\Omega)}\setminus\Omega)$ and there exists a $\varrho>0$ such that, either $\varphi\in W^{s,1}(\Omega_\varrho\setminus\Omega)$ or
$\varphi\in L^\infty(\Omega_\varrho\setminus\Omega)$, then
$\Tail_s(\varphi, \Omega_{\Theta \diam(\Omega)} \setminus \Omega;\,\cdot\,) \in L^1(\Omega)$.

The uniqueness of the minimizer follows from the strict convexity of $\F_s$. 
On the other hand, in order to prove the existence, we exploit the (unique) minimizers $u_M$ of the functionals $\F^M_s(\,\cdot\,,\Omega)$---considered within their natural domain. We exploit the hypothesis on the integrability of the tail, to prove a uniform estimate for the $W^{s,1}(\Omega)$ norm of the minimizers $u_M$, independently on $M\geq0$. Hence, up to subsequences, $u_M$ converges, as $M\to\infty$, to a limit function $u$, which is easily proved to minimize $\F_s$.

\smallskip

Moreover, we prove that if $u$ is a minimizer of $\F_s$ within $\W^s_\varphi(\Omega)$, then $u\in L^\infty_{\loc}(\Omega)$. Also, we show that if the exterior data $\varphi$ is bounded in a big enough neighborhood of $\Omega$, then $u\in L^\infty(\Omega)$, and we establish an apriori bound on the $L^\infty$ norm.

\smallskip

Let us go back to the relationship between the fractional area functional and the fractional perimeter.
We show that by appropriately rearranging a set $E$ in the vertical direction we decrease the $s$-perimeter. 
More precisely, given a set~$E \subseteq \R^{n + 1}$, we consider the function~$w_E: \R^n \to \R$ defined by
\bgs{
	w_E(x) := \lim_{R \rightarrow +\infty} \left( \int_{-R}^R \chi_{E}(x, t) \, dt - R \right)
}
for every~$x \in \R^n$.

Then, we have the following result.

\begin{theorem}\label{CH:0:Rearr_THM}
	Let~$E \subseteq \R^{n + 1}$ be such that~$E \setminus \Omega^\infty$ is a subgraph and
\[
	\Omega \times (-\infty, -M) \subseteq E \cap \Omega^\infty \subseteq \Omega \times (-\infty, M),
\]
	for some~$M > 0$. Then,
	\[
		\Per_s(\Sg(w_E), \Omega^M) \le \Per_s(E, \Omega^M).
	\]
	The inequality is strict unless~$\Sg(w_E)=E$.
\end{theorem}

Exploiting also the local boundedness of a minimizer, we prove that if $u:\R^n\to\R$
is a measurable function such that $u\in W^{s,1}(\Omega)$, then
\[
u\mbox{ minimizes }\F_s\mbox{ within }\W^s_u(\Omega)\quad\implies\quad
\Sg(u)\mbox{ is locally $s$-minimal in } \Omega^\infty.
\]
Theorem \ref{CH:0:Rearr_THM} extends to the fractional framework a well known result holding for the classical perimeter---see, e.g., \cite[Lemma 14.7]{Giusti}. However, notice that in the fractional framework, due to the nonlocal character of the functionals involved, we have to assume that the set $E$ is already a subgraph outside the cylinder $\Omega^\infty$.

We also observe that, since $u$ is locally bounded in $\Omega$ and its subgraph is locally $s$-minimal in the cylinder $\Omega^\infty$, by \cite[Theorem 1.1]{CaCo}
we have that $u\in C^\infty(\Omega)$---that is, minimizers of $\F_s$ are smooth.

\medskip

Let us now get to the Euler-Lagrange equation satisfied by minimizers.
We first introduce the notion of weak solutions.

Let $f\in C(\overline{\Omega})$. We say that a measurable function $u:\R^n\to\R$ is a weak solution of $\h_s u=f$ in $\Omega$ if
\[
\langle\h_s u,v\rangle=\int_\Omega fv\,dx,
\]
for every $v\in C^\infty_c(\Omega)$.

As a consequence of the convexity of $\F_s$, it is easy to prove that, given a measurable function $u:\R^n\to\R$ such that $u\in W^{s,1}(\Omega)$, it holds
\[
u \mbox{ is a minimizer of }\F_s\mbox{ in }\W^s_u(\Omega)\quad
\Longleftrightarrow\quad
u\mbox{ is a weak solution of }\h_s u=0\mbox{ in }\Omega.
\]

Another natural notion of solution for the equation $\h_s u=f$ is that of a viscosity solution---we refer to Section \ref{CH:4:ViscWeak_Sec} for the precise definition.
One of the main results of Chapter 
\ref{CH_Nonparametric} consists in proving that viscosity (sub)solutions are weak (sub)solutions. More precisely:

\begin{theorem}
	Let~$\Omega\subseteq\R^n$ be a bounded open set and let~$f\in C(\overline\Omega)$. Let~$u:\R^n\to\R$
	be such that~$u$ is locally integrable in~$\R^n$ and~$u$ is locally bounded in~$\Omega$. If~$u$
	is a viscosity subsolution,
	\bgs{
		\h_s u\le f\quad\mbox{in }\Omega,
	}
	then~$u$ is a weak subsolution,
	\bgs{
		\langle\h_s u,v\rangle\le\int_\Omega fv\,dx,\qquad\forall\,v\in C^\infty_c(\Omega)\mbox{ s.t. }v\ge0.
	}
\end{theorem}

\medskip

Combining the main results of Chapter \ref{CH_Nonparametric} and exploiting the interior regularity proved in~\cite{CaCo},
we obtain the following:

\begin{theorem}\label{CH:0:Equiv_Intro}
	Let~$u:\R^n\to\R$ be a measurable function such that $u\in W^{s,1}(\Omega)$.
	Then, the following are equivalent:
	\begin{itemize}
		\item[(i)] $u$ is a weak solution of~$\h_s u=0$ in~$\Omega$,
		\item[(ii)] $u$ minimizes~$\F_s$ in~$\W^s_u(\Omega)$,
		\item[(iii)] $u\in L^\infty_{\loc}(\Omega)$ and~$\Sg(u)$ is locally $s$-minimal in~$\Omega\times\R$,
		\item[(iv)] $u\in C^\infty(\Omega)$ and~$u$ is a pointwise solution of~$\h_s u=0$ in~$\Omega$.
	\end{itemize}
	Moreover, if~$u\in L^1_{\loc}(\R^n)\cap W^{s,1}(\Omega)$, then all of the above are equivalent to:
	\begin{itemize}
		\item[(v)] $u$ is a viscosity solution of~$\h_s u=0$ in~$\Omega$.
	\end{itemize}
\end{theorem}

We also point out the following global version of Theorem~\ref{CH:0:Equiv_Intro}:

\begin{corollary}\label{CH:0:Equiv_Intro_Global_Corollary}
	Let~$u\in W^{s,1}_{\loc}(\R^n)$. Then, the following are equivalent:
	\begin{itemize}
		\item[(i)] $u$ is a viscosity solution of~$\h_s u=0$ in~$\R^n$,
		\item[(ii)] $u$ is a weak solution of~$\h_s u=0$ in~$\R^n$,
		\item[(iii)] $u$ minimizes~$\F_s$ in~$\W^s_u(\Omega)$, for every open set~$\Omega\Subset\R^n$ with Lipschitz boundary,
		\item[(iv)] $u\in L^\infty_{\loc}(\R^n)$ and~$\Sg(u)$ is locally~$s$-minimal in~$\R^{n+1}$,
		\item[(v)] $u\in C^\infty(\R^n)$ and~$u$ is a pointwise solution of~$\h_s u=0$ in~$\R^n$.
	\end{itemize}
\end{corollary}

\medskip

Let us also mention that the functional framework introduced above, easily extends to the obstacle problem. Namely, besides imposing the exterior data condition $u=\varphi$ a.e. in $\Co\Omega$, we constrain the functions to lie above an obstacle, that is, given an open set $A\subseteq\Omega$ and an obstacle $\psi\in L^\infty(A)$, we restrict ourselves to consider those functions $u\in\W^s_\varphi(\Omega)$ such that $u\ge\psi$ a.e. in $A$.

In Chapter \ref{CH_Nonparametric} we briefly cover also this obstacle problem, proving the existence and uniqueness of a minimizer and its relationship with the geometric obstacle problem that involves the fractional perimeter.

\smallskip

Finally, in the last section of Chapter \ref{CH_Nonparametric}, we prove some approximation results for subgraphs having (locally) finite fractional perimeter. In particular, exploiting the surprising density result established in \cite{DSV17}, we show that $s$-minimal subgraphs can be appropriately approximated by subgraphs of $\sigma$-harmonic functions, for any fixed $\sigma\in(0,1)$.

\subsection{Rigidity results for nonlocal minimal graphs}

In Chapter \ref{CH_Bern_Mos_result} we prove a flatness result for entire nonlocal minimal graphs having some partial derivatives bounded from either above or below. This result generalizes fractional versions of classical theorems due to Bernstein and Moser.

Moreover, we show that entire graphs having constant fractional mean curvature are minimal, thus extending a celebrated result of Chern on classical CMC graphs.

\medskip

We are interested in subgraphs that locally minimize the $s$-perimeter in the whole space $\R^{n+1}$.
We recall that, as we have seen in Corollary \ref{CH:0:Equiv_Intro_Global_Corollary},
under very mild assumptions on the function $u:\R^n\to\R$,
the subgraph $\Sg(u)$ is locally $s$-minimal in $\R^{n+1}$ if and only if
$u$ satisfies the fractional mean curvature equation
\eqlab{ \label{CH:0:Hsu=0}
\h_su=0\quad\mbox{in }\R^n.
}
Moreover, again by Corollary \ref{CH:0:Equiv_Intro_Global_Corollary}, there are several equivalent notions of solution for the equation \eqref{CH:0:Hsu=0},
such as smooth solutions, viscosity solutions, and weak solutions.

In what follows, a solution of~\eqref{CH:0:Hsu=0} will always indicate a function~$u \in C^\infty(\R^n)$ that satisfies identity~\eqref{CH:0:Hsu=0} pointwise. We stress that no growth assumptions at infinity are made on~$u$.

The main contribution of Chapter \ref{CH_Bern_Mos_result} is the following result.

\begin{theorem} \label{CH:0:under_CRA_Pakmainthm}
	Let~$n \ge \ell \ge 1$ be integers,~$s \in (0, 1)$, and suppose that
	\begin{equation} \tag{$P_{s, \ell}$} \label{CH:0:SMINSINGCON}
	\mbox{there exist no singular~$s$-minimal cones in~$\R^\ell$.}
	\end{equation}
	Let~$u$ be a solution of $\h_su=0$ in $\R^n$, having~$n - \ell$ partial derivatives bounded on one side.
	Then,~$u$ is an affine function.
\end{theorem}

Characterizing the values of~$s$ and~$\ell$ for which~\eqref{CH:0:SMINSINGCON} is satisfied represents a challenging open problem.
Nevertheless, property~\eqref{CH:0:SMINSINGCON} is known to hold in the following cases:
\begin{itemize}
	\item when~$\ell = 1$ or~$\ell = 2$, for every~$s \in (0, 1)$;
	\item when~$3\le \ell \le 7$ and~$s \in (1 - \varepsilon_0, 1)$ for some small~$\varepsilon_0 \in (0,1]$ depending only on~$\ell$. 
\end{itemize}
Case~$\ell = 1$ holds by definition, while~$\ell = 2$ is the content of~\cite[Theorem~1]{SV13}. On the other hand, case~$3 \le \ell \le 7$ has been established in~\cite[Theorem~2]{regularity}.

As a consequence of Theorem~\ref{CH:0:under_CRA_Pakmainthm} and the last remarks, we immediately obtain the following result.

\begin{corollary} \label{CH:0:CORRigidity_mainthm}
	Let~$n \ge \ell \ge 1$ be integers and~$s \in (0, 1)$. Assume that either
	\begin{itemize}
		\item $\ell \in \{ 1, 2 \}$, or
		\item $3\le \ell \le 7$ and~$s \in (1 - \varepsilon_0, 1)$, with~$\varepsilon_0=\varepsilon_0(\ell) > 0$ as in~\cite[Theorem~2]{regularity}.
	\end{itemize}
	Let~$u$ be a solution of $\h_su=0$ in $\R^n$, having~$n-\ell$ partial derivatives bounded on one side.
	Then,~$u$ is an affine function.
\end{corollary}

We observe that
Theorem~\ref{CH:0:under_CRA_Pakmainthm} gives a new flatness result for~$s$-minimal graphs, under the assumption that~\eqref{CH:0:SMINSINGCON} holds true. It can be seen as a generalization of the fractional~De~Giorgi-type lemma contained in~\cite[Theorem~1.2]{FV17}, which is recovered here taking~$\ell = n$. In this case, we indeed provide an alternative proof of said result.

On the other hand, the choice~$\ell = 2$ gives an improvement of~\cite[Theorem~4]{FarV17}, when specialized to~$s$-minimal graphs. In light of these observations, Theorem~\ref{CH:0:under_CRA_Pakmainthm} and Corollary~\ref{CH:0:CORRigidity_mainthm} can be seen as a bridge between Bernstein-type theorems (flatness results in low dimensions) and Moser-type theorems (flatness results under global gradient bounds).

For classical minimal graphs, the counterpart of Corollary~\ref{CH:0:CORRigidity_mainthm} has been recently obtained by A. Farina in~\cite{F17}. In that case, the result is sharp and holds with~$\ell = \min \{ n, 7 \}$.
The proof of Theorem~\ref{CH:0:under_CRA_Pakmainthm} is based on the extension to the fractional framework of a strategy---which relies on a general splitting result for blow-downs of the subgraph $\Sg(u)$---devised by A. Farina for classical minimal graphs and previously unpublished. As a result, the ideas contained in Chapter \ref{CH_Bern_Mos_result} can be used to obtain a different, easier proof of~\cite[Theorem~1.1]{F17}

Let us also mention that, by using the same ideas that lead to Theorem~\ref{CH:0:under_CRA_Pakmainthm}, we can prove the following rigidity result for entire~$s$-minimal graphs that lie above a cone.

\begin{theorem}\label{CH:0:growth_THM}
	Let~$n\ge1$ be an integer and~$s\in(0,1)$. Let~$u$ be a solution of $\h_su=0$ in $\R^n$, and assume that there exists a constant~$C > 0$ for which
\[
	u(x)\ge - C (1+|x|)\quad\mbox{for every }x\in\R^n.
\]
	Then,~$u$ is an affine function.
\end{theorem}

We remark that in~\cite{CaCo} a rigidity result analogous to Theorem~\ref{CH:0:growth_THM} is deduced,  under the stronger, two-sided assumption
\[
|u(x)| \le C(1 + |x|)\quad\mbox{for every }x\in\R^n.
\]
Theorem~\ref{CH:0:growth_THM} thus improves~\cite[Theorem~1.5]{CaCo} directly. 

\medskip

Finally, we prove that if $u:\R^n\to\R$ is
such that
\[
\langle\h_su,v\rangle=h\int_{\R^n}v\,dx\quad\mbox{for every }v\in C^\infty_c(\R^n),
\]
for some constant $h\in\R$, then the constant must be $h=0$.

In particular, recalling Corollary \ref{CH:0:Equiv_Intro_Global_Corollary}, we see that if $u\in W^{s,1}_{\loc}(\R^n)$ is a weak solution of $\h_su=h$ in $\R^n$, then the subgraph of $u$ is locally $s$-minimal in $\R^{n+1}$.
This extends to the nonlocal framework a celebrated result of Chern, namely the Corollary of Theorem~1 in~\cite{C65}.

\subsection{A free boundary problem}

In Chapter \ref{CH_FreeBdary_CHPT} we study minimizers of the functional
\eqlab{\label{CH:0:NOnLOC_fReeBdA}
\Nl(u,\Omega)+\Per\big(\{u>0\},\Omega\big),
}
with $\Nl(u,\Omega)$ being, roughly speaking, the $\Omega$-contribution to the $H^s$ seminorm of a function $u:\R^n\to\R$, that is
\[
\Nl(u,\Omega):=\iint_{\R^{2n}\setminus(\Co\Omega)^2}\frac{|u(x)-u(y)|^2}{|x-y|^{n+2s}}\,dx\,dy,
\]
for some fixed index $s\in(0,1)$.

Similar functionals, defined as the superposition of an ``elastic energy" plus a ``surface tension" term, have already been considered in the following papers:
\begin{itemize}
	\item Dirichlet energy plus classical perimeter in \cite{ACKS},
	\item Dirichlet energy plus fractional perimeter in \cite{CSV},
	\item the nonlocal energy $\Nl$ plus the fractional perimeter in \cite{DSV}, and the corresponding one-phase problem in \cite{DV-onephase}.
\end{itemize}
Studying the functional defined in \eqref{CH:0:NOnLOC_fReeBdA} somehow completes this picture.

\smallskip

The main contributions of
Chapter \ref{CH_FreeBdary_CHPT}
consist in establishing a monotonicity formula for the minimizers of the functional \eqref{CH:0:NOnLOC_fReeBdA}, in exploiting it to investigate the properties of blow-up limits and in proving a dimension reduction result. Moreover, we show that, when $s<1/2$, the perimeter dominates---in some sense---over the nonlocal energy. As a consequence, we obtain a regularity result for the free boundary $\{u=0\}$.

\medskip

As a thechnical note, let us first observe that we can not directly work with the set $\{u>0\}$. Instead, we consider \emph{admissible pairs} $(u,E)$, with $u:\R^n\to\R$ a measurable function, and $E\subseteq\R^n$ such that
\[
u\ge0\quad\mbox{a.e. in }E\quad\mbox{and}\quad u\leq0\quad\mbox{a.e. in }\Co E.
\]
The set $E$ is usually referred to as the \emph{positivity set} of $u$. Then, given an index $s\in(0,1)$ and a bounded open set with Lipschitz boundary $\Omega\subseteq\R^n$, we define the functional
\[
\F_\Omega(u,E):=\Nl(u,\Omega)+\Per(E,\Omega),
\]
for every admissible pair $(u,E)$.

Let us now remark that if $u:\R^n\to\R$ is a measurable function, then
\eqlab{\label{CH:0:TAiL_MiN_ObsV}
\Nl(u,\Omega)<\infty\quad\implies\quad\int_{\R^n}\frac{|u(\xi)|^2}{1+|\xi|^{n+2s}}\,d\xi<\infty.
}
For a proof see, e.g., Lemma \ref{CH:APP:usef_ineq_tail}.
As a consequence, we also have that
$$
\int_{\R^n}\frac{|u(\xi)|}{1+|\xi|^{n+2s}}\,d\xi<\infty\quad\mbox{ and }\quad u\in L^2_{\loc}(\R^n).
$$

The notion of minimizers that we consider is the following:

\begin{defn}\label{CH:0:def_minim_FrEeBdAry}
	Given an admissible pair~$(u,E)$ such that $\F_\Omega(u,E)<\infty$, we say that a pair $(v,F)$ 
	is an
	\emph{admissible competitor} if \eqlab{\label{CH:0:compet_def}
	&\textrm{supp}(v-u)\Subset\Omega,\qquad F\Delta E\Subset\Omega,\\
	&v-u\in H^s(\R^n)\qquad\textrm{and}\qquad \Per(F,\Omega)<+\infty.
	}
	We say that the admissible pair $(u,E)$
	is \emph{minimizing} in $\Omega$ if $\F_\Omega(u,E)<\infty$ and
	\[
	\F_\Omega(u,E)\leq\F_\Omega(v,F),
	\]
	for every admissible competitor~$(v,F)$.
\end{defn}
Notice that the first line of \eqref{CH:0:compet_def} simply says that the pairs $(u,E)$ and $(v,F)$ are equal---in the measure theoretic sense---outside a compact subset of $\Omega$. Then, since $\F_\Omega(u,E)<\infty$, it is readily seen that the second line is equivalent to $\F_\Omega(v,F)<\infty$.

\smallskip

In particular we are interested in the following minimization problem, with respect to fixed ``exterior data''.
Given an admissible pair $(u_0,E_0)$ and a bounded open set $\Op\subseteq\R^n$ with Lipschitz boundary, such that
\eqlab{\label{CH:0:Dir_data}
\Omega\Subset\Op,\qquad \Nl(u_0,\Omega)<+\infty
\quad\textrm{ and }\quad \Per(E_0,\Op)<+\infty,
}
we want to find an admissible pair $(u,E)$ attaining the following infimum
\eqlab{\label{CH:0:Min_dir}
\inf\big\{\Nl(v,\Omega)+\Per(F,\Op)\,|\,(v,F)&\textrm{ admissible pair s.t. }v=u_0\textrm{ a.e. in }\Co\Omega\\
&
\quad
\textrm{ and }F\setminus\Omega=E_0\setminus\Omega\big\}.
}
Roughly speaking, as customary when dealing with minimization problems involving the classical perimeter, we are considering a (fixed) neighborhood $\Op$ of $\Omega$ (as small as we like)
in order to ``read'' the boundary data $\partial E_0\cap\partial\Omega$.

We prove that, fixed as exterior data any pair $(u_0,E_0)$ satisfying \eqref{CH:0:Dir_data}, there exists a pair $(u,E)$ realizing the infimum in \eqref{CH:0:Min_dir}. Moreover, we show that such a pair $(u,E)$
is also minimizing in the sense of Definition \ref{CH:0:def_minim_FrEeBdAry}.

\smallskip

A useful result consists in establishing a uniform bound for the energy of minimizing pairs.

\begin{theorem}\label{CH:0:TH:unif}
	Let~$(u,E)$ be a minimizing pair in~$B_2$. 
	Then
	\[ \iint_{\R^{2n}\setminus (\Co B_1)^2}\frac{|u(x)-u(y)|^2}{
		|x-y|^{n+2s}}\,dx\,dy + \Per(E,B_1)\le 
	C\left(1+\int_{\R^n}\frac{|u(y)|^2}{1+|y|^{n+2s}}\,dy\right),
	\]
	for some~$C=C(n,s)>0$.
\end{theorem}

In particular, Theorem \ref{CH:0:TH:unif} is exploited in the proof of the existence of a blow-up limit.
For this, we have first to introduce---through the extension technique of \cite{CS07}---the extended functional associated to the minimization of $\F_\Omega$. We write
$$
\R^{n+1}_+:=\{(x,z)\in\R^{n+1} {\mbox{ with }} x\in\R^n,\,z>0\}.
$$
Given a function~$u:\R^n\to\R$,
we consider the function~$\Ue:\R_+^{n+1}\to\R$
defined via the convolution with an appropriate Poisson kernel,
\bgs{
\Ue(\,\cdot\,,z)=u\ast\mathcal K_s(\,\cdot\,,z),\quad\textrm{where}\quad\K_s(x,z):=c_{n,s}\frac{z^{2s}}{(|x|^2+z^2)^{(n+2s)/2}},
}
and~$c_{n,s}>0$ is an appropriate normalizing constant. Such an extended function $\Ue$ is well defined---see, e.g., \cite{Extension}---provided $u:\R^n\to\R$ is such that
\[
\int_{\R^n}\frac{|u(\xi)|}{1+|\xi|^{n+2s}}\,d\xi<\infty.
\]
In light of \eqref{CH:0:TAiL_MiN_ObsV}, we can thus consider the extension function of a minimizer.

We use capital letters, like $X=(x,z)$, to
denote points in $\R^{n+1}$.
Given a set $\Omega\subseteq\R^{n+1}$, we write
\bgs{
\Omega_+:=\Omega\cap\{z>0\}\qquad\textrm{and}\qquad\Omega_0:=\Omega\cap\{z=0\}.
}
Moreover we identify the hyperplane $\{z=0\}\simeq\R^n$ via the projection function.

Given a bounded open set $\Omega\subseteq\R^{n+1}$ 
with Lipschitz boundary, such that $\Omega_0\not=\emptyset$, we define
\bgs{
\lf_\Omega(\Vf,F):=c_{n,s}'\int_{\Omega_+}|\nabla\Vf|^2z^{1-2s}\,dX+\Per(F,\Omega_0),
}
for $\Vf:\R^{n+1}_+\to\R$ and $F\subseteq\R^n\simeq\{z=0\}$
the positivity set of the trace of $\Vf$ on $\{z=0\}$, that is
\bgs{
\Vf\big|_{\{z=0\}}\geq0\quad\textrm{a.e. in }F\quad\textrm{and}\quad
\Vf\big|_{\{z=0\}}\leq0\quad\textrm{a.e. in }\Co F.
}
We call such a pair $(\Vf,F)$ an {\emph{admissible pair}}
for the extended functional. Then, we introduce the following notion of minimizer for the extended functional.

\begin{defn}\label{CH:0:adm-ext}
	Given an admissible pair~$(\Uc,E)$, such that $\lf_\Omega(\Uc,E)<\infty$, we say that a pair~$(\Vf,F)$
	is an \emph{admissible competitor}
	if~$\lf_\Omega(\Vf,F)<\infty$ and
	\begin{equation*}
	\textrm{supp}\,(\Vf-\Uc)\Subset\Omega\qquad
	\textrm{and}\qquad E\Delta F\Subset\Omega_0.
	\end{equation*}
	We say that an admssible pair $(\Uc,E)$ is \emph{minimal} in~$\Omega$
	if~$\lf_\Omega(\Uc,E)<\infty$ and
	\begin{equation*}
	\lf_\Omega(\Uc,E)\leq\lf_\Omega(\Vf,F),
	\end{equation*}
	for every admissible competitor $(\Vf,F)$.
\end{defn}

An important result consists in showing that an appropriate minimization problem involving the extended functionals is equivalent to the minimization of the original functional $\F_\Omega$. More precisely:

\begin{prop}\label{CH:0:Local_energy_prop}
	Let $(u,E)$ be an admissible pair for $\F$, s.t.~$\F_{B_R}(u.E)<+\infty$.
	Then, the pair $(u,E)$ is minimizing in $B_R$ if and only 
	if the pair $(\Ue,E)$ is minimal for $\lf_\Omega$,
	in every bounded open set $\Omega\subseteq\R^{n+1}$ 
	with Lipschitz boundary such that~$\emptyset\not=\Omega_0\Subset B_R$.
\end{prop}

One of the main reasons for introducing the extended functional, resides in the fact that it enables us to establish a Weiss-type monotonicity formula for minimizers.

We denote
\bgs{
\BaLL_r:=\{(x,z)\in\R^{n+1}\,|\,|x|^2+z^2<r^2\}\qquad\textrm{and}\qquad
\BaLL_r^+:=\BaLL_r\cap\{z>0\}.
}

\begin{theorem}[Weiss-type Monotonicity Formula]\label{CH:0:Monotonicity_teo}
	Let $(u,E)$ be a minimizing pair for $\F$ in $B_R$ and define the function $\Phi_u:(0,R)\to\R$ by
	\bgs{
	\Phi_u(r):=r^{1-n}&
	\left(c'_{n,s}\int_{\BaLL_r^+}|\nabla\Ue|^2z^{1-2s}\,dX+\Per(E,B_r)\right)\\
	&\qquad
	-c'_{n,s}\Big(s-\frac{1}{2}\Big)r^{-n}\int_{(\partial\BaLL_r)^+}\Ue^2z^{1-2s}\,d\Ha^n.
	}
	Then, the function $\Phi_u$ is increasing in $(0,R)$.	
	Moreover, $\Phi_u$ is constant in $(0,R)$ if and only if the
	extension $\Ue$ is homogeneous of degree $s-\frac{1}{2}$ in $\BaLL_R^+$
	and $E$ is a cone in $B_R$.
\end{theorem}

Here above, $(\partial\BaLL_r)^+:=\partial\BaLL_r\cap\{z>0\}$.
Let us now introduce the rescaled pairs $(u_\lambda,E_\lambda)$.
Given $u:\R^n\to\R$ and $E\subseteq\R^n$, we define
\bgs{
	u_\lambda(x):=\lambda^{\frac{1}{2}-s}u(\lambda x)\qquad\textrm{and}\qquad E_\lambda:=\frac{1}{\lambda}E,
}
for every $\lambda>0$.
We observe that---because of the scaling properties of $\F_\Omega$---a pair $(u,E)$ is minimal in $\Omega$ if and only if the rescaled pair $(u_\lambda,E_\lambda)$ is minimal in $\Omega_\lambda$ for every $\lambda>0$.

We prove the convergence of minimizing pairs under appropriate conditions and we exploit it---together with Theorem \ref{CH:0:TH:unif}---in the particularly important case of the blow-up sequence.

We say that the admissible pair~$(u,E)$ is a {\emph{minimizing cone}}
if it is a minimizing pair in $B_R$, for every $R>0$, and is such that~$u$ is homogeneous 
of degree~$s-\frac12$ and~$E$ is a cone 

\begin{theorem}\label{CH:0:TH:blow}
	Let $s>1/2$ and~$(u,E)$ be a minimizing pair in~$B_1$,
	with~$0\in\partial E$. Also assume that
	$u\in C^{s-\frac{1}{2}}(B_1)$.
	Then, there exist a minimizing cone~$(u_0,E_0)$
	and a sequence~$r_k\searrow 0$ such that~$u_{r_k}\to u_0$
	in~$L^\infty_{\loc}(\R^n)$ and~$E_{r_k}\xrightarrow{\loc} E_0$.
\end{theorem}
The homogeneity properties of the blow-up limit $(u_0,E_0)$ are a consequence of Theorem \ref{CH:0:Monotonicity_teo}.

We also point out that we establish appropriate estimates for the tail energies of the functions $u_r$, that allow us to weaken the assumptions of~\cite[Theorem~1.3]{DSV}, where the authors ask~$u$
to be~$C^{s-\frac12}$ in the whole of~$\R^n$.

\smallskip

We now mention the following dimensional reduction result.
Only in the following Theorem, let us redefine
\[\F_\Omega(u,E):=(c_{n,s}')^{-1}\Nl(u,\Omega)+\Per(E,\Omega).\]
We say that an admissible pair $(u,E)$ is minimizing in $\R^n$ if
it minimizes $\F_\Omega$ in any bounded open subset $\Omega\subseteq\R^n$ with Lipschitz boundary.

\begin{theorem}
	Let $(u,E)$ be an admissible pair and define
	\begin{equation*}
	u^\star(x,x_{n+1}):=u(x)\qquad\mbox{and}\qquad
	E^\star:=E\times\R.
	\end{equation*}
	Then, the pair $(u,E)$ is minimizing in $\R^n$
	if and only if the pair $(u^\star,E^\star)$ is minimizing in $\R^{n+1}$.
\end{theorem}

\smallskip

Finally, we observe that in the case $s<1/2$ the perimeter is, in some sense, the leading term of the functional $\F_\Omega$. As a consequence, we are able to prove the following regularity result:

\begin{theorem}
	Let $s\in(0,1/2)$ and let $(u,E)$ be a minimizing pair in $\Omega$. Assume that $u\in L^\infty_{\loc}(\Omega)$.
	Then, $E$ has almost minimal boundary in $\Omega$.
	More precisely, if $x_0\in\Omega$ and $d:=d(x_0,\Omega)/3$,
	then, for every $r\in(0,d]$ it holds
	\bgs{
	\Per(E,B_r(x_0))\le \Per(F,B_r(x_0))+ C\,r^{n-2s},\qquad\forall\,F\subseteq\R^n\mbox{ s.t. }E\Delta F\Subset B_r(x_0),
	}
	where
	\[
	C=C\left(s,x_0,d,\|u\|_{L^\infty(B_{2d}(x_0))},\int_{\R^n}\frac{|u(y)|}{1+|y|^{n+2s}}\,dy\right)>0.
	\]
	Therefore
	\begin{itemize}
		\item[(i)] $\partial^*E$ is locally $C^{1,\frac{1-2s}{2}}$ in $\Omega$,
		
		\item[(ii)] the singular set $\partial E\setminus\partial^*E$ is such that
		\[
		\Ha^\sigma\big((\partial E\setminus\partial^*E)\cap\Omega\big)=0,\qquad\mbox{for every }\sigma>n-8.
		\]
	\end{itemize}
\end{theorem}

\smallskip

We conclude by saying a few words about the one-phase problem, that corresponds to the case in which $u\ge0$ a.e. in $\R^n$. Even if these results are not included in this thesis, they will be part of the final version of the article on which Chapter \ref{CH_FreeBdary_CHPT} is based.
Following the arguments of \cite{DV-onephase}, we will prove that if $(u,E)$ is a minimizer of the one-phase problem in $B_2$, with $s>1/2$, and if $0\in\partial E$, then
$u\in C^{s-\frac{1}{2}}(B_{1/2})$. Notice in particular that, by Theorem \ref{CH:0:TH:blow}, this ensures the existence of a blow-up limit $(u_0,E_0)$.
Moreover, we will establish uniform density estimates for the positivity set $E$, from both sides.

\subsection{The Phillip Island penguin parade
	(a mathematical treatment)}

The goal of Chapter \ref{CH_PEngUInS} is to provide a simple, but rigorous, mathematical
model which describes the formation of groups of penguins
on the shore at sunset. 

\medskip

Penguins are flightless, so they are forced to walk while on land. 
In particular, they show rather specific behaviours in their homecoming,
which are interesting to
observe and to describe analytically.
We observed that penguins have the tendency
to waddle back and forth on the shore to create
a sufficiently large group and then walk home compactly together.
The mathematical framework that we introduce describes
this phenomenon, by taking into account
``natural parameters'', such as the eye-sight of the penguins and
their cruising speed. {The model that we propose
favours the formation of conglomerates of penguins that gather together,
but, on the other hand, it also allows the possibility of isolated and exposed
individuals}.
                           
The model that we propose is based on a set
of ordinary differential equations, with a number of degree
of freedom which is variable in time. Due to the discontinuous
behaviour of the speed of the penguins, the mathematical
treatment (to get existence and uniqueness of
the solution) is based on a ``stop-and-go'' procedure.
We use this setting to provide rigorous examples in which
at least some penguins manage to safely return home (there are also cases
in which some penguins {remain isolated}).
To facilitate the intuition of the model,
we also present some simple numerical simulations
that can be
compared with the actual movement of the penguins parade.

\begin{otherlanguage}{french}

\section{R\'esum\'e}

Cette th\`ese de doctorat est consacr\'ee \`a l'analyse de quelques probl\`emes de minimisation impliquant des fonctionnelles non locales. Les op\'erateurs non locaux ont fait l'objet d'une attention croissante au cours des derni\`eres ann\'ees, \`a la fois par leur int\'er\^et math\'ematique et par leurs applications---par exemple, pour mod\'eliser des processus de diffusion anormaux ou des transitions de phase \`a longue port\'ee.
Pour une introduction aux probl\`emes non locaux, le lecteur int\'eress\'e pourra consulter l'ouvrage \cite{bucval}.

\smallskip

Dans cette th\`ese, nous nous int\'eressons principalement au p\'erim\`etre $s$-fractionnaire---qui peut \^etre consid\'er\'e comme une version fractionnaire et non locale du p\'erim\`etre classique introduit par De Giorgi et Caccioppoli---et ses minimiseurs, les ensembles $s$-minimaux, qui ont \'et\'e consid\'er\'es dans \cite{CRS10} pour la pr\`emiere fois. Les fronti\`eres de ces ensembles $s$-minimaux sont g\'en\'eralement appel\'ees surfaces minimales non locales.
En particulier:
\begin{itemize}
	\item nous \'{e}tudions le comportement des ensembles ayant p\'{e}rim\`{e}tre fractionnaire (localement) fini, en prouvant la densit\'{e} des ensembles ouverts et lisses, un r\'{e}sultat asymptotique optimal pour $s \to 1^-$, et en \'{e}tudiant le lien existant entre le p\'{e}rim\`{e}tre fractionnaire et les ensembles ayant fronti\`{e}res fractales.
	\item Nous \'etablissons des r\'esultats d'existence et de compacit\'e pour les minimiseurs du p\'erim\`etre fractionnaire, qui sont une extension de ceux prouv\'es dans \cite{CRS10}.
	\item Nous \'etudions les ensembles $s$-minimaux dans des r\'egimes hautement non locaux, qui correspondent \`a de petites valeurs du param\`etre fractionnaire $s$. Nous montrons que, dans ce cas, les minimiseurs pr\'esentent un comportement compl\`etement diff\'erent de celui de leurs homologues locaux---les surfaces minimales (classiques).
	\item Nous introduisons un cadre fonctionnel pour \'etudier ces ensembles $s$-minimaux qui peuvent \^etre \'ecrits globalement en tant que sous-graphes. En particulier, nous prouvons des r\'esultats d'existence et d'unicit\'e pour les minimiseurs d'une version fractionnaire de la fonctionnelle d'aire  classique et une in\'egalit\'e de r\'earrangement impliquant que les sous-graphes de ces minimiseurs minimisent le p\'erim\`etre fractionnaire. Nous appelons les fronti\`eres de ces minimiseurs des graphes minimaux non locaux. De plus, nous montrons l'\'equivalence entre les minimiseurs et diverses notions de solution---\`a savoir, solutions faibles, solutions de viscosit\'e et solutions lisses ponctuelles---de l'\'equation de courbure moyenne fractionnaire.
	\item Nous montrons un r\'esultat de platitude pour des graphes minimaux non locaux entiers ayant des d\'eriv\'es partielles major\'ees ou minor\'ees---ainsi, en particulier, \'etendant au cadre fractionnaire des th\'eor\`emes classiques dus \`a Bernstein et Moser.
\end{itemize}

\smallskip

En outre, nous consid\'erons un probl\`eme \`a fronti\`ere libre, qui consiste en la minimisation d'une fonctionnelle d\'efinie comme la somme d'une \'energie non locale, plus le p\'erim\`etre classique de l'interface de s\'eparation entre les deux phases. 
Concernant ce probl\`eme:
\begin{itemize}
	\item nous prouvons l'existence de minimiseurs et introduisons un probl\`eme de minimisation \'equivalent, qui a une ``nature locale''---en exploitant la technique d'extension de \cite{CS07}.
	\item Nous \'etablissons des estimations d'\'energie uniformes et \'etudions la suite de blow-up d'un minimiseur. En particulier, nous prouvons une formule de monotonie qui implique que les limites de blow-up sont homog\`enes.
	\item Nous \'etudions la r\'egularit\'e de la fronti\`ere libre dans le cas o\`u le p\'erim\`etre a un r\^ole dominant sur l'\'energie non locale.
\end{itemize}

\smallskip

Nous mentionnons que le dernier chapitre de la th\`ese consiste en un article fournissant un mod\`ele math\'ematique d\'ecrivant la formation de groupes de manchots sur le rivage au coucher du soleil. \`A l'occasion d'un voyage de recherche \`a l'Universit\'e de Melbourne, nous avons vu le ``Phillip Island penguin parade'' et nous \'etions tellement fascin\'es par le comportement particulier des petits manchots que nous avons d\'ecid\'e de le d\'ecrire de mani\`ere math\'ematique.

\medskip

La th\`ese est divis\'ee en sept chapitres, chacun reposant sur l'un des articles de recherche suivants, que j'ai \'ecrit---seul ou en collaboration---au cours de mon doctorat:
\begin{enumerate}
	\item \emph{Fractional perimeters from a fractal perspective}, publi\'e dans Advanced Nonlinear Studies---voir \cite{Myfractal}. 
	\item \emph{Approximation of sets of finite fractional perimeter by smooth sets and comparison of local and global $s$-minimal surfaces}, publi\'e dans Interfaces and Free Boundaries---voir \cite{mine_cyl_stuff}.
	\item \emph{Complete stickiness of nonlocal minimal surfaces for small values of the fractional parameter}, co-auteur avec C. Bucur et E. Valdinoci, publi\'e dans Annales de l'Institut Henri Poincar\'e Analyse Non Lin\'eaire---voir \cite{BLV16}.
	\item \emph{On nonlocal minimal graphs}, co-auteur avec M. Cozzi, en cours de pr\'eparation.
	\item \emph{Bernstein-Moser-type results for nonlocal minimal graphs}, co-auteur avec M. Cozzi et A. Farina, soumis---voir \cite{CFL18}.
	\item Une version partielle et pr\'eliminaire de l'article \emph{A free boundary problem: superposition of nonlocal energy plus classical perimeter}, co-auteur avec S. Dipierro et E. Valdinoci, en cours de pr\'eparation.
	\item \emph{The Phillip Island penguin parade (a mathematical treatment)}, co-auteur avec S. Dipierro, P. Miraglio et E. Valdinoci, publi\'e dans ANZIAM Journal---voir \cite{DMLV16}.
\end{enumerate}

Les annexes contiennent des r\'esultats auxiliaires qui ont \'et\'e exploit\'es tout au long de la th\`ese.

\section{Une pr\'esentation plus d\'etaill\'ee}

Nous passons maintenant \`a une description d\'etaill\'ee du contenu et des principaux r\'esultats de cette th\`ese. Nous observons que chaque sujet a sa propre pr\'esentation, plus approfondie, au d\'ebut du chapitre correspondant. De plus, chaque chapitre a sa propre table des mati\`eres, pour aider le lecteur \`a naviguer entre les sections.

\subsection{Ensembles de p\'erim\`etre fractionnaire (localement) fini}

Le p\'erim\`etre $s$-fractionnaire et ses minimiseurs, les ensembles $s$-minimaux, ont \'et\'e introduits dans \cite{CRS10} en 2010, principalement motiv\'es par des applications aux probl\`emes de transition de phase en pr\`esence d'interactions \`a longue port\'ee.
Au cours des ann\'ees suivantes, ils ont suscit\'e un vif int\'er\^et, notamment en ce qui concerne la th\'eorie de la r\'egularit\'e et le comportement qualitatif des fronti\`eres des ensembles $s$-minimaux, qui sont les soi-disant surfaces minimales non locales. Nous invitons le lecteur int\'eress\'e \`a consulter \cite{V13} et~\cite[Chapter~6]{bucval}
pour une introduction, et \`a l'\'etude \cite{DV18} pour quelques d\'eveloppements r\'ecents.

En particulier, nous mentionnons que, m\^eme si la recherche de la r\'egularit\'e optimale des surfaces minimales non locales reste un probl\`eme ouvert et engageant, il est connu que les surfaces minimales non locales sont $(n-1)$-rectifiables.
Plus pr\'ecis\'ement, elles sont lisses,
sauf \'eventuellement pour un ensemble singulier de dimension de Hausdorff au plus \'egal \`a $n-3$ (voir \cite{CRS10}, \cite{SV13} et \cite{FV17}).
En cons\'equence, un ensemble $s$-minimal a p\'erim\`etre (au sens de De Giorgi et Caccioppoli) localement fini---et en fait, des estimations uniformes du p\'erim\`etre (classique) des ensembles $s$-minimaux sont disponibles (voir \cite{CSV16}).


\smallskip

D'autre part, la fronti\`ere d'un ensemble g\'en\'erique $E$ ayant $s$-p\'erim\`etre fini peut \^etre tr\`es irr\'eguli\`ere et peut m\^eme \^etre ``nulle part rectifiable'', comme dans le cas du flocon de neige de von Koch.

En fait, le $s$-p\'erim\`etre peut \^etre utilis\'e (en suivant l'article fondateur \cite{Visintin}) pour d\'efinir une ``dimension fractale'' pour la fronti\`ere, compris au sens de la th\'eorie de la mesure,
\begin{equation*}
\partial^-E:=\{x\in\R^n\,|\,0<|E\cap B_r(x)|<\omega_nr^n\textrm{ pour chaque }r>0\},
\end{equation*}
d'un ensemble $E\subseteq\R^n$.

Avant de continuer, nous rappelons la d\'efinition du $s$-p\'erim\`etre.
\'Etant donn\'e un param\`etre fractionnaire $s\in(0,1)$, nous d\'efinissons l'interaction
\begin{equation*}
\mathcal L_s(A,B):=\int_A\int_B\frac{1}{|x-y|^{n+s}}\,dx\,dy,
\end{equation*}
pour chaque couple d'ensembles disjoints $A,\,B\subseteq\mathbb R^n$.
Alors, le \emph{$s$-p\'erim\`etre} d'un ensemble $E\subseteq\mathbb R^n$ dans un ensemble ouvert $\Omega\subseteq\R^n$ est d\'efini comme
\begin{equation*}
\Per_s(E,\Omega):=\mathcal L_s(E\cap\Omega,\Co E\cap\Omega)+
\mathcal L_s(E\cap\Omega,\Co E\setminus\Omega)+
\mathcal L_s(E\setminus\Omega,\Co E\cap\Omega).
\end{equation*}
Nous \'ecrivons simplement $\Per_s(E):=\Per_s(E,\R^n)$.

On dit qu'un ensemble $E\subseteq\R^n$ a \emph{$s$-p\'erim\`etre localement fini} dans un ensemble ouvert $\Omega\subseteq\R^n$ si
\begin{equation*}
\Per_s(E,\Omega')<\infty\qquad\textrm{pour chaque ensemble ouvert }\Omega'\Subset\Omega.
\end{equation*}

Nous observons que nous pouvons r\'e\'ecrire le $s$-p\'erim\`etre comme
\begin{equation}\label{CH_fr_CH:1:func_per_form}
\Per_s(E,\Omega)=\frac{1}{2}\iint_{\R^{2n}\setminus(\Co\Omega)^2}\frac{|\chi_E(x)-\chi_E(y)|}{|x-y|^{n+s}}\,dx\,dy.
\end{equation}

La formule \eqref{CH_fr_CH:1:func_per_form} montre que le p\'erim\`etre fractionnaire est, approximativement, la $\Omega$-contribution
\`a la seminorme $W^{s,1}$ de la fonction caract\'eristique $\chi_E$.

Cette fonctionnelle est non locale, au sens qu'il faut conna\^itre l'ensemble $E$
dans tout $\R^n$, m\^eme pour calculer son $s$-p\'erim\`etre dans un petit domaine born\'e $\Omega$
(contrairement \`a ce qui se passe avec le p\'erim\`etre classique ou la
mesure $\Ha^{n-1}$, qui sont des fonctionnelles locales).
En plus, le $s$-p\'erim\`etre est ``fractionnaire'', dans le sens o\`u la seminorme $W^{s,1}$ mesure un ordre de r\'egularit\'e fractionnaire.

Nous observons que nous pouvons diviser le $s$-p\'erim\`etre comme
\bgs{
	\Per_s(E,\Omega)=\Per_s^L(E,\Omega)+\Per_s^{NL}(E,\Omega),
}
o\`u
\bgs{
	\Per_s^L(E,\Omega):=\Ll_s(E\cap\Omega,\Co E\cap\Omega)=\frac{1}{2}[\chi_E]_{W^{s,1}(\Omega)}
}
peut \^etre consid\'er\'e comme la ``partie locale'' du p\'erim\`etre fractionnaire, et
\bgs{
	\Per_s^{NL}(E,\Omega)&:=\Ll_s(E\cap\Omega,\Co E\setminus\Omega)+
	\Ll_s(E\setminus\Omega,\Co E\cap\Omega)\\
	&
	=\int_\Omega\int_{\Co\Omega}\frac{|\chi_E(x)-\chi_E(y)|}{|x-y|^{n+s}}\,dx\,dy,
}
qui peut \^etre consid\'er\'e comme la ``partie non locale''.

\subsubsection{\textbf{Fronti\`eres fractales}}\label{CH_fr_CH:0:Fractal_sec_Intro}



En 1991, dans l'article \cite{Visintin} l'auteur a sugg\'er\'e d'utiliser le param\`etre $s$ de la seminorme fractionnaire
$[\chi_E]_{W^{s,1}(\Omega)}$ (et de plus g\'en\'erales familles continues de fonctionnelles satisfaisant des opportunes formules de la co-aire g\'en\'eralis\'ees) comme un moyen de mesurer la codimension de la fronti\`ere comprise au sens de la th\'eorie de la mesure,
$\partial^-E$, d'un ensemble $E$ dans $\Omega$.
Il a prouv\'e que la dimension fractale obtenue de cette mani\`ere,
\begin{equation*}
\Dim_F(\partial^-E,\Omega):=n-\sup\{s\in(0,1)\,|\,[\chi_E]_{W^{s,1}(\Omega)}<\infty\},
\end{equation*}
est inf\'erieure ou \'egale \`a la dimension (sup\'erieure) de Minkowski.

La relation entre la dimension de Minkowski de la fronti\`ere d'un ensemble $E$
et la r\'egularit\'e fractionnaire (dans le sens des espaces de Besov) de la fonction caract\'eristique $\chi_E$ a \'et\'e \'etudi\'e aussi dans \cite{Sickel}, en 1999.
En particulier---voir \cite[Remark 3.10]{Sickel}---l'auteur a prouv\'e que
la dimension $\Dim_F$ du flocon de neige de von Koch $S$ co\"{i}ncide avec sa dimension de Minkowski, en exploitant le fait que
$S$ est un domaine de John.

La r\'egularit\'e de Sobolev d'une fonction caract\'eristique $\chi_E$ a \'et\'e approfondie dans \cite{FarRog}, en 2013,
o\`u les auteurs consid\`erent le cas dans lequel l'ensemble $E$ est une quasiball. Comme le flocon de neige de von Koch $S$ est un exemple typique de quasiball, les auteurs ont pu prouver que la dimension  $\Dim_F$ de $S$ co\"{i}ncide avec sa dimension de Minkowski.

\medskip

Dans le Chapitre \ref{CH_Fractals}, nous calculons la dimension  $\Dim_F$ du flocon de neige de von Koch $S$ de mani\`ere \'el\'ementaire, en utilisant uniquement l'invariance par roto-translation et la propri\'et\'e d'\'echelle du $s$-p\'erim\`etre, et la ``auto-similarit\'e'' de $S$.
Plus pr\'ecis\'ement, nous montrons que
\bgs{
	\Per_s(S)<\infty,\qquad\forall\,s\in\left(
	0,2-\frac{\log 4}{\log 3}\right),
}
et
\bgs{
	\Per_s(S)=\infty,\qquad\forall\,s\in\left[2-\frac{\log 4}{\log 3},1\right).
}
La d\'emonstration peut \^etre \'etendue de mani\`ere naturelle \`a tous les ensembles qui peuvent \^etre d\'efinis de mani\`ere r\'ecursive similaire \`a celle du flocon de von Koch. En cons\'equence, nous calculons la dimension  $\Dim_F$ de tous ces ensembles, sans avoir \`a les obliger \`a \^etre des domaines de John ou des quasiballs.

De plus, nous montrons que nous pouvons facilement obtenir beaucoup d'ensembles de ce type en modifiant de mani\`ere appropri\'ee des fractales auto-similaires bien connues, comme le flocon de neige de von Koch, le triangle de Sierpinski et l'\'eponge de Menger. Un exemple est illustr\'e dans la Figure \ref{CH_fr_CH:1:buffo_triang}.

\begin{figure}[htbp]
	\begin{center}
		\includegraphics[width=60mm]{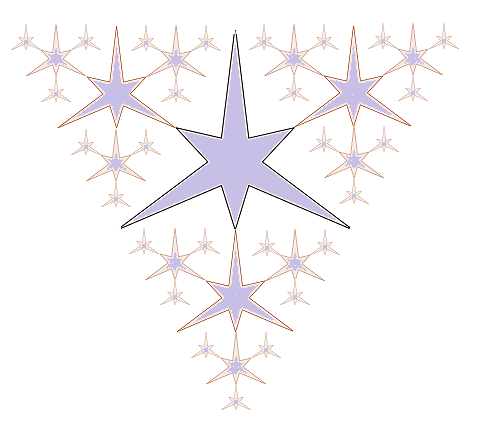}
		\caption{{\it Exemple d'un ensemble ``fractal'' construit en exploitant la structure du triangle de Sierpinski (visible \`a la quatri\`eme \'etape it\'erative).}}
		\label{CH_fr_CH:1:buffo_triang}
	\end{center}
\end{figure}

\subsubsection{\textbf{Asymptotique $s\to1^-$}}\label{CH_fr_CH:0:asympto1_subsec}

La discussion pr\'ec\'edente montre que
le $s$-p\'erim\`etre d'un ensemble $E$ ayant fronti\`ere irr\'eguli\`ere, \'eventuellement fractale, peut \^etre fini
pour $s$ sous un certain seuil, $s<\sigma$, et infini pour $s\in(\sigma,1)$.
D'autre part, il est bien connu que les ensembles avec une fronti\`ere r\'eguli\`ere ont $s$-p\'erim\`etre fini pour chaque $s$ et leur $s$-p\'erim\`etre converge, lorsque $s$ tend vers 1, au p\'erim\`etre classique, \`a la fois au sens classique (voir, par exemple, \cite{uniform}) et au sens de la $\Gamma$-convergence (voir, par exemple, \cite{Gamma}
et aussi \cite{Ponce} pour des r\'esultats connexes).

Dans le Chapitre \ref{CH_Fractals}, nous exploitons \cite[Theorem 1]{Davila} pour prouver une version optimale de cette propri\'et\'e asymptotique pour un ensemble $E$ ayant p\'erim\`etre classique fini dans un ensemble ouvert born\'e avec fronti\`ere de classe de Lipschitz. Plus pr\'ecis\'ement, nous prouvons que, si $E$ a p\'erim\`etre classique fini dans un voisinage de $\Omega$, alors
\begin{equation*}
\lim_{s\to1}(1-s)\Per_s(E,\Omega)=\omega_{n-1}\Per(E,\overline\Omega).
\end{equation*}

Nous observons que nous baissons la r\'egularit\'e demand\'ee dans \cite{uniform}, o\`u les auteurs ont exig\'e que le fronti\`ere $\partial E$ soit $C^{1,\alpha}$,
\`a la r\'egularit\'e optimale (demandent \`a $E$ seulement d'avoir p\'erim\`etre fini). En plus, nous n'avons pas \`a demander \`a $E$ de croiser $\partial\Omega$ ``transversalement'', c'est-\`a-dire que nous n'avons pas besoin
que
\begin{equation*}
\Ha^{n-1}(\partial^*E\cap\partial\Omega)=0,
\end{equation*}
o\`u $\partial^*E$ d\'enote la fronti\`ere r\'eduite de $E$.

En effet, nous prouvons que la partie non locale du $s$-p\'erim\`etre converge au p\'erim\`etre sur la fronti\`ere de $\Omega$, c'est-\`a-dire que nous prouvons que
\begin{equation*}
\lim_{s\to1}(1-s)\Per_s^{NL}(E,\Omega)=\omega_{n-1}\Ha^{n-1}(\partial^*E\cap\partial\Omega),
\end{equation*}
qui est, \`a la connaissance de l'auteur, un nouveau r\'esultat.

\subsubsection{\textbf{Approximation par ensembles ouverts lisses}}\label{CH_fr_CH:0:Appro_subsection}

Comme nous avons observ\'e dans la Section \ref{CH_fr_CH:0:Fractal_sec_Intro}, les ensembles ayant p\'erim\`etre fractionnaire fini peuvent avoir une fronti\`ere tr\`es rugueuse, qui peut en effet \^etre une fractale nulle part rectifiable (comme le flocon de neige de von Koch).\\
Cela repr\'esente une diff\'erence importante entre le p\'erim\`etre fractionnaire et le p\'erim\`etre classique, car les ensembles de Caccioppoli ont une partie "grande" de la fronti\`ere, dite fronti\`ere r\'eduite, qui est $(n-1)$-rectifiable (d'apr\`es le Th\'eor\`eme de structure de De Giorgi).

En tout cas, nous prouvons dans la premi\`ere partie du Chapitre \ref{CH_Appro_Min} qu'un ensemble a p\'erim\`etre fractionnaire (localement)  fini si et seulement si il peut \^etre approch\'e (de mani\`ere appropri\'ee) par des ensembles ouverts lisses.
Plus pr\'ecisément, nous prouvons ce qui suit:
\begin{frtheorem}\label{CH_fr_CH:0:density_smooth_thm}
	Soit $\Omega\subseteq\R^n$ un ensemble ouvert. Un ensemble $E\subseteq\R^n$ a $s$-p\'erim\`etre localement fini dans $\Omega$ si et seulement s'il existe une suite $E_h\subseteq\R^n$ de ensembles ouverts ayant fronti\`ere lisse et $\eps_h\to0^+$ tels que
	\begin{equation*}\begin{split}
	& (i)\quad E_h\xrightarrow{loc}E,\qquad\sup_{h\in\mathbb N}\Per_s(E_h,\Omega')<\infty\quad\textrm{pour chaque }\Omega'\Subset\Omega,\\
	& (ii)\quad\lim_{h\to\infty}\Per_s(E_h,\Omega')=\Per_s(E,\Omega')\quad\textrm{pour chaque }\Omega'\Subset\Omega,\\
	& (iii)\quad\partial E_h\subseteq N_{\eps_h}(\partial E).
	\end{split}
	\end{equation*}
	En outre, si $\Omega=\R^n$ et l'ensemble $E$ est tel que $|E|<\infty$ et $\Per_s(E)<\infty$, alors
	\begin{equation*}
	|E_h\Delta E|\to0,\qquad\quad\qquad \lim_{h\to\infty}\Per_s(E_h)=\Per_s(E),
	\end{equation*}
	et nous pouvons exiger que chaque $E_h$ soit born\'e (au lieu de demander $(iii)$).
\end{frtheorem}

Ci-dessus, $N_\delta(\partial E)$ d\'enote le $\delta$-voisinage tubulaire de $\partial E$.

Un tel r\'esultat est bien connu pour les ensembles de Caccioppoli (voir, par exemple,  \cite{Maggi}) et en effet, cette propri\'et\'e de densit\'e peut \^etre utilis\'ee pour d\'efinir la fonctionnelle de p\'erim\`etre (classique) comme \'etant la relaxation---par rapport \`a la convergence $L^1_{\loc}$---de la mesure $\Ha^{n-1}$ des fronti\`eres des ensembles ouverts lisses, c'est-\`a-dire
\begin{equation}\label{CH_fr_CH:0:liminfclassical}\begin{split}\Per(E,\Omega)=\inf\Big\{\liminf_{k\to\infty}\Ha^{n-1}(\partial E_h&\cap\Omega)\,\big|\,
E_h\subseteq\R^n\textrm{ ouvert ayant fronti\`ere}\\
&
\textrm{lisse, tel que }E_h\xrightarrow{loc}E\Big\}.
\end{split}
\end{equation}

Il est int\'eressant de noter que, dans \cite{DV18}, les auteurs ont prouv\'e, en exploitant le th\'eor\`eme de la divergence, que si $E\subseteq \R^n$ est un ensemble ouvert born\'e avec fronti\`ere lisse, alors
\begin{equation}\label{CH_fr_CH:0:liminf_fract_formula}
\Per_s(E)=c_{n,s}\int_{\partial E}\int_{\partial E}
\frac{2-|\nu_E(x)-\nu_E(y)|^2}{|x-y|^{n+s-2}}d\Ha^{n-1}_xd\Ha^{n-1}_y,
\end{equation}
o\`u $\nu_E$ d\'enote la normale externe de
$E$ et
\begin{equation*}
c_{n,s}:=\frac{1}{2s(n+s-2)}.
\end{equation*}
En exploitant la formule \eqref{CH_fr_CH:0:liminf_fract_formula}, la semicontinuit\'e inf\'erieure
du $s$-p\'erim\`etre et le Th\'eor\`eme \ref{CH_fr_CH:0:density_smooth_thm}, nous trouvons que, si
$E\subseteq\R^n$ est tel que $|E|<\infty$, alors
\begin{equation*}\begin{split}
\Per_s(E)&=\inf\bigg\{\liminf_{h\to\infty}c_{n,s}\int_{\partial E_h}\int_{\partial E_h}
\frac{2-|\nu_{E_h}(x)-\nu_{E_h}(y)|^2}{|x-y|^{n+s-2}}d\Ha^{n-1}_xd\Ha^{n-1}_y\,\big|\\
&
E_h\subseteq\R^n\textrm{ ensemble ouvert born\'e ayant fronti\`ere lisse, tel que }E_h\xrightarrow{loc}E\bigg\}.
\end{split}
\end{equation*}
Cela peut \^etre consid\'er\'e comme un analogue de \eqref{CH_fr_CH:0:liminfclassical} dans le cadre fractionnaire.

Nous mentionnons \'egalement que dans la Section \ref{CH:4:Appro_Section} nous allons prouver qu'un sous-graphe ayant $s$-p\'erim\`etre localement fini dans un cylindre  $\Omega \times\R$ peut \^etre approch\'e par les sous-graphes de fonctions lisses---et pas seulement par des ensembles ouverts lisses arbitraires.

\subsection{Surfaces minimales non locales}\label{CH_fr_CH:0:NMS_Sec}

La deuxi\`eme partie du Chapitre \ref{CH_Appro_Min} concerne les ensembles minimisant le p\'erim\`etre fractionnaire.
Les fronti\`eres de ces minimiseurs sont souvent appel\'es surfaces minimales non locales et apparaissent naturellement comme interfaces limites des mod\`eles de transition de phase \`a interaction \`a longue port\`ee. En particulier, dans les r\'egimes o\`u l'interaction \`a longue port\'ee est dominante, la fonctionnelle de Allen-Cahn non locale $\Gamma$-converge
au p\'erim\`etre fractionnaire (voir, par exemple, \cite{SV12})
et les interfaces minimales de l'\'equation de Allen-Cahn correspondante approchent localement de mani\`ere uniforme les surfaces minimales non locales (voir, par exemple, \cite{SV14}).

Nous rappelons maintenant la d\'efinition des ensembles minimisants introduite dans \cite{CRS10}.
\begin{frdefn}\label{CH_fr_CH:0:sMini_def}
	Soit $\Omega\subseteq\R^n$ un ensemble ouvert et soit $s\in(0,1)$. On dit qu'un ensemble $E\subseteq\R^n$ est \emph{$s$-minimal} dans $\Omega$ si $\Per_s(E,\Omega)<\infty$ et
	\[
	\Per_s(E,\Omega)\leq\Per_s(F,\Omega)\quad\mbox{ pour chaque }F\subseteq\R^n\mbox{ tel que }F\setminus\Omega=E\setminus\Omega.
	\]
\end{frdefn}
Parmi les nombreux r\'esultats, dans \cite{CRS10} les auteurs ont prouv\'e que,
si $\Omega\subseteq\R^n$ est un ensemble ouvert born\'e ayant fronti\`ere Lipschitz, alors pour chaque ensemble fix\'e
$E_0\subseteq\Co\Omega$ il existe un ensemble $E\subseteq\R^n$ qui est $s$-minimal dans $\Omega$ et tel que $E\setminus\Omega=E_0$.
L'ensemble $E_0$ est parfois appel\'e \emph{donné extérieur} et l'ensemble $E$ est dit \^etre $s$-minimal dans $\Omega$ par rapport \`a la donn\'ee ext\'erieure $E_0$.

Nous \'etendons le r\'esultat d'existence susmentionn\'e en prouvant que, dans un ensemble ouvert g\'en\'erique $\Omega$, il existe un ensemble $s$-minimal par rapport  \`a une certaine donn\'ee ext\'erieure $E_0\subseteq \Co\Omega$ fix\'ee, si et seulement si il existe un concurrent ayant $s$-p\'erim\`etre fini dans $\Omega$. Plus pr\'ecis\'ement:

\begin{frtheorem}\label{CH_fr_CH:0:Global_Existence}
	Soit $s\in(0,1)$, soit $\Omega\subseteq\R^n$ un ensemble ouvert et soit $E_0\subseteq\Co\Omega$. Alors, il existe un ensemble $E\subseteq\R^n$ qui est $s$-minimal dans $\Omega$ et tel que $E\setminus\Omega=E_0$, si et seulement si il existe un ensemble $F\subseteq\R^n$ tel que $F\setminus\Omega=E_0$ et $\Per_s(F,\Omega)<\infty$.
\end{frtheorem}

En cons\'equence, nous observons que, si $\Per_s(\Omega)<\infty$, alors il existe toujours un ensemble $s$-minimal par rapport \`a la donn\'ee ext\'erieure $E_0$, pour chaque ensemble $E_0\subseteq\Co\Omega$.

\medskip

Portons maintenant l'attention sur le cas dans lequel le domaine de minimisation n'est pas born\'e. Dans cette situation, il convient d'introduire la notion de minimiseur local.
\begin{frdefn}
	Soit $\Omega\subseteq\R^n$ un ensemble ouvert et soit $s\in(0,1)$. On dit qu'un ensemble $E\subseteq\R^n$ est \emph{localement $s$-minimal} dans $\Omega$ si $E$ est $s$-minimal dans chaque ensemble ouvert $\Omega'\Subset\Omega$.
\end{frdefn}

Notez en particulier que nous demandons \`a $E$ seulement d'avoir $s$-p\'erim\`etre localement fini dans $\Omega$ et pas d'avoir $s$-p\'erim\`etre fini dans tout le domaine. En effet, la principale raison de l'introduction des ensembles localement $s$-minimaux est donn\'ee par le fait qu'en g\'en\'eral, le $s$-p\'erim\`etre d'un ensemble n'est pas fini dans les domaines non born\'es.

Nous avons vu dans le Th\'eor\`eme \ref{CH_fr_CH:0:Global_Existence} que le seul obstacle \`a l'existence d'un ensemble $s$-minimal, par rapport  \`a une certaine donn\'ee ext\'erieure $E_0\subseteq \Co\Omega$ fix\'ee, est l'existence d'un concurrent ayant $s$-p\'erim\`etre fini. D'autre part, nous prouvons qu'un ensemble localement $s$-minimal existe toujours, peu importe ce que le domaine $\Omega$ et la donn\'ee ext\'erieure sont.

\begin{frtheorem}\label{CH_fr_CH:0:Local_Existence}
	Soit $s\in(0,1)$, soit $\Omega\subseteq\R^n$ un ensemble ouvert et soit $E_0\subseteq\Co\Omega$. Alors, il existe un ensemble $E\subseteq\R^n$ qui est localement $s$-minimal dans $\Omega$ et tel que $E\setminus\Omega=E_0$.
\end{frtheorem}
Quand $\Omega$ est un ensemble ouvert born\'e ayant fronti\`ere Lipschitz, nous montrons que les deux notions de minimiseur co\"incident. C'est-\`a-dire, si $\Omega\subseteq\R^n$ est un ensemble ouvert born\'e ayant fronti\`ere Lipschitz et $E\subseteq\R^n$, alors
\[
E\mbox{ est $s$-minimal dans }\Omega\quad\Longleftrightarrow\quad E\mbox{ est localement $s$-minimal dans }\Omega.
\]

Cependant, nous observons que cela n'est pas vrai dans un ensemble ouvert $\Omega$ arbitraire, car un ensemble $s$-minimal---au sens de la D\'efinition \ref{CH_fr_CH:0:sMini_def}---peut ne pas exister.

A titre d'exemple, nous consid\'erons la situation dans laquelle le domaine de minimisation est le cylindre
\[
\Omega^\infty:=\Omega\times\R\subseteq\R^{n+1},
\]
o\`u $\Omega\subseteq\R^n$ est un ensemble ouvert born\'e ayant fronti\`ere r\'eguli\`ere. Nous nous int\'eressons au cas o\`u la donn\'ee ext\'erieure est le sous-graphe d'une fonction mesurable $\varphi:\R^n\to\R$. C'est-\`a-dire, nous consid\'erons le sous-graphe
\[
\Sg(\varphi):=\left\{(x,x_{n+1})\in\R^{n+1}\,|\,x_{n+1}<\varphi(x)\right\},
\]
et nous voulons trouver un ensemble $E\subseteq\R^{n+1}$ qui minimise---dans un certain sens---le $s$-p\'erim\`etre dans le cylindre $\Omega^\infty$, par rapport \`a la donn\'ee ext\'erieure $E\setminus\Omega^\infty=\Sg(\varphi)\setminus\Omega^\infty$.

Une motivation pour consid\'erer un tel probl\`eme de minimisation est donn\'ee par le r\'ecent article \cite{graph}, o\`u les auteurs ont prouv\'e que si un tel ensemble de minimisation $E$ existe---et si $\varphi$ est une fonction continue---alors $E$ est en fait un sous-graphe global. Plus pr\'ecis\'ement, il existe une fonction $u:\R^n\to\R$, telle que $u=\varphi$ dans $\R^n\setminus\overline{\Omega}$ et $u\in C(\overline{\Omega})$, et telle que
\[
E=\Sg(u).
\]

On voit facilement que si une fonction $u:\R^n\to\R$ est assez r\'eguli\`ere dans $\Omega$, par exemple, si $u\in BV(\Omega)\cap L^\infty(\Omega)$, alors la partie locale du $s$-p\'erim\`etre du sous-graphe de $u$ est finie,
\[
\Per_s^L(\Sg(u),\Omega^\infty)<\infty.
\]
D'autre part, la partie non locale du $s$-p\'erim\`etre, en g\'en\'eral, est infinie, m\^eme pour des fonctions tr\`es r\'eguli\`eres $u$. En effet, nous prouvons que si $u\in L^\infty(\R^n)$, alors
\[
\Per_s^{NL}(\Sg(u),\Omega^\infty)=\infty.
\]

Une premi\`ere cons\'equence de cette observation---et de l'estimation a priori sur la ``variation verticale'' d'un ensemble de minimisation fourni par \cite[Lemma 3.3]{graph}---est le fait que, si $\varphi\in C(\R^n)\cap L^\infty(\R^n)$, alors il ne peut pas exister un ensemble $E$ qui est $s$-minimal dans $\Omega^\infty$---au sens de la D\'efinition \ref{CH_fr_CH:0:sMini_def}---par rapport \`a la donn\'ee ext\'erieure $\Sg(\varphi)\setminus\Omega^\infty$.

Toutefois, le Th\'eor\`eme \ref{CH_fr_CH:0:Local_Existence} garantit l'existence d'un ensemble $E\subseteq\R^{n+1}$ qui est localement $s$-minimal dans $\Omega^\infty$ et tel que $E\setminus\Omega^\infty=\Sg(\varphi)\setminus\Omega^\infty$. Donc, Th\'eor\`eme \ref{CH_fr_CH:0:Local_Existence} et \cite[Theorem 1.1]{graph} impliquent ensemble l'existence de sous-graphes minimisant (localement) le $s$-p\'erim\`etre, c'est-\`a-dire, des surfaces minimales non locales non param\'etriques.

Une deuxi\`eme cons\'equence consiste dans le fait que nous ne pouvons pas d\'efinir une version fractionnaire na\^ive de la fonctionnelle d'aire classique comme
\[
\mathscr{A}_s(u,\Omega):=\Per_s(\Sg(u),\Omega^\infty),
\]
puisque cela serait infinie m\^eme pour une fonction $u\in C^\infty_c(\R^n)$. Au Chapitre \ref{CH_Nonparametric} nous allons \'eviter ce probl\`eme en introduisant un cadre fonctionnel appropri\'e pour travailler avec des sous-graphes.

\subsection{Effets de stickiness pour les petits valeurs de $s$}

Le Chapitre \ref{Asympto0_CH_label} est consacr\'e \`a l'\'etude des ensembles $s$-minimaux dans des r\'egimes hautement non locaux, c'est-\`a-dire dans le cas o\`u le param\`etre fractionnaire $s\in(0,1)$ est tr\`es petit. Nous prouvons que, dans cette situation, le comportement des ensembles $s$-minimaux, d'une certaine mani\`ere, d\'eg\'en\`ere.

\medskip

Rappelons d'abord quelques r\'esultats connus concernant l'asymptotique $s\to1^-$.\\
Nous avons d\'ej\`a observ\'e dans la Section \ref{CH_fr_CH:0:asympto1_subsec} que le $s$-p\'erim\`etre converge vers le p\'erim\`etre classique lorsque $s\to1^-$.
De plus, quand $s\to1^-$, les ensembles $s$-minimaux convergent vers les minimiseurs du p\'erim\`etre classique, \`a la fois au ``sens uniforme'' (voir \cite{uniform, regularity}) et au sens de la
$\Gamma$-convergence (voir \cite{Gamma}).
En cons\'equence, on peut prouver (voir \cite{regularity}) que quand $s$ est suffisamment proche de 1, les surfaces minimales non locales ont la m\^eme r\'egularit\'e des surfaces minimales classiques.
Voir aussi \cite{DV18} pour une \'etude r\'ecente et assez compl\`ete des propri\'et\'es des ensembles $s$-minimaux lorsque $s$ est proche de 1. 

De plus, nous observons que la courbure moyenne fractionnaire converge \'egalement, comme $s\to1^-$, vers sa contrepartie classique. Pour \^etre plus pr\'ecis, rappelons d'abord que la courbure moyenne $s$-fractionnaire d'un ensemble $E$ en un point $q\in\partial E$ est d\'efinie comme l'int\'egrale au sens de la valeur principale
\[\I_s[E](q):=\PV\int_{\R^n}\frac{\chi_{\Co E}(y)-\chi_E(y)}{|y-q|^{n+s}}\,dy,\]
c'est-\`a-dire
\[\I_s[E](q):=\lim_{\varrho\to0^+}\I_s^\varrho[E](q),\qquad\textrm{o\`u}\qquad
\I_s^\varrho[E](q):=\int_{
	\Co B_\varrho(q)}\frac{\chi_{\Co E}(y)-\chi_E(y)}{|y-q|^{n+s}}\,dy.\] 
Remarquons qu'il est en effet n\'ecessaire d'interpr\'eter l'int\'egrale ci-dessus
au sens de la valeur principale, puisque l'int\'egrande est singuli\`ere et non int\'egrable dans un voisinage de $q$. D'autre part, s'il ya suffisamment d'annulation entre $E$ et $\Co E$ dans un voisinage de $q$---par exemple, si $\partial E$ est de classe $C^2$ autour de $q$---alors l'int\'egrale est bien d\'efinie au sens de la valeur principale.

La courbure moyenne fractionnaire a \'et\'e introduite dans \cite{CRS10}, o\`u les auteurs ont montr\'e qu'elle est l'op\'erateur d'Euler-Lagrange apparaissant dans la minimisation du $s$-p\'erim\`etre. En effet, si $E\subseteq\R^n$ est $s$-minimal dans un ensemble ouvert $\Omega$, alors 
\[
\I_s[E]=0\quad\mbox{sur } \partial E,
\]
dans un sens de viscosit\'e appropri\'e---pour plus de d\'etails voir, par exemple, l'Annexe \ref{CH:3:brr2}.

Il est connu (voir, par exemple, \cite[Theorem 12]{Abaty} et \cite{regularity}) que si $E\subseteq \Rn$ est un ensemble ayant fronti\`ere $C^2$, et $n\ge2$, alors pour tous $x\in \partial E$ on a que
\[ \lim_{s \to 1} (1-s)\I_s[E] (x) = \varpi_{n-1}H[E](x).\]
Ci-dessus $H$ d\'enote la courbure moyenne classique de $E$ au point $x$---selon la convention que nous prenons $H$ tel que la courbure de la boule est une quantit\'e positive---et
\[
\varpi_k:=\Ha^{k-1}(\{x\in\R^k\,|\,|x|=1\}),
\]
pour chaque $k\ge1$. Laissez-nous \'egalement d\'efinir $\varpi_0:=0$. Nous observons que pour $n=1$, nous avons
\[  \lim_{s \to 1} (1-s)\I_s[E] (x) = 0,\]
ce qui est compatible avec la notation $\varpi_0=0$---voir aussi Remarque \ref{CH:3:nuno}.

\medskip

Lorsque $s\to 0^+$, les asymptotiques sont plus compliqu\'es et pr\'esentent un comportement surprenant.  Cela est d\^u au fait que quand $s$ devient plus petit, la contribution non locale au compteur du $s$-p\'erim\`etre devient de plus en plus importante, tandis que la contribution locale perd de son influence. Quelques r\'esultats pr\'ecis \`a cet \'egard ont \'et\'e obtenus dans \cite{DFPV13}. L\`a, pour encoder le comportement \`a l'infini d'un ensemble, les auteurs ont introduit la quantit\'e
\[
\alpha(E)=\lim_{s\to 0^+} s\int_{\Co B_1} \frac{\chi_E(y)}{|y|^{n+s}}\, dy,
\]
qui appara\^it naturellement quand on regarde l'asymptotique pour $s\to0^+$ du p\'erim\`etre fractionnaire. En fait, dans \cite{DFPV13} les auteurs ont prouv\'e que, si $\Omega$ est un ensemble ouvert born\'e ayant fronti\`ere $C^{1,\gamma}$, pour quelque $\gamma\in (0,1]$, $E \subseteq \Rn$ a $s_0$-p\'erim\`etre fini dans $\Omega$, pour un certain $s_0\in (0,1)$, et $\alpha(E)$ existe, alors
\bgs{
	\lim_{s\to 0^+} s\Per_s(E,\Omega)=  \alpha(\Co E) |E\cap \Omega| + \alpha(E) |\Co E \cap \Omega|.
}

D'autre part, le comportement asymptotique lorsque $s\to 0^+$ de la courbure moyenne fractionnaire est \'etudi\'e au Chapitre \ref{Asympto0_CH_label}
(voit aussi \cite{DV18} pour le cas particulier dans lequel l'ensemble $E$ est born\'e).
Tout d'abord, puisque la quantit\'e $\alpha(E)$ peut ne pas exister---voir \cite[Example 2.8 et 2.9]{DFPV13}---nous d\'efinissons
\[
\overline \alpha (E):= \limsup_{s\to 0^+} s\int_{\Co B_1} \frac{\chi_E(y)}{|y|^{n+s}}\, dy
\quad\mbox{et} \quad
\underline \alpha(E) := \liminf_{s\to 0^+} s\int_{\Co B_1} \frac{\chi_E(y)}{|y|^{n+s}}\, dy.
\]

Nous prouvons que, lorsque $s\to 0^+$, la courbure moyenne $s$-fractionnaire devient compl\`etement indiff\'erente \`a la g\'eom\'etrie locale de la fronti\`ere $\partial E$, et en effet la valeur limite ne d\'epend que du comportement \`a l'infini de l'ensemble $E$. Plus pr\'ecis\'ement, si $E\subseteq\Rn$ et $p\in\partial E$ est tel que $\partial E$ est $C^{1,\gamma}$ autour de $p$,
pour un certain $\gamma\in(0,1]$, alors
\eqlab{\label{CH_fr_CH:0:liminf_curvat}
	\liminf_{s\to0^+} s\,\I_s[E](p) =\varpi_n -2 \overline \alpha(E),
}
et
\[
\limsup_{s\to0^+} s\,\I_s[E](p) =\varpi_n-2 \underline\alpha(E).
\]
Nous remarquons en particulier que si $E$ est born\'e, alors $\alpha(E)$ existe et $\alpha(E)=0$. Donc, si $E\subseteq\R^n$ est un ensemble ouvert born\'e ayant fronti\`ere $C^{1,\gamma}$, l'asymptotique est simplement
\[
\lim_{s\to0^+}s\,\I_s[E](p)=\varpi_n,
\]
pour chaque $p\in\partial E$---voir aussi \cite[Appendix~B]{DV18}.

Dans la Section \ref{CH:3:sectexamples} nous calculons la contribution \`a l'infini $\alpha(E)$ de quelques ensembles. Pour avoir quelques exemples en t\^ete, nous citons ici les cas suivants:
\begin{itemize}
	\item soit $S\subseteq\s^{n-1}$ et consid\`ere le c\^one
	\[
	C:=\{t\sigma\in\R^n\,|\,t\geq0,\,\sigma\in S\}.
	\]
	Alors, $\alpha(C)=\Ha^{n-1}(S)$.
	\item Si $u\in L^\infty(\R^n)$, alors
	$\alpha(\Sg(u))=\varpi_{n+1}/2$. Plus en g\'en\'eral, si $u:\R^n\to\R$ est telle que
	\[
	\lim_{|x|\to\infty}\frac{|u(x)|}{|x|}=0,
	\]
	alors $\alpha(\Sg(u))=\varpi_{n+1}/2$.
	\item Soit $u:\R^n\to\R$ telle que $u(x)\leq-|x|^2$, pour chaque $x\in\R^n\setminus B_R$, pour un certain $R>0$. Alors $\alpha(\Sg(u))=0$.
\end{itemize}
Approximativement, \`a partir des exemples ci-dessus, nous voyons que $\alpha(E)$ ne d\'epend pas de la g\'eom\'etrie locale ni de la r\'egularit\'e de $E$, mais seulement de son comportement \`a l'infini.

\smallskip

Maintenant, nous observons que, lorsque $s \to 0^+$, les ensembles $s$-minimaux pr\'esentent un comportement plut\^ot inattendu.

Par exemple, en \cite[Theorem 1.3]{boundary} il est prouv\'e que si nous considérons le premier quadrant du plan comme donn\'ee ext\'erieure, alors, assez \'etonnamment, si $s$ est assez petit, l'ensemble $s$-minimal dans $B_1\subseteq\R^2$ est vide dans $B_1$.
Les principaux r\'esultats du Chapitre \ref{Asympto0_CH_label} s'inspirent de ce r\'esultat.

Heuristiquement, afin de g\'en\'eraliser \cite[Theorem 1.3]{boundary} nous voulons prouver que, si $\Omega\subseteq\R^n$ est un ensemble ouvert born\'e et connexe ayant fronti\`ere lisse et si nous fixons comme donn\'ee ext\'erieure un ensemble $E_0\subseteq\Co\Omega$ tel que $\overline{\alpha}(E_0)<\varpi_n/2$, alors il y a une contradiction entre l'\'equation d'Euler-Lagrange d'un ensemble $s$-minimal et l'asymptotique de la courbure moyenne $s$-fractionnaire pour $s\to0^+$.

Pour motiver pourquoi nous attendons une telle contradiction, nous observons que l'asymptotique \eqref{CH_fr_CH:0:liminf_curvat} semble sugg\'erer que, si $s$ est assez petit, alors un ensemble $s$-minimal $E$ ayant donn\'ee ext\'erieure $E_0$ et tel que $\partial E\cap\Omega\not=\emptyset$ devrait avoir un point $p\in\partial E\cap\Omega$ tel que $\I_s[E](p)>0$---qui contredirait l'\'equation d'Euler-Lagrange. Pour \'eviter une telle contradiction, nous conclurions alors que $\partial E=\emptyset$ in $\Omega$, c'est-\`a-dire que soit $E\cap\Omega=\Omega$ ou $E\cap\Omega=\emptyset$. 

Afin de transformer cette id\'ee en argument rigoureux, nous montrons d'abord que nous pouvons minorer la courbure moyenne fractionnaire, uniform\'ement par rapport au rayon d'une boule tangente \`a $E$ ext\'erieurement.
Plus pr\'ecis\'ement:
\begin{frtheorem}\label{CH_fr_CH:0:positive_curvature}
	Soit $\Omega\subseteq\Rn$ un ensemble ouvert born\'e. 
	Soit $E_0\subseteq\Co\Omega$ tel que
	\bgs{
		\overline \alpha(E_0)<\frac{\varpi_n}2,
	}
	et soit
	\[
	\beta=\beta(E_0):=\frac{\varpi_n-2\overline \alpha(E_0)}4.
	\]
	Nous d\'efinissons
	\bgs{
		\delta_s=\delta_s(E_0):=e^{-\frac{1}{s}\log \frac{\varpi_n+2\beta}{\varpi_n+\beta}},
	}
	pour chaque $s\in(0,1)$.
	Alors, il existe $s_0=s_0(E_0,\Omega)\in(0,\frac{1}{2}]$ tel que, si $E\subseteq\Rn$ est tel que $E\setminus\Omega=E_0$
	et $E$ a une boule tangente ext\'erieurement de rayon (au moins) $\delta_\sigma$, pour un certain $\sigma\in(0,s_0)$, au point $q\in\partial E\cap\overline{\Omega}$, on a
	\bgs{
		\liminf_{\varrho\to0^+}\I_s^\varrho[E](q)\geq\frac{\beta}{s}>0,\qquad\forall\,s\in(0,\sigma].
	}
\end{frtheorem}

Introduisons maintenant la d\'efinition suivante.

\begin{frdefn}\label{CH_fr_CH:0:deltadense}
	Soit $\Omega\subseteq \Rn$ un ensemble ouvert born\'e.
	On dit qu'un ensemble $E$ est \emph{$\delta$-dense} dans $\Omega$, pour un certain $\delta>0$ fix\'e, si $|B_\delta(x)\cap E|>0$ pour chaque $x\in \Omega$ tel que $B_\delta(x)\Subset\Omega$.
\end{frdefn}

En exploitant un argument g\'eom\'etrique d\'elicat et le Th\'eor\`eme
\ref{CH_fr_CH:0:positive_curvature}, nous pouvons alors poursuivre l'id\'ee heuristique d\'ecrite ci-dessus et prouver le r\'esultat de classification suivant:

\begin{frtheorem}\label{CH_fr_CH:0:Main_THM}
	Soit $\Omega\subseteq\R^n$ un ensemble ouvert born\'e et connexe ayant fronti\`ere de classe $C^2$. Soit $E_0\subseteq \Co \Omega$ tel que
	\[
	\overline{\alpha}(E_0)<\frac{\varpi_n}{2}.
	\]  
	Alors, les deux r\'esultats suivants sont v\'erifi\'es.\\
	A)  Sont $s_0$ et $\delta_s$ comme dans le Th\'eor\`eme \ref{CH_fr_CH:0:positive_curvature}. Il existe $s_1=s_1(E_0,\Omega)\in (0,s_0]$ tel que, si $s<s_1$ et $E$ est un ensemble $s$-minimal dans $\Omega$ ayant donn\'ee ext\'erieure $E_0$, alors, soit
	\bgs{
		(A.1) \;  E\cap \Omega=\emptyset \quad  \mbox{ ou} \quad\; (A.2)\;  E \mbox{ est } \delta_s-\mbox{dense dans }\Omega.
	}
	\noindent
	B) Soit \\
	(B.1) il existe
	$\tilde s=\tilde s(E_0,\Omega)\in (0,1)$ tel que si $E$ est un ensemble $s$-minimal dans $\Omega$ ayant donn\'ee ext\'erieure $E_0$ et $s\in(0,\tilde s)$, alors
	\bgs{
		E\cap \Omega=\emptyset,
	}
	ou \\
	(B.2)    ils existent  $\delta_k \searrow 0$, $s_k \searrow 0$ et une suite d'ensembles $E_k$ tels que chaque $E_k$ est $s_k$-minimal dans $\Omega$ par rapport \`a la donn\'ee ext\'erieure $E_0$ et pour chaque $k$
	\bgs{
		\partial E_k \cap B_{\delta_k}(x) \neq \emptyset \quad \mbox{for every } B_{\delta_k}(x)\Subset \Omega.
	}
\end{frtheorem}

Approximativement, soit les ensembles $s$-minimaaux sont vides dans $\Omega$ quand $s$ est assez petit, ou nous pouvons trouver une suite $E_k$ d'ensembles $s_k$-minimaux, pour $s_k\searrow0$, dont les fronti\`eres ont tendance \`a remplir (topologiquement) le domaine $\Omega$ dans la limite $k\to\infty$.

Nous soulignons que le comportement typique consiste \`a \^etre vide. En fait, si la donn\'ee ext\'erieure $E_0\subseteq\Co\Omega$ n'entoure pas compl\`etement le domaine $\Omega$, nous avons le r\'esultat suivant:

\begin{frtheorem}\label{CH_fr_CH:0:Not_Surrounded}
	Soit $\Omega$ un ensemble ouvert born\'e et connexe ayant fronti\`ere $C^2$. Soit $E_0\subseteq \Co \Omega$ tel que
	\[
	\overline{\alpha}(E_0)<\frac{\varpi_n}{2},
	\]
	et soit $s_1$ comme dans le Th\'eor\`eme \ref{CH_fr_CH:0:Main_THM}. Supposons qu'ils existent $R>0$ et $x_0\in \partial \Omega$ tels que
	\[
	B_R(x_0)\setminus \Omega \subseteq \Co E_0.
	\]
	Alors, il existe $s_3=s_3(E_0,\Omega)\in(0,s_1]$ tel que, si $s<s_3$ et $E$ est un ensemble $s$-minimal dans $\Omega$ par rapport \`a la donn\'ee ext\'erieure $E_0$, alors 
	\[  E\cap \Omega=\emptyset .\]
\end{frtheorem}

Nous observons que la condition $\overline{\alpha}(E_0)<\varpi_n/2$
est en quelque sorte optimale. En effet,
lorsque $\alpha(E_0)$ existe et 
\[ \alpha(E_0)=\frac{\varpi_n}2,\]
plusieurs configurations peuvent se produire, selon la position de $\Omega$ par rapport \`a la donn\'ee ext\'erieure $E_0\setminus \Omega$---nous fournissons divers exemples au Chapitre \ref{Asympto0_CH_label}.

En outre, notez que lorsque $E$ est $s$-minimal dans $\Omega$ par rapport \`a $E_0$, alors $\Co E$ est $s$-minimal dans $\Omega$ par rapport \`a $\Co E_0$. En plus,
\[
\underline \alpha(E_0) >\frac{\varpi_n}{2} \qquad \implies \qquad \overline \alpha (\Co E_0)< \frac{\varpi_n}{2}.
\]
Ainsi, dans ce cas, nous pouvons appliquer les Th\'eor\`emes \ref{CH_fr_CH:0:positive_curvature}, \ref{CH_fr_CH:0:Main_THM} et \ref{CH_fr_CH:0:Not_Surrounded} \`a $\Co E$ par rapport \`a la donn\'ee ext\'erieure $\Co E_0$. Par exemple, si
$E$ est $s$-minimal dans $\Omega$ par rapport \`a la donn\'ee ext\'erieure $E_0$ tel que
\[ \underline \alpha(E_0) >\frac{\varpi_n}{2}, \]
et $s<s_1(\Co E_0, \Omega)$,
alors, soit
\[ E\cap \Omega=\Omega \qquad \mbox{ ou }  \qquad  \Co E \; \mbox{ est } \; \delta_s(\Co E_0)-\mbox{dense}.\]
Les analogues des Th\'eor\`emes mentionn\'es ci-dessus peuvent \^etre obtenus de la m\^eme mani\`ere.

Par cons\'equent, \`a partir de nos r\'esultats principaux et des observations ci-dessus, nous avons une classification compl\`ete des surfaces minimales non locales lorsque $s$ est petit, quand
\[ \alpha(E_0)\neq  \frac{\varpi_n}{2} .\] 

\medskip

Nous soulignons que les ph\'enom\`enes de stickiness d\'ecrits dans \cite{boundary}
et au Chapitre \ref{Asympto0_CH_label} sont sp\'ecifiques aux surfaces minimales non locales, car les surfaces minimales classiques traversent transversalement la fronti\`ere d'un domaine convexe.

Fait int\'eressant, ces ph\'enom\`enes de stickiness ne sont pas pr\'esents dans le cas du Laplacien fractionnaire, o\`u la donn\'ee du probl\`eme de Dirichlet est atteint de mani\`ere continue sous des hypoth\`eses plut\^ot g\'en\'erales, voir \cite{MR3168912}.
Cependant, les solutions des \'equations de $s$-Laplace ne sont g\'en\'eralement pas meilleures que $C^s$ \`a la fronti\`ere, donc la continuit\'e uniforme d\'eg\'en\`ere lorsque~$s\to0^+$.

D'autre part, dans le cas de fonctions harmoniques fractionnaires, une contrepartie partielle du ph\'enom\`ene de stickiness est, en un sens, donn\'ee par les solutions explosives \`a la fronti\`ere construites dans \cite{MR3393247,MR2985500}
(\`a savoir, dans ce cas, la fronti\`ere du sous-graphe de la fonction harmonique fractionnaire contient des murs verticaux).

Nous mentionnons aussi que des ph\'enom\`enes de stickiness pour sous-graphes minimaux non locaux---\'eventuellement en pr\'esence d'obstacles---seront \'etudi\'es dans le prochain article \cite{LuCla}.

\medskip

Dans la derni\`ere partie du Chapitre \ref{Asympto0_CH_label} nous prouvons que la courbure moyenne fractionnaire est continue pour toutes les variables.

Pour simplifier un peu la situation, supposons que $E_k,\,E\subseteq\R^n$ sont des ensembles ayant fronti\`eres $C^{1,\gamma}$, pour un certain $\gamma\in(0,1]$, tels que les fronti\`eres $\partial E_k$ convergent localement au sens $C^{1,\gamma}$ vers la fronti\`ere de $E$, pour $k\to\infty$.
Alors, nous prouvons que, si nous avons une s\'equence de points $x_k\in\partial E_k$ tels que $x_k\to x\in\partial E$ et une suite de param\`etres $s_k,\,s\in(0,\gamma)$ tels que $s_k\to s$, on a
\[
\lim_{k\to\infty}\I_{s_k}[E_k](x_k)=\I_s[E](x).
\]
En outre, nous \'etendons de mani\`ere appropri\'ee ce r\'esultat de convergence afin de couvrir \'egalement les cas dans lesquels $s_k\to1$ ou $s_k\to0$.

En particulier, consid\'erons un ensemble
$E\subseteq\R^n$ tel que $\alpha(E)$ existe et $\partial E$ est de classe $C^2$. Alors, si on d\'efinit
\sys[
\tilde \I_s  {[}E{]} (x):=]{ &s(1-s)\I_s[E](x),  & \mbox{ pour } &s\in (0,1) \\
	&{\varpi_{n-1}} H[E](x), &\mbox{ pour } &s=1\\
	&\varpi_n-2\alpha(E), &\mbox{ pour }
	&s=0,}
la fonction
\[
\tilde\I_{(\,\cdot\,)}[E](\,\cdot\,)
:[0,1]\times\partial E\longrightarrow\R,
\qquad(s,x)\longmapsto \tilde \I_s[E](x),
\]
est continue.
Il est int\'eressant de noter que la courbure moyenne fractionnaire en un point fix\'e $q\in\partial E$ peut changer de signe lorsque $s$ varie de 0 \`a 1. En outre---en cons\'equence de la continuit\'e dans le param\`etre fractionnaire $s$---dans un tel cas, il existe une valeur $\sigma\in(0,1)$ tel que $\I_\sigma[E](q)=0$.

\subsection{Cadre non param\'etrique}

Au Chapitre \ref{CH_Nonparametric}, nous introduisons un cadre fonctionnel pour \'etudier les minimiseurs du p\'erim\`etre fractionnaire qui peuvent \^etre \'ecrits globalement en tant que sous-graphes, c'est-\`a-dire
\[
\Sg(u)=\left\{(x,x_{n+1})\in\R^{n+1}\,|\,x_{n+1}<u(x)\right\},
\]
pour une fonction mesurable $u:\R^n\to\R$. Nous appelons les fronti\`eres de ces minimiseurs des \emph{graphes minimaux non locaux}.

Nous d\'efinissons une version fractionnaire de la fonctionnelle d'aire classique et nous \'etudions ses propri\'et\'es fonctionnelles et g\'eom\'etriques.
Ensuite, nous nous concentrons sur les minimiseurs et nous prouvons des r\'esultats d'existence et d'unicit\'e par rapport \`a une grande classe de donn\'ees ext\'erieures, qui inclut les fonctions localement born\'ees.
De plus, l'une des contributions principales du Chapitre \ref{CH_Nonparametric} consiste \`a prouver l'\'equivalence entre:
\begin{itemize}
	\item minimiseurs de la fonctionnelle d'aire fractionnaire,
	\item minimiseurs du p\'erim\`etre fractionnaire,
	\item solutions faibles de l'\'equation de courbure moyenne fractionnaire,
	\item solutions de viscosit\'e de l'\'equation de courbure moyenne fractionnaire,
	\item fonctions lisses r\'esolvant ponctuellement l'\'equation de courbure moyenne fractionnaire.
\end{itemize}

\medskip

Avant de donner un aper\c{c}u d\'etaill\'e des principaux r\'esultats, rappelons la d\'efinition de la fonctionnelle d'aire classique.
\'Etant donn\'e un ensemble ouvert born\'e
$\Omega\subseteq\R^n$ ayant fronti\`ere Lipschitz, la fonctionnelle d'aire est d\'efinie comme
\[
\mathscr A(u,\Omega):=\int_\Omega\sqrt{1+|\nabla u|^2}\,dx=\Ha^n\left(\left\{(x,u(x))\in\R^{n+1}\,|\,x\in\Omega\right\}\right),
\]
pour chaque fonction Lipschitz $u:\overline{\Omega}\to\R$. On \'etend alors cette fonctionnelle, en d\'efinissant la fonctionnelle d'aire relax\'ee d'une fonction $u\in L^1(\Omega)$ comme
\[
\mathscr A(u,\Omega):=\inf\left\{\liminf_{k\to\infty}\mathscr A(u_k,\Omega)\,|\,u_k\in C^1(\overline{\Omega}),\,\|u-u_k\|_{L^1(\Omega)}\to0\right\}.
\]
On voit bien que, si $u\in L^1(\Omega)$, alors
\eqlab{\label{CH_fr_CH:0:obs1_locarea}
	\mathscr A(u,\Omega)<\infty\quad\Longleftrightarrow\quad u\in BV(\Omega),
}
dans quel cas
\eqlab{\label{CH_fr_CH:0:obs2_locarea}
	\mathscr A(u,\Omega)=\Per\left(\Sg(u),\Omega\times\R\right).
}

Approximativement, les fonctions \`a variation born\'ee sont pr\'ecis\'ement les fonctions int\'egrables dont les sous-graphes ont p\'erim\`etre fini---pour les d\'etails, voir, par exemple, \cite{Giusti,GiaMart12}.

Nous pourrions donc \^etre tent\'es de d\'efinir une version fractionnaire de la fonctionnelle d'aire en consid\'erant le $s$-p\'erim\`etre \`a la place du p\'erim\`etre classique, d\'efinissant, pour une fonction mesurable $u:\R^n\to\R$,
\[
\mathscr A_s(u,\Omega):=\Per_s(\Sg(u),\Omega\times\R).
\]
Cependant, comme nous l'avons observ\'e \`a la fin de la Section \ref{CH_fr_CH:0:NMS_Sec}, une telle d\'efinition ne peut pas fonctionner, car
\[
\Per_s^{NL}(\Sg(u),\Omega\times\R)=\infty,
\]
m\^eme si $u\in C^\infty_c(\R^n)$.

Avant de poursuivre, quelques observations s'imposent.
M\^eme si la partie non locale du p\'erim\`etre fractionnaire dans le cylindre $\Omega^\infty:=\Omega\times\R$ est infinie, nous rappelons que nous savons---voir la fin de la Section \ref{CH_fr_CH:0:NMS_Sec}---que la partie locale est finie, si la fonction $u$ est assez r\'eguli\`ere dans $\Omega$.

Si la fonction $u$ est born\'ee dans $\Omega$, alors nous pouvons consid\'erer le p\'erim\`etre fractionnaire dans le ``cylindre tronqu\'e'' $\Omega^M:=\Omega\times(-M,M)$, o\`u $M\geq\|u\|_{L^\infty(\Omega)}$, au lieu du cylindre $\Omega^\infty$. Comme nous le verrons plus loin, en poursuivant cette id\'ee, nous obtenons une famille de fonctionnels d'aire fractionnaires $\F^M_s(\,\cdot\,,\Omega)$.

\smallskip

Par ailleurs, il existe une autre possibilit\'e de d\'efinir une fonctionnelle d'aire fractionnaire.
Dans \cite{regularity}, les auteurs ont observ\'e que lorsque~$E\subseteq\R^{n+1}$ est le sous-graphe d'une fonction~$u$,
sa courbure moyenne fractionnaire peut \^etre \'ecrite comme un op\'erateur int\'egrodiff\'erentiel agissant sur~$u$. Plus pr\'ecis\'ement, si~$u: \R^n \to \R$ est une fonction de classe~$C^{1, 1}$ dans un voisinage d'un point~$x \in \R^n$, nous avons
\[
H_s[\Sg(u)](x,u(x)) = \h_s u(x),
\]
o\`u
\[
\h_s u(x) := 2 \, \PV \int_{\R^n} G_s \left( \frac{u(x)-u(y)}{|x - y|} \right) \frac{dy}{|x - y|^{n+s}},
\]
et
\[
G_s(t):=\int_0^t g_s(\tau)\,d\tau,\qquad g_s(t):=\frac{1}{(1+t^2)^\frac{n+1+s}{2}} \quad \mbox{pour } t\in\R.
\]
Nous montrons maintenant que $\h_s$ est l'op\'erateur d'Euler-Lagrange associ\'e \`a une fonctionnelle (convexe) $\F_s(\,\cdot\,,\Omega)$, que nous consid\'ererons alors comme la fonctionnelle d'aire $s$-fractionnaire.

Commen\c{c}ons par remarquer que, lorsque $u$ n'est pas assez r\'egulier autour de $x$, la quantit\'e $\h_su(x)$ n'est g\'en\'eralement pas bien d\'efinie, en raison du manque d'annulation requise pour la valeur principale afin de converger. N\'eanmoins, nous pouvons comprendre l'op\'erateur $\h_s$ tel que d\'efini dans le sens faible (distributionnel) suivant. \'Etant donn\'ee une fonction mesurable~$u:\R^n \to \R$, nous d\'efinissons
\[
\langle \h_s u, v\rangle:=\int_{\Rn}\int_{\Rn} G_s \left( \frac{u(x)-u(y)}{|x-y|} \right) \big( v(x)-v(y) \big) \, \frac{dx\, dy}{|x-y|^{n+s}}
\]
pour chaque $v\in C^\infty_c(\R^n)$. Plus g\'en\'eralement, il est imm\'ediat de voir---en profitant du fait que~$G_s$ est born\'ee---que cette d\'efinition est bien pos\'ee pour chaque~$v\in W^{s,1}(\Rn)$. En effet, on a
\[
|\langle \h_s u ,v \rangle|\le\frac{\Lambda_{n,s}}{2} \, [v]_{W^{s,1}(\Rn)},
\]
o\`u
\[
\Lambda_{n,s}:=\int_\R g_s(t)\,dt<\infty.
\]
Partant,~$\h_s u$ peut \^etre interpr\'et\'ee comme une forme lin\'eaire et continue $\langle \h_s u,\,\cdot\, \rangle \in ( W^{s, 1}(\R^n))^*$.
Remarquablement, cela vaut pour chaque fonction mesurable~$u: \R^n \to \R$, quelle que soit sa r\'egularit\'e.

Nous d\'efinissons maintenant
\[
\G_s(t):=\int_0^t G_s(\tau)\,d\tau\quad\mbox{pour }t\in\R,
\]
et, \'etant donn\'e une fonction mesurable $u:\R^n\to\R$ et un ensemble ouvert $\Omega\subseteq\R^n$, nous d\'efinissons la \emph{fonctionnelle d'aire $s$-fractionnaire}
\[
\F_s(u,\Omega):=\iint_{\R^{2n}\setminus(\Co\Omega)^2}\G_s\left(\frac{u(x)-u(y)}{|x-y|}\right)\frac{dx\,dy}{|x-y|^{n-1+s}}.
\]
Ensuite, au moins formellement, nous avons
\[
\frac{d}{d\eps}\Big|_{\eps=0}\F_s(u+\eps v,\Omega)=\langle\h_s u,v\rangle\quad\mbox{pour chaque }v\in C^\infty_c(\Omega).
\]

\smallskip

Nous remarquons que dans le Chapitre \ref{CH_Nonparametric}, nous allons en fait consid\'erer des fonctionnelles plus g\'en\'erales du type aire fractionnaire---en prenant dans les d\'efinitions ci-dessus une fonction continue et paire $g:\R\to(0,1]$ satisfaisant une condition d'int\'egrabilit\'e appropri\'ee, et les fonctions correspondantes $G$ et $\G$, \`a la place de $g_s,\,G_s$ et $\G_s$ respectivement. Cependant, pour plus de simplicit\'e dans cette introduction, nous nous en tenons au ``cas g\'eom\'etrique''  correspondant au choix $g=g_s$.

\smallskip

Voyons maintenant les propri\'et\'es fonctionnelles de $\F_s(\,\cdot\,,\Omega)$ et sa relation avec le p\'erim\`etre fractionnaire.

\`A partir de maintenant, nous consid\'erons $n\ge1$, $s\in(0,1)$ et un ensemble ouvert born\'e $\Omega\subseteq\R^n$ ayant fronti\`ere Lipschitz.

Il est commode de scinder la fonctionnelle d'aire fractionnaire en tant que somme de sa partie locale et de sa partie non locale, c'est-\`a-dire
\[
\F_s(u,\Omega)=\A_s(u,\Omega)+\Nl_s(u,\Omega),
\]
o\`u
\[
\A_s(u,\Omega):=\int_\Omega\int_\Omega
\G_s\left(\frac{u(x)-u(y)}{|x-y|}\right)\frac{dx\,dy}{|x-y|^{n-1+s}}
\]
et
\[
\Nl_s(u,\Omega):=2\int_\Omega\int_{\Co\Omega}
\G_s\left(\frac{u(x)-u(y)}{|x-y|}\right)\frac{dx\,dy}{|x-y|^{n-1+s}}.
\]
Mentionnons tout d'abord l'observation int\'eressante suivante---voir, par exemple, Lemme \ref{CH:A:usef_ineq_hit}.
Si $u:\Omega\to\R$ est une fonction mesurable, alors
\[
[u]_{W^{s,1}(\Omega)}<\infty\quad\implies\quad\|u\|_{L^1(\Omega)}<\infty.
\]

En ce qui concerne la partie locale de la fonctionnelle d'aire fractionnaire, nous prouvons que, si $u:\Omega\to\R$ est une fonction mesurable, alors
\bgs{
	\A_s(u,\Omega)<\infty\quad&\Longleftrightarrow\quad u\in W^{s,1}(\Omega)\\
	&
	\Longleftrightarrow\quad
	\Per_s^L(\Sg(u),\Omega\times\R)<\infty.
}
En outre, si $u\in W^{s,1}(\Omega)$, alors
\[
\Per_s^L(\Sg(u),\Omega\times\R)
=\A_s(u,\Omega)+c,
\]
pour une certaine constante $c=c(n,s,\Omega)\geq0$.
Ces r\'esultats peuvent \^etre consid\'er\'es comme les contreparties fractionnaires de \eqref{CH_fr_CH:0:obs1_locarea} et \eqref{CH_fr_CH:0:obs2_locarea}.

D'autre part, pour que la partie non locale soit finie, nous devons imposer une condition d'int\'egrabilit\'e sur $u$ \`a l'infini, \`a savoir
\eqlab{\label{CH_fr_CH:0:Infty_int_cond}
	\int_\Omega\left(\int_{\Co\Omega}
	\frac{|u(y)|}{|x-y|^{n+s}}\,dy\right)dx<\infty.
}
Une telle condition est remplie, par exemple, si $u$ est globalement born\'ee dans $\R^n$ et, en g\'en\'eral, cela implique que la fonction $u$ doit avoir un comportement sous-lin\'eaire \`a l'infini. C'est donc une condition tr\`es restrictive.

En effet, on remarque que l'op\'erateur  $\h_su$ est bien d\'efini en un point $x$---\`a condition que $u$ soit assez r\'egulier dans un voisinage de $x$---sans avoir \`a imposer de conditions \`a $u$ \`a l'infini. De plus, comme nous l'avons observ\'e dans la Section \ref{CH_fr_CH:0:NMS_Sec}, en cons\'equence du Th\'eor\`eme \ref{CH_fr_CH:0:Local_Existence} et du \cite[Theorem 1.1]{graph} nous savons que, \'etant donn\'e toute fonction continue $\varphi:\R^n\to\R$, il existe une fonction $u:\R^n\to\R$ telle que $u=\varphi$ dans $\R^n\setminus\overline{\Omega}$, $u\in C(\overline{\Omega})$ et $\Sg(u)$ est localement $s$-minimal dans $\Omega^\infty$. Soulignons qu'aucune condition sur $\varphi$ \`a l'infini n'est requise.

Pour ces raisons, la condition \eqref{CH_fr_CH:0:Infty_int_cond} semble \^etre anormalement restrictive dans notre cadre---m\^eme si, \`a premi\`ere vue, elle semble n\'ecessaire, car elle est n\`ecessaire pour garantir que $\F_s$ soit bien d\'efini.

Afin d'\'eviter d'imposer la condition \eqref{CH_fr_CH:0:Infty_int_cond}, nous d\'efinissons---voir \eqref{CH:4:NMldef}---pour chaque $M\geq0$, la partie non locale ``tronqu\'ee'' $\Nl_s^M(u,\Omega)$ et la fonctionnelle d'aire fractionnaire tronqu\'ee
\[
\F_s^M(u,\Omega):=\A_s(u,\Omega)+\Nl_s^M(u,\Omega).
\]
Approximativement, l'id\'ee consiste \`a ajouter, \`a l'int\'erieur de la double int\'egrale d\'efinissant la partie non locale, un terme \'equilibrant la contribution venant de l'ext\'erieur de $\Omega$. Par exemple, dans le cas le plus simple $M=0$, on a
\[
\Nl_s^0(u,\Omega)=2\int_\Omega\left\{\int_{\Co\Omega}\left[\G_s\left(\frac{u(x)-u(y)}{|x-y|}\right)-\G_s\left(\frac{u(y)}{|x-y|}\right)\right]\frac{dy}{|x-y|^{n-1+s}}\right\}dx.
\]
Remarquablement, \'etant donn\'ee une fonction mesurable $u:\R^n\to\R$, on a
\[
|\Nl_s^M(u,\Omega)|<\infty\quad\mbox{si }u|_\Omega\in W^{s,1}(\Omega),
\]
quel que soit le comportement de $u$ dans $\Co\Omega$. D'autre part, nous remarquons qu'en g\'en\'eral, la partie non locale tronqu\'ee peut \^etre n\'egative, sauf si nous exigeons que $u$ soit born\'ee dans $\Omega$ et nous prenons $M\ge\|u\|_{L^\infty(\Omega)}$. D'un point de vue g\'eom\'etrique, les fonctionnelles d'aire fractionnaire tronquées correspondent \`a la prise en compte du p\'erim\`etre fractionnaire dans le cylindre tronqu\'e $\Omega^M$.

En fait, si $u:\R^n\to\R$ est une fonction mesurable telle que $u|_\Omega\in W^{s,1}(\Omega)\cap L^\infty(\Omega)$, et $M\ge\|u\|_{L^\infty(\Omega)}$, on a
\[
\F_s^M(u,\Omega)=\Per_s	\big(\Sg(u),\Omega\times(-M,M)\big)+c_M,
\]
pour une certaine constante $c_M=c_M(n,s,\Omega)\geq0$.

\medskip

Nous passons maintenant \`a l'\'etude des minimiseurs de la fonctionnelle d'aire fractionnaire.

\'Etant donn\'ee une fonction mesurable $\varphi:\Co\Omega\to\R$, nous d\'efinissons l'espace
\[
\W_\varphi^s(\Omega):=\left\{
u:\R^n\to\R\,|\,u|_\Omega\in W^{s,1}(\Omega)\mbox{ et }u=\varphi\mbox{ p.p. dans }\Co\Omega\right\},
\]
et on dit que $u\in\W^s_\varphi(\Omega)$ est un \emph{minimiseur} de $\F_s$ dans $\W^s_\varphi(\Omega)$, si
\[
\iint_{Q(\Omega)} \left\{ \G_s \left( \frac{u(x) - u(y)}{|x - y|} \right) - \G_s \left( \frac{v(x) - v(y)}{|x - y|} \right) \right\} \frac{dx\,dy}{|x - y|^{n - 1 + s}} \le 0
\]
pour chaque $v\in\W^s_\varphi(\Omega)$. Ci-dessus, nous avons utilis\'e la notation $Q(\Omega):=\R^{2n}\setminus(\Co\Omega)^2$.
Soulignons qu'une telle d\'efinition est bien pos\'ee sans devoir imposer de conditions \`a la \emph{donn\'ee ext\'erieure} $\varphi$, comme en effet---gr\^ace \`a l'in\'egalit\'e de type Hardy fractionnaire du Th\'eor\`eme \ref{CH:4:FHI}---nous avons
\[
\iint_{Q(\Omega)}\left|\G_s\left(\frac{u(x)-u(y)}{|x-y|}\right)
-\G_s\left(\frac{v(x)-v(y)}{|x-y|}\right)\right|\frac{dx\,dy}{|x-y|^{n-1+s}}
\le C\,\Lambda_{n,s}\|u-v\|_{W^{s,1}(\Omega)},
\]
pour chaques~$u,\,v\in\W^s_\varphi(\Omega)$,
pour une certaine constante $C=C(n,s,\Omega)>0$.

Nous prouvons l'existence de minimiseurs par rapport \`a des donn\'ees ext\'erieures satisfaisant une condition d'int\'egrabilit\'e appropri\'ee dans un voisinage du domaine $\Omega$. Plus pr\'ecis\'ement, \'etant donn\'e un ensemble ouvert~$\Op\subseteq\R^n$
tel que~$\Omega\Subset\Op$, nous d\'efinissons la \emph{queue tronqu\'ee} de~$\varphi:\Co\Omega\to\R$ au point~$x\in\Omega$ comme
\bgs{
	\Tail_s(\varphi, \Op\setminus\Omega;x):=\int_{\Op\setminus\Omega}\frac{|\varphi(y)|}{|x-y|^{n+s}}\,dy.
}
Nous utilisons la notation
\[
\Omega_\varrho:=\{x\in\R^n\,|\,d(x,\Omega)<\varrho\},
\]
pour $\varrho>0$, pour d\'enoter le $\varrho$-voisinage de $\Omega$.
Alors, nous prouvons ce qui suit:
\begin{frtheorem} 
	Il existe une constante~$\Theta > 1$, qui ne d\'epend que de~$n$ et~$s$, telle que, \'etant donn\'e toute fonction~$\varphi: \Co\Omega\to \R$
	avec~$\Tail_s(\varphi, \Omega_{\Theta \diam(\Omega)} \setminus \Omega;\,\cdot\,) \in L^1(\Omega)$, il existe un minimiseur unique~$u$ de~$\F_s$ dans~$\W^s_\varphi(\Omega)$. En plus,~$u$ satisfait
	\[
	\| u \|_{W^{s, 1}(\Omega)} \le C \left( \left\| \Tail_s(\varphi,\Omega_{\Theta \diam(\Omega)}\setminus \Omega;\,\cdot\,) \right\|_{L^1(\Omega)} + 1 \right),
	\]
	pour une certaine constante~$C=C(n,s,\Omega)>0$.
\end{frtheorem}
Nous observons que la condition sur l'intégrabilit\'e de la queue est beaucoup plus faible que \eqref{CH_fr_CH:0:Infty_int_cond}, puisque nous n'exigeons rien du comportement de $\varphi$ \`a l'ext\'erieur de $\Omega_{\Theta \diam(\Omega)}$. 

Nous mentionnons \'egalement que, approximativement, l'int\'egrabilit\'e de la queue \'equivaut \`a l'int\'egrabilit\'e de $\varphi$ plus certaines conditions de r\'egularit\'e pr\`es de la fronti\`ere $\partial\Omega$. Par exemple, si $\varphi\in L^1(\Omega_{\Theta \diam(\Omega)}\setminus\Omega)$ et il existe $\varrho>0$ tel que, soit $\varphi\in W^{s,1}(\Omega_\varrho\setminus\Omega)$ ou
$\varphi\in L^\infty(\Omega_\varrho\setminus\Omega)$, alors
$\Tail_s(\varphi, \Omega_{\Theta \diam(\Omega)} \setminus \Omega;\,\cdot\,) \in L^1(\Omega)$.

L'unicit\'e du minimiseur est une cons\'equence de la stricte convexit\'e de $\F_s$. 
D'autre part, afin de prouver l'existence, nous exploitons les (uniques) minimiseurs $u_M$ des fonctionnelles $\F^M_s(\,\cdot\,,\Omega)$---consid\'er\'es dans leur domaine naturel. Nous exploitons l'hypoth\`ese sur l'int\'egrabilit\'e de la queue pour prouver une estimation uniforme de la norme $W^{s,1}(\Omega)$ des minimiseurs $u_M$, ind\'ependamment de $M\geq0$. Donc, quitte \`a extraire des sous-suites, $u_M$ converge, lorsque $M\to\infty$, vers une fonction limite $u$, qui est facilement prouv\'e \^etre un minimiseur de $\F_s$.

\smallskip

En outre, nous prouvons que, si $u$ est un minimiseur de $\F_s$ dans $\W^s_\varphi(\Omega)$, alors $u\in L^\infty_{\loc}(\Omega)$. De plus, nous montrons que, si la donn\'ee ext\'erieure $\varphi$ est born\'ee dans un voisinage assez gros de $\Omega$, alors $u\in L^\infty(\Omega)$, et nous \'etablissons \'egalement une estimation a priori pour la norme $L^\infty$.

\smallskip

Revenons \`a la relation entre la fonctionnelle d'aire fractionnaire et le p\'erim\`etre fractionnaire.
Nous montrons qu'en r\'earrangeant correctement un ensemble $E$ dans la direction verticale, nous diminuons le $s$-p\'erim\`etre. Plus pr\'ecis\'ement, \`a partir d'un ensemble~$E \subseteq \R^{n + 1}$, nous consid\'erons la fonction~$w_E: \R^n \to \R$ d\'efinie comme
\bgs{
	w_E(x) := \lim_{R \rightarrow +\infty} \left( \int_{-R}^R \chi_{E}(x, t) \, dt - R \right)
}
pour chaque~$x \in \R^n$.

Alors, nous avons le r\'esultat suivant:

\begin{frtheorem}\label{CH_fr_CH:0:Rearr_THM}
	Soit~$E \subseteq \R^{n + 1}$ tel que~$E \setminus \Omega^\infty$ est un sous-graphe et
	\[
	\Omega \times (-\infty, -M) \subseteq E \cap \Omega^\infty \subseteq \Omega \times (-\infty, M),
	\]
	pour un certain~$M > 0$. Alors,
	\[
	\Per_s(\Sg(w_E), \Omega^M) \le \Per_s(E, \Omega^M).
	\]
	L'in\'galit\'e est stricte sauf si~$\Sg(w_E)=E$.
\end{frtheorem}

En exploitant \'egalement le fait qu'un minimiseur est localement born\'e, nous prouvons que, si $u:\R^n\to\R$
est une fonction mesurable telle que $u\in W^{s,1}(\Omega)$, alors
\[
u\mbox{ minimise }\F_s\mbox{ dans }\W^s_u(\Omega)\quad\implies\quad
\Sg(u)\mbox{ est localement $s$-minimal dans } \Omega^\infty.
\]
Le Th\'eor\`eme \ref{CH_fr_CH:0:Rearr_THM} \'etend au cadre fractionnaire un r\'esultat bien connu tenant pour le p\'erim\`etre classique---voir, par exemple, \cite[Lemma 14.7]{Giusti}. Cependant, notez que dans le cadre fractionnaire, en raison du caract\`ere non local des fonctionnelles impliqu\'ees, nous devons supposer que l'ensemble $E$ est d\'ej\`a un sous-graphe \`a l'ext\'erieur du cylindre $\Omega^\infty$.

Nous observons \'egalement que, puisque $u$ est localement born\'ee dans $\Omega$ et son sous-graphe est localement $s$-minimal dans le cylindre $\Omega^\infty$, gr\^ace \`a \cite[Theorem 1.1]{CaCo}
nous avons $u\in C^\infty(\Omega)$---c'est-\`a-dire, les minimiseurs de $\F_s$ sont lisses.

\medskip

Voyons maintenant l'\'equation d'Euler-Lagrange satisfaite par les minimiseurs.
Nous introduisons d'abord la notion de solutions faibles.

Soit $f\in C(\overline{\Omega})$. On dit qu'une fonction mesurable $u:\R^n\to\R$ est une solution faible de $\h_s u=f$ dans $\Omega$, si
\[
\langle\h_s u,v\rangle=\int_\Omega fv\,dx,
\]
pour chaque $v\in C^\infty_c(\Omega)$.

En cons\'equence de la convexit\'e de $\F_s$, il est facile de prouver que, \'etant donn\'ee une fonction mesurable $u:\R^n\to\R$ telle que $u\in W^{s,1}(\Omega)$, on a
\bgs{
	u \mbox{ minimise }\F_s\mbox{ dans }\W^s_u(\Omega)\quad
	\Longleftrightarrow\quad
	u\mbox{ est une solution faible de }\h_s u=0\mbox{ dans }\Omega.
}

Une autre notion naturelle de solution pour l'\'equation $\h_s u=f$ est celle d'une solution de viscosit\'e---nous nous r\'ef\'erons \`a la Section \ref{CH:4:ViscWeak_Sec} pour la d\'efinition pr\'ecise.
Un des principaux r\'esultats du Chapitre \ref{CH_Nonparametric} consiste \`a prouver que les (sous-)solutions de viscosit\'e sont des (sous-)solutions faibles. Plus pr\'ecis\'ement:

\begin{frtheorem}
	Soit~$\Omega\subseteq\R^n$ un ensemble ouvert born\'e et soit~$f\in C(\overline\Omega)$. Soit~$u:\R^n\to\R$
	localement int\'egrable et localement born\'e dans~$\Omega$. Si $u$ est une sous-solution de viscosit\'e,
	\bgs{
		\h_s u\le f\quad\mbox{dans }\Omega,
	}
	alors~$u$ est une sous-solution faible,
	\bgs{
		\langle\h_s u,v\rangle\le\int_\Omega fv\,dx,\qquad\forall\,v\in C^\infty_c(\Omega)\mbox{ telle que }v\ge0.
	}
\end{frtheorem}

\medskip

En combinant les principaux r\'esultats du Chapitre \ref{CH_Nonparametric} et en exploitant la r\'egularit\'e \`a la int\'erieure prouv\'ee dans \cite{CaCo}, on obtient ce qui suit:

\begin{frtheorem}\label{CH_fr_CH:0:Equiv_Intro}
	Soit~$u:\R^n\to\R$ une fonction mesurable telle que $u\in W^{s,1}(\Omega)$.
	Alors, les propositions suivantes sont \'equivalentes:
	\begin{itemize}
		\item[(i)] $u$ est une solution faible de~$\h_s u=0$ dans~$\Omega$,
		\item[(ii)] $u$ minimise~$\F_s$ dans~$\W^s_u(\Omega)$,
		\item[(iii)] $u\in L^\infty_{\loc}(\Omega)$ et~$\Sg(u)$ est localement $s$-minimal dans~$\Omega\times\R$,
		\item[(iv)] $u\in C^\infty(\Omega)$ et~$u$ est une solution ponctuelle de~$\h_s u=0$ dans~$\Omega$.
	\end{itemize}
	En plus, si~$u\in L^1_{\loc}(\R^n)\cap W^{s,1}(\Omega)$, alors les propositions ci-dessus sont \'equivalentes \`a:
	\begin{itemize}
		\item[(v)] $u$ est une solution de viscosit\'e de~$\h_s u=0$ dans~$\Omega$.
	\end{itemize}
\end{frtheorem}

Nous mentionnons \'egalement la version globale suivante du Th\'eor\`eme~\ref{CH_fr_CH:0:Equiv_Intro}:

\begin{frcorollary}\label{CH_fr_CH:0:Equiv_Intro_Global_Corollary}
	Soit~$u\in W^{s,1}_{\loc}(\R^n)$. Alors, les propositions suivantes sont \'equivalentes:
	\begin{itemize}
		\item[(i)] $u$ est une solution de viscosit\'e de~$\h_s u=0$ dans~$\R^n$,
		\item[(ii)] $u$ est une solution faible de~$\h_s u=0$ dans~$\R^n$,
		\item[(iii)] $u$ minimise~$\F_s$ dans~$\W^s_u(\Omega)$, pour chaque ensemble ouvert~$\Omega\Subset\R^n$ ayant fronti\`ere Lipschitz,
		\item[(iv)] $u\in L^\infty_{\loc}(\R^n)$ et~$\Sg(u)$ est localement~$s$-minimal dans~$\R^{n+1}$,
		\item[(v)] $u\in C^\infty(\R^n)$ et~$u$ est une solution ponctuelle de~$\h_s u=0$ dans~$\R^n$.
	\end{itemize}
\end{frcorollary}

\medskip

Nous signalons \'egalement que le cadre fonctionnel pr\'esent\'e ci-dessus s'\'etend facilement au probl\`eme avec obstacles. \`A savoir, en plus d'imposer la condition de la donn\'ee ext\'erieure $u=\varphi$ p.p. dans $\Co\Omega$, nous contraignons les fonctions \`a se trouver au-dessus d'un obstacle, c'est-\`a-dire, \`etant donn\'e un ensemble ouvert $A\subseteq\Omega$ et un obstacle $\psi\in L^\infty(A)$, nous nous bornons \`a consid\'erer ces fonctions $u\in\W^s_\varphi(\Omega)$ telles que $u\ge\psi$ p.p. dans $A$.

Au Chapitre \ref{CH_Nonparametric} nous examinons \'egalement bri\`evement ce probl\`eme d'obstacle, prouvant l'existence et l'unicit\'e d'un minimiseur, et sa relation avec le probl\`eme d'obstacle g\'eom\'etrique qui concerne le p\'erim\`etre fractionnaire.

\smallskip

Enfin, dans la derni\`ere Section du Chapitre \ref{CH_Nonparametric}, nous prouvons quelques r\'esultats d'approximation pour les sous-graphes ayant p\'erim\`etre fractionnaire (localement) fini. En particulier, en exploitant le r\'esultat surprenant de densit\'e \'etabli dans \cite{DSV17}, nous montrons que les sous-graphes $s$-minimaux peuvent \^etre approxim\'es de mani\`ere appropri\'ee par des sous-graphes de fonctions $\sigma$-harmoniques, pour chaque $\sigma\in(0,1)$ fix\'e.

\subsection{R\'esultats de rigidit\'e pour les graphes minimaux non locaux}

Au Chapitre \ref{CH_Bern_Mos_result} nous prouvons un r\'esultat de platitude pour des graphes minimaux non locaux entiers ayant des d\'eriv\'ees partielles minor\'es ou major\'es. Ce r\'esultat g\'en\'eralise au cadre fractionnaire des th\'eor\`emes classiques dues \`a Bernstein et Moser.

De plus, nous montrons que les graphes entiers ayant courbure moyenne fractionnaire constante sont minimales, \'etendant ainsi un r\'esultat c\'el\`ebre de Chern sur les graphes CMC classiques.

\medskip

Nous sommes int\'eress\'es par les sous-graphes qui minimisent localement le $s$-p\'erim\`etre dans tout l'espace $\R^{n+1}$.
Nous rappelons que, comme nous l'avons vu dans le Corollaire \ref{CH_fr_CH:0:Equiv_Intro_Global_Corollary}, sous des hypoth\`eses tr\`es faibles sur la fonction $u:\R^n\to\R$, le sous-graphe $\Sg(u)$ est localement $s$-minimal dans $\R^{n+1}$ si et seulement si $u$ satisfait \`a l'\'equation de courbure moyenne fractionnaire
\eqlab{ \label{CH_fr_CH:0:Hsu=0}
	\h_su=0\quad\mbox{dans }\R^n.
}
En outre, encore une fois gr\^ace au Corollaire \ref{CH_fr_CH:0:Equiv_Intro_Global_Corollary}, il existe plusieurs notions \'equivalentes de solution pour l'\'equation \eqref{CH_fr_CH:0:Hsu=0},
telles que solutions lisses, solutions de viscosit\'e et solutions faibles.

Dans ce qui suit, une solution de~\eqref{CH_fr_CH:0:Hsu=0} indiquera toujours une fonction~$u \in C^\infty(\R^n)$ qui satisfait l'identit\'e~\eqref{CH_fr_CH:0:Hsu=0} ponctuellement. Nous soulignons qu'aucune hypoth\`ese de croissance \`a l'infini n'est faite sur~$u$.

La contribution principale du Chapitre \ref{CH_Bern_Mos_result} est le r\'esultat suivant:

\begin{frtheorem} \label{CH_fr_CH:0:under_CRA_Pakmainthm}
	Soient~$n \ge \ell \ge 1$ des entiers,~$s \in (0, 1)$, et supposons que
	\begin{equation} \tag{$P_{s, \ell}$} \label{CH_fr_CH:0:SMINSINGCON}
	\mbox{il n'y a pas de c\^ones singuliers~$s$-minimaux dans~$\R^\ell$.}
	\end{equation}
	Soit~$u$ une solution de $\h_su=0$ dans $\R^n$, ayant~$n - \ell$ deriv\'ees partielles minor\'es ou major\'es.
	Alors,~$u$ est une fonction affine.
\end{frtheorem}

La caract\'erisation des valeurs de $s$ et $\ell$ pour lesquelles~\eqref{CH_fr_CH:0:SMINSINGCON} est satisfaite repr\'esente un probl\`eme ouvert difficile \`a r\'esoudre. N\'eanmoins, il est connu que la propri\'et\'e~\eqref{CH_fr_CH:0:SMINSINGCON} est vraie dans les cas suivants:
\begin{itemize}
	\item lorsque~$\ell = 1$ ou~$\ell = 2$, pour chaque~$s \in (0, 1)$;
	\item lorsque~$3\le \ell \le 7$ et~$s \in (1 - \varepsilon_0, 1)$ pour un certain~$\varepsilon_0 \in (0,1]$ ne d\'ependant que de~$\ell$. 
\end{itemize}
Le cas~$\ell = 1$ est vrai par d\'efinition, alors que le cas~$\ell = 2$ est le contenu de~\cite[Theorem~1]{SV13}. D'autre part, le cas~$3 \le \ell \le 7$ a \'et\'e \'etabli en~\cite[Theorem~2]{regularity}.

En cons\'equence du Th\'eor\`eme~\ref{CH_fr_CH:0:under_CRA_Pakmainthm} et des derni\`eres remarques, nous obtenons imm\'ediatement le r\'esultat suivant:

\begin{frcorollary} \label{CH_fr_CH:0:CORRigidity_mainthm}
	Soient~$n \ge \ell \ge 1$ des entiers et~$s \in (0, 1)$. Supposons que
	\begin{itemize}
		\item $\ell \in \{ 1, 2 \}$, ou
		\item $3\le \ell \le 7$ et~$s \in (1 - \varepsilon_0, 1)$, o\`u~$\varepsilon_0=\varepsilon_0(\ell) > 0$ est comme en~\cite[Theorem~2]{regularity}.
	\end{itemize}
	Soit~$u$ une solution de $\h_su=0$ dans $\R^n$, ayant~$n - \ell$ deriv\'ees partielles minor\'es ou major\'es.
	Alors,~$u$ est une fonction affine.
\end{frcorollary}

Nous observons que le Th\'eor\`eme~\ref{CH_fr_CH:0:under_CRA_Pakmainthm} est un nouveau r\'esultat de platitude pour les graphes~$s$-minimaux, en supposant que~\eqref{CH_fr_CH:0:SMINSINGCON} est vrai. Cela peut \^etre vu comme une g\'en\'eralisation du lemme de type De Giorgi fractionnaire contenu dans~\cite[Theorem~1.2]{FV17}, qui est r\'ecup\'er\'e ici en prenant~$\ell = n$. Dans ce cas, nous fournissons en effet une preuve alternative dudit r\'esultat.

D'autre part, le choix~$\ell = 2$ donne une am\'elioration de~\cite[Theorem~4]{FarV17}, quand sp\'ecialis\'e aux graphes~$s$-minimaux. \`A la lumi\`ere de ces observations, le Th\'eor\`eme~\ref{CH_fr_CH:0:under_CRA_Pakmainthm} et le Corollaraire~\ref{CH_fr_CH:0:CORRigidity_mainthm} peuvent \^etre vus comme un pont entre les th\'eor\`emes de type Bernstein (r\'esultats de platitude dans les dimensions basses) et les th\'eor\`emes de type Moser (r\'esultats de platitude en cons\'equence des estimations globales du gradient).

Pour les graphes minimaux classiques, la contrepartie de Corollaire~\ref{CH_fr_CH:0:CORRigidity_mainthm} a r\'ecemment \'et\'e obtenue par A. Farina dans~\cite{F17}. Dans ce cas, le r\'esultat est optimal et tient avec~$\ell = \min \{ n, 7 \}$.
La preuve du Th\'eor\`eme~\ref{CH_fr_CH:0:under_CRA_Pakmainthm} est bas\'ee sur l'extension au cadre fractionnaire d'une strat\'egie---qui repose sur un r\'esultat de splitting g\'en\'eral pour les blow-downs du sous-graphe $\Sg(u)$---con\c{c}u par A. Farina pour les graphes minimaux classiques et in\'edit. En cons\'equence, les id\'ees contenues dans le Chapitre \ref{CH_Bern_Mos_result} peuvent \^etre utilis\'ees pour obtenir une preuve diff\'erente, plus simple, de~\cite[Theorem~1.1]{F17}

Signalons \'egalement que, en utilisant les m\^emes id\'ees que celles qui conduisent au Th\'eor\`eme~\ref{CH_fr_CH:0:under_CRA_Pakmainthm}, nous pouvons prouver le r\'esultat de rigidit\'e suivant pour ces graphes $s$-mimimaux entiers qui sont situ\'es au-dessus d'un c\^one.

\begin{frtheorem}\label{CH_fr_CH:0:growth_THM}
	Sont~$n\ge1$ un entier et~$s\in(0,1)$. Soit~$u$ une solution de $\h_su=0$ dans $\R^n$, et supposons qu'il existe une constante~$C > 0$ telle que
	\[
	u(x)\ge - C (1+|x|)\quad\mbox{pour chaque }x\in\R^n.
	\]
	Alors,~$u$ est une fonction affine.
\end{frtheorem}

Nous remarquons que dans~\cite{CaCo} on en d\'eduit un r\'esultat de rigidit\'e analogue au Th\'eor\`eme~\ref{CH_fr_CH:0:growth_THM}, sous l'hypoth\`ese plus forte et bilat\'erale
\[
|u(x)| \le C(1 + |x|)\quad\mbox{pour chaque }x\in\R^n.
\]
Le Th\'eor\`eme~\ref{CH_fr_CH:0:growth_THM} am\'eliore donc~\cite[Theorem~1.5]{CaCo} directement.

\medskip

Enfin, nous prouvons que, si $u:\R^n\to\R$ est telle que
\[
\langle\h_su,v\rangle=h\int_{\R^n}v\,dx\quad\mbox{pour chaque }v\in C^\infty_c(\R^n),
\]
pour une certaine constante $h\in\R$, alors la constante doit \^etre $h=0$.

En particulier, en rappelant le Corollaire \ref{CH_fr_CH:0:Equiv_Intro_Global_Corollary}, on voit que, si $u\in W^{s,1}_{\loc}(\R^n)$ est une solution faible de $\h_su=h$ in $\R^n$, alors le sous-graphe de $u$ est localement $s$-minimal dans $\R^{n+1}$.
Cela \'etend au cadre non local un r\'esultat c\'el\`ebre de Chern, \`a savoir le corollaire du Th\'eor\`eme~1 de~\cite{C65}.

\subsection{Un probl\`eme \`a fronti\`ere libre}

Au Chapitre \ref{CH_FreeBdary_CHPT} nous \'etudions les minimiseurs de la fonctionnelle
\eqlab{\label{CH_fr_CH:0:NOnLOC_fReeBdA}
	\Nl(u,\Omega)+\Per\big(\{u>0\},\Omega\big),
}
o\`u $\Nl(u,\Omega)$ est, approximativement, la $\Omega$-contribution \`a la seminorme $H^s$ de la fonction $u:\R^n\to\R$, c'est-\`a-dire
\[
\Nl(u,\Omega):=\iint_{\R^{2n}\setminus(\Co\Omega)^2}\frac{|u(x)-u(y)|^2}{|x-y|^{n+2s}}\,dx\,dy,
\]
pour un certain param\`etre $s\in(0,1)$ fix\'e.

Des fonctionnelles similaires, d\'efinies comme la superposition d'un terme ``\'energie \'elastique'' et d'une ``tension de surface'', ont d\'ej\`a \'et\'e examin\'ees dans les articles suivants:
\begin{itemize}
	\item energie de Dirichlet plus p\'erim\`etre dans \cite{ACKS},
	\item energie de Dirichlet plus p\'erim\`etre fractionnaire dans \cite{CSV},
	\item l'energie non locale $\Nl$ plus le p\'erim\`etre dans \cite{DSV}, et le probl\`eme \`a une phase correspondant dans \cite{DV-onephase}.
\end{itemize}
L'\'etude de la fonctionnelle d\'efinie dans \eqref{CH_fr_CH:0:NOnLOC_fReeBdA} compl\`ete en quelque sorte cette situation.

\smallskip

Les contributions principales du Chapitre \ref{CH_FreeBdary_CHPT} consistent \`a \'etablir une formule de monotonie pour les minimiseurs de la fonctionnelle \eqref{CH_fr_CH:0:NOnLOC_fReeBdA}, \`a l'exploiter pour \'etudier les propri\'et\'es des limites de blow-up et \`a fournir un r\'esultat de r\'eduction de la dimension. De plus, nous montrons que, lorsque $s<1/2$, le p\'erim\`etre domine l'energie non locale. En cons\'equence, nous obtenons un r\'esultat de r\'egularit\'e pour la fronti\`ere libre $\{u=0\}$.

\medskip

En guise de note technique, observons d'abord que nous ne pouvons pas travailler directement avec l'ensemble $\{u>0\}$. Au lieu de cela, nous consid\'erons des \emph{paires admissibles} $(u,E)$, o\`u $u:\R^n\to\R$ est une fonction mesurable, et $E\subseteq\R^n$ est tel que
\[
u\ge0\quad\mbox{p.p. dans }E\quad\mbox{et}\quad u\leq0\quad\mbox{p.p. dans }\Co E.
\]
L'ensemble $E$ est g\'en\'eralement appel\'e \emph{ensemble de positivit\'e} de $u$. Alors, \'etant donn\'ee une valeur $s\in(0,1)$ et un ensemble ouvert ayant fronti\`ere Lipschitz $\Omega\subseteq\R^n$, nous d\'efinissons la fonctionnelle
\[
\F_\Omega(u,E):=\Nl(u,\Omega)+\Per(E,\Omega),
\]
pour chaque paire admissible $(u,E)$.

Remarquons maintenant que, si $u:\R^n\to\R$ est une fonction mesurable, alors
\eqlab{\label{CH_fr_CH:0:TAiL_MiN_ObsV}
	\Nl(u,\Omega)<\infty\quad\implies\quad\int_{\R^n}\frac{|u(\xi)|^2}{1+|\xi|^{n+2s}}\,d\xi<\infty.
}
Pour une preuve, voir par exemple, Lemme \ref{CH:APP:usef_ineq_tail}.
En cons\'equence, nous avons aussi
$$
\int_{\R^n}\frac{|u(\xi)|}{1+|\xi|^{n+2s}}\,d\xi<\infty\quad\mbox{ et }\quad u\in L^2_{\loc}(\R^n).
$$

La notion de minimiseurs que nous consid\'erons est la suivante:

\begin{frdefn}\label{CH_fr_CH:0:def_minim_FrEeBdAry}
	\'Etant donn\'ee une paire admissible~$(u,E)$ telle que $\F_\Omega(u,E)<\infty$, on dit que une paire $(v,F)$ 
	est un
	\emph{concurrent admissible} si \eqlab{\label{CH_fr_CH:0:compet_def}
		&\textrm{supp}(v-u)\Subset\Omega,\qquad F\Delta E\Subset\Omega,\\
		&v-u\in H^s(\R^n)\qquad\textrm{et}\qquad \Per(F,\Omega)<+\infty.
	}
	On dit que une paire admissible $(u,E)$
	est \emph{minimisante} dans $\Omega$ si $\F_\Omega(u,E)<\infty$ et
	\[
	\F_\Omega(u,E)\leq\F_\Omega(v,F),
	\]
	pour chaque concurrent admissible~$(v,F)$.
\end{frdefn}
Notez que la premi\`ere ligne de \eqref{CH_fr_CH:0:compet_def} dit simplement que les paires $(u,E)$ et $(v,F)$ sont \'egales---au sens th\'eorique de la mesure---en dehors d'un sous-ensemble compact de $\Omega$. Donc, puisque $\F_\Omega(u,E)<\infty$, on voit facilement que la deuxi\`eme ligne est \'equivalente \`a $\F_\Omega(v,F)<\infty$.

\smallskip

En particulier, nous nous int\'eressons au probl\`eme de minimisation suivant, par rapport \`a la ``donn\'ee ext\'erieure'' fix\'ee.
\'Etant donn\'ee une paire admissible $(u_0,E_0)$ et un ensemble ouvert born\'e $\Op\subseteq\R^n$ ayant fronti\`ere Lipschitz, tels que
\eqlab{\label{CH_fr_CH:0:Dir_data}
	\Omega\Subset\Op,\qquad \Nl(u_0,\Omega)<+\infty
	\quad\textrm{ et }\quad \Per(E_0,\Op)<+\infty,
}
nous voulons trouver une paire admissible $(u,E)$ atteignant l'infimum suivant
\eqlab{\label{CH_fr_CH:0:Min_dir}
	\inf\big\{\Nl(v,\Omega)+\Per(F,\Op)\,|\,(v,F)&\textrm{ paire admissible t.q. }v=u_0\textrm{ p.p. dans }\Co\Omega\\
	&
	\quad
	\textrm{ et }F\setminus\Omega=E_0\setminus\Omega\big\}.
}
Approximativement, comme d'habitude lorsqu'il s'agit de probl\`emes de minimisation impliquant le p\'erim\`etre classique, nous envisageons un voisinage (fixe)  $\Op$ de $\Omega$ (aussi petit que nous le souhaitons) afin de ``lire'' la donn\'ee sur la fronti\`ere, $\partial E_0\cap\partial\Omega$.

Nous prouvons que, \'etant fix\'ee une donn\'ee ext\'erieure $(u_0,E_0)$ satisfaisant \eqref{CH_fr_CH:0:Dir_data}, il existe une paire $(u,E)$ r\'ealisant l'infimum \eqref{CH_fr_CH:0:Min_dir}. De plus, nous montrons qu'une telle paire $(u,E)$ minimise aussi au sens de la D\'efinition \ref{CH_fr_CH:0:def_minim_FrEeBdAry}.

\smallskip

Un r\'esultat utile consiste \`a \'etablir une estimation uniforme de l'\'energie des paires minimisantes.

\begin{frtheorem}\label{CH_fr_CH:0:TH:unif}
	Soit~$(u,E)$ une paire minimisante dans~$B_2$. 
	Alors
	\[ \iint_{\R^{2n}\setminus (\Co B_1)^2}\frac{|u(x)-u(y)|^2}{
		|x-y|^{n+2s}}\,dx\,dy + \Per(E,B_1)\le 
	C\left(1+\int_{\R^n}\frac{|u(y)|^2}{1+|y|^{n+2s}}\,dy\right),
	\]
	pour une certaine~$C=C(n,s)>0$.
\end{frtheorem}

En particulier, le Th\'eor\`eme \ref{CH_fr_CH:0:TH:unif} est exploit\'e dans la preuve de l'existence d'une limite de blow-up.
Pour cela, nous devons d'abord introduire---par la technique d'extension de \cite{CS07}---la fonctionnelle \'etendue associ\'ee \`a la minimisation de $\F_\Omega$. Nous \'ecrivons
$$
\R^{n+1}_+:=\{(x,z)\in\R^{n+1} \,|\, x\in\R^n,\,z>0\}.
$$
\'Etant donn\'ee une fonction~$u:\R^n\to\R$,
nous consid\'erons la fonction~$\Ue:\R_+^{n+1}\to\R$
d\'efinie via la convolution avec un noyau de Poisson appropri\'e,
\bgs{
	\Ue(\,\cdot\,,z)=u\ast\mathcal K_s(\,\cdot\,,z),\quad\textrm{o\`u}\quad\K_s(x,z):=c_{n,s}\frac{z^{2s}}{(|x|^2+z^2)^{(n+2s)/2}},
}
et~$c_{n,s}>0$ est une constante de normalisation appropri\'ee. Une telle fonction \'etendue $\Ue$ est bien d\'efinie---voir, par exemple, \cite{Extension}---\`a condition que $u:\R^n\to\R$ est telle que
\[
\int_{\R^n}\frac{|u(\xi)|}{1+|\xi|^{n+2s}}\,d\xi<\infty.
\]
\`A la lumi\`ere de \eqref{CH_fr_CH:0:TAiL_MiN_ObsV}, nous pouvons donc consid\'erer la fonction d'extension d'un minimiseur.

Nous utilisons des lettres majuscules, comme $X=(x,z)$, pour d\'esigner les points dans $\R^{n+1}$.
\'Etant donn\'e un ensemble $\Omega\subseteq\R^{n+1}$, nous \'ecrivons
\bgs{
	\Omega_+:=\Omega\cap\{z>0\}\qquad\textrm{et}\qquad\Omega_0:=\Omega\cap\{z=0\}.
}
De plus, nous identifions l'hyperplan $\{z=0\} \simeq \R^n$ via la fonction de projection.
\'Etant donn\'e un ensemble ouvert born\'e $\Omega\subseteq\R^{n+1}$ 
ayant fronti\`ere Lipschitz, tel que $\Omega_0\not=\emptyset$, nous d\'efinissons
\bgs{
	\lf_\Omega(\Vf,F):=c_{n,s}'\int_{\Omega_+}|\nabla\Vf|^2z^{1-2s}\,dX+\Per(F,\Omega_0),
}
pour $\Vf:\R^{n+1}_+\to\R$ et $F\subseteq\R^n\simeq\{z=0\}$
l'ensemble de positivit\'e de la trace de $\Vf$ sur $\{z=0\}$, c'est-\`a-dire
\bgs{
	\Vf\big|_{\{z=0\}}\geq0\quad\textrm{p.p. dans }F\quad\textrm{et}\quad
	\Vf\big|_{\{z=0\}}\leq0\quad\textrm{p.p. dans }\Co F.
}
Nous appellons une telle paire $(\Vf,F)$ une {\emph{paire admissible}}
pour la fonctionnelle \'etendue. Alors, nous introduisons la notion suivante de minimiseur pour la fonctionnelle \'etendue.

\begin{frdefn}\label{CH_fr_CH:0:adm-ext}
	\'Etant donn\'ee une paire admissible~$(\Uc,E)$, telle que $\lf_\Omega(\Uc,E)<\infty$, on dit qu'une paire~$(\Vf,F)$
	est un \emph{concurrent admissible},
	si~$\lf_\Omega(\Vf,F)<\infty$ et
	\begin{equation*}
	\textrm{supp}\,(\Vf-\Uc)\Subset\Omega\qquad
	\textrm{et}\qquad E\Delta F\Subset\Omega_0.
	\end{equation*}
	On dit qu'une paire admissible $(\Uc,E)$ est \emph{minimale} dans~$\Omega$
	si~$\lf_\Omega(\Uc,E)<\infty$ et
	\begin{equation*}
	\lf_\Omega(\Uc,E)\leq\lf_\Omega(\Vf,F),
	\end{equation*}
	pour chaque concurrent admissible $(\Vf,F)$.
\end{frdefn}

Un r\'esultat important consiste \`a montrer qu'un probl\`eme de minimisation appropri\'e impliquant les fonctionnelles \'etendues \'equivaut \`a la minimisation de la fonctionnelle d'origine $\F_\Omega$. Plus pr\'ecis\'ement:

\begin{frprop}\label{CH_fr_CH:0:Local_energy_prop}
	Soit $(u,E)$ une paire admissible pour $\F$, telle que~$\F_{B_R}(u.E)<+\infty$.
	Alors, la paire $(u,E)$ est minimisante dans $B_R$ si et seulement si la paire $(\Ue,E)$ est minimale pour $\lf_\Omega$,
	dans chaque ensemble ouvert born\'e $\Omega\subseteq\R^{n+1}$ 
	ayant fronti\`ere Lipschitz tel que~$\emptyset\not=\Omega_0\Subset B_R$.
\end{frprop}

L'une des principales raisons d'introduire la fonctionnelle \'etendue r\'eside dans le fait qu'elle nous permet d'\'etablir une formule de monotonie de type Weiss pour les minimiseurs.

Nous notons
\bgs{
	\BaLL_r:=\{(x,z)\in\R^{n+1}\,|\,|x|^2+z^2<r^2\}\qquad\textrm{et}\qquad
	\BaLL_r^+:=\BaLL_r\cap\{z>0\}.
}

\begin{frtheorem}[Formule de Monotonie de type Weiss]\label{CH_fr_CH:0:Monotonicity_teo}
	Soit $(u,E)$ une paire minimisante pour $\F$ dans $B_R$ et d\'efinissons la fonction $\Phi_u:(0,R)\to\R$ comme
	\bgs{
		\Phi_u(r):=r^{1-n}&
		\left(c'_{n,s}\int_{\BaLL_r^+}|\nabla\Ue|^2z^{1-2s}\,dX+\Per(E,B_r)\right)\\
		&\qquad
		-c'_{n,s}\Big(s-\frac{1}{2}\Big)r^{-n}\int_{(\partial\BaLL_r)^+}\Ue^2z^{1-2s}\,d\Ha^n.
	}
	Alors, la fonction $\Phi_u$ est croissante dans $(0,R)$.	
	En outre, $\Phi_u$ est constante dans $(0,R)$ si et seulement si l'extension $\Ue$ est homog\`ene de degr\'e $s-\frac{1}{2}$ dans $\BaLL_R^+$
	et $E$ est un c\^one dans $B_R$.
\end{frtheorem}

Ci-dessus, $(\partial\BaLL_r)^+:=\partial\BaLL_r\cap\{z>0\}$.
Pr\'esentons maintenant les paires redimensionn\'ees $(u_\lambda,E_\lambda)$.
\'Etant donn\'e $u:\R^n\to\R$ et $E\subseteq\R^n$, nous d\'efinissons
\bgs{
	u_\lambda(x):=\lambda^{\frac{1}{2}-s}u(\lambda x)\qquad\textrm{et}\qquad E_\lambda:=\frac{1}{\lambda}E,
}
pour chaque $\lambda>0$.
Nous observons que---en raison des propri\'et\'es d'\'echelle de $\F_\Omega$---une paire $(u,E)$ est minimale dans $\Omega$ si et seulement si la paire redimensionn\'ee $(u_\lambda,E_\lambda)$ est minimale dans $\Omega_\lambda$, pour chaque $\lambda>0$.

Nous prouvons la convergence des paires minimisantes dans les conditions appropri\'ees et nous l'exploitons---en m\^eme temps que le Th\'eor\`eme \ref{CH_fr_CH:0:TH:unif}---dans le cas particuli\`erement important de la suite de blow-up.

On dit qu'une paire admissible~$(u,E)$ est un {\emph{c\^one minimisant}}
si elle est une paire minimisante dans $B_R$, pour chaque $R>0$, et elle est telle que~$u$ est homog\`ene de degr\'e $s-\frac12$ et~$E$ est un c\^one 

\begin{frtheorem}\label{CH_fr_CH:0:TH:blow}
	Soit $s>1/2$ et soit~$(u,E)$ une paire minimisante dans~$B_1$,
	avec~$0\in\partial E$. Supposons \'egalement que
	$u\in C^{s-\frac{1}{2}}(B_1)$.
	Alors, il existe un c\^one minimisant~$(u_0,E_0)$
	et une s\'equence~$r_k\searrow 0$ tels que~$u_{r_k}\to u_0$
	dans~$L^\infty_{\loc}(\R^n)$ et~$E_{r_k}\xrightarrow{\loc} E_0$.
\end{frtheorem}
Les propri\'et\'es d'homog\'en\'eit\'e de la limite de blow-up $(u_0,E_0)$ sont une cons\'equence du Th\'eor\`eme \ref{CH_fr_CH:0:Monotonicity_teo}.

Nous soulignons \'egalement que nous \'etablissons des estimations appropri\'ees pour les \'energies de queue des fonctions $u_r$, ce qui nous permet d'affaiblir les hypoth\`eses de~\cite[Theorem~1.3]{DSV}, o\`u les auteurs demandent \`a~$u$ d'\^etre~$C^{s-\frac12}$ dans tout~$\R^n$.

\smallskip

Nous mentionnons maintenant le r\'esultat de r\'eduction dimensionnelle suivant.
Seulement dans le Th\'eor\`eme suivant, red\'efinissons
\[\F_\Omega(u,E):=(c_{n,s}')^{-1}\Nl(u,\Omega)+\Per(E,\Omega).\]
On dit qu'une paire admissible $(u,E)$ est minimisante dans $\R^n$ si cela minimise $\F_\Omega$ dans chaque ensemble ouvert born\'e $\Omega\subseteq\R^n$ ayant fronti\`ere Lipschitz.

\begin{frtheorem}
	Soit $(u,E)$ une paire admissible et d\'efinissons
	\begin{equation*}
	u^\star(x,x_{n+1}):=u(x)\qquad\mbox{et}\qquad
	E^\star:=E\times\R.
	\end{equation*}
	Alors, la paire $(u,E)$ est minimisante dans $\R^n$
	si et seulement si la paire $(u^\star,E^\star)$ est minimisante dans $\R^{n+1}$.
\end{frtheorem}

\smallskip

Enfin, nous observons que dans le cas $s<1/2$, le p\'erim\`etre est en quelque sorte le terme principal de la fonctionnelle $\F_\Omega$. En cons\'equence, nous pouvons prouver le r\'esultat de r\'egularit\'e suivant:

\begin{frtheorem}
	Soit $s\in(0,1/2)$ et soit $(u,E)$ une paire minimisante dans $\Omega$. Suppose que $u\in L^\infty_{\loc}(\Omega)$.
	Alors, $E$ a fronti\`ere presque minimale dans $\Omega$.
	Plus pr\'ecis\'ement, si $x_0\in\Omega$ et $d:=d(x_0,\Omega)/3$,
	alors, pour chaque $r\in(0,d]$ on a
	\bgs{
		\Per(E,B_r(x_0))\le \Per(F,B_r(x_0))+ C\,r^{n-2s},\qquad\forall\,F\subseteq\R^n\mbox{ t.q. }E\Delta F\Subset B_r(x_0),
	}
	o\`u
	\[
	C=C\left(s,x_0,d,\|u\|_{L^\infty(B_{2d}(x_0))},\int_{\R^n}\frac{|u(y)|}{1+|y|^{n+2s}}\,dy\right)>0.
	\]
	Donc
	\begin{itemize}
		\item[(i)] $\partial^*E$ est localement $C^{1,\frac{1-2s}{2}}$ dans $\Omega$,
		
		\item[(ii)] l'ensemble singulier $\partial E\setminus\partial^*E$ est tel que
		\[
		\Ha^\sigma\big((\partial E\setminus\partial^*E)\cap\Omega\big)=0,\qquad\mbox{pour chaque }\sigma>n-8.
		\]
	\end{itemize}
\end{frtheorem}

\smallskip

Nous concluons en disant quelques mots sur le probl\`eme \`a une phase, qui correspond au cas dans lequel $u\ge0$ p.p. dans $\R^n$. M\^eme si ces r\'esultats ne sont pas inclus dans cette th\`ese, ils feront partie de la version finale de l'article sur lequel est bas\'e le Chapitre \ref{CH_FreeBdary_CHPT}.
En suivant les arguments de \cite{DV-onephase}, nous allons prouver que si $(u,E)$ est un minimiseur du probl\`eme \`a une phase dans $B_2$, pour $s>1/2$, et si $0\in\partial E$, alors
$u\in C^{s-\frac{1}{2}}(B_{1/2})$. Notez en particulier que, par le Th\'eor\`eme \ref{CH_fr_CH:0:TH:blow}, ceci garantit l'existence d'une limite de blow-up $(u_0,E_0)$.
De plus, nous \'etablirons des estimations de densit\'e uniforme pour l'ensemble de positivit\'e $E$, des deux c\^ot\'es.

\subsection{La parade de manchots \`a Phillip Island (traitement math\'ematique)}

Le Chapitre \ref{CH_PEngUInS} a pour but de fournir un mod\`ele math\'ematique simple, mais rigoureux, d\'ecrivant la formation de groupes de manchots sur le rivage au coucher du soleil.

\medskip

Les manchots sont incapables de voler, donc ils sont oblig\'es de marcher lorsqu'ils sont \`a terre.
En particulier, ils pr\'esentent des comportements assez sp\'ecifiques dans leur retour aux tani\`eres, qu'il est int\'eressant d'observer et de d\'ecrire analytiquement.
Nous avons observ\'e que les manchots ont tendance \`a se dandiner sur le rivage pour former un groupe suffisamment grand, puis \`a marcher de mani\`ere compacte chez eux.
Le cadre math\'ematique que nous introduisons d\'ecrit ce ph\'enom\`ene en prenant en compte des ``param\`etres naturels'', tels que la vue des manchots et leur vitesse de croisi\`ere.
Le mod\`ele que nous proposons favorise la formation de conglom\'erats de manchots qui se rassemblent, mais permet \'egalement des individus isol\'es et expos\'es.

Le mod\`ele que nous proposons repose sur un ensemble d'\'equations diff\'erentielles ordinaires, avec un nombre de degr\'es de libert\'e variable dans le temps. 
En raison du comportement discontinu de la vitesse des manchots, le traitement math\'ematique (pour obtenir l'existence et l'unicit\'e de la solution) est bas\'e sur une proc\'edure ``stop-and-go''.
Nous utilisons ce cadre pour fournir des exemples rigoureux dans lesquels au moins certains manchots parviennent \`a rentrer chez eux en toute s\'ecurit\'e (il existe aussi des cas dans lesquels certains manchots restent isol\'es).

Pour faciliter l'intuition du mod\`ele, nous pr\'esentons \'egalement quelques simples simulations num\'eriques, qui peuvent \^etre compar\'ees au mouvement r\'eel de la parade des manchots.

\end{otherlanguage}	

\end{chapter}

\begin{chapter}*{Notation and assumptions}
\adjustmtc

For the convenience of the reader, we collect some of the notation and assumptions used throughout the thesis.
	
	\begin{itemize}
		
		\item Unless otherwise stated, $\Omega$ and $\Omega'$ will always denote open sets.
				
		\item Given a set $A\subseteq\R^n$, we use the notation $\Co A$ to denote the complement of $A$ in $\R^n$, that is $\Co A:=\R^n\setminus A$.
		
		\item We write $\chi_E$ to denote the characteristic function of a set $E\subseteq\R^n$.
		
		\item We write $A\Subset B$ to mean that the closure of $A$ is compact in $\R^n$ and $\overline{A}\subseteq B$.
		
		\item In $\R^n$ we will usually write
		$|E|=\mathcal{L}^n(E)$ for the $n$-dimensional Lebesgue measure of a set $E\subseteq\R^n$.
		
		
		\item We write $\Ha^d$ for the $d$-dimensional Hausdorff measure, for any $d\geq0$.
		
		\item We define the dimensional constants
		\begin{equation*}
			\omega_d:=\frac{\pi^\frac{d}{2}}{\Gamma\big(\frac{d}{2}+1\big)},\qquad d\geq0.
		\end{equation*}
		In particular, we remark that $\omega_0=1$ and, if $k\in\mathbb N$, $k\ge1$, then $\omega_k=\mathcal{L}^k(B_1)$ is the volume of the $k$-dimensional unit ball $B_1\subseteq\R^k$
		and $k\,\omega_k=\Ha^{k-1}(\mathbb{S}^{k-1})$ is the surface area of the $(k-1)$-dimensional sphere
		\begin{equation*}
			\mathbb{S}^{k-1}:=\partial B_1=\{x\in\R^k\,|\,|x|=1\}.
		\end{equation*}
		Furthermore, in Chapter \ref{Asympto0_CH_label} we will make use of the notation
		\[
		\varpi_n:=\Ha^{n-1}(\s^{n-1})=n\,\omega_n\quad\textrm{ and }\quad\varpi_0:=0.
		\]
		
		\item By $A_h\xrightarrow{loc}A$ we mean that $\chi_{A_h}\to\chi_A$ in $L^1_{\loc}(\R^n)$,
		i.e. for every bounded open set $\Omega\subseteq\R^n$ we have $|(A_h\Delta A)\cap\Omega|\to0$.

		\item Since
		\[
			|E\Delta F|=0\quad\Longrightarrow\quad \Per(E,\Omega)=\Per(F,\Omega)\quad\textrm{and}\quad \Per_s(E,\Omega)=\Per_s(F,\Omega),
		\]
unless otherwise stated, we implicitly identify sets up to sets of negligible Lebesgue measure. 
		Moreover, whenever needed we will implicitly choose a particular representative for the class of $\chi_E$ in $L^1_{\loc}(\R^n)$,
		as in Remark \ref{CH:1:gmt_assumption}.\\
		We will not make this assumption in Section \ref{CH:1:Section_Irreg_bdary}, since the Minkowski content can be affected even by changes
		in sets of measure zero, that is, in general
		\begin{equation*}
			|\Gamma_1\Delta\Gamma_2|=0\quad\not\Rightarrow\quad
			\overline{\mathcal{M}}^r(\Gamma_1,\Omega)=\overline{\mathcal{M}}^r(\Gamma_2,\Omega)
		\end{equation*}
		(see Section \ref{CH:1:Section_Irreg_bdary} for a more detailed discussion).
		
		\item Given a set $F\subseteq\R^n$,
		the signed distance function $\bar{d}_F$ from $\partial F$, negative inside $F$, is defined as
		\begin{equation*}
		\bar{d}_F(x):=d(x,F)-d(x,\Co F)\qquad\mbox{for every }x\in\R^n,
		\end{equation*}
		where
		\[d(x,A)=\mbox{dist}(x,A):=\inf_{y\in A}|x-y|,\]
		denotes the usual distance from a set $A\subseteq\R^n$.
		For every $r\in\R$ we define the set
		\[
		F_r:=\left\{x\in\R^n\,|\,\bar{d}_F(x)<r\right\}.
		\]
		We also consider the open tubular $\varrho$-neighborhood of $\partial F$,
		\begin{equation*}
			N_\varrho(\partial F):=\{x\in\R^n\,|\,d(x,\partial F)<\varrho\}
			=\left\{|\bar{d}_F|<\varrho\right\},
		\end{equation*}
		for every $\varrho>0$. 
		Given a bounded open set $\Omega\subseteq\R^n$, the constant 
		\[
		r_0=r_0(\Omega)>0
		\]
		will have two different meanings, depending on the regularity of $\partial\Omega$:
		\begin{itemize}
			\item if $\Omega$ has Lipschitz boundary, then $r_0$ has the same meaning as in Proposition \ref{CH:1:bound_perimeter_unif}. Namely, for every $r\in(-r_0,r_0)$ the bounded open set $\Omega_r$ has Lipschitz boundary and the perimeters are uniformly bounded;
			\item if $\Omega$ has $C^2$ boundary, then $r_0$ has the same meaning as in Remark \ref{CH:3:ext_unif_omega}. Namely, the set $\Omega$ satisfies a strict interior and a strict exterior ball condition of radius $2 r_0$ at every point of the boundary.
		\end{itemize}
		For a more detailed discussion, see Appendix
		\ref{CH:1:Appendix_distance_function}
	\end{itemize}
	
	\begin{taggedtheorem}{MTA}[Measure theoretic assumption]\label{CH:1:gmt_assumption}
		Let $E\subseteq\R^n$ be a measurable set. Up to modifications in sets of Lebesgue measure zero, we can assume (see Appendix \ref{CH:1:Appendix_meas_th_bdary} for a detailed discussion) that
		$E$ contains its measure theoretic interior, it does not intersect its measure theoretic exterior and is such that the topological boundary coincides with the measure theoretic boundary. More precisely, we define
		\bgs{
	& E_{int}:=\left\{x\in\R^n\,|\,\exists\,r>0\mbox{ s.t. }|E\cap B_r(x)|=\omega_n r^n\right\},\\
	&
	E_{ext}:=\left\{x\in\R^n\,|\,\exists\,r>0\mbox{ s.t. }|E\cap B_r(x)|=0\right\},
	}
and the measure theoretic boundary
\bgs{
\partial^-E&:=\R^n\setminus\big(E_{int}\cup E_{ext}\big)\\
&
=\left\{x\in\R^n\,|\,0<|E\cap B_r(x)|<\omega_n r^n\mbox{ for every }r>0\right\}.
}
Then we assume that
\[
E_{int}\subseteq E,\qquad E\cap E_{ext}=\emptyset\qquad\textrm{and}\qquad\partial E=\partial^-E.
\]
As detailed in Appendix \ref{CH:1:Appendix_meas_th_bdary},
one way to do this consists in identifying the set $E$ with the set $E^{(1)}$ of points of density one.
	\end{taggedtheorem}

\end{chapter}

\mainmatter

\begin{chapter}{Fractional perimeters from a fractal perspective}\label{CH_Fractals}

\minitoc

%
%

\begin{section}{Introduction and main results}

The purpose of this chapter consists in better understanding the fractional nature of the nonlocal perimeters introduced in \cite{CRS10}. Following \cite{Visintin}, we exploit these fractional perimeters to introduce a definition of fractal dimension
for the measure theoretic boundary of a set.

We calculate the fractal dimension of sets which can be defined in a recursive way and
we give some examples of this kind of sets, explaining how to construct them starting from well known self-similar fractals.
In particular, we show that in the case of the von Koch snowflake $S\subseteq\R^2$ this fractal dimension
coincides with the Minkowski dimension.

We also obtain an optimal result for the asymptotics as $s\to1^-$ of the fractional perimeter
of a set having locally finite (classical) 
perimeter.

\smallskip

Now we give precise statements of the results obtained,
starting with the fractional analysis of fractal dimensions.

\begin{subsection}{Fractal boundaries}\label{CH:1:Section_intro_frac_bd}





We recall that we implicitly assume that all the sets we consider
contain their measure theoretic interior, do not intersect their measure theoretic exterior, and are such that their topological boundary coincides with their measure theoretic boundary---see Remark \ref{CH:1:gmt_assumption} and Appendix \ref{CH:1:Appendix_meas_th_bdary} 
for the details. We will not make this assumption in Section \ref{CH:1:Section_Irreg_bdary}, since the Minkowski content can be affected even by changes
in sets of measure zero.

\smallskip

We recall that we split the fractional perimeter as the sum
\begin{equation*}
\Per_s(E,\Omega)=\Per_s^L(E,\Omega)+\Per_s^{NL}(E,\Omega),
\end{equation*}
where
\begin{equation*}\begin{split}
&\Per_s^L(E,\Omega):=\mathcal L_s(E\cap\Omega,\Co E\cap\Omega)=\frac{1}{2}[\chi_E]_{W^{s,1}(\Omega)},\\
&
\Per_s^{NL}(E,\Omega):=\Ll_s(E\cap\Omega,\Co E\setminus\Omega)+\Ll_s(E\setminus\Omega,\Co E\cap\Omega).
\end{split}\end{equation*}
We can think of $\Per^L_s(E,\Omega)$ as the local part of the fractional perimeter, in the sense that if $|(E\Delta F)\cap\Omega|=0$,
then $\Per^L_s(F,\Omega)=\Per^L_s(E,\Omega)$.

We usually refer to $\Per_s^{NL}(E,\Omega)$ as the nonlocal part of the $s$-perimeter.

We say that a set $E$ has locally finite $s$-perimeter if it has finite $s$-perimeter in
every bounded open set $\Omega\subseteq\R^n$.

When $\Omega=\R^n$, we simply write
\begin{equation*}
\Per_s(E):=\Per_s(E,\R^n)=\frac{1}{2}[\chi_E]_{W^{s,1}(\R^n)}.
\end{equation*}

First of all, we prove in Section \ref{CH:1:Section_supp} that in some sense the measure theoretic boundary $\partial^-E$ is the ``right definition'' of boundary for
working with the $s$-perimeter.

To be more precise, we show that
\begin{equation*}
\partial^-E=\{x\in\R^n\,|\,\Per_s^L(E,B_r(x))>0,\,\forall\,r>0\},
\end{equation*}
and that if $\Omega$ is a connected open set, then
\begin{equation*}
\Per_s^L(E,\Omega)>0\quad\Longleftrightarrow\quad \partial^-E\cap\Omega\not=\emptyset.
\end{equation*}
This can be thought of as an analogue in the fractional framework of the fact that for a Caccioppoli set $E$ we have $\partial^-E=$ supp $|D\chi_E|$.

Now the idea of the definition of the fractal dimension consists in using the index $s$ of $\Per_s^L(E,\Omega)$ to measure the codimension of
$\partial^- E\cap\Omega$,
\begin{equation*}
\Dim_F(\partial^-E,\Omega):=n-\sup\{s\in(0,1)\,|\,\Per^L_s(E,\Omega)<\infty\}.
\end{equation*}

As shown in \cite{Visintin} (Proposition 11 and Proposition 13), the fractal dimension $\textrm{Dim}_F$ defined in this way is related to the (upper) Minkowski dimension (whose precise definition we recall in Definition \ref{CH:1:minkowski_def}) by
\begin{equation}\label{CH:1:intro_dim_ineq}
\Dim_F(\partial^-E,\Omega)\leq\overline{\Dim}_\mathcal M(\partial^-E,\Omega).
\end{equation}

For the convenience of the reader we provide a proof of inequality \eqref{CH:1:intro_dim_ineq} in 
Proposition \ref{CH:1:vis_prop}.

If $\Omega$ is a bounded open set with Lipschitz boundary, 
\eqref{CH:1:intro_dim_ineq} means that
\begin{equation*}
\Per_s(E,\Omega)<\infty\qquad\textrm{for every }s\in\big(0,n-\overline{\Dim}_\mathcal M(\partial^-E,\Omega)\big),
\end{equation*}
since the nonlocal part of the $s$-perimeter of any set $E\subseteq\R^n$ is
\begin{equation*}
\Per_s^{NL}(E,\Omega)\leq2\Per_s(\Omega)<\infty,\qquad\textrm{for every }s\in(0,1).
\end{equation*}

We show that for the von Koch snowflake \eqref{CH:1:intro_dim_ineq}
is actually an equality.

Namely, we prove the following:

\begin{theorem}[Fractal dimension of the von Koch snowflake]\label{CH:1:von_koch_snow}
Let $S\subseteq\R^2$ be the von Koch snowflake. Then
\begin{equation}\label{CH:1:koch1}
\Per_s(S)<\infty,\qquad\forall\,s\in\Big(0,2-\frac{\log4}{\log3}\Big),
\end{equation}
and
\begin{equation}\label{CH:1:koch2}
\Per_s(S)=\infty,\qquad\forall\,s\in\Big[2-\frac{\log4}{\log3},1\Big).
\end{equation}
Therefore
\begin{equation*}
\Dim_F(\partial S)=\Dim_\mathcal{M}(\partial S)=\frac{\log4}{\log3}.
\end{equation*}
\end{theorem}

Actually, exploiting the self-similarity of the von Koch curve, we have
\begin{equation*}
\Dim_F(\partial S,\Omega)=\frac{\log4}{\log3},
\end{equation*}
for every $\Omega$ such that $\partial S\cap\Omega\not=\emptyset$.
In particular, this is true for every $\Omega=B_r(p)$ with $p\in \partial S$ and $r>0$ as small as we want.

\smallskip

We remark that this represents a deep difference between the classical and the fractional perimeter.\\
Indeed, if a set $E$ has (locally) finite perimeter, then by De Giorgi's structure Theorem we know that its reduced boundary $\partial^*E$
is locally $(n-1)$-rectifiable. Moreover $\overline{\partial^*E}=\partial^-E$, so the reduced boundary is, in some sense,
a ``big'' portion of the measure theoretic boundary.

On the other hand, we have seen that there are (open) sets, like the von Koch snowflake, which have a ``nowhere rectifiable'' boundary
(meaning that $\partial^-E\cap B_r(p)$ is not $(n-1)$-rectifiable for every $p\in\partial^-E$ and $r>0$)
and still have finite $s$-perimeter for every $s\in(0,\sigma_0)$.

\subsubsection{Self-similar fractal boundaries}

Our argument for the von Koch snowflake is quite general and can be adapted to
compute the dimension $\Dim_F$ of
all sets which can be constructed in a similar recursive way.

To be more precise, 
we start with a bounded open set $T_0\subseteq\R^n$ with finite perimeter $\Per(T_0)<\infty$, which is, roughly speaking, our basic ``building block''.

Then we go on inductively by adding roto-translations of a scaling of the building block $T_0$,
i.e. sets of the form
\begin{equation*}
T_k^i=F_k^i(T_0):=\mathcal{R}_k^i\big(\lambda^{-k}T_0\big)+x_k^i,
\end{equation*}
where $\lambda>1$, $k\in\mathbb N$, $1\leq i\leq ab^{k-1}$, with $a,\,b\in\mathbb N$,
$\mathcal{R}_k^i\in SO(n)$ and $x_k^i\in\R^n$.
We ask that these sets do not overlap, i.e.
\begin{equation*}
|T^i_k\cap T^j_h|=0,\qquad\textrm{whenever }i\not=j\textrm{ or }k\not=h.
\end{equation*}
Then we define
\begin{equation}\label{CH:1:frac_ind_def}
T_k:=\bigcup_{i=1}^{ab^{k-1}}T_k^i\qquad\textrm{and}\qquad T:=\bigcup_{k=1}^\infty T_k.
\end{equation}
The final set $E$ is either
\begin{equation*}
E:=T_0\cup\bigcup_{k\geq1}\bigcup_{i=1}^{ab^{k-1}}T^i_k,\quad\textrm{or}\quad
E:=T_0\setminus\Big(\bigcup_{k\geq1}\bigcup_{i=1}^{ab^{k-1}}T^i_k\Big).
\end{equation*}

For example, the von Koch snowflake is obtained by adding pieces.\\
Examples obtained by removing the $T_k^i$'s
are the middle Cantor set $E\subseteq\R$, the Sierpinski triangle $E\subseteq\R^2$
and the Menger sponge $E\subseteq\R^3$.

We will consider just the set $T$ and exploit the same argument used for the von Koch snowflake
to compute the fractal dimension related to the $s$-perimeter.\\
However, we observe that the Cantor set, the Sierpinski triangle and the Menger sponge are such that $|E|=0$, i.e. $|T_0\Delta T|=0$.\\
Therefore neither the perimeter nor the $s$-perimeter can detect the fractal nature of the (topological) boundary of $T$
and indeed, since
\begin{equation*}
\Per(T)=\Per(T_0)<\infty,
\end{equation*}
we have $\Per_s(T)<\infty$ for every $s\in(0,1)$.

For example, in the case of the Sierpinski triangle, $T_0$ is an equilateral triangle
and $\partial^-T=\partial T_0$, even if $\partial T$ is a self-similar fractal.

The reason of this situation is that the fractal object is the topological boundary of $T$,
while the $s$-perimeter ``measures'' the measure theoretic boundary, which is regular.
Roughly speaking, the problem is that in these cases there is not room enough to find a small ball $B_k^i=F_k^i(B)\subseteq\Co T$
near each piece $T_k^i$.

Therefore, we will make the additional assumption that
\begin{equation}\label{CH:1:add_frac_self_hp}
\exists\,S_0\subseteq\Co T\quad\textrm{s.t. }|S_0|>0\quad\textrm{and }S_k^i:=F_k^i(S_0)\subseteq\Co T\quad\forall\,k,\,i.
\end{equation}
We remark that it is not necessary to ask that these sets do not overlap.

\begin{theorem}\label{CH:1:fractal_bdary_selfsim_dim}
Let $T\subseteq\R^n$ be a set which can be written as in \eqref{CH:1:frac_ind_def}.
If $\frac{\log b}{\log\lambda}\in(n-1,n)$ and \eqref{CH:1:add_frac_self_hp} holds true, then
\begin{equation*}
\Per_s(T)<\infty,\qquad\forall\,s\in\Big(0,n-\frac{\log b}{\log\lambda}\Big)
\end{equation*}
and
\begin{equation*}
\Per_s(T)=\infty,\qquad\forall\,s\in\Big[n-\frac{\log b}{\log\lambda},1\Big).
\end{equation*}
Thus
\begin{equation*}
\Dim_F(\partial^-T)=\frac{\log b}{\log\lambda}.
\end{equation*}
\end{theorem}

Furthermore, we show how to modify self-similar sets like the Sierpinski triangle, without altering their ``structure'', to obtain new sets which satisfy the hypothesis of Theorem \ref{CH:1:fractal_bdary_selfsim_dim} (see Remark \ref{CH:1:remark_general_recurs}
and the final part of Section \ref{CH:1:Section_self_sim_fr_bdar}).
An example is given in Figure \ref{CH:1:buffo_triang} above.

However, we also remark that the measure theoretic boundary of such
a new set will look quite different from the original fractal (topological) boundary and in general
it will be a mix of smooth parts and unrectifiable parts.

The most interesting examples of 
this kind of sets
are probably represented by bounded sets, because in this case the measure theoretic boundary
does indeed have, in some sense, a ``fractal nature'' (see Remark \ref{CH:1:self_sim_frac_bdry_nat_rmk}).\\
Indeed, if $T$ is bounded, then its boundary $\partial^-T$ is compact. Nevertheless, it has infinite (classical) perimeter
and actually $\partial^-T$ has Minkowski dimension strictly greater than $n-1$, thanks to \eqref{CH:1:intro_dim_ineq}.

However, even unbounded sets can have an interesting behavior. Indeed we obtain the following
\begin{prop}\label{CH:1:expl_farc_prop1}
Let $n\geq2$. For every $\sigma\in(0,1)$ there exists a Caccioppoli set $E\subseteq\R^n$
such that
\begin{equation*}
\Per_s(E)<\infty\qquad\forall\,s\in(0,\sigma)\quad\textrm{and}\quad \Per_s(E)=\infty\qquad\forall\,s\in[\sigma,1).
\end{equation*}
\end{prop}


\noindent
Roughly speaking, the interesting thing about this Proposition is the following. Since $E$ has locally finite perimeter, $\chi_E\in BV_{\loc}(\R^n)$,
it also has locally finite $s$-perimeter for every $s\in(0,1)$,
but 
the global perimeter $\Per_s(E)$ is finite if and only if $s<\sigma<1$.

\end{subsection}

\begin{subsection}{Asymptotics as $s\to1^-$}

In Section \ref{CH:1:Section_intro_frac_bd} we have shown that sets with an irregular, eventually fractal, boundary can have finite $s$-perimeter.

On the other hand, if the set $E$ is ``regular'', then it has finite $s$-perimeter for every $s\in(0,1)$.
Indeed, if $\Omega\subseteq\R^n$ is a bounded open set with Lipschitz boundary (or $\Omega=\R^n$),
then $BV(\Omega)\hookrightarrow W^{s,1}(\Omega)$. As a consequence of this embedding,
we find that
\begin{equation*}
\Per(E,\Omega)<\infty\qquad\Longrightarrow\qquad \Per_s(E,\Omega)<\infty\quad\textrm{for every }s\in(0,1).
\end{equation*}

Actually we can be more precise and obtain a sort of converse,
using only the local part of the $s$-perimeter and adding the condition
\begin{equation*}
\liminf_{s\to1^-}(1-s)\Per^L_s(E,\Omega)<\infty.
\end{equation*}

Indeed one has the following result, which is a combination of
\cite[Theorem 3']{BBM} and \cite[Theorem 1]{Davila}, restricted to characteristic functions:
\begin{theorem}\label{CH:1:Davila_conv_local}
Let $\Omega\subseteq\R^n$ be a bounded open set with Lipschitz boundary. Then
$E\subseteq\R^n$ has finite perimeter in $\Omega$ if and only if $\Per_s^L(E,\Omega)<\infty$
for every $s\in(0,1)$, and
\begin{equation}\label{CH:1:asymptotics_fin_cond}
\liminf_{s\to1}(1-s)\Per_s^L(E,\Omega)<\infty.
\end{equation}
In this case we have
\begin{equation}\label{CH:1:asymptotics_local_part}
\lim_{s\to1}(1-s)\Per_s^L(E,\Omega)=\frac{n\omega_n}{2}K_{1,n}\Per(E,\Omega).
\end{equation}
\end{theorem}
\noindent
We briefly show how to get this result (and in particular why the constant looks like that) from
the two Theorems cited above.
Then we compute the constant $K_{1,n}$ in an elementary way, proving that
\begin{equation*}
\frac{n\omega_n}{2}K_{1,n}=\omega_{n-1}.
\end{equation*}

Moreover we show the following:
\begin{remark}
Condition \eqref{CH:1:asymptotics_fin_cond} is necessary. Indeed, there exist bounded sets (see Example \ref{CH:1:inclusion_counterexample}) having finite $s$-perimeter for every $s\in(0,1)$ which do not have finite perimeter.
This also shows that in general the inclusion
\[
BV(\Omega)\subseteq\bigcap_{s\in(0,1)}W^{s,1}(\Omega)
\]
is strict.
\end{remark}

\begin{example}\label{CH:1:inclusion_counterexample}
Let $0<a<1$ and consider the open intervals $I_k:=(a^{k+1},a^k)$ for every $k\in\mathbb{N}$.
Define $E:=\bigcup_{k\in\mathbb{N}}I_{2k}$, which is a bounded (open) set.
Due to the infinite number of jumps $\chi_E\not\in BV(\mathbb{R})$. However it can be proved that
$E$ has finite $s$-perimeter for every $s\in(0,1)$. We postpone the proof to Section \ref{CH:1:Appendix_Example_proof}.
\end{example}

\begin{remark}
For completeness, we also mention a related result contained in \cite{DFPV13},
where the authors provide an example (Example 2.10) of a bounded set $E\subseteq\R$
which does not have finite $s$-perimeter for any $s\in(0,1)$. In particular,
this example proves that in general the inclusion
\[
\bigcup_{s\in(0,1)} W^{s,1}(\Omega)\subseteq L^1(\Omega)
\]
is strict.
\end{remark}

The main result of Section \ref{CH:1:Section_Asymptsto1} is the following Theorem, which extends the asymptotic convergence
of \eqref{CH:1:asymptotics_local_part} to the whole $s$-perimeter.

\begin{theorem}[Asymptotics]\label{CH:1:asymptotics_teo}
Let $\Omega\subseteq\R^n$ be an open set and let $E\subseteq\R^n$. Then, $E$ has locally finite perimeter in $\Omega$
if and only if $E$ has locally finite $s$-perimeter in $\Omega$ for every $s\in(0,1)$ and
\begin{equation*}
\liminf_{s\to1}(1-s)\Per_s^L(E,\Omega')<\infty,\qquad\forall\,\Omega'\Subset\Omega.
\end{equation*}
If $E$ has locally finite perimeter in $\Omega$, then
\begin{equation*}
\lim_{s\to1}(1-s)\Per_s(E,\Op)=\omega_{n-1}\Per(E,\overline\Op),
\end{equation*}
for every open set $\Op\Subset\Omega$ with Lipschitz boundary. More precisely,
\begin{equation*}
\lim_{s\to1}(1-s)\Per^L_s(E,\Op)=\omega_{n-1}\Per(E,\Op)
\end{equation*}
and
\begin{equation}\label{CH:1:nonlocal_per_asympt}
\lim_{s\to1}(1-s)\Per^{NL}_s(E,\Op)=\omega_{n-1}\Per(E,\partial\Op)=\omega_{n-1}\Ha^{n-1}(\partial^*E\cap\partial\Op).
\end{equation}
\end{theorem}

The proof of Theorem \ref{CH:1:asymptotics_teo} relies only on
\cite[Theorem 3']{BBM}, \cite[Theorem 1]{Davila} and on an appropriate estimate
of what happens in a neighborhood of $\partial\Op$.
The main improvement of the known asymptotics results is the convergence \eqref{CH:1:nonlocal_per_asympt}.

\end{subsection}

\end{section}

\begin{section}{Asymptotics as $s\to1^-$}\label{CH:1:Section_Asymptsto1}

We say that an open set $\Omega\subseteq\R^n$ is an extension domain if there exists a constant $C=C(n,s,\Omega)>0$ such that
for every $u\in W^{s,1}(\Omega)$ there exists $\tilde{u}\in W^{s,1}(\R^n)$ with $\tilde{u}_{|\Omega}=u$
and
\begin{equation*}
\|\tilde{u}\|_{W^{s,1}(\R^n)}\leq C\|u\|_{W^{s,1}(\Omega)}.
\end{equation*}
Every open set with bounded Lipschitz boundary is an extension domain (see \cite{HitGuide} for a proof).
By definition we consider $\R^n$ itself as an extension domain.

We begin with the following embedding.

\begin{prop}\label{CH:1:embedding_prop}
Let $\Omega\subseteq\R^n$ be an extension domain. Then
there exists a constant $C=C(n,s,\Omega)\geq 1$ such that
for every $u:\Omega\longrightarrow\R$
\begin{equation}\label{CH:1:embedding_ineq}
\|u\|_{W^{s,1}(\Omega)}\leq C\|u\|_{BV(\Omega)}.
\end{equation}
In particular we have the continuous embedding
\begin{equation*}
BV(\Omega)\hookrightarrow W^{s,1}(\Omega).
\end{equation*}
\end{prop}

\begin{proof}
The claim is trivially satisfied if the right hand side of \eqref{CH:1:embedding_ineq} is infinite, so let $u\in BV(\Omega)$.
Let $\{u_k\}\subseteq C^\infty(\Omega)\cap BV(\Omega)$ be an approximating sequence as in 
\cite[Theorem 1.17]{Giusti}, that is
\begin{equation*}
\|u-u_k\|_{L^1(\Omega)}\longrightarrow0\qquad\textrm{and}\qquad\lim_{k\to\infty}\int_\Omega|\nabla u_k|\,dx=|Du|(\Omega).
\end{equation*}
We only need to check that the $W^{s,1}$-seminorm of $u$ is bounded by its $BV$-norm.\\
Since $\Omega$ is an extension domain, we know (see \cite[Proposition 2.2]{HitGuide}) that
$\exists\, C(n,s)\geq1$ such that
\begin{equation*}
\|v\|_{W^{s,1}(\Omega)}\leq C\|v\|_{W^{1,1}(\Omega)}.
\end{equation*}
Then
\begin{equation*}
[u_k]_{W^{s,1}(\Omega)}\leq\|u_k\|_{W^{s,1}(\Omega)}\leq C\|u_k\|_{W^{1,1}(\Omega)}
=C\|u_k\|_{BV(\Omega)},
\end{equation*}
and hence, using Fatou's Lemma,
\begin{equation*}\begin{split}
[u]_{W^{s,1}(\Omega)}&\leq\liminf_{k\to\infty}[u_k]_{W^{s,1}(\Omega)}
\leq C\liminf_{k\to\infty}\|u_k\|_{BV(\Omega)}=C\lim_{k\to\infty}\|u_k\|_{BV(\Omega)}\\
&
=C\|u\|_{BV(\Omega)},
\end{split}\end{equation*}
proving \eqref{CH:1:embedding_ineq}.
\end{proof}

Given a set $E\subseteq\R^n$ and $r\in\R$, we denote
\begin{equation*}
E_r:=\{x\in\R^n\,|\,\bar{d}_E(x)<r\},
\end{equation*}
where $\bar{d}_E$ is the signed distance function from $E$ (see Appendix \ref{CH:1:Appendix_distance_function}).

\begin{corollary}\label{CH:1:embedding_fin_per_coroll}
\begin{itemize}
\item[(i)] If $E\subseteq\R^n$ has finite perimeter, i.e. $\chi_E\in BV(\R^n)$, then $E$ has also finite $s$-perimeter for every
$s\in(0,1)$.
\item[(ii)] Let $\Omega\subseteq\R^n$ be a bounded open set with Lipschitz boundary. Then there exists $r_0>0$
such that
\begin{equation}\label{CH:1:unif_bound_lip_frac_per}
\sup_{|r|<r_0}\Per_s(\Omega_r)<\infty.
\end{equation}
\item[(iii)] If $\Omega\subseteq\R^n$ is a bounded open set with Lipschitz boundary, then
\begin{equation*}
\Per_s^{NL}(E,\Omega)\leq 2\Per_s(\Omega)<\infty
\end{equation*}
for every $E\subseteq\R^n$.
\item[(iv)] Let $\Omega\subseteq\R^n$ be a bounded open set with Lipschitz boundary. Then
\begin{equation*}
\Per(E,\Omega)<\infty\qquad\Longrightarrow\qquad \Per_s(E,\Omega)<\infty\quad\textrm{for every }s\in(0,1).
\end{equation*}
\end{itemize}
\end{corollary}

\begin{proof}
Claim $(i)$ 
follows from
\begin{equation*}
\Per_s(E)=\frac{1}{2}[\chi_E]_{W^{s,1}(\R^n)}
\end{equation*}
and
Proposition \ref{CH:1:embedding_prop} with $\Omega=\R^n$.

$(ii)$ Let $r_0$ be as in
Proposition \ref{CH:1:bound_perimeter_unif}
and notice that
\begin{equation*}
\Per(\Omega_r)=\Ha^{n-1}\big(\{\bar{d}_\Omega=r\}\big),
\end{equation*}
so that
\begin{equation*}
\|\chi_{\Omega_r}\|_{BV(\R^n)}=|\Omega_r|+\Ha^{n-1}\big(\{\bar{d}_\Omega=r\}\big).
\end{equation*}
Thus
\begin{equation*}
\sup_{|r|<r_0}\Per_s(\Omega_r)\leq C\Big(|\Omega_{r_0}|+\sup_{|r|<r_0}\Ha^{n-1}\big(\{\bar{d}_\Omega=r\}\big)\Big)<\infty.
\end{equation*}

$(iii)$ Notice that
\begin{equation*}\begin{split}
&\Ll_s(E\cap\Omega,\Co E\setminus\Omega)\leq \Ll_s(\Omega,\Co\Omega)=\Per_s(\Omega),\\
&
\Ll_s(\Co E\cap\Omega,E\setminus\Omega)\leq \Ll_s(\Omega,\Co\Omega)=\Per_s(\Omega),
\end{split}
\end{equation*}
and use \eqref{CH:1:unif_bound_lip_frac_per} (with $\Omega_0=\Omega$).

$(iv)$ The nonlocal part of the $s$-perimeter is finite thanks to $(iii)$. As for the local part, recall that
\begin{equation*}
\Per(E,\Omega)=|D\chi_E|(\Omega)\qquad\textrm{and}\qquad \Per_s^L(E,\Omega)=\frac{1}{2}[\chi_E]_{W^{s,1}(\Omega)},
\end{equation*}
then use Proposition \ref{CH:1:embedding_prop}.
\end{proof}

\begin{subsection}{Asymptotics of the local part of the $s$-perimeter}
We recall the results of \cite{BBM} and \cite{Davila},
which straightforwardly give Theorem \ref{CH:1:Davila_conv_local}.

\begin{theorem}[Theorem 3' of \cite{BBM}]\label{CH:1:bb}
Let $\Omega\subseteq\R^n$ be a smooth bounded domain. Let $u\in L^1(\Omega)$. Then
$u\in BV(\Omega)$ if and only if
\begin{equation*}
\liminf_{n\to\infty}\int_\Omega\int_\Omega\frac{|u(x)-u(y)|}{|x-y|}\varrho_n(x-y)\,dxdy<\infty,
\end{equation*}
and then
\begin{equation*}
\begin{split}
C_1|Du|(\Omega)&\leq\liminf_{n\to\infty}\int_\Omega\int_\Omega\frac{|u(x)-u(y)|}{|x-y|}\varrho_n(x-y)\,dxdy\\
&
\leq\limsup_{n\to\infty}\int_\Omega\int_\Omega\frac{|u(x)-u(y)|}{|x-y|}\varrho_n(x-y)\,dxdy\leq C_2|Du|(\Omega),
\end{split}
\end{equation*}
for some constants $C_1$, $C_2$ depending only on $\Omega$.
\end{theorem}

This result was refined by D\'avila:

\begin{theorem}[Theorem 1 of \cite{Davila} ]
Let $\Omega\subseteq\R^n$ be a bounded open set with Lipschitz boundary. Let $u\in BV(\Omega)$. Then
\begin{equation*}
\lim_{k\to\infty}\int_\Omega\int_\Omega\frac{|u(x)-u(y)|}{|x-y|}\varrho_k(x-y)\,dxdy=K_{1,n}|Du|(\Omega),
\end{equation*}
where
\begin{equation*}
K_{1,n}=\frac{1}{n\omega_n}\int_{\mathbb{S}^{n-1}}|v\cdot e|\,d\sigma(v),
\end{equation*}
with $e\in\R^n$ any unit vector.
\end{theorem}

In the above Theorems $\varrho_k$ is any sequence of radial mollifiers i.e. of functions satisfying
\begin{equation}\label{CH:1:rule1}
\varrho_k(x)\geq0,\quad\varrho_k(x)=\varrho_k(|x|),\quad\int_{\R^n}\varrho_k(x)\,dx=1
\end{equation}
and
\begin{equation}\label{CH:1:rule2}
\lim_{k\to\infty}\int_\delta^\infty\varrho_k(r)r^{n-1}dr=0\quad\textrm{for all }\delta>0.
\end{equation}

In particular, for $R$ big enough, $R>$ diam$(\Omega)$, we can consider
\begin{equation*}
\varrho(x):=\chi_{[0,R]}(|x|)\frac{1}{|x|^{n-1}}
\end{equation*}
and define for any sequence $\{s_k\}\subseteq(0,1),\,s_k\nearrow1$,
\begin{equation*}
\varrho_k(x):=(1-s_k)\varrho(x)c_{s_k}\frac{1}{|x|^{s_k}},
\end{equation*}
where the $c_{s_k}$ are normalizing constants. Then
\begin{equation*}\begin{split}
\int_{\R^n}\varrho_k(x)\,dx&=(1-s_k)c_{s_k}n\omega_n\int_0^R\frac{1}{r^{n-1+s_k}}r^{n-1}\,dr\\
&
=(1-s_k)c_{s_k}n\omega_n\int_0^R\frac{1}{r^{s_k}}\,dr=c_{s_k}n\omega_nR^{1-s_k},
\end{split}
\end{equation*}
and hence taking $c_{s_k}:=\frac{1}{n\omega_n}R^{s_k-1}$ gives \eqref{CH:1:rule1}; notice that
$c_{s_k}\to\frac{1}{n\omega_n}$.\\
Also
\begin{equation*}\begin{split}
\lim_{k\to\infty}\int_\delta^\infty\varrho_k(r)r^{n-1}\,dr&=
\lim_{k\to\infty}(1-s_k)c_{s_k}\int_\delta^R\frac{1}{r^{s_k}}\,dr\\
&
=\lim_{k\to\infty}c_{s_k}(R^{1-s_k}-\delta^{1-s_k})=0,
\end{split}
\end{equation*}
giving \eqref{CH:1:rule2}.
With this choice we obtain
\begin{equation*}
\int_\Omega\int_\Omega\frac{|u(x)-u(y)|}{|x-y|}\varrho_k(x-y)\,dxdy=c_{s_k}(1-s_k)[u]_{W^{s_k,1}(\Omega)}.
\end{equation*}
Then, if $u\in BV(\Omega)$, D\'avila's Theorem gives
\begin{equation}\label{CH:1:limitperimeter}\begin{split}
\lim_{s\to1}(1-s)[u]_{W^{s,1}(\Omega)}&=\lim_{s\to1}\frac{1}{c_s}(c_s(1-s)[u]_{W^{s,1}(\Omega)})\\
&
=n\omega_nK_{1,n}|Du|(\Omega).
\end{split}
\end{equation}

\end{subsection}

\begin{subsection}{Proof of Theorem \ref{CH:1:asymptotics_teo}}

We split the proof of Theorem \ref{CH:1:asymptotics_teo} into several steps, which we believe are interesting on their own.

\begin{subsubsection}{The constant $\omega_{n-1}$}

We need to compute the constant $K_{1,n}$.
Notice that we can choose $e$ in such a way that $v\cdot e=v_n$.\\
Then using spheric coordinates for $\s^{n-1}$ we obtain $|v\cdot e|=|\cos\theta_{n-1}|$
and
\begin{equation*}
d\sigma=\sin\theta_2(\sin\theta_3)^2\ldots(\sin\theta_{n-1})^{n-2}d\theta_1\ldots d\theta_{n-1},
\end{equation*}
with $\theta_1\in[0,2\pi)$ and $\theta_j\in[0,\pi)$ for $j=2,\ldots,n-1$.
Notice that
\begin{equation*}\begin{split}
\Ha^k(\s^k)&=\int_0^{2\pi}\,d\theta_1\int_0^\pi\sin\theta_2\,d\theta_2\ldots
\int_0^\pi(\sin\theta_{k-1})^{k-2}\,d\theta_{k-1}\\
&
=\Ha^{k-1}(\s^{k-1})\int_0^\pi(\sin t)^{k-2}\,dt.
\end{split}
\end{equation*}
Then we get
\begin{equation*}
\begin{split}
\int_{\s^{n-1}}|v\cdot e|&\,d\sigma(v)=\Ha^{n-2}(\s^{n-2})\int_0^\pi(\sin t)^{n-2}|\cos t|\,dt\\
&
=\Ha^{n-2}(\s^{n-2})\Big(\int_0^\frac{\pi}{2}(\sin t)^{n-2}\cos t\,dt-\int_\frac{\pi}{2}^\pi(\sin t)^{n-2}\cos t\,dt\Big)\\
&
=\frac{\Ha^{n-2}(\s^{n-2})}{n-1}\Big(\int_0^\frac{\pi}{2}\frac{d}{dt}(\sin t)^{n-1}\,dt-\int_\frac{\pi}{2}^\pi\frac{d}{dt}(\sin t)^{n-1}\,dt\Big)\\
&
=\frac{2\Ha^{n-2}(\s^{n-2})}{n-1}.
\end{split}
\end{equation*}
Therefore
\begin{equation*}
n\omega_nK_{1,n}=2\frac{\Ha^{n-2}(\s^{n-2})}{n-1}=2\Ll^{n-1}(B_1(0))=2\omega_{n-1},
\end{equation*}
and hence \eqref{CH:1:limitperimeter} becomes
\begin{equation*}
\lim_{s\to1}(1-s)[u]_{W^{s,1}(\Omega)}=2\omega_{n-1}|Du|(\Omega),
\end{equation*}
for any $u\in BV(\Omega)$.


\end{subsubsection}

\begin{subsubsection}{Estimating the nonlocal part of the $s$-perimeter}

The aim of this subsection consists in proving that if $\Omega\subseteq\R^n$ is a bounded open set with Lipschitz boundary and
$E\subseteq\R^n$ has finite perimeter in $\Omega_\beta$, for some $\beta\in(0,r_0)$ and $r_0$ as in Proposition \ref{CH:1:bound_perimeter_unif}, then
\begin{equation}\label{CH:1:asymptotics_nonlocal_estimate}
\limsup_{s\to1}(1-s)\Per_s^{NL}(E,\Omega)
\leq2\omega_{n-1}\lim_{\varrho\to0^+}\Per(E,N_\varrho(\partial\Omega)).
\end{equation}
Actually, we prove something slightly more general than \eqref{CH:1:asymptotics_nonlocal_estimate}. Namely, that to estimate the nonlocal part of the $s$-perimeter
we do not necessarily need to use the sets $\Omega_\varrho$: any ``regular'' approximation of $\Omega$ will do.

More precisely, let
$A_k,\, D_k\subseteq\R^n$ be two sequences of
bounded open sets
with Lipschitz boundary strictly approximating $\Omega$ respectively from the inside and from the outside, that is
\begin{itemize}
\item[(i)] $A_k\subseteq A_{k+1}\Subset\Omega$ and $A_k\nearrow\Omega$, i.e. $\bigcup_k A_k=\Omega$,

\item[(ii)] $\Omega\Subset D_{k+1}\subseteq D_k$ and $D_k\searrow\overline{\Omega}$, i.e. $\bigcap_k D_k=\overline{\Omega}$.
\end{itemize}
We define for every $k$
\begin{equation*}\begin{split}
&\Omega_k^+:=D_k\setminus\overline{\Omega},\qquad\Omega_k^-:=\Omega\setminus\overline{A_k}
\qquad T_k:=\Omega_k^+\cup\partial\Omega\cup\Omega_k^-,\\
&\qquad\qquad d_k:=\min\{d(A_k,\partial\Omega),\,d(D_k,\partial\Omega)\}>0.
\end{split}
\end{equation*}
In particular, we observe that we can consider $\Omega_\varrho$ with $\varrho<0$ in place of $A_k$ and with $\varrho>0$ in place of $D_k$.
Then $T_k$ would be $N_\varrho(\partial\Omega)$ and $d_k=\varrho$.

\begin{prop}\label{CH:1:nonlocalpartappro}
Let $\Omega\subseteq\R^n$ be a bounded open set with Lipschitz boundary and let $E\subseteq\R^n$
be a set having finite perimeter in $D_1$.
Then
\begin{equation*}
\limsup_{s\to1}(1-s)\Per_s^{NL}(E,\Omega)\leq
2\omega_{n-1}\lim_{k\to\infty}\Per(E,T_k).
\end{equation*}
In particular, if $\Per(E,\partial\Omega)=0$, then
\begin{equation*}
\lim_{s\to1}(1-s)\Per_s(E,\Omega)=\omega_{n-1}\Per(E,\Omega).
\end{equation*}
\end{prop}
\begin{proof}
Since $\Omega$ is regular and $\Per(E,\Omega)<\infty$, we already know that
\begin{equation*}
\lim_{s\to1}(1-s)\Per_s^L(E,\Omega)=\omega_{n-1}\Per(E,\Omega).
\end{equation*}
Notice that, since $|D\chi_E|$ is a finite Radon measure on $D_1$ and
$T_k\searrow\partial\Omega$ as $k\nearrow\infty$, we have that
\begin{equation*}
\exists\lim_{k\to\infty}\Per(E,T_k)=\Per(E,\partial\Omega).
\end{equation*}
Consider the nonlocal part of the fractional perimeter,
\begin{equation*}
\Per_s^{NL}(E,\Omega)=\Ll_s(E\cap\Omega,\Co E\setminus\Omega)+\Ll_s(\Co E\cap\Omega,E\setminus\Omega),
\end{equation*}
and take any $k$. Then
\begin{equation*}\begin{split}
\Ll_s(E\cap\Omega,\Co E\setminus\Omega)&=\Ll_s(E\cap\Omega,\Co E\cap\Omega_k^+)+\Ll_s(E\cap\Omega,\Co E\cap(\Co\Omega\setminus D_k))\\
&
\leq\Ll_s(E\cap\Omega,\Co E\cap\Omega_k^+)+\frac{n\omega_n}{s}|\Omega|\frac{1}{d_k^s}\\
&
\leq\Ll_s(E\cap\Omega_k^-,\Co E\cap\Omega_k^+)+2\frac{n\omega_n}{s}|\Omega|\frac{1}{d_k^s}\\
&
\leq\Ll_s(E\cap(\Omega_k^-\cup\Omega_k^+),\Co E\cap(\Omega_k^-\cup\Omega_k^+))+2\frac{n\omega_n}{s}|\Omega|\frac{1}{d_k^s}\\
&
=\Per^L_s(E,T_k)+2\frac{n\omega_n}{s}|\Omega|\frac{1}{d_k^s}.
\end{split}\end{equation*}
Since we can bound the other term in the same way, we get
\begin{equation*}
\Per^{NL}_s(E,\Omega)\leq2\Per^L_s(E,T_k)+4\frac{n\omega_n}{s}|\Omega|\frac{1}{d_k^s}.
\end{equation*}
By hypothesis we know that $T_k$ is a bounded open set with Lipschitz boundary
\begin{equation*}
\partial T_k=\partial A_k\cup\partial D_k.
\end{equation*}
Therefore using \eqref{CH:1:asymptotics_local_part} we have
\begin{equation*}
\lim_{s\to1}(1-s)\Per^L_s(E,T_k)=\omega_{n-1}\Per(E,T_k),
\end{equation*}
and hence
\begin{equation*}
\limsup_{s\to1}(1-s)\Per_s^{NL}(E,\Omega)
\leq 2\omega_{n-1}\Per(E,T_k).
\end{equation*}
Since this holds true for any $k$, we get the claim.
\end{proof}

\end{subsubsection}

\begin{subsubsection}{Convergence in almost every $\Omega_\varrho$}

Having a ``continuous'' approximating sequence (the $\Omega_\varrho$) rather than numerable ones allows us to improve
Proposition \ref{CH:1:nonlocalpartappro}.

We first recall that if $E$ has finite perimeter, then De Giorgi's structure Theorem (see, e.g., \cite[Theorem 15.9]{Maggi}) guarantees in particular that
\begin{equation*}
|D\chi_E|=\Ha^{n-1}\llcorner\partial^*E
\end{equation*}
and hence
\begin{equation*}
\Per(E,B)=\Ha^{n-1}(\partial^*E\cap B)\qquad\textrm{for every Borel set }B\subseteq\R^n,
\end{equation*}
where $\partial^*E$ is the reduced boundary of $E$.

\begin{corollary}\label{CH:1:almost_asympt}
Let $\Omega\subseteq\R^n$ be a bounded open set with Lipschitz boundary and let $r_0$ be as in Proposition
\ref{CH:1:bound_perimeter_unif}. Let $E\subseteq\R^n$ be a set having finite perimeter in $\Omega_\beta$,
for some $\beta\in(0,r_0)$, and define
\begin{equation*}
S:=\left\{\delta\in(-r_0,\beta)\,|\,\Per(E,\partial\Omega_\delta)>0\right\}.
\end{equation*}
Then the set $S$ is at most countable. Moreover
\begin{equation}\label{CH:1:asymptotics_ae_convergence}
\lim_{s\to1}(1-s)\Per_s(E,\Omega_\delta)=\omega_{n-1}\Per(E,\Omega_\delta),
\end{equation}
for every $\delta\in(-r_0,\beta)\setminus S$.
\end{corollary}

\begin{proof}
We observe that
\begin{equation*}
\Per(E,\partial\Omega_\delta)=\Ha^{n-1}(\partial^*E\cap\{\bar{d}_\Omega=\delta\}),
\end{equation*}
for every $\delta\in(-r_0,\beta)$, and
\begin{equation}\label{CH:1:slice_everywhere}
M:=\Ha^{n-1}(\partial^*E\cap(\Omega_\beta\setminus\overline{\Omega_{-r_0}}))\leq \Per(E,\Omega_\beta)<\infty.
\end{equation}
Then we define the sets
\begin{equation*}
S_k:=\Big\{\delta\in(-r_0,\beta)\,|\,\Ha^{n-1}(\partial^*E\cap\{\bar{d}_\Omega=\delta\})>\frac{1}{k}\Big\},
\end{equation*}
for every $k\in\mathbb N$ and we remark that
\begin{equation*}
S=\bigcup_{k\in\mathbb N} S_k.
\end{equation*}
Since by \eqref{CH:1:slice_everywhere} we have
\begin{equation*}
\Ha^{n-1}\Big(\bigcup_{-r_0<\delta<\beta}(\partial^*E\cap\{\bar{d}_\Omega=\delta\})\Big)=M,
\end{equation*}
the number of elements in each $S_k$ is at most
\begin{equation*}
\sharp S_k\leq M\,k.
\end{equation*}
As a consequence the set $S$
is at most countable, as claimed.

Finally, since $\Omega_\delta$ is a bounded open set with Lipschitz boundary for every $\delta\in(-r_0,r_0)$
(see Proposition \ref{CH:1:bound_perimeter_unif}), we obtain \eqref{CH:1:asymptotics_ae_convergence}
by Proposition \ref{CH:1:nonlocalpartappro}.
\end{proof}

\end{subsubsection}

\begin{subsubsection}{Conclusion}

We are now ready to prove Theorem \ref{CH:1:asymptotics_teo}.

\begin{proof}[Proof of Theorem \ref{CH:1:asymptotics_teo}]
We begin by observing that if $E\subseteq\R^n$ and we have two open sets $\Op_1\subseteq\Op_2$,
then
\[
\Per_s(E,\Op_1)\le \Per_s(E,\Op_2).
\]
More precisely, we have
\begin{equation}\label{CH:1:subdomainssplit1}\begin{split}
\Per_s(E,\Op_2)&=\Per_s(E,\Op_1)+\Ll_s\big(E\cap(\Op_2\setminus\Op_1),\Co E\cap(\Op_2\setminus\Op_1)\big)\\
&
+\Ll_s\big(E\cap(\Op_2\setminus\Op_1),\Co E\setminus\Op_2\big)
+\Ll_s\big(\Co E\cap(\Op_2\setminus\Op_1), E\setminus\Op_2\big).
\end{split}\end{equation}
Moreover, we also have
\[
\Per_s^L(E,\Op_1)\le \Per_s(E,\Op_2)\quad\mbox{and}\quad \Per(E,\Op_1)\le \Per(E,\Op_2).
\]
Now suppose that $E$ has locally finite perimeter in $\Omega$ and let $\Omega'\Subset\Omega$. Notice that we can find a bounded open set $\Op$ with Lipschitz boundary, such that
\[
\Omega'\Subset\Op\Subset\Omega.
\]
Since $E$ has finite perimeter in $\Op$, by point $(iv)$ of Corollary \ref{CH:1:embedding_fin_per_coroll}, we know that $E$
has finite $s$-perimeter in $\Op$ (and hence also in $\Omega'\Subset\Op$) for every $s\in(0,1)$. Moreover, by
Theorem \ref{CH:1:Davila_conv_local} we obtain
\[
\liminf_{s\to1}(1-s)\Per_s^L(E,\Omega')\le\liminf_{s\to1}(1-s)\Per_s^L(E,\Op)<\infty.
\]
The converse implication is proved similarly.

\smallskip

Now suppose that $E$ has locally finite perimeter in $\Omega$ and let $\Op\Subset\Omega$ have Lipschitz boundary.
Let $r_0=r_0(\Op)>0$ be as in Proposition \ref{CH:1:bound_perimeter_unif}. Since $\Op\Subset\Omega$, we can find $\beta\in(0,r_0)$ small enough such that $\Op_\beta\Subset\Omega$. Moreover, since $E$ has locally finite perimeter in
$\Omega$, $E$ has finite perimeter in $\Op_\beta$.

Then, by Corollary \ref{CH:1:almost_asympt}, we can find $\delta\in(0,\beta)$ such that $\Per(E,\partial\Op_\delta)=0$ and we have
\begin{equation}\label{CH:1:tarta14}
\lim_{s\to1}(1-s)\Per_s(E,\Op_\delta)=\omega_{n-1}\Per(E,\Op_\delta).
\end{equation}
We also remark that, since $|\partial\Op|=0$, we can rewrite \eqref{CH:1:subdomainssplit1} as 
\begin{equation}\label{CH:1:subdomainssplit2}\begin{split}
\Per_s(E,\Op_\delta)&=\Per_s(E,\Op)+\Per^L_s(E,\Op_\delta\setminus\overline{\Op})\\
&
+\Ll_s\big(E\cap(\Op_\delta\setminus\overline{\Op}),\Co E\setminus\Op_\delta\big)
+\Ll_s\big(\Co E\cap(\Op_\delta\setminus\overline{\Op}), E\setminus\Op_\delta\big).
\end{split}\end{equation}
Let
\[
I_s:=\Ll_s\big(E\cap(\Op_\delta\setminus\overline{\Op}),\Co E\setminus\Op_\delta\big)
+\Ll_s\big(\Co E\cap(\Op_\delta\setminus\overline{\Op}), E\setminus\Op_\delta\big)
\]
and notice that
\begin{equation}\label{CH:1:tarta11}
I_s\le \Per^{NL}_s(E,\Op_\delta).
\end{equation}
Hence, since $\Per(E,\partial\Op_\delta)=0$, by \eqref{CH:1:tarta11} and Proposition \ref{CH:1:nonlocalpartappro} we obtain
\begin{equation}\label{CH:1:tarta12}
\lim_{s\to1}(1-s)I_s=0.
\end{equation}
Furthermore, since $E$ has finite perimeter in $\Op_\delta\setminus\overline{\Op}$, which is a bounded open set with Lipschitz boundary, by \eqref{CH:1:asymptotics_local_part} of Theorem \ref{CH:1:Davila_conv_local}, we find
\begin{equation}\label{CH:1:tarta13}
\lim_{s\to1}(1-s)\Per^L_s(E,\Op_\delta\setminus\overline{\Op})=\omega_{n-1}\Per(E,\Op_\delta\setminus\overline{\Op}).
\end{equation}
Therefore, by \eqref{CH:1:subdomainssplit2}, \eqref{CH:1:tarta14}, \eqref{CH:1:tarta12} and \eqref{CH:1:tarta13}, and exploiting the fact
that $\Per(E,\,\cdot\,)$ is a measure, we get
\begin{equation}\label{CH:1:tarta_fin}\begin{split}
\lim_{s\to1}(1-s)\Per(E,\Op)&=\omega_{n-1}\big(\Per(E,\Op_\delta)-\Per(E,\Op_\delta\setminus\overline{\Op})\big)\\
&
=\omega_{n-1}\Per(E,\overline{\Op}).
\end{split}\end{equation}
Finally, since by \eqref{CH:1:asymptotics_local_part} we know that
\begin{equation}\label{CH:1:tartina}
\lim_{s\to1}(1-s)\Per_s^L(E,\Op)=\omega_{n-1}\Per(E,\Op),
\end{equation}
by \eqref{CH:1:tarta_fin} and \eqref{CH:1:tartina} we obtain
\[
\lim_{s\to1}(1-s)\Per_s^{NL}(E,\Op)=\omega_{n-1}\Per(E,\partial\Op),
\]
concluding the proof of the Theorem.
\end{proof}

\end{subsubsection}

\end{subsection}

\end{section}

\begin{section}{Irregularity of the boundary}\label{CH:1:Section_Irreg_bdary}
\begin{subsection}{The measure theoretic boundary as ``support'' of the local part of the $s$-perimeter}\label{CH:1:Section_supp}

First of all we show that the (local part of the) $s$-perimeter does indeed measure
a quantity related to the measure theoretic boundary.
\begin{lemma}
Let $E\subseteq\R^n$ be a set of locally finite $s$-perimeter. Then
\begin{equation*}
\partial^-E=\{x\in\R^n\,|\,\Per_s^L(E,B_r(x))>0\textrm{ for every }r>0\}.
\end{equation*}
\end{lemma}
\begin{proof}
The claim follows from the following observation. Let $A,\,B\subseteq\R^n$ such that $A\cap B=\emptyset$; then
\begin{equation*}
\Ll_s(A,B)=0\quad\Longleftrightarrow\quad|A|=0\quad\textrm{or}\quad|B|=0.
\end{equation*}
Therefore
\begin{equation*}\begin{split}
x\in\partial^-E&\quad\Longleftrightarrow\quad
|E\cap B_r(x)|>0\textrm{  and }|\Co E\cap B_r(x)|>0\quad\forall\,r>0\\
&
\quad\Longleftrightarrow\quad
\Ll_s(E\cap B_r(x),\Co E\cap B_r(x))>0\quad\forall\,r>0,
\end{split}\end{equation*}
concluding the proof
\end{proof}

This characterization of $\partial^-E$ can be thought of as a fractional analogue
of \eqref{CH:1:support_perimeter}. However we can not really think of $\partial^-E$ as the support of
\begin{equation*}
\Per_s^L(E,\,\cdot\,):\Omega\longmapsto \Per_s^L(E,\Omega),
\end{equation*}
in the sense that, in general
\begin{equation*}
\partial^-E\cap\Omega=\emptyset\quad\not\Rightarrow\quad \Per_s^L(E,\Omega)=0.
\end{equation*}
For example, consider $E:=\{x_n\leq0\}\subseteq\R^n$ and notice that $\partial^-E=\{x_n=0\}$.
Let $\Omega:=B_1(2e_n)\cup B_1(-2e_n)$. Then $\partial^-E\cap\Omega=\emptyset$, but
\begin{equation*}
\Per_s^L(E,\Omega)=\Ll_s(B_1(2e_n),B_1(-2e_n))>0.
\end{equation*}

On the other hand, the only obstacle is the non connectedness of the set $\Omega$ and
indeed we obtain the following
\begin{prop}
Let $E\subseteq\R^n$ be a set of locally finite $s$-perimeter and let $\Omega\subseteq\R^n$ be an open set.
Then
\begin{equation*}
\partial^-E\cap\Omega\not=\emptyset\quad\Longrightarrow\quad \Per_s^L(E,\Omega)>0.
\end{equation*}
Moreover, if $\Omega$ is connected
\begin{equation*}
\partial^-E\cap\Omega=\emptyset\quad\Longrightarrow\quad \Per_s^L(E,\Omega)=0.
\end{equation*}
Therefore, if $\widehat{\mathcal O}(\R^n)$ denotes the family of bounded and connected open sets, then
$\partial^-E$ can be considered as the ``support'' of
\begin{equation*}\begin{split}
\Per_s^L(E,\,\cdot\,):\,&\widehat{\mathcal O}(\R^n)\longrightarrow [0,\infty)\\
&
\Omega\longmapsto \Per_s^L(E,\Omega),
\end{split}\end{equation*}
in the sense that, if $\Omega\in\widehat{\mathcal O}(\R^n)$, then
\begin{equation*}
\Per_s^L(E,\Omega)>0\quad\Longleftrightarrow\quad\partial^-E\cap\Omega\not=\emptyset.
\end{equation*}
\end{prop}
\begin{proof}
Let $x\in\partial^-E\cap\Omega$. Since $\Omega$ is open, we have $B_r(x)\subseteq\Omega$ for some $r>0$
and hence
\begin{equation*}
\Per_s^L(E,\Omega)\geq \Per_s^L(E,B_r(x))>0.
\end{equation*}
Let $\Omega$ be connected and suppose $\partial^-E\cap\Omega=\emptyset$.
Notice that we have the partition of $\R^n$ as $\R^n=E_{ext}\cup\partial^-E\cup E_{int}$ (see Appendix \ref{CH:1:Appendix_meas_th_bdary}). Thus we
can write $\Omega$ as the disjoint union
\begin{equation*}
\Omega=(E_{ext}\cap\Omega)\cup(E_{int}\cap\Omega).
\end{equation*}
However, since $\Omega$ is connected and both $E_{ext}$ and $E_{int}$ are open,
we must have $E_{ext}\cap\Omega=\emptyset$ or $E_{int}\cap\Omega=\emptyset$.
Now, if $E_{ext}\cap\Omega=\emptyset$ (the other case is analogous),
then $\Omega\subseteq E_{int}$ and hence $|\Co E\cap\Omega|=0$. Thus
\begin{equation*}
\Per_s^L(E,\Omega)=\Ll_s(E\cap\Omega,\Co E\cap\Omega)=0,
\end{equation*}
concluding the proof.
\end{proof}

\end{subsection}

\begin{subsection}{A notion of fractal dimension}

Let $\Omega\subseteq\R^n$ be an open set. Then
\begin{equation*}
t>s\qquad\Longrightarrow\qquad W^{t,1}(\Omega)\hookrightarrow W^{s,1}(\Omega),
\end{equation*}
(see, e.g., \cite[Proposition 2.1]{HitGuide}). As a consequence, for every $u\in L^1(\Omega)$ there exists a unique
$R(u)\in[0,1]$ such that
\begin{equation*}
[u]_{W^{s,1}(\Omega)}\quad\left\{\begin{array}{cc}
<\infty,& \forall\,s\in(0,R(u))\\
=\infty, &\forall\,s\in(R(u),1)
\end{array}\right.
\end{equation*}
that is
\begin{equation}\begin{split}\label{CH:1:frac_range}
R(u)&=\sup\left\{s\in(0,1)\,\big|\,[u]_{W^{s,1}(\Omega)}<\infty\right\}\\
&
=\inf\left\{s\in(0,1)\,\big|\,[u]_{W^{s,1}(\Omega)}=\infty\right\}.
\end{split}
\end{equation}

In particular, exploiting this result for characteristic functions, in \cite{Visintin} the author suggested the following definition of fractal dimension.
\begin{defn}
Let $\Omega\subseteq\R^n$ be an open set and let $E\subseteq\R^n$ such that $|E\cap\Omega|<\infty$. If $\partial^- E\cap\Omega\not=\emptyset$, we define
\begin{equation*}
\Dim_F(\partial^- E,\Omega):=n-R(\chi_E),
\end{equation*}
the fractal dimension of $\partial^- E$ in $\Omega$, relative to the fractional perimeter.
If $\Omega=\R^n$, we drop it in the formulas.
\end{defn}

Notice that in the case of sets \eqref{CH:1:frac_range} becomes
\begin{equation}\label{CH:1:frac_range_sets}
\begin{split}
R(\chi_E)&=\sup\left\{s\in(0,1)\,\big|\,\Per_s^L(E,\Omega)<\infty\right\}\\
&
=\inf\left\{s\in(0,1)\,\big|\,\Per_s^L(E,\Omega)=\infty\right\}.
\end{split}
\end{equation}
We observe that, since $\Per^L_s(\Co E,\Omega)=\Per^L_s(E,\Omega)$, in order to define the fractal dimension of $\partial^-E$
in $\Omega$, it is actually enough to require that either $|E\cap\Omega|<\infty$ or $|\Co E\cap\Omega|<\infty$. Clearly,
when the open set $\Omega$ is bounded, such assumptions are trivially satisfied.

In particular we can consider $\Omega$ to be the whole of $\R^n$, or a bounded open set with Lipschitz boundary.
In the first case the local part of the fractional perimeter coincides with the whole fractional perimeter, while in the second case we know that we can bound the nonlocal part with $2\Per_s(\Omega)<\infty$ for every $s\in(0,1)$. Therefore, in both cases in
\eqref{CH:1:frac_range_sets}
we can as well take the whole fractional perimeter $\Per_s(E,\Omega)$ instead of just the local part.

\smallskip

Now we recall the definition of Minkowski dimension, given in terms of the Minkowski contents.
For equivalent definitions of the Minkowski dimension and for the main properties, we refer to \cite{Mattila} and \cite{Falconer} and the references cited therein.

For simplicity, given $\Gamma\subseteq\R^n$ we set
\begin{equation*}
\bar{N}_\varrho^\Omega(\Gamma):=\overline{N_\varrho(\Gamma)}\cap\Omega
=\{x\in\Omega\,|\,d(x,\Gamma)\leq\varrho\},
\end{equation*}
for any $\varrho>0$.

\begin{defn}\label{CH:1:minkowski_def}
Let $\Omega\subseteq\R^n$ be an open set. For any $\Gamma\subseteq\R^n$ and $r\in[0,n]$ we define the inferior and superior $r$-dimensional Minkowski contents of $\Gamma$ relative to the set $\Omega$ as, respectively
\begin{equation*}
\underline{\mathcal{M}}^r(\Gamma,\Omega):=\liminf_{\varrho\to0}\frac{|\bar{N}_\varrho^\Omega(\Gamma)|}{\varrho^{n-r}},\qquad
\overline{\mathcal{M}}^r(\Gamma,\Omega):=\limsup_{\varrho\to0}\frac{|\bar{N}_\varrho^\Omega(\Gamma)|}{\varrho^{n-r}}.
\end{equation*}
Then we define the lower and upper Minkowski dimensions of $\Gamma$ in $\Omega$ as
\begin{equation*}\begin{split}
\underline{\Dim}_\mathcal{M}(\Gamma,\Omega)&:=\inf\big\{r\in[0,n]\,|\,\underline{\mathcal{M}}^r(\Gamma,\Omega)=0\big\}\\
&
=n-\sup\big\{r\in[0,n]\,|\,\underline{\mathcal{M}}^{n-r}(\Gamma,\Omega)=0\big\},
\end{split}\end{equation*}
\begin{equation*}\begin{split}
\overline{\Dim}_\mathcal{M}(\Gamma,\Omega)&:=\sup\big\{r\in[0,n]\,|\,\overline{\mathcal{M}}^r(\Gamma,\Omega)=\infty\big\}\\
&
=n-\inf\big\{r\in[0,n]\,|\,\overline{\mathcal{M}}^{n-r}(\Gamma,\Omega)=\infty\big\}.
\end{split}
\end{equation*}
If they agree, we write
\begin{equation*}
\Dim_\mathcal{M}(\Gamma,\Omega)
\end{equation*}
for the common value and call it the Minkowski dimension of $\Gamma$ in $\Omega$.
If $\Omega=\R^n$ or $\Gamma\Subset\Omega$, we drop the $\Omega$ in the formulas.
\end{defn}

\begin{remark}
Let $\Dim_\mathcal{H}$ denote the Hausdorff dimension. In general one has
\begin{equation*}
\Dim_\mathcal{H}(\Gamma)\leq\underline{\Dim}_\mathcal{M}(\Gamma)\leq\overline{\Dim}_\mathcal{M}(\Gamma),
\end{equation*}
and all the inequalities might be strict (for some examples, see, e.g., \cite[Section 5.3]{Mattila}). However for some sets, like self-similar sets which satisfy appropriate symmetric and regularity conditions, they are all equal (see, e.g., \cite[Corollary 5.8]{Mattila}).
\end{remark}

Now we give a proof of the relation \eqref{CH:1:intro_dim_ineq} (obtained in \cite{Visintin}).
For related results, see also \cite{Sickel} and \cite{FarRog}.

\begin{prop}\label{CH:1:vis_prop}
Let $\Omega\subseteq\R^n$ be a bounded open set. Then for every $E\subseteq\R^n$ such that $\partial^- E\cap\Omega\not=\emptyset$ and $\overline{\Dim}_\mathcal{M}(\partial^-E,\Omega)\geq n-1$ we have
\begin{equation*}
\Dim_F(\partial^-E,\Omega)\leq\overline{\Dim}_\mathcal{M}(\partial^-E,\Omega).
\end{equation*}
\end{prop}
\begin{proof}
By hypothesis we have
\begin{equation*}
\overline{\Dim}_\mathcal{M}(\partial^-E,\Omega)=n-\inf\big\{r\in(0,1)\,|\,\overline{\mathcal{M}}^{n-r}(\partial^-E,\Omega)=\infty\big\},
\end{equation*}
and we need to show that
\begin{equation*}
\inf\big\{r\in(0,1)\,|\,\overline{\mathcal{M}}^{n-r}(\partial^-E,\Omega)=\infty\big\}
\leq
\sup\{s\in(0,1)\,|\,\Per_s^L(E,\Omega)<\infty\}.
\end{equation*}
Up to modifying $E$ on a set of Lebesgue measure zero
we can suppose that $\partial E=\partial^-E$, as in Remark \ref{CH:1:gmt_assumption}. Notice that this does not affect the
$s$-perimeter.

Now for any $s\in(0,1)$
\begin{equation*}\begin{split}
2\Per_s^L(E,\Omega)&=\int_\Omega\,dx\int_\Omega\frac{|\chi_E(x)-\chi_E(y)|}{|x-y|^{n+s}}\,dy\\
&
=\int_\Omega dx\int_0^\infty d\varrho\int_{\partial B_\varrho(x)\cap\Omega}\frac{|\chi_E(x)-\chi_E(y)|}{|x-y|^{n+s}}\,d\Ha^{n-1}(y)\\
&
=\int_\Omega dx\int_0^\infty\frac{d\varrho}{\varrho^{n+s}}\int_{\partial B_\varrho(x)\cap\Omega}|\chi_E(x)-\chi_E(y)|\,d\Ha^{n-1}(y).
\end{split}
\end{equation*}
Notice that
\begin{equation*}
d(x,\partial E)>\varrho\quad\Longrightarrow\quad\chi_E(y)=\chi_E(x),\quad\forall\,y\in\overline{B_\varrho(x)},
\end{equation*}
and hence
\begin{equation*}\begin{split}
\int_{\partial B_\varrho(x)\cap\Omega}|\chi_E(x)-\chi_E(y)|\,d\Ha^{n-1}(y)&
\leq\int_{\partial B_\varrho(x)\cap\Omega}\chi_{\bar{N}_\varrho(\partial E)}(x)\,d\Ha^{n-1}(y)\\
&
\leq n\omega_n\varrho^{n-1}\chi_{\bar{N}_\varrho(\partial E)}(x).
\end{split}
\end{equation*}
Therefore
\begin{equation*}
2\Per_s^L(E,\Omega)\leq n\omega_n\int_0^\infty\frac{d\varrho}{\varrho^{1+s}}\int_\Omega
\chi_{\bar{N}_\varrho(\partial E)}(x)
=n\omega_n\int_0^\infty\frac{|\bar{N}^\Omega_\varrho(\partial E)|}{\varrho^{1+s}}\,d\varrho.
\end{equation*}
We claim that
\begin{equation}\label{CH:1:visintin_proof}
\overline{\mathcal{M}}^{n-r}(\partial E,\Omega)<\infty\quad\Longrightarrow\quad \Per_s^L(E,\Omega)<\infty,\quad\forall\,s\in(0,r).
\end{equation}
Indeed
\begin{equation*}
\limsup_{\varrho\to0}\frac{|\bar{N}^\Omega_\varrho(\partial E)|}{\varrho^r}<\infty\quad\Longrightarrow\quad\exists\,C>0\textrm{ s.t. }
\sup_{\varrho\in(0,C]}\frac{|\bar{N}^\Omega_\varrho(\partial E)|}{\varrho^r}\leq M<\infty.
\end{equation*}
Hence
\begin{equation*}\begin{split}
2\Per_s^L(E,\Omega)&\leq n\omega_n\Big\{\int_0^C\frac{|\bar{N}^\Omega_\varrho(\partial E)|}{\varrho^{1-(r-s)+r}}\,d\varrho
+\int_C^\infty\frac{|\bar{N}^\Omega_\varrho(\partial E)|}{\varrho^{1+s}}\,d\varrho\Big\}\\
&
\leq n\omega_n\Big\{
M\int_0^C\frac{1}{\varrho^{1-(r-s)}}\,d\varrho+|\Omega|\int_C^\infty\frac{1}{\varrho^{1+s}}\,d\varrho
\Big\}\\
&
=n\omega_n\Big\{
\frac{M}{r-s}C^{r-s}+\frac{|\Omega|}{sC^s}
\Big\}<\infty,
\end{split}\end{equation*}
proving \eqref{CH:1:visintin_proof}.
This implies that
\begin{equation*}
r\leq\sup\{s\in(0,1)\,|\,\Per_s^L(E,\Omega)<\infty\},
\end{equation*}
for every $r\in(0,1)$ such that $\overline{\mathcal{M}}^{n-r}(\partial E,\Omega)<\infty$.

Thus, for $\eps>0$ very small, we have
\begin{equation*}
\inf\big\{r\in(0,1)\,|\,\overline{\mathcal{M}}^{n-r}(\partial^-E,\Omega)=\infty\big\}-\eps
\leq\sup\{s\in(0,1)\,|\,\Per_s^L(E,\Omega)<\infty\}.
\end{equation*}
Letting $\eps$ tend to zero, we conclude the proof.
\end{proof}

In particular, if $\Omega$ has Lipschitz boundary we obtain:
\begin{corollary}\label{CH:1:fractal_dim_coroll}
Let $\Omega\subseteq\R^n$ be a bounded open set with Lipschitz boundary. Let $E\subseteq\R^n$ such that $\partial^-E\cap\Omega\not=\emptyset$ and
$\overline{\Dim}_\mathcal{M}(\partial^-E,\Omega)\in[n-1,n)$. Then
\begin{equation*}
\Per_s(E,\Omega)<\infty\qquad\textrm{for every }s\in\left(0,n-\overline{\Dim}_\mathcal{M}(\partial^-E,\Omega)\right).
\end{equation*}
\end{corollary}

\begin{remark}\label{CH:1:fractal_dim_rmk}
Actually, Proposition \ref{CH:1:vis_prop} and Corollary \ref{CH:1:fractal_dim_coroll} still remain true when $\Omega=\R^n$, provided the set $E$ we are considering is bounded.
Indeed, if $E$ is bounded, we can apply the previous results with $\Omega=B_R$ such that $E\Subset\Omega$. Moreover, since $\Omega$
has a regular boundary, as remarked above we can take the whole $s$-perimeter in
\eqref{CH:1:frac_range_sets}, instead of just the local part. But then, since $\Per_s(E,\Omega)=\Per_s(E)$, we see that
\begin{equation*}
\Dim_F(\partial^-E,\Omega)=\Dim_F(\partial^-E,\R^n).
\end{equation*}
\end{remark}

\begin{subsubsection}{Remarks about the Minkowski content of $\partial^-E$}

In the beginning of the proof of Proposition \ref{CH:1:vis_prop}
we chose a particular representative for the class of $E$ in order to have $\partial E=\partial^-E$.
This can be done since it does not affect the $s$-perimeter
and we are already considering the Minkowski dimension of $\partial^-E$.

On the other hand, if we consider a set $F$ such that $|E\Delta F|=0$,
we can use the same proof to obtain the inequality
\begin{equation*}
\Dim_F(\partial^-E,\Omega)\leq\overline{\Dim}_\mathcal{M}(\partial F,\Omega).
\end{equation*}
It is then natural to ask whether we can find a ``better'' representative $F$,
whose (topological) boundary $\partial F$ has Minkowski dimension strictly smaller than that
of $\partial^-E$.

First of all, we remark that the Minkowski content can be influenced by changes in sets of measure zero.
Roughly speaking, this is because the Minkowski content is not a purely measure theoretic notion,
but rather a combination of metric and measure.

For example, let $\Gamma\subseteq\R^n$ and define $\Gamma':=\Gamma\cup\mathbb Q^n$.
Then $|\Gamma\Delta\Gamma'|=0$, but
$N_\delta(\Gamma')=\R^n$ for every $\delta>0$.

In particular, considering different representatives for $E$
we will get different topological boundaries and hence different Minkowski dimensions.

However, since the measure theoretic boundary minimizes the size of the topological boundary, that is
\begin{equation*}
\partial^-E=\bigcap_{|F\Delta E|=0}\partial F,
\end{equation*}
(see Appendix \ref{CH:1:Appendix_meas_th_bdary}), it minimizes also the Minkowski dimension.\\
Indeed, for every $F$ such that $|F\Delta E|=0$ we have
\begin{equation*}\begin{split}
\partial^-E\subseteq\partial F&\quad\Longrightarrow\quad \bar N_\varrho^\Omega(\partial^-E)\subseteq
\bar N_\varrho^\Omega(\partial F)\\
&
\quad\Longrightarrow\quad
\overline{\mathcal M}^r(\partial^-E,\Omega)\leq
\overline{\mathcal M}^r(\partial F,\Omega)\\
&
\quad\Longrightarrow\quad
\overline{\Dim}_\mathcal{M}(\partial^-E,\Omega)\leq
\overline{\Dim}_\mathcal{M}(\partial F,\Omega).
\end{split}
\end{equation*}

\end{subsubsection}

\end{subsection}

\begin{subsection}{Fractal dimension of the von Koch snowflake}

The von Koch snowflake $S\subseteq\R^2$ is an example of a bounded open set with fractal boundary,
for which the Minkowski dimension and the fractal dimension introduced above coincide.

Moreover its boundary is ``nowhere rectifiable'', in the sense that
$\partial S\cap B_r(p)$ is not $(n-1)$-rectifiable for any $r>0$
and $p\in\partial S$.

\smallskip

First of all we recall how to construct the von Koch curve. Then the snowflake is made of three
von Koch curves.

Let $\Gamma_0$ be a line segment of unit length. 
The set $\Gamma_1$ consists of the four segments obtained by removing the middle third of $\Gamma_0$ and replacing it by the other two sides of the equilateral triangle based on the removed segment.\\
We construct $\Gamma_2$ by applying the same procedure to each of the segments in $\Gamma_1$ and so on.
Thus $\Gamma_k$ comes from replacing the middle third of each straight line segment of $\Gamma_{k-1}$ by the other two sides of an equilateral triangle.

As $k$ tends to infinity, the sequence of polygonal curves $\Gamma_k$ approaches a limiting curve $\Gamma$, called the von Koch curve.\\
If we start with an equilateral triangle with unit length side and perform the same construction on all three sides, we obtain the von Koch snowflake $\Sigma$ (see Figure \ref{CH:1:Figure_Koch_flake}).
Let $S$ be the bounded region enclosed by $\Sigma$, so that $S$ is open and $\partial S=\Sigma$. We still
call $S$ the von Koch snowflake.

It can be shown (see, e.g., \cite{Falconer}) that the Hausdorff dimension of the von Koch snowflake is equal to its Minkowski dimension
and
\begin{equation*}
\Dim_\Ha(\Sigma)=\Dim_\mathcal M(\Sigma)=\frac{\log4}{\log3}
\end{equation*}

Now we explain how to construct $S$ in a recursive way and we observe that
\begin{equation*}
\partial^-S=\partial S=\Sigma.
\end{equation*}

\begin{figure}[htbp]
\begin{center}
\includegraphics[width=100mm]{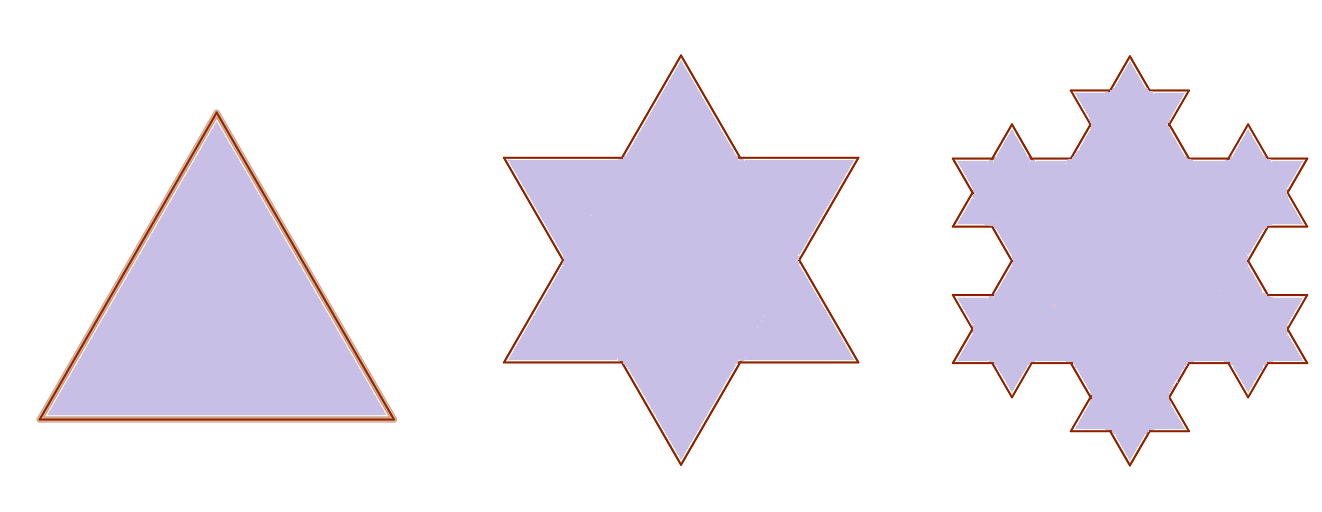}
\caption{{\it The first three steps of the construction of the von Koch snowflake}}
\label{CH:1:Figure_Koch_flake}
\end{center}
\end{figure}

As starting point for the snowflake take the
equilateral triangle $T$ of side 1, with barycenter in the origin and a vertex on the $y$-axis, $P=(0,t)$ with $t>0$.\\
Then $T_1$ is made of three triangles of side $1/3$, $T_2$ of $3\cdot4$ triangles of side $1/3^2$ and so on.\\
In general $T_k$ is made of $3\cdot4^{k-1}$ triangles of side $1/3^k$, call them $T_k^1,\ldots,T_k^{3\cdot4^{k-1}}$.
Let $x^i_k$ be the baricenter of $T_k^i$ and $\Per_k^i$ the vertex which does not touch $T_{k-1}$.

Then $S=T\cup\bigcup T_k$. Also notice that $T_k$ and $T_{k-1}$ touch only on a set of measure zero.

For each triangle
$T^i_k$ there exists a rotation $\mathcal{R}_k^i\in SO(n)$ such that
\begin{equation*}
T_k^i=F_k^i(T):=\mathcal{R}_k^i\Big(\frac{1}{3^k}T\Big)+x_k^i.
\end{equation*}
We choose the rotations so that $F_k^i(P)=\Per_k^i$.

Notice that for each triangle $T_k^i$ we can find a small ball which is contained in the complementary
of the snowflake,
$B_k^i\subseteq\Co S$,
and touches the triangle in the vertex $\Per_k^i$. Actually these balls can be obtained as the images
of the affine transformations $F_k^i$
of a fixed
ball $B$.

To be more precise, fix a small ball contained in the
complementary of $T$, which has the center on the $y$-axis
and touches $T$ in the vertex $P$, say $B:=B_{1/1000}(0,t+1/1000)$. Then
\begin{equation}\label{CH:1:koch3}
B_k^i:=F_k^i(B)\subseteq\Co S
\end{equation}
for every $i,\,k$. To see this, imagine constructing the snowflake $S$ using the same affine transformations $F_k^i$
but starting with $T\cup B$ in place of $T$.

We know that $\partial^-S\subseteq\partial S$ (see Appendix \ref{CH:1:Appendix_meas_th_bdary}).\\
On the other hand, let $p\in\partial S$. Then
every ball $B_\delta(p)$ contains at least a triangle $T^i_k\subseteq S$ and its corresponding ball
$B^i_k\subseteq\Co S$ (and actually infinitely many). Therefore
\[
0<|B_\delta(p)\cap S|<\omega_n\delta^n
\]
for every $\delta>0$ and hence $p\in\partial^-S$.

\begin{proof}[Proof of Theorem \ref{CH:1:von_koch_snow}]
Since $S$ is bounded, its boundary is $\partial^-S=\Sigma$, and
$\Dim_\mathcal M(\Sigma)=\frac{\log4}{\log3}$,
we obtain \eqref{CH:1:koch1} from Corollary \ref{CH:1:fractal_dim_coroll}
and Remark \ref{CH:1:fractal_dim_rmk}.

Exploiting the construction of $S$ given above and \eqref{CH:1:koch3} we prove \eqref{CH:1:koch2}.\\
We have
\begin{equation*}\begin{split}
\Per_s(S)&=\Ll_s(S,\Co S)=\Ll_s(T,\Co S)+\sum_{k=1}^\infty\Ll_s(T_k,\Co S)\\
&
=\Ll_s(T,\Co S)+\sum_{k=1}^\infty\sum_{i=1}^{3\cdot4^{k-1}}\Ll_s(T_k^i,\Co S)
\geq\sum_{k=1}^\infty\sum_{i=1}^{3\cdot4^{k-1}}\Ll_s(T_k^i,\Co S)\\
&
\geq\sum_{k=1}^\infty\sum_{i=1}^{3\cdot4^{k-1}}\Ll_s(T_k^i,B_k^i)\qquad\textrm{(by }\eqref{CH:1:koch3})\\
&
=\sum_{k=1}^\infty\sum_{i=1}^{3\cdot4^{k-1}}\Ll_s(F_k^i(T),F_k^i(B))\\
&
=\sum_{k=1}^\infty\sum_{i=1}^{3\cdot4^{k-1}}\Big(\frac{1}{3^k}\Big)^{2-s}\Ll_s(T,B)\qquad\textrm{(by Proposition }\ref{CH:1:elementary_properties})\\
&
=\frac{3}{3^{2-s}}\Ll_s(T,B)\sum_{k=0}^\infty\Big(\frac{4}{3^{2-s}}\Big)^k.
\end{split}\end{equation*}
We remark that
\begin{equation*}
\Ll_s(T,B)\leq\Ll_s(T,\Co T)=\Per_s(T)<\infty,
\end{equation*}
for every $s\in(0,1)$.

To conclude, notice that the last series is divergent if $s\geq2-\frac{\log4}{\log3}$.
\end{proof}

Exploiting the self-similarity of the von Koch curve, we show that the fractal dimension
of $S$ is the same in every open set which contains a point of $\partial S$.

\begin{corollary}\label{CH:1:koch_coroll}
Let $S\subseteq\R^2$ be the von Koch snowflake. Then
\begin{equation*}
\Dim_F(\partial S,\Omega)=\frac{\log4}{\log3}
\end{equation*}
for every open set $\Omega$ such that $\partial S\cap\Omega\not=\emptyset$.
\end{corollary}

\begin{proof}
Since $\Per_s(S,\Omega)\leq \Per_s(S)$, we have
\begin{equation*}
\Per_s(S,\Omega)<\infty,\qquad\forall\,s\in\Big(0,2-\frac{\log4}{\log3}\Big).
\end{equation*}
On the other hand, if $p\in\partial S\cap\Omega$, then $B_r(p)\subseteq\Omega$
for some $r>0$.
Now notice that $B_r(p)$ contains a rescaled version of the von Koch curve, including
all the triangles $T_k^i$ which constitute it and the relative
balls $B_k^i$. We can thus repeat the argument above to obtain
\begin{equation*}
\Per_s(S,\Omega)\geq \Per_s(S,B_r(p))=\infty,\qquad\forall\,s\in\Big[2-\frac{\log4}{\log3},1\Big),
\end{equation*}
concluding the proof.
\end{proof}

\end{subsection}

\begin{subsection}{Self-similar fractal boundaries}\label{CH:1:Section_self_sim_fr_bdar}

\begin{proof}[Proof of Theorem \ref{CH:1:fractal_bdary_selfsim_dim}]
Arguing as we did with the von Koch snowflake, we show that $\Per_s(T)$ is
bounded both from above and from below by the series
\begin{equation*}
\sum_{k=0}^\infty\Big(\frac{b}{\lambda^{n-s}}\Big)^k,
\end{equation*}
which converges if and only if $s<n-\frac{\log b}{\log\lambda}$.

Indeed
\begin{equation*}\begin{split}
\Per_s(T)&=\Ll_s(T,\Co T)=\sum_{k=1}^\infty\sum_{i=1}^{ab^{k-1}}\Ll_s(T_k^i,\Co T)\\
&
\leq
\sum_{k=1}^\infty\sum_{i=1}^{ab^{k-1}}\Ll_s(T_k^i,\Co T_k^i)
=
\sum_{k=1}^\infty\sum_{i=1}^{ab^{k-1}}\Ll_s(F_k^i(T_0),F_k^i(\Co T_0))\\
&
=\frac{a}{\lambda^{n-s}}\Ll_s(T_0,\Co T_0)\sum_{k=0}^\infty\Big(\frac{b}{\lambda^{n-s}}\Big)^k,
\end{split}\end{equation*}
and
\begin{equation*}\begin{split}
\Per_s(T)&=\Ll_s(T,\Co T)=\sum_{k=1}^\infty\sum_{i=1}^{ab^{k-1}}\Ll_s(T_k^i,\Co T)\\
&
\geq
\sum_{k=1}^\infty\sum_{i=1}^{ab^{k-1}}\Ll_s(T_k^i,S_k^i)
=
\sum_{k=1}^\infty\sum_{i=1}^{ab^{k-1}}\Ll_s(F_k^i(T_0),F_k^i(S_0))\\
&
=\frac{a}{\lambda^{n-s}}\Ll_s(T_0,S_0)\sum_{k=0}^\infty\Big(\frac{b}{\lambda^{n-s}}\Big)^k.
\end{split}\end{equation*}
Also notice that, since $\Per(T_0)<\infty$, we have
\begin{equation*}
\Ll_s(T_0,S_0)\leq\Ll_s(T_0,\Co T_0)=\Per_s(T_0)<\infty,
\end{equation*}
for every $s\in(0,1)$.
\end{proof}

Now suppose that $T$ does not satisfy \eqref{CH:1:add_frac_self_hp}.
Then we can obtain a set $T'$ which does, simply by removing
a portion $S_0$ from the building block $T_0$.\\
To be more precise, let $S_0\subseteq T_0$ be such that
\begin{equation*}
|S_0|>0,\quad|T_0\setminus S_0|>0\quad\textrm{ and }\quad \Per(T_0\setminus S_0)<\infty.
\end{equation*}
Then define a new
building block $T'_0:=T_0\setminus S_0$ and the set
\begin{equation*}
T':=\bigcup_{k=1}^\infty\bigcup_{i=1}^{ab^{k-1}}F_k^i(T'_0).
\end{equation*}
This new set has exactly the same structure of $T$, since we are using the same
collection $\{F_k^i\}$ of affine maps. 

Notice that
\begin{equation*}
S_0\subseteq T_0\quad\Longrightarrow\quad F_k^i(S_0)\subseteq F_k^i(T_0),
\end{equation*}
and
\begin{equation*}
F_k^i(T'_0)=F_k^i(T_0)\setminus F_k^i(S_0),
\end{equation*}
for every $k,\,i$.
Thus
\begin{equation*}
T'=T\setminus\Big(\bigcup_{k=1}^\infty\bigcup_{i=1}^{ab^{k-1}}F_k^i(S_0)\Big)
\end{equation*}
satisfies \eqref{CH:1:add_frac_self_hp}.

\begin{remark}\label{CH:1:remark_general_recurs}
Roughly speaking, what matters in order to obtain a set which satisfies the hypothesis of Theorem  \ref{CH:1:fractal_bdary_selfsim_dim}
is that there exists a bounded open set $T_0$ such that
\begin{equation*}
|F_k^i(T_0)\cap F_h^j(T_0)|=0,\qquad\textrm{if }i\not=j\textrm{ or }k\not=h.
\end{equation*}
This can be thought of as a compatibility criterion for
the family of affine maps $\{F_k^i\}$.
We also need to ask that the ratio of the logarithms
of the growth factor and the scaling factor is $\frac{\log b}{\log\lambda}\in(n-1,n)$.
Then we are free to choose
as building block any set $T'_0\subseteq T_0$ such that
\begin{equation*}
|T'_0|>0,\qquad|T_0\setminus T'_0|>0\qquad\textrm{and }\Per(T'_0)<\infty,
\end{equation*}
and the set
\begin{equation*}
T':=\bigcup_{k=1}^\infty\bigcup_{i=1}^{ab^{k-1}}F_k^i(T'_0).
\end{equation*}
satisfies the hypothesis of Theorem \ref{CH:1:fractal_bdary_selfsim_dim}.
\end{remark}

Therefore, even if the Sierpinski triangle and the Menger sponge
do not satisfy \eqref{CH:1:add_frac_self_hp},
we can exploit their structure to construct new sets which do.

However, we remark that
the new boundary $\partial^-T'$ will look very different from the original fractal. Actually, in general
it will be a mix of unrectifiable pieces and smooth pieces. In particular, we can not hope to get an
analogue of Corollary \ref{CH:1:koch_coroll}.
Still, the following Remark shows that the new (measure theoretic) boundary retains at least some of the ``fractal nature'' of the original set.

\begin{remark}\label{CH:1:self_sim_frac_bdry_nat_rmk}
If the set $T$ of Theorem \ref{CH:1:fractal_bdary_selfsim_dim} is bounded, exploiting Proposition \ref{CH:1:vis_prop} and Remark \ref{CH:1:fractal_dim_rmk}
we obtain
\begin{equation*}
\overline{\Dim}_{\mathcal M}(\partial^-T)\geq\frac{\log b}{\log\lambda}>n-1.
\end{equation*}
Moreover, notice that if $\Omega$ is a bounded open set with Lipschitz boundary, then
\begin{equation*}
\Per(E,\Omega)<\infty\quad\Longrightarrow\quad\Dim_F(E,\Omega)=n-1.
\end{equation*}
Therefore, if $T\Subset B_R$, then
\begin{equation*}
\Per(T)=\Per(T,B_R)=\infty,
\end{equation*}
even if $T$ is bounded (and hence $\partial^-T$ is compact).
\end{remark}

\begin{subsubsection}{Sponge-like sets}

The simplest way to construct the set $T'$ consists in simply removing a small ball
$S_0:=B\Subset T_0$ from $T_0$.

In particular, suppose that $|T_0\Delta T|=0$, as with the Sierpinski triangle.\\
Define
\begin{equation*}
S:=\bigcup_{k=1}^\infty\bigcup_{i=1}^{ab^{k-1}}F_k^i(B)
\quad\textrm{and}\quad
T':=\bigcup_{k=1}^\infty\bigcup_{i=1}^{ab^{k-1}}F_k^i(T_0\setminus B)=T\setminus S.
\end{equation*}
Then
\begin{equation}\label{CH:1:fractal_spazz_end}
|T_0\Delta T|=0\quad\Longrightarrow\quad |T'\Delta (T_0\setminus S)|=0.
\end{equation}
Now the set $E:=T_0\setminus S$ looks like a sponge, in the sense that it
is a bounded open set with an infinite number of holes (each one at a positive, but non-fixed distance from the others).

From \eqref{CH:1:fractal_spazz_end} we get $\Per_s(E)=\Per_s(T')$. Thus, since $T'$ satisfies the hypothesis of
Theorem \ref{CH:1:fractal_bdary_selfsim_dim}, we obtain
\begin{equation*}
\Dim_F(\partial^-E)=\frac{\log b}{\log\lambda}.
\end{equation*}

\end{subsubsection}

\begin{subsubsection}{Dendrite-like sets}

Depending on the form of the set $T_0$ and on the affine maps $\{F_k^i\}$,
we can define more intricate sets $T'$.

As an example we consider the Sierpinski triangle $E\subseteq\R^2$.\\
It is of the form $E=T_0\setminus T$,
where the building block $T_0$ is an equilateral triangle, say with side length one,
a vertex on the $y$-axis and baricenter in 0.
The pieces $T_k^i$ are obtained with a scaling factor $\lambda=2$
and the growth factor is $b=3$ (see, e.g., \cite{Falconer} for the construction).
As usual, we consider the set
\begin{equation*}
T=\bigcup_{k=1}^\infty\bigcup_{i=1}^{3^{k-1}}T_k^i.
\end{equation*}
However, as remarked above, we have $|T\Delta T_0|=0$.

Starting from $k=2$ each triangle $T_k^i$ touches with (at least) a vertex (at least) another triangle $T_h^j$.
Moreover, each triangle $T_k^i$ gets touched in the middle point of each side (and actually it gets touched in infinitely many points).

Exploiting this situation, we can remove from $T_0$ six smaller triangles, so that the new building block
$T'_0$ is a star polygon centered in 0, with six vertices, one in each vertex of $T_0$ and one in each middle point of the sides of $T_0$.

\begin{figure}[htbp]
\begin{center}
\includegraphics[width=90mm]{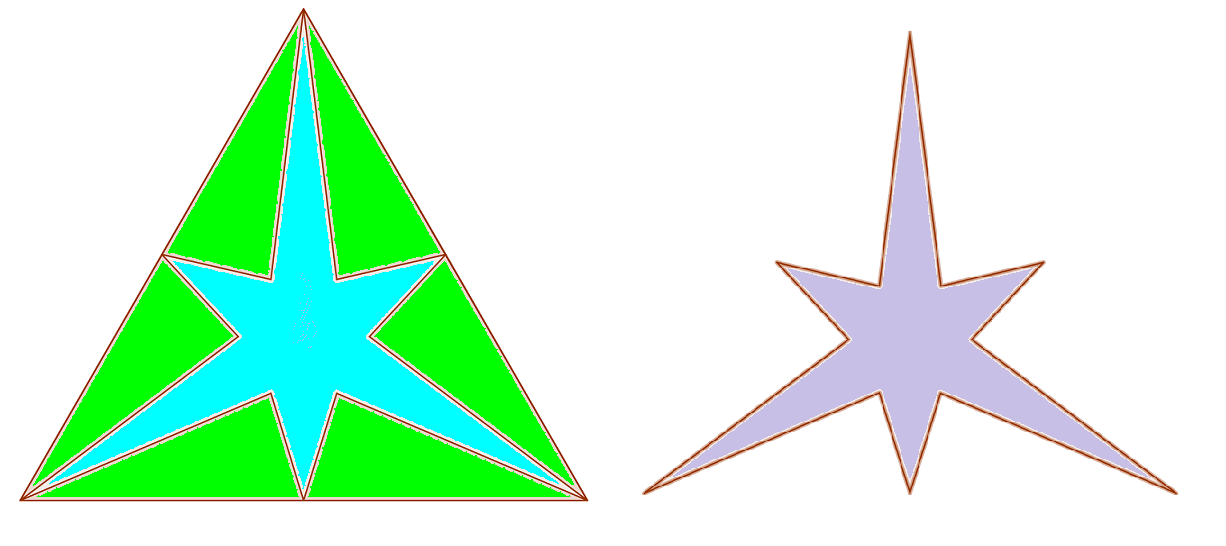}
\caption{{\it Removing the six triangles (in green) to obtain the new ``building block'' $T'_0$ (on the right)}}
\end{center}
\end{figure}

The resulting set
\begin{equation*}
T'=\bigcup_{k=1}^\infty\bigcup_{i=1}^{3^{k-1}}F_k^i(T'_0)
\end{equation*}
will have an infinite number of ramifications.

\begin{figure}[htbp]
\begin{center}
\includegraphics[width=110mm]{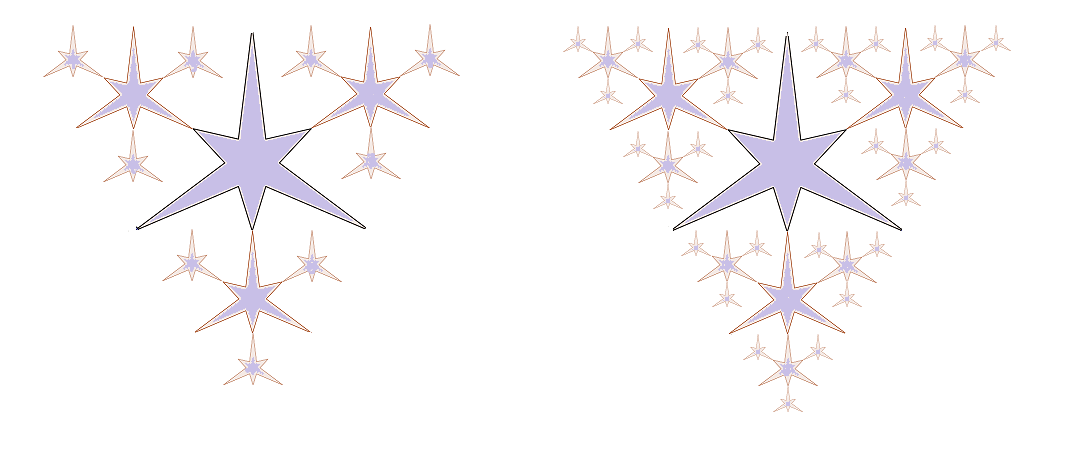}
\caption{{\it The third and fourth steps of the iterative construction of the set $T'$}}
\end{center}
\end{figure}

Since $T'$ satisfies the hypothesis of
Theorem \ref{CH:1:fractal_bdary_selfsim_dim}, we obtain
\begin{equation*}
\Dim_F(\partial^-T')=\frac{\log 3}{\log2}.
\end{equation*}

\end{subsubsection}

\begin{subsubsection}{``Exploded'' fractals}

In all the previous examples, the sets $T_k^i$ are accumulated in a bounded region.

On the other hand, imagine making a fractal like the von Koch snowflake or the Sierpinski triangle
``explode'' and then rearrange the pieces $T_k^i$ in such a way that $d(T_k^i,T_h^j)\geq d$,
for some fixed $d>0$.

Since the shape of the building block is not important, we can consider $T_0:=B_{1/4}(0)\subseteq\R^n$,
with $n\geq2$. Moreover, since the parameter $a$ does not influence the dimension, we can fix $a=1$.

Then we rearrange the pieces obtaining
\begin{equation}\label{CH:1:exploded_frac_def}
E:=\bigcup_{k=1}^\infty\bigcup_{i=1}^{b^{k-1}}B_\frac{1}{4\lambda^k}(k,0,\ldots,0,i).
\end{equation}
Define for simplicity
\begin{equation*}
B_k^i:=B_\frac{1}{4\lambda^k}(k,0,\ldots,0,i)\quad\textrm{and}\quad x_k^i:=k\,e_1+i\,e_n,
\end{equation*}
and notice that
\begin{equation*}
B_k^i=\lambda^{-k}B_\frac{1}{4}(0)+x_k^i.
\end{equation*}
Since for every $k,\,h$ and every $i\not=j$ we have
\begin{equation*}
d(B_k^i,B_h^j)\geq\frac{1}{2},
\end{equation*}
the boundary of the set $E$ is the disjoint union of $(n-1)$-dimensional spheres
\begin{equation*}
\partial^-E=\partial E=\bigcup_{k=1}^\infty\bigcup_{i=1}^{b^{k-1}}\partial B_k^i,
\end{equation*}
and in particular is smooth.

The (global) perimeter of $E$ is
\begin{equation*}
\Per(E)=\sum_{k=1}^\infty\sum_{i=1}^{b^{k-1}}\Per(B_k^i)=\frac{1}{\lambda}\Per(B_{1/4}(0))\sum_{k=0}^\infty
\Big(\frac{b}{\lambda^{n-1}}\Big)^k=+\infty,
\end{equation*}
since $\frac{\log b}{\log\lambda}>n-1$.

However $E$ has locally finite perimeter, since its boundary is smooth and every ball $B_R$ intersects
only finitely many $B_k^i$'s,
\begin{equation*}
\Per(E,B_R)<\infty,\qquad\forall\,R>0.
\end{equation*}
Therefore it also has locally finite $s$-perimeter for every $s\in(0,1)$
\begin{equation*}
\Per_s(E,B_R)<\infty,\qquad\forall\,R>0,\qquad\forall\,s\in(0,1).
\end{equation*}

What is interesting is that the set $E$ satisfies the hypothesis of Theorem \ref{CH:1:fractal_bdary_selfsim_dim} and hence it also has finite global $s$-perimeter
for every $s<\sigma_0:=n-\frac{\log b}{\log\lambda}$,
\begin{equation*}
\Per_s(E)<\infty\qquad\forall\,s\in(0,\sigma_0)\quad\textrm{and}\quad \Per_s(E)=\infty\qquad\forall\,s\in[\sigma_0,1).
\end{equation*}

Thus we obtain Proposition \ref{CH:1:expl_farc_prop1}.

\begin{proof}[Proof of Proposition \ref{CH:1:expl_farc_prop1}]
It is enough to choose a natural number $b\geq2$ and take $\lambda:=b^\frac{1}{n-\sigma}$. Notice that
$\lambda>1$ and
\begin{equation*}
\frac{\log b}{\log\lambda}=n-\sigma\in(n-1,n).
\end{equation*}
Then we can define $E$ as in \eqref{CH:1:exploded_frac_def} and we are done.
\end{proof}

\end{subsubsection}

\end{subsection}

\begin{subsection}{Elementary properties of the $s$-perimeter}

In the following Proposition we collect some elementary but useful properties of the fractional perimeter
which we have exploited throughout the chapter.

\begin{prop}\label{CH:1:elementary_properties}
Let $\Omega\subseteq\R^n$ be an open set.
\begin{itemize}
\item[(i)] (Subadditivity) Let $E,\,F\subseteq\R^n$ be such that $|E\cap F|=0$. Then
\begin{equation*}
\Per_s(E\cup F,\Omega)\leq \Per_s(E,\Omega)+\Per_s(F,\Omega).
\end{equation*}

\item[(ii)] (Translation invariance) Let $E\subseteq\R^n$ and $x\in\R^n$. Then
\begin{equation*}
\Per_s(E+x,\Omega+x)=\Per_s(E,\Omega).
\end{equation*}

\item[(iii)] (Rotation invariance) Let $E\subseteq\R^n$ and $\mathcal{R}\in SO(n)$ a rotation. Then
\begin{equation*}
\Per_s(\mathcal{R}E,\mathcal{R}\Omega)=\Per_s(E,\Omega).
\end{equation*}

\item[(iv)] (Scaling) Let $E\subseteq\R^n$ and $\lambda>0$. Then
\begin{equation*}
\Per_s(\lambda E,\lambda\Omega)=\lambda^{n-s}\Per_s(E,\Omega).
\end{equation*}
\end{itemize}
\end{prop}

\begin{proof}
(i) follows from the following observations. Let $A_1,\,A_2,\,B\subseteq\R^n$. If $|A_1\cap A_2|=0$, then
\begin{equation*}
\Ll_s(A_1\cup A_2,B)
=\Ll_s(A_1,B)+\Ll_s(A_2,B).
\end{equation*}
Moreover
\begin{equation*}
A_1\subseteq A_2\quad\Longrightarrow\quad\Ll_s(A_1,B)\leq\Ll_s(A_2,B),
\end{equation*}
and
\begin{equation*}
\Ll_s(A,B)=\Ll_s(B,A).
\end{equation*}
Therefore
\begin{equation*}\begin{split}
\Per_s(E\cup F,\Omega)&=\Ll_s((E\cup F)\cap\Omega,\Co(E\cup F))+\Ll_s((E\cup F)\setminus\Omega,\Co(E\cup F)\cap\Omega)\\
&
=\Ll_s(E\cap\Omega,\Co(E\cup F))+\Ll_s(F\cap\Omega,\Co(E\cup F))\\
&
\qquad+\Ll_s(E\setminus\Omega,\Co(E\cup F)\cap\Omega)+\Ll_s(F\setminus\Omega,\Co(E\cup F)\cap\Omega)\\
&
\leq\Ll_s(E\cap\Omega,\Co E)+\Ll_s(F\cap\Omega,\Co F)\\
&
\qquad+\Ll_s(E\setminus\Omega,\Co E\cap\Omega)+\Ll_s(F\setminus\Omega,\Co F\cap\Omega)\\
&
=\Per_s(E,\Omega)+\Per_s(F,\Omega).
\end{split}\end{equation*}

(ii), (iii) and (iv) follow simply by changing variables in $\Ll_s$ and the following observations:
\begin{equation*}\begin{split}
&(x+A_1)\cap(x+A_2)=x+A_1\cap A_2,\qquad x+\Co A=\Co(x+A),\\
&
\mathcal{R}A_1\cap\mathcal{R}A_2=\mathcal{R}(A_1\cap A_2),\qquad\mathcal{R}(\Co A)=\Co(\mathcal{R}A),\\
&
(\lambda A_1)\cap(\lambda A_2)=\lambda(A_1\cap A_2),\qquad\lambda(\Co A)=\Co(\lambda A).
\end{split}\end{equation*}
For example, for claim (iv) we have
\begin{equation*}\begin{split}
\Ll_s(\lambda A,\lambda B)&=\int_{\lambda A}\int_{\lambda B}\frac{dx\,dy}{|x-y|^{n+s}}
=\int_A\lambda^n\,dx\int_B\frac{\lambda^n\,dy}{\lambda^{n+s}|x-y|^{n+s}}\\
&
=\lambda^{n-s}\Ll_s(A,B).
\end{split}
\end{equation*}
Then
\begin{equation*}\begin{split}
\Per_s(\lambda E,\lambda\Omega)&=\Ll_s(\lambda E\cap\lambda\Omega,\Co(\lambda E))+
\Ll_s(\lambda E\cap\Co(\lambda\Omega),\Co(\lambda E)\cap\lambda\Omega)\\
&
=\Ll_s(\lambda(E\cap\Omega),\lambda\Co E)+\Ll_s(\lambda(E\setminus\Omega),\lambda(\Co E\cap\Omega))\\
&
=\lambda^{n-s}\left(\Ll_s(E\cap\Omega,\Co E)+\Ll_s(E\setminus\Omega,\Co E\cap\Omega)\right)\\
&
=\lambda^{n-s}\Per_s(E,\Omega).
\end{split}\end{equation*}
This concludes the proof of the Proposition.
\end{proof}

\end{subsection}

\end{section}

\section{Proof of Example \ref{CH:1:inclusion_counterexample}}\label{CH:1:Appendix_Example_proof}

Note that $E\subseteq (0,a^2]$.
Let $\Omega:=(-1,1)\subseteq\mathbb{R}$. Then $E\Subset\Omega$ and $\textrm{dist}(E,\partial\Omega)=1-a^2=:d>0$.
Now
\begin{equation*}
\Per_s(E)=\int_E\int_{\Co E\cap\Omega}\frac{dxdy}{|x-y|^{1+s}}+
\int_E\int_{\Co\Omega}\frac{dxdy}{|x-y|^{1+s}}.
\end{equation*}
As for the second term, we have
\begin{equation*}
\int_E\int_{\Co\Omega}\frac{dxdy}{|x-y|^{1+s}}\leq\frac{2|E|}{sd^s}<\infty.
\end{equation*}
We split the first term into three pieces
\begin{equation*}\begin{split}
\int_E&\int_{\Co E\cap\Omega}\frac{dxdy}{|x-y|^{1+s}}\\
&
=\int_E\int_{-1}^0\frac{dxdy}{|x-y|^{1+s}}
+\int_E\int_{\Co E\cap(0,a)}\frac{dxdy}{|x-y|^{1+s}}+\int_E\int_a^1\frac{dxdy}{|x-y|^{1+s}}\\
&
=\mathcal{I}_1+\mathcal{I}_2+\mathcal{I}_3.
\end{split}
\end{equation*}
Note that $\Co E\cap(0,a)=\bigcup_{k\in\mathbb{N}}I_{2k-1}=\bigcup_{k\in\mathbb{N}}(a^{2k},a^{2k-1})$.\\
A simple calculation shows that, if $a<b\leq c<d$, then
\begin{equation}\label{CH:1:rectangle_integral}\begin{split}
\int_a^b&\int_c^d\frac{dxdy}{|x-y|^{1+s}}=\\
&
\frac{1}{s(1-s)}\big[(c-a)^{1-s}+(d-b)^{1-s}-(c-b)^{1-s}-(d-a)^{1-s}\big].
\end{split}
\end{equation}
Also note that, if $n>m\geq1$, then
\begin{equation}\label{CH:1:derivative_bound}\begin{split}
(1-a^n)^{1-s}-(1-a^m)^{1-s}&=\int_m^n\frac{d}{dt}(1-a^t)^{1-s}\,dt\\
&
=(s-1)\log a\int_m^n\frac{a^t}{(1-a^t)^s}\,dt\\
&
\leq a^m (s-1)\log a\int_m^n\frac{1}{(1-a^t)^s}\,dt\\
&
\leq(n-m)a^m\frac{(s-1)\log a}{(1-a)^s}.
\end{split}
\end{equation}
Now consider the first term
\begin{equation*}
\mathcal{I}_1=\sum_{k=1}^\infty\int_{a^{2k+1}}^{a^{2k}}\int_{-1}^0\frac{dxdy}{|x-y|^{1+s}}.
\end{equation*}
Use \eqref{CH:1:rectangle_integral} and notice that $(c-a)^{1-s}-(d-a)^{1-s}\leq0$ to get
\begin{equation*}
\int_{-1}^0\int_{a^{2k+1}}^{a^{2k}}\frac{dxdy}{|x-y|^{1+s}}
\leq\frac{1}{s(1-s)}\big[(a^{2k})^{1-s}-(a^{2k+1})^{1-s}\big]\leq\frac{1}{s(1-s)}(a^{2(1-s)})^k.
\end{equation*}
Then, as $a^{2(1-s)}<1$ we get
\begin{equation*}
\mathcal{I}_1\leq\frac{1}{s(1-s)}\sum_{k=1}^\infty(a^{2(1-s)})^k<\infty.
\end{equation*}
As for the last term
\begin{equation*}
\mathcal{I}_3=\sum_{k=1}^\infty\int_{a^{2k+1}}^{a^{2k}}\int_a^1\frac{dxdy}{|x-y|^{1+s}},
\end{equation*}
use \eqref{CH:1:rectangle_integral} and notice that $(d-b)^{1-s}-(d-a)^{1-s}\leq0$ to get
\begin{equation*}\begin{split}
\int_{a^{2k+1}}^{a^{2k}}\int_a^1\frac{dxdy}{|x-y|^{1+s}}&
\leq\frac{1}{s(1-s)}\big[(1-a^{2k+1})^{1-s}-(1-a^{2k})^{1-s}\big]\\
&
\leq\frac{-\log a}{s(1-a)^s}a^{2k}\quad\textrm{by }\eqref{CH:1:derivative_bound}.
\end{split}
\end{equation*}
Thus
\begin{equation*}
\mathcal{I}_3\leq\frac{-\log a}{s(1-a)^s}\sum_{k=1}^\infty(a^2)^k<\infty.
\end{equation*}
Finally we split the second term
\begin{equation*}
\mathcal{I}_2=\sum_{k=1}^\infty\sum_{j=1}^\infty\int_{a^{2k+1}}^{a^{2k}}\int_{a^{2j}}^{a^{2j-1}}
\frac{dxdy}{|x-y|^{1+s}}
\end{equation*}
into three pieces according to the cases $j>k$, $j=k$ and $j<k$.

If $j=k$, using \eqref{CH:1:rectangle_integral} we get
\begin{equation*}\begin{split}
\int_{a^{2k+1}}^{a^{2k}}&\int_{a^{2k}}^{a^{2k-1}}
\frac{dxdy}{|x-y|^{1+s}}=\\
&
=\frac{1}{s(1-s)}\big[(a^{2k}-a^{2k+1})^{1-s}+(a^{2k-1}-a^{2k})^{1-s}-(a^{2k-1}-a^{2k+1})^{1-s}\big]\\
&
=\frac{1}{s(1-s)}\big[a^{2k(1-s)}(1-a)^{1-s}+a^{(2k-1)(1-s)}(1-a)^{1-s}\\
&
\quad\quad\quad\quad\quad-a^{(2k-1)(1-s)}(1-a^2)^{1-s}\big]\\
&
=\frac{1}{s(1-s)}(a^{2(1-s)})^k\Big[(1-a)^{1-s}+\frac{(1-a)^{1-s}}{a^{1-s}}-\frac{(1-a^2)^{1-s}}{a^{1-s}}\Big].
\end{split}
\end{equation*}
Summing over $k\in\mathbb{N}$ we get
\begin{equation*}\begin{split}
\sum_{k=1}^\infty&\int_{a^{2k+1}}^{a^{2k}}\int_{a^{2k}}^{a^{2k-1}}
\frac{dxdy}{|x-y|^{1+s}}=\\
&
=\frac{1}{s(1-s)}\frac{a^{2(1-s)}}{1-a^{2(1-s)}}\Big[(1-a)^{1-s}+\frac{(1-a)^{1-s}}{a^{1-s}}-\frac{(1-a^2)^{1-s}}{a^{1-s}}\Big]<\infty.
\end{split}
\end{equation*}
In particular note that
\begin{equation*}\begin{split}
(1-s)&\Per_s(E)\geq(1-s)\mathcal{I}_2\\
&
\geq\frac{1}{s(1-a^{2(1-s)})}\big[a^{2(1-s)}(1-a)^{1-s}+a^{1-s}(1-a)^{1-s}-a^{1-s}(1-a^2)^{1-s}\big],
\end{split}
\end{equation*}
which tends to $+\infty$ when $s\to1$. This shows that $E$ cannot have finite perimeter.

To conclude let $j>k$, the case $j<k$ being similar, and consider
\begin{equation*}
\sum_{k=1}^\infty\sum_{j=k+1}^\infty\int_{a^{2j}}^{a^{2j-1}}\int_{a^{2k+1}}^{a^{2k}}
\frac{dxdy}{|x-y|^{1+s}}.
\end{equation*}
Again, using \eqref{CH:1:rectangle_integral} and $(d-b)^{1-s}-(d-a)^{1-s}\leq0$, we get
\begin{equation*}\begin{split}
\int_{a^{2j}}^{a^{2j-1}}&\int_{a^{2k+1}}^{a^{2k}}
\frac{dxdy}{|x-y|^{1+s}}\\
&
\leq\frac{1}{s(1-s)}\big[(a^{2k+1}-a^{2j})^{1-s}-(a^{2k+1}-a^{2j-1})^{1-s}\big]\\
&
=\frac{a^{1-s}}{s(1-s)}(a^{2(1-s)})^k\big[(1-a^{2(j-k)-1})^{1-s}-(1-a^{2(j-k)-2})^{1-s}\big]\\
&
\leq\frac{a^{1-s}}{s(1-s)}(a^{2(1-s)})^k\frac{(s-1)\log a}{(1-a)^s}a^{2(j-k)-2}\quad\quad\textrm{by }\eqref{CH:1:derivative_bound}\\
&
=\frac{-\log a}{s(1-a^s)a^{s+1}}(a^{2(1-s)})^k(a^2)^{j-k},
\end{split}
\end{equation*}
for $j\geq k+2$. Then
\begin{equation*}
\begin{split}
\sum_{k=1}^\infty&\sum_{j=k+2}^\infty\int_{a^{2j}}^{a^{2j-1}}\int_{a^{2k+1}}^{a^{2k}}
\frac{dxdy}{|x-y|^{1+s}}\\
&
\leq\frac{-\log a}{s(1-a^s)a^{s+1}}\sum_{k=1}^\infty(a^{2(1-s)})^k\sum_{h=2}^\infty(a^2)^h<\infty.
\end{split}
\end{equation*}
If $j=k+1$ we get
\begin{equation*}\begin{split}
\sum_{k=1}^\infty\int_{a^{2k+2}}^{a^{2k+1}}\int_{a^{2k+1}}^{a^{2k}}\frac{dxdy}{|x-y|^{1+s}}&
\leq\frac{1}{s(1-s)}\sum_{k=1}^\infty(a^{2k+1}-a^{2k+2})^{1-s}\\
&
=\frac{a^{1-s}(1-a)^{1-s}}{s(1-s)}\sum_{k=1}^\infty(a^{2(1-s)})^k<\infty.
\end{split}
\end{equation*}
This shows that also $\mathcal{I}_2<\infty$, so that $\Per_s(E)<\infty$ for every $s\in(0,1)$ as claimed.

\end{chapter}

\begin{chapter}[Approximation of sets of finite fractional perimeter]{Approximation of sets of finite fractional perimeter by smooth sets
and comparison of local and global $s$-minimal surfaces}\label{CH_Appro_Min}

\minitoc

%
%
%
%

\section{Introduction and main results}

%

This chapter is divided in two parts.

In the first part we prove that a set has (locally) finite fractional perimeter if and only if it can be approximated (in an appropriate way) by smooth open sets.
To be more precise, we show that a set $E$ has locally finite $s$-perimeter if and only if we can find a sequence of smooth open sets
which converge in measure to $E$,
whose boundaries converge to that of $E$ in a uniform sense, and whose $s$-perimeters converge to that of $E$ in every bounded open set.

%
%


\smallskip

The second part of this chapter is concerned with sets minimizing the fractional perimeter.

We recall that, given a set $A$ and an open set $\Omega$, we will write $A\Subset\Omega$ to mean that the closure $\overline A$ of $A$ is compact and $\overline{A}\subseteq \Omega$. In particular, notice that if $A\Subset \Omega$, then $A$ must be bounded.

We consider sets which are locally $s$-minimal in an open set $\Omega\subseteq\R^n$,
namely sets which
minimize the $s$-perimeter in every open subset $\Omega'\Subset\Omega$,
and we prove existence and compactness results which extend those of \cite{CRS10}.

We also compare this definition of local $s$-minimal set with
the definition of $s$-minimal set introduced in \cite{CRS10}, proving that they coincide
when the domain $\Omega$ is a bounded open set with Lipschitz boundary (see Theorem \ref{CH:2:confront_min_teo}).

In particular, the following existence results are proven:
\begin{itemize}
\item if $\Omega$ is an open set and $E_0$ is a fixed set, then there exists a set $E$ which is locally $s$-minimal in $\Omega$
and such that $E\setminus\Omega=E_0\setminus\Omega$;
\item there exist minimizers in the class of subgraphs, namely nonlocal nonparametric minimal surfaces (see Theorem \ref{CH:2:nonparametric_exist_teo}
for a precise statement);
\item if $\Omega$ is an open set which has finite $s$-perimeter, then for every fixed set $E_0$ there exists
a set $E$ which is $s$-minimal in $\Omega$ and such that $E\setminus\Omega=E_0\setminus\Omega$.
\end{itemize}

On the other hand, we show that when the domain $\Omega$ is unbounded the nonlocal part
of the $s$-perimeter can be infinite, thus preventing the existence of
competitors having finite $s$-perimeter in $\Omega$ and hence also of ``global'' $s$-minimal sets.
In particular, we study this situation in a cylinder $\Omega^\infty:=\Omega\times\R\subseteq\R^{n+1}$,
considering as exterior data the subgraph of a (locally) bounded function.

\medskip

In the following subsections we present the precise statements of the main results of this chapter.

\subsection{Sets having (locally) finite $s$-perimeter}

We recall that we implicitly assume that all the sets we consider
contain their measure theoretic interior, do not intersect their measure theoretic exterior, and are such that their topological boundary coincides with their measure theoretic boundary---see Remark \ref{CH:1:gmt_assumption} and Appendix \ref{CH:1:Appendix_meas_th_bdary} 
for the details.

\smallskip


We recall that we say that a set $E\subseteq\R^n$ has locally finite $s$-perimeter in an open set $\Omega\subseteq\R^n$ if
\begin{equation*}
\Per_s(E,\Omega')<\infty\qquad\textrm{for every open set }\Omega'\Subset\Omega.
\end{equation*}

We remark that the family of sets having finite $s$-perimeter in $\Omega$ need not coincide with the family of sets of locally
finite $s$-perimeter in $\Omega$,
not even when $\Omega$
is ``nice'' (say bounded and with Lipschitz boundary). To be more precise, since
\begin{equation}\label{CH:2:sup_for_per_eq}
\Per_s(E,\Omega)=\sup_{\Omega'\Subset\Omega}\Per_s(E,\Omega'),
\end{equation}
(see Proposition \ref{CH:2:continuity_in_open_set_seq} and Remark \ref{CH:2:continuity_in_open_set_seq_rmk}),
a set which has finite $s$-perimeter in $\Omega$ has also locally finite $s$-perimeter.
However the converse, in general, is false.\\
When $\Omega$ is not bounded it is clear that also for sets of locally finite $s$-perimeter
the sup in \eqref{CH:2:sup_for_per_eq}
may be infinite (consider, e.g., $\Omega=\R^n$ and $E=\{x_n\leq0\}$).

Actually, as shown in Remark \ref{CH:2:counterex_loc_fin_per}, this may happen even when $\Omega$ is bounded and has Lipschitz boundary.
Roughly speaking, this is because the set $E$ might oscillate more and more as it approaches the boundary $\partial\Omega$.


\subsection{Approximation by smooth open sets}
We denote by $N_\varrho(\Gamma)$ the $\varrho$-neighborhood of a set $\Gamma\subseteq\R^n$, that is
\begin{equation*}
N_\varrho(\Gamma):=\{x\in\R^n\,|\,d(x,\Gamma)<\varrho\}.
\end{equation*}

The main approximation result is the following. In particular it shows that open sets with smooth boundary are dense in the family of sets of locally finite $s$-perimeter.

\begin{theorem}\label{CH:2:density_smooth_teo}
Let $\Omega\subseteq\R^n$ be an open set. A set $E\subseteq\R^n$ has locally finite $s$-perimeter in $\Omega$
if and only if there exists a sequence $E_h\subseteq\R^n$ of open sets with smooth boundary and $\eps_h\longrightarrow0^+$ such that
\begin{equation*}\begin{split}
& (i)\quad E_h\xrightarrow{loc}E,\qquad\sup_{h\in\mathbb N}\Per_s(E_h,\Omega')<\infty\quad\textrm{for every }\Omega'\Subset\Omega,\\
& (ii)\quad\lim_{h\to\infty}\Per_s(E_h,\Omega')=\Per_s(E,\Omega')\quad\textrm{for every }\Omega'\Subset\Omega,\\
& (iii)\quad\partial E_h\subseteq N_{\eps_h}(\partial E).
\end{split}
\end{equation*}
Moreover, if $\Omega=\R^n$ and the set $E$ is such that $|E|<\infty$ and $\Per_s(E)<\infty$, then
\begin{equation}\label{CH:2:last_eq_teo_app}
E_h\longrightarrow E,\qquad\quad\qquad \lim_{h\to\infty}\Per_s(E_h)=\Per_s(E),
\end{equation}
and we can require each set $E_h$ to be bounded (instead of asking $(iii)$).
\end{theorem}

We recall that, as we have observed in Section \ref{CH:0:Appro_subsection}, such a result is well known for Caccioppoli sets (see, e.g.,  \cite{Maggi}) and indeed this density property can be used to define the (classical) perimeter functional as the relaxation (with respect to $L^1_{\loc}$ convergence) of the $\Ha^{n-1}$ measure of boundaries of smooth open sets, that is
\begin{equation}\label{CH:2:liminfclassical}\begin{split}\Per(E,\Omega)=\inf\Big\{\liminf_{k\to\infty}\Ha^{n-1}(\partial E_h&\cap\Omega)\,\big|\,
E_h\subseteq\R^n\textrm{ open with smooth}\\
&
\textrm{boundary, s.t. }E_h\xrightarrow{loc}E\Big\}.
\end{split}
\end{equation}

The scheme of the proof of Theorem \ref{CH:2:density_smooth_teo} is the following.

First of all, in Section \ref{CH:2:Approx_Section} we prove appropriate approximation results for the functional
\begin{equation*}
\Fc(u,\Omega)=\frac{1}{2}\iint_{\R^{2n}\setminus(\Co\Omega)^2}\frac{|u(x)-u(y)|}{|x-y|^{n+s}}\,dx\,dy,
\end{equation*}
which we believe might be interesting on their own.

Then we exploit the generalized coarea formula
\begin{equation*}
\Fc(u,\Omega)=\int_{-\infty}^\infty \Per_s(\{u>t\},\Omega)\,dt,
\end{equation*}
and Sard's Theorem to obtain the approximation of the set $E$ by superlevel sets of smooth functions which approximate $\chi_E$.

Finally, a diagonal argument guarantees the convergence of the $s$-perimeters
in every open set $\Omega'\Subset\Omega$.

\begin{remark}
Let $\Omega\subseteq\R^n$ be a bounded open set with Lipschitz boundary and consider a set $E$ which has
finite $s$-perimeter in $\Omega$. Notice that if we apply Theorem \ref{CH:2:density_smooth_teo},
in point $(ii)$ we do not get the convergence of the $s$-perimeters in $\Omega$, but only
in every $\Omega'\Subset\Omega$.
On the other hand, if we can find an open set $\mathcal O$ such that $\Omega\Subset\mathcal O$
and
\[\Per_s(E,\mathcal O)<\infty,\]
then we can apply Theorem \ref{CH:2:density_smooth_teo} in $\mathcal O$. In particular, since $\Omega\Subset\mathcal O$,
by point $(ii)$ we obtain
\begin{equation}\label{CH:2:conv_whole_Omega}
\lim_{h\to\infty}\Per_s(E_h,\Omega)=\Per_s(E,\Omega).
\end{equation}
\end{remark}

Still, when $\Omega$ is a bounded open set with Lipschitz boundary, we can always obtain the convergence \eqref{CH:2:conv_whole_Omega}
at the cost of weakening a little our request on the uniform convergence of the boundaries.

\begin{theorem}\label{CH:2:appro_in_bded_open}
Let $\Omega\subseteq\R^n$ be a bounded open set with Lipschitz boundary. A set $E\subseteq\R^n$ has finite $s$-perimeter
in $\Omega$ if and only if there exists a sequence $\{E_h\}$ of open sets with smooth boundary
and $\eps_h\longrightarrow0^+$ such that
\begin{equation*}\begin{split}
& (i)\quad E_h\xrightarrow{loc}E,\qquad\sup_{h\in\mathbb N}\Per_s(E_h,\Omega)<\infty,\\
& (ii)\quad\lim_{h\to\infty}\Per_s(E_h,\Omega)=\Per_s(E,\Omega),\\
& (iii)\quad\partial E_h\setminus N_{\eps_h}(\partial\Omega)\subseteq N_{\eps_h}(\partial E).
\end{split}
\end{equation*}
\end{theorem}

Notice that in point $(iii)$ we do not ask the convergence of the boundaries in the whole of $\R^n$
but only in $\R^n\setminus N_\delta(\partial\Omega)$ (for any fixed $\delta>0$).
Since
$N_{\eps_h}(\partial\Omega)\searrow\partial\Omega$,
roughly speaking, the convergence holds in $\R^n$ ``in the limit''.

Moreover, we remark that point $(ii)$ in Theorem \ref{CH:2:appro_in_bded_open}
guarantees the convergence of the $s$-perimeters also in every $\Omega'\Subset\Omega$ (see Remark \ref{CH:2:rmk_conv_every_cpt_subopen}).

Finally, from the lower semicontinuity of the $s$-perimeter and Theorem \ref{CH:2:appro_in_bded_open}, we obtain
\begin{corollary}
Let $\Omega\subseteq\R^n$ be a bounded open set with Lipschitz boundary and let $E\subseteq\R^n$. Then
\begin{equation*}\begin{split}
\Per_s(E,\Omega)=\inf\Big\{\liminf_{h\to\infty}\Per_s(E_h&,\Omega)\,\big|\,E_h\subseteq\R^n\textrm{ open with smooth}\\
&
\textrm{boundary, s.t. }E_h\xrightarrow{loc}E\Big\}.
\end{split}
\end{equation*}

\end{corollary}

For similar approximation results see also \cite{CSV16} and \cite{CDV17}.

It is interesting to observe that in \cite{DV18} the authors have proved, by exploiting the divergence Theorem, that
if $E\subseteq\R^n$ is a bounded open set with smooth boundary, then
\begin{equation}\label{CH:2:liminffract}
\Per_s(E)=c_{n,s}\int_{\partial E}\int_{\partial E}
\frac{2-|\nu_E(x)-\nu_E(y)|^2}{|x-y|^{n+s-2}}d\Ha^{n-1}_xd\Ha^{n-1}_y,
\end{equation}
where $\nu_E$ denotes the external normal of $E$ and
\begin{equation*}
c_{n,s}:=\frac{1}{2s(n+s-2)}.
\end{equation*}

Notice that in order to consider the right hand side of \eqref{CH:2:liminffract},
we need the boundary of the set $E$ to be at least locally $(n-1)$-rectifiable, so that the Hausdorff dimension of $\partial E$
is $n-1$ and $E$ has a well defined normal vector at $\Ha^{n-1}$-a.e. $x\in\partial E$.
Therefore, the equality \eqref{CH:2:liminffract} cannot hold true for a generic set $E$ having finite $s$-perimeter,
since, as remarked in the beginning of the Introduction, such a set could have a nowhere rectifiable boundary.

Nevertheless, as a consequence of the equality \eqref{CH:2:liminffract}, of the lower semicontinuity of the $s$-perimeter
and of Theorem \ref{CH:2:density_smooth_teo}, we obtain the following Corollary, which can be thought of
as an analogue of \eqref{CH:2:liminfclassical} in the fractional setting.

\begin{corollary}
Let $E\subseteq\R^n$ be such that $|E|<\infty$. Then
\begin{equation*}\begin{split}
\Per_s(E)=\inf\Big\{&\liminf_{h\to\infty}c_{n,s}\int_{\partial E_h}\int_{\partial E_h}
\frac{2-|\nu_{E_h}(x)-\nu_{E_h}(y)|^2}{|x-y|^{n+s-2}}d\Ha^{n-1}_xd\Ha^{n-1}_y\,\big|\\
&\quad
E_h\subseteq\R^n\textrm{ bounded open set with smooth
boundary, s.t. }E_h\xrightarrow{loc}E\Big\}.
\end{split}
\end{equation*}
\end{corollary}


\subsection{Nonlocal minimal surfaces}

First of all we recall the definition of (locally) $s$-minimal sets.

\begin{defn}
Let $\Omega\subseteq\R^n$ be an open set and let $s\in(0,1)$. We say that a set $E\subseteq\R^n$
is \emph{$s$-minimal} in $\Omega$ if $\Per_s(E,\Omega)<\infty$ and
\begin{equation*}
F\setminus\Omega=E\setminus\Omega\quad\Longrightarrow\quad \Per_s(E,\Omega)\leq \Per_s(F,\Omega).
\end{equation*}
We say that a set $E\subseteq\R^n$ is \emph{locally $s$-minimal} in $\Omega$ if it is $s$-minimal in every open subset
$\Omega'\Subset\Omega$.
\end{defn}

When the open set $\Omega\subseteq\R^n$ is bounded and has Lipschitz boundary, the notions of $s$-minimal set and locally $s$-minimal set coincide.

\begin{theorem}\label{CH:2:confront_min_teo}
Let $\Omega\subseteq\R^n$ be a bounded open set with Lipschitz boundary and let $E\subseteq\R^n$. The following are equivalent:
\begin{itemize}
\item[(i)] $E$ is $s$-minimal in $\Omega$;

\item[(ii)] $\Per_s(E,\Omega)<\infty$ and
\begin{equation*}
\Per_s(E,\Omega)\leq \Per_s(F,\Omega)\qquad
\textrm{for every }F\subseteq\R^n\quad\textrm{s.t.}\quad E\Delta F\Subset\Omega;
\end{equation*}

\item[(iii)] $E$ is locally $s$-minimal in $\Omega$.
\end{itemize}
\end{theorem}

 We remark that a set as in $(ii)$ is called a local minimizer for $\Per_s(-,\Omega)$ in \cite{Gamma} and a ``nonlocal area minimizing surface'' in $\Omega$ in \cite{DdPW18}.

\begin{remark}
The implications $(i)\Longrightarrow(ii)\Longrightarrow(iii)$ actually hold in any open set $\Omega\subseteq\R^n$.
\end{remark}

In \cite{CRS10} the authors proved that if $\Omega$ is a bounded open set with Lipschitz boundary,
then given any fixed set $E_0\subseteq\R^n$ we can find a set $E$ which is $s$-minimal in $\Omega$
and such that $E\setminus\Omega=E_0\setminus\Omega$.

This is because
\begin{equation*}
\Per_s(E_0\setminus\Omega,\Omega)\leq \Per_s(\Omega)<\infty,
\end{equation*}
so the exterior datum $E_0\setminus\Omega$ is itself an admissible competitor with finite $s$-perimeter in $\Omega$ and we can use
the direct method of the Calculus of Variations to obtain a minimizer.

In Section \ref{CH:2:Compactness_Section} we prove a compactness property which we use in Section \ref{CH:2:Existence_Section} to
prove the following existence results, which extend that of \cite{CRS10}.

\begin{theorem}\label{CH:2:glob_min_exist}
Let $\Omega\subseteq\R^n$ be an open set and let $E_0\subseteq\R^n$. Then there exists
a set $E\subseteq\R^n$ $s$-minimal in $\Omega$, with $E\setminus\Omega=E_0\setminus\Omega$, if and only if
there exists a set $F\subseteq\R^n$, with $F\setminus\Omega=E_0\setminus\Omega$ and such that
$
\Per_s(F,\Omega)<\infty.
$
\end{theorem}

An immediate consequence of this Theorem is the existence of $s$-minimal sets in open sets having finite $s$-perimeter.
\begin{corollary}
Let $s\in(0,1)$ and let $\Omega\subseteq\R^n$ be an open set such that
\[\Per_s(\Omega)<\infty.\]
Then for every $E_0\subseteq\R^n$
there exists a set $E\subseteq\R^n$ $s$-minimal in $\Omega$, with $E\setminus\Omega=E_0\setminus\Omega$.
\end{corollary}

Even if we cannot find a competitor with finite $s$-perimeter, we can always find a locally $s$-minimal set.

\begin{corollary}\label{CH:2:loc_min_set_cor}
Let $\Omega\subseteq\R^n$ be an open set and let $E_0\subseteq\R^n$.
Then there exists a set $E\subseteq\R^n$ locally $s$-minimal in $\Omega$, with $E\setminus\Omega=E_0\setminus\Omega$.
\end{corollary}

In Section \ref{CH:2:Min_Compactness_Section} we also prove compactness results for (locally) $s$-minimal sets
(by slightly modifying the proof of \cite[Theorem 3.3]{CRS10}, which proved compactness for $s$-minimal sets in a ball).
Namely, we prove that every limit set of a sequence of (locally) $s$-minimal sets is itself (locally) $s$-minimal.

\begin{theorem}\label{CH:2:minimal_comp}
Let $\Omega\subseteq\R^n$ be a bounded open set with Lipschitz boundary. Let $\{E_k\}$ be a sequence of $s$-minimal sets in $\Omega$, with $E_k\xrightarrow{loc}E$. Then $E$ is $s$-minimal in $\Omega$ and
\begin{equation}\label{CH:2:conv_perimeter}
\Per_s(E,\Omega)=\lim_{k\to\infty}\Per_s(E_k,\Omega).
\end{equation}

\end{theorem}

\begin{corollary}\label{CH:2:local_minima_comp}
Let $\Omega\subseteq\R^n$ be an open set. Let $\{E_h\}$ be a sequence of sets locally $s$-minimal
in $\Omega$, with $E_h\xrightarrow{loc}E$. Then $E$ is locally $s$-minimal in $\Omega$ and
\begin{equation}\label{CH:2:conv_perimeter_locally}
\Per_s(E,\Omega')=\lim_{h\to\infty}\Per_s(E_h,\Omega'),\qquad\textrm{for every }\Omega'\Subset\Omega.
\end{equation}

\end{corollary}

\subsubsection{Minimal sets in cylinders}

We have seen in Corollary \ref{CH:2:loc_min_set_cor} that a locally $s$-minimal set always exists, no matter what the domain $\Omega$
or the exterior data $E_0\setminus\Omega$ are.

On the other hand, by Theorem \ref{CH:2:glob_min_exist} we know that the only requirement needed for the existence of an $s$-minimal set
is the existence of a competitor with finite $s$-perimeter.\\
We show that even in the case of a regular domain, like the cylinder $\Omega^\infty:=\Omega\times\R$,
with $\Omega\subseteq\R^n$ bounded with $C^{1,1}$ boundary,
such a competitor might not exist. Roughly speaking,
this is a consequence of the unboundedness of the domain $\Omega^\infty$, which forces the nonlocal part of the $s$-perimeter
to be infinite.

In Section \ref{CH:2:Min_Sets_Cyl_Section} we study (locally) $s$-minimal sets in $\Omega^\infty$, with respect to the exterior data
given by the subgraph of a function $v$, that is
\[\Sg(v):=\left\{(x,t)\in\R^{n+1}\,|\,t<v(x)\right\}.\]

In particular, we consider sets which are $s$-minimal in the ``truncated'' cylinders
\[\Omega^k:=\Omega\times(-k,k),\]
showing that if the function $v$ is locally bounded, then these $s$-minimal sets cannot ``oscillate'' too much.
Namely their boundaries are constrained in a cylinder $\Omega\times(-M,M)$
independently on $k$.
As a consequence, we can find $k_0$ big enough such that a set $E$ is locally $s$-minimal in $\Omega^\infty$
if and only if it is $s$-minimal in $\Omega^{k_0}$ (see Lemma \ref{CH:2:bded_cyl_prop} and Proposition \ref{CH:2:bded_cyl_coroll} for the precise statements).

However, in general a set $s$-minimal in $\Omega^\infty$ does not exist. As an example we prove that
there cannot exist an $s$-minimal set having as exterior data the subgraph of a bounded function.

\smallskip

Frst of all, we recall that we can write the fractional perimeter as the sum
\begin{equation*}
\Per_s(E,\Omega)=\Per_s^L(E,\Omega)+\Per_s^{NL}(E,\Omega),
\end{equation*}
where
\begin{equation*}\begin{split}
&\Per_s^L(E,\Omega):=\mathcal L_s(E\cap\Omega,\Co E\cap\Omega)=\frac{1}{2}[\chi_E]_{W^{s,1}(\Omega)},\\
&
\Per_s^{NL}(E,\Omega):=\Ll_s(E\cap\Omega,\Co E\setminus\Omega)+\Ll_s(E\setminus\Omega,\Co E\cap\Omega).
\end{split}\end{equation*}
We can think of $\Per^L_s(E,\Omega)$ as the local part of the fractional perimeter, in the sense that if $|(E\Delta F)\cap\Omega|=0$,
then $\Per^L_s(F,\Omega)=\Per^L_s(E,\Omega)$.

\smallskip

The main result of Section \ref{CH:2:Min_Sets_Cyl_Section} is the following:

\begin{theorem}\label{CH:2:bound_unbound_per_cyl_prop}
Let $\Omega\subseteq\R^n$ be a bounded open set. Let $E\subseteq\R^{n+1}$ be such that
\begin{equation}\label{CH:2:bound_hp_forml_subgraph}
\Omega\times(-\infty,-k]\subseteq E\cap\Omega^\infty\subseteq\Omega\times(-\infty,k],
\end{equation}
for some $k\in\mathbb N$, and suppose that $\Per_s(E,\Omega^{k+1})<\infty$. Then
\begin{equation*}
\Per_s^L(E,\Omega^\infty)<\infty.
\end{equation*}
On the other hand, if
\begin{equation}\label{CH:2:bound_hp_forml_subgraph2}
\{x_{n+1}\leq-k\}\subseteq E\subseteq\{x_{n+1}\leq k\},
\end{equation}
then
\begin{equation*}
\Per^{NL}_s(E,\Omega^\infty)=\infty.
\end{equation*}
In particular, if $\Omega$ has $C^{1,1}$ boundary and $v\in L^\infty(\R^n)$, there cannot exist an $s$-minimal set in $\Omega^\infty$ with exterior data
\[
\Sg(v)\setminus\Omega^\infty=\{(x,t)\in\R^{n+1}\,|\,x\in\Co\Omega,\quad t<v(x)\}.\]
\end{theorem}

\begin{remark}
From Theorem \ref{CH:2:glob_min_exist} we see that if $v\in L^\infty(\R^n)$,
there cannot exist a set $E\subseteq\R^{n+1}$ such that $E\setminus\Omega^\infty=\Sg(v)\setminus\Omega^\infty$
and $\Per_s(E,\Omega^\infty)<\infty$.
\end{remark}

As a consequence of the computations developed in the proof of Theorem \ref{CH:2:bound_unbound_per_cyl_prop},
in the end of Section \ref{CH:2:Min_Sets_Cyl_Section} we also show that we cannot define a ``naive'' fractional nonlocal version
of the area functional as
\begin{equation*}
\mathscr A_s(u,\Omega):=\Per_s(\Sg(u),\Omega^\infty),
\end{equation*}
since this would be infinite even for very regular functions.

\smallskip

To conclude, we remark that as an immediate consequence of Corollary \ref{CH:2:loc_min_set_cor}
and \cite[Theorem 1.1]{graph}, we obtain an existence result for the Plateau's problem in the class of subgraphs.

\begin{theorem}\label{CH:2:nonparametric_exist_teo}
Let $\Omega\subseteq\R^n$ be a bounded open set with $C^{1,1}$ boundary.
For every function $v\in C(\R^n)$ there exists a function
$u\in C(\overline{\Omega})$ such that, if
\begin{equation*}
\tilde{u}:=\chi_\Omega u+(1-\chi_\Omega)v,
\end{equation*}
then $\Sg(\tilde{u})$ is locally $s$-minimal in $\Omega^\infty$.
\end{theorem}

Notice that, as remarked in \cite{graph}, the function $\tilde{u}$ need not be continuous. Indeed, because of boundary stickiness effects
of $s$-minimal surfaces (see, e.g., \cite{boundary}), in general we might have
\[u_{|_{\partial\Omega}}\not=v_{|_{\partial\Omega}}.\]

\section{Tools}

We collect here some auxiliary results that we will exploit in the following sections.

We begin by pointing out the following easy but useful result.
\begin{prop}\label{CH:2:subopensets}
Let $\Omega'\subseteq\Omega\subseteq\R^n$ be open sets and let $E\subseteq\R^n$. Then
\begin{equation*}\begin{split}
\Per_s(E,\Omega)=\Per_s(E,\Omega')&+\Ll_s\big(E\cap(\Omega\setminus\Omega'),\Co E\setminus\Omega\big)
+\Ll_s\big(E\setminus\Omega,\Co E\cap(\Omega\setminus\Omega')\big)\\
&
\qquad
+\Ll_s\big(E\cap(\Omega\setminus\Omega'),\Co E\cap(\Omega\setminus\Omega')\big).
\end{split}
\end{equation*}

As a consequence,

$(i)\quad$ if $E\subseteq\Omega$, then
\begin{equation*}
\Per_s(E,\Omega)=\Per_s(E),
\end{equation*}

$(ii)\quad$ if $E,\,F\subseteq\R^n$ have finite $s$-perimeter in $\Omega$
and $E\Delta F\subseteq\Omega'\subseteq\Omega$, then
\begin{equation*}
\Per_s(E,\Omega)-\Per_s(F,\Omega)=\Per_s(E,\Omega')-\Per_s(F,\Omega').
\end{equation*}

\end{prop}

\begin{remark}\label{CH:2:bded_set_frac_per}
In particular, if $E$ has finite $s$-perimeter in $\Omega$, then it has finite $s$-perimeter also in every open set $\Omega'\subseteq\Omega$.

\end{remark}

\subsection{Bounded open sets with Lipschitz boundary}

It is convenient to recall here some notation and results concerning the signed distance function, since we will make extensive use of such results in the subsequent sections.

Given a set $E\subseteq\R^n$, with $E\not=\emptyset$, the distance function from $E$ is defined as
\begin{equation*}
d_E(x)=d(x,E):=\inf_{y\in E}|x-y|,\qquad\textrm{for }x\in\R^n.
\end{equation*}
The signed distance function from $\partial E$, negative inside $E$, is then defined as
\begin{equation*}
\bar{d}_E(x)=\bar{d}(x,E):=d(x,E)-d(x,\Co E).
\end{equation*}

We also define for every $r\in\R$ the sets
\begin{equation*}
E_r:=\{x\in\R^n\,|\,\bar{d}_E(x)<r\}.
\end{equation*}
Notice that if $\varrho>0$, then
\begin{equation*}
N_\varrho(\partial\Omega)=\{|\bar{d}_\Omega|<\varrho\}=\Omega_\varrho\setminus\overline{\Omega_{-\varrho}}
\end{equation*}
is the $\varrho$-tubular neighborhood of $\partial\Omega$.

Let $\Omega\subseteq\R^n$ be a bounded open set with Lipschitz boundary.
It is well known (see, e.g.,
\cite[Theorem 4.1]{Doktor}) that also the bounded open sets $\Omega_r$ have
Lipschitz boundary, when $r$ is small enough, say $|r|<r_0$.
Also notice that
\begin{equation*}
\partial\Omega_r=\{\bar{d}_\Omega=r\}.
\end{equation*}

Moreover the perimeter of $\Omega_r$
can be bounded uniformly in $r\in(-r_0,r_0)$ (see also Appendix \ref{CH:1:Appendix_distance_function} for a more detailed discussion)

\begin{prop}\label{CH:2:bound_perimeter_unif}
Let $\Omega\subseteq\R^n$ be a bounded open set with Lipschitz boundary. Then there exists $r_0>0$ such that
$\Omega_r$ is a bounded open set with Lipschitz boundary for every $r\in(-r_0,r_0)$ and
\begin{equation}\label{CH:2:bound_perimeter_unif_eq}
\sup_{|r|<r_0}\Ha^{n-1}(\{\bar{d}_\Omega=r\})<\infty.
\end{equation}
\end{prop}

As a consequence, exploiting the embedding $BV(\R^n)\hookrightarrow W^{s,1}(\R^n)$ we obtain
a uniform bound for the (global) $s$-perimeters
of the sets $\Omega_r$ (see Corollary \ref{CH:1:embedding_fin_per_coroll}).
\begin{corollary}
Let $\Omega\subseteq\R^n$ be a bounded open set with Lipschitz boundary. Then there exists $r_0>0$
such that
\begin{equation}\label{CH:2:unif_bound_lip_frac_per}
\sup_{|r|<r_0}\Per_s(\Omega_r)<\infty.
\end{equation}
\end{corollary}

\subsubsection{Increasing sequences}

In particular, Proposition \ref{CH:2:bound_perimeter_unif} shows that if $\Omega$ is a bounded open set with Lpschitz boundary,
then we can approximate it strictly from the inside with a sequence of bounded open sets
$\Omega_k:=\Omega_{-1/k}\Subset\Omega$. Moreover,
\eqref{CH:2:bound_perimeter_unif_eq} gives a uniform bound on the measure of the boundaries
of the approximating sets.

Now we prove that any open set $\Omega\not=\emptyset$ can be approximated strictly from the inside with
a sequence of bounded open sets with smooth boundaries.

\begin{prop}\label{CH:2:first_approx_prop}
Let $\Omega\subseteq\R^n$ be a bounded open set. For every $\eps>0$ there exists a bounded open set $\mathcal O_\eps\subseteq\R^n$ with smooth boundary, such that
\begin{equation*}
\mathcal O_\eps\Subset\Omega\qquad\textrm{and}\qquad\partial\mathcal O_\eps\subseteq N_\eps(\partial\Omega).
\end{equation*}
\end{prop}

\begin{proof}

We show that we can approximate the set $\Omega_{-\eps/2}$
with a bounded open set $\mathcal O_\eps$ with smooth boundary such that $\partial\mathcal O_\eps\subseteq N_{\eps/4}(\partial\Omega_{-\eps/2})$.\\
In general $\mathcal O_\eps\not\subseteq\Omega_{-\eps/2}$.
However
\begin{equation}\label{CH:2:eq_app_op}
\mathcal O_\eps\subseteq N_{\eps/4}(\Omega_{-\eps/2})\Subset\Omega
\quad\textrm{and indeed}\quad\Omega_{-3\eps/4}\subseteq\mathcal O_\eps\subseteq\Omega_{-\eps/4},
\end{equation}
proving the claim.

Let $u:=\chi_{\Omega_{-\eps/2}}$ and consider the regularized function
\begin{equation*}
v:=u_{\eps/4}=u\ast\eta_{\eps/4}
\end{equation*}
 (see Section \ref{CH:2:Approx_Section} for the details about the mollifier $\eta_\eps$).
Since $v\in C^\infty(\R^n)$, we know from Sard's Theorem that the superlevel set $\{v>t\}$ is an open set with smooth boundary for a.e. $t\in(0,1)$.
Moreover notice that $0\leq v\leq1$, with
\begin{equation*}
\textrm{supp }v\subseteq N_{\eps/4}(\textrm{supp }u)=N_{\eps/4}(\Omega_{-\eps/2})\subseteq\Omega_{-\eps/4},
\end{equation*}
and
\begin{equation*}
v(x)=1\qquad\textrm{for every }x\in\Big\{y\in\Omega_{-\eps/2}\,\big|\,d(y,\partial\Omega_{-\eps/2})>\frac{\eps}{4}\Big\}
\supseteq\Omega_{-\frac{3}{4}\eps}.
\end{equation*}
This shows that $\mathcal O_\eps:=\{v>t\}$ (for any ``regular'' $t$) satisfies
\eqref{CH:2:eq_app_op}.
\end{proof}

\begin{corollary}\label{CH:2:regular_approx_open_sets_coroll}
Let $\Omega\subseteq\R^n$ be an open set. Then there exists a sequence $\{\Omega_k\}$ of bounded open sets with smooth boundary such that $\Omega_k\nearrow\Omega$ strictly, i.e.
\begin{equation*}
\Omega_k\Subset\Omega_{k+1}\Subset\Omega\qquad\textrm{and}\qquad\bigcup_{k\in\mathbb N}\Omega_k=\Omega.
\end{equation*}
In particular $\Omega_k\xrightarrow{loc}\Omega$.
\end{corollary}

\begin{proof}

It is enough to notice that we can approximate $\Omega$ strictly from the inside with bounded open sets $\mathcal O_k\subseteq\R^n$, that is
\begin{equation*}
\mathcal O_k\Subset\mathcal O_{k+1}\Subset\Omega\qquad\textrm{and}\qquad\bigcup_{k\in\mathbb N}\mathcal O_k=\Omega.
\end{equation*}

Then we can exploit Proposition \ref{CH:2:first_approx_prop}, and in particular \eqref{CH:2:eq_app_op}, to find bounded open sets $\Omega_k\subseteq\R^n$ with smooth boundary such that
\begin{equation*}
\mathcal O_k\Subset\Omega_k\Subset\mathcal O_{k+1}.
\end{equation*}
Indeed we can take as $\Omega_k$ a set $\mathcal O_\eps$ corresponding to $\mathcal O_{k+1}$, with $\eps$
small enough to guarantee $\mathcal O_k\Subset\mathcal O_\eps$.\\
As for the sets $\mathcal O_k$, if $\Omega$ is bounded we can simply take
$
\mathcal O_k:=\Omega_{-2^{-k}}.
$
If $\Omega$ is not bounded, we can consider the sets $\Omega\cap B_{2^k}$ and define
\begin{equation*}
\mathcal O_k:=\big\{x\in\Omega\cap B_{2^k}\,|\,d\big(x,\partial(\Omega\cap B_{2^k})\big)>2^{-k}\big\}.
\end{equation*}
To conclude, notice that we have $\chi_{\Omega_k}\longrightarrow\chi_\Omega$ pointwise everywhere in $\R^n$, which implies the convergence in $L^1_{\loc}(\R^n)$.
\end{proof}

\subsubsection{Some uniform estimates for $\varrho$-neighborhoods}

The uniform bound \eqref{CH:2:bound_perimeter_unif_eq} on the perimeters of the sets $\Omega_\delta$
allows us to obtain the following estimates, which will be used in the sequel.

\begin{lemma}\label{CH:2:UniFNeiGHboEstiLeMmA}
Let $\Omega\subseteq\R^n$ be a bounded open set with Lipschitz boundary. Let $\delta\in(0,r_0)$. Then
\begin{equation}\label{CH:2:uniform_bound_strips}\begin{split}
&(i)\quad\Ll_s(\Omega_{-\delta},\Omega\setminus\Omega_{-\delta})\leq C\,\delta^{1-s},\\
&
(ii)\quad\Ll_s(\Omega,\Omega_\delta\setminus\Omega)\leq C\,\delta^{1-s}
\quad\textrm{and}
\quad\Ll_s(\Omega\setminus\Omega_{-\delta},\Co\Omega)\leq C\,\delta^{1-s},
\end{split}\end{equation}
where the constant $C$ is
\begin{equation*}
C:=\frac{n\omega_n}{s(1-s)}\,\sup_{|r|<r_0}\Ha^{n-1}(\{\bar{d}_\Omega=r\}).
\end{equation*}
\end{lemma}

\begin{proof}

By using the coarea formula for $\bar{d}_\Omega$ and exploiting \eqref{CH:2:bound_perimeter_unif_eq}, we get
\begin{equation*}\begin{split}
\Ll_s(\Omega_{-\delta},\Omega\setminus\Omega_{-\delta})&
=\int_{-\delta}^0\Big(\int_{\{\bar{d}_\Omega=\varrho\}}\Big(\int_{\Omega_{-\delta}}\frac{dx}{|x-y|^{n+s}}\Big)d\Ha^{n-1}_y\Big)d\varrho\\
&
\leq\int_{-\delta}^0\Big(\int_{\{\bar{d}_\Omega=\varrho\}}
\Big(\int_{\Co B_{\varrho+\delta}(y)}\frac{dx}{|x-y|^{n+s}}\Big)d\Ha^{n-1}_y\Big)d\varrho\\
&
=\frac{n\omega_n}{s}\int_{-\delta}^0\frac{\Ha^{n-1}(\{\bar{d}_\Omega=\varrho\})}{(\varrho+\delta)^s}\,d\varrho\\
&
\leq M \frac{n\omega_n}{s(1-s)}\int_{-\delta}^0\frac{d}{d\varrho}(\varrho+\delta)^{1-s}\,d\varrho=M \frac{n\omega_n}{s(1-s)}\,\delta^{1-s}.
\end{split}\end{equation*}

In the same way we obtain point $(ii)$,
\begin{equation*}\begin{split}
\Ll_s(\Omega_\delta\setminus\Omega,\Omega)&
=\int^\delta_0\Big(\int_{\{\bar{d}_\Omega=\varrho\}}\Big(\int_\Omega\frac{dx}{|x-y|^{n+s}}\Big)d\Ha^{n-1}_y\Big)d\varrho\\
&
\leq\int^\delta_0\Big(\int_{\{\bar{d}_\Omega=\varrho\}}\Big(\int_{\Co B_\varrho(y)}\frac{dx}{|x-y|^{n+s}}\Big)d\Ha^{n-1}_y\Big)d\varrho\\
&
=\frac{n\omega_n}{s}\int^\delta_0\frac{\Ha^{n-1}(\{\bar{d}_\Omega=\varrho\})}{\varrho^s}\,d\varrho\\
&
\leq M \frac{n\omega_n}{s(1-s)}\int^\delta_0\frac{d}{d\varrho}\varrho^{1-s}\,d\varrho=M \frac{n\omega_n}{s(1-s)}\,\delta^{1-s},
\end{split}\end{equation*}
(the other estimate in point $(ii)$ is analogous).
\end{proof}

\subsection{(Semi)continuity of the $s$-perimeter}

As shown in \cite[Theorem 3.1]{CRS10}, Fatou's Lemma gives the lower semicontinuity of the functional $\Ll_s$.
\begin{prop}\label{CH:2:semicont_first_prop}
Suppose
\begin{equation*}
A_k\xrightarrow{loc}A\qquad\textrm{and}\qquad B_k\xrightarrow{loc}B.
\end{equation*}
Then
\begin{equation}\label{CH:2:semicontinuity_first}
\Ll_s(A,B)\leq\liminf_{k\to\infty}\Ll_s(A_k,B_k).
\end{equation}
In particular, if
\begin{equation*}
E_k\xrightarrow{loc}E\qquad\textrm{and}\qquad \Omega_k\xrightarrow{loc}\Omega,
\end{equation*}
then
\begin{equation*}
\Per_s(E,\Omega)\leq\liminf_{k\to\infty}\Per_s(E_k,\Omega_k).
\end{equation*}
\end{prop}

\begin{proof}
If the right hand side of \eqref{CH:2:semicontinuity_first} is infinite, we have nothing to prove, so we can suppose that it is finite.
By definition of the liminf, we can find $k_i\nearrow\infty$ such that
\begin{equation*}
\lim_{i\to\infty}\Ll_s(A_{k_i},B_{k_i})=\liminf_{k\to\infty}\Ll_s(A_k,B_k)=:I.
\end{equation*}
Since $\chi_{A_{k_i}}\to\chi_A$ and $\chi_{B_{k_i}}\to\chi_B$ in
$L^1_{\loc}(\R^n)$, up to passing to a subsequence we can suppose that
\begin{equation*}
\chi_{A_{k_i}}\longrightarrow\chi_A\qquad\textrm{and}\qquad\chi_{B_{k_i}}\longrightarrow\chi_B\qquad\textrm{a.e. in }\R^n.
\end{equation*}
Then, since
\begin{equation*}
\Ll_s(A_{k_i},B_{k_i})=\int_{\R^n}\int_{\R^n}\frac{1}{|x-y|^{n+s}}\chi_{A_{k_i}}(x)\chi_{B_{k_i}}(y)\,dx\,dy,
\end{equation*}
Fatou's Lemma gives
\begin{equation*}
\Ll_s(A,B)\leq\liminf_{i\to\infty}\Ll_s(A_{k_i},B_{k_i})=I,
\end{equation*}
proving \eqref{CH:2:semicontinuity_first}.

The second inequality follows just by summing the contributions defining the fractional perimeter.
\end{proof}

\noindent
Keeping $\Omega$ fixed we obtain \cite[Theorem 3.1]{CRS10}.

On the other hand, if we keep the set $E$ fixed and approximate the open set $\Omega$ with a
sequence of open subsets $\Omega_k\subseteq\Omega$, we get a continuity property.

\begin{prop}\label{CH:2:continuity_in_open_set_seq}
Let $\Omega\subseteq\R^n$ be an open set and let $\{\Omega_k\}$ be any sequence of open sets such that
$\Omega_k\xrightarrow{loc}\Omega$.
Then for every set $E\subseteq\R^n$
\begin{equation*}
\Per_s(E,\Omega)\leq\liminf_{k\to\infty}\Per_s(E,\Omega_k).
\end{equation*}
Moreover, if $\Omega_k\subseteq\Omega$ for every $k$, then
\begin{equation}\label{CH:2:limit_sub_open}
\Per_s(E,\Omega)=\lim_{k\to\infty}\Per_s(E,\Omega_k),
\end{equation}
(whether it is finite or not).
\end{prop}

\begin{proof}
Since $\Omega_k\xrightarrow{loc}\Omega$,
Proposition \ref{CH:2:semicont_first_prop} gives the first statement.
Now notice that if $\Omega_k\subseteq\Omega$, Proposition \ref{CH:2:subopensets} implies
\begin{equation*}
\Per_s(E,\Omega_k)\leq \Per_s(E,\Omega),
\end{equation*}
and hence
\begin{equation*}
\limsup_{k\to\infty}\Per_s(E,\Omega_k)\leq \Per_s(E,\Omega),
\end{equation*}
concluding the proof.
\end{proof}

\begin{remark}\label{CH:2:continuity_in_open_set_seq_rmk}
As a consequence, exploiting Corollary \ref{CH:2:regular_approx_open_sets_coroll}, we get
\begin{equation*}
\Per_s(E,\Omega)=\sup_{\Omega'\subsetneq\Omega}\Per_s(E,\Omega')
=\sup_{\Omega'\Subset\Omega}\Per_s(E,\Omega').
\end{equation*}
\end{remark}

\begin{remark}\label{CH:2:counterex_loc_fin_per}
Consider the set $E\subseteq\R$ constructed in the proof of \cite[Example 2.10]{DFPV13}.
That is, let $\beta_k>0$ be a decreasing sequence such that
\[M:=\sum_{k=1}^\infty\beta_k<\infty\quad\textrm{and}
\quad\sum_{k=1}^\infty\beta_{2k}^{1-s}=\infty,\quad\forall\,s\in(0,1).\]
Then define
\[\sigma_m:=\sum_{k=1}^m\beta_k,\qquad I_m:=(\sigma_m,\sigma_{m+1}),\qquad E:=\bigcup_{j=1}^\infty I_{2j},\]
and let $\Omega:=(0,M)$. As shown in \cite{DFPV13},
\[\Per_s(E,\Omega)=\infty,\qquad\forall\,s\in(0,1).\]
On the other hand
\[\Per(E,\Omega')<\infty,\qquad\forall\,\Omega'\Subset\Omega,\]
hence $E$ has locally finite $s$-perimeter in $\Omega$, for every $s\in(0,1)$.

Indeed, notice that the intervals $I_{2j}$ accumulate near $M$. Thus, for every $\eps>0$,
all but a finite number of the intervals $I_{2j}$'s fall outside of the open set
$\mathcal O_\eps:=(\eps,M-\eps)$. Therefore
$\Per(E,\mathcal O_\eps)<\infty$ and hence
\[\Per_s(E,\mathcal O_\eps)<\infty,\quad\forall\,s\in(0,1).\]
Since $\mathcal O_\eps\nearrow\Omega$ as $\eps\to0^+$, the set $E$ has locally finite $s$-perimeter in $\Omega$
for every $s\in(0,1)$.
\end{remark}

\begin{prop}\label{CH:2:subcont_lem_approx}
Let $\Omega\subseteq\R^n$ be an open set and let $\{E_h\}$ be a sequence of sets such that
\begin{equation*}
E_h\xrightarrow{loc}E\qquad\textrm{and}\qquad \lim_{h\to\infty}\Per_s(E_h,\Omega)=\Per_s(E,\Omega)<\infty.
\end{equation*}
Then
\begin{equation*}
\lim_{h\to\infty}\Per_s(E_h,\Omega')=\Per_s(E,\Omega')\qquad\textrm{for every open set }\Omega'\subseteq\Omega.
\end{equation*}
\end{prop}

\begin{proof}
The claim follows from classical properties of limits of sequences.

Indeed, let
\begin{equation*}
\begin{split}
&\qquad\qquad\qquad\qquad\qquad a_h:=\Per_s(E_h,\Omega'),\\
&
b_h:=\Ll_s\big(E_h\cap(\Omega\setminus\Omega'),\Co E_h\setminus\Omega\big)
+\Ll_s\big(E_h\setminus\Omega,\Co E_h\cap(\Omega\setminus\Omega')\big)\\
&
\qquad\qquad\quad+\Ll_s\big(E_h\cap(\Omega\setminus\Omega'),\Co E_h\cap(\Omega\setminus\Omega')\big),
\end{split}\end{equation*}
and let $a$ and $b$ be the corresponding terms for $E$.\\
Notice that, by Proposition \ref{CH:2:subopensets}, we have
\begin{equation*}
\Per_s(E_h,\Omega)=a_h+b_h\qquad\textrm{and}\qquad \Per_s(E,\Omega)=a+b.
\end{equation*}
From Proposition \ref{CH:2:semicont_first_prop} we have
\begin{equation*}
a\leq\liminf_{h\to\infty}a_h\qquad\textrm{and}\qquad b\leq\liminf_{h\to\infty}b_h,
\end{equation*}
and by hypothesis we know that
\begin{equation*}
\lim_{h\to\infty}(a_h+b_h)=a+b.
\end{equation*}
Therefore
\begin{equation*}
a+b\leq\liminf_{h\to\infty}a_h+\liminf_{h\to\infty}b_h\leq\liminf_{h\to\infty}(a_h+b_h)=a+b,
\end{equation*}
and hence
\begin{equation*}
0\leq\liminf_{h\to\infty}b_h-b=a-\liminf_{h\to\infty}a_h\leq0,
\end{equation*}
so that
\begin{equation*}
a=\liminf_{h\to\infty}a_h\qquad\textrm{and}\qquad b=\liminf_{h\to\infty}b_h.
\end{equation*}
Then, since
\begin{equation*}
\limsup_{h\to\infty}a_h+\liminf_{h\to\infty}b_h\leq\limsup_{h\to\infty}(a_h+b_h)=a+b,
\end{equation*}
we obtain
\begin{equation*}
a=\liminf_{h\to\infty}a_h\leq\limsup_{h\to\infty}a_h\leq a,
\end{equation*}
concluding the proof.
\end{proof}

\subsection{Compactness}\label{CH:2:Compactness_Section}


\begin{prop}[Compactness]\label{CH:2:compact_prop}
Let $\Omega\subseteq\R^n$ be an open set.
If $\{E_h\}$ is a sequence of sets such that
\begin{equation}\label{CH:2:compact_seq_diag_eq}
\limsup_{h\to\infty}\Per_s^L(E_h,\Omega')\leq c(\Omega')<\infty,\quad\forall\,\Omega'\Subset\Omega,
\end{equation}
then there exists a subsequence $\{E_{h_i}\}$ and $E\subseteq\R^n$ such that
\begin{equation*}
E_{h_i}\cap\Omega\xrightarrow{loc}E\cap\Omega.
\end{equation*}
\end{prop}

\begin{proof}
We want to use a compact Sobolev embedding (see, e.g., \cite[Corollary 7.2]{HitGuide}) to construct a limit set via a diagonal argument.

Thanks to Corollary \ref{CH:2:regular_approx_open_sets_coroll} we know that we can find an increasing sequence of bounded open sets $\{\Omega_k\}$ with smooth boundary such that
\begin{equation*}
\Omega_k\Subset\Omega_{k+1}\Subset\Omega\qquad\textrm{and}\qquad\bigcup_{k\in\mathbb N}\Omega_k=\Omega.
\end{equation*}
Moreover, hypothesis \eqref{CH:2:compact_seq_diag_eq} guarantees that
\begin{equation}\label{CH:2:compact_seq_diag_eq2}
\forall k\quad\exists h(k)\textrm{ s.t.}\quad
\Per_s^L(E_h,\Omega_k)\leq c_k<\infty,\quad\forall h\geq h(k).
\end{equation}
Clearly
\begin{equation*}
\|\chi_{E_h}\|_{L^1(\Omega_k)}\leq|\Omega_k|<\infty,
\end{equation*}
and hence, since $[\chi_{E_h}]_{W^{s,1}(\Omega_k)}=2\Per_s^L(E_h,\Omega_k)$,  we have
\begin{equation*}
\|\chi_{E_h}\|_{W^{s,1}(\Omega_k)}\leq c'_k,\quad\forall h\geq h(k).
\end{equation*}
Therefore \cite[Corollary 7.2]{HitGuide} (notice that each $\Omega_k$ is an extension domain) guarantees for every fixed $k$ the existence of a subsequence $h_i\nearrow\infty$ (with $h_1\geq
h(k)$) such that
\begin{equation*}
E_{h_i}\cap\Omega_k\xrightarrow{i\to\infty} E^k
\end{equation*}
in measure, for some set $E^k\subseteq\Omega_k$.

Applying this argument for $k=1$ we get a subsequence $\{h_i^1\}$ with
\begin{equation*}
E_{h_i^1}\cap\Omega_1\xrightarrow{i\to\infty}E^1.
\end{equation*}
Applying again this argument in $\Omega_2$, with $\{E_{h_i^1}\}$ in place of $\{E_h\}$, we get a subsequence
$\{h_i^2\}$ of $\{h_i^1\}$ with
\begin{equation*}
E_{h_i^2}\cap\Omega_2\xrightarrow{i\to\infty}E^2.
\end{equation*}
Notice that, since $\Omega_1\subseteq\Omega_2$, we must have $E^2\cap\Omega_1=E^1$ in measure (by the uniqueness of the limit in $\Omega_1$). We can also suppose that $h_1^2>h_1^1$.\\
Proceeding inductively in this way we get an increasing subsequence $\{h_1^k\}$ such that
\begin{equation*}
E_{h_1^i}\cap\Omega_k\xrightarrow{i\to\infty}E^k,\qquad\textrm{for every }k\in\mathbb{N},
\end{equation*}
with $E^{k+1}\cap\Omega_k=E^k$. Therefore if we define $E:=\bigcup_kE^k$, since $\bigcup_k\Omega_k=\Omega$, we get
\begin{equation*}
E_{h_1^i}\cap\Omega\xrightarrow{loc}E,
\end{equation*}
concluding the proof.
\end{proof}

\begin{remark}\label{CH:2:min_app_seq_rmk}
If $E_h$ is $s$-minimal in $\Omega_k$ for every $h\geq h(k)$, then by minimality we get
\begin{equation*}
\Per_s^L(E_h,\Omega_k)\leq \Per_s(E_h,\Omega_k)\leq \Per_s(E_h\setminus\Omega_k,\Omega_k)\leq \Per_s(\Omega_k)=:c_k<\infty,
\end{equation*}
since $\Omega_k$ is bounded and has Lipschitz boundary. Therefore $\{E_h\}$ satisfies the hypothesis of Proposition \ref{CH:2:compact_prop} and we can find a convergent subsequence.
\end{remark}

\section{Generalized coarea and approximation by smooth sets}\label{CH:2:GCAABSS_SEC}

We begin by showing that the $s$-perimeter satisfies a generalized coarea formula (see also \cite{Visintin} and \cite[Lemma 10]{Gamma}).
In the end of this section we will exploit this formula to prove that a set $E$ of locally finite $s$-perimeter can be approximated by smooth sets whose
$s$-perimeter converges to that of $E$.

\smallskip

Let $\Omega\subseteq\R^n$ be an open set. Given a function $u:\R^n\longrightarrow\R$,
we define the functional
\begin{equation}\label{CH:2:def_ext_func_coarea}
\Fc(u,\Omega):=\frac{1}{2}\int_\Omega\int_\Omega\frac{|u(x)-u(y)|}{|x-y|^{n+s}}dx\,dy
+\int_\Omega\int_{\Co\Omega}\frac{|u(x)-u(y)|}{|x-y|^{n+s}}dx\,dy,
\end{equation}
that is, half the ``$\Omega$-contribution'' to the $W^{s,1}$-seminorm of $u$.\\
Notice that
\begin{equation*}
\Fc(\chi_E,\Omega)=\Per_s(E,\Omega)
\end{equation*}
and, clearly
\begin{equation*}
\Fc(u,\R^n)=\frac{1}{2}[u]_{W^{s,1}(\R^n)}.
\end{equation*}

\begin{prop}[Coarea]\label{CH:2:coarea_prop}
Let $\Omega\subseteq\R^n$ be an open set and let $u:\R^n\longrightarrow\R$.
Then
\begin{equation}\label{CH:2:coarea_formula}
\Fc(u,\Omega)=\int_{-\infty}^\infty \Per_s(\{u>t\},\Omega)\,dt.
\end{equation}
In particular
\begin{equation*}
\frac{1}{2}[u]_{W^{s,1}(\Omega)}=\int_{-\infty}^\infty \Per_s^L(\{u>t\},\Omega)\,dt.
\end{equation*}
\end{prop}

\begin{proof}
Notice that for every $x,\, y\in\R^n$ we have
\begin{equation}\label{CH:2:coarea_first_formula}
|u(x)-u(y)|=\int_{-\infty}^\infty|\chi_{\{u>t\}}(x)-\chi_{\{u>t\}}(y)|\,dt.
\end{equation}
Indeed, the function $t\longmapsto|\chi_{\{u>t\}}(x)-\chi_{\{u>t\}}(y)|$ takes only the values $\{0,1\}$
and it is different from 0 precisely in the interval having $u(x)$ and $u(y)$ as extremes.
Therefore, if we plug \eqref{CH:2:coarea_first_formula} into \eqref{CH:2:def_ext_func_coarea} and use Fubini's Theorem, we get
\begin{equation*}
\Fc(u,\Omega)=\int_{-\infty}^\infty\Fc(\chi_{\{u>t\}},\Omega)\,dt=\int_{-\infty}^\infty \Per_s(\{u>t\},\Omega)\,dt,
\end{equation*}
as wanted.
\end{proof}

\subsection{Approximation results for the functional $\Fc$}\label{CH:2:Approx_Section}

In this section we prove the approximation properties for the functional $\Fc$ which we
need for the proofs of Theorem \ref{CH:2:density_smooth_teo} and Theorem \ref{CH:2:appro_in_bded_open}.
To this end we consider a (symmetric) smooth function $\eta$ such that
\begin{equation*}
\eta\in C^\infty_c(\R^n),\quad\textrm{supp }\eta\subseteq B_1,\quad\eta\geq0,\quad\eta(-x)=\eta(x),\quad\int_{\R^n}\eta\,dx=1,
\end{equation*}
and we define the mollifier
\begin{equation*}
\eta_\eps(x):=\frac{1}{\eps^n}\eta\Big(\frac{x}{\eps}\Big),
\end{equation*}
for every $\eps\in(0,1)$.
Notice that supp $\eta_\eps\subseteq B_\eps$ and $\int_{\R^n}\eta_\eps=1$.

Given $u\in L^1_{\loc}(\R^n)$, we define the $\eps$-regularization of $u$ as the convolution
\begin{equation*}
u_\eps(x):=(u\ast\eta_\eps)(x)=\int_{\R^n}u(x-\xi)\eta_\eps(\xi)\,d\xi,\quad\textrm{for every }x\in\R^n.
\end{equation*}
It is well known that $u_\eps\in C^\infty(\R^n)$ and
\begin{equation*}
u_\eps\longrightarrow u\qquad\textrm{in }L^1_{\loc}(\R^n).
\end{equation*}
Moreover, if $u=\chi_E$, then
\begin{equation}\label{CH:2:forml1_smoothing}
0\leq u_\eps\leq1\qquad\textrm{and}\quad u_\eps(x)=\left\{\begin{array}{cc}1, & \textrm{if }|B_\eps(x)\setminus E|=0\\
0, & \textrm{if }|B_\eps(x)\cap E|=0\end{array}\right.,
\end{equation}
(see, e.g., \cite[Section 12.3]{Maggi}).

\begin{lemma}\label{CH:2:dens_lemma}
$(i)\quad$ Let $u\in L^1_{\loc}(\R^n)$ and let $\Omega\subseteq\R^n$ be an open set. Then
\begin{equation}\label{CH:2:forml_lemma_coarea}
\Fc(u,\Omega)<\infty\quad\Longrightarrow\quad\lim_{\eps\to0^+}\Fc(u_\eps,\Omega')=\Fc(u,\Omega')\qquad\forall\,\Omega'\Subset\Omega.
\end{equation}
$(ii)\quad$ Let $u\in W^{s,1}(\R^n)$. Then
\begin{equation*}
\lim_{\eps\to0}[u_\eps]_{W^{s,1}(\R^n)}=[u]_{W^{s,1}(\R^n)}.
\end{equation*}
$(iii)\quad$ Let $u\in W^{s,1}(\R^n)$. Then there exists $\{u_k\}\subseteq C_c^\infty(\R^n)$ such that
\begin{equation*}
\|u-u_k\|_{L^1(\R^n)}\longrightarrow0\quad\textrm{and}\quad\lim_{k\to\infty}[u_k]_{W^{s,1}(\R^n)}=[u]_{W^{s,1}(\R^n)}.
\end{equation*}
Moreover, if $u=\chi_E$, then $0\leq u_k\leq1$.
\end{lemma}

\begin{proof}

$(i)\quad$
Given $\mathcal O\subseteq\R^n$, let $Q(\mathcal O):=\R^{2n}\setminus(\Co\mathcal O)^2$, so that
\begin{equation*}
\Fc(u,\mathcal O)=\frac{1}{2}\iint_{Q(\mathcal O)}\frac{|u(x)-u(y)|}{|x-y|^{n+s}}\,dx\,dy.
\end{equation*}
 Notice that if $\mathcal O\subseteq\Omega$, then
$Q(\mathcal O)\subseteq Q(\Omega)$
and hence
\begin{equation}\label{CH:2:forml1_coarea}
\Fc(u,\mathcal O)\leq\Fc(u,\Omega).
\end{equation}
Now let $\Omega'\Subset\Omega$ and notice that for $\eps$ small enough we have
\begin{equation}\label{CH:2:forml2_coarea}
Q(\Omega'-\eps\xi)\subseteq Q(\Omega)\qquad\textrm{for every }\xi\in B_1.
\end{equation}
As a consequence
\begin{equation}\label{CH:2:forml4_coarea}
\Fc(u_\eps,\Omega')\leq\int_{B_1}\Fc(u,\Omega'-\eps\xi)\eta(\xi)\,d\xi\leq\Fc(u,\Omega).
\end{equation}
The second inequality follows from \eqref{CH:2:forml2_coarea}, \eqref{CH:2:forml1_coarea} and $\int_{B_1}\eta=1$.\\
As for the first inequality, we have
\begin{equation*}\begin{split}
\iint_{Q(\Omega')}&\frac{|u_\eps(x)-u_\eps(y)|}{|x-y|^{n+s}}dx\,dy\\
&
=\iint_{Q(\Omega')}\Big|\int_{\R^n}\big(u(x-\xi)-u(y-\xi)\big)\frac{1}{\eps^n}\eta\Big(\frac{\xi}{\eps}\Big)\,d\xi\Big|\frac{dx\,dy}{|x-y|^{n+s}}\\
&
=\iint_{Q(\Omega')}\Big|\int_{B_1}\big(u(x-\eps\xi)-u(y-\eps\xi)\big)\eta(\xi)\,d\xi\Big|\frac{dx\,dy}{|x-y|^{n+s}}\\
&
\leq\int_{B_1}\Big(\iint_{Q(\Omega')}\frac{|u(x-\eps\xi)-u(y-\eps\xi)|}{|x-y|^{n+s}}dx\,dy\Big)\eta(\xi)\,d\xi\\
&
=\int_{B_1}\Big(\iint_{Q(\Omega'-\eps\xi)}\frac{|u(x)-u(y)|}{|x-y|^{n+s}}dx\,dy\Big)\eta(\xi)\,d\xi.
\end{split}\end{equation*}


We prove something stronger than the claim, that is
\begin{equation}\label{CH:2:forml_app_98}
\lim_{\eps\to0^+}\Fc(u_\eps-u,\Omega')=0.
\end{equation}
Indeed, notice that
\begin{equation*}
|\Fc(u_\eps,\Omega')-\Fc(u,\Omega')|\leq\Fc(u_\eps-u,\Omega').
\end{equation*}
Let $\psi:\R^{2n}\longrightarrow\R$ be defined as
\begin{equation*}
\psi(x,y):=\frac{u(x)-u(y)}{|x-y|^{n+s}}.
\end{equation*}
Moreover, for every $\eps>0$ and $\xi\in B_1$, we consider the left translation by $\eps(\xi,\xi)$
in $\R^{2n}$, that is
\begin{equation*}
(L_{\eps\xi}f)(x,y):=f(x-\eps\xi,y-\eps\xi),
\end{equation*}
for every $f:\R^{2n}\longrightarrow\R$.\\
Since $\psi\in L^1(Q(\Omega))$, for every $\delta>0$ there exists $\Psi\in C_c^1(Q(\Omega))$ such that
\begin{equation*}
\|\psi-\Psi\|_{L^1(Q(\Omega))}\leq\frac{\delta}{2}.
\end{equation*}
We have
\begin{equation*}\begin{split}
\Fc(u_\eps&-u,\Omega')
=\iint_{Q(\Omega')}\frac{|u_\eps(x)-u_\eps(y)-u(x)+u(y)|}{|x-y|^{n+s}}dx\,dy\\
&
\leq\int_{B_1}\Big(\iint_{Q(\Omega')}\frac{|u(x-\eps\xi)-u(y-\eps\xi)-u(x)+u(y)|}{|x-y|^{n+s}}dx\,dy\Big)\eta(\xi)\,d\xi\\
&
=\int_{B_1}\|L_{\eps\xi}\psi-\psi\|_{L^1(Q(\Omega'))}\eta(\xi)\,d\xi\\
&
\leq\int_{B_1}\Big(\|L_{\eps\xi}\psi-L_{\eps\xi}\Psi\|_{L^1(Q(\Omega'))}
+\|L_{\eps\xi}\Psi-\Psi\|_{L^1(Q(\Omega'))}\\
&
\qquad\qquad\qquad\qquad
+\|\Psi-\psi\|_{L^1(Q(\Omega'))}\Big)\eta(\xi)\,d\xi.
\end{split}\end{equation*}
Notice that
\begin{equation*}
\|L_{\eps\xi}\psi-L_{\eps\xi}\Psi\|_{L^1(Q(\Omega'))}
=\|\psi-\Psi\|_{L^1(Q(\Omega'-\eps\xi))}
\leq\|\psi-\Psi\|_{L^1(Q(\Omega))}
\end{equation*}
and hence
\begin{equation*}
\Fc(u_\eps-u,\Omega')
\leq\delta+\int_{B_1}\|L_{\eps\xi}\Psi-\Psi\|_{L^1(Q(\Omega'))}\eta(\xi)\,d\xi.
\end{equation*}
For $\eps>0$ small enough we have
\begin{equation*}
\textrm{supp}(L_{\eps\xi}\Psi-\Psi)\subseteq N_1(\textrm{supp }\Psi)=:K\Subset\R^{2n},
\end{equation*}
and
\begin{equation*}
|\Psi(x-\eps\xi,y-\eps\xi)-\Psi(x,y)|\leq2\max_{\textrm{supp }\Psi}|\nabla\Psi|\,\eps.
\end{equation*}
Thus
\begin{equation*}
\int_{B_1}\|L_{\eps\xi}\Psi-\Psi\|_{L^1(Q(\Omega'))}\eta(\xi)\,d\xi
\leq2|K|\max_{\textrm{supp }\Psi}|\nabla\Psi|\,\eps.
\end{equation*}
Passing to the limit as $\eps\to0^+$ then gives
\begin{equation*}
\limsup_{\eps\to0^+}\Fc(u_\eps-u,\Omega')\leq\delta.
\end{equation*}
Since $\delta$ is arbitrary, we get \eqref{CH:2:forml_app_98}.

$(ii)\quad$ Reasoning as above we obtain
\begin{equation*}\begin{split}
\int_{\R^n}&\int_{\R^n}\frac{|u_\eps(x)-u_\eps(y)|}{|x-y|^{n+s}}dx\,dy\\
&
\leq\int_{B_1}\Big(\int_{\R^n}\int_{\R^n}\frac{|u(x-\eps\xi)-u(y-\eps\xi)|}{|x-y|^{n+s}}dx\,dy\Big)\eta(\xi)\,d\xi\\
&
=\int_{B_1}\Big(\int_{\R^n}\int_{\R^n}\frac{|u(x)-u(y)|}{|x-y|^{n+s}}dx\,dy\Big)\eta(\xi)\,d\xi\\
&
=[u]_{W^{s,1}(\R^n)}\int_{B_1}\eta(\xi)\,d\xi,
\end{split}\end{equation*}
that is
\begin{equation*}
[u_\eps]_{W^{s,1}(\R^n)}\leq[u]_{W^{s,1}(\R^n)}.
\end{equation*}
This and Fatou's Lemma give
\begin{equation*}
[u]_{W^{s,1}(\R^n)}\leq\liminf_{\eps\to0}[u_\eps]_{W^{s,1}(\R^n)}\leq\limsup_{\eps\to0}[u_\eps]_{W^{s,1}(\R^n)}
\leq[u]_{W^{s,1}(\R^n)},
\end{equation*}
concluding the proof.

$(iii)\quad$ The proof is a classical cut-off argument. We consider a sequence of cut-off functions $\psi_k\in C_c^\infty(\R^n)$
such that
\begin{equation*}
0\leq\psi_k\leq1,\quad\textrm{supp }\psi_k\subseteq B_{k+1}\quad\textrm{and}\quad\psi_k\equiv1\quad\textrm{in }B_k.
\end{equation*}
We can also assume that
\begin{equation*}
\sup_{k\in\mathbb N}|\nabla\psi_k|\leq M_0<\infty.
\end{equation*}
It is enough to show that
\begin{equation}\label{CH:2:forml_lemma_approx1}
\lim_{k\to\infty}\|u-\psi_k u\|_{L^1(\R^n)}=0\quad\textrm{and}\quad\lim_{k\to\infty}[\psi_k u]_{W^{s,1}(\R^n)}=
[u]_{W^{s,1}(\R^n)}.
\end{equation}

Indeed then we can use $(ii)$ to approximate each $\psi_k u$ with a smooth function $u_k:=(u\psi_k)\ast\eta_{\eps_k}$,
for $\eps_k$ small enough to have
\begin{equation*}
\|\psi_k u-u_k\|_{L^1(\R^n)}<2^{-k}\quad\textrm{and}\quad|[\psi_k u]_{W^{s,1}(\R^n)}-[u_k]_{W^{s,1}(\R^n)}|<2^{-k}.
\end{equation*}
Therefore
\begin{equation*}
\|u-u_k\|_{L^1(\R^n)}
\leq\|u-\psi_k u\|_{L^1(\R^n)}+2^{-k}\longrightarrow0
\end{equation*}
and
\begin{equation*}
|[u]_{W^{s,1}(\R^n)}-[u_k]_{W^{s,1}(\R^n)}|\leq|[u]_{W^{s,1}(\R^n)}-[\psi_k u]_{W^{s,1}(\R^n)}|+2^{-k}\longrightarrow0.
\end{equation*}
Also notice that
\begin{equation*}
\textrm{supp }u_k\subseteq N_{\eps_k}(\textrm{supp }\psi_k u)\subseteq B_{k+2}
\end{equation*}
so that $u_k\in C_c^\infty(\R^n)$ for every $k$. Moreover, from the definition of $u_k$ it follows that if $u=\chi_E$, then
$0\leq u_k\leq1$.

For a proof of \eqref{CH:2:forml_lemma_approx1} see, e.g., \cite[Lemma 12]{FSV15}.
\end{proof}


Now we show that if $\Omega$ is a bounded open set with Lipschitz boundary and if $u=\chi_E$,
then we can find smooth functions $u_h$ such that
\[\Fc(u_h,\Omega)\longrightarrow\Fc(u,\Omega).\]

We first need the following two results.

\begin{lemma}\label{CH:2:stuff_forml1}
Let $\Omega\subseteq\R^n$ be a bounded open set with Lipschitz boundary. Let $u\in L^\infty(\R^n)$ be such that
$\Fc(u,\Omega)<\infty$. For every $\delta\in(0,r_0)$ let
\begin{equation*}
\varphi_\delta:=1-\chi_{\{|\bar{d}_\Omega|<\delta\}}.
\end{equation*}
Then
\begin{equation}\label{CH:2:global_bded_conv_l1}
u\varphi_\delta\xrightarrow{\delta\to0}u\quad\textrm{in }L^1(\R^n),
\end{equation}
and
\begin{equation*}
\lim_{\delta\searrow0^+}\Fc(u\varphi_\delta,\Omega)=\Fc(u,\Omega).
\end{equation*}
\end{lemma}

\begin{proof}
First of all, notice that
\begin{equation*}
\int_{\R^n}|u\varphi_\delta-u|\,dx=
\int_{\{|\bar{d}_\Omega|<\delta\}}|u|\,dx\leq\|u\|_{L^\infty(\R^n)}\,|\{|\bar{d}_\Omega|<\delta\}|\xrightarrow{\delta\to0}0.
\end{equation*}
Now
\begin{equation*}\begin{split}
\int_\Omega&\int_\Omega\frac{|(u\varphi_\delta)(x)-(u\varphi_\delta)(y)|}{|x-y|^{n+s}}\,dx\,dy\\
&
=
\int_{\Omega_{-\delta}}\int_{\Omega_{-\delta}}\frac{|u(x)-u(y)|}{|x-y|^{n+s}}\,dx\,dy
+
2\int_{\Omega_{-\delta}}\Big(\int_{\Omega\setminus\Omega_{-\delta}}\frac{|u(x)|}{|x-y|^{n+s}}\,dy\Big)dx.
\end{split}\end{equation*}
Since $\Omega_{-\delta}\subseteq\Omega$, we have
\begin{equation*}
\int_{\Omega_{-\delta}}\int_{\Omega_{-\delta}}\frac{|u(x)-u(y)|}{|x-y|^{n+s}}\,dx\,dy
\leq
\int_\Omega\int_\Omega\frac{|u(x)-u(y)|}{|x-y|^{n+s}}\,dx\,dy.
\end{equation*}
On the other hand, since $|\Omega\setminus\Omega_{-\delta}|\longrightarrow0$, we get
\begin{equation*}
\frac{|u(x)-u(y)|}{|x-y|^{n+s}}\chi_{\Omega_{-\delta}}(x)\chi_{\Omega_{-\delta}}(y)\xrightarrow{\delta\to0}
\frac{|u(x)-u(y)|}{|x-y|^{n+s}}\chi_\Omega(x)\chi_\Omega(y),
\end{equation*}
for a.e. $(x,y)\in\R^n\times\R^n$.

Therefore, by Fatou's Lemma we obtain
\begin{equation}\label{CH:2:stuff_formul1}
[u]_{W^{s,1}(\Omega)}\leq\liminf_{\delta\searrow0}[u]_{W^{s,1}(\Omega_{-\delta})}
\leq\limsup_{\delta\searrow0}[u]_{W^{s,1}(\Omega_{-\delta})}\leq[u]_{W^{s,1}(\Omega)}.
\end{equation}
Moreover, by point $(i)$ of \eqref{CH:2:uniform_bound_strips} we get
\begin{equation*}\begin{split}
2\int_{\Omega_{-\delta}}\Big(\int_{\Omega\setminus\Omega_{-\delta}}\frac{|u(x)|}{|x-y|^{n+s}}\,dy\Big)dx&
\leq2\|u\|_{L^\infty(\R^n)}\Ll_s(\Omega_{-\delta},\Omega\setminus\Omega_{-\delta})\\
&
\leq 2C\|u\|_{L^\infty(\R^n)}\,\delta^{1-s}.
\end{split}\end{equation*}
Therefore we find
\begin{equation*}
\lim_{\delta\searrow0}[u\varphi_\delta]_{W^{s,1}(\Omega)}=[u]_{W^{s,1}(\Omega)}.
\end{equation*}

Now
\begin{equation*}\begin{split}
\int_\Omega&\int_{\Co\Omega}\frac{|(u\varphi_\delta)(x)-(u\varphi_\delta)(y)|}{|x-y|^{n+s}}\,dx\,dy\\
&
=
\int_{\Omega_{-\delta}}\int_{\Co\Omega_\delta}\frac{|u(x)-u(y)|}{|x-y|^{n+s}}\,dx\,dy
+
\int_{\Omega_{-\delta}}\Big(\int_{\Omega_\delta\setminus\Omega}\frac{|u(x)|}{|x-y|^{n+s}}\,dy\Big)dx\\
&
\qquad\qquad+
\int_{\Omega\setminus\Omega_{-\delta}}\Big(\int_{\Co\Omega_\delta}\frac{|u(x)|}{|x-y|^{n+s}}\,dy\Big)dx.
\end{split}\end{equation*}
Since $\Omega_{-\delta}\subseteq\Omega$ and $\Co\Omega_\delta\subseteq\Co\Omega$, we have
\begin{equation*}
\int_{\Omega_{-\delta}}\int_{\Co\Omega_\delta}\frac{|u(x)-u(y)|}{|x-y|^{n+s}}\,dx\,dy
\leq
\int_\Omega\int_{\Co\Omega}\frac{|u(x)-u(y)|}{|x-y|^{n+s}}\,dx\,dy.
\end{equation*}
Moreover, since both $|\Omega\setminus\Omega_{-\delta}|\longrightarrow0$ and
$|\Co\Omega\setminus\Co\Omega_\delta|\longrightarrow0$, we have
\begin{equation*}
\frac{|u(x)-u(y)|}{|x-y|^{n+s}}\chi_{\Omega_{-\delta}}(x)\chi_{\Co\Omega_\delta}(y)\xrightarrow{\delta\to0}
\frac{|u(x)-u(y)|}{|x-y|^{n+s}}\chi_\Omega(x)\chi_{\Co\Omega}(y),
\end{equation*}
for a.e. $(x,y)\in\R^n\times\R^n$.

Therefore, again by Fatou's Lemma we obtain
\begin{equation*}
\lim_{\delta\searrow0}\int_{\Omega_{-\delta}}\int_{\Co\Omega_\delta}\frac{|u(x)-u(y)|}{|x-y|^{n+s}}\,dx\,dy
=
\int_\Omega\int_{\Co\Omega}\frac{|u(x)-u(y)|}{|x-y|^{n+s}}\,dx\,dy.
\end{equation*}

Furthermore, by point $(ii)$ of \eqref{CH:2:uniform_bound_strips} we get
\begin{equation*}\begin{split}
\int_{\Omega_{-\delta}}\Big(&\int_{\Omega_\delta\setminus\Omega}\frac{|u(x)|}{|x-y|^{n+s}}\,dy\Big)dx
\leq\|u\|_{L^\infty(\R^n)}\Ll_s(\Omega_{-\delta},\Omega_\delta\setminus\Omega)\\
&
\leq\|u\|_{L^\infty(\R^n)}\Ll_s(\Omega,\Omega_\delta\setminus\Omega)
\leq C\|u\|_{L^\infty(\R^n)}\delta^{1-s}
\end{split}\end{equation*}
and also
\begin{equation*}
\int_{\Omega\setminus\Omega_{-\delta}}\Big(\int_{\Co\Omega_\delta}\frac{|u(x)|}{|x-y|^{n+s}}\,dy\Big)dx
\leq C\|u\|_{L^\infty(\R^n)}\delta^{1-s}.
\end{equation*}
Thus
\begin{equation*}
\lim_{\delta\searrow0}\int_\Omega\int_{\Co\Omega}\frac{|(u\varphi_\delta)(x)-(u\varphi_\delta)(y)|}{|x-y|^{n+s}}\,dx\,dy
=\int_\Omega\int_{\Co\Omega}\frac{|u(x)-u(y)|}{|x-y|^{n+s}}\,dx\,dy,
\end{equation*}
concluding the proof.
\end{proof}

\begin{lemma}\label{CH:2:second_lemma_for_func_appro}
Let $\Omega\subseteq\R^n$ be a bounded open set with Lipschitz boundary. Let $v\in L^\infty(\R^n)$ be
such that $\Fc(v,\Omega)<\infty$ and
\begin{equation*}
v\equiv0\quad\textrm{in }\{|\bar{d}_\Omega|<\delta/2\},
\end{equation*}
for some $\delta\in(0,r_0)$. Then
\begin{equation*}
\big|\Fc(v,\Omega)-\Fc(v,\Omega_{-\delta/2})\big|\leq C\|v\|_{L^\infty(\R^n)}\delta^{1-s},
\end{equation*}
where $C=C(n,s,\Omega)>0$ does not depend on $v$.
\end{lemma}

\begin{proof}
Since
\begin{equation*}
v\equiv0\quad\textrm{in }\{|\bar{d}_\Omega|<\delta/2\},
\end{equation*}
we have
\begin{equation*}
\Fc(v,\Omega)=\Fc(v,\Omega_{-\delta/2})
+2\int_{\Omega\setminus\Omega_{-\delta/2}}\Big(\int_{\Co\Omega_{\delta/2}}\frac{|v(y)|}{|x-y|^{n+s}}\,dy\Big)dx.
\end{equation*}
Now, by point $(ii)$ of \eqref{CH:2:uniform_bound_strips} we have
\begin{equation*}\begin{split}
\int_{\Omega\setminus\Omega_{-\delta/2}}\Big(\int_{\Co\Omega_{\delta/2}}\frac{|v(y)|}{|x-y|^{n+s}}\,dy\Big)&
\leq \|v\|_{L^\infty(\R^n)}\Ll_s(\Omega\setminus\Omega_{-\delta/2},\Co\Omega)\\
&
\leq 2^{s-1}C\|v\|_{L^\infty(\R^n)}\,\delta^{1-s}.
\end{split}
\end{equation*}
\end{proof}

\begin{prop}\label{CH:2:density_bded_reg_set}
Let $\Omega\subseteq\R^n$ be a bounded open set with Lipschitz boundary. Let $u\in L^\infty(\R^n)$ be such that
$\Fc(u,\Omega)<\infty$.
Then there exists a sequence $\{u_h\}\subseteq C^\infty(\R^n)$ such that
\begin{equation*}\begin{split}
&(i)\quad \|u_h\|_{L^\infty(\R^n)}\leq\|u\|_{L^\infty(\R^n)},\quad\textrm{and}\quad0\leq u_h\leq1\quad\textrm{if}\quad
0\leq u\leq1,\\
&
(ii)\quad
u_h\xrightarrow{h\to\infty}u\quad\textrm{in }L^1_{\loc}(\R^n),\\
&
(iii)\quad\lim_{h\to\infty}\Fc(u_h,\Omega)=\Fc(u,\Omega).
\end{split}\end{equation*}
\end{prop}

\begin{proof}


By Lemma \ref{CH:2:stuff_forml1} we know that for every $h\in\mathbb N$ we can find $\delta_h$ small enough such that
\begin{equation}\label{CH:2:forml_eqtn}
\|u-u\varphi_{\delta_h}\|_{L^1(\R^n)}<2^{-h}\quad\textrm{and}
\quad\big|\Fc(u,\Omega)-\Fc(u\varphi_{\delta_h},\Omega)\big|<2^{-h}.
\end{equation}
We can assume that $\delta_h\searrow0$.

By point $(i)$ of Lemma \ref{CH:2:dens_lemma} we know that for every $h$ we can find $\eps_h$ small enough such that
\begin{equation}\label{CH:2:more_formul2}
\|(u\varphi_{\delta_h})\ast\eta_{\eps_h}-u\varphi_{\delta_h}\|_{L^1(B_h)}<2^{-h}
\end{equation}
and
\begin{equation}\label{CH:2:more_formul1}
\big|\Fc(u\varphi_{\delta_h},\Omega_{-\delta_h/2})-\Fc((u\varphi_{\delta_h})\ast\eta_{\eps_h},\Omega_{-\delta_h/2})\big|<2^{-h}.
\end{equation}
Taking $\eps_h$ small enough, we can also assume that
\begin{equation}\label{CH:2:more_formul3}
(u\varphi_{\delta_h})\ast\eta_{\eps_h}\equiv0\qquad\textrm{in }\{|\bar{d}_\Omega|<\delta_h/2\},
\end{equation}
since the $\eps$-convolution enlarges the support at most to an $\eps$-neighborhood of the original support.

Let $u_h:=(u\varphi_{\delta_h})\ast\eta_{\eps_h}$. Since we are taking the $\eps_h$-regularization of
the function $u\varphi_{\delta_h}$, which is just
the product of $u$ with a characteristic function, point $(i)$ of our claim is immediate.

By
\eqref{CH:2:more_formul2}
and the first part of \eqref{CH:2:forml_eqtn} we get point $(ii)$.

As for point $(iii)$, exploiting \eqref{CH:2:more_formul3} and Lemma \ref{CH:2:second_lemma_for_func_appro}, we obtain
\begin{equation*}\begin{split}
\big|\Fc(u,\Omega)&-\Fc(u_h,\Omega)\big|\\
&
\leq
\big|\Fc(u,\Omega)-\Fc(u\varphi_{\delta_h},\Omega)\big|+
\big|\Fc(u\varphi_{\delta_h},\Omega)-\Fc(u\varphi_{\delta_h},\Omega_{-\delta_h/2})\big|\\
&
\qquad+\big|\Fc(u\varphi_{\delta_h},\Omega_{-\delta_h/2})-\Fc(u_h,\Omega_{-\delta_h/2})\big|\\
&
\qquad\qquad+\big|\Fc(u_h,\Omega_{-\delta_h/2})-\Fc(u_h,\Omega)\big|\\
&
\leq 2^{-h}+2^sC\|u\|_{L^\infty(\R^n)}\delta_h^{1-s}+2^{-h},
\end{split}\end{equation*}
which goes to 0 as $h\longrightarrow\infty$.
\end{proof}

\subsection{Proofs of Theorem \ref{CH:2:density_smooth_teo} and Theorem \ref{CH:2:appro_in_bded_open}}

Exploiting Lemma \ref{CH:2:dens_lemma} and the coarea formula, we can now prove Theorem \ref{CH:2:density_smooth_teo}.

\begin{proof}[Proof of Theorem \ref{CH:2:density_smooth_teo}]
The ``if part'' is trivial. Indeed, just from point $(i)$ and the lower semicontinuity of the $s$-perimeter we get
\begin{equation*}
\Per_s(E,\Omega')\leq\liminf_{h\to\infty}\Per_s(E_h,\Omega')<\infty,
\end{equation*}
for every $\Omega'\Subset\Omega$.

Now suppose that $E$ has locally finite $s$-perimeter in $\Omega$.\\
The scheme of the proof is similar to that of the classical case (see, e.g., the proof of \cite[Theorem 13.8]{Maggi}).

Given a sequence $\eps_h\searrow0^+$ we consider the $\eps_h$-regularization of $u:=\chi_E$ and define
the sets
\begin{equation*}
E_h^t:=\{u_{\eps_h}>t\}\quad\textrm{with }t\in(0,1).
\end{equation*}
Sard's Theorem guarantees that for a.e. $t\in(0,1)$ the sequence $\{E_h^t\}_h$ is made of open sets with smooth boundary.
We will get our sets $E_h$ by opportunely choosing $t$.

Since $u_{\eps_h}\longrightarrow\chi_E$ in $L^1_{\loc}(\R^n)$, it is readily seen that for a.e. $t\in(0,1)$
\begin{equation*}
E_h^t\xrightarrow{loc}E,
\end{equation*}
and hence the lower semicontinuity of the $s$-perimeter gives
\begin{equation}\label{CH:2:forml1_pf_coarea}
\Per_s(E,\mathcal O)\leq\liminf_{h\to\infty}\Per_s(E^t_h,\mathcal O),
\end{equation}
for every open set $\mathcal O\subseteq\R^n$.

Moreover from \eqref{CH:2:forml1_smoothing} we have
\begin{equation*}
\{0<u_\eps<1\}\subseteq N_\eps(\partial E)\qquad\forall\,\eps>0,
\end{equation*}
and hence, since $\partial E^t_h\subseteq\{u_{\eps_h}=t\}$, we obtain
\begin{equation}\label{CH:2:forml6_pf_coarea}
\partial E_h^t\subseteq N_{\eps_h}(\partial E),
\end{equation}
which will give $(iii)$ once we choose our $t$.

We improve \eqref{CH:2:forml1_pf_coarea} by showing that, if $\Omega'\Subset\Omega$ is a fixed bounded open set,
then for a.e. $t\in(0,1)$ (with the set of exceptional values of $t$ possibly depending on $\Omega'$),
\begin{equation}\label{CH:2:forml2_pf_coarea}
\Per_s(E,\Omega')=\liminf_{h\to\infty}\Per_s(E^t_h,\Omega').
\end{equation}
 By \eqref{CH:2:forml1_pf_coarea} and Fatou's Lemma, we have
\begin{equation}\label{CH:2:forml5_pf_coarea}
\Per_s(E,\Omega')\leq\int_0^1\liminf_{h\to\infty}\Per_s(E_h^t,\Omega')\,dt\leq\liminf_{h\to\infty}\int_0^1 \Per_s(E_h^t,\Omega')\,dt.
\end{equation}
Let $\mathcal O$ be a bounded open set such that $\Omega'\Subset\mathcal O\Subset\Omega$.
Since $E$ has locally finite $s$-perimeter in $\Omega$, we have $\Per_s(E,\mathcal O)<\infty$.
Then, since $\Omega'\Subset\mathcal O$, point $(i)$ of Lemma \ref{CH:2:dens_lemma}
(with $\mathcal O$ in the place of $\Omega$) implies
\begin{equation}\label{CH:2:forml3_pf_coarea}
\lim_{h\to\infty}\Fc(u_{\eps_h},\Omega')=\Fc(\chi_E,\Omega')=\Per_s(E,\Omega').
\end{equation}
Since $0\leq u_{\eps_h}\leq1$, we have $E^t_h=\R^n$ if $t<0$ and $E^t_h=\emptyset$ if $t>1$, and hence
rewriting \eqref{CH:2:forml3_pf_coarea} exploiting the coarea formula,
\begin{equation*}
\lim_{h\to\infty}\int_0^1\Per_s(E_h^t,\Omega')\,dt=\Per_s(E,\Omega').
\end{equation*}
This and \eqref{CH:2:forml5_pf_coarea} give
\begin{equation*}
\int_0^1\liminf_{h\to\infty}\Per_s(E_h^t,\Omega')\,dt=\Per_s(E,\Omega')=\int_0^1 \Per_s(E,\Omega')\,dt,
\end{equation*}
which implies
\begin{equation}\label{CH:2:end_approx_proof_eq}
\Per_s(E,\Omega')=\liminf_{h\to\infty}\Per_s(E_h^t,\Omega'),\quad\textrm{for a.e. }t\in(0,1),
\end{equation}
as claimed.

Now let the sets $\Omega_k\Subset\Omega$ be as in Corollary \ref{CH:2:regular_approx_open_sets_coroll}.
From \eqref{CH:2:end_approx_proof_eq} we deduce that for a.e. $t\in(0,1)$ we have
\begin{equation}\label{CH:2:end_approx_proof_eq2}
\Per_s(E,\Omega_k)=\liminf_{h\to\infty}\Per_s(E_h^t,\Omega_k),\qquad\forall\,k\in\mathbb N.
\end{equation}

Therefore, combining all we wrote so far, we find that for a.e. $t\in(0,1)$ the sequence $\{E_h^t\}_h$ is made of open sets with smooth boundary such that $E_h^t\xrightarrow{loc}E$ and both
\eqref{CH:2:forml6_pf_coarea} and \eqref{CH:2:end_approx_proof_eq2} hold true.

To conclude, by a diagonal argument we can find $t_0\in(0,1)$ and $h_i\nearrow\infty$ such that,
if we define $E_i:=E^{t_0}_{h_i}$, then
$\{E_i\}$ is a sequence of open sets with smooth boundary
such that $E_i\xrightarrow{loc}E$, with $\partial E_i\subseteq N_{\eps_{h_i}}(\partial E)$, and
\begin{equation}\label{CH:2:forml7_pf_coarea}
\Per_s(E,\Omega_k)=\lim_{i\to\infty}\Per_s(E_i,\Omega_k),\qquad\forall\,k\in\mathbb N.
\end{equation}

Now notice that if $\Omega'\Subset\Omega$, then there exists a $k$ such that $\Omega'\Subset\Omega_k$.
Therefore by \eqref{CH:2:forml7_pf_coarea} and Proposition \ref{CH:2:subcont_lem_approx} we get $(ii)$.

This concludes the proof of the first part of the claim.

\smallskip

Now suppose that $\Omega=\R^n$ and $|E|,\,\Per_s(E)<\infty$.

Since $|E|<\infty$, we know that $u_\eps\longrightarrow\chi_E$ in $L^1(\R^n)$. Therefore we obtain
$E_h^t\longrightarrow E$ for a.e. $t\in(0,1)$.\\
Moreover, from point $(ii)$ of Lemma \ref{CH:2:dens_lemma} we know that
\begin{equation*}
\Fc(u,\R^n)<\infty\qquad\Longrightarrow\qquad\lim_{\eps\to0}\Fc(u_\eps,\R^n)=\Fc(u,\R^n).
\end{equation*}
We can thus repeat the proof above and obtain
\begin{equation*}
\Per_s(E)=\liminf_{h\to\infty}\Per_s(E_h^t),
\end{equation*}
for a.e. $t\in(0,1)$. For any fixed ``good'' $t_0\in(0,1)$ this directly implies, with no need of a diagonal argument,
the existence of a subsequence $h_i\nearrow\infty$ such that
\begin{equation*}
\Per_s(E)=\lim_{i\to\infty}\Per_s(E_{h_i}^{t_0}).
\end{equation*}
We are left to show that in this case we can take the sets $E_h$ to be bounded.

To this end, it is enough to replace the functions $u_{\eps_k}$ with the functions $u_k$
obtained in point $(iii)$ of Lemma \ref{CH:2:dens_lemma}.\\
Indeed, since $u_k$ has compact support, for each $t\in(0,1)$ the set
\begin{equation*}
E_k^t:=\{u_k>t\}
\end{equation*}
is bounded. Since $u_k\longrightarrow u$ in $L^1(\R^n)$ we still find
\begin{equation*}
E_k^t\xrightarrow{loc}E\quad\textrm{for a.e. }t\in(0,1),
\end{equation*}
and, since $0\leq u_k\leq1$ and
\begin{equation*}
\lim_{k\to\infty}\Fc(u_k,\R^n)=\Per_s(E),
\end{equation*}
we can use again the coarea formula to conclude as above.
\end{proof}



\begin{proof}[Proof of Theorem \ref{CH:2:appro_in_bded_open}]
Exploiting the approximating sequence obtained in Proposition \ref{CH:2:density_bded_reg_set}, we can now prove Theorem
 \ref{CH:2:appro_in_bded_open} exactly as above.
 
 As for point $(iii)$, recall that the functions $u_h$ of Proposition \ref{CH:2:density_bded_reg_set} are defined as
 \[u_h=(\chi_E\varphi_{\delta_h})\ast\eta_{\eps_h}.\]
 Notice that, since we can suppose that $\eps_h<\delta_h/2$, we have
 \[u_h=\chi_E\ast\eta_{\eps_h},\qquad\textrm{in }\R^n\setminus N_{2\delta_h}(\partial\Omega).\]
 Therefore, for every $t\in(0,1)$ we find
 \[\partial\{u_h>t\}\subseteq N_{\eps_h}(\partial E)\subseteq N_{2\delta_h}(\partial E),\qquad
 \textrm{in }\R^n\setminus N_{2\delta_h}(\partial\Omega).\]
 
 This gives point $(iii)$ once we choose an appropriate $t$, as in the proof of Theorem \ref{CH:2:density_smooth_teo}.
\end{proof}

\begin{remark}\label{CH:2:rmk_conv_every_cpt_subopen}
We remark that by Proposition \ref{CH:2:subcont_lem_approx} we have also
\begin{equation*}
\lim_{h\to\infty}\Per_s(E_h,\Omega')=\Per_s(E,\Omega'),\qquad\textrm{for every }\Omega'\Subset\Omega.
\end{equation*}
\end{remark}

\section{Existence and compactness of (locally) $s$-minimal sets}

\subsection{Proof of Theorem \ref{CH:2:confront_min_teo}}

\begin{proof}[Proof of Theorem \ref{CH:2:confront_min_teo}]
$(i)\Longrightarrow(ii)\quad$ is obvious.

$(ii)\Longrightarrow(iii)\quad$ Let $\Omega'\Subset\Omega$ and let $F\subseteq\R^n$ be such that $F\setminus\Omega'=E\setminus\Omega'$.\\
Since $E\Delta F\subseteq\Omega'\Subset\Omega$, we have
\begin{equation*}
\Per_s(E,\Omega)\leq \Per_s(F,\Omega).
\end{equation*}
Then, since $F\setminus\Omega'=E\setminus\Omega'$, by Proposition \ref{CH:2:subopensets} we get
\begin{equation*}
\Per_s(E,\Omega')\leq \Per_s(F,\Omega').
\end{equation*}

$(iii)\Longrightarrow(i)\quad$ Let $E$ be locally $s$-minimal in $\Omega$.

First of all we prove that $\Per_s(E,\Omega)<\infty$.\\
Indeed, since $E$ is locally $s$-minimal in $\Omega$, in particular it is $s$-minimal in every $\Omega_r$, with $r\in(-r_0,0)$. Thus,
by minimality and \eqref{CH:2:unif_bound_lip_frac_per}, we get
\begin{equation*}
\Per_s(E,\Omega_r)\leq \Per_s(E\setminus\Omega_r,\Omega_r)\leq \Per_s(\Omega_r)\leq M<\infty,
\end{equation*}
for every $r\in(-r_0,0)$.
Therefore 
by \eqref{CH:2:limit_sub_open} we obtain $\Per_s(E,\Omega)\leq M$.

Now let $F\subseteq\R^n$ be such that $F\setminus\Omega=E\setminus\Omega$. Take a sequence $\{r_k\}\subseteq(-r_0,0)$ such that $r_k\nearrow0$,
let $\Omega_k:=\Omega_{r_k}$, and define
\begin{equation*}
F_k:=(F\cap\Omega_k)\cu\Per(E\setminus\Omega_k).
\end{equation*}
The local minimality of $E$ gives
\begin{equation*}
\Per_s(E,\Omega_k)\leq \Per_s(F_k,\Omega_k),\qquad\textrm{for every }k\in\mathbb N,
\end{equation*}
and by \eqref{CH:2:limit_sub_open} we know that
\begin{equation*}
\Per_s(E,\Omega)=\lim_{k\to\infty}\Per_s(E,\Omega_k).
\end{equation*}
Since $F_k=F$ outside $\Omega\setminus\Omega_k$, and $F_k=E$ in $\Omega\setminus\Omega_k$, we obtain
\begin{equation*}\begin{split}
\Per_s(F,\Omega_k)&-\Per_s(F_k,\Omega_k)=\Ll_s(F\cap\Omega_k,\Co F\cap(\Omega\setminus\Omega_k))\\
&
+\Ll_s(\Co F\cap\Omega_k,F\cap(\Omega\setminus\Omega_k))
-\Ll_s(F\cap\Omega_k,\Co E\cap(\Omega\setminus\Omega_k))\\
&\qquad\quad -\Ll_s(\Co F\cap\Omega_k,E\cap(\Omega\setminus\Omega_k)).
\end{split}
\end{equation*}
Notice that each of the four terms in the right hand side is less or equal than $\Ll_s(\Omega_k,\Omega\setminus\Omega_k)$.
Thus
\begin{equation*}
a_k:=|\Per_s(F,\Omega_k)-\Per_s(F_k,\Omega_k)|\leq4\,\Ll_s(\Omega_k,\Omega\setminus\Omega_k).
\end{equation*}
Notice that from point $(i)$ of \eqref{CH:2:uniform_bound_strips} we have $a_k\longrightarrow0$.

Now
\begin{equation*}
\Per_s(F,\Omega)+a_k\geq \Per_s(F,\Omega_k)+a_k\geq \Per_s(F_k,\Omega_k)\geq \Per_s(E,\Omega_k),
\end{equation*}
and hence, passing to the limit $k\to\infty$, we get
\begin{equation*}
\Per_s(F,\Omega)\geq \Per_s(E,\Omega).
\end{equation*}
Since $F$ was an arbitrary competitor for $E$, we see that $E$ is $s$-minimal in $\Omega$.
\end{proof}

\subsection{Proofs of Theorem \ref{CH:2:minimal_comp} and Corollary \ref{CH:2:local_minima_comp}}\label{CH:2:Min_Compactness_Section}

We slightly modify the proof of \cite[Theorem 3.3]{CRS10} to show that the conclusion remains true in any bounded open set $\Omega$ with Lipschitz boundary.

\begin{proof}[Proof of Theorem \ref{CH:2:minimal_comp}]

Assume $F=E$ outside $\Omega$ and let
\begin{equation*}
F_k:=(F\cap\Omega)\cu\Per(E_k\setminus\Omega).
\end{equation*}
Since $F_k=E_k$ outside $\Omega$ and $E_k$ is $s$-minimal in $\Omega$, we have
\begin{equation*}
\Per_s(F_k,\Omega)\geq \Per_s(E_k,\Omega).
\end{equation*}
On the other hand, since $F_k=F$ inside $\Omega$, we have
\begin{equation*}
|\Per_s(F_k,\Omega)-\Per_s(F,\Omega)|\leq\Ll_s(\Omega,(F_k\Delta F)\setminus\Omega)=
\Ll_s(\Omega,(E_k\Delta E)\setminus\Omega)=:b_k.
\end{equation*}
Thus
\begin{equation*}
\Per_s(F,\Omega)+b_k\geq \Per_s(F_k,\Omega)\geq \Per_s(E_k,\Omega).
\end{equation*}
If we prove that $b_k\longrightarrow0$, then by lower semicontinuty of the fractional perimeter
\begin{equation}\label{CH:2:inequality1}
\Per_s(F,\Omega)\geq\limsup_{k\to\infty}\Per_s(E_k,\Omega)\geq\liminf_{k\to\infty}\Per_s(E_k,\Omega)\geq \Per_s(E,\Omega).
\end{equation}
This shows that $E$ is $s$-minimal in $\Omega$.
Moreover, \eqref{CH:2:conv_perimeter} follows from \eqref{CH:2:inequality1} by taking $F=E$.

We are left to show $b_k\longrightarrow0$.\\
Let $r_0$ be as in Proposition \ref{CH:2:bound_perimeter_unif} and let $R>r_0$. In the end we will let $R\longrightarrow\infty$.
Define
\begin{equation*}
a_k(r):=\Ha^{n-1}\big((E_k\Delta E)\cap\{\bar{d}_\Omega=r\})\big)
\end{equation*}
for every $r\in[0,r_0)$.\\
We split $b_k$ as the sum
\begin{equation*}\begin{split}
b_k&=\Ll_s\big(\Omega,(E_k\Delta E)\cap (\Omega_{r_0}\setminus\Omega)\big)
+\Ll_s\big(\Omega,(E_k\Delta E)\cap (\Omega_R\setminus\Omega_{r_0})\big)\\
&
\qquad\qquad\qquad+\Ll_s\big(\Omega,(E_k\Delta E)\setminus\Omega_R\big).
\end{split}\end{equation*}
Notice that if $x\in\Omega$ and $y\in(\Omega_R\setminus\Omega_{r_0})$, then 
$|x-y|\geq r_0$, and hence
\begin{equation*}\begin{split}
\Ll_s\big(\Omega,(E_k\Delta E)\cap (\Omega_R\setminus\Omega_{r_0})\big)&
=\int_{\Omega_R\setminus\Omega_{r_0}}\chi_{E_k\Delta E}(y)\,dy\int_\Omega\frac{1}{|x-y|^{n+s}}dx\\
&
\leq \frac{|\Omega|}{r_0^{n+s}}|(E_k\Delta E)\cap (\Omega_R\setminus\Omega_{r_0})|.
\end{split}\end{equation*}
Since $E_k\xrightarrow{loc}E$ and $\Omega_R\setminus\Omega_{r_0}$ is bounded, for every fixed $R$ we find
\begin{equation*}
\lim_{k\to\infty}\Ll_s\big(\Omega,(E_k\Delta E)\cap (\Omega_R\setminus\Omega_{r_0})\big)=0.
\end{equation*}
As for the last term, we have
\begin{equation*}
\Ll_s\big(\Omega,(E_k\Delta E)\setminus\Omega_R\big)\leq\Ll_s(\Omega,\Co\Omega_R)\leq
\int_\Omega dx\int_{\Co B_R(x)}\frac{dy}{|x-y|^{n+s}}=\frac{n\omega_n}{s\,R^s}|\Omega|.
\end{equation*}
We are left to estimate the first term. By using the coarea formula, we obtain
\begin{equation*}\begin{split}
\Ll_s\big(\Omega,(E_k&\Delta E)\cap (\Omega_{r_0}\setminus\Omega)\big)\\
&
=\int_0^{r_0}\Big(\int_{\{\bar{d}_\Omega=r\}}\chi_{E_k\Delta E}(y)\Big(\int_\Omega\frac{dx}{|x-y|^{n+s}}\Big)d\Ha^{n-1}_y\Big)dr\\
&
\leq
\int_0^{r_0}\Big(\int_{\{\bar{d}_\Omega=r\}}\chi_{E_k\Delta E}(y)\Big(\int_{\Co B_r(y)}\frac{dx}{|x-y|^{n+s}}\Big)d\Ha^{n-1}_y\Big)dr\\
&
=\frac{n\omega_n}{s}\int_0^{r_0}\frac{a_k(r)}{r^s}\,dr.
\end{split}
\end{equation*}
Notice that
\begin{equation*}
\int_0^{r_0}a_k(r)\,dr=|(E_k\Delta E)\cap(\Omega_{r_0}\setminus\Omega)|\xrightarrow{k\to\infty}0,
\end{equation*}
so that
\begin{equation*}
a_k(r)\xrightarrow{k\to\infty}0\qquad\textrm{for a.e. }r\in[0,r_0).
\end{equation*}
Moreover, exploiting \eqref{CH:2:bound_perimeter_unif_eq} we get
\begin{equation*}
\int_0^{r_0}\frac{a_k(r)}{r^s}\,dr\leq M\int_0^{r_0}\frac{1}{r^s}\,dr=\frac{M}{1-s}r_0^{1-s},
\end{equation*}
and hence, by dominated convergence, we obtain
\begin{equation*}
\lim_{k\to\infty}\int_0^{r_0}\frac{a_k(r)}{r^s}\,dr=0.
\end{equation*}
Therefore
\begin{equation*}
\limsup_{k\to\infty}b_k\leq\frac{n\omega_n}{s}|\Omega|\,R^{-s}.
\end{equation*}
Letting $R\longrightarrow\infty$, we obtain $b_k\longrightarrow0$, concluding the proof.
\end{proof}

\begin{proof}[Proof of Corollary \ref{CH:2:local_minima_comp}]
Let the sets $\Omega_k\Subset\Omega$ be as in Corollary
\ref{CH:2:regular_approx_open_sets_coroll}.
By Theorem \ref{CH:2:minimal_comp} we see that $E$ is $s$-minimal in each $\Omega_k$. Moreover \eqref{CH:2:conv_perimeter} gives
\begin{equation*}
\Per_s(E,\Omega_k)=\lim_{h\to\infty}\Per_s(E_h,\Omega_k),
\end{equation*}
for every $k$. Now if $\Omega'\Subset\Omega$, then $\Omega'\subseteq\Omega_k$ for some $k$.
Thus $E$ is $s$-minimal in $\Omega'$
and we obtain \eqref{CH:2:conv_perimeter_locally} by Proposition \ref{CH:2:subcont_lem_approx}.
\end{proof}

\subsection{Proofs of Theorem \ref{CH:2:glob_min_exist} and Corollary \ref{CH:2:loc_min_set_cor}}\label{CH:2:Existence_Section}


We can exploit Proposition \ref{CH:2:compact_prop} to extend the existence result \cite[Theorem 3.2]{CRS10}
to any open set $\Omega$, provided a competitor with finite fractional perimeter exists.

\begin{proof}[Proof of Theorem \ref{CH:2:glob_min_exist}]
The ``only if'' part is trivial. Now suppose there exists a competitor for $E_0$ with finite $s$-perimeter in $\Omega$. Then
\begin{equation*}
\inf\{\Per_s(E,\Omega)\,|\,E\setminus\Omega=E_0\setminus\Omega\}<\infty
\end{equation*}
and we can find a minimizing sequence, that is $\{E_h\}$ with $E_h\setminus\Omega=E_0\setminus\Omega$ and
\begin{equation*}
\lim_{h\to\infty}\Per_s(E_h,\Omega)=\inf\{\Per_s(E,\Omega)\,|\,E\setminus\Omega=E_0\setminus\Omega\}.
\end{equation*}
Let $\Omega'\Subset\Omega$. Since, for every $h\in\mathbb N$ we have
\begin{equation*}
\Per_s(E_h,\Omega')\leq \Per_s(E_h,\Omega)\leq M<\infty,
\end{equation*}
we can use Proposition \ref{CH:2:compact_prop} to find a set $E'\subseteq\Omega$ such that
\begin{equation*}
E_h\cap\Omega\xrightarrow{loc}E'
\end{equation*}
(up to subsequence). Since $E_h\setminus\Omega=E_0\setminus\Omega$ for every $h$, if we set $E:=E'\cu\Per(E_0\setminus\Omega)$, then
\begin{equation*}
E_h\xrightarrow{loc}E.
\end{equation*}
The semicontinuity of the fractional perimeter concludes the proof.
\end{proof}

\begin{remark}\label{CH:2:rmk_crs_existence}
In particular, if $\Omega$ is a bounded open set with Lipschitz boundary, then (as already proved in \cite{CRS10}) we can always find an $s$-minimal set for every $s\in(0,1)$,
no matter what the external data $E_0\setminus\Omega$ is. Indeed in this case
\begin{equation*}
\Per_s(E_0\setminus\Omega,\Omega)\leq \Per_s(\Omega)<\infty.
\end{equation*}
Actually, in order to have the existence of $s$-minimal sets for some fixed $s\in(0,1)$, the open set $\Omega$ need not be bounded nor have a regular boundary. It is enough to have
\[\Per_s(\Omega)<\infty.\]
Then $E_0\setminus\Omega$ has finite $s$-perimeter in $\Omega$ and we can apply Theorem
\ref{CH:2:glob_min_exist}.
\end{remark}


Now we prove that a locally $s$-minimal set always exists, without having to assume the existence of a competitor having finite
fractional perimeter.

\begin{proof}[Proof of Corollary \ref{CH:2:loc_min_set_cor}]
Let the sets $\Omega_k$ be as in Corollary
\ref{CH:2:regular_approx_open_sets_coroll}.\\
From Theorem \ref{CH:2:glob_min_exist} and Remark \ref{CH:2:rmk_crs_existence} we know that for every $k$ we can find a set $E_k$ which is $s$-minimal in $\Omega_k$
and such that $E_k\setminus\Omega_k=E_0\setminus\Omega_k$.\\
Notice that, since the sequence $\Omega_k$ is increasing, the set $E_h$ is $s$-minimal in $\Omega_k$ for every $h\geq k$.\\
This gives us a sequence $\{E_h\}$ satisfying the hypothesis of Proposition \ref{CH:2:compact_prop} (see Remark\ref{CH:2:min_app_seq_rmk}),
and hence (up to a subsequence)
\begin{equation*}
E_h\cap\Omega\xrightarrow{loc}F,
\end{equation*}
for some $F\subseteq\Omega$. Since $E_h\setminus\Omega=E_0\setminus\Omega$ for every $h$,
if we set $E:=F\cu\Per(E_0\setminus\Omega)$, we obtain
\begin{equation*}
E_h\xrightarrow{loc}E.
\end{equation*}
Theorem \ref{CH:2:minimal_comp} guarantees that $E$ is $s$-minimal in every $\Omega_k$
and hence also locally $s$-minimal in $\Omega$. Indeed, if $\Omega'\Subset\Omega$,
then for some $k$ big enough we have $\Omega'\subseteq\Omega_k$. Now, since $E$ is $s$-minimal in $\Omega_k$, it is
$s$-minimal also in $\Omega'$.
\end{proof}


\section{Locally $s$-minimal sets in cylinders}\label{CH:2:Min_Sets_Cyl_Section}

Given a bounded open set $\Omega\subseteq\R^n$, 
we consider
the cylinders
\begin{equation*}
\Omega^k:=\Omega\times(-k,k),\qquad\Omega^\infty:=\Omega\times\R.
\end{equation*}
We recall that, given any set $E_0\subseteq\R^{n+1}$, by Corollary \ref{CH:2:loc_min_set_cor} we can find a set $E\subseteq\R^{n+1}$ which is locally $s$-minimal in
$\Omega^\infty$
 and such that $E\setminus\Omega^\infty=E_0\setminus\Omega^\infty$.

\begin{remark}\label{CH:2:rmk_from_compact_to_any_subset}
Actually, if $\Omega$ has Lipschitz boundary then $E$ is $s$-minimal in every cylinder $\mathcal O=\Omega\times(a,b)$ of finite height (notice that $\mathcal O$
is not compactly contained in $\Omega^\infty$).
Indeed, $\mathcal O$ is a bounded open set with Lipschitz boundary and $E$ is locally $s$-minimal in $\mathcal O$.
Thus, by Theorem \ref{CH:2:confront_min_teo}, $E$ is $s$-minimal in $\mathcal O$.\\
As a consequence, $E$ is $s$-minimal in every bounded open subset $\Omega'\subseteq\Omega^\infty$.
\end{remark}

We are going to consider as exterior data the subgraph
\begin{equation*}
E_0=\Sg(v):=\{(x,t)\in\R^{n+1}\,|\,t<v(x)\},
\end{equation*}
of a function $v:\R^n\longrightarrow\R$, which is locally bounded, i.e.
\begin{equation}\label{CH:2:locally_bounded_assumption}
M_r:=\sup_{|x|\leq r}|v(x)|<\infty,\qquad\textrm{for every }r>0.
\end{equation}

The following result is an immediate consequence of (the proof of)  \cite[Lemma 3.3]{graph}.

\begin{lemma}\label{CH:2:bded_cyl_prop}
Let $\Omega\subseteq\R^n$ be a bounded open set with $C^{1,1}$ boundary and let $v:\R^n\longrightarrow\R$
be locally bounded. There exists a constant $M=M(n,s,\Omega,v)>0$
such that if $E\subseteq\R^{n+1}$ is locally $s$-minimal in $\Omega^\infty$,
with $E\setminus\Omega^\infty=\Sg(v)\setminus\Omega^\infty$,
then
\begin{equation*}
\Omega\times(-\infty,-M]\subseteq E\cap\Omega^\infty\subseteq\Omega\times(-\infty,M].
\end{equation*}
As a consequence
\begin{equation}\label{CH:2:formu_trivia2}
E\setminus\big(\Omega\times[-M,M]\big)=\Sg(v)\setminus\big(\Omega\times[-M,M]\big).
\end{equation}
\end{lemma}
\begin{proof}
By Remark \ref{CH:2:rmk_from_compact_to_any_subset}, the set $E$ is $s$-minimal in $\Omega^\infty$ in the sense considered
in \cite{graph}.
Thus, \cite[Lemma 3.3]{graph} guarantees that
\[E\cap\Omega^\infty\subseteq\Omega\times(-\infty,M].\]
Moreover, the same argument used in the proof shows also that
\[\Co E\cap\Omega^\infty\subseteq\Omega\times[-M,\infty),\]
(up to considering a bigger $M$).

Since $M>M_{R_0}$, where $R_0$ is such that $\Omega\Subset B_{R_0}$,
we get \eqref{CH:2:formu_trivia2}, concluding the proof.
\end{proof}

Roughly speaking, Lemma \ref{CH:2:bded_cyl_prop} gives an a priori bound on the variation of
$\partial E$ in the ``vertical'' direction.
In particular, from \eqref{CH:2:formu_trivia2} we see that it is enough to look for a locally $s$-minimal set among sets which coincide with
$\Sg(v)$ out of $\Omega\times[-M,M]$.

As a consequence, we can prove that a set is locally $s$-minimal in $\Omega^\infty$
if and only if it is $s$-minimal in $\Omega\times[-M,M]$.

\begin{prop}\label{CH:2:bded_cyl_coroll}
Let $\Omega\subseteq\R^n$ be a bounded open set with $C^{1,1}$ boundary and let $v:\R^n\longrightarrow\R$ be locally bounded. Let $M$ be as in Lemma \ref{CH:2:bded_cyl_prop} and
let $k_0$ be the smallest integer $k_0>M$.
Let
$F\subseteq\R^{n+1}$ be $s$-minimal in $\Omega^{k_0}$,
with respect to the exterior data
\begin{equation}\label{CH:2:ext_data_eq_cyl}
F\setminus\Omega^{k_0}=\Sg(v)\setminus\Omega^{k_0}.
\end{equation}
Then $F$ is $s$-minimal in $\Omega^k$ for every $k\geq k_0$,
hence is locally $s$-minimal in $\Omega^\infty$.
\end{prop}
\begin{proof}
Let $E\subseteq\R^{n+1}$ be locally $s$-minimal in $\Omega^\infty$, with respect to the exterior data 
\[E\setminus\Omega^\infty=\Sg(v)\setminus\Omega^\infty.\]
Recall that by Remark \ref{CH:2:rmk_from_compact_to_any_subset} the set $E$ is $s$-minimal in $\Omega^k$
for every $k$. In particular
\[\Per_s(E,\Omega^k)<\infty\qquad\forall\,k\in\mathbb N.\]

To prove the Proposition, it is enough to show that
\begin{equation}\label{CH:2:bded_to_unbded_cyl}
\Per_s(F,\Omega^k)=\Per_s(E,\Omega^k),\qquad\textrm{for every }k\geq k_0.
\end{equation}
Indeed, notice that by \eqref{CH:2:ext_data_eq_cyl} and \eqref{CH:2:formu_trivia2} we have
\begin{equation}\label{CH:2:triv_fin1}
F\setminus\Omega^{k_0}=\Sg(v)\setminus\Omega^{k_0}=E\setminus\Omega^{k_0},
\end{equation}
hence, clearly,
\[F\setminus\Omega^k=E\setminus\Omega^k,\qquad\forall\,k\geq k_0.\]
Then, since $E$ is $s$-minimal in $\Omega^k$, from \eqref{CH:2:bded_to_unbded_cyl}
we conclude that also $F$ is $s$-minimal in $\Omega^k$, for every $k\geq k_0$. In turn, this implies that $F$
is locally $s$-minimal in $\Omega^\infty$.

Exploiting
Proposition \ref{CH:2:subopensets}, by \eqref{CH:2:triv_fin1} we obtain that for every $k\geq k_0$
\begin{equation}\label{CH:2:fin_triv2}
\Per_s(F,\Omega^k)=\Per_s(F,\Omega^{k_0})+c_k,\qquad \Per_s(E,\Omega^k)=\Per_s(E,\Omega^{k_0})+c_k,
\end{equation}
where
\begin{equation*}\begin{split}
c_k=\Ll_s\big(\Sg(v)&\cap(\Omega^k\setminus\Omega^{k_0}),\Co\Sg(v)\setminus\Omega^k\big)
+
\Ll_s\big(\Sg(v)\setminus\Omega^k,\Co\Sg(v)\cap(\Omega^k\setminus\Omega^{k_0})\big)\\
&
\qquad
+\Ll_s\big(\Sg(v)\cap(\Omega^k\setminus\Omega^{k_0}),\Co\Sg(v)\cap(\Omega^k\setminus\Omega^{k_0})\big),
\end{split}
\end{equation*}
which is finite and does not depend on $E$ nor $F$. To see that $c_k$ is finite, simply notice that
\begin{equation*}
c_k\leq \Per_s(E,\Omega^k)<\infty.
\end{equation*}
Now, by \eqref{CH:2:triv_fin1} and the minimality of $F$ we have
\begin{equation*}
\Per_s(F,\Omega^{k_0})\leq \Per_s(E,\Omega^{k_0}).
\end{equation*}
On the other hand, since also the set $E$ is $s$-minimal in $\Omega^{k_0}$,
again by \eqref{CH:2:triv_fin1} we get
\begin{equation*}
\Per_s(E,\Omega^{k_0})\leq \Per_s(F,\Omega^{k_0}).
\end{equation*}
This and \eqref{CH:2:fin_triv2} give
\begin{equation*}
\Per_s(F,\Omega^k)=\Per_s(F,\Omega^{k_0})+c_k=\Per_s(E,\Omega^k),
\end{equation*}
proving \eqref{CH:2:bded_to_unbded_cyl} and concluding the proof.
\end{proof}

It is now natural to wonder whether the set $F$ is actually $s$-minimal in $\Omega^\infty$. The answer, in general, is no.
Indeed, Theorem \ref{CH:2:bound_unbound_per_cyl_prop} shows
that in general we cannot hope to find an $s$-minimal set in $\Omega^\infty$.

\begin{proof}[Proof of Theorem \ref{CH:2:bound_unbound_per_cyl_prop}]

Notice that by \eqref{CH:2:bound_hp_forml_subgraph} we have
\begin{equation*}\begin{split}
&E\cap(\Omega^\infty\setminus\Omega^{k+1})=\Omega\times(-\infty,-k-1),\\
&
\Co E\cap(\Omega^\infty\setminus\Omega^{k+1})=\Omega\times(k+1,\infty),
\end{split}\end{equation*}
and
\begin{equation*}
E\cap\Omega^{k+1}\subseteq\Omega\times(-k-1,k),\qquad\qquad\Co E\cap\Omega^{k+1}
\subseteq\Omega\times(-k,k+1).
\end{equation*}
Thus
\begin{equation*}\begin{split}
\Per_s^L(E,\Omega^\infty)&=
\Per_s^L(E,\Omega^{k+1})+\Ll_s(E\cap(\Omega^\infty\setminus\Omega^{k+1}),\Co E\cap\Omega^{k+1})\\
&\qquad+\Ll_s(\Co E\cap(\Omega^\infty\setminus\Omega^{k+1}),E\cap\Omega^{k+1})
+\Per_s^L(E,\Omega^\infty\setminus\Omega^{k+1})\\
&
\leq \Per_s^L(E,\Omega^{k+1})+2\Ll_s(\Omega\times(-\infty,-k-1),\Omega\times(-k,k+1))\\
&\qquad+\Ll_s(\Omega\times(-\infty,-k-1),\Omega\times(k+1,\infty)).
\end{split}
\end{equation*}
Since $d(\Omega\times(-\infty,-k-1),\Omega\times(-k,k+1))=1$, we get
\begin{equation*}\begin{split}
\Ll_s(\Omega\times(-\infty&,-k-1),\Omega\times(-k,k+1))\\
&
\leq\int_{\Omega\times(-k,k+1)}\Big(\int_{\Co B_1(X)}\frac{dY}{|X-Y|^{n+1+s}}\Big)\,dX\\
&
=\frac{(n+1)\omega_{n+1}}{s}(2k+1)|\Omega|.
\end{split}\end{equation*}
As for the last term, since $n+1\geq2$, we have
\begin{equation*}\begin{split}
\Ll_s(\Omega\times(-\infty&,-k-1),\Omega\times(k+1,\infty))\\
&
=\int_\Omega\int_\Omega\Big(\int_{-\infty}^{-k-1}\int_{k+1}^\infty\frac{dt\,d\tau}{(|x-y|^2+(t-\tau)^2)^\frac{n+1+s}{2}}
\Big)dx\,dy\\
&
\leq|\Omega|^2\int_{-\infty}^{-k-1}\Big(\int_{k+1}^\infty\frac{dt}{(t-\tau)^{n+1+s}}\Big)d\tau\\
&
=\frac{|\Omega|^2}{n+s}\int_{-\infty}^{-k-1}\frac{d\tau}{(k+1-\tau)^{n+s}}\\
&
=\frac{|\Omega|^2}{(n+s)(n-1+s)}\,\frac{1}{(2k+2)^{n-1+s}}.
\end{split}\end{equation*}

This shows that $\Per_s^L(E,\Omega^\infty)<\infty$.\\
Now suppose that $E\subseteq\R^{n+1}$ satisfies \eqref{CH:2:bound_hp_forml_subgraph2}. Then
\begin{equation*}
\Per_s^{NL}(E,\Omega^\infty)\geq2\Ll_s(\Omega\times(-\infty,-k),\Co\Omega\times(k,\infty)).
\end{equation*}
Since $\Omega$ is bounded, we can take $R>0$ big enough such that $\Omega\Subset B_R$. For every
$T>T_0:=\max\{k,R\}$ we have
\begin{equation*}
\Omega\times(-\infty,-T)\subseteq\Omega\times(-\infty,-k)\quad\textrm{and}\quad (B_T\setminus B_R)\times(T,\infty)
\subseteq\Co\Omega\times(k,\infty).
\end{equation*}
Thus for every $T>T_0$
\begin{equation*}\begin{split}
\Ll_s(\Omega&\times(-\infty,-k),\Co\Omega\times(k,\infty))
\geq\Ll_s(\Omega\times(-\infty,-T),(B_T\setminus B_R)\times(T,\infty))\\
&
=\int_\Omega dx\int_{B_T\setminus B_R}dy\int_{-\infty}^{-T}dt
\int_T^\infty\frac{d\tau}{(|x-y|^2+(\tau-t)^2)^\frac{n+1+s}{2}}=:a_T.
\end{split}\end{equation*}
Notice that for every $x\in\Omega,\,y\in B_T\setminus B_R,\,t\in(-\infty,-T)$ and $\tau\in(T,\infty)$, we have
\begin{equation*}
|x-y|\leq|x|+|y|\leq R+T\leq 2T\leq \tau-t,
\end{equation*}
and hence
\begin{equation*}\begin{split}
a_T&\geq\frac{1}{2^\frac{n+1+s}{2}}
\int_\Omega dx\int_{B_T\setminus B_R}dy\int_{-\infty}^{-T}dt
\int_T^\infty\frac{d\tau}{(\tau-t)^{n+1+s}}\\
&
=\frac{|\Omega|}{2^\frac{n+1+s}{2}(n+s)(n-1+s)}\frac{|B_T\setminus B_R|}{(2T)^{n-1+s}}.
\end{split}\end{equation*}
Since $|B_T\setminus B_R|\sim T^n$ as $T\to\infty$, we get $a_T\longrightarrow\infty$. Therefore,
since
\begin{equation*}
\Per^{NL}_s(E,\Omega^\infty)\geq 2a_T\qquad\textrm{for every }T>T_0,
\end{equation*}
we obtain $\Per^{NL}_s(E,\Omega^\infty)=\infty$.

To conclude, let $\Omega$ be bounded, with $C^{1,1}$ boundary, and let $v\in L^\infty(\R^n)$.\\
Suppose that there exists a set $E\subseteq\R^{n+1}$ which is $s$-minimal in $\Omega^\infty$
with respect to the exterior data $E\setminus\Omega^\infty=\Sg(v)\setminus\Omega^\infty$.\\
Then, thanks to Lemma \ref{CH:2:bded_cyl_prop}, we can find $k$ big enough such that $E$ satisfies
\eqref{CH:2:bound_hp_forml_subgraph2}. 
Since this implies $\Per_s(E,\Omega^\infty)=\infty$, we reach a contradiction concluding the proof.
\end{proof}

\begin{corollary}\label{CH:2:non_well_def_frac_area}
In particular
\begin{equation}\label{CH:2:forml_eq1}
u\in BV_{\loc}(\R^n)\cap L^\infty_{\loc}(\R^n)\quad\Longrightarrow\quad \Per_s^L(\Sg(u),\Omega^\infty)<\infty,
\end{equation}
and
\begin{equation}\label{CH:2:forml_eq2}
u\in
L^\infty(\R^n)\quad\Longrightarrow\quad \Per_s^{NL}(\Sg(u),\Omega^\infty)=\infty,
\end{equation}
for every bounded open set $\Omega\subseteq\R^n$.

Furthermore, if $|u|\leq M$ in $\Omega$ and there exists $\Sigma\subseteq\mathbb S^{n-1}$ with $\Ha^{n-1}(\Sigma)>0$
such that either
\begin{equation*}
u(r\omega)\leq M\quad\textrm{or}\quad u(r\omega)\geq-M\qquad\textrm{for every }\omega\in\Sigma\quad\textrm{and}\quad r\geq r_0,
\end{equation*}
then $\Per^{NL}_s(\Sg(u),\Omega^\infty)=\infty$.
\end{corollary}

\begin{proof}
Both \eqref{CH:2:forml_eq1} and \eqref{CH:2:forml_eq2} are immediate from Theorem \ref{CH:2:bound_unbound_per_cyl_prop}, so
we only need to prove the last claim.

Since $\Omega$ is bounded, we can find $R>0$ such that $\Omega\Subset B_R$.\\
For every $T>T_0:=\max\{M,R,r_0\}$ define
\begin{equation*}
\mathcal S(T):=\{x=r\omega\in\R^n\,|\,r\in(T_0,T),\,\omega\in\Sigma\}.
\end{equation*}
Notice that $\mathcal S(T)\subseteq B_T$ and
\begin{equation*}\begin{split}
|\mathcal S(T)|&
=\int_{T_0}^T\Big(\int_{\partial B_r}\chi_{\mathcal S(T)}\,d\Ha^{n-1}\Big)dr
=\int_{T_0}^T\Ha^{n-1}(r\Sigma)\,dr\\
&
=\frac{\Ha^{n-1}(\Sigma)}{n}(T^n-T_0^n).
\end{split}
\end{equation*}
Suppose that $u(r\omega)\leq M$ for every $r\geq r_0$ and $\omega\in\Sigma$. Then, arguing as in the second part of the proof
of Theorem \ref{CH:2:bound_unbound_per_cyl_prop}, we obtain
\begin{equation*}\begin{split}
\Per^{NL}_s(\Sg(u),\Omega^\infty)&\geq\Ll_s(\Sg(u)\cap\Omega^\infty,\Co \Sg(u)\setminus\Omega^\infty)\\
&
\geq\Ll_s(\Omega\times(-\infty,-T),\mathcal S(T)\times(T,\infty))\\
&
\geq\frac{|\Omega|}{2^\frac{n+1+s}{2}(n+s)(n-1+s)}\frac{|\mathcal S(T)|}{(2T)^{n-1+s}},
\end{split}
\end{equation*}
for every $T>T_0$.
Since
\begin{equation*}
\frac{|\mathcal S(T)|}{(2T)^{n-1+s}}\sim T^{1-s},
\end{equation*}
which tends to $\infty$ as $T\to\infty$, we get our claim.
\end{proof}


In the classical framework, the area functional of a function $u\in C^{0,1}(\R^n)$ is defined as
\begin{equation*}
\mathcal A(u,\Omega):=\int_\Omega\sqrt{1+|\nabla u|^2}\,dx=\Ha^n\big(\{(x,u(x))\in\R^{n+1}\,|\,x\in\Omega\}\big),
\end{equation*}
for any bounded open set $\Omega\subseteq\R^n$. Exploiting the subgraph of $u$ one then defines the relaxed area
functional of a function $u\in BV_{\loc}(\R^n)$ as
\begin{equation}\label{CH:2:relaxed_area}
\mathcal A(u,\Omega):=\Per(\Sg(u),\Omega^\infty).
\end{equation}
Notice that when $u$ is Lipschitz the two definitions coincide.

One might then be tempted to define a nonlocal fractional version of the area functional by replacing the 
classical perimeter in \eqref{CH:2:relaxed_area} with the $s$-perimeter, that is
\begin{equation*}
\mathcal A_s(u,\Omega):=\Per_s(\Sg(u),\Omega^\infty).
\end{equation*}
However Corollary \ref{CH:2:non_well_def_frac_area}
shows that this definition is ill-posed even for regular functions $u$.\\
On the other hand, it is worth remarking that one could use just the local part of the $s$-perimeter,
but then the resulting functional
\begin{equation*}
\mathcal A_s^L(u,\Omega):=\Per_s^L(\Sg(u),\Omega^\infty)=\frac{1}{2}[\chi_{\Sg(u)}]_{W^{s,1}(\Omega^\infty)}
\end{equation*}
has a local nature.

Exploiting \cite[Theorem 1]{Davila}, we obtain the following:
\begin{lemma}
Let $\Omega\subseteq\R^n$ be a bounded open set with Lipschitz boundary and let $u\in BV(\Omega)\cap L^\infty(\Omega)$.
Then
\begin{equation*}
\lim_{s\to1^-}(1-s)\mathcal A^L_s(u,\Omega)=\omega_n\mathcal A(u,\Omega).
\end{equation*}
\end{lemma}

\begin{proof}
Let $k$ be such that $|u|\leq k$. Then $E=\Sg(u)$ satisfies \eqref{CH:2:bound_hp_forml_subgraph} and hence, arguing as in the beginning of the proof of Theorem \ref{CH:2:bound_unbound_per_cyl_prop}, we
get
\begin{equation*}
\mathcal A_s^L(u,\Omega)=\Per_s^L(\Sg(u),\Omega^{k+1})+O(1),
\end{equation*}
as $s\to1$.
Since $\Sg(u)$ has finite perimeter in $\Omega^{k+1}$, which is a bounded open set with Lipschitz boundary,
we conclude using \cite[Theorem 1]{Davila} (see also Theorem \ref{CH:1:asymptotics_teo} for the asymptotics as $s\to1$ of the $s$-perimeter).\\
Indeed, notice that since $|u|\leq k$, we have
\begin{equation*}
\Per(\Sg(u),\Omega^{k+1})=\Per(\Sg(u),\Omega^\infty)=\mathcal A(u,\Omega).
\qedhere
\end{equation*}
\end{proof}

\end{chapter}

\begin{chapter}[Complete stickiness of nonlocal minimal surfaces]{Complete stickiness of nonlocal minimal surfaces
for small values of the fractional parameter}\label{Asympto0_CH_label}

\minitoc

%
%

\section{Introduction and main results}



In this chapter, we deal with the behavior of $s$-minimal sets when the fractional parameter $s\in (0,1)$ is small. In particular
\begin{itemize}
\item we give the asymptotic behavior of the fractional mean curvature as $s\to 0^+$, 
\item we classify the behavior of $s$-minimal surfaces, in dependence of the exterior data at infinity.
\end{itemize}  
Moreover, we prove the continuity of the fractional mean curvature in all variables for $s\in [0,1]$.

\smallskip

It is convenient to recall the definition of the $s$-fractional mean curvature of a set $E$ at a point $q\in\partial E$ (which is the fractional counterpart of the classical mean curvature). It is defined as the principal value integral
\[\I_s[E](q):=\PV\int_{\R^n}\frac{\chi_{\Co E}(y)-\chi_E(y)}{|y-q|^{n+s}}\,dy,\]
that is
\[\I_s[E](q):=\lim_{\varrho\to0^+}\I_s^\varrho[E](q),\qquad\textrm{where}\qquad
\I_s^\varrho[E](q):=\int_{
\Co B_\varrho(q)}\frac{\chi_{\Co E}(y)-\chi_E(y)}{|y-q|^{n+s}}\,dy.\] 
For the main properties of the fractional mean curvature, we refer, e.g., to \cite{Abaty}.

Let us also recall here the notation for the area of the $(n-1)$-dimensional sphere as 
 \[ \varpi_n:=\mathcal H^{n-1}\left(\left\{x\in\R^n\,|\,|x|=1\right\}\right) ,\]
 where  $\mathcal H^{n-1}
 $ is the $(n-1)$-dimensional Hausdorff  measure. The volume of the $n$-dimensional unit ball is then
  \[\omega_n=|B_1|=\frac{\varpi_n}n.\]
Moreover, we set $\varpi_0:=0$.

This chapter is organized as follows. We set some notations and recall some known results in the following Subsection \ref{CH:3:defnsknwown}. Also, we give some  preliminary results on the contribution from infinity of sets in Section \ref{CH:3:contr_infty}.

In Section \ref{CH:3:classify}, we consider exterior data ``occupying at infinity'' in measure, with respect to an appropriate weight, less than an half-space. To be precise
\bgs{\label{CH:3:mainhyp} \alpha(E_0)<\frac{\varpi_n}{2}. }
In this hypothesis:
\begin{itemize}
\item In Subsection \ref{CH:3:sectnotfull} we give some asymptotic estimates of the density, in particular showing that when $s$ is small
enough $s$-minimal sets cannot fill their domain. 
\item In Subsection \ref{CH:3:estimatecurvature} we give some estimates on the fractional mean curvature. In particular we show that if a set $E$ has an exterior tangent ball of radius $\delta$ at some point $p\in\partial E$, then the $s$-fractional mean curvature of $E$ in $p$ is strictly positive for every $s<s_\delta$. 
\item In Subsection \ref{CH:3:alternative} we prove that when the fractional parameter is small and the exterior data at infinity occupies (in measure, with respect to the weight) less than half the space, then $s$-minimal sets completely stick at the boundary (that is, they are empty inside the domain), or become ``topologically dense'' in their domain. A similar result, which says that $s$-minimal sets fill the domain or their complementaries become dense, can be obtained in the same way, when the exterior data occupies in the appropriate sense more than half the space (so this threshold is somehow optimal).
\item Subsection \ref{CH:3:sticky} narrows the set of minimal sets that become dense in the domain for $s$ small. As a matter of fact, if the exterior data does not completely surround the domain, $s$-minimal sets completely stick at the boundary. 
\end{itemize}
In Section \ref{CH:3:sectexamples}, we provide some examples in which we are able to explicitly compute the contribution from infinity of sets. 
Section \ref{CH:3:cont} contains the continuity of the fractional mean curvature operator in all its variables for $s\in[0,1]$. As a corollary, we show that for $s\to 0^+$ the fractional mean curvature at a regular point of the boundary of a set, takes into account only the behavior of that set at infinity.   The continuity property implies that the mean curvature at a regular point on the boundary of a set may change sign, as $s$ varies,
depending on the signs of the two asymptotics as $s\to1^-$ and $s\to0^+$.

In Appendix \ref{CH:3:appendicite} and Appendix \ref{CH:3:appendicite2} we collect some useful results that we use in the present chapter. Worth mentioning are Appendixes \ref{CH:3:brr2} and \ref{CH:3:appendicite3}. The first of the two gathers some known results on the regularity of $s$-minimal surfaces, so as to state the Euler-Lagrange equation pointwisely in the interior of $\Omega$. In the latter we prove that the Euler-Lagrange equation holds (at least as a inequality) at $\partial E \cap \partial \Omega$, as long as the two boundaries do not intersect ``transversally''.


\subsection {Statements of the main results}

 We remark that the quantity $\alpha$,
\eqlab{\label{CH:3:alpha} \alpha(E)=\lim_{s\to 0^+} s\int_{\Co B_1} \frac{\chi_E(y)}{|y|^{n+s}}\, dy,  }
may not exist---see \cite[Example 2.8 and 2.9]{DFPV13}. For this reason, we define
\eqlab{\label{CH:3:baralpha1} \overline \alpha (E):= \limsup_{s\to 0^+} s\int_{\Co B_1} \frac{\chi_E(y)}{|y|^{n+s}}\, dy ,\quad \quad \underline \alpha(E) := \liminf_{s\to 0^+} s\int_{\Co B_1} \frac{\chi_E(y)}{|y|^{n+s}}\, dy.}

This set parameter plays an important role in describing
the asymptotic behavior of the fractional mean curvature as~$s\to0^+$
for unbounded sets. As a matter of fact, the limit as~$s\to0^+$
of the fractional mean curvature for a {\em bounded} set
is a positive, universal constant (independent of the set),
see, e.g., \cite[Appendix~B]{DV18}).
On the other hand, this asymptotic behavior changes for {\em unbounded}
sets, due to the set function $\alpha(E)$, as described explicitly
in the following result:
  \begin{theorem}\label{CH:3:asympts}[Proof in Section \ref{CH:3:cont}]
Let $E\subseteq\Rn$ and let $p\in\partial E$ be such that $\partial E$ is $C^{1,\gamma}$ near $p$,
for some $\gamma\in(0,1]$. Then
\bgs{& \liminf_{s\to0^+} s\,\I_s[E](p) =\varpi_n -2 \overline \alpha(E)\\
& \limsup_{s\to0^+} s\,\I_s[E](p) =\varpi_n-2 \underline\alpha(E).}
\end{theorem}
We notice that if~$E$ is bounded, then $\underline\alpha(E)
=\overline\alpha(E)=\alpha(E)=0$, hence Theorem~\ref{CH:3:asympts}
reduces in this case to formula~(B.1) in  \cite{DV18}.
Actually, we can estimate the fractional mean curvature from below (above) uniformly with respect to the radius of the exterior
(interior) tangent ball to $E$. To be more precise, if there exists an exterior tangent ball at $p\in \partial E$ of radius $\delta>0$, then for every $s<s_\delta$ we have 
\[ \liminf_{\varrho \to 0^+} s\,\I^\varrho_s[E](p)\geq \frac{\varpi_n -2\overline \alpha(E)}4.\]    
More explicitly, we have the following result:
   \begin{theorem}\label{CH:3:positivecurvature}[Proof in Section \ref{CH:3:estimatecurvature}]
Let $\Omega\subseteq\Rn$ be a bounded open set. 
Let $E_0\subseteq\Co\Omega$ be such that
\eqlab{\label{CH:3:weak_hp_beta}\overline \alpha(E_0)<\frac{\varpi_n}2,}
and let
\[\beta=\beta(E_0):=\frac{\varpi_n-2\overline \alpha(E_0)}4.\] We define 
\eqlab{\label{CH:3:delta_wild_index_def}
\delta_s=\delta_s(E_0):=e^{-\frac{1}{s}\log \frac{\varpi_n+2\beta}{\varpi_n+\beta}} ,}
for every $s\in(0,1)$.
Then, there exists $s_0=s_0(E_0,\Omega)\in(0,\frac{1}{2}]$ such that, if $E\subseteq\Rn$ is such that $E\setminus\Omega=E_0$
and $E$ has an exterior tangent ball
of radius (at least) $\delta_\sigma$, for some $\sigma\in(0,s_0)$, at some point $q\in\partial E\cap\overline{\Omega}$, then
\eqlab{\label{CH:3:unif_pos_curv}\liminf_{\varrho\to0^+}\I_s^\varrho[E](q)\geq\frac{\beta}{s}>0,\qquad\forall\,s\in(0,\sigma].}
 \end{theorem}

Given an open set $\Omega\subseteq\R^n$ and $\delta\in\R$, we consider the open set
\[\Omega_\delta:=\{x\in\R^n\,|\,\bar{d}_\Omega(x)<\delta\},\]
where $\bar{d}_\Omega$ denotes the signed distance function from $\partial\Omega$, negative inside $\Omega$.

It is well known (see, e.g., \cite{GilTru,Ambrosio})
that if $\Omega$ is bounded and $\partial \Omega$ is of class $C^2$, then the distance function is also of class $C^2$
in a neighborhood of $\partial\Omega$. Namely, there exists $r_0>0$
such that
\[\bar{d}_\Omega\in C^2(N_{2r_0}(\partial\Omega)),\quad
\mbox{ where }\quad N_{2r_0}(\partial\Omega):=\{x\in\R^n\,|\,|\bar{d}_\Omega(x)|<2r_0\}.\]
As a consequence, since $|\nabla\bar{d}_\Omega|=1$,
the open set $\Omega_\delta$ has $C^2$ boundary for every $|\delta|<2r_0$.
For a more detailed discussion, see Appendix \ref{CH:3:A2} and the references cited therein.

The constant $r_0$ will have the above meaning throughout this whole chapter.

\smallskip

We give the next definition.
\begin{defn}\label{CH:3:wild}
   Let $\Omega\subseteq \Rn$ be a bounded open set.
   We say that a set $E$ is $\delta$-{dense} in $\Omega$ for some fixed $\delta>0$ if $|B_\delta(x)\cap E|>0$ for any $x\in \Omega$ for which $B_\delta(x)\Subset\Omega$.
  \end{defn}
  
\noindent Notice that if $E$ is $\delta$-dense  then $E$ cannot have an exterior tangent ball of radius greater or equal than $\delta$ at any point $p\in \partial E\cap \Omega_{-\delta}$.

\noindent We observe that the notion for a set of
being $\delta$-dense is a ``topological'' notion,
rather than a measure theoretic one. 
Indeed, $\delta$-dense sets need not be ``irregular'' nor ``dense'' in the measure theoretic sense (see Remark \ref{CH:3:deltadance}).

\smallskip
With this definition and using Theorem \ref{CH:3:positivecurvature} we obtain the following classification.
\begin{theorem}\label{CH:3:THM}[Proof in Section \ref{CH:3:alternative}]
  Let $\Omega$ be a bounded and  connected open set with $C^2$ boundary. Let $E_0\subseteq \Co \Omega$ such that\[\overline \alpha(E_0)<\frac{\varpi_n}{2}.\]  
 Then the following two results hold.\\
  A)  Let $s_0$ and $\delta_s$ be as in Theorem \ref{CH:3:positivecurvature}. There exists
  $s_1=s_1(E_0,\Omega)\in (0,s_0]$ such that if $s<s_1$ and $E$ is an $s$-minimal set in $\Omega$ with exterior data $E_0$, then either
     \bgs{(A.1) \;  E\cap \Omega=\emptyset \quad  \mbox{ or} \quad\; (A.2)\;  E \mbox{ is } \delta_s-\mbox{dense}.}
 \noindent
 B) Either \\
(B.1) there exists
  $\tilde s=\tilde s(E_0,\Omega)\in (0,1)$ such that if $E$ is an $s$-minimal set in $\Omega$ with exterior data $E_0$ and $s\in(0,\tilde s)$, then
     \bgs{  E\cap \Omega=\emptyset,}
     or \\
    (B.2)    there exist  $\delta_k \searrow 0$, $s_k \searrow 0$ and a sequence of sets  $E_k$ such that each $E_k$ is $s_k$-minimal in $\Omega$ with exterior data $E_0$ and for every $k$
     \bgs{ \partial E_k \cap B_{\delta_k}(x) \neq \emptyset \quad \forall \; B_{\delta_k}(x)\Subset \Omega.}
     \end{theorem}

We remark here that Definition \ref{CH:3:wild} allows the $s$-minimal set 
to completely fill $\Omega$. The next theorem states that for $s$ small enough (and $\overline \alpha(E)<\varpi_n/2$) we can exclude  this possibility.
\begin{theorem}\label{CH:3:notfull}[Proof in Section \ref{CH:3:sectnotfull}]
Let $\Omega\subseteq\R^n$ be a bounded open set of finite classical perimeter and let $E_0\subseteq\Co\Omega$ be such that
\[\overline{\alpha}(E_0)<\frac{\varpi_n}{2}.\]
For every $\delta>0$ and every $\gamma\in(0,1)$ there exists $\sigma_{\delta,\gamma}=\sigma_{\delta,\gamma}(E_0,\Omega)\in(0,\frac{1}{2}]$ such that if $E\subseteq\R^n$ is $s$-minimal in $\Omega$, with exterior data $E_0$ and $s<\sigma_{\delta,\gamma}$, then
\eqlab{\label{CH:3:1666}\big|(\Omega\cap B_\delta(x))\setminus E\big|\ge\gamma\, \frac{\varpi_n-2\overline{\alpha}(E_0)}{\varpi_n-\overline{\alpha}(E_0)}\big|\Omega\cap B_\delta(x)\big|,\qquad\forall\,x\in\overline{\Omega}.}
\end{theorem}

\begin{remark}
 Let $\Omega$ and $ E_0$ be as in Theorem \ref{CH:3:notfull} and fix $\gamma=\frac{1}{2}$.
 \begin{enumerate}
\item Notice that we can find $\bar \delta >0$ and  $\bar x \in \Omega$ such that
\[ B_{2\bar\delta} (\bar x ) \subseteq \Omega.\]
Now if $s<\sigma_{\bar \delta,\frac{1}{2}}$ and $E$ is $s$-minimal in $\Omega$ with respect to $E_0$, \eqref{CH:3:1666} says that
\[ |B_{\bar\delta} (\bar x ) \cap \Co E|>0.\] 
Then (since the ball is connected), either $B_{\bar\delta} (\bar x ) \subseteq \Co E$ or there exists a point
\[x_0\in\partial E\cap \overline B_{\bar\delta} (\bar x ).\]
In this case, since $d(x_0, \partial \Omega )\ge\bar \delta$, \cite[Corollary 4.3]{CRS10} implies that
\[B_{\bar\delta c_s}(z)\subseteq\Co E\cap B_{\bar \delta}(x_0)\subseteq\Co E\cap\Omega\] for some $z$, where $c_s\in(0,1]$ denotes the constant of the clean ball condition (as introduced in \cite[Corollary 4.3]{CRS10}) and depends only on $s$ (and $n$). In both cases, there exists a ball of radius $\bar \delta c_s$ contained in $\Co E \cap \Omega$. 
\item If $s<\sigma_{\bar \delta,\frac{1}{2}}$ and $E$ is $s$-minimal and $\delta_s$-dense, then 
we have that
\[\delta_s>c_s\bar \delta.\]
On the other hand, we have an explicit expression for $\delta_s$, given in \eqref{CH:3:delta_wild_index_def}. Therefore, if one could prove that $c_s$ goes to zero slower than $\delta_s$, one could exclude the existence of $s$-minimal sets that are $\delta_s$-dense (for all sufficiently small $s$). 
\end{enumerate}
\end{remark}
\smallskip 

     
     An interesting result is related to $s$-minimal sets whose exterior data does not completely surround $\Omega$. In this case, the $s$-minimal set, for small values of $s$, is always
empty in $\Omega$. More precisely:

 \begin{theorem}\label{CH:3:boundedset}[Proof in Section \ref{CH:3:sticky}]
  Let $\Omega$ be a  bounded and  connected open set with $C^2$ boundary. Let $E_0\subseteq \Co \Omega$ such that  \[\overline \alpha(E_0)<\frac{\varpi_n}{2},\]
 and let $s_1$ be as in Theorem \ref{CH:3:THM}. Suppose that there exists $R>0$ and $x_0\in \partial \Omega$ such that \[B_R(x_0)\setminus \Omega \subseteq \Co E_0.\]  Then, there exists $s_3=s_3(E_0,\Omega)\in(0,s_1]$ such that if $s<s_3$ and $E$ is an $s$-minimal set in $\Omega$ with exterior data $E_0$, then 
    \[  E\cap \Omega=\emptyset .\]
     \end{theorem}

We notice that Theorem \ref{CH:3:boundedset} prevents the existence 
of $s$-minimal sets that are $\delta$-dense (for any $\delta$).    

\begin{remark}
The indexes $s_1$ and $s_3$ are defined as follows
\[s_1:=\sup\{s\in(0,s_0)\,|\,\delta_s<r_0\}\]
and
\[s_3:=\sup\Big\{s\in(0,s_0)\,\big|\,\delta_s<\frac{1}{2}\min\{r_0,R\}\Big\}.\]
Clearly, $s_3\leq s_1\leq s_0$.
\end{remark}

\begin{remark} We point out that condition \eqref{CH:3:weak_hp_beta}
is somehow optimal. Indeed,
when $\alpha(E_0)$ exists and 
\[ \alpha(E_0)=\frac{\varpi_n}2,\]
several configurations may occur, depending on the position of $\Omega$ with respect to the exterior data $E_0\setminus \Omega$. As an example, take 
\[ \mathfrak P =\{ (x',x_n) \; \big| \; x_n> 0\}.\] Then, for any $\Omega\subseteq \Rn$  a bounded open set with $C^2$ boundary, the only $s$-minimal set with exterior data given by $\mathfrak P \setminus \Omega$ is $\mathfrak P$ itself. So, if $E$ is $s$-minimal with respect to $\mathfrak P\setminus \Omega$ then
\bgs{&\Omega\subseteq \mathfrak P & \quad \implies \quad &E\cap \Omega=\Omega\\
 & \Omega\subseteq \Rn \setminus \mathfrak P &\quad \implies \quad &E\cap \Omega=\emptyset.} 
 On the other hand, if one takes $\Omega= B_1$, then 
 \[  E\cap B_1 = \mathfrak P  \cap B_1.\] 
 
 As a further example, we consider the supergraph
 \[ E_0:=\{(x',x_n) \; \big| \; x_n > \tanh x_1\},\] for which we have that (see Example \ref{CH:3:tanh})
 \[\alpha(E_0)=\frac{\varpi_n}2.\]  Then for every $s$-minimal set in $\Omega$ with exterior data $E_0\setminus \Omega$, we have that
 \bgs{&\Omega\subseteq \{ (x',x_n) \; \big| \; x_n> 1\} & \quad \implies \quad &E\cap \Omega=\Omega\\
 & \Omega\subseteq \{ (x',x_n) \; \big| \; x_n<-1\} &\quad \implies \quad &E\cap \Omega=\emptyset.} 
Taking $\Omega=B_2$, we have by the maximum principle in Proposition \ref{CH:3:maximum_principle}  that every set $E$ which is $s$-minimal in $B_2$, with respect to $E_0\setminus B_2$, satisfies
\bgs{  B_2\cap  \{ (x',x_n) \; \big| \; x_n> 1\}\subseteq E,   \qquad 
  B_2 \cap \{ (x',x_n) \; \big| \; x_n<-1\} \subseteq \Co E .} 
 On the other hand, we are not able to establish what happens in $B_2\cap \{ (x',x_n) \; \big| \; -1<x_n< 1\} $.
\end{remark}

\begin{remark}
We notice that when $E$ is $s$-minimal in $\Omega$ with respect to $E_0$, then $\Co E$ is $s$-minimal in $\Omega$ with respect to $\Co E_0$. Moreover
\[ \underline \alpha(E_0) >\frac{\varpi_n}{2} \qquad \implies \qquad \overline \alpha (\Co E_0)< \frac{\varpi_n}{2}.\]
So in this case we can apply Theorems \ref{CH:3:positivecurvature}, \ref{CH:3:THM}, \ref{CH:3:notfull} and \ref{CH:3:boundedset} to $\Co E$ with respect to the exterior data $\Co E_0$. For instance, if
$E$ is $s$-minimal in $\Omega$ with exterior data $E_0$ with
\[ \underline \alpha(E_0) >\frac{\varpi_n}{2}, \]
and $s<s_1(\Co E_0, \Omega)$,
 then either
\[ E\cap \Omega=\Omega \qquad \mbox{ or }  \qquad  \Co E \; \mbox{ is } \; \delta_s(\Co E_0)-\mbox{dense}.\]
 The analogues of the just mentioned Theorems can be obtained similarly.
\end{remark}

We point out that from our main results and the last two remarks, we have a complete classification of nonlocal minimal surfaces when $s$ is small whenever
\[ \alpha(E_0)\neq  \frac{\varpi_n}{2} .\] 

In the last section of the chapter, we prove the continuity of the fractional mean curvature in all variables (see Theorem \ref{CH:3:everything_converges} and Proposition \ref{CH:3:propsto0}). As a consequence, we have the following result. 

\begin{prop}\label{CH:3:rsdfyish}
Let $E\subseteq\R^n$ and let $p\in\partial E$ such that $\partial E$ is $C^{1,\alpha}$ in $B_R(p)$ for some
$R>0$ and $\alpha\in(0,1]$. Then the function
\[\I_{(\,\cdot\,)}[E](\,\cdot\,):(0,\alpha)\times(\partial E\cap B_R(p))\longrightarrow\R,
\qquad(s,x)\longmapsto\I_s[E](x)\]
is continuous.\\
Moreover, if $\partial E\cap B_R(p)$ is $C^2$ and for every $x\in \partial E\cap B_R(p)$ we define
\sys[ \tilde \I_s  {[}E{]} (x):=]{ &s(1-s)\I_s[E](x),  & \mbox{ for } &s\in (0,1) \\	
						&{\varpi_{n-1}} H[E](x), &\mbox{ for } &s=1,}
then the function
\[\tilde \I_{(\,\cdot\,)}[E](\,\cdot\,):(0,1]\times(\partial E\cap B_R(p))\longrightarrow\R,
\qquad(s,x)\longmapsto \tilde \I_s[E](x)\]
is continuous.\\
Finally, if $\partial E\cap B_R(p)$ is $C^{1,\alpha}$ and $\alpha(E)$ exists, and if for every $x\in \partial E\cap B_R(p)$ we denote
\[\tilde\I_0[E](x):=\varpi_n-2\alpha(E),\]
then the function
\[\tilde \I_{(\,\cdot\,)}[E](\,\cdot\,):[0,\alpha)\times(\partial E\cap B_R(p))\longrightarrow\R,
\qquad(s,x)\longmapsto \tilde \I_s[E](x)\]
is continuous.
\end{prop}

As a consequence of the continuity of the fractional mean curvature and the asymptotic result in  Theorem \ref{CH:3:asympts} we
establish that, by varying the fractional parameter $s$,
the nonlocal mean curvature may change sign at a point
where the classical mean curvature is negative, as one can observe in Theorem \ref{CH:3:changeyoursign}.


\subsection{Definitions, known facts and notations}\label{CH:3:defnsknwown}
We recall here some basic facts on $s$-minimal sets and surfaces, on the fractional mean curvature operator, and some notations, that we will use in the course of this chapter.

\subsubsection{Measure theoretic assumption}
We recall the following notations and measure theoretic assumptions, which are assumed throughout the chapter.

Let $E\subseteq\R^n$ be a measurable set. Up to modifications in sets of measure zero, we can assume (see Remark \ref{CH:1:gmt_assumption} and Appendix \ref{CH:1:Appendix_meas_th_bdary})
that $E$ contains the measure theoretic interior
\begin{equation*}
E_{int}:=\Big\{x\in\R^n\,|\,\exists\,r>0\textrm{ s.t. }|E\cap B_r(x)|=\frac{\varpi_n}n r^n\Big\}\subseteq E,
\end{equation*}
the complementary $\Co E$ contains the measure theoretic exterior
\begin{equation*}
E_{ext}:=\{x\in\R^n\,|\,\exists\,r>0\textrm{ s.t. }|E\cap B_r(x)|=0\}\subseteq\Co E,
\end{equation*}
and the topological boundary of $E$ coincides with its measure theoretic boundary, $\partial E=\partial^-E$,
where
\begin{equation*}\begin{split}
\partial^-E&:=\R^n\setminus(E_{int}\cup E_{ext})\\
&
=\Big\{x\in\R^n\,|\,0<|E\cap B_r(x)|<\frac{\varpi_n}{n}r^n\textrm{ for every }r>0\Big\}.
\end{split}
\end{equation*}
In particular, we remark that both $E_{int}$ and $E_{ext}$ are open sets.


\subsubsection{H\"{o}lder continuous functions}
We will use the following notation for the class of H\"{o}lder continuous functions.

Let $\alpha\in (0,1]$, let $S\subseteq\R^n$ and let $v:S\longrightarrow\R^m$. The $\alpha$-H\"{o}lder semi-norm of $v$ in $S$
is defined as
\[[v]_{C^{0,\alpha}(S,\R^m)}:=\sup_{x\neq y\in S}\frac{|v(x)-v(y)|}{|x-y|^\alpha}.\]
With a slight abuse of notation, we will omit the $\R^m$ in the formulas.
We also define
\[\|v\|_{C^0(S)}:=\sup_{x\in S}|v(x)|\quad\textrm{and}\quad\|v\|_{C^{0,\alpha}(S)}
:=\|v\|_{C^0(S)}+[v]_{C^{0,\alpha}(S)}.\]

Given an open set $\Omega\subseteq\R^n$, we define the space of uniformly H\"{o}lder continuous functions
$C^{0,\alpha}(\overline \Omega,\R^m)$ as
\[C^{0,\alpha}(\overline \Omega,\R^m):=\{v\in C^0(\overline{\Omega},\R^m)\,|\,
\|v\|_{C^{0,\alpha}(\overline{\Omega})}<\infty\}.\] 

Recall that $C^1(\overline{\Omega})$ is the space of those functions $u:\overline{\Omega}\longrightarrow\R$ such that
$u\in C^0(\overline{\Omega})\cap C^1(\Omega)$ and
such that $\nabla u$ can be continuously extended to $\overline{\Omega}$.
For every $S\subseteq\overline{\Omega}$ we write
\[\|u\|_{C^{1,\alpha}(S)}:=\|u\|_{C^0(S)}+\|\nabla u\|_{C^{0,\alpha}(S)},\]
and we define
\[C^{1,\alpha}(\overline \Omega):=\{u\in C^1(\overline{\Omega})\,|\,
\|u\|_{C^{1,\alpha}(\overline{\Omega})}<\infty\}.\]

We will usually consider the local versions of the above spaces. Given an open set $\Omega\subseteq\R^n$,
the space of locally H\"{o}lder continuous functions $C^{k,\alpha}(\Omega)$, with $k\in\{0,1\}$, is defined as
\[C^{k,\alpha}(\Omega):=\{u\in C^k(\Omega)\,|\,\|u\|_{C^{k,\alpha}(\mathcal O)}<\infty \mbox{ for every } \mathcal O\Subset \Omega\}.\] 

\subsubsection{The Euler-Lagrange equation}
 
We recall that the fractional mean curvature gives the Euler-Lagrange equation of an $s$-minimal set.
To be more precise, if $E$ is $s$-minimal in $\Omega$, then
\[\I_s[E]=0,\qquad\textrm{on}\quad\partial E\cap\Omega,\]
in an appropriate viscosity sense (see \cite[Theorem 5.1]{CRS10}).

Actually, by exploiting the interior regularity theory of $s$-minimal sets, the equation is satisfied in the classical sense
in a neighborhood of every ``viscosity point'' (see Appendix \ref{CH:3:brr2}). 
That is, if $E$ has at $p\in\partial E\cap\Omega$ a tangent ball (either interior or exterior), then $\partial E$ is $C^\infty$ in $B_r(p)$, for some $r>0$ small enough, and
\[\I_s[E](x)=0,\qquad\forall\,x\in\partial E\cap B_r(p).\] 
Moreover, if the boundary of $\Omega$ is of class $C^{1,1}$, then the Euler-Lagrange equation (at least as an inequality) holds also at a point $p\in\partial E\cap\partial\Omega$,
provided that the boundary $\partial E$ and the boundary $\partial\Omega$ do not intersect ``transversally'' in $p$ (see Theorem \ref{CH:3:EL_boundary_coroll}).


\section{Contribution to the mean curvature coming from infinity}\label{CH:3:contr_infty}
In this section, we study in detail the quantities $\alpha(E)$, $\overline\alpha(E),\underline \alpha(E)$) as defined in \eqref{CH:3:alpha}, \eqref{CH:3:baralpha1}. As a first remark, notice that these definitions  are independent on the radius of the ball (see
\cite[Observation 3 in Subsection 3.3]{DFPV13}) so we have that for any $R>0$
\eqlab{ \label{CH:3:baralpha}  \overline \alpha (E)= \limsup_{s\to 0^+} s\int_{\Co B_R} \frac{\chi_E(y)}{|y|^{n+s}}\, dy, \quad \underline \alpha(E) := \liminf_{s\to 0^+} s\int_{\Co B_R} \frac{\chi_E(y)}{|y|^{n+s}}\, dy .}
Notice that
\[ \overline \alpha(E) = \varpi_n -\underline \alpha(\Co E),  \quad \underline \alpha(E) = \varpi_n - \overline \alpha(\Co E).\] 
We define
\[ \alpha_s(q,r,E):=\int_{\Co B_r(q)} \frac{\chi_{E}(y) }{|q-y|^{n+s}} \, dy .\]
Then, the quantity $\alpha_s(q,r,E)$
somehow ``stabilizes'' for small $s$ independently on how large or where we take the ball, as rigorously given by the following result:

\begin{prop}\label{CH:3:unifrq}
Let $K\subseteq \Rn$ be a compact set and $[a,b]\subseteq \R$ be a closed interval, with $0<a<b$. Then 
\bgs{\label{CH:3:name1} \lim_{s\to 0^+}s|\alpha_s(q,r,E)-\alpha_s(0,1,E)| =0\quad \mbox{ uniformly in } q \in K, r\in [a,b].
}
Moreover,
for any bounded open set $\Omega\subseteq \Rn$ and any fixed $r>0$,
we have that \eqlab{\label{CH:3:claimalpha}\limsup_{s\to 0^+} s\inf_{q\in \overline \Omega} \alpha_s(q,r,E)= \limsup_{s\to 0^+} s\sup_{q\in \overline \Omega} \alpha_s(q,r,E)=\overline \alpha(E).}

\end{prop}
\begin{proof}
Let us fix $r\in [a,b]$ and $q\in K$, and $R>0$ such that $K\subseteq \overline B_{ R}$. Let also $\eps\in (0,1)$ be a fixed positive small quantity (that we will take arbitrarily small further on), such that 
\[ R>(\eps b)/(1-\eps).\]
We notice that if $x\in B_r(q)$, we have that 
$|x|<r+|q|<{R}/{\eps},$
hence $B_r(q)\subseteq B_{R/\eps}$.
We write that
\[\alpha_s(q,R,E)= \int_{\Co B_r(q)}\frac{\chi_E(y)}{|q-y|^{n+s}}\, dy = \int_{\Co B_{ R/\eps}} \frac{\chi_E(y)}{|q-y|^{n+s}}\, dy + \int_{ B_{ R/\eps}\setminus B_r(q)} \frac{\chi_E(y)}{|q-y|^{n+s}}\, dy.\]
Now for $y\in \Co B_{R/\eps}$ we have that $|y-q|\geq |y|-|q| \geq (1-\eps)|y|$, thus for any $q\in \overline B_R$ 
\eqlab{\label{CH:3:secondaaa1} \int_{ \Co B_{ R/\eps}} \frac{\chi_E(y)}{|q-y|^{n+s}} \, dy \leq &\; (1-\eps)^{-n-s} \int_{\Co B_{ R/\eps}} \frac{\chi_E(y)}{|y|^{n+s}}\, dy =(1-\eps)^{-n-s} \alpha_s(0,R/\eps,E)
. }
Moreover 
\eqlab{\label{CH:3:secondaaa11} \int_{ B_{ R/\eps}\setminus B_r(q)} \frac{\chi_E(y)}{|q-y|^{n+s}}\, dy \leq &\; \int_{B_{ R/\eps}\setminus B_r(q)} \frac{dy}{|q-y|^{n+s}} \leq \varpi_n \int_r^{R/\eps+ R} t^{-s-1}\, dt\\ = &\; \varpi_n \frac{r^{-s}- R^{-s}\eps^s(1+\eps)^{-s} }{s} \leq \varpi_n \frac{a^{-s}- R^{-s}\eps^s(1+\eps)^{-s} }{s} . }
Notice also that
since $B_r(q)\subseteq B_{R/\eps}$ and $|q-y|\leq |q|+|y|\leq (\eps+1)|y|$ for any $y\in \Co B_{R/\eps}$,
we obtain that
\eqlab{\label{CH:3:primariga}\int_{\Co B_r(q)}\frac{\chi_E(y) }{|q-y|^{n+s}}\, dy \geq &\; \int_{\Co B_{ R/\eps}} \frac{\chi_E(y) }{|q-y|^{n+s}}\, dy  \geq (1+\eps)^{-n-s} \int_{\Co B_{ R/\eps}} \frac{\chi_E(y)}{|y|^{n+s}}\, dy.
}
Putting\eqref{CH:3:secondaaa1}, \eqref{CH:3:secondaaa11} and \eqref{CH:3:primariga} together, we get that
\bgs{  0\leq \alpha_s(q,r,E) -(1+\eps)^{-n-s}  \alpha_s(0,R/\eps,E)
 \leq &\;\alpha_s(0, R/\eps, E) \left((1-\eps)^{-n-s}-(1+\eps)^{-n-s}\right)  \\
 &\;+ \varpi_n \frac{a^{-s}- R^{-s}\eps^s(1+\eps)^{-s} }{s}.}
Now we have that
\bgs{| \alpha_s(0, R/\eps,E)
-\alpha_s(0,1,E)| \leq \left|\int_{B_{R/\eps}\setminus B_1} \frac{dy}{|y|^{n+s}} \right| \leq \varpi_n \frac{ |1-R^{-s}\eps^{s}|}{s}. 
}
So by the triangle inequality we obtain 
\bgs{ |\alpha_s(q,r,E)-&(1+\eps)^{-n-s}\alpha_s(0,1,E)|
\leq \alpha_s(0, R/\eps, E)  \left((1-\eps)^{-n-s}-(1+\eps)^{-n-s}\right)
 \\
  &\;+ \frac{\varpi_n}s \big[ a^{-s}- R^{-s}\eps^s(1+\eps)^{-s} +  (1+\eps)^{-n-s}   { |1-R^{-s}\eps^{s}|}\big] . }
Hence,
it holds that
\[ \limsup_{s\to 0^+}s |\alpha_s(q,r,E)-(1+\eps)^{-n}\alpha_s(0,1,E)| \leq  \left((1-\eps)^{-n}-(1+\eps)^{-n}\right) \overline \alpha(E) ,\]
uniformly in $q\in K$ and in $r\in [a,b]$.\\
Letting $\eps \to 0^+$, 
we conclude that
 \[\lim_{s\to 0^+}s|\alpha_s(q,r,E)-\alpha_s(0,1,E)| =0,\]
uniformly in $q\in K$ and in $r\in [a,b]$.

Now, we consider $K$ such that $K=\overline \Omega$.
 Using the inequalities \eqref{CH:3:secondaaa1}, \eqref{CH:3:secondaaa11} and \eqref{CH:3:primariga} we have that for any $q\in \overline \Omega$
\bgs{ 
(1+\eps)^{-n-s} &\int_{\Co B_{ R/\eps}} \frac{\chi_E(y)}{|y|^{n+s}} \,dy
\leq \int_{\Co B_r(q)}\frac{\chi_E(y)}{|q-y|^{n+s}}\, dy\\
&
\le (1-\eps)^{-n-s} \int_{\Co B_{R/\eps}} \frac{\chi_E(y)}{|y|^{n+s}}\, dy
+ \varpi_n \frac{a^{-s}-R^{-s}\eps^s(1+\eps)^{-s}}s.
}
Passing to limsup 
it follows that
\bgs{ &(1+\eps)^{-n} \overline \alpha(E) \leq \limsup_{s\to 0^+} s\inf_{q\in \overline \Omega} \int_{\Co B_r(q)}\frac{\chi_E(y)}{|q-y|^{n+s}} \, dy\\& \leq \limsup_{s\to 0^+} s\sup_{q\in \overline \Omega}\int_{\Co B_r(q)}\frac{\chi_E(y)}{|q-y|^{n+s}} \, dy  \leq  (1-\eps)^{-n} \overline\alpha(E).}
Sending $\eps \to 0$  we obtain the conclusion.
\end{proof}

\begin{remark}\label{CH:3:finmeas}
Let $E\subseteq \Rn$ be such that $|E|<\infty$. Then
\[ \alpha(E)=0.\]
Indeed,
\[ |\alpha_s(0,1,E)|\leq |E|,\]
hence
\[ \limsup_{s\to 0} s|\alpha_s(0,1,E)|=0.\]
 \end{remark}
 
Now, we discuss some useful properties of $\overline\alpha$. 
Roughly speaking, the quantity $\overline\alpha$ takes into account
the ``largest possible asymptotic opening'' of a set, and so it
possesses nice geometric features such as monotonicity, additivity
and geometric invariances. The detailed list of these properties is
the following:

\begin{prop}\label{CH:3:subsetssmin} \quad \\ 
(i) (Monotonicity) Let $E,F\subseteq \Rn$ be such that for some $r>0$ and $q\in \Rn$\[ E\setminus B_r(q)\subseteq  F\setminus B_r(q).\] Then
\[\overline \alpha(E)\leq \overline \alpha(F).\]
(ii) (Additivity) Let $E,F\subseteq \Rn$ be such that for some $r>0$ and $q\in \Rn$ \[ (E\cap F)\setminus B_r(q)= \emptyset.\] Then
\[ \overline \alpha (E\cup F)\leq \overline \alpha(E)+\overline \alpha(F).\]
Moreover, if $\alpha(E), \alpha(F)$ exist, then $\alpha(E\cup F)$ exists and
\[ \alpha(E\cup F)= \alpha(E)+\alpha(F).\]
(iii) (Invariance with respect
to rigid motions) Let $E\subseteq \Rn$, $x\in \Rn$ and $\mathcal R \in \mathcal {SO}(n)$ be a rotation. Then
\[ \overline \alpha(E+x)=\overline \alpha(E) \quad{\mbox{ and }}\quad \overline \alpha( \mathcal R E)=\overline \alpha(E).\]
(iv) (Scaling) Let $E\subseteq \Rn$ and $\lambda >0$. Then for some $r>0$ and $q\in \Rn$
\[ \alpha_s(q,r,\lambda E) =  \lambda^{-s} 
\alpha_s\left(\frac{q}{\lambda},\frac{r}{\lambda},E\right)
\quad{\mbox{ and }}\quad
\overline \alpha(\lambda E) =\overline \alpha(E).\] 
(v) (Symmetric difference)
Let  $E, F\subseteq \Rn$. Then for every $r>0$ and $q\in \Rn$
\[ |\alpha_s(q,r,E)-\alpha_s(q,r,F)|\leq \alpha_s(q,r,E\Delta F).\]
As a consequence, if $|E\Delta F|<\infty$ and $\alpha(E)$  exists, then $\alpha(F)$ exists and  
\[\alpha(E)=\alpha(F). \]
\end{prop}

\begin{proof} (i) It is enough to notice that for every $s\in (0,1)$
\[ \alpha_s(q, r,E) \leq \alpha_s(q, r,F) .\]
Then, passing to limsup and recalling \eqref{CH:3:claimalpha} we conclude that
\[\overline \alpha(E)\leq \overline \alpha(F).\]
(ii)  We notice that for every $s\in (0,1)$ \[ \alpha_s(q,r,E\cup F) = \alpha_s(q,r,E)+\alpha_s(q,r, F) \]  
and passing to limsup and liminf as $s\to 0^+$  we obtain the desired claim.\\
(iii) By a change of variables,
we have that
\[ \alpha_s(0,1,E+x)= \int_{\Co B_1} \frac{\chi_{E+x}(y)}{|y|^{n+s}}\, dy = \int_{\Co B_1(-x)} \frac{\chi_{E}(y)}{|x+y|^{n+s}}\, dy= \alpha_s(-x,1,E).\]
Accordingly, the invariance by translation 
follows after passing to limsup and using \eqref{CH:3:claimalpha}. \\
In addition, the invariance by rotations is obvious, using a change of variables.\\
(iv) Changing the variable $y=\lambda x$ we deduce that 
\bgs{
\alpha_s(q,r,\lambda E)&
=\int_{\Co B_r(q)}\frac{\chi_{\lambda E} (y)}{|q-y|^{n+s}}\, dy
=\lambda ^{-s}\int_{\Co B_{\frac{r}{\lambda}}(\frac{q}{\lambda})} \frac{\chi_E(x)}{|\frac{q}{\lambda}-x|^{n+s}}\, dx\\
&
=\lambda^{-s} \alpha_s\left(\frac{q}{\lambda},\frac{r}{\lambda},E\right).
}
Hence, the claim follows by passing to limsup as $s\to 0^+$.\\
(v) We have that
\bgs{
|\alpha_s(q,r,E)-\alpha_s(q,r,F)| &\leq \int_{\Co B_r(q)} \frac{ |\chi_{E}(y)-\chi_F(y)|}{|y-q|^{n+s}}\, dy
= \int_{\Co B_r(q)} \frac{ \chi_{E\Delta F}(y)}{|y-q|^{n+s}}\, dy	\\
&
= \alpha_s(q,r,E\Delta F).
}
The second part of the claim follows applying the Remark \ref{CH:3:finmeas}.
\end{proof}

We recall the definition (see (3.1) in  \cite{DFPV13})
\[\mu(E):=\lim_{s \to 0^+} s \Per_s(E,\Omega),\]
where $\Omega$ is a bounded open set with $C^2$ boundary.
Moreover, we define 
\[ \overline \mu(E)= \limsup_{s\to 0^+} s\Per_s(E,\Omega)\]
and give the following result:
\begin{prop}\label{CH:3:barmubaral} Let $\Omega\subseteq\Rn$ be a bounded open set with finite classical perimeter
and let $E_0\subseteq\Co \Omega$. Then
\[ \overline \mu(E_0)= \overline \alpha(E_0) |\Omega|.\]
\end{prop}
\begin{proof}
Let $R>0$ be fixed such that $\Omega \subseteq B_R$, $y\in \Omega$ be any fixed point and $\eps\in(0,1)$ be small enough such that $R/\eps> R+1$. This choice of $\eps$ assures that $B_1(y)\subseteq B_{R/\eps}$. 
We have that
\bgs{
\int_{\Rn} \frac{\chi_{E_0}(x)}{|x-y|^{n+s}}\, dx &=  \int_{\Co B_{R/\eps} } \frac{\chi_{E_0}(x)}{|x-y|^{n+s}}\, dx 
+ \int_{ B_{R/\eps}\setminus B_{1}(y)} \frac{\chi_{E_0}(x)}{|x-y|^{n+s}}\, dx\\
&\qquad\quad
+ \int_{  B_{1}(y)} \frac{\chi_{E_0}(x)}{|x-y|^{n+s}}\, dx.
}
Since $|x-y|\geq (1-\eps)|x|$ whenever $x\in \Co B_{R/\eps}$, we get
\[ \int_{\Co B_{R/\eps} } \frac{\chi_{E_0}(x)}{|x-y|^{n+s}}\, dx\leq (1-\eps)^{-n-s}\int_{\Co B_{R/\eps} } \frac{\chi_{E_0}(x)}{|x|^{n+s}}\, dx.\]
Also we have that
\[ \int_{ B_{R/\eps}\setminus B_{1}(y)} \frac{\chi_{E_0}(x)}{|x-y|^{n+s}}\, dx\leq \varpi_n\int_1^{R/\eps+R} \varrho^{-s-1}\, d\varrho \leq \varpi_n\frac{1- \left(\frac{R}{\eps}+R\right)^{-s}}s.\]
Also, we can assume that $s<1/2$ 
(since we are interested in what happens for $s\to 0$). In this way,
if $|x-y|<1$ we have that $|x-y|^{-n-s}\leq|x-y|^{-n-\frac{1}{2}}$, and so
\[\int_{  B_1(y)} \frac{\chi_{E_0}(x)}{|x-y|^{n+s}}\, dx\leq \int_{  B_1(y)} \frac{\chi_{E_0}(x)}{|x-y|^{n+\frac{1}2}}\, dx.\] 
Also, since $E_0\subseteq\Co\Omega$, we have that 
\[\int_{  B_1(y)} \frac{\chi_{E_0}(x)}{|x-y|^{n+\frac{1}2}}\, dx  \leq \int_{  B_1(y)\setminus \Omega} \frac{dx}{|x-y|^{n+\frac{1}2}}
\leq \int_{\Co \Omega} \frac{dx}{|x-y|^{n+\frac{1}2}}.  \] 
This means that
\[\int_{\Omega} \int_{  B_1(y)} \frac{\chi_{E_0}(x)}{|x-y|^{n+s}}\, dx\,dy\leq \int_{\Omega} \int_{\Co \Omega} \frac{dx}{|x-y|^{n+\frac{1}2}}=\Per_{\frac{1}2}(\Omega)=c<\infty,\]
since $\Omega$ has a finite classical perimeter.
In this way, it follows that
\eqlab{\label{CH:3:mumu1}
s&\Per_s(E_0,\Omega) = \int_{\Omega} \int_{\Rn}\frac{\chi_{E_0}(x) }{|x-y|^{n+s}}dx \, dy\\
&\leq s (1-\eps)^{-n-s} |\Omega|\int_{\Co B_{R/\eps} } \frac{\chi_{E_0}(x)}{|x|^{n+s}} dx
+\varpi_n\left(1- \left(\frac{R}{\eps}+R\right)^{-s}\right)|\Omega| +sc.
}
Furthermore, notice that if $x\in B_{R/\eps}$ we have that $|x-y|\leq (1+\eps)|x|$, hence
\[ \int_{\Rn}\frac{\chi_{E_0}(x)}{|x-y|^{n+s}}\, dx\geq \int_{\Co B_{{R}/{\eps}}} \frac{\chi_{E_0}(x)}{|x-y|^{n+s}} \, dx \geq (1+\eps)^{-n-s} \int_{\Co B_{{R}/{\eps}}} \frac{\chi_{E_0}(x)}{|x|^{n+s}} \, dx.  \]
Thus for any $\eps>0$
\[ s\Per_s(E_0,\Omega) \geq s |\Omega|  (1+\eps)^{-n-s}  \int_{\Co B_{{R}/{\eps}}} \frac{\chi_{E_0}(x)}{|x|^{n+s}} \, dx .\] 
Passing to  limsup as $s\to 0^+$ here above and in \eqref{CH:3:mumu1} it follows that
\[ (1+\eps)^{-n} \overline \alpha(E_0)\,  |\Omega| \leq \overline \mu(E_0) \leq (1-\eps)^{-n} \overline \alpha(E_0)\, |\Omega| .\]
Sending $\eps \to 0$, we obtain the desired conclusion.
\end{proof}

\section{Classification of nonlocal minimal surfaces for small $s$}\label{CH:3:classify}

\subsection{Asymptotic estimates of the density (Theorem \ref{CH:3:notfull})}\label{CH:3:sectnotfull}

 
The importance of Theorem \ref{CH:3:notfull} is threefold:
\begin{itemize}
\item first of all, it is an interesting result in itself, by stating (in the usual hypothesis in which the contribution from infinity of the exterior data $E_0$ is less than that of a half-space) that any ball of fixed radius, centered at some $x\in\overline\Omega$, contains at least a portion of the complement of an $s$-minimal set $E$, when $s$ is small enough. We further observe that Theorem \ref{CH:3:notfull} actually provides a ``uniform'' measure theoretic estimate of how big this portion is, purely in terms of the fixed datum
$\overline{\alpha}(E_0)$.

\item Moreover, we point out that Definition \ref{CH:3:wild} does not exlude apriori ``full'' sets, i.e. sets $E$ such that $E\cap\Omega=\Omega$.  Hence, in the situation of point $(A)$ of Theorem \ref{CH:3:THM}, one may wonder whether an $s$-minimal set $E$, which is
$\delta_s$-dense, can actually completely cover $\Omega$.  The answer is no: Theorem \ref{CH:3:notfull}
proves in particular that the contribution from infinity forces the domain $\Omega$, for $s$ small enough, to contain at least
a non-trivial portion of the complement of $E$.

\item Finally, the density estimate of Theorem \ref{CH:3:notfull} serves as an auxiliary result for the proof of part (B) of our main Theorem \ref{CH:3:THM}.
\end{itemize}

 \begin{proof}[Proof of Theorem \ref{CH:3:notfull}]
We begin with two easy but useful preliminary remarks. We
observe that,
given a set $F\subseteq\R^n$ and two open sets $\Omega'\subseteq\Omega$, we have
\eqlab{\label{CH:3:monny}
\Per_s(F,\Omega')\leq \Per_s(F,\Omega).
}
Also, we point out that,
given an open set $\mathcal O\subseteq\R^n$ and a set $F\subseteq\R^n$, then by the definition
of the fractional perimeter,
it holds
\eqlab{\label{CH:3:Yamamoto}
F\cap\Omega=\emptyset\quad\implies\quad
\Per_s(F,\mathcal O)=\int_F\int_{\mathcal O}\frac{dx\,dy}{|x-y|^{n+s}}.
}

With these observations at hand, we are ready to proceed with the proof of the Theorem.
We argue by contradiction.

Suppose that there exists $\delta>0$ and $\gamma\in(0,1)$ for which we can find a sequence $s_k\searrow0$,
a sequence of sets $\{E_k\}$ such that each $E_k$ is $s_k$-minimal in $\Omega$ with exterior data $E_0$, and a sequence
of points $\{x_k\}\subseteq\overline{\Omega}$ such that
\eqlab{\label{CH:3:density_contrad_proof}
\big|(\Omega\cap B_\delta(x_k))\setminus E_k\big|< \gamma\,\frac{\varpi_n-2\overline{\alpha}(E_0)}{\varpi_n-\overline{\alpha}(E_0)}\big|\Omega\cap B_\delta(x_k)\big|.
}
As a first step, we are going to exploit \eqref{CH:3:density_contrad_proof} in order to obtain a
bound from below for the limit as $k\to\infty$ of
$s_k \Per_{s_k}(E_k,\Omega\cap B_\delta(x_k))$ (see the forthcoming inequality \eqref{CH:3:dolph1}).

First of all we remark that, since $\overline{\Omega}$ is compact, up to passing to subsequences we can suppose that $x_k\longrightarrow x_0$, for some $x_0\in\overline{\Omega}$.
Now we observe that from \eqref{CH:3:density_contrad_proof} it follows that 
\bgs{
|E_k\cap(\Omega\cap B_\delta(x_k))\big|&=|\Omega\cap B_\delta(x_k)|-\big|(\Omega\cap B_\delta(x_k))\setminus E_k\big|\\
&
>\frac{(1-\gamma)\varpi_n-(1-2\gamma)\overline{\alpha}(E_0)}{\varpi_n-\overline{\alpha}(E_0)}\,|\Omega\cap B_\delta(x_k)|,
}
and hence, since $x_k\longrightarrow x_0$,
\eqlab{\label{CH:3:Esti_pf_eqn1}
\liminf_{k\to\infty}|E_k\cap(\Omega\cap B_\delta(x_k))\big|
\geq\frac{(1-\gamma)\varpi_n-(1-2\gamma)\overline{\alpha}(E_0)}{\varpi_n-\overline{\alpha}(E_0)}\,|\Omega\cap B_\delta(x_0)|.
}
Notice that, since $\Omega$ is bounded, we can find $R>0$ such that $\Omega\Subset B_R(q)$ for every $q\in\overline{\Omega}$.
Then
we obtain that
\bgs{
\Per_{s_k}(E_k,\Omega\cap B_\delta(x_k))&\ge
\int_{E_k\cap(\Omega\cap B_\delta(x_k))}
\Big(\int_{\Co E_k\setminus(\Omega\cap B_\delta(x_k))}\frac{dz}{|y-z|^{n+s_k}}\Big)dy\\
&
\ge\int_{E_k\cap(\Omega\cap B_\delta(x_k))}\Big(\int_{\Co \Omega}\frac{\chi_{\Co E_0}(z)}{|y-z|^{n+s_k}}\,dz\Big)dy\\
&\ge \int_{E_k\cap(\Omega\cap B_\delta(x_k))}\Big(\inf_{q\in \overline \Omega } \int_{\Co \Omega}\frac{\chi_{\Co E_0}(z)}{|q-z|^{n+s_k}}\,dz\Big)dy\\
 &
\ge\big|E_k\cap(\Omega\cap B_\delta(x_k))\big|\inf_{q\in\overline{\Omega}}\int_{\Co B_R(q)}
\frac{\chi_{\Co E_0}(z)}{|q-z|^{n+s_k}}\,dz.
}
So, thanks to Proposition \ref{CH:3:unifrq} and recalling \eqref{CH:3:Esti_pf_eqn1}, we find
\eqlab{\label{CH:3:dolph1}
\liminf_{k\to\infty}s_k \Per_{s_k}&(E_k,\Omega\cap B_\delta(x_k))\\
&
\ge\Big(\liminf_{k\to\infty}|E_k\cap(\Omega\cap B_\delta(x_k))\big|\Big)\Big(\liminf_{k\to\infty}s_k\,
\inf_{q\in\overline{\Omega}}\int_{\Co B_R(q)}\frac{\chi_{\Co E_0}(z)}{|q-z|^{n+s_k}}\,dz\Big)\\
&
=\big(\varpi_n-\overline{\alpha}(E_0)\big)\Big(\liminf_{k\to\infty}|E_k\cap(\Omega\cap B_\delta(x_k))\big|\Big)\\
&
\geq\big(\varpi_n-\overline{\alpha}(E_0)\big)
\frac{(1-\gamma)\varpi_n-(1-2\gamma)\overline{\alpha}(E_0)}{\varpi_n-\overline{\alpha}(E_0)}\,|\Omega\cap B_\delta(x_0)|.
}

On the other hand, as a second step we claim that
\eqlab{\label{CH:3:contrad_dens_limsup}
\limsup_{k\to\infty}s_k \Per_{s_k}(E_k,\Omega\cap B_\delta(x_k))\le\overline{\alpha}(E_0)\big|\Omega\cap B_\delta(x_0)\big|.
}
We point out that obtaining the inequality \eqref{CH:3:contrad_dens_limsup} is a crucial step of the proof. Indeed,
exploiting both \eqref{CH:3:contrad_dens_limsup} and \eqref{CH:3:dolph1}, we obtain
\eqlab{\label{CH:3:contrad_dens_liminf}
\overline{\alpha}(E_0)\,|\Omega\cap B_\delta(x_0)|&\ge\liminf_{k\to\infty}s_k \Per_{s_k}(E_k,\Omega\cap B_\delta(x_k))\\
&
\ge
\big((1-\gamma)\varpi_n-(1-2\gamma)\overline{\alpha}(E_0)\big)|\Omega\cap B_\delta(x_0)|.}
Then, since $x_0\in\overline{\Omega}$ implies that
\[
|\Omega\cap B_\delta(x_0)|>0,
\]
by \eqref{CH:3:contrad_dens_liminf} we get
\[
\overline{\alpha}(E_0)\ge(1-\gamma)\varpi_n-(1-2\gamma)\overline{\alpha}(E_0)
\quad\textrm{ that is }\quad
(1-\gamma)\overline{\alpha}(E_0)\ge(1-\gamma)\frac{\varpi_n}{2}.
\]
Therefore, since $\gamma\in(0,1)$ and by hypothesis $\overline{\alpha}(E_0)<\frac{\varpi_n}{2}$, we reach a contradiction, concluding the proof.

We are left to prove \eqref{CH:3:dolph1}. For this, we exploit the minimality of the sets $E_k$ in order to compare the $s_k$-perimeter
of $E_k$ with the $s_k$-perimeter of appropriate competitors $F_k$.

We first remark that, since $x_k\longrightarrow x_0$, for every $\eps>0$
there exists $\tilde k_\eps$ such that
\eqlab{\label{CH:3:subset_balls_eps}
\Omega\cap B_\delta (x_k)\subseteq\Omega\cap B_{\delta+\eps}(x_0),\qquad
\forall\,k\ge\tilde k_\eps.}
We fix a small $\eps>0$. We will let $\eps\to0$ later on.

We also observe that, since $E_k$ is $s_k$-minimal in $\Omega$, it is $s_k$-minimal also in every
$\Omega'\subseteq\Omega$, hence in particular in $\Omega \cap B_{\delta+\eps}(x_0)$.
Now we proceed to define the sets
\eqlab{\label{CH:3:density_competitor_def}
F_k:=E_0\cup(E_k\cap(\Omega\setminus B_{\delta+\eps}(x_0)))
=E_k\setminus \big(\Omega\cap B_{\delta+\eps}(x_0)\big).
}
Then, by \eqref{CH:3:monny}, \eqref{CH:3:subset_balls_eps}, \eqref{CH:3:density_competitor_def} and by the minimality of $E_k$
in $\Omega \cap B_{\delta+\eps}(x_0)$, for every $k\ge\tilde k_\eps$ we find that
\bgs{
\Per_{s_k}(E_k,\Omega\cap B_\delta(x_k))\le \Per_{s_k}(E_k,\Omega\cap B_{\delta+\eps}(x_0))
\le \Per_{s_k}(F_k,\Omega\cap B_{\delta+\eps}(x_0)).
}
We observe that by the definition \eqref{CH:3:density_competitor_def} we have that
\[
F_k\cap\big(\Omega\cap B_{\delta+\eps}(x_0)\big)=\emptyset.
\] 
Therefore, recalling \eqref{CH:3:Yamamoto}
and the definition \eqref{CH:3:density_competitor_def} of the sets $F_k$, we obtain that
\bgs{
\Per_{s_k}&(F_k,\Omega\cap B_{\delta+\eps}(x_0))
=\int_{E_0\cup(E_k\cap(\Omega\setminus B_{\delta+\eps}(x_0)))}\int_{\Omega\cap B_{\delta+\eps}(x_0)}
\frac{dy\,dz}{|y-z|^{n+s_k}}\\
&
=\int_{E_0}\int_{\Omega\cap B_{\delta+\eps}(x_0)}\frac{dy\,dz}{|y-z|^{n+s_k}}
+
\int_{E_k\cap(\Omega\setminus B_{\delta+\eps}(x_0))}\int_{\Omega\cap B_{\delta+\eps}(x_0)}\frac{dy\,dz}{|y-z|^{n+s_k}}\\
&
\le\int_{E_0}\int_{\Omega\cap B_{\delta+\eps}(x_0)}\frac{dy\,dz}{|y-z|^{n+s_k}}
+\int_{\Omega\setminus B_{\delta+\eps}(x_0)}\int_{\Omega\cap B_{\delta+\eps}(x_0)}\frac{dy\,dz}{|y-z|^{n+s_k}}\\
&
=:I^1_k+I^2_k.
}
Furthermore, again by \eqref{CH:3:Yamamoto}, we have that 
\eqlab{\label{CH:3:dens_asympt_pf}
I^1_k=\Per_{s_k}(E_0,\Omega\cap B_{\delta+\eps}(x_0))
\quad\mbox{and}\quad
I^2_k=\Per_{s_k}(\Omega\setminus B_{\delta+\eps}(x_0),\Omega\cap B_{\delta+\eps}(x_0)).
}
We observe that the open set $\Omega\cap B_{\delta+\eps}(x_0)$ has finite classical perimeter.
Thus, we can exploit the equalities \eqref{CH:3:dens_asympt_pf} and apply Proposition \ref{CH:3:barmubaral} twice, obtaining
\[
\limsup_{k\to\infty}s_k I^1_k\le
\overline{\alpha}(E_0)\big|\Omega\cap B_{\delta+\eps}(x_0)\big|,
\]
and
\eqlab{\label{CH:3:dens_asympt_pf1}
\limsup_{k\to\infty}s_k I^2_k\le\overline{\alpha}(\Omega\setminus B_{\delta+\eps}(x_0))\big|\Omega\cap B_{\delta+\eps}(x_0)\big|,
}
for every $\eps>0$. Also notice that, since $\Omega$ is bounded, by Remark \ref{CH:3:finmeas} we have
\[
\overline{\alpha}(\Omega\setminus B_{\delta+\eps}(x_0))=\alpha(\Omega\setminus B_{\delta+\eps}(x_0))=0,
\]
and hence, by \eqref{CH:3:dens_asympt_pf1},
\[
\lim_{k\to\infty}s_kI^2_k=0.
\]
Therefore, combining these computations we find that
\bgs{
\limsup_{k\to\infty}s_k \Per_{s_k}(E_k,\Omega\cap B_\delta(x_k))\le
\limsup_{k\to\infty}s_k I^1_k\le
\overline{\alpha}(E_0)\big|\Omega\cap B_{\delta+\eps}(x_0)\big|,
}
for every $\eps>0$ small.
To conclude, we let $\eps\to0$ and we obtain \eqref{CH:3:contrad_dens_limsup}.
\end{proof}

It is interesting to observe that, as a straightforward consequence of Theorem \ref{CH:3:notfull}, when $\alpha(E_0)=0$ we
know that any sequence of $s$-minimal sets is asymptotically empty inside $\Omega$, as $s\to0^+$. More precisely 
\begin{corollary}
Let $\Omega\subseteq\R^n$ be a bounded open set of finite classical perimeter and let $E_0\subseteq\Co\Omega$ be such that
$\alpha(E_0)=0$. Let $s_k\in(0,1)$ be such that $s_k\searrow0$ and let $\{E_k\}$
be a sequence of sets such that each $E_k$ is $s_k$-minimal in $\Omega$ with exterior data $E_0$. Then
\[\lim_{k\to\infty}|E_k\cap\Omega|=0.\]
\end{corollary}

\begin{proof}
Fix $\delta>0$. Since $\overline{\Omega}$ is compact, we can find a finite number of points $x_1,\dots,x_m\in\overline{\Omega}$
such that
\[\overline{\Omega}\subseteq\bigcup_{i=1}^mB_\delta(x_i).\]
By Theorem \ref{CH:3:notfull} (by using the fact that $\alpha(E_0)=0$) we know that for every $\gamma\in(0,1)$ we can find a $k(\gamma)$ big enough such that
\bgs{
\big|(\Omega\cap B_\delta(x_i))\setminus E _k\big|\ge\gamma\, \big|\Omega\cap B_\delta(x_i)\big|.
}
Then, 
\bgs{
		\big|E_k\cap(\Omega\cap B_\delta(x_i))\big|= \big|\Omega\cap B_\delta(x_i)\big|- 	\big|(\Omega\cap B_\delta(x_i))\setminus E_k\big|
				\le 
(1-\gamma)|\Omega\cap B_\delta(x_i)|,
}
for every $i=1,\dots,m$ and every $k\ge k(\gamma)$. Thus
\[|E_k\cap\Omega|\le(1-\gamma)\sum_{i=1}^m|\Omega\cap B_\delta(x_i)|,\]
for every $k\ge k(\gamma)$, and hence
\[\limsup_{k\to\infty}|E_k\cap\Omega|\le(1-\gamma)\sum_{i=1}^m|\Omega\cap B_\delta(x_i)|,\]
for every $\gamma\in(0,1)$. Letting $\gamma\longrightarrow1^-$ concludes the proof.
\end{proof}

We recall here that any set $E_0$ of finite measure has $\alpha(E_0)=0$ (check Remark \ref{CH:3:finmeas}).

\subsection{Estimating the fractional mean curvature (Theorem \ref{CH:3:positivecurvature})}\label{CH:3:estimatecurvature}

Thanks to the previous preliminary work, we are now in the position
of completing the proof of Theorem \ref{CH:3:positivecurvature}.

\begin{proof}[Proof of Theorem \ref{CH:3:positivecurvature}]
Let $R:=2\,\max\{1,\textrm{diam}(\Omega)\}$. First of all, \eqref{CH:3:claimalpha} implies that
\bgs{ \liminf_{s\to 0^+} \bigg(\varpi_n R^{-s} -2s\sup_{q\in \overline \Omega} \int_{\Co B_R(q)}\frac{\chi_E(y)}{|q-y|^{n+s}} \, dy \bigg) = {\varpi_n-2\overline \alpha(E_0)}=4 \beta.} Notice that by \eqref{CH:3:weak_hp_beta}, $\beta>0$. Hence for every $s$ small enough, say $s<s'\leq\frac{1}{2}$ with $s'=s'(E_0,\Omega)$, we have that
\eqlab{\label{CH:3:tildebeta}\varpi_n R^{-s} - 2s\sup_{q\in \overline \Omega} \int_{\Co B_R(q)}\frac{\chi_E(y)}{|q-y|^{n+s}} \, dy \geq \frac {7}2 \beta.}

Now, let $E\subseteq\Rn$ be such that $E\setminus\Omega=E_0$, suppose that $E$ has an exterior tangent ball
of radius $\delta<R/2$ at $q\in\partial E\cap\overline{\Omega}$, that is
\[B_\delta(p)\subseteq\Co E\quad\textrm{and}\quad q\in\partial B_\delta(p),\]
and let $s<s'$.
Then for $\varrho$ small enough (say $\varrho<\delta/2$) we conclude that
\bgs{
	\I^\varrho_s[E](q)= \int_{B_R(q)\setminus B_{\varrho}(q)} \frac{\chi_{ \Co E} (y)-\chi_{E}(y)}{|q-y|^{n+s}}\, dy + \int_{\Co B_R(q)} \frac{\chi_{ \Co E} (y)-\chi_{E}(y)}{|q-y|^{n+s}}\, dy.
}

 Let $D_\delta= B_\delta(p) \cap B_\delta(p')$, where $p'$ is the symmetric of $p$ with respect to $q$, i.e. the ball $B_\delta(p')$ is the ball tangent to  $B_\delta(p)$ in $q$. Let also $K_\delta$ be the convex hull of $D_\delta$ and let $\Per_\delta:=K_\delta-D_\delta$. Notice that $B_\varrho(q)\subseteq K_\delta \subseteq B_R(q)$ . Then
\bgs{
\int_{B_R(q)\setminus B_{\varrho}(q)} &\frac{\chi_{ \Co E} (y)-\chi_{E}(y)}{|q-y|^{n+s}}\, dy
= \int_{D_\delta\setminus B_{\varrho}(q)} \frac{\chi_{ \Co E} (y)-\chi_{E}(y)}{|q-y|^{n+s}}\, dy\\
&
+\int_{\Per_\delta\setminus B_{\varrho}(q) } \frac{\chi_{ \Co E} (y)-\chi_{E}(y)}{|q-y|^{n+s}}\, dy
+\int_{B_R(q) \setminus K_\delta} \frac{\chi_{ \Co E} (y)-\chi_{E}(y)}{|q-y|^{n+s}}\, dy.
}
Since $B_\delta(p) \subseteq \Co E$, by symmetry we obtain that
\bgs{
\int_{D_\delta \setminus B_{\varrho}(q)} &\frac{\chi_{ \Co E} (y)-\chi_{E}(y)}{|q-y|^{n+s}}\, dy\\
& =\int_{B_\delta(p)\setminus B_{\varrho}(q)} \frac{dy}{|q-y|^{n+s}}
+\int_{B_\delta(p')\setminus B_{\varrho}(q) }\frac{\chi_{ \Co E} (y)-\chi_{E}(y)}{|q-y|^{n+s}}\, dy  \geq 0.
}
Moreover, from \cite[Lemma 3.1]{graph} (here applied with $\lambda =1$) we have that
\[\bigg| \int_{\Per_\delta \setminus B_{\varrho}(q)} \frac{\chi_{ \Co E} (y)-\chi_{E}(y)}{|q-y|^{n+s}}\, dy \bigg|\leq  \int_{\Per_\delta} \frac{dy}{|q-y|^{n+s}}\leq  \frac{C_0}{1-s} {\delta^{-s} }
 ,\] 
with $C_0=C_0(n)>0$.
Notice that $B_\delta(q)\subseteq K_\delta$  so
\bgs{\bigg | \int_{B_R(q)\setminus K_\delta } \frac{\chi_{ \Co E} (y)-\chi_{E}(y)}{|q-y|^{n+s}}\, dy\bigg | \leq &\;\int_{B_R(q)\setminus B_\delta(q)} \frac{dy}{|q-y|^{n+s}} = \varpi_n  \frac{\delta^{-s}- R^{-s}}{s} .} 
Therefore for every $\varrho <\delta/2$
one has that
\bgs{\label{CH:3:due}\int_{B_R(q)\setminus B_\varrho(q)} \frac{\chi_{ \Co E} (y)-\chi_{E}(y)}{|q-y|^{n+s}}\, dy \geq -\frac{C_0}{1-s}\delta^{-s} - \frac{\varpi_n}{s} \delta^{-s} + \frac{\varpi_n}{s} R^{-s}  .}
Thus, using \eqref{CH:3:tildebeta}
\bgs{
\I^\varrho_s&[E](q)=\int_{B_R(q)\setminus B_\varrho(q)} \frac{\chi_{ \Co E} (y)-\chi_{E}(y)}{|q-y|^{n+s}}\, dy
+ \int_{\Co B_R(q)} \frac{\chi_{ \Co E} (y)-\chi_{E}(y)}{|q-y|^{n+s}}\, dy \\
&\geq  -\frac{C_0}{1-s}\delta^{-s} - \frac{\varpi_n}{s} \delta^{-s} + \frac{\varpi_n}{s} R^{-s}
+ \int_{\Co B_R(q)} \frac{dy}{|q-y|^{n+s}} -  2\int_{ \Co B_R(q)} \frac{\chi_E(y)}{|q-y|^{n+s}}\,dy \\ 
&\geq - \delta^{-s} \Big(\frac{C_0}{1-s} +\frac{\varpi_n} s\Big) +  \frac{\varpi_n}s R^{-s}
+ \bigg(\frac{\varpi_n}s R^{-s}- 2\sup_{q\in \overline \Omega}
 \int_{ \Co B_R(q)} \frac{\chi_E(y)}{|q-y|^{n+s}}\,dy\bigg)\\
 &\geq - \delta^{-s} \Big(\frac{C_0}{1-s} +\frac{\varpi_n} s\Big) +  \frac{\varpi_n}s R^{-s}  +\frac{7\beta}{2s}\\
 &\geq  - \delta^{-s} \Big(2C_0 +\frac{\varpi_n} s\Big) +  \frac{\varpi_n}s R^{-s}  +\frac{7\beta}{2s},
 }
where we also exploited that $s<s'\leq 1/2$.
Since $R>1$, we have
\[R^{-s}\to1^-,\qquad\textrm{as }s\to0^+.\]
Therefore we can find $s''=s''(E_0,\Omega)$ small enough such that
\[\varpi_nR^{-s}\geq\varpi_n-\frac{\beta}{2},\qquad\forall s<s''.\]

Now let
\[s_0=s_0(E_0,\Omega):=\min\Big\{s',s'',\frac{\beta}{2C_0}\Big\}.\]
Then, for every $s<s_0$ we have
\eqlab{\label{CH:3:important_estimate_curv}\I_s^\varrho[E](q)&\geq\frac{1}{s}\Big\{-\delta^{-s}\big((2C_0)s+\varpi_n\big)+\varpi_n R^{-s}+\frac{7}{2}\beta\Big\}\\
&
\geq\frac{1}{s}\big\{-\delta^{-s}(\varpi_n+\beta)+\varpi_n+3\beta\big\},}
for every $\varrho\in(0,\delta/2)$.

Notice that if we fix $s\in(0,s_0)$, then for every
\[\delta\geq e^{-\frac{1}{s}\log\frac{\varpi_n+2\beta}{\varpi_n+\beta}}=:\delta_s(E_0),\]
we have that
\[-\delta^{-s}(\varpi_n+\beta)+\varpi_n+3\beta\geq\beta>0.\]
To conclude, we
let $\sigma\in(0,s_0)$ and suppose that $E$ has an exterior tangent ball of radius $\delta_\sigma$
at $q\in\partial E\cap\overline{\Omega}$.
Notice that, since $\delta_\sigma<1$, we have
\[-(\delta_\sigma)^{-s}(\varpi_n+\beta)+\varpi_n+3\beta\geq
-(\delta_\sigma)^{-\sigma}(\varpi_n+\beta)+\varpi_n+3\beta=\beta,\qquad\forall\,s\in(0,\sigma].\]
Then \eqref{CH:3:important_estimate_curv} gives that
\[\liminf_{\varrho\to0^+}\I_s^\varrho[E](q)\geq\frac{\beta}{s}>0,\qquad\forall\,s\in(0,\sigma],\]
which concludes the proof.
\end{proof}

\begin{remark}
We remark that
\[\log\frac{\varpi_n+2\beta}{\varpi_n+\beta}>0,\]
thus
\[\delta_s\to0^+\qquad\textrm{as }s\to 0^+.\]
\end{remark}


As a consequence of Theorem \ref{CH:3:positivecurvature}, we have that, as $s\to0^+$,
the $s$-minimal sets with small mass at infinity have small mass in $\Omega$.
The precise result goes as follows:

\begin{corollary}
Let $\Omega\subseteq\Rn$ be a bounded open set,
let $E\subseteq\Rn$ be such that
\[\overline \alpha(E) <\frac{\varpi_n}2,\]
and suppose that $\partial E$ is of class $C^2$ in $\Omega$.
Then, for every $\Omega'\Subset\Omega$ there exists
$\tilde s=\tilde s(E\cap\overline{\Omega'})\in(0,s_0)$ such that for every $s\in(0,\tilde s]$
\eqlab{\I_s[E](q)\geq\frac{\varpi_n -2\overline \alpha(E)}{4s}>0,\qquad\forall\,q\in\partial E\cap\overline{\Omega'}.}
\end{corollary}

\begin{proof}
Since $\partial E$ is of class $C^2$ in $\Omega$ and $\Omega'\Subset\Omega$, the set $E$ satisfies a uniform exterior
ball condition of radius $\tilde\delta=\tilde\delta(E\cap\overline{\Omega'})$ in $\overline{\Omega'}$, meaning that
$E$ has an exterior tangent ball of radius at least $\tilde\delta$ at every point $q\in\partial E\cap\overline{\Omega'}$.

Now, since $\delta_s\to0^+$ as $s\to0^+$, we can find $\tilde s=\tilde s(E\cap\overline{\Omega'})
<s_0(E\setminus\Omega,\Omega)$,
small enough such that $\delta_s<\tilde\delta$ for every $s\in(0,\tilde s]$. Then we can conclude by applying Theorem \ref{CH:3:positivecurvature}.
\end{proof}

\subsection{Classification of $s$-minimal surfaces (Theorem \ref{CH:3:THM})}\label{CH:3:alternative}

To classify the behavior of the $s$-minimal surfaces when $s$ is small,
we need to take into account the ``worst case scenario'',
that is the one in which the set behaves very badly in terms
of oscillations and lack of regularity. To this aim,
we make an observation about $\delta$-dense sets.

   \begin{center}
\begin{figure}[htpb]
	\hspace{0.79cm}
	\begin{minipage}[b]{0.79\linewidth}
	\centering
	\includegraphics[width=0.79\textwidth]{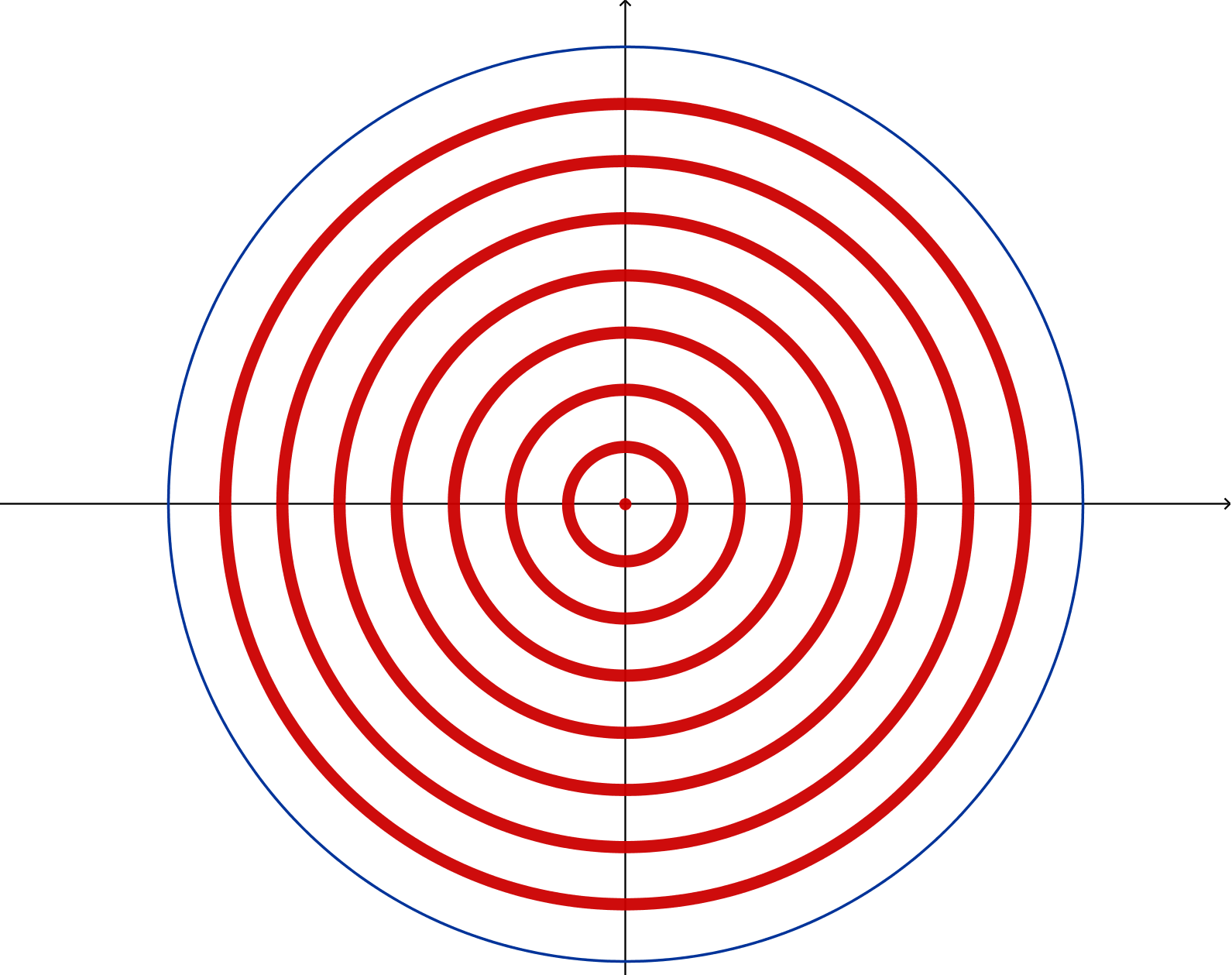}
	\caption{A $\delta$-dense set of measure $<\eps$}   
	\end{minipage}
\end{figure} 
	\end{center}
  \begin{remark}\label{CH:3:deltadance}
  For every $k\ge1$ and every $\eps<2^{-k-1}$, we define the sets
\[\Gamma_k^\eps:=B_\eps\cup\bigcup_{i=1}^{2^k-1}\Big\{x\in\R^n\,\big|\,\frac{i}{2^k}-\eps<|x|<\frac{i}{2^k}+\eps\Big\}
\quad\textrm{and}\quad\Gamma_k:=\{0\}\cup\bigcup_{i=1}^{2^k-1}\partial B_\frac{i}{2^k}.\]
%
\noindent Notice that for every $\delta>0$ there exists $\tilde k=\tilde k(\delta)$ such that for every $k\geq\tilde k$ we have
\[B_\delta(x)\cap\Gamma_k\not=\emptyset,\qquad\forall\,B_\delta(x)\subseteq B_1.\]
Thus, for every $k\geq\tilde k(\delta)$ and $\eps<2^{-k-1}$, the set $\Gamma_k^\eps$ is $\delta$-dense in $B_1$.
\noindent Moreover, notice that
\[\Gamma_k=\bigcap_{\eps\in(0,2^{-k-1})}\Gamma_k^\eps\quad\textrm{ and }\quad\lim_{\eps\to0^+}|\Gamma_k^\eps|=0.\]
It is also worth remarking that the sets $\Gamma_k^\eps$ have smooth boundary.
\noindent In particular, for every $\delta>0$ and every $\eps>0$ small, we can find a set $E\subseteq B_1$ which is $\delta$-dense in $B_1$ and whose measure is $|E|<\eps$.
This means that we can find an open set $E$ with smooth boundary, whose measure is arbitrarily small
and which is ``topologically arbitrarily dense'' in $B_1$.
\end{remark}

We introduce the following useful geometric observation. 

\begin{prop}\label{CH:3:tgball}
Let $\Omega\subseteq \Rn$ be a bounded and connected open set with $C^2$ boundary and let $\delta\in(0,r_0)$,
for $r_0$ given in \eqref{CH:3:r01}. 
If $E$ is not $\delta$-{dense} in $\Omega$ 
and $|E\cap\Omega|>0$, 
then there exists a point $q\in\partial E\cap\Omega$ such that $E$ has an exterior tangent ball at $q$ of radius $\delta$ (contained in $\Omega$),
i.e.
there exist $p\in \Co E\cap \Omega$ such that
\[ B_\delta(p)\Subset\Omega,\qquad q \in \partial B_{\delta}(p) \cap \partial E
\quad \mbox{ and }  \quad B_{ \delta}(p) \subseteq \Co E.\] 
\end{prop}
\begin{proof}
Using Definition \ref{CH:3:wild}, we have that there exists  $x\in \Omega$ for which $B_\delta(x)\Subset \Omega$ and $|B_\delta(x)\cap E|=0$, so $B_\delta(x)\subseteq E_{ext}$. If $B_\delta(x)$ is tangent to $\partial E$ then we are done.

Notice that
\[B_\delta(x)\Subset\Omega\quad\Longrightarrow\quad d(x,\partial\Omega)>\delta,\]
and let
\[\delta':=\min\{r_0,d(x,\partial\Omega)\}\in(\delta,r_0].\]
Now we consider the open set $\Omega_{-\delta'}\subseteq \Omega$
\[ \Omega_{-\delta'} :=\{ \bar d_{\Omega}<-\delta'\},\]  so $x\in \Omega_{-\delta'}$. According to Remark \ref{CH:3:c21} and Lemma \ref{CH:3:geomlem} we have that $\Omega_{-\delta'}$ has $C^2$ boundary and that 
\eqlab{\label{CH:3:unifball} \Omega_{-\delta'} \mbox{ satisfies the uniform interior ball condition of radius at least } r_0.} 
 
 We have two possibilities:
 \eqlab{\label{CH:3:inside} &\mbox{ i) } && \overline  E\cap \Omega_{-\delta'}\neq \emptyset \\ 
 				  &\mbox{ ii) }&&\emptyset \neq  \overline{E}\cap\Omega \subseteq \Omega \setminus \Omega_{-\delta'}.} 
 
 If i) happens, we pick any point $y\in \overline{E}\cap \Omega_{-\delta'}$.  
The set $\overline{\Omega_{-\delta'}}$ is path connected (see Proposition \ref{CH:3:retract}), so there exists a path $c:[0,1]\longrightarrow
\R^n$ that connects $x$ to $y$ and that stays inside $\overline{\Omega_{-\delta'}}$, that is
\[c(0)=x,\qquad c(1)=y\quad\textrm{ and }\quad c(t)\in\overline{\Omega_{-\delta'}},\quad\forall\,t\in[0,1].\]
Moreover, since $\delta<\delta'$, we have
\[B_\delta\big(c(t)\big)\Subset\Omega\qquad\forall\,t\in[0,1].\]
Hence,
we can ``slide the ball'' $B_{\delta}(x)$ along the path and we obtain the desired
claim thanks to Lemma \ref{CH:3:slidetheballs}.

%

Now, if we are in the case ii) of \eqref{CH:3:inside}, then $\Omega_{-\delta'}\subseteq E_{ext}$, so we dilate $\Omega_{-\delta'}$ until we first touch 
$\overline E$. That is, we consider 
\[\tilde \varrho:=\inf\{\varrho\in[0,\delta'] \; \big| \; \Omega_{-\varrho}\subseteq E_{ext}\}.\]
Notice that by hypothesis $\tilde \varrho>0$. Then 
\[\overline{\Omega_{-\tilde \varrho}} \subseteq \overline {E_{ext}}=E_{ext}\cup\partial E. \]
If
\[ \partial \Omega_{-\tilde\varrho} \cap \partial E=\emptyset\quad  \mbox { then }\quad  \overline{ \Omega_{-\tilde \varrho}} \subseteq  E_{ext},\] 
hence we have that
\[ d= d\left(  \overline E\cap \Omega\setminus \Omega_{-\delta'}, \overline {\Omega_{-\tilde \varrho}}\right)\in(0,\tilde \varrho),\]
therefore
\[ \Omega_{-\tilde \varrho} \subseteq \Omega_{-(\tilde \varrho-d)}\subseteq E_{ext}.\]
This is in contradiction with the definition of $\tilde \varrho$. Hence, 
 there exists $q \in \partial \Omega_{-\tilde \varrho} \cap \partial E$. 

Recall that, by definition of $\tilde\varrho$, we have
$\Omega_{-\tilde \varrho} \subseteq \Co E.$ Thanks to \eqref{CH:3:unifball}, there exists a tangent ball at $q$ interior to $\Omega_{-\tilde \varrho}$, hence a tangent ball at $q$ exterior to $E$, of radius at least $r_0>\delta$. 
This concludes the proof of the lemma.
\end{proof}

We observe that part $(A)$ of Theorem \ref{CH:3:THM}
is essentially a consequence of Theorem \ref{CH:3:positivecurvature}.
Indeed, if an $s$-minimal set $E$ is not $\delta_s$-dense and it is not empty in $\Omega$, then
by Proposition \ref{CH:3:tgball} we can find a point $q\in\partial E\cap\Omega$ at which $E$ has an exterior tangent ball of radius $\delta_s$. Then Theorem \ref{CH:3:positivecurvature} implies that the $s$-fractional mean curvature of $E$ in $q$ is strictly positive,
contradicting the Euler-Lagrange equation.

On the other hand, part $(B)$ of Theorem \ref{CH:3:THM} follows from a careful asymptotic use
of the density estimates provided by Theorem \ref{CH:3:notfull}.
For the reader's facility, we also recall that $r_0$ 
has the same meaning here and across the chapter, 
as clarified in Appendix \ref{CH:3:A2}.
We now proceed with the precise arguments of the proof.

\begin{proof}[Proof of Theorem \ref{CH:3:THM}]
We begin by proving part $(A)$.\\
First of all, since $\delta_s\to0^+$, we can find $s_1=s_1(E_0,\Omega)\in(0,s_0]$ such that $\delta_s<r_0$
for every $s\in(0,s_1)$.

Now let $s\in(0,s_1)$ and let $E$ be $s$-minimal in $\Omega$, with exterior data $E_0$.

We suppose that  $E\cap \Omega\neq \emptyset$ and prove that $E$ has to be $\delta_s$-dense. 

Suppose by contradiction that $E$ is not $\delta_s$-dense. Then, in view of 
Proposition \ref{CH:3:tgball}, there exists $p\in \Co E\cap \Omega$ such that
\[ q \in \partial B_{\delta_s}(p) \cap (\partial E\cap\Omega) \quad \mbox{ and }  \quad B_{ \delta_s}(p) \subseteq \Co E.\] 
Hence we use the Euler-Lagrange theorem at $q$, i.e.
\[\I_s[E](q) \leq 0,\] 
to obtain a contradiction with Theorem \ref{CH:3:positivecurvature}. This says that $E$ is not $\delta_s$-dense and concludes the proof of part $(A)$ of Theorem \ref{CH:3:THM}.

Now we prove the part $(B)$ of the Theorem.\\
Suppose that point $(B.1)$ does not hold true. Then we can find a sequence $s_k\searrow0$ and a sequence of sets $E_k$
such that each $E_k$ is $s_k$-minimal in $\Omega$ with exterior data $E_0$ and
\[E_k\cap\Omega\not=\emptyset.\]
We can assume that $s_k<s_1(E_0,\Omega)$ for every $k$. Then part $(A)$ implies that each $E_k$ is $\delta_{s_k}$-dense,
that is
\[|E_k\cap B_{\delta_{s_k}}(x)|>0\quad\forall\,B_{\delta_{s_k}}(x)\Subset\Omega.\]
Fix $\gamma=\frac{1}{2}$, take a sequence $\delta_h\searrow0$ and let $\sigma_{\delta_h,\frac{1}{2}}$ be as in Theorem \ref{CH:3:notfull}. Recall that $\delta_s\searrow0$ as $s\searrow0$.
Thus for every $h$ we can find $k_h$ big enough such that
\eqlab{\label{CH:3:koala}s_{k_h}<\sigma_{\delta_h,\frac{1}{2}}\qquad\textrm{and}\qquad\delta_{s_{k_h}}<\delta_h.}
In particular, this implies
\eqlab{\label{CH:3:koala1}|E_{k_h}\cap B_{\delta_h}(x)|\ge|E_k\cap B_{\delta_{s_{k_h}}}(x)|
>0\quad\forall\,B_{\delta_h}(x)\Subset\Omega,}
for every $h$. On the other hand, by \eqref{CH:3:koala} and Theorem \ref{CH:3:notfull}, we also have that
\eqlab{\label{CH:3:koala2}|\Co E_{k_h}\cap B_{\delta_h}(x)|>0\quad\forall\,B_{\delta_h}(x)\Subset\Omega.}
This concludes the proof of part $(B)$. Indeed, notice that since $B_{\delta_h}(x)$ is connected, \eqref{CH:3:koala1} and \eqref{CH:3:koala2}
together imply that
\[\partial E_{k_h}\cap B_{\delta_h}(x)\neq\emptyset\quad\forall\,B_{\delta_h}(x)\Subset\Omega.\]
\end{proof}

\subsection{Stickiness to the boundary is a typical behavior (Theorem \ref{CH:3:boundedset})}\label{CH:3:sticky}

Now we show that the ``typical behavior''
of the nonlocal minimal surfaces
is to stick at the boundary whenever
they are allowed to do it,
in 
the precise
sense given 
by Theorem \ref{CH:3:boundedset}.

\begin{proof}[Proof of Theorem \ref{CH:3:boundedset}]
Let
\[\delta:=\frac{1}{2}\min\{r_0,R\},\]
and notice that (see Remark \ref{CH:3:ext_unif_omega})
\[B_\delta(x_0+\delta\nu_\Omega(x_0))\subseteq B_R(x_0)\setminus\Omega\subseteq\Co E_0.\]
Since $\delta_s\to0^+$, we can find $s_3=s_3(E_0,\Omega)\in(0,s_0]$ such that $\delta_s<\delta$ for every $s\in(0,s_3)$.

Now let $s\in(0,s_3)$ and let $E$ be $s$-minimal in $\Omega$, with exterior data $E_0$.

We claim that
\eqlab{\label{CH:3:pf_bdedset_eq1}B_\delta(x_0-r_0\nu_\Omega(x_0))\subseteq E_{ext}.}
We observe that this is indeed a crucial step to
prove Theorem \ref{CH:3:boundedset}. Indeed, once this is established,
by Remark \ref{CH:3:ext_unif_omega} we obtain that
\[B_\delta(x_0-r_0\nu_\Omega(x_0))\Subset\Omega.\]
Hence, since $\delta_s<\delta$, we deduce from \eqref{CH:3:pf_bdedset_eq1}
that $E$ is not $\delta_s$-dense.
Thus, since $s<s_3\leq s_1$, Theorem \ref{CH:3:THM} 
implies that $E\cap\Omega=\emptyset$, which concludes the proof of Theorem \ref{CH:3:boundedset}.

This, we are left to prove \eqref{CH:3:pf_bdedset_eq1}. Suppose by contradiction that
\[\overline{E}\cap B_\delta(x_0-r_0\nu_\Omega(x_0))\not=\emptyset,\]
and consider the segment $c:[0,1]\longrightarrow\Rn$,
\[c(t):=x_0+\big((1-t)\delta-t\,r_0\big)\nu_\Omega(x_0).\]
Notice that
\[B_\delta\big(c(0)\big)\subseteq E_{ext}\quad\mbox{ and }\quad B_\delta\big(c(1)\big)\cap\overline{E}\not=\emptyset,\]
so
\[t_0:=\sup\Big\{\tau\in[0,1]\,\big|\,\bigcup_{t\in[0,\tau]}B_\delta\big(c(t)\big)\subseteq E_{ext}\Big\}<1.\]
Arguing as in Lemma \ref{CH:3:slidetheballs}, we conclude that
\[B_\delta\big(c(t_0)\big)\subseteq E_{ext}\quad\mbox{ and }\quad \exists\,q\in\partial B_\delta\big(c(t_0)\big)\cap\partial E.\]
By definition of $c$, we have that either $q\in\Omega$ or
\[q\in\partial\Omega\cap B_R(x_0).\]
In both cases (see \cite[Theorem 5.1]{CRS10} and Theorem \eqref{CH:3:EL_boundary_coroll}) we have
\[\I_s[E](q)\leq0,\]
which gives a contradiction with Theorem \ref{CH:3:positivecurvature} and concludes the proof.
\end{proof}

 
   \section{The contribution from infinity of some supergraphs} \label{CH:3:sectexamples} We compute  in this Subsection the contribution from infinity of some particular supergraphs.
 \begin{example}[The cone] \label{CH:3:THECONE} Let $ S \subseteq  \mathbb S^{n-1}$ be a portion of the unit sphere, $\mathfrak o:=\Ha^{n-1}(S)$ and
  \[ C:=\{ t\sigma  \; \big| \; t\geq 0, \;\sigma\in  S)\}.\] 
  Then the contribution from infinity is given by the opening of the cone,
  \eqlab{\label{CH:3:cony} \alpha(C) = \mathfrak o.}
  Indeed,
  \[\alpha_s(0,1,C)= \int_{\Co B_1} \frac{\chi_C(y)}{|y|^{n+s}}\, dy = \Ha^{n-1}(S) \int_1^{\infty} t^{-s-1}\, dt=
  \frac{\mathfrak o}s,  \]
and we obtain the claim by passing to the limit. Notice that this says in particular that the contribution from infinity of a half-space is $\varpi_n/2$.
  \end{example} 
  
  \begin{center}
\begin{figure}[htpb]
	\hspace{0.99cm}
	\begin{minipage}[b]{0.99\linewidth}
	\centering
	\includegraphics[width=0.99\textwidth]{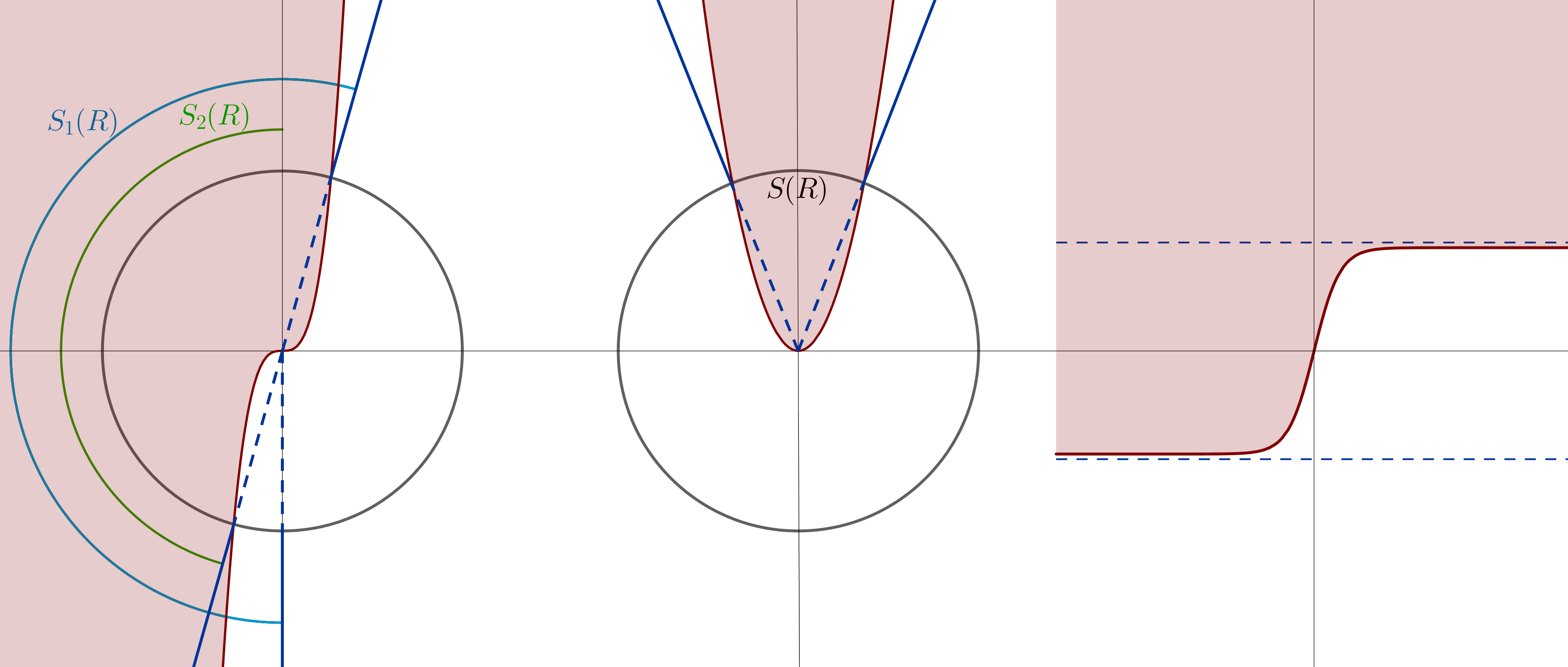}
	\caption{The contribution from infinity of $x^3$, $x^2$ and $\tanh x$}   
	\label{CH:3:x3}
	\end{minipage}
\end{figure} 
	\end{center}
	  \begin{example}[The parabola]
  We consider the supergraph 
  \[ E:=\{ (x',x_n) \; \big| \; x_n\geq |x'|^2\},\] and we show that, in this case, \[ \alpha(E)=0.\]
  In order to see this, we take any $R>0$, intersect the ball $B_R$ with the parabola and build a cone on this intersection (see the second picture in Figure \ref{CH:3:x3}), i.e. we define
  \[ S(R):=\partial B_R \cap E , \quad \quad C_R =\{ t\sigma  \; \big| \; t\geq 0, \;\sigma\in  S(R)\}.\] 
  We can explicitly compute the opening of this cone, that is
  \[ \mathfrak o (R)= \Bigg( \arcsin \frac{\left(\sqrt{4R^2+1}-1\right)^{1/2}}{R\, \sqrt 2 }\Bigg)\frac{\varpi_n}{\pi} .\]
  Since $E\subseteq C_R$ outside of $B_R$, thanks to the monotonicity property in Proposition \ref{CH:3:subsetssmin} and to \eqref{CH:3:cony}, we have that
  \[  \overline \alpha (E) \leq \overline \alpha(C_R)  =\mathfrak o(R).\]
  Sending $R\to \infty$, we find that
\[  \overline \alpha (E) =0, \quad \mbox{ thus } \quad \alpha(E)=0.\]
  \end{example}
  More generally, if we consider for any given $c, \eps>0$ a function $u$ such that
  \[ u(x')>c|x'|^{1+\eps}, \quad \mbox{ for  any }|x'|>R \mbox{ for  some } R>0\]
  and\[   E:=\{ (x',x_n) \; \big| \; x_n\geq u(x')\},\]
  then
\[ \alpha(E)=0.\]
On the other hand, if we consider a function that is not rotation invariant, things can go differently, as we see in the next example. 
  
   \begin{example}[The supergraph of $x^3$]
 We consider the supergraph 
  \[ E:=\{ (x,y) \; \big| \; y\geq x^3\}.\] 
In this case, we show that \[\alpha(E) = \pi .\]
For this, given $R>0$, we intersect  $\partial B_R$ with $E$  and denote by $S_1(R)$ and $S_2(R)$ the arcs on the circle as the first picture in Figure \ref{CH:3:x3}. We consider
the cones\[ C^1_R :=\{ t\sigma  \; \big| \; t\geq 0, \;\sigma\in  S_1(R)\}\, \quad  C^2_R :=\{ t\sigma  \; \big| \; t\geq 0, \;\sigma\in  S_2(R)\} \] and
notice that outside of $B_R$, it holds that $C_R^2\subseteq E\subseteq C_R^1$. 
Let  $\overline x_R $ be the solution of \[ x^6+x^2=R^2,\] that is the $x$-coordinate in absolute value of the intersection points $\partial B_R \cap \partial E$. Since $f(x)=x^6+x^2$ is increasing on $(0,\infty)$ and $ R^2=f(\overline x_R)<f(R^{1/3}),$ we have that  $\overline x_R< R^{1/3}$. Hence  
\[ \mathfrak o^1(R)= \pi + \arcsin \frac{\overline x_R}R  \leq \pi+  \arcsin \frac{R^{1/3}}R  ,\quad   \mathfrak o^2(R)  \geq \pi- \arcsin \frac{R^{1/3}}R .\] 
Thanks to the monotonicity property in Proposition \ref{CH:3:subsetssmin} and to \eqref{CH:3:cony} we have that
  \[  \overline \alpha (E) \leq\alpha(C_R^1)=  \mathfrak o^1(R), \quad \underline \alpha(E) \geq \alpha(C_R^2)= \mathfrak o^2(R)  \]
  and sending  $R\to \infty$ we obtain that
\[ \overline \alpha(E) \leq  \pi , \quad \underline \alpha(E) \geq \pi.\] Thus $\alpha(E)$ exists and we obtain the desired conclusion.
 \end{example}

 \begin{example}[The supergraph of a bounded function]\label{CH:3:tanh}
 We consider the supergraph 
  \[ E:=\{ (x',x_n) \; \big| \; x_n\geq u(x') \}, \quad \quad{\mbox{with}}
\quad \quad \|u\|_{L^\infty(\Rn)} <M.\] 
We show that, in this case,
\[\alpha(E) = \frac{\varpi_n}2 .\]
To this aim, let \bgs{ &\mathfrak{P}_1:=\{ (x',x_n)\; \big| \; x_n>M\} \\
	&\mathfrak{P}_2:=\{ (x',x_n)\; \big| \; x_n<-M\}.}
	We have that 
 \[   \mathfrak{P}_1 \subseteq E , \quad \quad \mathfrak{P}_2 \subseteq \Co E . \] 
  Hence by Proposition \ref{CH:3:subsetssmin}
  \[ \underline \alpha(E) \geq \overline \alpha (\mathfrak{P}_1)=\frac{\varpi_n}2,\quad  \quad  
  \underline \alpha(\Co E)\geq \overline\alpha (\mathfrak{P}_2)=\frac{\varpi_n}2 .\]
  Since $\underline \alpha(CE) =\varpi_n - \overline \alpha(E)$ we find that
   \[ \overline \alpha(E)\leq \frac{\varpi_n}2,\]
  thus the conclusion. An example of this type is depicted in Figure \ref{CH:3:x3} (more generally, the result holds for the supergraph in $\Rn$  
 $ \{ (x',x_n) \; \big| \; x_n\geq \tanh x_1\}$).
 \end{example}

 \begin{example}[The supergraph of a sublinear graph]\label{CH:3:candygr}
   More generally, we can take the supergraph of a function that grows sublinearly at infinity, i.e.  
  \[  E:=\{(x',x_n)\;\big|\; x_n>u(x')\}, \qquad{\mbox{with}}\qquad 
\lim_{|x'|\to +\infty} \frac{|u(x')|}{|x'|}=0 .\]
In this case, we show that \[ \alpha(E)= \frac{\varpi_n}2.\]
Indeed, for any $\eps>0$ we have that there exists $R=R(\eps)>0$ such that
\[ |u(x')|<\eps |x'|, \quad \forall \;|x'|>R.\]
We denote
\[ S_1(R):= \partial B_R \cap \{(x',x_n)\;\big|\; x_n>\eps|x'|\}, \qquad S_2(R):= \partial B_R \cap \{(x',x_n)\;\big|\; x_n<-\eps|x'|\}\]
and
\[ C_R^i=\{ t\sigma\; \big| \; t\geq 0, \; \sigma \in S_i(R)\}, \quad \mbox{for }\; i=1,2.\]
We have that outside of $B_R$
\[ C_R^1\subseteq E, \qquad C_R^2\subseteq \Co E ,\] and 
\[ \alpha (C_R^1)= \alpha (C_R^2)= \frac{ \varpi_n}{\pi} \left( \frac{\pi}2- \arctan \eps \right).\]
We use Proposition \ref{CH:3:subsetssmin}, (i), and letting $\eps$ go to zero, we obtain that $\alpha(E)$ exists and 
\[ \alpha (E)= \frac{\varpi_n}2.\]
 A particular example of this type is given by
 \[ E:=\{(x',x_n)\;\big|\; x_n>c|x'|^{1-\eps} \}, \quad \mbox{ when } |x'|>R\; \;  \mbox{ for some } \eps\in(0,1],\,c\in\R,\, R>0     .\]
 \end{example}
 In particular using the additivity property in Proposition \ref{CH:3:subsetssmin} we can compute $\alpha$ for sets that lie between two graphs. 

  \begin{center}
\begin{figure}[htpb]
	\hspace{0.89cm}
	\begin{minipage}[b]{0.89\linewidth}
	\centering
	\includegraphics[width=0.89\textwidth]{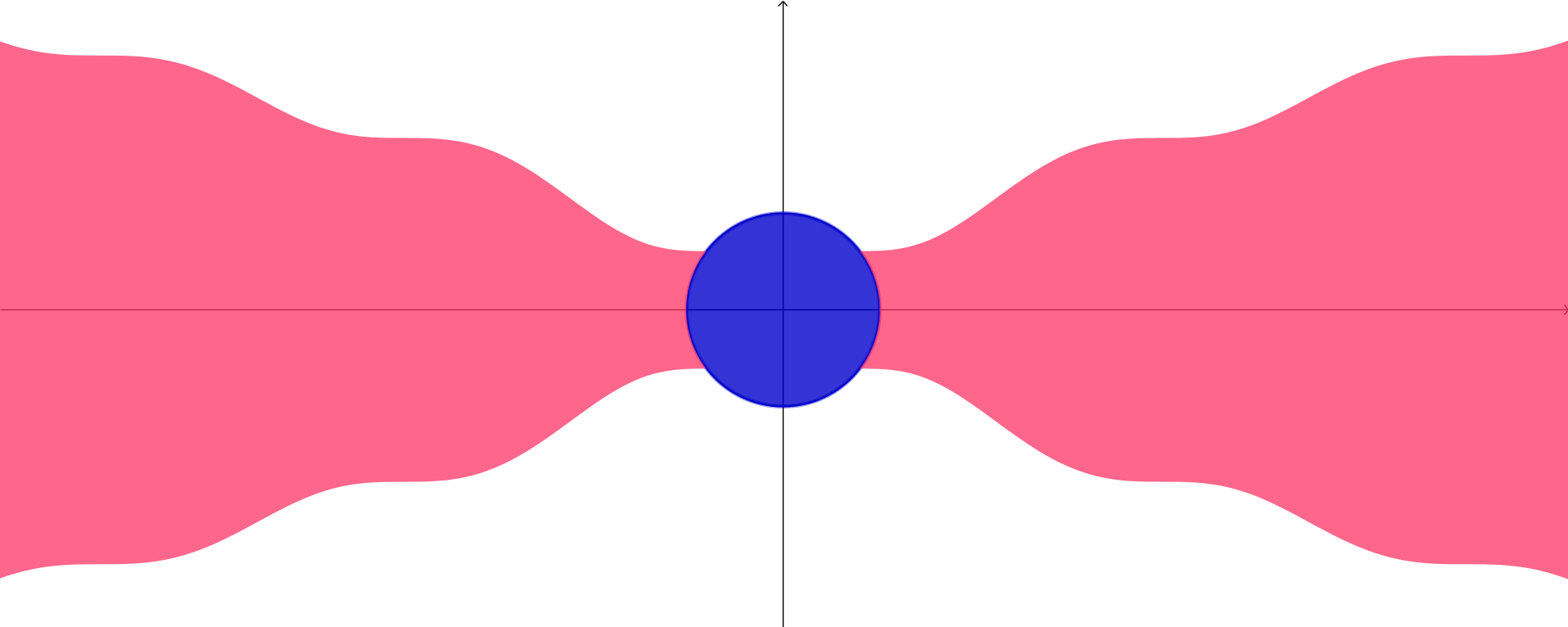}
	\caption{The ``butterscotch hard candy'' graph}   
	\label{CH:3:candy_pic}
	\end{minipage}
\end{figure} 
	\end{center}
	 \begin{example}[The ``butterscotch hard candy'']
Let~$E\subseteq\R^n$ be such that
\bgs{
E\cap\{|x'|>R\}\subseteq\{(x',x_n)\;\big|\; |x'|>R\;,\;\,|x_n|<c|x'|^{1-\eps} \},
}
for some $\eps\in(0,1],\,c>0$ and $R>0$
(an example of such a set $E$ is given in Figure \ref{CH:3:candy_pic}).
In this case, we have that \[\alpha(E)=0.\]
Indeed, we can write $E_1:=E\cap\{|x'|>R\}$ and $E_2:=
E\cap\{|x'|\le R\}$. Then,
using
 the computations in Example \ref{CH:3:candygr}, we have by
the monotonicity and the additivity properties in Proposition \ref{CH:3:subsetssmin} that
\[ \overline{\alpha}(E_1) 
\le \alpha\big(\{x_n>-c|x'|^{1-\eps}\}\big)-\alpha\big(\{x_n>c|x'|^{1-\eps}\}\big)=0.\]
Moreover, $E_2$ lies inside $\{|x_1|\le R\}$. Hence, again by Proposition \ref{CH:3:subsetssmin} and by Example \ref{CH:3:THECONE},
we find
\[\overline{\alpha}(E_2)\le\alpha\big(\{|x_1|\le R\}\big)=\alpha\big(\{x_1\le R\}\big)-\alpha\big(\{x_1<-R\}\big)=0.\]
Consequently, using again the additivity property in Proposition \ref{CH:3:subsetssmin}, we obtain that
$$ \overline\alpha(E)\le \overline\alpha(E_1)+\overline\alpha(E_2)
=0,$$
that is the desired result.
\end{example}
We can also compute $\alpha$ for sets that have different growth ratios in different directions. For this, we have the following example.
   \begin{example}[The supergraph of a superlinear function on a small cone]
  We consider a set lying in the half-space, deprived of a set that grows linearly at infinity.  We denote by $\tilde S$ the portion of the sphere given by 
 \bgs{ \tilde{S}:=\Big\{\sigma \in \mathbb S^{n-2}\,\Big|\,\sigma=(&\cos\sigma_1, \sin \sigma_1\cos\sigma_2,\dots, \sin \sigma_1\dots\sin \sigma_{n-2}),\\ &  \mbox{ with } \sigma_i\in \lr{\frac{\pi}2-\bar \eps, \frac{\pi}2+\bar \eps},\;  i=1,\dots,n-2\Big\} ,} where $\overline \eps \in (0,{\pi}/{2})$.
 For $x_0\in \Rn$ and $k>0$ we define  the supergraph $E\subseteq \Rn$ as
\bgs{\label{CH:3:sigma}  E:=\big\{ (x',x_n)\in \Rn \; &\big| \; x_n\geq u(x')\big \} \quad   \mbox{ where }  \quad
 u(x')=\alig{ & \, k|x'-x_0'| &\mbox{ for } &x'\in X, \\
		 & \, 0 &\mbox{ for } &x'\notin X,  }\\
	  &X= \{ x'   \in \R^{n-1}  \mbox{ s.t. }x' = t\sigma +x_0',\,  \sigma \in  \tilde S\}.   }
We remark that $X\subseteq \{x_n=0\}$ is the cone ``generated'' by $\tilde S$ 
and centered at $x_0$.	 Then
\eqlab{\label{CH:3:alphasigma} \alpha(E) = \frac{\varpi_n}2 -  \mathcal{H}^{n-2}(\tilde S)  \int_{0}^k \frac{dt}{(1+t^2)^{\frac{n}2}}.}
Let 
\[ \mathfrak{P}_+
:=\{ (x',x_n) \; \big| \; x_n>0\}, \quad \quad \mathfrak{P}_-:=\{ (x',x_n) \; \big| \; x_n<0\}\] and we consider
the subgraph
\[ F:= \big\{ (x',x_n)\; \big| \; 0<x_n< u(x')\big \}.\]
Then \[ E \cup F=\mathfrak{P}_+, \quad \quad \mathfrak{P}_-\cup F = \Co E.\] 
Using the additivity property in Proposition \ref{CH:3:subsetssmin}, we see that
\eqlab{\label{CH:3:ah123}  \overline \alpha(E)\geq \frac{\varpi_n}2-\overline \alpha(F), \quad \quad  \varpi_n - \underline \alpha(E) = \overline \alpha(\Co E)\leq \frac{\varpi_n}2 +\overline \alpha(F).} 
Let $R>0$ be arbitrary. We get that
\[\alpha_s(x_0,R,F) \leq  \int_{\left(B'_R(x'_0)\times \R\right)\cap \Co B_R(x_0)}\frac{\chi_{F}(y)}{|y-x_0|^{n+s}}\, dy + \int_{\Co \left( B'_R(x'_0)\times \R\right)}\frac{\chi_{F}(y)}{|y-x_0|^{n+s}}\, dy\]
so
\eqlab{\label{CH:3:i123}  \alpha_s(x_0,R,F) \leq &\;  \int_{B'_R(x'_0)}\frac{dy'}{|y'-x_0'|^{n-1+s}}\int_{\frac{\sqrt{R^2-|y'-x_0'|^2}}{|y'-x_0'|}}^{\infty} \frac{dt}{(1+t^2)^{\frac{n+s}2}}
   \\ &\; + 
   \int_{\Co B'_R(x'_0)\cap X}\frac{dy'}{|y'-x_0'|^{n-1+s}}\int_{0}^{k} \frac{dt}{(1+t^2)^{\frac{n+s}2}}\\=&\; I_1+I_2.} 
   Using that $1+t^2\geq \max\{1,t^2\}$ and passing to polar coordinates, we obtain that
   \bgs{ I_1= &\; \int_{B'_R(x_0')}\frac{dy'}{|y'-x_0'|^{n-1+s}}\bigg(\int_{\frac{\sqrt{R^2-|y'-x_0'|^2}}{|y'-x_0'|}}^{\frac{R}{|y'-x_0'|}} \frac{dt}{(1+t^2)^{\frac{n+s}2}} + \int_{\frac{R}{|y'-x_0'|}}^{\infty}\frac{dt}{(1+t^2)^{\frac{n+s}2}}\bigg)\\
   \leq&\; \varpi_{n-1} \bigg( \int_0^R \tau^{-s-2}   \left(R- \sqrt{R^2-\varrho^2}\right)\, d \varrho  +  \frac{R^{-n-s+1}}{n+s-1} \int_0^R \varrho^{n-2}\, d \varrho    \bigg)\\
   =&\;   \varpi_{n-1} \bigg( R^{-s} \int_0^1 \tau^{-s-2}   \left(1- \sqrt{1-\tau^2}\right)\, d \tau  +  \frac{R^{-s}}{(n+s-1)(n-1)}  \bigg)  .}
   Also, for any $\tau \in (0,1)$ we have that
\[1 -  \sqrt{1-\tau^2 }\leq c \tau^2,\] for some positive constant $c$, independent on $n,s$.
Therefore 
\bgs{ I_1\leq  \frac{ c \varpi_{n-1} R^{-s} }{1-s} + \frac{ \varpi_{n-1} R^{-s} }{(n-1)(n+s-1)} .}
Moreover, 
\[ I_2 =   \mathcal{H}^{n-2}(\tilde S) \frac{R^{-s}}s \int_{0}^k \frac{dt}{(1+t^2)^{\frac{n+s}2}}.\]
So
passing to limsup and liminf as $s\to 0^+$ in \eqref{CH:3:i123} and using Fatou's lemma we obtain that
\[\overline \alpha(F)\leq   \mathcal{H}^{n-2}(\tilde S)\int_{0}^k \frac{dt}{(1+t^2)^{\frac{n}2}},\quad \quad \underline \alpha(F) \geq  \mathcal{H}^{n-2}(\tilde S)\int_{0}^k \frac{dt}{(1+t^2)^{\frac{n}2}}.\] In particular $\alpha(F)$ exists, and from \eqref{CH:3:ah123} we get that 
\[ \frac{\varpi_n}2 -\alpha(F) \leq \underline \alpha(E)\leq \overline \alpha(E)\leq\frac{\varpi_n}2 -\alpha(F).\] 
Therefore, $\alpha(E)$ exists and 
\[  \alpha(E) =  \frac{\varpi_n}2 -  \mathcal{H}^{n-2}(\tilde S)  \int_{0}^k \frac{dt}{(1+t^2)^{\frac{n}2}} .\] 
 \end{example}

\section[Continuity of the fractional mean curvature]{Continuity of the fractional mean curvature and a sign changing
property of the nonlocal mean curvature}\label{CH:3:cont}

We use a formula proved in \cite{regularity} to show that the $s$-fractional mean curvature is continuous
with respect to $C^{1,\alpha}$ convergence of sets, for any $s<\alpha$ and with respect to $C^2$ convergence of sets, for $s$ close to 1.

By $C^{1,\alpha}$ convergence of sets we mean that our sets locally converge in measure and can locally be described as the supergraphs
of functions which converge in $C^{1,\alpha}$. 

\begin{defn}\label{CH:3:convofsets}
Let $E\subseteq\R^n$ and let $q\in\partial E$ such that $\partial E$ is $C^{1,\alpha}$ near $q$, for some $\alpha\in(0,1]$. We say that the sequence $E_k\subseteq\R^n$ converges to $E$ in a $C^{1,\alpha}$ sense (and write $E_k\xrightarrow{C^{1,\alpha}}E$) in a neighborhood of $q$ if:\\
(i)  the sets $E_k$ locally converge in measure to $E$, i.e.
\[|(E_k\Delta E)\cap B_r|\xrightarrow{k\to\infty}0 \quad \mbox{ for any } r>0\]
and \\
(ii) the boundaries $\partial E_k$ converge to $\partial E$ in $C^{1,\alpha}$ sense in a neighborhood of $q$.\\
We define in a similar way the $C^2$ convergence of sets.
\end{defn}

More precisely, 
we denote  
\[Q_{r,h}(x):=B'_r(x')\times(x_n-h,x_n+h),\]
for $x\in\R^n$, $r,h>0$. If $x=0$, we drop it in formulas and simply write $Q_{r,h}:=Q_{r,h}(0)$. Notice that up to a translation and a rotation, we can suppose that $q=0$ and
\eqlab{\label{CH:3:opossum1}E\cap Q_{2r,2h}=\{(x',x_n)\in\R^n\,|\,x'\in B'_{2r},\,u(x')<x_n<2h\},}
for some $r,h>0$ small enough and $u\in C^{1,\alpha}(\overline B'_{2r})$ such that $u(0)=0$. 
Then, point $(ii)$ means that we can write
\eqlab{\label{CH:3:opossum}E_k\cap Q_{2r,2h}=\{(x',x_n)\in\R^n\,|\,x'\in B'_{2r},\,u_k(x')<x_n<2h\},}
for some functions $u_k\in C^{1,\alpha}(\overline B_{2r}')$ such that
\eqlab{\label{CH:3:norm_graph_conv}\lim_{k\to\infty}\|u_k-u\|_{C^{1,\alpha}(\overline B'_{2r})}=0.}
\smallskip
We remark that, by the continuity of $u$,
up to considering a smaller $r$, we can suppose that
\eqlab{\label{CH:3:bded_graph_hp}|u(x')|<\frac{h}{2},\qquad\forall\,x'\in B'_{2r}.}

We have the following result.

\begin{theorem}\label{CH:3:everything_converges}
Let $E_k\xrightarrow{C^{1,\alpha}}E$ in a neighborhood of $q\in \partial E$. Let $q_k\in\partial E_k$ be such that $
q_k\longrightarrow q$ and let $s,s_k\in(0,\alpha)$ be such that $s_k\xrightarrow{k\to\infty} s$.
Then
\[\lim_{k\to\infty}\I_{s_k}[E_k](q_k)=\I_s[E](q).\]

Let $E_k\xrightarrow{C^2}E$ in a neighborhood of $q\in \partial E$. Let $q_k\in\partial E_k$ be such that $q_k\longrightarrow q$ and let $s_k\in(0,1)$ be such that $s_k\xrightarrow{k\to\infty} 1$. Then
\[\lim_{k\to\infty}(1-s_k)\I_{s_k}[E_k](q_k)=\varpi_{n-1}H[E](q).\]
%
\end{theorem}
\medskip 

A similar problem is studied also in \cite{mattheorem}, where the author estimates the
difference between the fractional mean curvature of a set $E$ with $C^{1,\alpha}$
boundary and that of the set
$\Phi(E)$, where $\Phi$ is a $C^{1,\alpha}$ diffeomorphism of $\R^n$,
%
in terms of the $C^{0,\alpha}$ norm of the Jacobian of the diffeomorphism $\Phi$.


\smallskip

When $s\to 0^+$ we do not need the $C^{1,\alpha}$ convergence of sets, but only the uniform boundedness  of the $C^{1,\alpha}$ norms of the functions defining the boundary of $E_k$ in a neighborhood of the boundary points. However, we have to require that the measure of the symmetric difference is uniformly bounded. More precisely:

\begin{prop}\label{CH:3:propsto0}
Let $ E\subseteq \Rn$
be such that $\alpha(E)$ exists.  Let 
$q \in \partial E$ be such that 
\bgs{ E\cap Q_{r,h}(q)=\{(x',x_n)\in\R^n\,|\, x'\in B'_{r}(q'),\,u(x')<x_n<h+q_n\},}
for some $r,h>0$ small enough and $u\in C^{1,\alpha}(\overline B'_{r}(q'))$ such that $u(q')=q_n$. 
Let $E_k\subseteq \Rn$ be such that
\[ |E_k\Delta E|<C_1  \] 
for some $C_1>0$. Let $q_k\in \partial E_k \cap B_d$, for some $d>0$, such that 
 \bgs{E_k\cap Q_{r,h}(q_k)=\{(x',x_n)\in\R^n\,|\,x'\in B'_{r}(q_k'),\,u_k(x')<x_n<h+q_{k,n}\} }
for some functions $u_k\in C^{1,\alpha}(\overline B_{r}'(q_k'))$ such that $u_k(q_k')=q_{k,n}$ and
\[ \|u_k\|_{C^{1,\alpha}(\overline B'_{r}(q'_k))} <C_2 \] 
for some $C_2>0$. Let $s_k\in(0,\alpha)$ be such that $s_k\xrightarrow{k\to\infty} 0$. Then
\bgs{\lim_{k\to\infty}s_k\I_{s_k}[E_k](q_k)=\varpi_n-2\alpha(E).}
\end{prop}

In particular, fixing $E_k=E$ in Theorem \ref{CH:3:everything_converges} and Proposition \ref{CH:3:propsto0} we obtain Proposition \ref{CH:3:rsdfyish} stated in the Introduction. 

To prove
Theorem \ref{CH:3:everything_converges}
we prove at first the following preliminary result. 

\begin{lemma}\label{CH:3:supergraph_hp_for_proof}
Let  $E_k\xrightarrow{C^{1,\alpha}}E$ in a neighborhood of $0\in \partial E$. Let $q_k\in\partial E_k$ be such that $
q_k\longrightarrow 0$.  Then \[ E_k-q_k\xrightarrow{C^{1,\beta}}E  \quad \mbox{ in a neighborhood of $0$},\] 
for every $\beta \in(0,\alpha)$.\\
Moreover, if $E_k\xrightarrow{C^{2}}E$ in a neighborhood of $0\in \partial E$, $q_k\in\partial E_k$ are such that $
q_k\longrightarrow 0$ and $\mathcal{R}_k\in  SO(n)$ are such that \[\lim_{k\to \infty} |\mathcal R_k -\mbox{Id}|=0,\] then
\[ \mathcal R_k (E_k-q_k) \xrightarrow{C^{2}}E \quad \mbox{ in a neighborhood of $0$ }.\]
%
%
%
\end{lemma}

\begin{proof}
First of all, notice that since $q_k\longrightarrow0$, for $k$ big enough we have
\[|q'_k|<\frac{1}{2}r\qquad\textrm{and}\qquad|q_{k,n}|=|u_k(q'_k)|<\frac{1}{8}h.\]
 By \eqref{CH:3:bded_graph_hp} and \eqref{CH:3:norm_graph_conv}, we
see that for $k$ big enough
\[|u_k(x')|\leq\frac{3}{4}h,\qquad\forall\,x'\in B_{2r}'.\]
Therefore
\[|u_k(x')-q_{k,n}|<\frac{7}{8}h<h,\qquad\forall\,x'\in B_{2r}'.\]
If we define
\[\tilde{u}_k(x'):=u_k(x'+q'_k),\qquad x'\in \overline{B}'_r,\]
for every $k$ big enough we have 
\eqlab{\label{CH:3:graphs_for_the_proof_eq}
(E_k-q_k) \cap Q_{r,h} =\{(x',x_n)\in\R^n\,|\,x'\in B'_r,\, \tilde u_k(x')<x_n<h\}.}
It is easy to check that the sequence $E_k-q_k$ locally converges in measure to $E$. We claim that
\eqlab{\label{CH:3:conv_transl_graph}
\lim_{k\to\infty}\|\tilde{u}_k-u\|_{C^{1,\beta}(\overline{B}'_r)}=0.}
Indeed, let 
\[ \tau_k u(x'):= u(x'+q_k').\] 
We have that
\[ \|\tilde u_k-\tau_k u\|_{C^{1}(\overline B_r')}  \leq \|u_k-u\|_{C^{1}\big(\overline B'_{\frac{3}2r}\big)}  \]
and that
\[ \|\tau_k u-u\|_{C^{1}(\overline B_r')} \leq  \|\nabla u\|_{C^0\big(\overline B'_{\frac{3}2r}\big)}  |q_k'| + \|u\|_{C^{1,\alpha}\big( \overline B'_{\frac{3r}2}\big)}|q'_k|^{\alpha} .\] Thus by the triangular inequality 
\[  \lim_{k \to \infty} \|\tilde u_k -u\|_{C^{1}(\overline B_r')} =0,\]
thanks to \eqref{CH:3:norm_graph_conv} and the fact that $q_k\to 0$.

Now, notice that $\nabla (\tilde u_k) =\tau_k (\nabla u_k)$, so 
\[ [\nabla \tilde u_k -\nabla u]_{C^{0,\beta}(\overline B'_r)} \leq  [\tau_k (\nabla  u_k -\nabla u)]_{C^{0,\beta}(\overline B'_r)}+ [\tau_k(\nabla  u) -\nabla u)]_{C^{0,\beta}(\overline B'_r)}.\] 
Therefore
\[[\tau_k (\nabla  u_k -\nabla u)]_{C^{0,\beta}(\overline B'_r)} \leq  [\nabla  u_k -\nabla u]_{C^{0,\beta}\big(\overline B'_{\frac{3r}2}\big)}\]
and for every $\delta >0$ we obtain
\[  [\tau_k( \nabla  u) -\nabla u]_{C^{0,\beta}(\overline B'_r)} \leq \frac{2}{\delta^\beta} \|\tau_k (\nabla u) -\nabla u\|_{C^{0}\big(\overline B'_{\frac{3r}2}\big)} + 2[\nabla u]_{C^{0,\alpha}(\overline B'_r)} \delta^{\alpha-\beta}.\]
Sending $k\to \infty$ we find that
\[ \limsup_{k\to \infty}  [\tau_k (\nabla  u) -\nabla u)]_{C^{0,\beta}(\overline B'_r)}  \leq 2[\nabla u]_{C^{0,\alpha}(\overline B'_r)} \delta^{\alpha-\beta}\]
for every $\delta>0$, hence
\[ \lim_{k \to \infty}  [\nabla \tilde u_k -\nabla u]_{C^{0,\beta}(\overline B'_r)} =0.\]
This concludes the proof of the first part of the Lemma.\\
As for the second part, the $C^2$ convergence of sets in a neighborhood of $0$ can be proved similarly. Some care must be taken when considering rotations, since one needs to use the implicit function theorem.
\end{proof}

\begin{proof}[Proof of Theorem \ref{CH:3:everything_converges}]

Up to a translation and a rotation, we can suppose that $q=0$ and $\nu_E(0)=0$. Then we can find $r,h>0$ small enough
and $u\in C^{1,\alpha}(\overline B'_r)$ such that we can write $E\cap Q_{2r,2h}$ as in \eqref{CH:3:opossum1}.

Since $s_k\to s\in(0,\alpha)$ for $k$ large enough we can suppose that $s_k,s \in[\sigma_0,\sigma_1]$ for $0<\sigma_0<\sigma_1<\beta<\alpha$.
Notice that there exists $\delta>0$ such that
\eqlab{\label{CH:3:continuity_eq2}
B_\delta\Subset Q_{r,h}.}
We take an arbitrary $R>1$ as large as we want and define the sets
\[F_k:= (E_k\cap B_R) -q_k.\]
From Lemma \ref{CH:3:supergraph_hp_for_proof} we have that in a neighborhood of $0$
\[ F_k\xrightarrow{C^{1,\beta}} E\cap B_R.\] 
In other words,
\eqlab{\label{CH:3:convvvvv1}\lim_{k\to \infty} |F_k \Delta (E\cap B_R)| =0.}
Moreover, if $u_k$ is a function defining $E_k$ as a supergraph in a neighborhood of $0$ as in \eqref{CH:3:opossum},
denoting $\tilde u_k(x')=u_k(x'+q_k')$ we have that
\[F_k\cap Q_{r,h}=\{(x',x_n)\in\Rn\,|\,x'\in B'_r,\,\tilde{u}_k(x')<x_n<h\}\]
and that
\eqlab{\label{CH:3:convvvvv} \lim_{k\to \infty} \|\tilde u_k -u\|_{C^{1,\beta}(\overline B'_r)} =0, \quad  \quad  \|\tilde u_k\|_{C^{1,\beta}(\overline B'_r)}\leq M \; \mbox{ for some } \; M>0.}
We also remark that, by \eqref{CH:3:bded_graph_hp} we can write
\[E\cap Q_{r,h}=\{(x',x_n)\in\R^n\,|\,x'\in B'_r,\,u(x')<x_n<h\}.\]

Exploiting \eqref{CH:3:graphs_for_the_proof_eq} we can write the fractional mean curvature of $F_k$ in $0$
by using formula \eqref{CH:3:complete_curv_formula}, that is
\begin{equation}\label{CH:3:5rs6ydbfd}
\begin{split}
\I_{s_k}[F_k](0)&=2\int_{B'_r}\Big\{G_{s_k}\Big(\frac{\tilde{u}_k(y')-\tilde{u}_k(0)}{|y'|}\Big)
-G_{s_k}\Big(\nabla \tilde{u}_k(0)\cdot\frac{y'}{|y'|}\Big)\Big\}\frac{dy'}{|y'|^{n-1+s_k}}\\
&
\qquad\qquad+\int_{\Rn}\frac{\chi_{\Co F_k}(y)-\chi_{F_k}(y)}{|y|^{n+s_k}}\chi_{\Co Q_{r,h}}(y)\,dy.\end{split}\end{equation}
Now, we denote as in \eqref{CH:3:mathcalg}
\[\mathcal G(s_k,\tilde{u}_k,y'):=\mathcal G(s_k,\tilde{u}_k,0,y')= G_{s_k}\Big(\frac{\tilde{u}_k(y')-\tilde{u}_k(0)}{|y'|}\Big)
-G_{s_k}\Big(\nabla \tilde{u}_k(0)\cdot\frac{y'}{|y'|}\Big)\] and we rewrite the
identity in \eqref{CH:3:5rs6ydbfd} as
\bgs{\I_{s_k}[F_k](0)&=2\int_{B'_r}\mathcal G(s_k,\tilde{u}_k,y') \frac{dy'}{|y'|^{n-1+s_k}}+\int_{\R^n}\frac{\chi_{\Co F_k}(y)-\chi_{F_k}(y)}{|y|^{n+s_k}}\chi_{\Co Q_{r,h}}(y)\,dy.} 
Also, with this notation and by formula \eqref{CH:3:complete_curv_formula} we have for $E$
\[\I_s[E\cap B_R](0)=2\int_{B'_r}\mathcal G(s,u,y')\frac{dy'}{|y'|^{n-1+s}}
+\int_{\R^n}\frac{\chi_{\Co (E\cap B_R)}(y)-\chi_{E\cap B_R}(y)}{|y|^{n+s}}\chi_{\Co Q_{r,h}}(y)\,dy.\]
%
We can  suppose that $r<1$. We begin by showing that for every $y'\in B_r'\setminus\{0\}$ we have
\eqlab{\label{CH:3:pwise_conv}\lim_{k\to\infty}\mathcal G(s_k,\tilde{u}_k,y')=\mathcal G(s,u,y').}
First of all, we observe that
\[|\mathcal G(s_k,\tilde{u}_k,y')-\mathcal G(s,u,y')|
\leq|\mathcal G(s_k,\tilde{u}_k,y')-\mathcal G(s,\tilde{u}_k,y')|+|\mathcal G(s,\tilde{u}_k,y')-\mathcal G(s,u,y')|.\]
Then
\bgs{|\mathcal G(s_k,\tilde{u}_k,y')-\mathcal G(s,\tilde{u}_k,y')|&
=\Big|\int_{\nabla \tilde{u}_k(0)\cdot\frac{y'}{|y'|}}^{\frac{\tilde{u}_k(y')-\tilde{u}_k(0)}{|y'|}}(g_{s_k}(t)-g_s(t))\,dt\Big|\\
&
\leq2\int_0^{+\infty}|g_{s_k}(t)-g_s(t)|\,dt.}
Notice that for every $t\in\R$
\[\lim_{k\to\infty}|g_{s_k}(t)-g_s(t)|=0,\qquad\textrm{and}\qquad |g_{s_k}(t)-g_s(t)|\leq2 g_{\sigma_0}(t),\quad\forall\,k\in\mathbb N.\]
Since $g_{\sigma_0}\in L^1(\R)$, by the Dominated Convergence Theorem we obtain that
\[\lim_{k\to\infty}|\mathcal G(s_k,\tilde{u}_k,y')-\mathcal G(s,\tilde{u}_k,y')|=0.\]
We estimate
\bgs{|\mathcal G&(s,\tilde{u}_k,y')-\mathcal G(s,u,y')|\leq\Big|G_s\Big(\frac{\tilde{u}_k(y')-\tilde{u}_k(0)}{|y'|}\Big)
-G_s\Big(\frac{u(y')-u(0)}{|y'|}\Big)\Big|\\
&
\qquad\qquad\qquad\qquad
+\Big|G_s\Big(\nabla\tilde{u}_k(0)\cdot\frac{y'}{|y'|}\Big)-G_s\Big(\nabla u(0)\cdot\frac{y'}{|y'|}\Big)\Big|\\
&
\leq\Big|\frac{\tilde{u}_k(y')-\tilde{u}_k(0)}{|y'|}-\frac{u(y')-u(0)}{|y'|}\Big|
+|\nabla\tilde{u}_k(0)-\nabla u(0)|\\
&
=\Big|\nabla(\tilde{u}_k-u)(\xi)\cdot\frac{y'}{|y'|}\Big|+|\nabla\tilde{u}_k(0)-\nabla u(0)|\\
&
\leq2\|\nabla\tilde{u}_k-\nabla u\|_{C^0(\overline B'_r)},
}
which, by \eqref{CH:3:conv_transl_graph}, tends to 0 as $k\to\infty$. This proves the pointwise convergence claimed in \eqref{CH:3:pwise_conv}.\\
Therefore, for every $y'\in B'_r\setminus\{0\}$,
\[\lim_{k\to\infty}\frac{\mathcal G(s_k,\tilde{u}_k,y')}{|y'|^{n-1+s_k}}=
\frac{\mathcal G(s,u,y')}{|y'|^{n-1+s}}.\]
Thus, by \eqref{CH:3:Holder_useful} we obtain that
\[\Big|\frac{\mathcal G(s_k,\tilde{u}_k,y')}{|y'|^{n-1+s_k}}\Big|
\leq\|\tilde{u}_k\|_{C^{1,\beta}(\overline B'_r)}\frac{1}{|y'|^{n-1-(\beta-s_k)}}
\leq \frac{M}{|y'|^{n-1-(\beta-\sigma_1)}}\in L^1_{\loc}(\R^{n-1}),
\]
given \eqref{CH:3:convvvvv}. 
The Dominated Convergence Theorem then implies that
\eqlab{\label{CH:3:first_piece_conv}
\lim_{k\to\infty}\int_{B'_r}\mathcal G(s_k,\tilde{u}_k,y')\frac{dy'}{|y'|^{n-1+s_k}}=
\int_{B'_r}\mathcal G(s,u,y')\frac{dy'}{|y'|^{n-1+s}}.
}

Now, we show that 
\eqlab{\label{CH:3:second_piece} \lim_{k \to \infty}\int_{\Rn}\frac{\chi_{\Co F_k}(y)-\chi_{F_k}(y)}{|y|^{n+s_k}} \chi_{\Co Q_{r,h}}(y)\, dy = \int_{\Rn}\frac{\chi_{\Co (E\cap B_R)}(y)-\chi_{E\cap B_R}(y)}{|y|^{n+s}}\chi_{\Co Q_{r,h}}(y) \, dy.}  
For this, we observe that
\bgs{\Big|\int_{\Co Q_{r,h}}&(\chi_{\Co (E\cap B_R)}(y)-\chi_{E\cap B_R}(y))\Big(\frac{1}{|y|^{n+s_k}}-\frac{1}{|y|^{n+s}}\Big)
dy\Big|
\leq\int_{\Co B_\delta}\Big|\frac{1}{|y|^{n+s_k}}-\frac{1}{|y|^{n+s}}\Big|dy,
}where we have used
\eqref{CH:3:continuity_eq2}  in the last inequality. 
For $y\in \Co B_1$ 
\bgs{ \Big|\frac{1}{|y|^{n+s_k}}-\frac{1}{|y|^{n+s}}\Big| \leq \frac{2}{|y|^{n+\sigma_0} }\in L^1(\Co B_1) }
and  for  $y\in  B_1\setminus B_\delta$ 
\bgs{ \Big|\frac{1}{|y|^{n+s_k}}-\frac{1}{|y|^{n+s}}\Big| \leq \frac{2}{|y|^{n+\sigma_1}} \in L^1(B_1\setminus B_\delta).} 
We use then the Dominated Convergence Theorem and get that
\bgs{ \lim_{k \to \infty} \int_{\Co Q_{r,h}}(\chi_{\Co (E\cap B_R)}(y)-\chi_{E\cap B_R}(y))\Big(\frac{1}{|y|^{n+s_k}}-\frac{1}{|y|^{n+s}}\Big)
dy =0.} 
Now
\bgs{ \bigg|\int_{\Co Q_{r,h}} & \;\frac{ \chi_{\Co F_k}(y)-\chi_F{_k}(y) -\left(\chi_{\Co (E\cap B_R)(y) }- \chi_{E\cap B_R} (y)\right) }{|y|^{n+s_k}}\, dy\bigg| =2\int_{\Co Q_{r,h}} \frac{\chi_{F_k \Delta (E\cap B_R)} (y)}{|y|^{n+s_k}}\, dy \\ \leq& \;2 \frac{ |F_k \Delta (E\cap B_R)|}{\delta^{n+\sigma_1}} \xrightarrow{k\to \infty} 0,}
according to \eqref{CH:3:convvvvv1}. The last two limits prove \eqref{CH:3:second_piece}. Recalling \eqref{CH:3:first_piece_conv}, we obtain that
\[ \lim_{k\to \infty} \I_{s_k}[F_k](0) = \I_s[E\cap B_R](0).\] 
We have that $\I_{s_k} [F_k](0)= \I_{s_k} [E_k\cap B_R](q_k)$, so
\bgs{
|\I_{s_k}[E_k](q_k)&-\I_s[E](0)|\leq |\I_{s_k}[E_k](q_k) - \I_{s_k}[E_k\cap B_R](q_k)|\\
&+ |\I_{s_k}[F_k](0)- \I_s[E\cap B_R](0)| + |\I_s[E\cap B_R](0)- \I_s[E](0)|.
}
Since
\bgs{\label{CH:3:ekbounded} | \I_{s_k}[E_k](q_k) - \I_{s_k} [E_k\cap B_R] (q_k) | +| \I_{s}[E](0) - \I_{s} [E\cap B_R] (0) | \leq  \frac{4 \varpi_n }{\sigma_0} R^{-\sigma_0},    }
sending $R\to \infty$
\[ \lim_{k \to \infty} \I_{s_k}[E_k](q_k)=\I_s[E](0).\]
This concludes the proof of the first part of the Theorem.

\bigskip

In order to prove the second part of Theorem \ref{CH:3:everything_converges}, we fix $R>1$ and we denote
\[F_k:=\mathcal R_k\big((E_k\cap B_R)-q_k\big),\]
where $\mathcal R_k\in SO(n)$ is a rotation such that
\[\mathcal R_k:\nu_{E_k}(0)\longmapsto\nu_E(0)=-e_n\quad\mbox{ and }\quad
\lim_{k\to\infty}|\mathcal R_k-\mbox{Id}|=0.\]
Thus, by Lemma \ref{CH:3:supergraph_hp_for_proof} we know that $F_k\xrightarrow{C^2}E$ in a neighborhood of $0$.\\
To be more precise,
\eqlab{\label{CH:3:convvvvv2}\lim_{k\to \infty} |F_k \Delta (E\cap B_R)| =0.}
Moreover, there exist $r,h>0$ small enough and $v_k,u\in C^2(\overline B'_r)$ such that
\bgs{&F_k\cap Q_{r,h}=\{(x',x_n)\in\Rn\,|\,x'\in B'_r,\,v_k(x')<x_n<h\},\\
&
E\cap Q_{r,h}=\{(x',x_n)\in\Rn\,|\,x'\in B'_r,\,u(x')<x_n<h\}}
and that
\eqlab{\label{CH:3:convvvvv3} \lim_{k\to \infty} \|v_k -u\|_{C^2(\overline B'_r)} =0.}
Notice that $0\in\partial F_k$ and $\nu_{F_k}(0)=e_n$ for every $k$, that is,
\eqlab{\label{CH:3:opossum4}v_k(0)=u(0)=0,\quad\nabla v_k(0)=\nabla u(0)=0.}

We claim that
\eqlab{\label{CH:3:opossum3}\lim_{k\to\infty}(1-s_k)\big|\I_{s_k}[F_k](0)-\I_{s_k}[E\cap B_R](0)\big|=0.}
By \eqref{CH:3:opossum4} 
%
and formula \eqref{CH:3:complete_curv_formula} we have that
\bgs{\label{CH:3:opossum5}
\I_{s_k}[F_k](0)&=2\int_{B'_r}\frac{dy'}{|y'|^{n+s_k-1}} \int_0^{\frac{v_k(y')}{|y'|} }\frac{dt}{(1+t^2)^{\frac{n+s_k}2}}
+\int_{\Co Q_{r,h}}\frac{\chi_{\Co F_k}(y)-\chi_{F_k}(y)}{|y|^{n+s_k}}\,dy\\
&= \I^{loc}_{s_k}[F_k](0) +\int_{\Co Q_{r,h}}\frac{\chi_{\Co F_k}(y)-\chi_{F_k}(y)}{|y|^{n+s_k}}\,dy.}
We use the same formula for $E\cap B_R$ and prove at first that
\bgs{
\bigg|\int_{\Co Q_{r,h}}\frac{\chi_{\Co F_k}(y)-\chi_{F_k}(y)-\chi_{\Co(E\cap B_R)}(y)+\chi_{E\cap B_R}(y)}{|y|^{n+s_k}}\,dy\bigg|&\le\frac{|F_k\Delta(E\cap B_R)|}{\delta^{n+s_k}}\\
&\le\frac{|F_k\Delta(E\cap B_R)|}{\delta^{n+1}},
}
(where we have used \eqref{CH:3:continuity_eq2}), which tends to 0 as $k\to\infty$, by \eqref{CH:3:convvvvv2}.

Moreover, notice that by the Mean Value Theorem and \eqref{CH:3:opossum4} we have
\[|(v_k-u)(y')|\le\frac{1}{2}|D^2(v_k-u)(\xi')||y'|^2\le\frac{\|v_k-u\|_{C^2(\overline B'_r)}}{2}|y'|^2.\]
Thus
\bgs{ &\big| \I^{loc}_{s_k}[F_k](0)  - \I^{loc}_{s_k}[E\cap B_R](0) |
\le2\int_{B'_r} \frac{dy'}{|y'|^{n+s_k-1}} \bigg|\int_{\frac{u(y')}{|y'|}}^{\frac{v_k(y')}{|y'|}}\frac{dt}{(1+t^2)^\frac{n+s_k}{2}}\bigg|\\
&
\le2\int_{B'_r}|y'|^{-n-s_k}|(v_k-u)(y')|\,d y'
\le\frac{\varpi_{n-1}\,\|v_k-u\|_{C^2(\overline B'_r)}}{1-s_k}r^{1-s_k},
}
hence by \eqref{CH:3:convvvvv3} we obtain
\eqlab{\label{CH:3:opossum7}
\lim_{k\to\infty}(1-s_k)\big|\I^{loc}_{s_k}[F_k](0)  - \I^{loc}_{s_k}[E\cap B_R](0)|=0.
}
This concludes the proof of claim \eqref{CH:3:opossum3}.

Now we use the triangle inequality and have that
\bgs{
\big|(1&-s_k)\I_{s_k}[E_k](q_k)-H[E](0)\big|\le(1-s_k)\big|\I_{s_k}[E_k](q_k)-\I_{s_k}[F_k](0)\big|\\
&
+(1-s_k)\big|\I_{s_k}[F_k](0)-\I_{s_k}[E\cap B_R](0)\big|
+\big|(1-s_k)\I_{s_k}[E\cap B_R](0)-H[E](0)\big|.
}
The last term in the right hand side converges by \cite[Theorem 12]{Abaty}. As for the first term,
notice that
\[\I_{s_k}[F_k](0)=\I_{s_k}[E_k\cap B_R](q_k),\]
hence
\[\lim_{k\to\infty}(1-s_k)\big|\I_{s_k}[E_k\cap B_R](q_k)-\I_{s_k}[E_k](q_k)\big|\le\limsup_{k\to\infty}(1-s_k)\frac{2\varpi_n}{s_k}R^{-s_k}=0.\]
Sending $k\to\infty$ in the triangle inequality above, we conclude the proof of the second part of Theorem \ref{CH:3:everything_converges}.
\end{proof}

\begin{remark}
In relation to the second part of the proof, we point out that using the directional fractional mean curvature defined in \cite[ Definition 6, Theorem 8]{Abaty}, we can write
\bgs{ \I^{loc}_{s_k}[F_k](0)= &\;2\int_{\mathbb S^{n-2}}\bigg[\int_0^r\varrho^{n-2}\bigg(\int_0^{v_k(\varrho e)}\frac{dt}{(\varrho^2+t^2)^\frac{n+s_k}{2}}\bigg)d\varrho\bigg]d\mathcal H^{n-2}_e \\
=&\;2\int_{\mathbb S^{n-2}} \overline K_{s_k,e} d\mathcal H^{n-2}_e.}
One is then actually able to prove that
\bgs{\lim_{k\to\infty}(1-s_k)\overline K_{s_k,e}[E_k-q_k](0)=H_e[E](0),}
uniformly in $e\in\mathbb S^{n-2}$, by using formula \eqref{CH:3:opossum7} and the first claim of \cite[Theorem 12]{Abaty}.
\end{remark}

\begin{remark}\label{CH:3:nuno}
	The proof of Theorem \ref{CH:3:everything_converges}, as well as the proof of the next Proposition \ref{CH:3:propsto0}, settles the case in which $n\geq 2$. For $n=1$, the proof follows in the same way, after observing that the local contribution to the fractional mean curvature is equal to zero because of symmetry. As a matter of fact, the formula in~\eqref{CH:3:complete_curv_formula} for the fractional mean curvature (which has no meaning for $n=1$) is not required.\\
	We remark also that in our notation $\varpi_0=0$. This gives consistency to the second claim of Theorem \ref{CH:3:everything_converges} also for $n=1$.
\end{remark}

\bigskip
We prove now the continuity of the fractional mean curvature as $s\to 0$.

\begin{proof}[Proof of Proposition \ref{CH:3:propsto0}]
Up to a translation, we can take $q=0$ and $u(0)=0$. \\
For $R>2\max\{r,h\}$, we write
\bgs{  \I_{s_k}[E_k](q_k) =&\;
\PV \int_{ Q_{r,h}(q_k)} \frac{\chi_{\Co E_k}(y) -\chi_{E_k}(y)}{|y-q_k|^{n+s_k}}\, dy
+\int_{\Co Q_{r,h}(q_k)} \frac{\chi_{\Co E_k}(y)-\chi_{E_k}(y)}{|y-q_k|^{n+s_k}}\, dy 
\\
=&\; \PV \int_{ Q_{r,h}(q_k)} \frac{\chi_{\Co E_k}(y)-\chi_{E_k}(y)}{|y-q_k|^{n+s_k}}\, dy  + \int_{ B_R(q_k)\setminus Q_{r,h}(q_k)}\frac{\chi_{\Co E_k}(y)-\chi_{E_k}(y)}{|y-q_k|^{n+s_k}}\, dy
\\
&\; +\int_{\Co B_R(q_k)} \frac{\chi_{\Co E_k}(y)-\chi_{E_k}(y)}{|y-q_k|^{n+s_k}}\, dy  
 \\
=&\;I_1(k)+I_2(k)+I_3(k).}
Now  using  \eqref{CH:3:complete_curv_formula}, \eqref{CH:3:mathcalg} and \eqref{CH:3:Holder_useful}  we have that 
\bgs{|I_1(k)|\leq &\;  2 \int_{B'_r(q_k')} \frac{| \mathcal G(s_k,u_k,q_k', y')|}{|y'-q_k'|^{n+s_k-1}} \,dy'
 \leq 2  \| u_k \|_{C^{1,\alpha}(\overline{B}'_r(q_k'))}\int_{B'_r(q_k')} \frac{  |y'-q_k'|^{\alpha}}{|y'-q_k'|^{n+s_k-1}} \,dy'
 \\
 \leq  &\; 2 C_2 \varpi_{n-1} \frac{r^{\alpha-s_k}}{\alpha-s_k}.
}
Using \eqref{CH:3:continuity_eq2} we also have that
\[|I_2(k)| \leq  \int_{ B_R(q_k)\setminus B_\delta(q_k)}\frac{dy}{|y-q_k|^{n+s_k}} 
= \varpi_n\frac{\delta^{-s_k}-R^{-s_k}}{s_k}. \]
Thus
\eqlab{\label{CH:3:kangaroo}
\lim_{k\to\infty}s_k\big(|I_1(k)|+|I_2(k)|\big)=0.}
Furthermore
\bgs{ \big | s_k& I_3(k)- \big(\varpi_n-2 s_k \alpha_{s_k}(0,R,E) \big)\big|
\\
\leq &\; \bigg|s_k \int_{\Co B_R(q_k)} \frac{dy}{|y-q_k|^{n+s_k}}  -2 s_k \int_{\Co B_R(q_k)} \frac{\chi_{E_k}(y)}{|y-q_k|^{n+s_k}}\, dy   - \varpi_n+2s_k \alpha_{s_k}(q_k,R,E))\bigg|
\\
&\;+2s_k |\alpha_{s_k}(q_k,R,E)- \alpha_{s_k}(0,R,E)|
\\
\leq& \; |\varpi_n R^{-s_k}-\varpi_n| + 2 s_k\bigg| \int_{\Co B_R(q_k)}\frac{\chi_{E_k}(y)}{|y-q_k|^{n+s_k}}\, dy  - \int_{\Co B_R(q_k)} \frac{\chi_{E}(y)}{|y-q_k|^{n+s_k}}\, dy \bigg| 
\\
&\;+ 2s_k |\alpha_{s_k}(q_k,R,E)- \alpha_{s_k}(0,R,E)|
\\
\leq &\; |\varpi_n R^{-s_k}-\varpi_n| + 2 s_k \int_{\Co B_R(q_k)}\frac{\chi_{E_k\Delta E}(y)}{|y-q_k|^{n+s_k}}\, dy  
+ 2s_k| \alpha_{s_k}(q_k,R,E)- \alpha_{s_k}(0,R,E)|
\\
\leq &\; |\varpi_n R^{-s_k}-\varpi_n| + 2 C_1 s_k R^{-n-s_k} +2s_k | \alpha_{s_k}(q_k,R,E)- \alpha_{s_k}(0,R,E)|,
 }
 where we have used that $|E_k\Delta E|<C_1$.
 
 Therefore, since $q_k\in B_d$ for every $k$, as a consequence of Proposition \ref{CH:3:unifrq} it follows that
\eqlab{\label{CH:3:kangaroo1}\lim_{k\to \infty}\big | s_k I_3(k)- &\big(\varpi_n-2 s_k \alpha_{s_k}(0,R,E) \big)\big|=0.}

Hence, by \eqref{CH:3:kangaroo} and \eqref{CH:3:kangaroo1}, we get that
\bgs{ \lim_{k \to \infty} s_k \I_{s_k}[E_k](q_k)=\varpi_n-2\lim_{k\to\infty}s_k\alpha_{s_k}(0,R,E)=\varpi_n-2\alpha(E),}
concluding the proof.
\end{proof} 

\begin{proof}[Proof of Theorem \ref{CH:3:asympts}]
Arguing as in the proof of Proposition \ref{CH:3:propsto0}, by keeping fixed $E_k=E$ and $q_k=p$, we obtain
\bgs{ \liminf_{s\to0} s\, \I_s[E](p)=\varpi_n-2\limsup_{s\to0}s\,\alpha_s(0,R,E)=\varpi_n-2\overline{\alpha}(E),}
and similarly for the limsup.
\end{proof}

As a corollary of Theorem \ref{CH:3:everything_converges} and Theorem \ref{CH:3:asympts}, we have the following result.
\begin{theorem}\label{CH:3:changeyoursign}
Let $E\subseteq\R^n$ and let $p\in\partial E$ be such that $\partial E\cap B_r(p)$ is $C^2$ for some $r>0$.
Suppose that the classical mean curvature of $E$ in $p$ is $H(p)<0$. Also assume that 
\[\overline \alpha(E) < \frac{\varpi_n}2.\] Then there exist $\sigma_0<\tilde{s}<\sigma_1$ in $(0,1)$ such that

$(i)\quad\I_s[E](p)>0$ for every $s\in(0,\sigma_0]$, and actually
\[\liminf_{s\to0^+}s \;\I_s[E](p)=\varpi_n- 2\overline \alpha(E),\]

$(ii)\quad\I_{\tilde{s}}[E](p)=0,$

$(iii)\quad\I_s[E](p)<0$ for every $s\in[\sigma_1,1)$,
and actually
\[ \lim_{s\to 1} (1-s)\;\I_s[E](p)= \varpi_{n-1}H[E](p).\]

\end{theorem}

\end{chapter}

\begin{chapter}{On nonlocal minimal graphs}\label{CH_Nonparametric}

\minitoc

\section{Introduction}


The aim of this chapter consists in introducing a functional framework for studying minimizers of the fractional perimeter that can be globally written as subgraphs.

More precisely, we define a functional~$\F_s$, which can be considered as a fractional and nonlocal version of the area functional,
and we exploit it to study nonlocal minimal graphs.

One of the main difficulties in defining a fractional and nonlocal version of the classical area functional is that,
as observed in Chapter \ref{CH_Appro_Min},
\bgs{
\Per_s(\{x_{n+1}<u(x)\},\Omega\times\R)=\infty,
}
independently of the regularity of~$u$---see Theorem \ref{CH:2:bound_unbound_per_cyl_prop} and Corollary \ref{CH:2:non_well_def_frac_area}.
Nevertheless, this problem can be avoided by working in the ``truncated cylinders''~$\Omega\times(-M,M)$.
In the functional setting that we introduce, this leads us to consider a family of functionals~$\F^M_s$, instead of
only the global functional~$\F_s$.

Exploiting these approximating functionals, we prove existence and uniqueness results for minimizers of the functional~$\F_s$---and actually of more general functionals---for a large class of exterior data which includes locally bounded functions.

Moreover, one of the main contributions of this chapter consists in proving the equivalence of:
\begin{itemize}
\item minimizers of the functional~$\F_s$,
\item minimizers ot the fractional perimeter,
\item weak solutions of the fractional mean curvature equation,
\item viscosity solutions of the fractional mean curvature equation,
\item smooth pointwise solutions of the fractional mean curvature equation,
\end{itemize}
(see Theorem \ref{CH:4:equiv_intro}).

We observe that the functional framework introduced in this chapter easily adapts to the obstacle problem.
Hence we prove also existence and uniqueness results for the nonlocal Plateau problem with (eventually discontinuous)
obstacles.

\smallskip

Now we proceed to give the definitions and the precise statements of the main results of the chapter.

\subsection{Definitions and main results}

Let~$g: \R \to \R$ be a continuous function satisfying
\bgs{
g(t) = g(-t)  \quad \mbox{for every } t \in \R, \qquad 0 < g \le 1 \quad \mbox{in } \R,
}
and
\bgs{
\lambda := \int_{0}^{+\infty} g(t) t \, dt < \infty.
}
Then, we define
\begin{equation*}
G(t) := \int_0^t g(\tau) \, d\tau \quad \mbox{and} \quad \G(t) := \int_0^t G(\tau) \, d\tau = \int_0^t \left( \int_0^\tau g(\sigma) \, d\sigma \right) d\tau.
\end{equation*}

Given any function~$u: \R^n \to \R$, we also formally set
\begin{equation} \label{CH:4:Fcdef}
\F(u, \Omega) := \iint_{Q(\Omega)} \G \left( \frac{u(x) - u(y)}{|x - y|} \right) \frac{dx\,dy}{|x - y|^{n - 1 + s}},
\end{equation}
where
$$
Q(\Omega):=\R^{2n}\setminus(\Co\Omega)^2.
$$

A particularly important example of function~$g$ is given by
\eqlab{\label{CH:4:gdef}
g_s(t) := \frac{1}{(1 + t^2)^{\frac{n + 1 + s}{2}}}.
}
We indicate with~$G_s$ and~$\G_s$ respectively the first and second integrals of~$g_s$ as in~\eqref{CH:4:Gdefs}. Furthermore,~$\F_s$ denotes the functional corresponding to~$\G_s$ in light of definition~\ref{CH:4:Fcdef}.

We will consider the following space
\bgs{
\W^s(\Omega):=\{u:\R^n\to\R\,|\,u|_\Omega\in W^{s,1}(\Omega)\}.
}
Given a function~$\varphi:\Co\Omega\to\R$ we also define the space
\eqlab{\label{CH:4:domain_w_data}
\W^s_\varphi(\Omega):=\{v\in\W^s(\Omega)\,|\,v=\varphi\mbox{ a.e. in }\Co\Omega\}.
}

Our aim will be that of minimizing the functional~$\F$ in~$\W^s_\varphi(\Omega)$, given
a fixed function~$\varphi:\Co\Omega\to\R$ as exterior data.

However, we remark that the functional~$\F$ is not well defined on functions~$u\in\W^s_\varphi(\Omega)$, unless the function~$\varphi$ has a suitable growth at infinity, namely
\eqlab{\label{CH:4:tartainfty}
\int_\Omega\Big(\int_{\Co\Omega}\frac{|\varphi(y)|}{|x-y|^{n+s}}dy\Big)dx<\infty,
}
which is a quite restrictive condition.

Nevertheless, as ensured by Lemma~\ref{CH:4:tartariccio}---exploiting the fractional Hardy-type inequality of Theorem \ref{CH:4:FHI}---the following definition of minimizer is well posed:

\begin{defn}\label{CH:4:minim_tarta_def}
Let~$\Omega \subseteq \R^n$ be a bounded open set
with Lipschitz boundary. A function~$u\in\W^s(\Omega)$ is a \emph{minimizer} of~$\F$ in~$\Omega$ if
$$
\iint_{Q(\Omega)} \left\{ \G \left( \frac{u(x) - u(y)}{|x - y|} \right) - \G \left( \frac{v(x) - v(y)}{|x - y|} \right) \right\} \frac{dx\,dy}{|x - y|^{n - 1 + s}} \le 0
$$
for every~$v\in\W^s(\Omega)$ such that~$v=u$ almost everywhere in~$\Co\Omega$.
\end{defn}

Fixed~$\varphi:\Co\Omega\to\R$, we consider the problem of finding a function~$u\in\W^s_\varphi(\Omega)$
which is a minimizer for the functional~$\F$ in the sense of Definition~\ref{CH:4:minim_tarta_def}.

One of the main difficulties of this chapter will be that of finding such a minimizer without imposing the global condition~\eqref{CH:4:minim_tarta_def} on the exterior data. This will be done by asking a suitable weaker condition
on the exterior data~$\varphi$ and by exploiting a ``truncation procedure'' for the functional~$\F$.

\begin{defn}
Let~$\Omega$ be a bounded open set and let~$u:\Co\Omega\to\R$. Given an open set~$\Op\subseteq\R^n$
such that~$\Omega\subseteq\Op$, we define
the ``truncated tail'' of~$u$ at a point~$x\in\Omega$ as
\bgs{
\Tail_s(u, \Op\setminus\Omega;x):=\int_{\Op\setminus\Omega}\frac{|u(y)|}{|x-y|^{n+s}}\,dy.
}
\end{defn}


It is convenient to recall that, given a set $F\subseteq\R^n$, we denote
\[
F_r:=\left\{x\in\R^n\,|\,\bar{d}_F(x)<r\right\},
\]
for any $r\in\R$, with $\bar{d}_F$ denoting the signed distance function from $\partial F$, negative inside $F$. In particular, if $\Omega\subseteq\R^n$ is a bounded open set and $\varrho>0$, then
\[
\Omega_{-\varrho}\Subset\Omega\Subset\Omega_\varrho.
\]
We will make extensive use of this notation in the present chapter.

One of the main results of this chapter consists in proving the existence and uniqueness of a minimizer~$u\in\W^s_\varphi(\Omega)$
for exterior data~$\varphi$ whose tail is integrable in a large enough neighborhood of~$\Omega$.

\begin{theorem} \label{CH:4:Dirichlet}
Let~$n \ge 1$,~$s \in (0, 1)$, and~$\Omega \subseteq \R^n$ be a bounded open set with Lipschitz boundary.
Then, there is a constant~$\Theta > 1$, depending only on~$n$ and~$s$ and~$g$, such that,
given any function~$\varphi: \Co\Omega\to \R$
with~$\Tail_s(\varphi, \Omega_{\Theta \diam(\Omega)} \setminus \Omega;\,\cdot\,) \in L^1(\Omega)$, there exists a unique minimizer~$u$ of~$\F$ within~$\W^s_\varphi(\Omega)$. Moreover,~$u$ satisfies
\begin{equation} \label{CH:4:Ws1estformin}
\| u \|_{W^{s, 1}(\Omega)} \le C \left( \left\| \Tail_s(\varphi,\Omega_{\Theta \diam(\Omega)}\setminus \Omega;\,\cdot\,) \right\|_{L^1(\Omega)} + 1 \right),
\end{equation}
for some constant~$C > 0$ depending only on~$n$,~$s$,~$g$ and~$\Omega$.
\end{theorem}

We remark that asking
\bgs{
\left\| \Tail_s(\varphi,\Op\setminus \Omega;\,\cdot\,) \right\|_{L^1(\Omega)}
=\int_\Omega\Big(\int_{\Op\setminus\Omega}\frac{|\varphi(y)|}{|x-y|^{n+s}}dy\Big)dx<\infty,
}
is a much weaker requirement than asking~\eqref{CH:4:minim_tarta_def}, since
we impose no conditions on~$\varphi$ in~$\Co\Op$.

The proof of Theorem \ref{CH:4:Dirichlet} is the content of Section \ref{CH:4:Proof_of_Dir_Sec}. The argument exploits the minimizers of appropriate truncated functionals $\F^M(\,\cdot\,,\Omega)$, considered within their natural domain, and an apriori bound on the $W^{s,1}(\Omega)$ norm, which gives \eqref{CH:4:Ws1estformin}. These topics are studied in Section \ref{CH:4:TRUNC_FUN_SEC}.

See also Section \ref{CH:4:FTPOTFAF_SEC} for the definition of the functionals $\F^M(\,\cdot\,,\Omega)$ and for their main functional properties, and Section \ref{CH:4:areageom} for the relationship existing between $\F^M(\,\cdot\,,\Omega)$ and the $s$-perimeter---in the geometric case $g=g_s$.

We also observe that if~$\varphi$ is bounded in~$\Op\setminus\Omega$, then, since~$\Omega$ is bounded and has Lipschitz boundary, we have $\Tail_s(\varphi,\Op\setminus \Omega;\,\cdot\,)\in L^1(\Omega)$---for more details about the integrability of the truncated tail, we refer to Lemma~\ref{CH:4:tail_equiv_cond_Lemma}.

Hence, the boundedness of~$\varphi$ in a large enough neighborhood of~$\Omega$ is enough to guarantee the existence of a unique minimizer of~$\F$.
Furthermore, in this case
we prove that the minimizer is bounded also in~$\Omega$. More precisely:

\begin{theorem} \label{CH:4:minareboundedthm}
Let~$n \ge 1$,~$s \in (0, 1)$,~$\Omega \subseteq \R^n$ be a bounded open set with Lipschitz boundary,
and~$R_0 > 0$ be such that~$\Omega \subseteq B_{R_0}$. There exists a large constant~$\Theta > 1$,
depending only on~$n$,~$s$ and~$g$,
such that if~$u\in\W^s(\Omega)$ is a minimizer of~$\F$ in~$\Omega$,
bounded in~$B_{\Theta R_0} \setminus \Omega$, then~$u$ is also bounded in~$\Omega$ and
\bgs{
\| u \|_{L^\infty(\Omega)} \le R_0 + \| u \|_{L^\infty(B_{\Theta R_0} \setminus \Omega)}.
}
\end{theorem}

We observe that, even when the exterior data~$\varphi$ is not bounded in a neighborhood of~$\Omega$,
we are nevertheless able to prove that the minimizer of~$\F$ in~$\W^s_\varphi(\Omega)$,
if it exists, is locally bounded inside~$\Omega$ (see Proposition~\ref{CH:4:Linftylocprop}).

Moreover, we point out that in order to obtain the global boundedness of the minimizer~$u\in\W^s_\varphi(\Omega)$
inside~$\Omega$, it is actually enough to require the function~$\varphi$ to be bounded only in a neighborhood~$\Omega_r\setminus\Omega$, with~$r>0$ as small as we want. However, we remark that in this case the apriori~$L^\infty$ bound is not as clean as the one of Theorem~\ref{CH:4:minareboundedthm} (see Theorem~\ref{CH:4:Bdary_Bdedness_Thm} for the precise statement).

\smallskip

Let us also mention that in Section \ref{CH:4:NPPWO_SEC} we will partially extend the above results to the obstacle problem. More precisely, we will prove the existence and uniqueness of a minimizer, in the case of locally bounded exterior data only, and we will establish an apriori bound on the $L^\infty(\Omega)$ of the minimizer.

\smallskip
 
The Euler-Lagrange operator associated to the minimization of~$\F$ is
\bgs{
\h u(x):=2\,\PV\int_{\R^n}G\Big(\frac{u(x)-u(y)}{|x-y|}\Big)\frac{dy}{|x-y|^{n+s}}.
}
We remark that in order for~$\h u(x)$ to be well defined, the function $u$ must be regular enough (e.g.~$C^{1,\alpha}$
for some~$\alpha>s$) in a neighborhood of the point~$x$.

On the other hand, we can always define~$\h u$ in the distributional sense, for any measurable~$u:\R^n\to\R$, as
the linear functional
\bgs{
\langle\h u,v\rangle:=\int_{\R^n}\int_{\R^n}G\Big(\frac{u(x)-u(y)}{|x-y|}\Big)\big(v(x)-v(y)\big)\frac{dx\,dy}{|x-y|^{n+s}},
}
for every~$v\in W^{s,1}(\R^n)$.

This observation prompts us to give the following definition of weak solution:

\begin{defn}
Let~$\Omega\subseteq\R^n$ be a bounded open set and let~$f\in C(\overline{\Omega})$. We say that a function~$u:\R^n\to\R$
is a weak solution of~$\h u=f$ in~$\Omega$ if
\bgs{
\langle \h u,v\rangle=\int_\Omega fv\,dx,\qquad\forall\,v\in C_c^\infty(\Omega).
}
\end{defn}
Some elementary properties of the operator $\h$ are studied in Section \ref{CH:4:SFATELO_SEC}.

Exploiting the convexity of the functional~$\F$, it is easy to verify---see Lemma \ref{CH:4:weak_implies_min_lemma}---that if we add the requirement that~$u\in\W^s(\Omega)$, then
\bgs{
u \mbox{ is a minimizer of }\F\mbox{ in }\Omega\quad
\Longleftrightarrow\quad
u\mbox{ is a weak solution of }\h u=0\mbox{ in }\Omega.
}

Besides distributional solutions, another natural notion of solutions to consider for
the problem
\[
\left\{\begin{array}{cc}
\h u=f & \mbox{in }\Omega,\\
u=\varphi & \mbox{in }\Co\Omega.
\end{array}\right.
\]
is that of viscosity solutions. We will use~$C^{1,1}$ functions as test functions.

\begin{defn}
Let~$\Omega\subseteq\R^n$ be a bounded open set and let~$f\in C(\overline{\Omega})$.
We say that a function~$u:\R^n\to\R$ is a (viscosity) subsolution of~$\h u=f$ in~$\Omega$, and we write
\[
\h u\leq f\qquad\textrm{in }\Omega,
\]
if~$u$ is upper semicontinuous in~$\Omega$ and whenever the following happens:
\begin{itemize}
\item[(i)] $x_0\in\Omega$,
\item[(ii)] $v\in C^{1,1}(B_r(x_0)),$ for some~$r<d(x_0,\partial\Omega)$,
\item[(iii)] $v(x_0)=u(x_0)$ and~$v(y)\geq u(y)$ for every~$y\in B_r(x_0)$,
\end{itemize}
then if we define
\[
\tilde v(x):=\left\{\begin{array}{cc} v(x) & \textrm{if }x\in B_r(x_0),\\
u(x) & \textrm{if }x\in\R^n\setminus B_r(x_0),\end{array}\right.
\]
we have
\[
\h\tilde v(x_0)\leq f(x_0).
\]
A supersolution is defined similarly. A (viscosity) solution is a function~$u:\R^n\to\R$ which is continuous in~$\Omega$
and which is both a subsolution and a supersolution.
\end{defn}

We remark that in the definition of a viscosity subsolution we do not ask~$u$ to be upper semicontinuous in~$\overline\Omega$
but only in~$\Omega$. Furthermore, we do not ask~$u$ to belong to the functional space~$\W^s(\Omega)$.

Another important result of this chapter consists in proving that viscosity (sub)solutions are also weak (sub)solutions.

\begin{theorem}\label{CH:4:Gen_viscweak}
Let~$\Omega\subseteq\R^n$ be a bounded open set and let~$f\in C(\overline\Omega)$. Let~$u:\R^n\to\R$
be such that~$u$ is locally integrable in~$\R^n$ and~$u$ is locally bounded in~$\Omega$. If~$u$
is a viscosity subsolution,
\bgs{
\h u\le f\quad\mbox{in }\Omega,
}
then~$u$ is a weak subsolution,
\bgs{
\langle\h u,v\rangle\le\int_\Omega fv\,dx,\qquad\forall\,v\in C^\infty_c(\Omega)\mbox{ s.t. }v\ge0.
}
\end{theorem}

It is worth to mention also a global version, for viscosity solutions, of this Theorem. Given a continuous function~$f\in C(\R^n)$,
we say that a function~$u:\R^n\to\R$ is a viscosity solution of~$\h u=f$ in~$\R^n$ if~$u\in C(\R^n)$ and~$u$ is
a viscosity solution in every bounded open set~$\Omega\subseteq\R^n$.

\begin{corollary}\label{CH:4:Visco_GLOB_CoRoLl}
Let~$f\in C(\R^n)$ and let~$u:\R^n\to\R$. If~$u$ is a viscosity solution of~$\h u=f$ in~$\R^n$, then~$u$
is a weak solution,
\bgs{
\langle\h u,v\rangle=\int_{\R^n} fv\,dx,\qquad\forall\,v\in C^\infty_c(\R^n).
}
\end{corollary}

The study of viscosity (sub)solutions and the proof of Theorem \ref{CH:4:Gen_viscweak} are carried out in Section \ref{CH:4:ViscWeak_Sec}.

\subsubsection{Geometric case}
The case in which~$g=g_s$ is particularly important, because it is connected with the nonlocal minimal surfaces.
In particular,
\bgs{
\h_su(x)=\I_s[\Sg(u)](x,u(x)),
}
is the~$s$-fractional mean curvature of the subgraph of~$u$,
\bgs{
\Sg(u):=\{X=(x,x_{n+1})\in\R^{n+1}\,|\,x_{n+1}<u(x)\},
}
at the point~$(x,u(x))\in\partial\Sg(u)$ (provided~$u$ is regular enough near~$x$).

Therefore, the equation
\bgs{
\h_su=0
}
is, at least formally, the Euler-Lagrange equation of an~$s$-minimal set which can be globally written as a subgraph.

Before going on, we recall that the~$s$-fractional perimeter of a set~$E\subseteq\R^{n+1}$
in an open set~$\Op\subseteq\R^{n+1}$
is defined as
\bgs{
\Per_s(E,\Op)=\Ll_s(E\cap\Op,\Co E\cap\Op)+\Ll_s(E\cap\Op,\Co E\setminus\Op)+\Ll_s(E\setminus\Op,\Co E\cap\Op),
}
where
\bgs{
\Ll_s(A,B):=\int_A\int_B\dKers,
}
for every couple of disjoint sets~$A,B\subseteq\R^{n+1}$. We also observe that we can rewrite the $s$-perimeter as
\bgs{
\Per_s(E,\Op)=\frac{1}{2}\iint_{\R^{2(n+1)}\setminus(\Co\Op)^2}\frac{|\chi_E(X)-\chi_E(Y)|}{|X-Y|^{n+1+s}}dX\,dY.
}
The $s$-fractional mean curvature of~$E$
at~$X\in\partial E$ is the principal value integral
\bgs{
\I_s[E](X):=\PV\int_{\R^{n+1}}\frac{\chi_{\Co E}(Y)-\chi_E(Y)}{|X-Y|^{n+1+s}}dy.
}

\begin{defn}
Let~$\Op\subseteq\R^{n+1}$ be an open set and let~$E\subseteq\R^{n+1}$. We say that~$E$ is $s$-minimal in~$\Op$
if~$\Per_s(E,\Op)<\infty$ and
\bgs{
F\setminus\Op=E\setminus\Op\qquad\implies\qquad \Per_s(E,\Op)\le \Per_s(F,\Op).
}
We say that~$E$ is locally~$s$-minimal in~$\Op$ if it is $s$-minimal in every~$\Op'\Subset\Op$.
\end{defn}

In this chapter we are interested in the case where the domain is a cylinder,~$\Op=\Omega\times\R$.
For simplicity, we introduce the following notation:
\bgs{
\Omega^M:=\Omega\times(-M,M),\quad\forall\,M\ge0\qquad\textrm{and}\qquad\Omega^\infty:=\Omega\times\R.
}

We remark that when~$\Omega$ is bounded and has Lipschitz boundary, then a set~$E$ is locally $s$-minimal in~$\Omega^\infty$ if and only if it is $s$-minimal in~$\Omega^M$, for every~$M>0$---see Remark~\ref{CH:2:rmk_from_compact_to_any_subset}.

We show that appropriately rearranging a set $E$ in the vertical direction we decrease the $s$-perimeter.

More precisely, given a set~$E \subseteq \R^{n + 1}$, we consider the function~$w_E: \R^n \to \R$ defined by
\eqlab{\label{CH:4:rearr_func_def}
w_E(x) := \lim_{R \rightarrow +\infty} \left( \int_{-R}^R \chi_{E}(x, t) \, dt - R \right)
}
for every~$x \in \R^n$, together with its subgraph~$E_\star := \Sg(w_E)$.

Then we have the following result.

\begin{theorem} \label{CH:4:Persdecreases}
Let~$n \ge 1$,~$s \in (0, 1)$, and~$\Omega \subseteq \R^n$ be an open set with Lipschitz boundary. Let~$E \subseteq \R^{n + 1}$ be such that~$E \setminus \Omega^\infty$ is a subgraph and
\begin{equation} \label{CH:4:EboundedinOmega}
\Omega \times (-\infty, -M) \subseteq E \cap \Omega^\infty \subseteq \Omega \times (-\infty, M),
\end{equation}
for some~$M > 0$. Then,
\begin{equation} \label{CH:4:sPerdecreases}
\Per_s(E_\star, \Omega^M) \le \Per_s(E, \Omega^M).
\end{equation}
The inequality is strict unless~$E_\star=E$.
\end{theorem}

The proof of Theorem \ref{CH:4:Persdecreases} can be found in Section \ref{CH:4:REARR_SECTI} and is based on a rearrangement inequality that
we establish for rather general 1-dimensional integral set functions.

Combining the main results of this chapter and exploiting the interior regularity proved in~\cite{CaCo},
we obtain the following Theorem---whose proof is in Section \ref{CH:4:Geom_min_proof_Section}.

\begin{theorem}\label{CH:4:equiv_intro}
Let~$n \ge 1$,~$s \in (0, 1)$,~$\Omega\subseteq\R^n$ be a bounded open set with Lipschitz boundary, and
let~$u\in\W^s(\Omega)$.
Then, the following are equivalent:
\begin{itemize}
\item[(i)] $u$ is a weak solution of~$\h_s u=0$ in~$\Omega$,
\item[(ii)] $u$ is a minimizer of~$\F_s$ in~$\Omega$,
\item[(iii)] $u\in L^\infty_{\loc}(\Omega)$ and~$\Sg(u)$ is locally $s$-minimal in~$\Omega\times\R$,
\item[(iv)] $u\in C^\infty(\Omega)$ and~$u$ is a pointwise solution of~$\h_s u=0$ in~$\Omega$.
\end{itemize}
Moreover, if~$u\in\W^s(\Omega)\cap L^1_{\loc}(\R^n)$, then all of the above are equivalent to:
\begin{itemize}
\item[(v)] $u$ is a viscosity solution of~$\h_s u=0$ in~$\Omega$.
\end{itemize}
\end{theorem}

We also point out the following global version of Theorem~\ref{CH:4:equiv_intro}:
\begin{corollary}\label{CH:4:equiv_intro_global_corollary}
Let~$u\in W^{s,1}_{\loc}(\R^n)$. Then, the following are equivalent:
\begin{itemize}
\item[(i)] $u$ is a viscosity solution of~$\h_s u=0$ in~$\R^n$,
\item[(ii)] $u$ is a weak solution of~$\h_s u=0$ in~$\R^n$,
\item[(iii)] $u$ is a minimizer of~$\F_s$ in~$\Omega$, for every bounded open set~$\Omega\subseteq\R^n$,
\item[(iv)] $u\in L^\infty_{\loc}(\R^n)$ and~$\Sg(u)$ is locally~$s$-minimal in~$\R^{n+1}$,
\item[(v)] $u\in C^\infty(\R^n)$ and~$u$ is a pointwise solution of~$\h_s u=0$ in~$\R^n$.
\end{itemize}
\end{corollary}

In~\cite{graph} the authors observed that if a set~$E$ is locally $s$-minimal in~$\Omega^\infty$,
with~$\Omega\subseteq\R^n$ a bounded open set with~$C^2$ boundary, and~$E=\Sg(\varphi)$
in~$\Co\Omega^\infty$, with~$\varphi\in L^\infty(B_{\tilde R}\setminus\Omega)$
for some~$\tilde R=\tilde R(n,s,\Omega)>0$ big enough, then
\bgs{
\Omega \times (-\infty, -M_0) \subseteq E \cap \Omega^\infty \subseteq \Omega \times (-\infty, M_0),
}
for some~$M_0=M_0(n,s,\Omega,\varphi)>0$. Roughly speaking, this is an a priori bound on the ``vertical
variation'' of the nonlocal minimal surface~$\partial E$ in terms of the exterior data~$\varphi$
and can be thought of as the geometric counterpart of
Theorem~\ref{CH:4:minareboundedthm}.

Exploiting this observation, Theorem~\ref{CH:4:Persdecreases} and Theorem~\ref{CH:4:equiv_intro},
we conclude that there exists a unique locally $s$-minimal set~$E\subseteq\R^{n+1}$ in~$\Omega^\infty$ having as exterior data the subgraph
of~$\varphi$ and ~$E$ is the subgraph of the function~$u\in\W^s_\varphi(\Omega)$
that minimizes~$\F_s$.


\begin{theorem}\label{CH:4:uniqueness_teo}
Let~$\Omega\subseteq\R^n$ be a bounded open set with~$C^2$ boundary and let~$\tilde{R}(n,s,\Omega)$ be as
defined above.
Let~$\varphi:\R^n\to\R$ be such that~$\varphi\in L^\infty(B_{\tilde{R}}\setminus\Omega)$.
If~$E\subseteq\R^{n+1}$ is locally $s$-minimal in~$\Omega^\infty$ and~$E\setminus\Omega^\infty=\Sg(\varphi)\setminus\Omega^\infty$, then~$E=\Sg(u)$, for some~$u\in\B_{M_0}\W^s_\varphi(\Omega)$,
with~$M_0(n,s,\Omega,\varphi)>0$ defined as above. Moreover,~$u$ is the unique minimizer of~$\F_s$
in~$\W^s_\varphi(\Omega)$.
\end{theorem}

We point out that the existence of a locally $s$-minimal set as in Theorem~\ref{CH:4:uniqueness_teo}
is ensured by Corollary \ref{CH:2:loc_min_set_cor}.

In particular, Theorem~\ref{CH:4:uniqueness_teo} extends the result obtained in~\cite{graph}
to a much wider family of exterior data~$\varphi$.
Moreover, it is interesting to observe that, to the best of the authors' knowledge, this also provides the only uniqueness result available for (locally) $s$-minimal sets,
besides the trivial case where the exterior data is an half-space.

\smallskip

We conclude the Introduction with some observations concerning the regularity of the minimizers of the functional~$\F_s$.

Thanks to the interior regularity results proven in~\cite{CaCo} and the fact that
the subgraph of a minimizer~$u$ of~$\F_s$ is locally $s$-minimal,
we know that~$u\in C^\infty(\Omega)$.

On the other hand, we point out that a minimizer~$u$ of~$\F_s$ need not be continuous across the boundary of~$\Omega$
and indeed, in general the subgraph~$\Sg(u)$ sticks to the boundary of the cylinder~$\Omega^\infty$.
For examples of this typically nonlocal phenomenon, we refer in particular to~\cite[Theorems~1.2 and~1.4]{boundary}.
Indeed, the exterior data considered in~\cite[Theorem 1.2]{boundary} is the subgraph of
the function~$\varphi:\R\setminus(-1,1)\to\R$ defined as
\[
\varphi(t):=-M\quad\mbox{if }t\le-1\quad\mbox{and}\quad\varphi(t):=M\quad\mbox{if }t\ge1.
\]
Hence, by Theorem~\ref{CH:4:uniqueness_teo}, we know that there exists a unique locally $s$-minimal set~$E$ with exterior data~$\Sg(\varphi)$,
which is given by~$E=\Sg(u)$,
where~$u\in \W_\varphi^s(-1,1)$ is the minimizer of~$\F_s$. Then~\cite[Theorem 1.2]{boundary}
says that~$u$ does not attain the exterior data~$\varphi$, which is smooth and globally bounded,
in a continuous way, but rather ``sticks'' to the boundary of the cylinder~$(-1,1)\times\R$.

The same behavior is observed in~\cite[Theorem 1.4]{boundary} where the exterior data
can be chosen to be a small, smooth and compactly supported bump function. Again, by Theorem~\ref{CH:4:uniqueness_teo}, we know that the locally $s$-minimal set is given by the subgraph of the function~$u$ which minimizes~$\F_s$. Furthermore, we remark that in this case the
exterior data can be taken to be ``arbitrarily close'' to the constant function 0, so in some sense this kind of phenomenon
is the typical boundary behavior of minimizers of~$\F_s$.

Furthermore, we mention the forthcoming paper~\cite{LuCla}, where this behavior is investigated in the case where the fractional parameter $s$ is small, also in the presence of obstacles.

Nevertheless, even if in general the minimizer of~$\F_s$ is not continuous across the boundary of the domain,
not even when the exterior data is smooth and globally bounded,
we point out that no gap phenomenon occurs, as shown by Proposition~\ref{CH:4:NoLavrentievgap_prop}.

\medskip

Finally, we mention that in Section \ref{CH:4:Appro_Section} we prove some approximation results for subgraphs having (locally) finite fractional perimeter. In particular, exploiting the surprising density result established in \cite{DSV17}, we show that $s$-minimal subgraphs can be appropriately approximated by subgraphs of $\sigma$-harmonic functions, for any fixed $\sigma\in(0,1)$---see Theorem~\ref{CH:4:appro_nonlomin_fracharm_theorem}.

%
%

\section{Preliminary results}

\subsection{Elementary properties of the functions~$g$,~$G$, and~$\G$}

We begin by recalling the following definitions given in the introduction. 
We consider a continuous function~$g: \R \to \R$ satisfying
\begin{align}
\label{CH:4:geven}
g(t) = g(-t) & \quad \mbox{for every } t \in \R, \\
\label{CH:4:gbounds}
0 < g \le 1 & \quad \mbox{in } \R,
\end{align}
and
\begin{equation} \label{CH:4:intglelambda}
\lambda := \int_{0}^{+\infty} t g(t) \, dt < \infty.
\end{equation}
We also observe that
\eqlab{\label{CH:4:gintegr}
\Lambda := \int_\R g(t) \, dt \le 2(\lambda+1)<\infty.
}
As remarked in the Introduction, it is easily seen that the function~$g_s$ defined in~\eqref{CH:4:gdef}
satisfies these assumptions. When considering~$g_s$, we will denote
\bgs{
\Lambda_{n,s}:=\int_\R g_s(t)\,dt.
}

Then, we define
\begin{equation} \label{CH:4:Gdefs}
G(t) := \int_0^t g(\tau) \, d\tau, \qquad \G(t) := \int_0^t G(\tau) \, d\tau = \int_0^t \left( \int_0^\tau g(\sigma) \, d\sigma \right) d\tau
\end{equation}
and
\begin{equation} \label{CH:4:Gbardef}
\overline{G}(t) := \int_{-\infty}^t g(\tau) \, d\tau = \int_{-t}^{+\infty} g(\tau) \, d\tau,
\end{equation}
for every~$t \in \R$. Notice that
\begin{equation} \label{CH:4:GbarG}
\overline{G}(t) = \frac{\Lambda}{2} + G(t) \quad \mbox{for every } t \in \R.
\end{equation}

The following lemma collects the main properties of these functions that will be used in the forthcoming sections.

\begin{lemma} \label{CH:4:gsprop}
The functions~$G$ and~$\G$ are respectively of class~$C^1$ and~$C^2$. Furthermore, the following facts hold true.
\begin{enumerate}[label=$(\alph*)$,leftmargin=*]
\item The function~$G$ is odd, increasing, satisfies~$G(0) = 0$ and
\begin{equation} \label{CH:4:Gbounds}
c_\star \min \{ 1, |t| \} \le |G(t)| \le \min \left\{ \frac{\Lambda}{2}, |t| \right\} \quad \mbox{for every } t \in \R,
\end{equation}
where
\eqlab{\label{CH:4:cstardef}
c_\star=c_\star(g):=\inf_{t\in[0,1]}g(t) > 0.
}
Moreover,
\eqlab{\label{CH:4:lip_G_bla}
|G(t)-G(\tau)|\le|t-\tau|\quad \mbox{for every } t,\tau \in \R.
}
\item The function~$\G$ is even, increasing on~$[0,\infty)$, strictly convex and such that~$\G(0) = 0$. It satisfies
\begin{align}
\label{CH:4:GGbounds}
\frac{c_\star}{2} \min \left\{ |t|, t^2 \right\} & \le \G(t) \le \frac{t^2}{2}, \\
\label{CH:4:GGbetterbound}
\frac{\Lambda}{2} |t| - \lambda & \le \G(t) \le \frac{\Lambda}{2} |t|,
\end{align}
for every~$t \in \R$, and
\eqlab{\label{CH:4:Lip_Gcal}
|\G(t)-\G(\tau)|\le\frac{\Lambda}{2} \, |t-\tau|\quad \mbox{for every } t,\tau \in \R.
}
\end{enumerate}
\end{lemma}
\begin{proof}
Almost all the statements follow immediately from definitions~\eqref{CH:4:Gdefs} and~\eqref{CH:4:Gbardef}. The only properties that require an explicit proof are the lower bounds on~$|G|$ and~$\G$

To obtain the left-hand inequality in~\eqref{CH:4:Gbounds} we assume without loss of generality that~$t \ge 0$ and distinguish between the cases~$t > 1$ and~$t \in [0, 1]$. In the first situation, the claim simply follows by~\eqref{CH:4:Gdefs} along with the monotonicity of~$G$ and~\eqref{CH:4:cstardef}, as indeed
$$
G(t) \ge G(1) =\int_0^1g(t)\,dt \ge c_\star.
$$
Conversely, when~$t \in [0, 1]$ we have
$$
G(t) = \int_0^t g(\tau) \, d\tau \ge c_\star \, t,
$$
thanks again to~\eqref{CH:4:cstardef}.

To get the lower bound in~\eqref{CH:4:GGbounds}, we first notice that we can restrict ourselves to~$t \ge 1$, since the case~$t \in [0, 1]$ can be deduced straight-away from~\eqref{CH:4:Gbounds} and the definition of $\G$. For~$t \ge 1$ we apply~\eqref{CH:4:Gbounds} to compute
\bgs{
\G(t) = \int_0^1 G(\tau) \, d\tau + \int_1^t G(\tau) \, d\tau
\ge c_\star \left( \int_0^1 \tau \, d\tau + \int_1^t d\tau \right)
= \frac{c_\star}{2} \big( 1 + 2 (t - 1) \big) \ge \frac{c_\star}{2} \, t.
}

Finally, to establish the first inequality in~\eqref{CH:4:GGbetterbound}, we recall definitions~\eqref{CH:4:intglelambda}-\eqref{CH:4:Gdefs} and compute, for~$t \ge 0$,
\begin{align*}
\G(t) - \frac{\Lambda}{2} t & = \int_0^t \left( \int_0^\tau g(\sigma) \, d\sigma \right) d\tau - \left( \int_0^{+\infty} g(\sigma) \, d\sigma \right) t = - \int_0^t \left( \int_{\tau}^{+\infty} g(\sigma) \, d\sigma \right) d\tau \\
& = - \int_0^t \left( \int_0^\sigma g(\sigma) \, d\tau \right) d\sigma - \int_t^{+\infty} \left( \int_0^t g(\sigma) \, d\tau \right) d\sigma \\
& = - \int_0^t \sigma g(\sigma) \, d\sigma - t \int_t^{+\infty} g(\sigma) \, d\sigma = - \lambda + \int_t^{+\infty} (\sigma - t) g(\sigma) \, d\sigma \ge - \lambda.
\end{align*}
Note that the third identity follows by Fubini's theorem. The proof of the lemma is thus complete.
\end{proof}

We stress that hypothesis~\eqref{CH:4:intglelambda} has only been used to deduce the left-hand inequality in~\eqref{CH:4:GGbetterbound}. If one drops it, the weaker lower bound
$$
\G(t) \ge \frac{c_\star}{2} |t| - \frac{c_\star}{2} \quad \mbox{for every } t \in \R
$$
can still be easily deduced from~\eqref{CH:4:GGbounds}. This estimate is indeed sufficient for most of the applications presented in the remainder of this chapter. However, we will make crucial use of the finer bound~\eqref{CH:4:GGbetterbound} at some point in the proof of Proposition~\ref{CH:4:Linftylocprop}. Therefore, such result and all those that rely on it need assumption~\eqref{CH:4:intglelambda} to hold.

Note that the function~$g(t) = 1/(1 + t^2)$ fulfills hypotheses~\eqref{CH:4:geven},~\eqref{CH:4:gbounds},~\eqref{CH:4:gintegr}, but not~\eqref{CH:4:intglelambda}. Also, the corresponding second antiderivative~$\G$ does not satisfies the lower bound in~\eqref{CH:4:GGbetterbound} or any bound of the form~$\G(t) \ge \Lambda |t| / 2 - C$ for some constant~$C > 0$.

\subsection{Functional theoretic properties of the fractional area functionals}\label{CH:4:FTPOTFAF_SEC}

In this subsection we introduce the area-type functionals~$\F^M$ and determine some basic properties of the local part~$\A$ and nonlocal part~$\Nl^M$.

First of all, we observe that we can split the functional~$\F$ defined in~\eqref{CH:4:Fcdef} into the two components
\bgs{
\F(u,\Omega)=\A(u,\Omega)+\Nl(u,\Omega),
}
with
\eqlab{\label{CH:4:Adef}
\A(u,\Omega):=\int_\Omega\int_\Omega\G\left(\frac{u(x)-u(y)}{|x-y|}\right)\frac{dx\,dy}{|x-y|^{n-1+s}}
}
and
\bgs{
\Nl(u,\Omega):=2\int_\Omega\int_{\Co\Omega}\G\left(\frac{u(x)-u(y)}{|x-y|}\right)\frac{dx\,dy}{|x-y|^{n-1+s}}.
}
As shown in Lemma~\ref{CH:4:Adomainlem}, in order for the local part~$\A(u,\Omega)$ to be well defined, it is 
necessary and sufficient that~$u\in W^{s,1}(\Omega)$. On the other hand,
for the nonlocal part~$\Nl(u,\Omega)$ to be well defined, we would have to impose some restrictive condition
on the behavior of~$u$ in the whole~$\R^n$---namely~\eqref{CH:4:tartainfty}.

For this reason, given any real number~$M \ge 0$ we define for a function~$u: \R^n \to \R$
the ``truncated'' nonlocal part
\eqlab{
\label{CH:4:NMldef}
\Nl^M(u,\Omega):= \int_\Omega \left\{ \int_{\Co\Omega} \left[ \int_{\frac{-M-u(y)}{|x-y|}}^{\frac{u(x)-u(y)}{|x-y|}}
\overline{G}(t) \, dt + \int^{\frac{M-u(y)}{|x-y|}}_{\frac{u(x)-u(y)}{|x-y|}}
\overline{G}(-t) \, dt \right] \frac{dy}{\kers} \right\} dx,
}
and we introduce the functional
\begin{equation} \label{CH:4:FMdef}
\F^M(u,\Omega) := \A(u,\Omega) + \Nl^M(u,\Omega).
\end{equation}
When~$g=g_s$ we will add the subscript~$s$ to the functionals, that is, we will write~$\A_s,\,\Nl_s^M$
and~$\F_s^M$.

As a motivation for introducing the functionals~$\Nl^M$, we observe that in the geometric situation---that is, when~$g=g_s$---considering the functional~$\Nl_s^M$
in place of~$\Nl_s$ amounts, roughly speaking, to considering the nonlocal part of the fractional perimeter of the subgraph~$\Sg(u)$ in the ``truncated''
cylinder~$\Omega^M=\Omega\times(-M,M)$
instead of considering the nonlocal part of the fractional perimeter in the whole cylinder~$\Omega\times\R$---which would be infinite.
This relationship with the fractional perimeter will be made precise in the forthcoming Subsection~\ref{CH:4:areageom}.

From the functional point of view,
as proved in Lemma~\ref{CH:4:NMdomainlem}, the advantage of considering~$\Nl^M$ instead of~$\Nl$ consists in that we do not need to impose any condition on the function~$u$ outside of the domain~$\Omega$
for~$\Nl^M(u,\Omega)$ to be well defined.

We now proceed to establish the natural domain of definition of the local part~$\A$.
Notice that for the integral defining it to be meaningful (albeit possibly infinite) one only needs~$u$ to be defined in~$\Omega$.

\begin{lemma} \label{CH:4:Adomainlem}
Let~$n \ge 1$,~$s\in(0,1)$,~$\Omega \subseteq \R^n$ be a bounded open set, and let~$u: \Omega \rightarrow \R$
be a measurable function. Then
\begin{equation} \label{CH:4:AGags}
\frac{c_\star}{2} \Big( [u]_{W^{s,1}(\Omega)} - c_s(\Omega) \Big) \le \A(u,\Omega) \le \frac{\Lambda}{2} \, [u]_{W^{s,1}(\Omega)},
\end{equation}
where~$c_\star > 0$ is the constant defined in~\eqref{CH:4:cstardef} and
\begin{equation} \label{CH:4:csOmegadef}
c_s(\Omega) := \frac{\Ha^{n-1}(\mathbb S^{n-1})}{1-s}\,|\Omega| \diam (\Omega)^{1-s}.
\end{equation}
Therefore,
\bgs{
u\in W^{s,1}(\Omega) \quad \mbox{if and only if} \quad \A(u,\Omega) < \infty.
}
\end{lemma}
\begin{proof}
The upper bound in~\eqref{CH:4:AGags} immediately follows by observing that~$\G(t) \le \Lambda |t| / 2$ for every~$t \in \R$, thanks to the right-hand inequality in formula~\eqref{CH:4:GGbounds} of Lemma~\ref{CH:4:gsprop}. To get the lower bound, we recall the left-hand side of~\eqref{CH:4:GGbounds} and compute
$$
\A(u, \Omega) \ge \frac{c_\star}{2} \left( \int_\Omega \int_\Omega \frac{|u(x)-u(y)|}{|x-y|^{n + s}} \, dx\,dy - \int_\Omega \int_\Omega \dkers \right).
$$
The conclusion follows now by Lemma~\ref{CH:4:dumb_kernel_lemma} in Appendix~\ref{CH:4:app}.
Finally, we observe that if~$u$ is a measurable function such that~$[u]_{W^{s,1}(\Omega)}<\infty$,
then~$u\in L^1(\Omega)$ by Lemma~\ref{CH:A:usef_ineq_hit}.
\end{proof}

In the following result we present an equivalent representation for~$\Nl^M(u, \Omega)$, given in terms of the function~$\G$.
We also establish its finiteness when the restriction of~$u$ to~$\Omega$ belongs to the space~$W^{s, 1}(\Omega)$.
Interestingly, no assumption on the behavior of~$u$ outside of~$\Omega$ is needed.

\begin{lemma} \label{CH:4:NMdomainlem}
Let~$n \ge 1$,~$s \in (0, 1)$,~$M \ge 0$,~$\Omega\subseteq\Rn$ be a bounded open set with Lipschitz boundary and
let~$u: \R^n \to \R$ be a measurable function. Then,
\begin{equation} \label{CH:4:Nlbound}
\left| \Nl^M(u,\Omega) \right| \le C\,\Lambda \big( \| u \|_{W^{s, 1}(\Omega)} + M \big),
\end{equation}
where~$\Lambda$ is the positive constant defined in~\eqref{CH:4:gintegr}
and~$C > 0$ is a constant depending only on~$n$,~$s$ and~$\Omega$. Hence,
$$
\left| \Nl^M(u,\Omega) \right| < \infty \quad \mbox{if} \quad u|_\Omega \in W^{s, 1}(\Omega).
$$
Furthermore, we have the identity
\begin{equation} \label{CH:4:nonlocal_explicit}
\begin{aligned}
\Nl^M(u,\Omega) & = \int_{\Omega} \left\{ \int_{\Co \Omega} \left[ 2 \, \G \left( \frac{u(x)-u(y)}{|x-y|} \right) -
\G \left( \frac{M+u(y)}{|x-y|} \right) \right. \right. \\
& \quad \left. \left. -\G \left( \frac{M-u(y)}{|x-y|} \right) \right] \frac{dy}{\kers} \right\} dx + M \Lambda \int_{\Omega} \int_{\Co \Omega} \frac{dx \, dy}{|x-y|^{n+s}}.
\end{aligned}
\end{equation}
\end{lemma}
\begin{proof}
We can assume that~$u|_\Omega\in W^{s,1}(\Omega)$, otherwise~\eqref{CH:4:Nlbound} is trivially satisfied.
Taking advantage of~\eqref{CH:4:GbarG} and of the right-hand inequality in~\eqref{CH:4:Gbounds}, we get that
$$
\left| \Nl^M(u,\Omega) \right| \le 2\Lambda \left[ \int_\Omega \left( |u(x)| \int_{\Co\Omega} \frac{dy}{|x - y|^{n + s}} \right) dx + M \int_{\Omega} \int_{\Co \Omega} \frac{dx \, dy}{|x-y|^{n+s}} \right].
$$
We remark that the last double integral in the previous formula is the $s$-fractional perimeter of~$\Omega$ in~$\R^n$,
which is finite, since~$\Omega$ is bounded and has Lipschitz boundary.
Then~\eqref{CH:4:Nlbound} follows by Corollary~\ref{CH:4:FHI_corollary}.

On the other hand, identity~\eqref{CH:4:nonlocal_explicit} is a simple consequence of definition~\eqref{CH:4:NMldef}, formula~\eqref{CH:4:GbarG} and the symmetry properties of~$G$ and~$\G$.
\end{proof}

We stress that, in order to have~$\Nl^M(u, \Omega)$ finite, the requirement~$u|_\Omega \in W^{s, 1}(\Omega)$ is far from being optimal. In fact, as the previous proof showed, it suffices that~$u|_\Omega$ lies in a suitable weighted~$L^1$ space over~$\Omega$---that contains for instance~$L^\infty(\Omega)$. Nevertheless, such a requirement does not limit our analysis, since it is needed to have~$\A(u, \Omega)$ finite, according to Lemma~\ref{CH:4:Adomainlem}. We inform the interested reader that a more precise result on the natural domain of definition of~$\Nl_s^M$ will be provided by Lemma~\ref{CH:4:per_of_subgraph_lem2} in the forthcoming Subsection~\ref{CH:4:areageom}.

Furthermore, we observe that if~$u:\R^n\to\R$ is such that~$u\in L^\infty(\Omega)$ and~$M\ge\|u\|_{L^\infty(\Omega)}$,
then~$\Nl^M(u,\Omega)\ge0$---since the integrand inside the square brackets in~\eqref{CH:4:NMldef} is non-negative.
On the other hand, we remark that in general the nonlocal part~$\Nl^M(\,\cdot\,,\Omega)$ can assume also negative values,
as proved in the following Example~\ref{CH:4:Exe_neg_nonloc}.

\begin{example}\label{CH:4:Exe_neg_nonloc}
Let~$n \ge 1$,~$s \in (0, 1)$,~$M \ge 0$,~$\Omega\subseteq\Rn$ be a bounded open set with Lipschitz boundary.
There exists a positive
constant~$C=C(n,s,\Omega,g,M)>0$ big enough such that, if~$u:\R^n\to\R$ is the constant function~$u\equiv T$,
for some~$T\ge C$, then
\bgs{
\F^M(u,\Omega)=\Nl^M(u,\Omega)<0.
}
\begin{proof}
Let us fix~$R>0$ such that~$\Omega\Subset B_R$. By identity~\eqref{CH:4:nonlocal_explicit}
and recalling that~$\G\ge0$, we obtain
\bgs{
	\Nl^M(u,\Omega) & = -\int_{\Omega} \left\{ \int_{\Co \Omega} \left[
	\G \left( \frac{M+T}{|x-y|} \right) +\G \left( \frac{M-T}{|x-y|} \right) \right] \frac{dy}{\kers} \right\} dx \\
	&\qquad\qquad\qquad
	+ M \Lambda \int_{\Omega} \int_{\Co \Omega} \frac{dx \, dy}{|x-y|^{n+s}}\\
	&
	\le-\int_{\Omega} \int_{\Co \Omega}
	\G \left( \frac{M+T}{|x-y|} \right)\frac{dx\,dy}{\kers}
	+ M \Lambda \int_{\Omega} \int_{\Co \Omega} \frac{dx \, dy}{|x-y|^{n+s}}\\
	&
	\le-\int_{\Omega} \int_{B_R\setminus \Omega}
	\G \left( \frac{M+T}{|x-y|} \right)\frac{dx\,dy}{\kers}
	+ M \Lambda \int_{\Omega} \int_{\Co \Omega} \frac{dx \, dy}{|x-y|^{n+s}}.
}
By exploiting~\eqref{CH:4:GGbounds}, the fact that~$\Omega$ is bounded and has Lipschitz boundary---hence it has finite $s$-perimeter---and
Lemma~\ref{CH:4:dumb_kernel_lemma}, we find that
\bgs{
\int_{\Omega} \int_{B_R\setminus \Omega}
\G \left( \frac{M+T}{|x-y|} \right)\frac{dx\,dy}{\kers}&\ge
\frac{c_\star}{2}
\int_\Omega\int_{B_R\setminus \Omega}
\left[\frac{M+T}{|x-y|}-1\right]\frac{dx\,dy}{\kers}\\
&
=C_1(M+T)-C_2,
}
with~$C_1,\,C_2>0$ depending only on~$n,\,s,\,\Omega$ and~$g$. Therefore,
\bgs{
\Nl^M(u,\Omega)\le-C_1(M+T)+C_2+M \Lambda \int_{\Omega} \int_{\Co \Omega} \frac{dx \, dy}{|x-y|^{n+s}},
}
which is negative, provided we take~$T>0$
big enough. This concludes the proof.
\end{proof}
\end{example}

We collect the results of Lemmas~\ref{CH:4:Adomainlem} and~\ref{CH:4:NMdomainlem} in the following unifying statement.

\begin{lemma} \label{CH:4:FMdomainlem}
Let~$n \ge 1$,~$s \in (0, 1)$,~$M \ge 0$,~$\Omega\subseteq\Rn$ be a bounded open set with Lipschitz boundary
and let~$u\in\W^s(\Omega)$. Then,~$\F^M(u, \Omega)$ is finite and it holds
$$
\left|\F^M(u, \Omega)\right| \le C\,\Lambda \left( \| u \|_{W^{s, 1}(\Omega)} + M \right),
$$
for some constant~$C > 0$ depending only on~$n$,~$s$ and~$\Omega$.
\end{lemma}

We conclude this subsection by specifying the convexity properties enjoyed by
the functionals~$\A$,~$\Nl^M$, and~$\F^M$.

\begin{lemma}\label{CH:4:conv_func}
Let~$s \in (0, 1)$ and let~$\Omega\subseteq\Rn$ be a bounded open set with Lipschitz boundary. The following facts hold true:
\begin{enumerate}[label=$(\roman*)$,leftmargin=*]
\item The functional~$\A(\,\cdot\,,\Omega)$ is convex on~$W^{s, 1}(\Omega)$.
\item Given any~$M \ge 0$ and measurable function~$\varphi: \Co\Omega \to \R$,
the functionals~$\Nl^M(\,\cdot\,, \Omega)$ and~$\F^M(\,\cdot\,,\Omega)$
are strictly convex on the space~$\W^s_\varphi(\Omega)$
defined in~\eqref{CH:4:domain_w_data}.
\end{enumerate}
\end{lemma}
\begin{proof}
The convexity of the functionals is an immediate consequence of the (strict) convexity of~$\G$ warranted by Lemma~\ref{CH:4:gsprop}. We point out that the convexity of~$\Nl^M(\,\cdot\,, \Omega)$ is due also to the fact that the second and third summands
appearing inside square brackets in the representation~\eqref{CH:4:nonlocal_explicit} are constant on~$\W_\varphi^s(\Omega)$.
Indeed, given~$u,\,v\in\W^s_\varphi(\Omega)$ and~$t\in(0,1)$, we have the identity
\eqlab{\label{CH:4:conv_tarta}
\Nl^M&(tu+(1-t)v,\Omega)-t\,\Nl^M(u,\Omega)-(1-t)\Nl^M(v,\Omega)\\
&
=2\int_\Omega\bigg\{\int_{\Co\Omega}\bigg[\G\left(t\frac{u(x)-\varphi(y)}{|x-y|}+(1-t)\frac{v(x)-\varphi(y)}{|x-y|}\right)
-t\,\G\left(\frac{u(x)-\varphi(y)}{|x-y|}\right)\\
&
\qquad\qquad\qquad-(1-t)\G\left(\frac{v(x)-\varphi(y)}{|x-y|}\right)\bigg]\frac{dy}{|x-y|^{n-1+s}}\bigg\}dx,
}
and the convexity of~$\G$ guarantees that the integrand in the double integral above is nonpositive.
Furthermore, the strict convexity of~$\G$ implies that the quantity in~\eqref{CH:4:conv_tarta} is equal to zero if and only if
\bgs{
\frac{u(x)-\varphi(y)}{|x-y|}=\frac{v(x)-\varphi(y)}{|x-y|}\quad\mbox{for a.e. }(x,y)\in\Omega\times\Co\Omega,
}
i.e. if and only if~$u=v$ almost everywhere in~$\Omega$---and hence in~$\R^n$.
\end{proof}

\subsection{Geometric properties of the fractional area functionals} \label{CH:4:areageom}

This subsection is devoted to the description of a few geometric properties enjoyed by~$\A_s$,~$\Nl_s^M$ and~$\F_s^M$. More specifically, in this subsection we consider the case~$g=g_s$ and we show the connection existing between
the fractional perimeter~$\Per_s$ and these
functionals, that ultimately motivates their introduction.

First of all, we remark that we can split the $s$-perimeter into its local and nonlocal parts, as
\bgs{
\Per_s(E,\Op)=\Per_s^L(E,\Op)+\Per_s^{NL}(E,\Op),
}
with
\bgs{
\Per_s^L(E,\Op):=\Ll_s(E\cap\Op,\Co E\cap\Op)=\frac{1}{2}[\chi_E]_{W^{s,1}(\Op)}.
}

We begin with a result that deals with the local part~$\A_s$.

Before going on, we recall that by Lemma~\ref{CH:A:usef_ineq_hit} we know that
a function having finite~$W^{s,1}(\Omega)$-seminorm also belongs to~$L^1(\Omega)$.

\begin{lemma} \label{CH:4:per_of_subgraph_lem1}
Let~$s \in (0, 1)$,~$\Omega \subseteq \R^n$ be a bounded open set and let~$u:\Omega\to\R$ be a measurable function. Then,
\bgs{
u \in W^{s, 1}(\Omega) \quad \mbox{if and only if} \quad \Per_s^L \left( \Sg(u), \Omega^\infty \right) < \infty.
}
In particular, it holds
\begin{equation} \label{CH:4:PerAid}
\Per_s^L \left( \Sg(u),\Omega^\infty \right)=\A_s(u,\Omega) + \Per_s^L \left( \{x_{n+1} < 0\},\Omega^\infty \right).
\end{equation}
\end{lemma}
\begin{proof}
Using Lebesgue's monotone convergence theorem, we write
$$
\Per_s^L \left( \Sg(u), \Omega^\infty \right) = \lim_{\delta \rightarrow 0^+} \iint_{\Omega^2 \cap \{ |x - y| > \delta \}} dx\,dy \int_{-\infty}^{u(x)} dx_{n + 1} \int_{u(y)}^{+\infty} \frac{dy_{n + 1}}{|X - Y|^{n + 1 + s}}.
$$
Fix any small~$\delta > 0$ and let~$(x, y) \in \Omega^2 \cap \{ |x - y| > \delta \}$. Shifting variables, we see that
$$
\int_{-\infty}^{u(x)} dx_{n + 1} \int_{u(y)}^{+\infty} \frac{dy_{n + 1}}{|X - Y|^{n + 1 + s}} = \int_{-\infty}^{u(x) - u(y)} dx_{n + 1} \int_{0}^{+\infty} \frac{dy_{n + 1}}{|X - Y|^{n + 1 + s}},
$$
so that
\begin{align*}
\Per_s^L \left( \Sg(u), \Omega^\infty \right) & = \lim_{\delta \rightarrow 0^+} \iint_{\Omega^2 \cap \{ |x - y| > \delta \}} dx\,dy \int_{0}^{u(x) - u(y)} dx_{n + 1} \int_{0}^{+\infty} \frac{dy_{n + 1}}{|X - Y|^{n + 1 + s}} \\
& \quad + \Per_s^L \left( \{ x_{n + 1} < 0 \}, \Omega^\infty \right).
\end{align*}
After a renormalization of both variables~$x_{n + 1}$ and~$y_{n + 1}$, we have
$$
\int_{0}^{u(x) - u(y)} dx_{n + 1} \int_{0}^{+\infty} \frac{dy_{n + 1}}{|X - Y|^{n + 1 + s}} = \frac{1}{|x - y|^{n - 1 + s}} \int_{0}^{\frac{u(x) - u(y)}{|x - y|}} dt \int_{0}^{+\infty} \frac{dr}{\left[ 1 + (r - t)^2 \right]^{\frac{n + 1 + s}{2}}}.
$$
Changing coordinates once again and recalling definition~\eqref{CH:4:gdef}, we obtain that
\begin{align*}
\int_{0}^{u(x) - u(y)} dx_{n + 1} \int_{0}^{+\infty} \frac{dy_{n + 1}}{|X - Y|^{n + 1 + s}} & = \frac{1}{|x - y|^{n - 1 + s}} \int_{0}^{\frac{u(x) - u(y)}{|x - y|}} \left( \int_{-t}^{+\infty} \frac{d\tau'}{\left[ 1 + (\tau')^2 \right]^{\frac{n + 1 + s}{2}}} \right) dt \\
& = \frac{1}{|x - y|^{n - 1 + s}} \int_{0}^{\frac{u(x) - u(y)}{|x - y|}} \left( \int_{-\infty}^{t} g_s(\tau) \, d\tau \right) dt.
\end{align*}
By~\eqref{CH:4:Gdefs} and~\eqref{CH:4:gintegr}, we get
\begin{align*}
\int_{0}^{\frac{u(x) - u(y)}{|x - y|}} \left( \int_{-\infty}^{t} g_s(\tau) \, d\tau \right) dt = \frac{\Lambda_{n, s}}{2} \frac{u(x) - u(y)}{|x - y|} + \G_s \left( \frac{u(x) - u(y)}{|x - y|} \right).
\end{align*}
Since, by symmetry,
$$
\iint_{\Omega^2 \cap \{ |x - y| > \delta \}} \frac{u(x) - u(y)}{|x - y|} \frac{dx\,dy}{ |x - y|^{n - 1 + s}} = 0,
$$
we conclude that
$$
\Per_s^L \left( \Sg(u), \Omega^\infty \right) = \lim_{\delta \rightarrow 0^+} \iint_{\Omega^2 \cap \{ |x - y| > \delta \}} \G_s \left( \frac{u(x) - u(y)}{|x - y|} \right) dx\,dy + \Per_s^L \left( \{ x_{n + 1} < 0 \}, \Omega^\infty \right).
$$
The claim of the lemma now follows by taking advantage once again of Lebesgue's monotone convergence theorem and recalling definition~\eqref{CH:4:Adef}.
\end{proof}

\begin{lemma} \label{CH:4:per_of_subgraph_lem2}
Let~$s \in (0, 1)$,~$\Omega \subseteq \R^n$ be a bounded open set with Lipschitz boundary and let~$u: \R^n \to \R$ be
such that~$u|_\Omega \in L^\infty(\Omega)$. Then, for any~$M \ge \| u \|_{L^\infty(\Omega)}$, the quantity~$\Nl_s^M(u, \Omega)$ is finite and it holds
\begin{equation} \label{CH:4:NMid}
\Nl_s^M(u,\Omega)=\Ll_s \left( \Sg(u)\cap\Omega^M, \Co\Sg(u)\setminus\Omega^\infty \right)
+\Ll_s \left( \Co\Sg(u)\cap\Omega^M, \Sg(u)\setminus\Omega^\infty \right).
\end{equation}
\end{lemma}
\begin{proof}
Thanks to the fact that~$M \ge \| u \|_{L^\infty(\Omega)}$, we write
\begin{align*}
\Ll_s \left( \Sg(u)\cap\Omega^M, \Co\Sg(u)\setminus\Omega^\infty \right) & = \int_{\Omega} dx \int_{\Co \Omega} dy \int_{-M}^{u(x)} dx_{n + 1} \int_{u(y)}^{+\infty} \frac{dy_{n + 1}}{|X - Y|^{n + 1 + s}}, \\
\Ll_s \left( \Co\Sg(u)\cap\Omega^M, \Sg(u)\setminus\Omega^\infty \right) & = \int_{\Omega} dx \int_{\Co \Omega} dy \int_{u(x)}^{M} dx_{n + 1} \int_{-\infty}^{u(y)} \frac{dy_{n + 1}}{|X - Y|^{n + 1 + s}}.
\end{align*}
By arguing as in the proof of Lemma~\ref{CH:4:per_of_subgraph_lem1} and recalling definitions~\eqref{CH:4:gdef},~\eqref{CH:4:Gdefs}
and~\eqref{CH:4:Gbardef}, we then have
\begin{align*}
\int_{-M}^{u(x)} dx_{n + 1} \int_{u(y)}^{+\infty} \frac{dy_{n + 1}}{|X - Y|^{n + 1 + s}} & = \int_{-M - u(y)}^{u(x) - u(y)} dx_{n + 1} \int_{0}^{+\infty} \frac{dy_{n + 1}}{|X - Y|^{n + 1 + s}} \\
& = \frac{1}{|x - y|^{n - 1 + s}} \int_{\frac{-M - u(y)}{|x - y|}}^{\frac{u(x) - u(y)}{|x - y|}} dx_{n + 1} \int_{- x_{n + 1}}^{+\infty} \frac{d\tau}{(1 + \tau^2)^{\frac{n + 1 + s}{2}}} \\
& = \frac{1}{|x - y|^{n - 1 + s}} \int_{\frac{-M - u(y)}{|x - y|}}^{\frac{u(x) - u(y)}{|x - y|}} \overline{G}_s(t) \, dt
\end{align*}
for every~$x \in \Omega$ and~$y \in \Co \Omega$. Hence,
$$
\Ll_s \left( \Sg(u)\cap\Omega^M, \Co\Sg(u)\setminus\Omega^\infty \right) = \int_{\Omega} dx \int_{\Co \Omega} dy \left( \frac{1}{|x - y|^{n - 1 + s}} \int_{\frac{-M - u(y)}{|x - y|}}^{\frac{u(x) - u(y)}{|x - y|}} \overline{G}_s(t) \, dt \right).
$$
Similarly,
$$
\Ll_s \left( \Co\Sg(u)\cap\Omega^M, \Sg(u)\setminus\Omega^\infty \right) = \int_{\Omega} dx \int_{\Co \Omega} dy \left( \frac{1}{|x - y|^{n + 1 + s}} \int_{\frac{u(x) - u(y)}{|x - y|}}^{\frac{M - u(y)}{|x - y|}} \overline{G}_s(-t) \, dt \right).
$$
By combining the last two identities and recalling definition~\eqref{CH:4:NMldef}, we are led to~\eqref{CH:4:NMid}.
\end{proof}

\begin{prop} \label{CH:4:per_of_subgraph_prop}
Let~$s \in (0, 1)$ and~$\Omega \subseteq \R^n$ be a bounded open set with Lipschitz boundary. Let~$u: \R^n \to \R$ be such that~$u|_\Omega \in L^\infty(\Omega)$ and take~$M \ge \| u \|_{L^\infty(\Omega)}$. Then,
\begin{equation} \label{CH:4:u_iff_Sgu}
u|_\Omega \in W^{s,1}(\Omega) \quad \mbox{if and only if} \quad \Per_s \left( \Sg(u), \Omega^M \right) < \infty.
\end{equation}
In particular, it holds
\begin{equation} \label{CH:4:per_of_subgraph}
\Per_s \left( \Sg(u), \Omega^M \right) = \F^M_s(u,\Omega) + \kappa_{\Omega, M},
\end{equation}
where~$\kappa_{\Omega, M}$ is the non-negative constant given by
\begin{equation} \label{CH:4:kappadef}
\kappa_{\Omega, M} := \Per_s^L \left( \{ x_{n + 1} < 0 \}, \Omega^\infty \right) - \Per_s^L \left( \{ x_{n + 1} < 0 \}, \Omega^\infty \setminus \Omega^M \right).
\end{equation}
\end{prop}
\begin{proof}
The proposition is an almost immediate consequence of Lemmas~\ref{CH:4:per_of_subgraph_lem1} and~\ref{CH:4:per_of_subgraph_lem2}. First, we observe that the following identities are true:
\begin{align*}
\Ll_s \left( \Sg(u) \cap \Omega^M, \Omega^M \setminus \Sg(u) \right) & = \int_{\Omega} dx \int_{\Omega} dy \int_{-M}^{u(x)} dx_{n + 1} \int_{u(y)}^M \frac{dy_{n + 1}}{|X - Y|^{n + 1 + s}}, \\
\Ll_s \left( \Sg(u) \cap \Omega^M, \Co \Sg(u) \setminus \Omega^M \right) & = \int_{\Omega} dx \int_{\Omega} dy \int_{-M}^{u(x)} dx_{n + 1} \int_{M}^{+\infty} \frac{dy_{n + 1}}{|X - Y|^{n + 1 + s}} \\
& \quad + \int_{\Omega} dx \int_{\Co \Omega} dy \int_{-M}^{u(x)} dx_{n + 1} \int_{u(y)}^{+\infty} \frac{dy_{n + 1}}{|X - Y|^{n + 1 + s}}, \\
\Ll_s \left( \Omega^M \setminus \Sg(u), \Sg(u) \setminus \Omega^M \right) & = \int_{\Omega} dx \int_{\Omega} dy \int_{u(x)}^M dx_{n + 1} \int_{-\infty}^{-M} \frac{dy_{n + 1}}{|X - Y|^{n + 1 + s}} \\
& \quad + \int_{\Omega} dx \int_{\Co \Omega} dy \int_{u(x)}^M dx_{n + 1} \int_{-\infty}^{u(y)} \frac{dy_{n + 1}}{|X - Y|^{n + 1 + s}}.
\end{align*}
Note that we took advantage of the fact that~$M \ge \| u \|_{L^\infty(\Omega)}$ in order to obtain the above formulas. In light of this, it is not hard to see that
\begin{align*}
\Per_s \left( \Sg(u), \Omega^M \right) & = \Per_s^L \left( \Sg(u),\Omega^\infty \right) - \Per_s^L \left( \{ x_{n + 1} < 0 \}, \Omega^\infty \setminus \Omega^M \right) \\
& \quad + \Ll_s \left( \Sg(u) \cap \Omega^M, \Co \Sg(u) \setminus \Omega^\infty \right) + \Ll_s \left( \Co \Sg(u) \cap \Omega^M, \Sg(u) \setminus \Omega^\infty \right).
\end{align*}
Identity~\eqref{CH:4:per_of_subgraph} follows by recalling definition~\eqref{CH:4:FMdef} and applying~\eqref{CH:4:PerAid} and~\eqref{CH:4:NMid}.
\end{proof}

\subsection{Some facts about the Euler-Lagrange operator}\label{CH:4:SFATELO_SEC}

We collect here some observations about the nonlocal integrodifferential operator~$\h$, which is formally defined on a function~$u: \R^n \to \R$ at a point~$x \in \R^n$ by
$$
\h u(x) := 2 \, \PV \int_{\R^n} G \left( \frac{u(x)-u(y)}{|x-y|} \right) \frac{dy}{|x-y|^{n+s}}.
$$
We begin by introducing the following useful notation
\bgs{
\delta_g(u,x;\xi):=G\Big(\frac{u(x)-u(x+\xi)}{|\xi|}\Big)-G\Big(\frac{u(x-\xi)-u(x)}{|\xi|}\Big),
}
and we observe that by symmetry we can write
\eqlab{\label{CH:4:delta_op}
\h u(x)=\PV\int_{\R^n}\frac{\delta_g(u,x;\xi)}{|\xi|^{n+s}}\,d\xi.
}
From now on, unless otherwise stated, we will always consider~$\h u(x)$ as written in~\eqref{CH:4:delta_op}.

We remark that when~$g=g_s$, we will write~$\h_s$ for the corresponding nonlocal operator.
From a geometric standpoint, the quantity~$\h_s u$ describes the~$s$-mean curvature of the subgraph of~$u$. Indeed, it holds
\eqlab{\label{CH:4:Geom_ELop_curvature}
\I_s[\Sg(u)](x,u(x)) = \h_s u(x)
}
for every~$x \in \R^n$ at which~$u$ is of class~$C^{1, \alpha}$, for some~$\alpha>s$---see~\cite[Appendix~B.1]{BLV16} for the details of this computation.

We also define
\bgs{
\h^{\geq r}u(x):=\int_{\R^n\setminus B_r}\frac{\delta_g(u,x;\xi)}{|\xi|^{n+s}}\,d\xi,
\qquad\forall\,r>0,
}
so that
\[\h u(x)=\lim_{r\to0^+}\h^{\geq r}u(x).\]

\begin{remark}\label{CH:4:rmk_tail}
Let~$u:\R^n\to\R$. Then~$\h^{\geq r}u(x)$ is finite for every~$x\in\R^n$ and~$r>0$.
Indeed, exploiting the boundedness of~$G$ we find
\[\Big|\frac{\delta_g(u,x;\xi)}{|\xi|^{n+s}}\Big|\leq\frac{\Lambda}{|\xi|^{n+s}},\]
which is summable in~$\R^n\setminus B_r$. In particular
\[|\h^{\geq r}u(x)|\leq\frac{n\omega_n}{s}\Lambda\,r^{-s}.\]
\end{remark}

One of the main advantages of writing the nonlocal operator~$\h u(x)$ as in~\eqref{CH:4:delta_op}
is that the integral is well defined in the classical sense, provided~$u$
is regular enough around~$x$.

\begin{lemma}\label{CH:4:classical_form_reg_func}
Let~$s\in(0,1)$ and let~$u:\R^n\to\R$ such that~$u\in C^{1,\gamma}(B_r(x))$, for some~$x\in\R^n$, $r>0$
and~$\gamma\in(s,1]$. Then
\[\h^{<\varrho}u(x):=\int_{B_\varrho}\frac{\delta_g(u,x;\xi)}{|\xi|^{n+s}}\,d\xi\]
is well defined for every~$\varrho>0$ and
\eqlab{\label{CH:4:formula_regular_func}
\h u(x)=\h^{<\varrho}u(x)+\h^{\geq\varrho}u(x)=\int_{\R^n}\frac{\delta_g(u,x;\xi)}{|\xi|^{n+s}}\,d\xi.
}
\end{lemma}
\begin{proof}
We begin by proving that
\eqlab{\label{CH:4:local_C1gamma}
\Big|\frac{u(x+\xi)+u(x-\xi)-2u(x)}{|\xi|}\Big|
\leq2^\gamma\|u\|_{C^{1,\gamma}(B_r(x))}|\xi|^\gamma,\quad\forall\,\xi\in B_r\setminus\{0\}.
}
Indeed, by the Mean Value Theorem we have
\[u(x+\xi)-u(x)=\nabla u(x+t\xi)\cdot\xi\quad\textrm{and}\quad u(x-\xi)-u(x)=\nabla u(x-\tau\xi)\cdot(-\xi),\]
for some~$t,\,\tau\in[0,1]$. Thus
\bgs{
&\Big|\frac{u(x+\xi)+u(x-\xi)-2u(x)}{|\xi|}\Big|=\Big|\frac{\nabla u(x+t\xi)\cdot\xi-\nabla u(x-\tau\xi)\cdot\xi}{|\xi|}\Big|\\
&\qquad\qquad
\leq|\nabla u(x+t\xi)-\nabla u(x-\tau\xi)|
\leq[\nabla u]_{C^\gamma(B_r(x))}|(t+\tau)\xi|^\gamma\\
&\qquad\qquad
\leq2^\gamma\|u\|_{C^{1,\gamma}(B_r(x))}|\xi|^\gamma,
}
as claimed.

Now we remark that, thanks to Remark~\ref{CH:4:rmk_tail},
in order to prove the lemma it is enough to show that~$\delta_g(u,x;\xi)|\xi|^{-n-s}$ is summable in~$B_\varrho$,
for~$\varrho>0$ small enough. For this, just notice that by~\eqref{CH:4:lip_G_bla}
we have
\[|\delta_g(u,x;\xi)|\leq\Big|\frac{u(x+\xi)+u(x-\xi)-2u(x)}{|\xi|}\Big|.\]
Then the conclusion follows from~\eqref{CH:4:local_C1gamma}.
\end{proof}

We stress that the right hand side of~\eqref{CH:4:formula_regular_func} is defined in the classical sense,
not as a principal value. Also notice that, thanks to Remark~\ref{CH:4:rmk_tail}, we need not ask any growth condition for~$u$ at infinity.

When~$u$ is not regular enough around~$x$, the quantity~$\h u(x)$ is in general not well-defined, due to the lack of cancellation required for the principal value to converge. Nevertheless, as already observed in the Introduction,
we can understand the operator~$\h$ as defined in the following weak (distributional) sense. Given a function~$u:\R^n \to \R$, we set
\eqlab{\label{CH:4:weak_opossum_curv}
\langle \h u, v\rangle:=\int_{\Rn}\int_{\Rn} G \left( \frac{u(x)-u(y)}{|x-y|} \right) \big( v(x)-v(y) \big) \, \frac{dx\,dy}{|x-y|^{n+s}}
}
for every $v\in C^\infty_c(\Rn)$. More generally, it is immediate to see that~\eqref{CH:4:weak_opossum_curv} is well-defined for every~$v:\Rn \to \R$ such that~$[v]_{W^{s,1}(\Rn)}<\infty$. Indeed, taking advantage of the boundedness of~$G$, one has that
\begin{equation} \label{CH:4:hsucont}
|\langle \h u ,v \rangle|\le\frac{\Lambda}{2} \, [v]_{W^{s,1}(\Rn)},
\end{equation}
with~$\Lambda$ as in~\eqref{CH:4:gintegr}. Hence,~$\h u$ induces a continuous linear functional on~$W^{s, 1}(\R^n)$, that is
$$
\langle \h u, - \rangle \in \left( W^{s, 1}(\R^n) \right)^*.
$$
Remarkably, this holds for every measurable function~$u: \R^n \to \R$, regardless of its regularity.

%
%

Estimate~\eqref{CH:4:hsucont} says that the pairing~$(u, v) \mapsto \langle \h u, v \rangle$ is continuous in the second component~$v$, with respect to the~$W^{s, 1}(\R^n)$ topology. The next lemma shows that we also have continuity in~$u$ with respect to convergence a.e.~in~$\R^n$.


\begin{lemma}\label{CH:4:easylemma}
Let~$u_k,\,u:\R^n\to\R$ be such that~$u_k\to u$ almost everywhere in~$\R^n$ and let~$v\in W^{s,1}(\R^n)$.
Then
\bgs{
\lim_{k\to\infty}\langle\h u_k,v\rangle=\langle\h u,v\rangle.
}
\end{lemma}

Lemma~\ref{CH:4:easylemma} is a simple consequence of Lebesgue's dominated convergence theorem, thanks to the fact that~$\| G \|_{L^\infty(\R)} = \Lambda / 2$.

The next result shows that the nonlocal mean curvature operator~$\h$ naturally arises when computing the Euler-Lagrange equation associated to the fractional area functional.

\begin{lemma} \label{CH:4:1varlem}
Let~$n \ge 1$,~$s \in (0, 1)$,~$M\ge0$,~$\Omega \subseteq \R^n$ be a bounded open set with Lipschitz boundary, and~$u \in \W^s(\Omega)$. Then,
\begin{equation} \label{CH:4:MEL}
\left. \frac{d}{d\eps} \right|_{\eps=0} \F^M(u+\eps v, \Omega) = \langle \h u, v \rangle \quad \mbox{for every } v \in \W^s_0(\Omega).
\end{equation}
\end{lemma}
\begin{proof}
First, notice that~$u + \varepsilon v \in \W^s(\Omega)$ for every~$\varepsilon \in \R$. Hence, by Lemma~\ref{CH:4:FMdomainlem}, both~$\F^M(u, \Omega)$ and~$\F^M(u + \varepsilon v, \Omega)$ are finite. Now, by Lagrange's mean value theorem, there exists a function~$\tilde{\tau}_\varepsilon: \R \times \R \to [-|\varepsilon|, |\varepsilon|]$ such that
$$
\G \left( A + \varepsilon B \right) - \G \left( A \right) = \varepsilon \, G \left( A + \tilde{\tau}_\varepsilon(A, B) B \right) B
$$
for every~$A, B \in \R$. As~$v = 0$ in~$\Co \Omega$, calling
$$
\tau_\varepsilon(x, y) := \tilde{\tau}_\varepsilon \left( \frac{u(x)- u(y)}{|x - y|}, \frac{v(x)- v(y)}{|x - y|} \right) \quad \mbox{for every } x, y \in \R^n, 
$$
we have
\begin{align*}
& \F^M(u+\eps v,\Omega) - \F^M(u,\Omega) \\
& \hspace{20pt} = \varepsilon \int_{\R^n} \int_{\R^n} G \left( \frac{ u(x)- u(y) }{|x - y|} + \tau_\varepsilon(x, y) \frac{ v(x)- v(y) }{|x - y|} \right) \big( v(x)- v(y) \big) \, \frac{dx\,dy}{|x - y|^{n + s}}.
\end{align*}
Since~$G$ is bounded,~$v \in \W_0^s(\Omega)$, and~$|\tau_\varepsilon| \le \varepsilon$, we may conclude the proof using Lebesgue's dominated convergence theorem.
\end{proof}

For more details about the Euler-Lagrange equation of minimizers, we refer to
Lemma~\ref{CH:4:weak_implies_min_lemma}.

%
%

\section{Viscosity implies weak}
\label{CH:4:ViscWeak_Sec}

\subsection{Viscosity (sub)solutions}

We are interested in viscosity solutions of the equation
\[
\left\{\begin{array}{cc}
\h u=f & \mbox{in }\Omega,\\
u=\varphi & \mbox{in }\R^n\setminus\Omega.
\end{array}\right.
\]
We will use $C^{1,1}$ functions as test functions.
First we point out the following easy remark.
\begin{remark}\label{CH:4:mah}
Let $u,\,v:\R^n\to\R$ be such that
\[u(x_0)=v(x_0)\qquad\textrm{and}\qquad u(x)\leq v(x)\quad\forall\,x\in\R^n.\]
Then
\[\delta_g(u,x_0;\xi)\geq\delta_g(v,x_0;\xi)\qquad\forall\,\xi\in\R^n,\]
hence also
\[\h u(x_0)\geq\h v(x_0).\]
Indeed, it is enough to notice that
\[\delta_g(u,x_0;\xi)=G\Big(\frac{u(x_0)-u(x_0+\xi)}{|\xi|}\Big)+G\Big(\frac{u(x_0)-u(x_0-\xi)}{|\xi|}\Big),\]
and recall that $G$ is increasing.
\end{remark}

\begin{defn}\label{CH:4:VISCO_DEF}
Let $\Omega\subseteq\R^n$ be a bounded open set and let $f\in C(\overline{\Omega})$.
We say that a function $u:\R^n\to\R$ is a (viscosity) subsolution of $\h u=f$ in $\Omega$, and we write
\[\h u\leq f\qquad\textrm{in }\Omega,\]
if $u$ is upper semicontinuous in $\Omega$ and whenever the following happens:
\begin{itemize}
\item[(i)] $x_0\in\Omega$,
\item[(ii)] $v\in C^{1,1}(B_r(x_0)),$ for some $r<d(x_0,\partial\Omega)$,
\item[(iii)] $v(x_0)=u(x_0)$ and $v(y)\geq u(y)$ for every $y\in B_r(x_0)$,
\end{itemize}
then if we define
\[\tilde v(x):=\left\{\begin{array}{cc} v(x) & \textrm{if }x\in B_r(x_0),\\
u(x) & \textrm{if }x\in\R^n\setminus B_r(x_0),\end{array}\right.\]
we have
\[\h\tilde v(x_0)\leq f(x_0).\]
A supersolution is defined similarly. A (viscosity) solution is a function $u:\R^n\to\R$ which is continuous in $\Omega$
and which is both a subsolution and a supersolution.
\end{defn}

From now on, we will concentrate on viscosity subsolutions, the corresponding statements for supersolutions being obtained by considering $-u$ in place of $u$.\\
Unless otherwise stated, $f$ will always be supposed to be continuous in the closure of $\Omega$.\\

A crucial observation is the following.\\
Roughly speaking, for a function $u$ to be touched from above
at some point $x_0$ by a $C^{1,1}$ function means that $u$ is $C^{1,1}$ ``from above'' at $x_0$. From the geometric point
of view, the subgraph of $u$ has an exterior tangent paraboloid at the point $(x_0,u(x_0))$.

This regularity of $u$ ``from above'' at a point $x_0$, coupled with the property of being a viscosity subsolution,
is enough to guarantee that $\h u(x_0)$ is well defined.\\
As a consequence, a viscosity subsolution is a classical subsolution in every ``viscosity point'', i.e. in every point which can
be touched from above by a $C^{1,1}$ function. More precisely:

\begin{prop}\label{CH:4:pwise_prop}
Let
\[\h u\le f\qquad\textrm{in }\Omega,\]
and let $x_0\in\Omega$. Suppose that there exists a function $ v\in C^{1,1}(B_r(x_0))$ that touches $u$ from above at $x_0$, that is
\[ v(x_0)=u(x_0)\qquad\textrm{and}\qquad v(y)\geq u(y)\quad\forall\,y\in B_r(x_0).\]
Then $\h u(x_0)$ is defined in the classical sense and
\[\h u(x_0)\leq f(x_0).\]
\end{prop}
\begin{proof}
We begin by showing that
$\delta_g(u,x_0;\xi)|\xi|^{-n-s}$ is integrable in $\R^n$, so that $\h u(x_0)$ is well defined in the classical sense.

For the argument we follow \cite[Proposition 1]{Lin16}.
We consider the functions
\[ v_\varrho(y):=\left\{\begin{array}{cc} v(y) & \textrm{if }y\in B_\varrho(x_0),\\
u(y) & \textrm{if }y\in\R^n\setminus B_\varrho(x_0),\end{array}\right.\]
for every $\varrho\in(0,r]$ and we denote
\[\delta_g^+( v_\varrho,x_0;\xi):=\max\{\delta_g( v_\varrho,x_0;\xi),0\}\qquad\textrm{and}\qquad
\delta_g^-( v_\varrho,x_0;\xi):=\max\{-\delta_g( v_\varrho,x_0;\xi),0\}.\]
We remark that, since $ v\in C^{1,1}(B_r(x_0))$, the function $\delta_g( v_\varrho,x_0;\xi)|\xi|^{-n-s}$ is integrable
in $\R^n$, that is
\[\int_{\R^n}\frac{\delta_g^+( v_\varrho,x_0;\xi)+\delta_g^-( v_\varrho,x_0;\xi)}{|\xi|^{n+s}}\,d\xi=
\int_{\R^n}\frac{|\delta_g( v_\varrho,x_0;\xi)|}{|\xi|^{n+s}}\,d\xi<+\infty.\]
Moreover, notice that
\[\delta_g(u,x_0;\xi)\geq\delta_g( v_{\varrho_1},x_0;\xi)\geq\delta_g( v_{\varrho_2},x_0;\xi),\qquad\textrm{for every }
0<\varrho_1\leq\varrho_2\leq r.\]
Therefore, in particular
\eqlab{\label{CH:4:delta-}
\int_{\R^n}\frac{\delta_g^-(u,x_0;\xi)}{|\xi|^{n+s}}\,d\xi
\le\int_{\R^n}\frac{|\delta_g(v_r,x_0;\xi)|}{|\xi|^{n+s}}\,d\xi<+\infty.}

Now, since $u$ is a subsolution, we have
\[\int_{\R^n}\frac{\delta_g( v_\varrho,x_0;\xi)}{|\xi|^{n+s}}\,d\xi\leq f(x_0),\]
that is
\[\int_{\R^n}\frac{\delta_g^+( v_\varrho,x_0;\xi)}{|\xi|^{n+s}}\,d\xi
\leq\int_{\R^n}\frac{\delta_g^-( v_\varrho,x_0;\xi)}{|\xi|^{n+s}}\,d\xi +f(x_0).\]
Since
\[\delta_g^+( v_\varrho,x_0;\xi)\nearrow\delta_g^+(u,x_0;\xi)\qquad\textrm{as }\varrho\searrow0,\]
the monotone convergence Theorem gives
\[\lim_{\varrho\to0^+}\int_{\R^n}\frac{\delta_g^+( v_\varrho,x_0;\xi)}{|\xi|^{n+s}}\,d\xi
=\int_{\R^n}\frac{\delta_g^+(u,x_0;\xi)}{|\xi|^{n+s}}\,d\xi.\]
Moreover
\bgs{
\int_{\R^n}\frac{\delta_g^+( v_{\varrho_1},x_0;\xi)}{|\xi|^{n+s}}\,d\xi&
\leq\int_{\R^n}\frac{\delta_g^-( v_{\varrho_1},x_0;\xi)}{|\xi|^{n+s}}\,d\xi+f(x_0)\\
&
\leq\int_{\R^n}\frac{\delta_g^-( v_{\varrho_2},x_0;\xi)}{|\xi|^{n+s}}\,d\xi+f(x_0)<+\infty,}
for every $0<\varrho_1\leq \varrho_2\leq r$. Thus
\begin{equation}\label{CH:4:pwise_def}
\int_{\R^n}\frac{\delta_g^+(u,x_0;\xi)}{|\xi|^{n+s}}\,d\xi\leq
\int_{\R^n}\frac{\delta_g^-( v_\varrho,x_0;\xi)}{|\xi|^{n+s}}\,d\xi+f(x_0)<+\infty,
\end{equation}
for every $\varrho\in(0,r]$.
By \eqref{CH:4:delta-} and \eqref{CH:4:pwise_def}, we see that
$\delta_g(u,x_0;\xi)|\xi|^{-n-s}$ is integrable in $\R^n$ and hence $\h u(x_0)$ is well defined.

Finally, since for every $\varrho\in(0,r]$ we have
\[\frac{\delta_g^-( v_\varrho,x_0;\xi)}{|\xi|^{n+s}}\leq\frac{\delta_g^-( v_r,x_0;\xi)}{|\xi|^{n+s}},\]
which is integrable in $\R^n$, by Lebesgue's dominated convergence Theorem we can pass to the limit $\varrho\to0$ in the right hand side of \eqref{CH:4:pwise_def}, obtaining
\[\int_{\R^n}\frac{\delta_g^+(u,x_0;\xi)}{|\xi|^{n+s}}\,d\xi\leq
\int_{\R^n}\frac{\delta_g^-(u,x_0;\xi)}{|\xi|^{n+s}}\,d\xi+f(x_0),\]
that is
\[\h u(x_0)\le f(x_0),\]
as claimed.
\end{proof}

For later use, it is convenient to introduce the following definition.
\begin{defn}\label{CH:4:defC11}
Let $u:\R^n\to\R$ and let $x_0\in\R^n$. The function $u$ is $C^{1,1}$ at $x_0$, and we write
$u\in C^{1,1}(x_0)$, if there exist $\ell\in\R^n$ and $M,\,r>0$ such that
\eqlab{\label{CH:4:localC11}
|u(x_0+\xi)-u(x_0)-\ell\cdot\xi|\le M|\xi|^2,\qquad\forall\,\xi\in B_r.
}
\end{defn}

We remark that we clearly have
\bgs{
u\in C^{1,1}(B_R(x_0))\quad\implies\quad u\in C^{1,1}(x_0).
}
Roughly speaking, being $C^{1,1}$ at $x_0$ means that there exist both an interior and an exterior tangent paraboloid
to the subgraph of $u$ at the point $(x_0,u(x_0))$.

As a consequence of Proposition \ref{CH:4:pwise_prop}, we obtain the following Corollary:
\begin{corollary}\label{CH:4:pwise_cor}
Let
\[\h u\le f\qquad\textrm{in }\Omega,\]
and let $x_0\in\Omega$. If $u\in C^{1,1}(x_0)$, then $\h u(x_0)$ is well defined and
\bgs{
\h u(x_0)\le f(x_0).
}
\end{corollary}

\begin{proof}
Consider the paraboloid
\bgs{
q(x):=u(x_0)+\ell\cdot(x-x_0)+M|x-x_0|^2,\qquad\forall\,x\in B_r(x_0),
}
with $\ell,\,M$ and $r$ as in Definition \ref{CH:4:defC11}. Then $q\in C^{1,1}(B_r(x_0))$ and
by \eqref{CH:4:localC11} we know that $q$ touches $u$ from above
at $x_0$.
Thus the conclusion follows from Proposition \ref{CH:4:pwise_prop}.
\end{proof}

\subsection{Sup-convolutions}
In this subsection we introduce and study the sup-convolutions $u^\eps$ of a viscosity subsolution $u$.
These provide a sequence of subsolutions which converge to $u$ and which enjoy nice regularity properties,
since they are semiconvex functions.\\

We will consider only globally bounded subsolutions.

%

\begin{defn}\label{CH:4:supconvolution}
Let $u:\R^n\to\R$ be a bounded function. We define the sup-convolution $u^\eps$ of $u$ as
\bgs{
u^\eps(x):=\sup_{y\in\R^n}\left\{u(y)-\frac{1}{\eps}|y-x|^2\right\}\qquad\forall\,x\in\R^n,
}
for every $\eps>0$.
\end{defn}

Now we point out some easy properties of sup-convolutions.\\
We begin by remarking that, by definition,
\bgs{
u^\eps\ge u\quad\mbox{in }\R^n.
}
Moreover, notice that if we denote
\bgs{
\sup_{\R^n}|u|=:M<+\infty,
}
then
\eqlab{\label{CH:4:silence}
u^\eps(x)=\sup_{|y-x|\le\sqrt{2M\eps}}\left\{u(y)-\frac{1}{\eps}|y-x|^2\right\}\qquad\forall\,x\in\R^n.
}
Indeed, if $|y-x|>\sqrt{2M\eps}$, then
\bgs{
u(y)-\frac{1}{\eps}|y-x|^2<-M\le u(x),
}
but we know that
\bgs{
u^\eps(x)\ge u(x).
}

\begin{remark}\label{CH:4:attained_conv}
Given an open set $\Omega\subseteq\R^n$, we denote
\eqlab{\label{CH:4:open_visc}
\Omega^\eps:=\left\{x\in\Omega\,|\,d(x,\partial\Omega)>2\sqrt{2M\eps}\right\}.
}
If $u$ is upper semicontinuous in an open set $\Omega$, then for every $x\in\Omega^\eps$ there exists $y_0\in\Omega$ such that
\bgs{
u^\eps(x)=u(y_0)-\frac{1}{\eps}|y_0-x|^2=\max_{|y-x|\le\sqrt{2M\eps}}\left\{u(y)-\frac{1}{\eps}|y-x|^2\right\}.
}
This follows straightforwardly from \eqref{CH:4:silence} and the upper semicontinuity of $u$.
\end{remark}

In the following Theorem we collect some important properties of sup-convolutions
which can be found in \cite[Chapter 1]{Ambrosio}.

We first recall the definition of semiconvex functions.

\begin{defn}
Let $\Omega\subseteq\R^n$ be an open set and let $u:\Omega\to\R$. We say that $u$ is semiconvex in $\Omega$ if there
exists a constant $c\ge0$ such that
\bgs{
x\mapsto u(x)+\frac{c}{2}|x|^2
}
is convex in any ball $B\subseteq\Omega$. The smallest constant $c\ge0$ for which this happens is called the
semiconvexity constant of $u$ and denoted $sc(u,\Omega)$.
\end{defn}

\begin{theorem}\label{CH:4:supconvo_teo}
Let $u:\R^n\to\R$ be a bounded function. Then
$u^\eps$ is semiconvex in $\R^n$, with semiconvexity constant
\bgs{
sc(u^\eps,\R^n)\le\frac{2}{\eps}.
}
Therefore $u^\eps\in W^{1,\infty}_{\loc}(\R^n)$ and $\nabla u^\eps\in BV_{\loc}(\R^n,\R^n)$. Moreover
$u^\eps\in C^{1,1}(x)$ for almost every $x\in\R^n$.

If $u$ is upper semicontinuous in an open set $\Omega\subseteq\R^n$, then
\bgs{
u^\eps(x)\searrow u(x)\quad\mbox{as }\eps\searrow0,\qquad\forall\,x\in\Omega.
}
The convergence is locally uniform if $u$ is continuous in $\Omega$.
\end{theorem}

\begin{proof}
The semiconvexity of $u^\eps$ follows by \cite[Proposition 4, $(i)$]{Ambrosio}. Then, by \cite[Theorem 15]{Ambrosio} this implies that
$u^\eps\in W^{1,\infty}_{\loc}(\R^n)$ and by \cite[Theorem 16]{Ambrosio} that $\nabla u^\eps\in BV_{\loc}(\R^n,\R^n)$.
That $u^\eps\in C^{1,1}(x)$ for almost every $x\in\R^n$ follows from the Taylor expansion in point \cite[Theorem 16, $(ii)$]{Ambrosio}.
Finally, the convergence of $u^\eps$ to $u$ is obtained by arguing as in the proof of \cite[Proposition 4, $(ii)$]{Ambrosio}.
\end{proof}

One of the most important features of sup-convolutions consists in preserving the viscosity subsolution property
(eventually up to a small error).
More precisely:
\begin{prop}\label{CH:4:convo_visc}
Let $\Omega\subseteq\R^n$ be a bounded open set, let $f\in C(\overline{\Omega})$ and let $u:\R^n\to\R$ be a bounded function,
\bgs{
M:=\sup_{\R^n}|u|<+\infty,
}
such that
\bgs{
\h u\le f\quad\mbox{in }\Omega.
}
Then
\bgs{
\h u^\eps\le f+c_\eps\quad\mbox{in }\Omega^\eps,
}
where $\Omega^\eps$ is defined in \eqref{CH:4:open_visc} and the constant $c_\eps\ge0$ depends
only on $\eps,\,M$ and the modulus of continuity of $f$. More precisely,
\bgs{
c_\eps:=\sup_{\substack{x,\,y\in\overline\Omega \\ |x-y|\le\sqrt{2M\eps}}}|f(x)-f(y)|.
}
In particular
\eqlab{\label{CH:4:errorproperty}
c_\eps\searrow0\quad\mbox{as }\eps\searrow0\qquad\mbox{and}\qquad c_\eps=0\quad{if}\quad f\mbox{ is constant.}
}
\end{prop}

\begin{proof}
Let $x_0\in\Omega^\eps$ and suppose that there exists $v\in C^{1,1}(B_r(x_0))$ such that
\bgs{
v(x_0)=u^\eps(x_0)\quad\mbox{and}\quad v(x)\ge u^\eps(x),\quad\forall\,x\in B_r(x_0).
}
We need to show that
\bgs{
\h\tilde v(x_0)\le f(x_0)+c_\eps.
}
By Remark \ref{CH:4:attained_conv} we know that we can find $y_0\in\Omega$ such that $|y_0-x_0|\le\sqrt{2M\eps}$ and
\bgs{
u^\eps(x_0)=u(y_0)-\frac{1}{\eps}|y_0-x_0|^2,
}
Then we define
\bgs{
\psi(x):=v(x+x_0-y_0)+\frac{1}{\eps}|y_0-x_0|^2,\qquad\forall\,x\in B_r(y_0),
}
and we remark that clearly $\psi\in C^{1,1}(B_r(y_0))$. Moreover
\bgs{
\psi(y_0)=v(x_0)+\frac{1}{\eps}|y_0-x_0|^2
=u^\eps(x_0)+\frac{1}{\eps}|y_0-x_0|^2=u(y_0).
}
Then notice that, since $v\ge u^\eps$ in $B_r(x_0)$, by definition of $u^\eps$ we obtain
\bgs{
u(y)-\frac{1}{\eps}|y-x|^2\le u^\eps(x)\le v(x),\qquad\forall\,y\in\R^n\mbox{ and }x\in B_r(x_0).
}
Taking $y\in B_r(y_0)$ and $x:=y+x_0-y_0$ gives
\bgs{
u(y)\le \psi(y),\qquad\forall\,y\in B_r(y_0).
}
Thus $\psi$ touches $u$ from above at $y_0$ and hence
\bgs{
\h\tilde\psi(y_0)\le f(y_0).
}
Now notice that by changing variables we find
\bgs{
\h^{<r}\tilde\psi(y_0)&
=2\PV\int_{B_r(y_0)}G\Big(\frac{\psi(y_0)-\psi(y)}{|y-y_0|}\Big)\frac{dy}{|y-y_0|^{n+s}}\\
&
=2\PV\int_{B_r(x_0)}G\Big(\frac{v(x_0)-v(x)}{|x-x_0|}\Big)\frac{dx}{|x-x_0|^{n+s}}\\
&
=\h^{<r}\tilde v(x_0).
}
On the other hand,
\bgs{
\h^{\ge r}\tilde\psi(y_0)&
=2\int_{\R^n\setminus B_r(y_0)}G\Big(\frac{u^\eps(x_0)+\frac{1}{\eps}|y_0-x_0|^2-u(y)}
{|y-y_0|}\Big)\frac{dy}{|y-y_0|^{n+s}}\\
&
=2\int_{\R^n\setminus B_r(x_0)}G\Big(\frac{u^\eps(x_0)+\frac{1}{\eps}|y_0-x_0|^2-u(x+y_0-x_0)}
{|x-x_0|}\Big)\frac{dx}{|x-x_0|^{n+s}}.
}
We remark that by taking $y:=x+y_0-x_0$ in the definition of $u^\eps(x)$ as a sup, we get
\bgs{
\frac{1}{\eps}|y_0-x_0|^2-u(x+y_0-x_0)\ge-u^\eps(x),\qquad\forall\,x\in\R^n\setminus B_r(x_0).
}
Hence
\bgs{
\h^{\ge r}\tilde\psi(y_0)\ge2\int_{\R^n\setminus B_r(x_0)}G\Big(\frac{u^\eps(x_0)-u^\eps(x)}
{|x-x_0|}\Big)\frac{dx}{|x-x_0|^{n+s}}=\h^{\ge r}\tilde v(x_0).
}
This implies that
\bgs{
\h\tilde v(x_0)\le\h\tilde\psi(y_0)\le f(y_0)\le f(x_0)+c_\eps,
}
as claimed. To conclude the proof, notice that \eqref{CH:4:errorproperty} follows from the
definition of $c_\eps$ and the uniform continuity of $f$.
\end{proof}

As a consequence, exploiting the regularity of $u^\eps$ we find that $u^\eps$ is
a classical subsolution almost everywhere in $\Omega^\eps$.

\begin{corollary}\label{CH:4:convo_pwise}
Let $\Omega\subseteq\R^n$ be a bounded open set, let $f\in C(\overline{\Omega})$ and let $u:\R^n\to\R$ be a bounded function such that
\bgs{
\h u\le f\quad\mbox{in }\Omega.
}
Then for almost every $x\in\Omega^\eps$ we have that
$\h u^\eps(x)$ is well defined and
\bgs{
\h u^\eps(x)\le f(x)+c_\eps.
}
\end{corollary}

\begin{proof}
By Theorem \ref{CH:4:supconvo_teo} we know that $u^\eps\in C^{1,1}(x)$ for almost every $x\in\R^n$.
Then the conclusion follows from Proposition \ref{CH:4:convo_visc} and Corollary \ref{CH:4:pwise_cor}.
\end{proof}

\subsection{Weak (sub)solutions}

Given a function $u:\R^n\to\R$, we define
\bgs{
\langle\h u,v\rangle:=\int_{\R^n}\int_{\R^n}G\Big(\frac{u(x)-u(y)}{|x-y|}\Big)\big(v(x)-v(y)\big)\frac{dx\,dy}{|x-y|^{n+s}},
}
for every $v\in W^{s,1}(\R^n)$. In particular, this is well defined for every $v\in C_c^\infty(\Omega)$,
where we understand that $v$ is extended by zero outside $\Omega$.

\begin{defn}
Let $\Omega\subseteq\R^n$ be a bounded open set and let $f\in C(\overline\Omega)$.
We say that a function $u:\R^n\to\R$
is a weak subsolution in $\Omega$ if
\bgs{
\langle\h u,v\rangle\le\int_\Omega fv\,dx,\qquad\forall\,v\in C_c^\infty(\Omega)\mbox{ s.t. }v\ge0.
}
\end{defn}

We want to pass from a function $u$ which is a subsolution almost everywhere
to a weak subsolution.
In order to do this, it is enough to ask $u$ to have a gradient with bounded variation.
More precisely, we introduce the space
\bgs{
BH(\Omega)&:=\left\{u\in W^{1,1}(\Omega)\,|\,\nabla u\in BV(\Omega,\R^n)\right\}\\
&
=\left\{u\in W^{1,1}(\Omega)\,|\,\partial_j u\in BV(\Omega),\,\forall\,j=1,\dots,n\right\},
}
endowed with the norm
\bgs{
\|u\|_{BH(\Omega)}:=\|u\|_{W^{1,1}(\Omega)}+|D^2u|(\Omega).
}

For the properties of the space $BH(\Omega)$ of functions of bounded Hessian, we refer the interested reader to
\cite{Demengel}.
We remark that sometimes the notation $BV^2(\Omega)=BH(\Omega)$ is also used in the literature.

We only recall the following ``density'' property, \cite[Proposition 1.4]{Demengel}:
\begin{prop}\label{CH:4:bh_approx}
Let $\Omega\subseteq\R^n$ be a bounded open set with $C^2$ boundary and let $u\in BH(\Omega)$.
Then there exist $u_k\in C^2(\Omega)\cap W^{2,1}(\Omega)$ such that
\bgs{
\lim_{k\to\infty}\left\{\|u-u_k\|_{W^{1,1}(\Omega)}+\big||D^2u|(\Omega)-|D^2u_k|(\Omega)\big|\right\}=0.
}
\end{prop}

Exploiting this density property, we can prove the following:

\begin{lemma}\label{CH:4:bv2_prop}
Let $\Omega\subseteq\R^n$ be a bounded open set and let $u\in BH(\Omega)$.
Let $\Omega'\Subset\Omega$ and let $d:=dist(\Omega',\partial\Omega)/2$. Then
\eqlab{\label{CH:4:bv2esti}
\int_{\Omega'}|u(x+\xi)+u(x-\xi)-2u(x)|\,dx\le 2|\xi|^2|D^2u|(\Omega),\qquad\forall\,\xi\in B_d.
}
\end{lemma}

\begin{proof}
Let $\Op\subseteq\Omega$ be a bounded open set with $C^2$ boundary such that
\eqlab{\label{CH:4:op_appro_BH}
\Omega'\Subset\Op,\qquad\mbox{with}\quad dist(\Omega',\partial\Op)>d.
}
By Proposition \ref{CH:4:bh_approx} we can find $u_k\in C^2(\Op)\cap W^{2,1}(\Op)$ such that
\eqlab{\label{CH:4:appro_BH_conv}
\lim_{k\to\infty}\left\{\|u-u_k\|_{W^{1,1}(\Op)}+\big||D^2u|(\Op)-|D^2u_k|(\Op)\big|\right\}=0.
}
Now notice that
\bgs{
|u_k(x+\xi)+u_k(x-\xi)-2u_k(x)|\le|u_k(x+\xi)&-u_k(x)-\nabla u_k(x)\cdot\xi|\\
&+
|u_k(x-\xi)-u_k(x)-\nabla u_k(x)\cdot(-\xi)|.
}
Then by Taylor's formula with integral remainder we have
\bgs{
|u_k(x+\xi)-u_k(x)-\nabla u_k(x)\cdot\xi|\le|\xi|^2\int_0^1|D^2u_k(x+t\xi)|\,dt,
}
and similarly for $-\xi$. Integrating in $x$ over $\Omega'$ and switching the order of integration
gives
\bgs{
\int_{\Omega'}|u_k(x+\xi)&+u_k(x-\xi)-2u_k(x)|\,dx\le|\xi|^2\int_{\Omega'}\Big(\int_{-1}^1|D^2u_k(x+t\xi)|\,dt\Big)dx\\
&
=|\xi|^2\int_{-1}^1\Big(\int_{\Omega'}|D^2u_k(x+t\xi)|\,dx\Big)dt\le2|\xi|^2|D^2u_k|(\Op),
}
since $|\xi|<d$ and $\Op$ satisfies \eqref{CH:4:op_appro_BH}.

Then by Fatou's Lemma and \eqref{CH:4:appro_BH_conv} we obtain
\bgs{
\int_{\Omega'}|u(x+\xi)+u(x-\xi)-2u(x)|\,dx&\le2|\xi|^2\liminf_{k\to\infty}|D^2u_k|(\Op)=2|\xi|^2|D^2u|(\Op)\\
&
\le2|\xi|^2|D^2u|(\Omega),
}
proving \eqref{CH:4:bv2esti} and concluding the proof of the Lemma.
\end{proof}

\begin{prop}\label{CH:4:BHprop_curvature}
Let $\Omega\subseteq\R^n$ be a bounded open set and let $u\in BH(\Omega)$.
Then
\bgs{
\h u\in L^1_{\loc}(\Omega),
}
and
\bgs{
\langle\h u,v\rangle=\int_{\Omega}\h u(x)v(x)\,dx,\qquad\forall\,v\in C^\infty_c(\Omega).
}
\end{prop}

\begin{proof}
Let $\Omega'\Subset\Omega$ and let $d:=dist(\Omega',\partial\Omega)/2$. We recall that
\bgs{
|\delta_g(u,x;\xi)|\le\frac{|u(x+\xi)+u(x-\xi)-2u(x)|}{|\xi|}.
}
Therefore, by Remark \ref{CH:4:rmk_tail} and \eqref{CH:4:bv2esti} we obtain
\eqlab{\label{CH:4:fundaBHesti}
\int_{\Omega'}&|\h u(x)|\,dx\le\int_{\Omega'}dx\int_{\R^n}\frac{|\delta_g(u,x;\xi)|}{|\xi|^{n+s}}\,d\xi\\
&
\le |\Omega'|\frac{n\omega_n}{s}\Lambda\,d^{-s}
+\int_{B_d}\Big(\int_{\Omega'}|u(x+\xi)+u(x-\xi)-2u(x)|\,dx\Big)\frac{d\xi}{|\xi|^{n+1+s}}\\
&
\le |\Omega'|\frac{n\omega_n}{s}\Lambda\,d^{-s}
+2|D^2u|(\Omega)\frac{n\omega_n}{1-s}d^{1-s}<+\infty.
}
This proves that $\h u\in L^1_{\loc}(\Omega)$.

As a consequence of \eqref{CH:4:fundaBHesti}, since
\bgs{
|\h^{\ge\varrho}u(x)|\le\int_{\R^n}\frac{|\delta_g(u,x;\xi)|}{|\xi|^{n+s}}\,d\xi,\qquad\forall\,\varrho>0,
}
given $v\in C^\infty_c(\Omega)$ we can apply Lebesgue's dominated convergence Theorem to obtain
\bgs{
\lim_{\varrho\to0^+}\int_{\R^n} \h^{\ge\varrho}u(x)v(x)\,dx=\int_{\R^n} \h u(x)v(x)\,dx.
}
Now notice that by symmetry
\bgs{
\int_{\R^n} \h^{\ge\varrho}u(x)v(x)\,dx=\int_{\R^n}\int_{\R^n}G\Big(\frac{u(x)-u(y)}{|x-y|}\Big)\big(v(x)-v(y)\big)
\big(1-\chi_{B_\varrho}(x-y)\big)\frac{dx\,dy}{|x-y|^{n+s}}.
}
Finally, since $v\in C^\infty_c(\Omega)\subseteq W^{s,1}(\R^n)$, we can apply again Lebesgue's dominated convergence Theorem to obtain
\bgs{
\lim_{\varrho\to0^+}\int_{\R^n}\int_{\R^n}G\Big(\frac{u(x)-u(y)}{|x-y|}\Big)\big(v(x)-v(y)\big)
\big(1-\chi_{B_\varrho}(x-y)\big)\frac{dx\,dy}{|x-y|^{n+s}}=\langle\h u,v\rangle,
}
concluding the proof.
\end{proof}

Then exploiting Proposition \ref{CH:4:BHprop_curvature}, Theorem \ref{CH:4:supconvo_teo} and Corollary \ref{CH:4:convo_pwise} we immediately obtain the following:

\begin{corollary}\label{CH:4:convo_weak}
Let $\Omega\subseteq\R^n$ be a bounded open set, let $f\in C(\overline\Omega)$ and let $u:\R^n\to\R$ be a bounded function such that
\bgs{
\h u\le f\quad\mbox{in }\Omega.
}
Then
\eqlab{\label{CH:4:supconvoweak}
\langle \h u^\eps,v\rangle\le\int_{\Omega}\big(f(x)+c_\eps\big)v(x)\,dx,\qquad\forall\,v\in C^\infty_c(\Omega^\eps)\mbox{ s.t. }v\ge0.
}
\end{corollary}

\subsection{Viscosity implies weak}

As a consequence of Lemma \ref{CH:4:easylemma}, exploiting the fact that supconvolutions are weak subsolutions we obtain the following:

\begin{theorem}\label{CH:4:ViscWeak_reg}
Let $\Omega\subseteq\R^n$ be a bounded open set, let $f\in C(\overline\Omega)$ and
let $u:\R^n\to\R$ be a viscosity subsolution,
\bgs{
\h u\le f\quad\mbox{in }\Omega.
}
Suppose that $u$ is bounded and
assume also that there exists a closed set $S\subseteq\R^n\setminus\Omega$
such that $|S|=0$ and $u$ is upper semicontinuous in $\R^n\setminus S$.
Then $u$ is a weak subsolution in $\Omega$,
\bgs{
\langle\h u,v\rangle\le\int_\Omega fv\,dx,\qquad\forall\,v\in C^\infty_c(\Omega)\mbox{ s.t. }v\ge0.
}
\end{theorem}

\begin{proof}
The hypothesis on $u$ and Theorem \ref{CH:4:supconvo_teo} imply that
\bgs{
u^\eps(x)\to u(x),\qquad\forall\,x\in\R^n\setminus S,
}
and hence almost everywhere in $\R^n$. Let $v\in C^\infty_c(\Omega)$ be such that $v\ge0$. Notice that
\bgs{
supp\,v\subseteq\Omega^\eps,
}
for every $\eps$ small enough. Thus, by \eqref{CH:4:supconvoweak} and recalling \eqref{CH:4:errorproperty},
we obtain
\bgs{
\langle\h u,v\rangle=\lim_{\eps\to0^+}\langle\h u^\eps,v\rangle\le
\lim_{\eps\to0^+}\int_\Omega(f+c_\eps)v\,dx=\int_\Omega fv\,dx,
}
concluding the proof.
\end{proof}

In particular, if $|\partial\Omega|=0$, we allow for $\partial\Omega\subseteq S$, so we are not asking $u$
to be continuous across $\partial\Omega$.

\medskip

We are now going to use an approximation procedure to extend
Theorem \ref{CH:4:ViscWeak_reg} to the case of arbitrary exterior data.

The crucial point consists in the following observation, that follows essentially from the fact that
$\h^{\ge d}u(x)$ can be bounded in terms of only $d$, independently both of $u$ or $x$ (see Remark \ref{CH:4:rmk_tail}).

\begin{theorem}\label{CH:4:approtrick}
Let $\Omega\subseteq\R^n$ be a bounded open set and let $f\in C(\overline\Omega)$. Let $u:\R^n\to\R$ be locally integrable
in $\R^n$ and suppose that
\bgs{
\h u\le f\quad\mbox{in }\Omega.
}
Let $u_k:\R^n\to\R$ be such that $u_k\to u$ in $L^1_{\loc}(\R^n)$. Given two open sets
\bgs{
\Omega'\Subset\Op\subseteq\Omega,
}
we define
\begin{equation*}
\bar u_k(x):=\left\{\begin{array}{cc}
u(x) & \mbox{if }x\in\Op,\\
u_k(x) & \mbox{if }x\in\R^n\setminus\Op.
\end{array}\right.
\end{equation*}
Then for every $k\in\N$ there exists
a constant $e_k\ge0$ such that $e_k\to0$ and
\bgs{
\h \bar u_k\le f+e_k\quad\mbox{in }\Omega'.
}
\end{theorem}

\begin{proof}
We denote $d:=dist(\Omega',\partial\Op)>0$ and we remark that for every $x\in\Omega'$ we have
\eqlab{\label{CH:4:approdel2}
\delta_g(\bar u_k,x;\xi)=\delta_g(u,x;\xi),\qquad\forall\,\xi\in B_d.
}
On the other hand, let
\bgs{
\omega_k(x):=\h^{\ge d}\bar u_k(x)-\h^{\ge d}u(x),\qquad\forall\,x\in\Omega',
}
and let $R_0>0$ be such that $\Omega\subseteq B_{R_0}$. Then for every $x\in\Omega'$ and $R\ge R_0+d$ we have
\bgs{
|\omega_k(x)|&\le2\int_{\R^n\setminus B_d(x)}\Big|G\Big(\frac{u(x)-\bar u_k(y)}{|x-y|}\Big)-G\Big(\frac{u(x)-u(y)}{|x-y|}\Big)\Big|
\frac{dy}{|x-y|^{n+s}}\\
&
\le2\int_{B_R(x)\setminus B_d(x)}\frac{|\bar u_k(y)-u(y)|}{|x-y|^{n+1+s}}\,dy
+2\Lambda\int_{\R^n\setminus B_R(x)}\frac{dy}{|x-y|^{n+s}}\\
&
\le\frac{2}{d^{n+1+s}}\|u_k-u\|_{L^1(B_{R+R_0})}+\frac{2\Lambda n\omega_n}{s}R^{-s}.
}
Hence for every $k\in\N$ we obtain
\bgs{
\sup_{x\in\Omega'}|\omega_k(x)|\le\frac{2}{d^{n+1+s}}\|u_k-u\|_{L^1(B_{R+R_0})}+\frac{2\Lambda n\omega_n}{s}R^{-s},
\qquad\forall\,R\ge R_0+d.
}
Thus, if we define
\bgs{
e_k:=\inf_{R\ge R_0+d}\left(\frac{2}{d^{n+1+s}}\|u_k-u\|_{L^1(B_{R+R_0})}+\frac{2\Lambda n\omega_n}{s}R^{-s}\right),
}
we get
\eqlab{\label{CH:4:unif_est}
\sup_{x\in\Omega'}|\omega_k(x)|\le e_k,\qquad\forall\,k\in\N.
}
Now notice that, since $u_k\to u$ in $L^1_{\loc}(\R^n)$, we have
\bgs{
\limsup_{k\to\infty}e_k\le\limsup_{k\to\infty}\left(\frac{2}{d^{n+1+s}}\|u_k-u\|_{L^1(B_{R+R_0})}+\frac{2\Lambda n\omega_n}{s}R^{-s}\right)
=\frac{2\Lambda n\omega_n}{s}R^{-s},
}
for every $R\ge R_0+d$. Letting $R\nearrow+\infty$ proves that $e_k\to0$.

\smallskip

Now let $x_0\in\Omega'$ be such that there exists $v\in C^{1,1}(B_r(x_0))$ with $r<dist(x_0,\Omega')$ and
\bgs{
v(x_0)=\bar u_k(x_0)=u(x_0)\quad\mbox{and}\quad v(x)\ge \bar u_k(x)=u(x)\quad\forall\,x\in B_r(x_0).
}
By Proposition \ref{CH:4:pwise_prop} we obtain
\bgs{
\h u(x_0)\le f(x_0).
}
Hence, by \eqref{CH:4:approdel2} and \eqref{CH:4:unif_est} we get
\bgs{
\h \bar u_k(x_0)&=\h^{<d}u(x_0)+\h^{\ge d}\bar u_k(x_0)=\h u(x_0)+\omega_k(x_0)\\
&
\le f(x_0)+|\omega_k(x_0)|\le f(x_0)+e_k.
}
Finally, notice that if we set
\[\tilde v_k(x):=\left\{\begin{array}{cc} v(x) & \textrm{if }x\in B_r(x_0),\\
\bar u_k(x) & \textrm{if }x\in\R^n\setminus B_r(x_0),\end{array}\right.\]
then by Remark \ref{CH:4:mah} we obtain
\bgs{
\h\tilde v_k(x_0)\le\h \bar u_k(x_0)\le f(x_0)+e_k,
}
concluding the proof.
\end{proof}

With this fundamental approximation tool at hand, we are ready to prove the general ``viscosity implies weak'' Theorem.

\begin{proof}[Proof of Theorem \ref{CH:4:Gen_viscweak}]
Let $v\in C^\infty_c(\Omega)$ such that $v\ge0$. Then we can find two open sets such that
\bgs{
supp\,v\Subset\Omega'\Subset\Op\Subset\Omega,
}
and such that $|\partial\Op|=0$.

Since $u$ is locally integrable in $\R^n$, we can find a sequence of functions $u_k\in C(\R^n)\cap L^\infty(\R^n)$
such that
\bgs{
u_k\to u\quad\mbox{in }L^1_{\loc}(\R^n)\quad\mbox{and}\quad\mbox{a.e. in }\R^n.
}

Now let $\bar u_k$ and $e_k$ be as defined in Theorem \ref{CH:4:approtrick}. Notice that since $u$ is locally bounded in $\Omega$,
it is bounded in $\Op$, and hence the functions $\bar u_k$ are bounded in $\R^n$.

Moreover, since $u$ is upper semicontinuous in $\Omega$ and $u_k$ is continuous in $\R^n$,
the functions $\bar u_k$ are upper semicontinuous in $\R^n\setminus\partial\Op$.

We can thus apply Theorem \ref{CH:4:approtrick} and Theorem \ref{CH:4:ViscWeak_reg} to obtain
\bgs{
\langle \h\bar u_k,v\rangle\le\int_{\Omega'} (f+e_k)v\,dx,\qquad\forall\,k\in\N.
}
Then, since $\bar u_k\to u$ almost everywhere in $\R^n$ and $e_k\to 0$,
by Lemma \ref{CH:4:easylemma} we get
\bgs{
\langle \h u,v\rangle=\lim_{k\to\infty}\langle \h\bar u_k,v\rangle\le\lim_{k\to\infty}\int_{\Omega'} (f+e_k)v\,dx
=\int_\Omega fv\,dx.
}
This concludes the proof of the Theorem.
\end{proof}

\section{Minimizers of~$\F_s^M$ versus minimizers of~$\Per_s$}\label{CH:4:REARR_SECTI}

Here, we bring forward our analysis of the geometric properties enjoyed the functional~$\F_s^M$, and in particular of its relation with the~$s$-perimeter.

%
We will show that the~$s$-perimeter of a set~$E$, which is a subgraph outside $\Omega^\infty$ and whose boundary is trapped inside a strip of finite height inside $\Omega^\infty$, decreases under a vertical rearrangement that transforms~$E$ into a global subgraph. This fact will be a consequence of a new rearrangement inequality for rather general~$1$-dimensional integral set functions, that we establish in the following subsection.

\subsection{A one-dimensional rearrangement inequality}

Let~$K: \R \to \R$ be a non-negative function. Given two sets~$A, B \subseteq \R$, we define
\begin{equation} \label{CH:4:IKdef}
\I_K(A, B) := \int_A \int_B d\mu, \quad \mbox{where } \, d\mu = d\mu_K(x - y) := K(x - y) \, dx\,dy,
\end{equation}
whenever this quantity is finite.

Fix two real numbers~$\alpha, \beta$ and consider two sets~$A, B$ satisfying
\begin{equation} \label{CH:4:ABconstr}
(-\infty, \alpha) \subseteq A \quad \mbox{and} \quad (\beta, +\infty) \subseteq B.
\end{equation}
We define the~\emph{decreasing rearrangement}~$A_*$ of~$A$ as
\begin{equation} \label{CH:4:decrrearr}
A_* := (-\infty, a_* ), \quad \mbox{with } \, a_* := \lim_{R \rightarrow +\infty} \left( \int_{-R}^{R} \chi_A(t) \, dt - R \right).
\end{equation}
Similarly, we define the~\emph{increasing rearrangement}~$B^*$ of~$B$ as
\begin{equation} \label{CH:4:incrrearr}
B^* := (b^*, +\infty), \quad \mbox{with } \, b^* := \lim_{R \rightarrow +\infty} \left( R - \int_{-R}^R \chi_B(t) \, dt \right).
\end{equation}
Notice that, up to a set of vanishing measure---actually, a point---it holds
\begin{equation} \label{CH:4:decrincr}
B^* = \Co (\Co B)_*.
\end{equation}

The next result shows that the value of~$\I_K$ decreases when their arguments are appropriately rearranged.

\begin{prop} \label{CH:4:iotadecreasesprop}
Let~$A, B \subseteq \R$ be two sets satisfying
$$
(-\infty, \ubar{\alpha}) \subseteq A \subseteq (-\infty, \bar{\alpha}) \quad \mbox{and} \quad (\bar{\beta}, +\infty) \subseteq B \subseteq (\ubar{\beta}, +\infty),
$$
for some real numbers~$\ubar{\alpha} < \bar{\alpha}$ and~$\ubar{\beta} < \bar{\beta}$. Let~$K: \R \to \R$ be a non-negative function and suppose that
\begin{equation} \label{CH:4:Kerint}
\I_K \! \left( (-\infty, \bar{\alpha}), (\ubar{\beta}, +\infty) \right) < \infty.
\end{equation}
Then,
\begin{equation} \label{CH:4:iotadecreases}
\I_K(A_*, B^*) \le \I_K(A, B).
\end{equation}
\end{prop}
\begin{proof}
First of all, we observe that we can restrict ourselves to assume that~$A$ and~$B$ are both open sets. Indeed, if~$A$ and~$B$ are merely measurable, by the outer regularity of the Lebesgue measure there exist two sequences of open sets~$\{ A_k \}, \, \{ B_k \}$ with~$A \subseteq A_k \subseteq (-\infty, \bar{\alpha})$ and~$B \subseteq B_k \subseteq (\ubar{\beta}, +\infty)$ for any~$k \in \N$, and such that~$|A_k \setminus A|, |B_k \setminus B| \rightarrow 0$ as~$k \rightarrow +\infty$. Suppose now that~\eqref{CH:4:iotadecreases} holds with~$A_k$ and~$B_k$ respectively in place of~$A$ and~$B$. By this and the fact that, by definitions~\eqref{CH:4:decrrearr}-\eqref{CH:4:incrrearr}, it clearly holds~$A_* \subseteq (A_k)_*$ and~$B^* \subseteq (B_k)^*$ for any~$k$, we deduce that
$$
\I_K(A_*, B^*) \le \lim_{k \rightarrow +\infty} \I_K((A_k)_*, (B_k)^*) \le \lim_{k \rightarrow +\infty} \I_K(A_k, B_k) = \I_K(A, B).
$$
The last identity follows from Lebesgue's dominated convergence theorem, which can be used thanks to~\eqref{CH:4:Kerint}. In light of this, it suffices to prove~\eqref{CH:4:iotadecreases} when~$A$ and~$B$ are open sets.

Next, we recall that each open subset of the real line can be written as the union of countably many disjoint open intervals. In our setting, we have
$$
A = \bigcup_{k = 0}^{+\infty} A^{(k)}, \mbox{ with } \, A^{(k)} := \bigcup_{i = 0}^k A_i,
$$
and
$$
B = \bigcup_{k = 0}^{+\infty} B^{(k)}, \mbox{ with } \, B^{(k)} := \bigcup_{j = 0}^k B_j,
$$
for two sequences~$\{ A_i \}$,~$\{ B_j \}$ of open intervals satisfying~$A_{i_1} \cap A_{i_2} = \varnothing$ for every~$i_1 \ne i_2$ and~$B_{j_1} \cap B_{j_2} = \varnothing$ for every~$j_1 \ne j_2$, and such that~$(-\infty, \ubar{\alpha}) \subseteq A_0$ and~$(\bar{\beta}, +\infty) \subseteq B_0$, Suppose now that~\eqref{CH:4:iotadecreases} holds when~$A$ and~$B$ are the unions of finitely many disjoint open intervals. In particular,~\eqref{CH:4:iotadecreases} is true with~$A^{(k)}$ and~$B^{(k)}$ in place of~$A$ and~$B$, respectively. Hence,
\begin{equation} \label{CH:4:finiteimpliesnum}
\I_K((A^{(k)})_*, (B^{(k)})^*) \le \I_K(A^{(k)}, B^{(k)}) \le \I_K(A, B)
\end{equation}
for every~$k \in \N$. On the other hand, it is easy to see that
$$
(-\infty, \ubar{\alpha}) \subseteq (A^{(k - 1)})_* \subseteq (A^{(k)})_* \subseteq A_* \quad \mbox{and} \quad (\bar{\beta}, +\infty) \subseteq (B^{(k - 1)})^* \subseteq (B^{(k)})^* \subseteq B^*
$$
for every~$k \in \N$. Therefore, both~$|A_* \setminus (A^{(k)})_*|$ and~$|B^* \setminus (B^{(k)})^*|$ go to~$0$ as~$k \rightarrow +\infty$. Lebesgue's monotone convergence theorem then yields that
$$
\I_K(A_*, B^*) = \lim_{k \rightarrow +\infty} \I_K((A^{(k)})_*, (B^{(k)})^*). 
$$
The combination of this and~\eqref{CH:4:finiteimpliesnum} gives~\eqref{CH:4:iotadecreases}.

In light of the considerations that we just made, we are left to prove~\eqref{CH:4:iotadecreases} when~$A$ and~$B$ are unions of finitely many disjoint open intervals. Thus, we fix~$M, N \in \N \cup \{ 0 \}$ and assume that
$$
A = \bigcup_{i = 0}^M A_i \quad \mbox{and} \quad B = \bigcup_{j = 0}^{N} B_j,
$$
with
\begin{align*}
A_0 := (-\infty, a_0) \quad \mbox{and} \quad A_i := (a_{2 i - 1}, a_{2 i}) & \quad \mbox{for } i = 1, \ldots, M, \\
B_0 := (b_0, +\infty) \quad \mbox{and} \quad B_j := (b_{2 j}, b_{2 j - 1}) & \quad \mbox{for } j = 1, \ldots, N,
\end{align*}
where~$\{ a_i \}_{i = 0}^{2 M}, \{ b_j \}_{j = 0}^{2 N} \subseteq \R$ are two sets of points satisfying~$a_{i - 1} < a_{i}$ and~$b_{j - 1} < b_j$, for every~$i = 1, \ldots, 2 M$ and~$j = 1, \ldots, 2 N$. In this framework, inequality~\eqref{CH:4:iotadecreases} takes the form
\begin{equation} \label{CH:4:ID2tech0}
\sum_{\substack{i = 0, \ldots, M \\ j = 0, \ldots, N}} \int_{A_i} \int_{B_j} d\mu \ge \int_{A_*} \int_{B^*} d\mu.
\end{equation}

Clearly, when~$M = N = 0$ there is nothing to prove, as it holds~$A_* = A$ and~$B^* = B$. In case either~$M = 0$ or~$N = 0$, the verification of~\eqref{CH:4:ID2tech0} is also simple. Indeed, suppose for instance that~$N = 0$ and~$M \ge 1$. Then,~$B^* = B = (b_0, +\infty)$ and~$A_* = (-\infty, a_*)$ for some~$a_* \in \R$. Up to a set of measure zero we may write~$A_*$ as the union of the~$M + 1$ disjoint intervals~$\{ C_i \}_{i = 0}^M$ given by~$C_i = A_i - \bar{a}_i$, with~$\bar{a}_i \ge 0$ for every~$i$. Accordingly, 
$$
\int_{A_*} \int_{B^*} d\mu = \sum_{i = 0}^M \int_{C_i} \int_{b_0}^{+\infty} d\mu = \sum_{i = 0}^M \int_{A_i} \int_{b_0 + \bar{a}_i}^{+\infty} d\mu \le \sum_{i = 0}^M \int_{A_i} \int_{b_0}^{+\infty} d\mu = \int_A \int_B d\mu,
$$
that is~\eqref{CH:4:ID2tech0}. Note that the second identity follows by adding to both of the variables of the double integral the same quantity~$\bar{a}_i$. That is, we applied the change of coordinates~$x = w - \bar{a}_i$,~$y = z - \bar{a}_i$ and got
$$
\int_{C_i} \int_{b_0}^{+\infty} d\mu = \int_{C_i} \int_{b_0}^{+\infty} K(x - y) \, dx\,dy = \int_{A_i} \int_{b_0 + \bar{a}_i}^{+\infty} K(w - z) \, dw dz = \int_{A_i} \int_{b_0 + \bar{a}_i}^{+\infty} d\mu,
$$
exploiting the fact that~$K$ is translation-invariant.

As the case~$M = 0$,~$N \ge 1$ is completely analogous, we can now address the validity of~\eqref{CH:4:ID2tech0} when~$M, N \ge 1$. Recalling definitions~\eqref{CH:4:decrrearr}-\eqref{CH:4:incrrearr}, it is immediate to see that
\begin{align*}
A_* = \left( -\infty, a_* \right), & \quad \mbox{with } \, a_* = a_0 + \sum_{\ell = 1}^M |A_\ell| = a_0 + \sum_{\ell = 1}^M (a_{2 \ell} - a_{2 \ell - 1}),\\
B^* = \left( b^*, +\infty \right), & \quad \mbox{with } \, b^* = b_0 - \sum_{\ell = 1}^N |B_j| = b_0 - \sum_{\ell = 1}^N (b_{2 \ell - 1} - b_{2 \ell}).
\end{align*}
Set
\begin{align}
\label{CH:4:Cidef} C_i & := A_i - \bar{a}_i, \quad \mbox{with } \, \bar{a}_i := \sum_{\ell = 0}^{i - 1} (a_{2 \ell + 1} - a_{2 \ell}) \, \mbox{ for } i = 1, \ldots, M \mbox{ and } \bar{a}_0 := 0, \\
\label{CH:4:Djdef} D_j & := B_j + \bar{b}_j, \quad \mbox{with } \, \bar{b}_j := \sum_{\ell = 0}^{j - 1} (b_{2 \ell} - b_{2 \ell + 1}) \, \hspace{2pt} \mbox{ for } j = 1, \ldots, N \mbox{ and } \bar{b}_0 := 0.
\end{align}
The families~$\{ C_i \}_{i = 0}^M$ and~$\{ D_j \}_{j = 0}^{N}$ are both made up of consecutive open intervals. Moreover, up to sets of measure zero, we have
\begin{equation} \label{CH:4:CDfill}
A_* = \bigcup_{i = 0}^M C_i \quad \mbox{and} \quad B^* = \bigcup_{j = 0}^{N} D_j.
\end{equation}
Consequently, we can equivalently express~\eqref{CH:4:ID2tech0} as
\begin{equation} \label{CH:4:ID2tech1}
\sum_{\substack{ i = 0, \ldots, M \\ j = 0, \ldots, N}} \int_{A_i} \int_{B_j} d\mu \ge \sum_{\substack{ i = 0, \ldots, M \\ j = 0, \ldots, N}} \int_{C_i} \int_{D_j} d\mu.
\end{equation}

Fix any~$j = 1, \ldots, N$. We compute
$$
\int_{A_0} \int_{B_j} d\mu = \int_{C_0} \int_{D_j - \bar{b}_j} d\mu = \int_{C_0 + \bar{b}_j} \int_{D_j} d\mu = \int_{\left( C_0 + \bar{b}_j \right) \setminus C_0} \int_{D_j} d\mu + \int_{C_0} \int_{D_j} d\mu.
$$
Notice that the first identity follows from definitions~\eqref{CH:4:Cidef}-\eqref{CH:4:Djdef}, the second by applying to both variables of the double integral a shift of length~$\bar{b}_j$, and the third since~$C_0 \subseteq C_0 + \bar{b}_j$. Similarly,
$$
\int_{A_i} \int_{B_0} d\mu = \int_{C_i} \int_{\left( D_0 - \bar{a}_i \right) \setminus D_0} d\mu + \int_{C_i} \int_{D_0} d\mu.
$$
for every~$i = 1, \ldots, M$. Furthermore, by a translation of size~$\bar{b}_j - \bar{a}_i$, we may also write
$$
\int_{A_i} \int_{B_j} d\mu = \int_{C_i + \bar{a}_i} \int_{D_j - \bar{b}_j} d\mu = \int_{C_i + \bar{b}_j} \int_{D_j - \bar{a}_i} d\mu
$$
for every~$i = 1, \ldots, M$ and~$j = 1, \ldots, N$. Finally, again by~\eqref{CH:4:Cidef}-\eqref{CH:4:Djdef}---with~$i = j = 0$---we have
$$
\int_{A_0} \int_{B_0} d\mu = \int_{C_0} \int_{D_0} d\mu.
$$
Applying the last four identities together with~\eqref{CH:4:CDfill}, formula~\eqref{CH:4:ID2tech1} becomes
\begin{equation} \label{CH:4:ID2tech2}
\sum_{\substack{ i = 0, \ldots, M \\ j = 0, \ldots, N}} \int_{E_{i; j}} \int_{F_{j; i}} \, d\mu \ge \int_{a_0}^{a_*} \int_{b^*}^{b_0} d\mu,
\end{equation}
where we put
\begin{equation} \label{CH:4:EFdef}
\begin{aligned}
E_{0; 0} & := \{ a_0 \}, & \, F_{0; 0} & := \{ b_0 \}, && \\
E_{i; 0} & := C_i, & \, F_{0; i} & := \left( D_0 - \bar{a}_i \right) \setminus D_0, \, && \, \mbox{for } i = 1, \ldots, M, \\
E_{0; j} & := \left( C_0 + \bar{b}_j \right) \setminus C_0, & \, F_{j; 0} & := D_j, \, && \, \mbox{for } j = 1, \ldots, N, \\
E_{i; j} & := C_i + \bar{b}_j, & \, F_{j; i} & := D_j - \bar{a}_i, \, && \, \mbox{for } i = 1, \ldots, M, \, j = 1, \ldots, N.
\end{aligned}
\end{equation}

We now claim that
\begin{equation} \label{CH:4:ID2claim}
[a_0, a_*] \times [b^*, b_0] \subseteq \bigcup_{\substack{i = 0, \ldots, M \\ j = 0, \ldots, N}} \overline{E_{i; j}} \times \overline{F_{j; i}}.
\end{equation}
Observe that~\eqref{CH:4:ID2claim} is stronger than~\eqref{CH:4:ID2tech2}, and therefore that its validity would lead us to the conclusion of the proof.

Before showing that~\eqref{CH:4:ID2claim} is true, we make some considerations on the intervals~$E_{i; j}$'s and~$F_{j; i}$'s. Given a bounded non-empty interval~$I \subseteq \R$, we indicate with~$\ell(I)$ and~$r(I)$ its left and right endpoint, respectively. We have that
\begin{alignat}{2}
\label{CH:4:rlE}
r(E_{i - 1; j}) & = \ell(E_{i; j}), && \qquad \mbox{for } i = 1, \ldots, M, \, j = 0, \ldots, N, \\
\label{CH:4:rlF}
r(F_{j; i}) & = \ell(F_{j - 1; i}), && \qquad \mbox{for } i = 0, \ldots, M, \, j = 1, \ldots, N, \\
\label{CH:4:Eends}
r(E_{M; j}) & \ge a_*, && \qquad \mbox{for } j = 0, \ldots, N, \\
\label{CH:4:Fends}
\ell(F_{N; i}) & \le b^*, && \qquad \mbox{for } i = 0, \ldots, M.
\end{alignat}
To check~\eqref{CH:4:rlE}, we recall definitions~\eqref{CH:4:EFdef},~\eqref{CH:4:Cidef},~\eqref{CH:4:Djdef}, and notice that
\begin{align*}
r(E_{i - 1; j}) & = r(A_{i - 1}) - \bar{a}_{i - 1} + \bar{b}_j = a_{2 i - 2} - \bar{a}_{i} + (a_{2 i - 1} - a_{2 i - 2}) + \bar{b}_j \\
& = \ell(A_{i}) - \bar{a}_{i} + \bar{b}_j = \ell(E_{i; j})
\end{align*}
for every~$i = 1, \ldots, M$ and~$j = 0, \ldots, N$. On the other hand, it holds
\begin{align*}
r(E_{M; j}) & = r(A_M) - \bar{a}_M + \bar{b}_j = a_{2 M} - \sum_{\ell = 0}^{M - 1} (a_{2 \ell + 1} - a_{2 \ell}) + \bar{b}_j \\
& = a_0 + \sum_{\ell = 1}^{M} (a_{2 \ell} - a_{2 \ell - 1}) + \bar{b}_j \ge a_*,
\end{align*}
which gives~\eqref{CH:4:Eends}. Items~\eqref{CH:4:rlF} and~\eqref{CH:4:Fends} follow analogously.

In view of formulas~\eqref{CH:4:rlE},~\eqref{CH:4:Eends} and~\eqref{CH:4:rlF},~\eqref{CH:4:Fends}, we immediately deduce that
\begin{equation} \label{CH:4:Ecover}
[a_0, a_*] \subseteq \bigcup_{i = 0}^M \overline{E_{i ; j}} \quad \mbox{for any } j = 0, \ldots, N
\end{equation}
and
\begin{equation} \label{CH:4:Fcover}
[b^*, b_0] \subseteq \bigcup_{j = 0}^N \overline{F_{j;i}} \quad \mbox{for any } i = 0, \ldots, M,
\end{equation}
respectively---recall that~$\ell(E_{0; j}) = a_0$ and~$r(F_{0; i}) = b_0$ for any such~$j$ and~$i$.

On top of the previous facts, we also claim that
\begin{equation} \label{CH:4:ID2claim2pre}
\ell(E_{i; j}) > \ell(E_{i; j - 1}) \quad \mbox{for every } i = 1, \ldots, M, \, j = 1, \ldots, N
\end{equation}
and
\begin{equation} \label{CH:4:ID2claim2}
r(F_{j; i}) < r(F_{j; i - 1}) \quad \mbox{for every } i = 1, \ldots, M, \, j = 1, \ldots, N.
\end{equation}
Indeed, for~$i = 1, \ldots, M$ and~$j = 1, \ldots, N$ we have
$$
r(F_{j; i}) = r(D_j) - \bar{a}_i = r(D_j) - \bar{a}_{i - 1} - (a_{2 i - 1} - a_{2 i - 2}) < r(D_j) - \bar{a}_{i - 1} = r(F_{j; i}).
$$
This proves~\eqref{CH:4:ID2claim2}, while~\eqref{CH:4:ID2claim2pre} can be checked in a similar fashion.

Thanks to the previous remarks, we can now address the proof of~\eqref{CH:4:ID2claim}. Let
\begin{equation} \label{CH:4:pdef}
p = (x, y) \in [a_0, a_*] \times [b^*, b_0]
\end{equation}
and suppose by contradiction that~$p$ does not belong to the right-hand side of~\eqref{CH:4:ID2claim}. I.e.,
\begin{equation} \label{CH:4:pcontradict}
p \notin \overline{E_{i; j}} \times \overline{F_{j; i}} \quad \mbox{for every } i = 0, \ldots, M \mbox{ and } j = 0, \ldots, N.
\end{equation}
In light of~\eqref{CH:4:Ecover}, we know that in correspondence to every~$j = 0, \ldots, N$ we can pick an~$i_j \in \{ 0, \ldots, M \}$ in such a way that
\begin{equation} \label{CH:4:xinE}
x \in \overline{E_{i_j; j}}.
\end{equation}
We claim that
\begin{equation} \label{CH:4:ijnonincr}
\{i_j\}_{j = 0}^N \mbox{ is non-increasing}.
\end{equation}
Indeed, suppose that we have constructed the (finite) sequence~$\{ i_\ell \}$ up to the index~$\ell = j - 1$, with~$j \in \{ 1, \ldots, N \}$. Of course, when~$i_{j - 1} = M$ we necessarily have~$i_{j} \le i_{j - 1}$. On the other hand, if~$i_{j - 1} \le M - 1$, using~\eqref{CH:4:ID2claim2pre} and~\eqref{CH:4:rlE}, we infer that
$$
\ell(E_{i_{j - 1} + 1; j}) > \ell(E_{i_{j - 1} + 1; j - 1}) = r(E_{i_{j - 1}; j - 1}) \ge x.
$$
Hence, also in this case~$i_j$ falls within the set~$\{ 0, \ldots, i_{j - 1} \}$ and~\eqref{CH:4:ijnonincr} is established.

Next, by comparing~\eqref{CH:4:xinE} and~\eqref{CH:4:pcontradict}, we notice that~$y \notin \overline{F_{j; i_j}}$. This amounts to say that, for every index~$j = 0, \ldots, N$,
\begin{equation} \label{CH:4:y<lory>r}
\mbox{ either } \, y < \ell(F_{j; i_j}) \, \mbox{ or } \, y > r(F_{j; i_j}).
\end{equation}
We now claim that the latter possibility cannot occur, i.e.,~that
\begin{equation} \label{CH:4:ID2claimcontr}
y < \ell(F_{j; i_j})
\end{equation}
for every~$j = 0, \ldots, N$. Note that~\eqref{CH:4:ID2claimcontr} would lead us to a contradiction. Indeed, by using it with~$j = N$ and in combination with~\eqref{CH:4:pdef} and~\eqref{CH:4:Fends}, we would get
$$
b^* \le y < \ell(F_{N; i_N}) \le b^*,
$$
which is clearly impossible. Therefore, to finish the proof we are only left to show that~\eqref{CH:4:ID2claimcontr} holds true for every~$j = 0, \ldots, N$. To achieve this, we argue inductively. First, we check that~\eqref{CH:4:ID2claimcontr} is verified for~$j = 0$. Indeed, by~\eqref{CH:4:pdef} and~\eqref{CH:4:EFdef},
$$
y \le b_0 = r(F_{0; i_0}),
$$
and thus~\eqref{CH:4:y<lory>r} yields that~$y < \ell(F_{0; i_0})$. Secondly, we pick any~$j \in \{ 1, \ldots, N\}$ and assume that~\eqref{CH:4:ID2claimcontr} is valid with~$j - 1$ in place of~$j$. Then, recalling~\eqref{CH:4:rlF},~\eqref{CH:4:ijnonincr} and possibly~\eqref{CH:4:ID2claim2} (applied iteratively), we get that
$$
y < \ell(F_{j - 1; i_{j - 1}}) = r(F_{j; i_{j - 1}}) \le r(F_{j; i_j}).
$$
By comparing this with~\eqref{CH:4:y<lory>r}, we finally deduce that claim~\eqref{CH:4:ID2claimcontr} holds true. Thus, the proof is complete.
\end{proof}

\subsection{Vertical rearrangements and the~$s$-perimeter} \label{CH:4:rearrangsub}

We now take advantage of Proposition~\ref{CH:4:iotadecreasesprop} to show that~$\Per_s$ decreases under vertical rearrangements. Given a set~$E \subseteq \R^{n + 1}$, we consider the function~$w_E: \R^n \to \R$ defined by
$$
w_E(x) := \lim_{R \rightarrow +\infty} \left( \int_{-R}^R \chi_{E}(x, t) \, dt - R \right)
$$
for every~$x \in \R^n$, together with its subgraph~$E_\star := \Sg(w_E)$. We have the following result.

\begin{proof}[Proof of Theorem \ref{CH:4:Persdecreases}]
Denote with~$G$ either the set~$E$ or its rearrangement~$E_\star$. Observe that~$E$ and~$E_\star$ coincide outside of~$\Omega^\infty$, and are both given by the subgraph of the same function~$v: \Co \Omega \to \R$. Hence,
\begin{equation} \label{CH:4:GoutOmega}
G \setminus \Omega^\infty = \Big\{ (x, t) \in \left( \Co \Omega \right) \times \R \,|\, t < v(x) \Big\}.
\end{equation}
It is also clear that~$E_\star$ satisfies~\eqref{CH:4:EboundedinOmega}. Accordingly,
\begin{equation} \label{CH:4:GboundedinOmega}
\Omega \times (-\infty, -M) \subseteq G \cap \Omega^\infty \subseteq \Omega \times (-\infty, M).
\end{equation}

We compute
\begin{align*}
\Per_s(G, \Omega^M) & = \Ll_s(G \cap \Omega^M,\Co G \cap \Omega^M) + \Ll_s(G \cap \Omega^M, \Co G \setminus \Omega^M) + \Ll_s(G \setminus \Omega^M, \Co G \cap \Omega^M) \\
& = \Ll_s(G \cap \Omega^M,\Co G \cap \Omega^M) \\
& \quad + \Ll_s(G \cap \Omega^M, \Co G \cap (\Omega^\infty \setminus \Omega^M)) + \Ll_s(G \cap \Omega^M, \Co G \setminus \Omega^\infty) \\
& \quad + \Ll_s(G \cap (\Omega^\infty \setminus \Omega^M), \Co G \cap \Omega^M) + \Ll_s(G \setminus \Omega^\infty, \Co G \cap \Omega^M) \\
& = \Ll_s(G \cap \Omega^\infty,\Co G \cap \Omega^\infty) - \Ll_s(G \cap (\Omega^\infty \setminus \Omega^M), \Co G \cap (\Omega^\infty \setminus \Omega^M)) \\
& \quad + \Ll_s(G \cap \Omega^M, \Co G \setminus \Omega^\infty) + \Ll_s(G \setminus \Omega^\infty, \Co G \cap \Omega^M).
\end{align*}
Thanks to~\eqref{CH:4:GboundedinOmega}, we may write
$$
\Ll_s(G \cap (\Omega^\infty \setminus \Omega^M), \Co G \cap (\Omega^\infty \setminus \Omega^M)) = \Ll_s(\Omega \times (-\infty, -M), \Omega \times (M, +\infty)) =: C_M.
$$
Note that~$C_M$ is a constant depending only on~$n$,~$s$,~$\Omega$, and~$M$. Moreover, using~\eqref{CH:4:GoutOmega} and again~\eqref{CH:4:GboundedinOmega}, we have
$$
\Ll_s(G \cap \Omega^M, \Co G \setminus \Omega^\infty) = \Ll_s(G \cap \Omega^\infty, \Co G \setminus \Omega^\infty) - D^{(1)}_M
$$
and
$$
\Ll_s(G \setminus \Omega^\infty, \Co G \cap \Omega^M) = \Ll_s(G \setminus \Omega^\infty, \Co G \cap \Omega^\infty) - D^{(2)}_M,
$$
where~$D^{(1)}_M := \Ll_s(\Omega \times (-\infty, -M), \Co \Sg(v) \setminus \Omega^\infty)$ and~$D^{(2)}_M := \Ll_s(\Sg(v) \setminus \Omega^\infty, \Omega \times (M, +\infty))$ are constants depending only on~$n$,~$s$,~$\Omega$,~$M$, and~$v$. Putting together the last four identities, we find that
\begin{align*}
\Per_s(G, \Omega^M) & = \Ll_s(G \cap \Omega^\infty,\Co G \cap \Omega^\infty) + \Ll_s(G \cap \Omega^\infty, \Co G \setminus \Omega^\infty) + \Ll_s(G \setminus \Omega^\infty, \Co G \cap \Omega^\infty) \\
& \quad - C_M - D^{(1)}_M - D^{(2)}_M.
\end{align*}
In particular, inequality~\eqref{CH:4:Persdecreases} will be verified if we prove that
\begin{equation} \label{CH:4:LsleLsequiv}
\begin{aligned}
\Ll_s(E_\star \cap \Omega^\infty, \Co E_\star \cap \Omega^\infty) & \le \Ll_s(E \cap \Omega^\infty,\Co E \cap \Omega^\infty), \\
\Ll_s(E_\star \cap \Omega^\infty, \Co E_\star \setminus \Omega^\infty) & \le \Ll_s(E \cap \Omega^\infty, \Co E \setminus \Omega^\infty), \\
\Ll_s(E_\star \setminus \Omega^\infty, \Co E_\star \cap \Omega^\infty) & \le \Ll_s(E \setminus \Omega^\infty, \Co E \cap \Omega^\infty).
\end{aligned}
\end{equation}

Set
$$
G(x) := \Big\{ t \in \R \,|\, (x, t) \in G \Big\} \quad \mbox{for } x \in \R^n
$$
and
$$
K_{a}(t) := \frac{1}{\left( a^2 + t^2 \right)^{\frac{n + 1 + s}{2}}} \quad \mbox{for } a, t \in \R.
$$
Using the notation of~\eqref{CH:4:IKdef} and Fubini's theorem, we may write
\begin{equation} \label{CH:4:LsasIK}
\begin{aligned}
\Ll_s(G \cap \Omega^\infty,\Co G \cap \Omega^\infty) & = \int_{\Omega} \int_{\Omega} \I_{K_{|x - y|}} \! \left( G(x), \Co G(y) \right) dx\,dy,\\
\Ll_s(G \cap \Omega^\infty, \Co G \setminus \Omega^\infty) & = \int_{\Omega} \int_{\Co \Omega} \I_{K_{|x - y|}} \! \left( G(x), \Co G(y) \right) dx\,dy,\\
\Ll_s(G \setminus \Omega^\infty, \Co G \cap \Omega^\infty) & = \int_{\Co \Omega} \int_{\Omega} \I_{K_{|x - y|}} \! \left( G(x), \Co G(y) \right) dx\,dy.
\end{aligned}
\end{equation}
Recalling the definition of decreasing rearrangement of a subset of the real line introduced in~\eqref{CH:4:decrrearr}, we observe that~$E(x)_* = (-\infty, w_E(x)) = E_\star(x)$ for every~$x \in \R^n$. Also notice that~$\I_{K_a} \! \left( (-\infty, \alpha), (\beta, +\infty) \right) < \infty$ for every~$\alpha, \beta \in \R$ and~$a \ne 0$. By this, we are allowed to apply Proposition~\ref{CH:4:iotadecreasesprop} and deduce that
$$
\I_{K_{|x - y|}} \! \left( E_\star(x), \Co E_\star(y) \right) \le \I_{K_{|x - y|}} \! \left( E(x), \Co E(y) \right) \quad \mbox{for a.e.~} x, y \in \R^n,
$$
where we also took advantage of property~\eqref{CH:4:decrincr}. In view of~\eqref{CH:4:LsasIK}, this last inequality ensures the validity of~\eqref{CH:4:LsleLsequiv}. The proof is thus finished.
\end{proof}

\section{Minimizers} \label{CH:4:Linftysec}

This section is devoted to the study of the minimizers of~$\F$.
As observed in the Introduction, we will prove the existence of minimizers with the aid of an appropriate
approximation procedure, which makes use of the ``truncated functionals''~$\F^M$ and of their own minimizers.
For this reason, we introduce straight away the following auxiliary functional spaces.
Given a bounded open set~$\Omega\subseteq\R^n$,~$s\in(0,1)$ and~$M\ge0$, we define
\bgs{
&\qquad\B\W^s(\Omega):=\{u\in\W^s(\Omega)\,|\,u|_\Omega\in L^\infty(\Omega)\}\\
&\textrm{and}\quad\B_M\W^s(\Omega):=\{u\in\B\W^s(\Omega)\,|\,\|u\|_{L^\infty(\Omega)}\le M\}.
}
Moreover, given a function~$\varphi:\Co\Omega\to\R$, we define
\bgs{
&\qquad\B\W^s_\varphi(\Omega):=\{u\in\B\W^s(\Omega)\,|\,u=\varphi\textrm{ a.e. in }\Co\Omega\}\\
&\textrm{and}\quad\B_M\W^s_\varphi(\Omega):=\{u\in\B_M\W^s(\Omega)\,|\,u=\varphi\textrm{ a.e. in }\Co\Omega\}.
}

We begin by recalling the definition of minimizer in the context of the Dirichlet problem. Let~$\Omega\subseteq\R^n$
be a bounded open set with Lipschitz boundary,~$s\in(0,1)$ and let~$\varphi:\Co\Omega\to\R$.
We say that a function~$u\in\W^s_\varphi$ is a minimizer of~$\F$ in~$\W^s_\varphi(\Omega)$ if
\bgs{
\iint_{Q(\Omega)}\left\{\G\left(\frac{u(x)-u(y)}{|x-y|}\right)-\G\left(\frac{v(x)-v(y)}{|x-y|}\right)\right\}
\frac{dx\,dy}{|x-y|^{n-1+s}}\le0,
}
for every~$v\in\W^s_\varphi(\Omega)$.

It is now convenient to point out the following useful result, which is easily obtained by arguing as in the proof of
Lemma~\ref{CH:4:NMdomainlem}, exploiting formula~\eqref{CH:4:nonlocal_explicit} and
the global Lipschitzianity of~$\G$---see~\eqref{CH:4:Lip_Gcal}.

\begin{lemma} \label{CH:4:tartariccio}
Let~$n \ge 1$,~$s \in (0, 1)$,~$M \ge 0$,~$\Omega\subseteq\Rn$ be a bounded open set with Lipschitz boundary
and let~$\varphi:\Co\Omega\to\R$.
There exists a constant~$C>0$, depending only on~$n$,~$s$ and~$\Omega$, such that
\bgs{
\iint_{Q(\Omega)}\left|\G\left(\frac{u(x)-u(y)}{|x-y|}\right)
-\G\left(\frac{v(x)-v(y)}{|x-y|}\right)\right|\frac{dx\,dy}{|x-y|^{n-1+s}}
\le C\,\Lambda\|u-v\|_{W^{s,1}(\Omega)},
}
for every~$u,\,v\in\W^s_\varphi(\Omega)$, with~$\Lambda$ as defined in~\eqref{CH:4:gintegr}.
Moreover, we have the identity
\eqlab{\label{CH:4:Lem_tartariccio_eqn1}
\F^M(u,\Omega)-\F^M(v,\Omega)=
\iint_{Q(\Omega)} \left\{ \G\left(\frac{u(x)-u(y)}{|x-y|}\right)
-\G\left(\frac{v(x)-v(y)}{|x-y|}\right) \right\} \frac{dx\,dy}{|x-y|^{n-1+s}}.
}
As a consequence, if~$u, u_k\in\W^s_\varphi(\Omega)$ are such that~$\|u-u_k\|_{W^{s,1}(\Omega)}\to0$ as~$k\to\infty$,
then
\bgs{
\lim_{k\to\infty}\F^M(u_k,\Omega)=\F^M(u,\Omega).
}
\end{lemma}

\begin{remark}\label{CH:4:tartapower_Remark}
In this Remark we collect the following straightforward but important
consequences of Lemma~\ref{CH:4:tartariccio}:
\begin{itemize}
\item[(i)] it guarantees that the definition of minimizer is well posed;
\item[(ii)] it provides an equivalent characterization of a minimizer of~$\F$ in~$\W^s_\varphi(\Omega)$
as a function~$u\in\W^s_\varphi(\Omega)$ that
minimizes~$\F^M(\,\cdot\,,\Omega)$, i.e. such that
\bgs{
\F^M(u,\Omega)=\inf\left\{\F^M(v,\Omega)\,|\,v\in\W^s_\varphi(\Omega)\right\};
}
\item[(iii)] by point~$(ii)$ and by the strict convexity of~$\F^M$---see point~$(ii)$ of Lemma~\ref{CH:4:conv_func}---we obtain that a minimizer of~$\F$ in~$\W^s_\varphi(\Omega)$,
if it exists, is unique;
\item[(iv)] as a consequence of the density of the spaces~$C^\infty_c(\Omega)$
and~$W^{s,1}(\Omega)\cap L^\infty(\Omega)$ in the fractional Sobolev
space~$W^{s,1}(\Omega)$---see, e.g., Appendix~\ref{CH:4:Density_appendix}---Lemma~\ref{CH:4:tartariccio} implies that to verify the minimality of~$u\in\W^s_\varphi(\Omega)$ we can limit ourselves to consider
either competitors~$v\in\W^s_\varphi(\Omega)$ such that~$v|_\Omega\in C^\infty_c(\Omega)$,
or~$v\in\B\W^s_\varphi(\Omega)$.
\end{itemize}
\end{remark}

In light of point~$(ii)$ of Remark~\ref{CH:4:tartapower_Remark}, we could have considered as definition of minimizer
just that of a function~$u\in\W^s_\varphi(\Omega)$ that minimizes the functional~$\F^0$---or
the functional~$\F^M$, for some fixed~$M>0$---in~$\W^s_\varphi(\Omega)$.
However, we remark that such a definition is not very helpful when trying to prove existence results and indeed it presents some difficulties, first of all the fact that the functional~$\F^M$ in general is not non-negative in the space~$\W^s_\varphi(\Omega)$ and may indeed change sign---see Example~\ref{CH:4:Exe_neg_nonloc}. Hence,
lower semicontinuity and compactness properties are not straightforward.

\smallskip

Now we turn our attention to the Euler-Lagrange equation satisfied by minimizers.

We recall that, given a bounded open set~$\Omega\subseteq\R^n$ and~$s\in(0,1)$, we say that
a measurable function~$u:\R^n\to\R$
is a weak solution of~$\h u=0$ in~$\Omega$ if
\eqlab{\label{CH:4:sanzo}
\langle \h u,v\rangle=0\quad\mbox{for every }v\in C^\infty_c(\Omega).
}

\begin{remark}\label{CH:4:crocodile_dens_rmk}
Notice that, if~$\Omega$ has Lipschitz boundary, then, by density, in~\eqref{CH:4:sanzo}
we can as well consider~$v\in\W^s_0(\Omega)$ as test function.
Indeed, in light of Corollary~\ref{CH:4:FHI_corollary} we have that
\bgs{
\|v\|_{W^{s,1}(\R^n)}\le C(n,s,\Omega)\|v\|_{W^{s,1}(\Omega)}\quad\mbox{for every }v\in\W^s_0(\Omega).
}
Hence, since~$\langle\h u,\,\cdot\,\rangle\in\big(W^{s,1}(\R^n)\big)^*$, by the density of~$C^\infty_c(\Omega)$
in~$W^{s,1}(\Omega)$ we find that~\eqref{CH:4:sanzo} implies that
\bgs{
\langle\h u,v\rangle=0\quad\mbox{for every }v\in\W^s_0(\Omega).
}
\end{remark}

Exploiting the convexity of the functionals~$\F^M$, we can prove the equivalence between weak solutions (with ``finite energy'')
and minimizers.

\begin{lemma}\label{CH:4:weak_implies_min_lemma}
Let~$n \ge 1$,~$s \in (0, 1)$,~$\Omega \subseteq \R^n$ be a bounded open set with Lipschitz boundary, and~$u \in \W^s(\Omega)$. Then,~$u$ is a weak solution of~$\h u = 0$ in~$\Omega$ if and only if~$u$ is a minimizer of~$\F$ in~$\Omega$.
\end{lemma}

\begin{proof}
Suppose that~$u$ is a weak solution, let~$v\in\W^s_u(\Omega)$ and define~$w:=v-u$. Notice that,
since~$w\in\W^s_0(\Omega)$, by Remark~\ref{CH:4:crocodile_dens_rmk} we have
\bgs{
\langle\h u,w\rangle=0.
}
Now we observe that the convexity of~$\G$ implies that
\bgs{
\G(t)-\G(\tau)\ge G(\tau)(t-\tau)\quad\mbox{for every }t,\,\tau\in\R.
}
Thus, by~\eqref{CH:4:Lem_tartariccio_eqn1} we obtain
\bgs{
\F^M(v,\Omega)-\F^M(u,\Omega)\ge\langle\h u,w\rangle=0.
}
Since~$v\in\W^s_u(\Omega)$ is arbitrary, the function~$u$ minimizes~$\F^M(\,\cdot\,,\Omega)$ in~$\W^s_u(\Omega)$,
and hence---by point~$(ii)$ of Remark~\ref{CH:4:tartapower_Remark}---$u$ is a minimizer of~$\F$ in~$\Omega$,
in the sense of Definition~\ref{CH:4:minim_tarta_def}.

To conclude the proof of the Lemma,
the converse implication follows by
point~$(ii)$ of Remark~\ref{CH:4:tartapower_Remark} and Lemma~\ref{CH:4:1varlem}.
\end{proof}

It is interesting to observe that Lemmas~\ref{CH:4:easylemma} and~\ref{CH:4:weak_implies_min_lemma} imply straight away that
the set of minimizers of~$\F$ is closed in~$\W^s(\Omega)$, with respect to almost everywhere convergence.

\begin{prop}
Let~$n \ge 1$,~$s \in (0, 1)$ and let~$\Omega \subseteq \R^n$ be a bounded open set with Lipschitz boundary.
Let~$\{u_k\}\subseteq \W^s(\Omega)$ be such
that each~$u_k$ is a minimizer of~$\F$ in~$\Omega$.
If~$u_k\to u$ a.e. in~$\R^n$, for some function~$u\in\W^s(\Omega)$, then~$u$ is a minimizer of~$\F$ in~$\Omega$.
\end{prop}

Before going on, we briefly explain why we consider condition~\eqref{CH:4:tartainfty} to be too restrictive in our
framework---even if at first glance it seems to be necessary, since it is required in order to guarantee that~$\F$ is well defined
on~$\W^s_\varphi(\Omega)$---and why it makes sense to expect the existence of minimizers even when the
exterior data~$\varphi:\R^n\to\R$ does not satisfy~\eqref{CH:4:tartainfty}.

First of all, we observe that if~$\Omega\subseteq\R^n$ is a bounded open set with Lipschitz boundary
and~$\varphi:\R^n\to\R$ is bounded in a neighborhood of~$\Omega$, then it is readily seen
that~$\varphi$ satisfies~\eqref{CH:4:tartainfty}
if and only if
\eqlab{\label{CH:4:tail_cond_bla}
\int_{\R^n}\frac{|\varphi(y)|}{1+|y|^{n+s}}\,dy<\infty.
}
 We remark in particular that~\eqref{CH:4:tail_cond_bla} forces~$\varphi$ to grow sublinearly at infinity.

Let now~$u:\R^n\to\R$ be such that~$u=\varphi$ almost everywhere in~$\Co\Omega$ and suppose that $u\in C^2(B_r(x))$,
for some~$x\in\Omega$ and~$r>0$. Then, the condition~\eqref{CH:4:tail_cond_bla} is the same condition needed in order to guarantee the well definiteness of the fractional Laplacian
\bgs{
(-\Delta)^\frac{s}{2}u(x)=\frac{1}{2}\int_{\R^n}\frac{2u(x)-u(x+y)+u(x-y)}{|y|^{n+s}}\,dy.
}
On the other hand, as observed in Lemma~\ref{CH:4:classical_form_reg_func},
the operator~$\h u$ is well defined at~$x$ just thanks to the local regularity of~$u$,
with no need of assumptions about the growth of~$u$ at infinity.
We further mention that condition~\eqref{CH:4:tail_cond_bla} is needed in order to define
the fractional $\frac{s}{2}$-Laplacian of a function as a tempered distribution. Contrarily,
we can always define the operator~$\h u$ in the distributional sense of~\eqref{CH:4:weak_opossum_curv},
without having to make any assumption on the function~$u$, besides measurability.

Also, we recall that we have a definition of minimizer of~$\F$, namely Definition~\ref{CH:4:minim_tarta_def},
which---as ensured by Lemma~\ref{CH:4:tartariccio}---makes sense without having to impose any
restriction on the exterior data.
 
Thus, differently to what happens in the context of the fractional Laplacian, where condition~\eqref{CH:4:tail_cond_bla}
is totally natural, in our framework it seems to be unnecessarily restrictive.

Let us now switch our attention to the geometric situation, which corresponds to the choice~$g=g_s$.
Let~$\Omega\subseteq\R^n$ be a bounded open set with~$C^2$ boundary.
We consider as exterior data a continuous function~$\varphi\in C(\R^n)$, but we make no assumption on the behavior
of~$\varphi$ at infinity. Then, we know that there exists a function~$u:\R^n\to\R$ such
that~$u\in C^\infty(\Omega)\cap C(\overline{\Omega})$ and~$u=\varphi$ in~$\R^n\setminus\overline{\Omega}$,
whose subgraph~$\Sg(u)$ is locally $s$-minimal in the cylinder~$\Omega^\infty$.
The existence follows from~\cite[Theorem~1.1]{graph} and~Theorem \ref{CH:2:nonparametric_exist_teo},
while the interior smoothness is guaranteed by~\cite[Theorem~1.1]{CaCo}.
Thus, we know that in this case the ``geometric problem'' of (locally) minimizing the $s$-perimeter in~$\Omega^\infty$
with respect to the exterior data~$\Sg(\varphi)\setminus\Omega^\infty$ has a solution, which is given by the subgraph of
a function~$u$, even if~$\varphi$ does not satisfy~\eqref{CH:4:tartainfty}.
To go one step further, we now observe that the function~$u$ is actually the minimizer of~$\F_s$ in~$\W^s_\varphi(\Omega)$.
Indeed, thanks to the smoothness of~$u$ and the minimality of~$\Sg(u)$---see~\cite[Theorem~5.1]{CRS10}---we have that
\bgs{
\h_su(x)=\I_s[\Sg(u)](x,u(x))=0\qquad\mbox{for every }x\in\Omega,
}
and hence, by Proposition~\ref{CH:4:BHprop_curvature},
\bgs{
\langle \h_su,v\rangle=0\qquad\mbox{for every }v\in C^\infty_c(\Omega).
}
Then, by Lemma~\ref{CH:4:weak_implies_min_lemma}, we conclude that~$u$ minimizes~$\F_s$ in~$\W^s_\varphi(\Omega)$.

For a more detailed discussion about the equivalence between stationary functions, minimizers of~$\F$ and ``geometric minimizers'', in a more general situation, we refer to the forthcoming proof of Theorem~\ref{CH:4:equiv_intro}
in Section~\ref{CH:4:Geom_min_proof_Section}.

\subsection{Minimizers of the truncated functionals $\F^M$}\label{CH:4:TRUNC_FUN_SEC}
As we have just anticipated, we are going to prove the existence of minimizers of~$\F$ by making use of the minimizers of
the truncated functionals~$\F^M$.
In order to motivate why we should expect this strategy to work, let us indulge a little longer in the discussion about the geometric 
situation.

Again, we consider a bounded open set~$\Omega\subseteq\R^n$ with~$C^2$ boundary and we fix as exterior data a continuous function~$\varphi\in C(\R^n)$. As a first step, we observe that~\cite[Theorem~1.1]{graph}
says that if~$E\subseteq\R^{n+1}$ is a set which is locally $s$-minimal in the cylinder~$\Omega^\infty$
and~$E\setminus\Omega^\infty=\Sg(\varphi)\setminus\Omega^\infty$, then~$E$ is globally a subgraph, that is,~$E=\Sg(u)$,
for some function~$u:\R^n\to\R$ such that~$u\in C(\overline{\Omega})$
and~$u=\varphi$ in~$\R^n\setminus\overline{\Omega}$.

Therefore, we are reduced to prove the existence of
a set~$E$ which is locally $s$-minimal in~$\Omega^\infty$, with exterior data~$\Sg(\varphi)\setminus\Omega^\infty$.

We recall that, in order to do this, the argument exploited in the proof of Corollary \ref{CH:2:loc_min_set_cor} is the following. We first consider the minimization problem in the truncated cylinders~$\Omega^k$, that is, we take a set~$E_k\subseteq\R^{n+1}$ which is
$s$-minimal in~$\Omega^k$ and such that~$E_k\setminus\Omega^k=\Sg(\varphi)\setminus\Omega^k$. The existence of such sets
is guaranteed by~\cite[Theorem~3.2]{CRS10}, since~$\Omega^k$ is a bounded open set with Lipschitz boundary.
Then, a compactness argument which exploits uniform perimeter estimates for $s$-minimal sets guarantees the existence of
a set~$E$ such that~$\chi_{E_k}\to\chi_E$ in~$L^1_{\loc}(\R^{n+1})$, up to subsequences. Notice that we have~$E\setminus\Omega^\infty=\Sg(\varphi)\setminus\Omega^\infty$. Finally, the $s$-minimality
of the approximating sets~$E_k$ implies that the limit set~$E$ is locally $s$-minimal in~$\Omega^\infty$---we refer
the interested reader to Chapter \ref{CH_Appro_Min} for the rigorous details of the argument.

Now we recall that, when restricted to the functional space~$\B_M\W^s(\Omega)$,
the functional~$\F^M_s$ corresponds to the $s$-fractional perimeter in the truncated cylinder~$\Omega^M$---by
Proposition~\ref{CH:4:per_of_subgraph_prop}.
Hence, the problem of finding a set~$E_k\subseteq\R^{n+1}$ which is $s$-minimal in~$\Omega^k$,
with respect to the exterior data~$\Sg(\varphi)\setminus\Omega^k$ corresponds,
when~$k\ge\|\varphi\|_{L^\infty(\Omega)}$, to the functional problem of minimizing~$\F^k_s$
in the space~$\B_k\W^s_\varphi(\Omega)$---see Proposition~\ref{CH:4:rearrangement_truncated_Prop}.
As we are going to prove in a moment, by making use of the direct method of the Calculus of Variations and exploiting the convexity of~$\F^k_s$, this minimizing problem has a unique solution~$u_k$.
Then, if we want to follow the same strategy exploited in the geometric situation, we should aim to prove that~$u_k\to u$
almost everywhere in~$\R^n$, up to subsequences. This step is quite simple when working with sets, thanks to universal perimeter estimates. On the other hand, in the functional setting the situation is a little trickier and the existence of
a limit function~$u$ is ensured by the uniform estimates provided by Proposition~\ref{CH:4:Ws1prop}.
Finally, we can exploit the minimality of the functions~$u_k$ in the space~$\B_k\W^s_\varphi(\Omega)$
to obtain the minimality of~$u$ in~$\W^s_\varphi(\Omega)$.

\smallskip

Let us now get to the proofs of the aforementioned results.

We begin by observing that~$\F^M$ is lower semicontinuous in~$\B_M\W^s(\Omega)$ with respect to pointwise convergence
almost everywhere.

\begin{lemma}[\bf Semicontinuity] \label{CH:4:semicontlem}
Let~$n \ge 1$,~$s \in (0, 1)$,~$M > 0$, and~$\Omega \subseteq \R^n$ be an open set.
Let~$\{ u_k \}_{k \in \N}\subseteq\B_M\W^s(\Omega)$ be a sequence of functions
converging to some~$u: \R^n \to \R$ a.e.~in~$\R^n$.
Then,
\bgs{
\F^M(u,\Omega)\le\liminf_{k\to\infty}\F^M(u_k,\Omega).
}
\end{lemma}

\begin{proof}
The proof is a consequence of Fatou's lemma, applied separately to the functionals~$\A$ and~$\Nl^M$.
Notice that, in order to use this result with~$\Nl^M$, 
the uniform bound
\bgs{
\|u_k\|_{L^\infty(\Omega)}\le M
}
is fundamental to guarantee that the quantity inside square brackets in~\eqref{CH:4:NMldef}
is non-negative---recall that~$\overline{G} \ge 0$ by definition~\eqref{CH:4:Gbardef}.
\end{proof}


Next is a compactness result for sequences uniformly bounded with respect to~$\A$. 

\begin{lemma}[\bf Compactness] \label{CH:4:complem}
Let~$n \ge 1$,~$s \in (0, 1)$, and~$\Omega \subseteq \R^n$ be a bounded open set with Lipschitz boundary. Let~$\{ u_k \}_{k \in \N}$ be a sequence functions~$u_k: \Omega \to \R$ satisfying
$$
\sup_{k \in \N} \left(\|u_k\|_{L^1(\Omega)}+\A(u_k, \Omega)\right) < \infty.
$$
Then, there exists a function~$u \in W^{s, 1}(\Omega)$ such that~$\{ u_k \}$ converges to~$u$ a.e.~in~$\Omega$,
up to a subsequence.
\end{lemma}

Lemma~\ref{CH:4:complem} follows at once from the compact embedding~$W^{s, 1}(\Omega) \hookrightarrow \hookrightarrow L^1(\Omega)$---see,~e.g.,~\cite[Theorem~7.1]{HitGuide}---and Lemma~\ref{CH:4:Adomainlem}.

By combining the last two results, we easily obtain the existence of a (unique) minimizer~$u_M$ of~$\F^M(\,\cdot\,, \Omega)$ among all functions in~$\B_M \W^s(\Omega)$ with fixed values outside of~$\Omega$.

\begin{prop}\label{CH:4:ertyui}
Let~$n\ge1$,~$s \in (0, 1)$,~$\Omega \subseteq \R^n$ be a bounded open set with Lipschitz boundary, and~$\varphi: \Co \Omega \to \R$ be a given function. For every~$M > 0$, there exists a unique minimizer~$u_M$ of~$\F^M(\,\cdot\,, \Omega)$ in~$\B_M\W^s_\varphi(\Omega)$, i.e., there exists a unique~$u_M \in \B_M \W^s_\varphi(\Omega)$ for which
\begin{equation} \label{CH:4:uMmin}
\F^M(u_M,\Omega) = \inf \Big\{ \F^M(v, \Omega) \,|\, v \in \B_M\W_\varphi^s(\Omega) \Big\}.
\end{equation}
\end{prop}
\begin{proof}
Since~$\B_M\W^s_\varphi(\Omega)$ is a convex subset of~$\W^s_\varphi(\Omega)$,
the uniqueness of the minimizer of~$\F^M(\,\cdot\,, \Omega)$ within~$\B_M\W^s_\varphi(\Omega)$
is a consequence of the strict convexity of~$\F^M(\,\cdot\,, \Omega)$---see point~(ii) of Lemma~\ref{CH:4:conv_func}.
Therefore, we are only left to establish its existence.

Let~$\{ u^{(k)} \} \subseteq \B_M\W^s_\varphi(\Omega)$ be a minimizing sequence, that is
$$
\lim_{k\to\infty}\F^M(u^{(k)},\Omega)=\inf \Big\{ \F^M(v, \Omega) \,|\, v \in \B_M\W_\varphi^s(\Omega) \Big\} =: m.
$$
Clearly,~$\F^M(u^{(k)},\Omega) \le 2m$ for~$k$ large enough. Now, since~$\| u^{(k)} \|_{L^\infty(\Omega)} \le M$, we know that~$\Nl^M(u^{(k)}, \Omega) \ge 0$---recall definitions~\eqref{CH:4:NMldef} and~\eqref{CH:4:Gbardef}---and therefore~$\A(u^{(k)}, \Omega) \le 2 m$ for~$k$ large. In light of Lemma~\ref{CH:4:complem}, we then deduce that~$\{ u^{(k)} \}$ converges (up to a subsequence) to a function~$u_M \in \B_M\W^s_\varphi(\Omega)$ a.e.~in~$\R^n$. Identity~\eqref{CH:4:uMmin} follows by applying Lemma~\ref{CH:4:semicontlem}.
\end{proof}

We briefly mention here that if for some~$M_0>0$ we have~$\|u_{M_0}\|_{L^\infty(\Omega)}<M_0$, then---as a
consequence of the strict convexity of~$\F^M$---we obtain that~$u_M=u_{M_0}$
for every~$M\ge M_0$. It is readily seen that this implies that the function~$u_{M_0}$
minimizes~$\F$ in~$\W^s_\varphi(\Omega)$.
Therefore, in order to guarantee the existence of a minimizer,
it is enough to prove an a priori $L^\infty$ estimate. Depending on the exterior data, this is indeed possible---see
Theorem~\ref{CH:4:minareboundedthm} and Section~\ref{CH:4:Bded_proofs_Section}.

We will not pursue this strategy here, but we will exploit it to prove the existence of a solution to the obstacle problem. For more details we thus refer to the proof of Theorem~\ref{CH:4:Dirichlet_obstacle}.

Instead, we now prove the following a priori estimate on the~$W^{s,1}$ norm.


%

\begin{prop} \label{CH:4:Ws1prop}
Let~$n \ge 1$,~$s \in (0, 1)$,~$M \ge 0$, and~$\Omega \subseteq \R^n$ be a bounded open set with Lipschitz boundary.
Let~$\varphi:\Co\Omega\to\R$
with~$\Tail_s(\varphi, \Omega_{\Theta \diam(\Omega)} \setminus \Omega;\,\cdot\,) \in L^1(\Omega)$.
If~$u\in\W^s_\varphi(\Omega)$ is such that
\bgs{
\F^M(u,\Omega)\le\F^M(v,\Omega)\quad\mbox{for every }v\in\B_M\W^s_\varphi(\Omega),
}
then
\bgs{
\diam(\Omega)^{- s} \| u \|_{L^1(\Omega)} + [u]_{W^{s, 1}(\Omega)}
\le C \left(  \left\| \Tail_s(\varphi, \Omega_{\Theta \diam(\Omega)} \setminus \Omega;\,\cdot\,) \right\|_{L^1(\Omega)}
+ \diam(\Omega)^{1 - s} |\Omega| \right),
}
for two constants~$\Theta, C > 1$, depending only on~$n$,~$s$ and~$g$.
\end{prop}

We observe that Proposition~\ref{CH:4:Ws1prop} applies in particular to the minimizers~$u_M$,
but we stress that, in general, in the hypothesis we are not assuming~$u$ to be bounded.

\begin{proof}[Proof of Proposition~\ref{CH:4:Ws1prop}]
We use the function~$v := \chi_{\Co \Omega} u \in \B_M \W_\varphi^s(\Omega)$ as a competitor for~$u$. We get
\begin{equation} \label{CH:4:uMlev}
0 \le \F^M(v, \Omega) - \F^M(u, \Omega) = - \A(u, \Omega) + 2 \int_{\Omega} \int_{\Co \Omega} \frac{\I(x, y)}{|x - y|^{n - 1 + s}} \, dx\,dy,
\end{equation}
with
$$
\I(x, y) := \G \left( \frac{v(x) - v(y)}{|x - y|} \right) - \G \left( \frac{u(x) - u(y)}{|x - y|} \right).
$$
Write~$d := \diam(\Omega)$. On the one hand, by Lemma~\ref{CH:4:Adomainlem},
\begin{equation} \label{CH:4:AsuM}
\A(u, \Omega) \ge \frac{c_\star}{2} [u]_{W^{s, 1}(\Omega)} 
-\frac{c_\star\,\Ha^{n-1}(\mathbb S^{n-1})}{2(1-s)} |\Omega| d^{1 - s},
\end{equation}
with~$c_\star>0$ as defined in~\eqref{CH:4:cstardef}. On the other hand, let~$R := \Theta d$, with~$\Theta \ge 1$ to be chosen later. Recalling the definition of~$v$ and taking advantage of point~$(b)$ of Lemma~\ref{CH:4:gsprop}, we obtain
$$
\I(x, y) \le \frac{\Lambda}{2} \frac{|u(y)|}{|x - y|} +\frac{c_\star}{2}- \frac{c_\star}{2} \frac{|u(x) - u(y)|}{|x - y|} \quad \mbox{for every } x \in \Omega, \, y \in \Omega_R \setminus \Omega
$$
and
$$
\I(x, y) \le \frac{\Lambda}{2} \frac{|u(x)|}{|x - y|} \quad \mbox{for every } x \in \Omega, \, y \in \Co \Omega_R.
$$
Hence, exploiting Lemma~\ref{CH:4:dumb_kernel_lemma} and observing that~$c_\star\le\Lambda$, we get
\begin{align*}
2\int_{\Omega} \int_{\Co \Omega} \frac{\I(x, y)}{|x - y|^{n - 1 + s}} \, dx\,dy &\le \int_{\Omega} \left( \int_{\Omega_R \setminus \Omega} \frac{\Lambda |u(y)| - c_\star |u(x) - u(y)|}{|x - y|^{n + s}} \, dy \right) dx \\
&\quad+ c_\star \int_{\Omega} \int_{\Omega_R \setminus \Omega} \frac{dx\,dy}{|x - y|^{n - 1 + s}}\\
&\quad+ \Lambda \int_{\Omega} |u(x)| \left( \int_{\Co \Omega_R} \frac{dy}{|x - y|^{n + s}} \right) dx  \\
& \le \Lambda\, \bigg( \! \left\| \Tail_s(u, \Omega_{\Theta d} \setminus \Omega; \,\cdot\,) \right\|_{L^1(\Omega)}
+\frac{\Ha^{n-1}(\mathbb S^{n-1})}{1-s} \Theta^{1 - s} d^{1 - s} |\Omega| \\
& \quad + \Theta^{-s} d^{-s} \| u \|_{L^1(\Omega)} \bigg) - c_\star \int_{\Omega} \int_{\Omega_{\Theta d} \setminus \Omega} \frac{|u(x) - u(y)|}{|x - y|^{n + s}} \, dx\,dy.
\end{align*}
Putting together this estimate with~\eqref{CH:4:uMlev} and~\eqref{CH:4:AsuM}, and recalling that~$\Theta\ge1$, we find that
\begin{equation} \label{CH:4:semiest}
\begin{aligned}
\int_{\Omega} \int_{\Omega_{\Theta d}} \frac{|u(x) - u(y)|}{|x - y|^{n + s}} \, dx\,dy & \le \frac{\Lambda}{c_\star} \, \Bigg( \! \left\| \Tail_s(u, \Omega_{\Theta d} \setminus \Omega; \,\cdot\,) \right\|_{L^1(\Omega)} \\
& \quad +2\frac{\Ha^{n-1}(\mathbb S^{n-1})}{1-s} \Theta^{1 - s} d^{1 - s} |\Omega| + \Theta^{-s} d^{-s} \| u \|_{L^1(\Omega)} \Bigg).
\end{aligned}
\end{equation}
Now we observe that
\eqlab{\label{CH:4:esti_diam_meas}
\diam(\Omega_d)=3d\qquad\mbox{and}\qquad|\Omega_d\setminus\Omega|\ge c_nd^n,
}
for some dimensional constant~$c_n>0$ depending only on~$n$. Indeed, the equality is an immediate consequence of the definition
of~$\Omega_d$, while the measure estimate follows by observing that if we take a point~$x_0\in\partial\Omega_{d/2}$,
then~$B_{d/2}(x_0)\subseteq\Omega_d\setminus\Omega$, and hence
\bgs{
|\Omega_d\setminus\Omega|\ge|B_{d/2}(x_0)|=\frac{|B_1|}{2^n}\,d^n.
}
Since~$v = 0$ in~$\Omega$ and~$v = u$ outside of~$\Omega$, using Lemma~\ref{CH:4:FPI} and exploiting~\eqref{CH:4:esti_diam_meas}, we may now estimate
\begin{align*}
\| u \|_{L^1(\Omega)} & = \| u - v \|_{L^1(\Omega)} \le \frac{\diam(\Omega_d)^{n + s}}{|\Omega_d \setminus \Omega|} \int_{\Omega} |u(x)| \left( \int_{\Omega_d \setminus \Omega} \frac{dy}{|x - y|^{n + s}} \right) dx \\
& \le C d^s \left( \int_{\Omega} \int_{\Omega_{d}\setminus \Omega} \frac{|u(x) - u(y)|}{|x - y|^{n + s}} \, dx\,dy + \left\| \Tail_s(u, \Omega_{\Theta d} \setminus \Omega;\,\cdot\,) \right\|_{L^1(\Omega)} \right),
\end{align*}
with~$C>0$ depending only on~$n$ and~$s$.
Using this estimate together with~\eqref{CH:4:semiest} and recalling that~$\Theta \ge 1$, we get
$$
\| u \|_{L^1(\Omega)} \le C \left(  d^s\left\| \Tail_s(u, \Omega_{\Theta d} \setminus \Omega; \,\cdot\,) \right\|_{L^1(\Omega)} + \Theta^{1 - s} d |\Omega| + \Theta^{-s} \| u \|_{L^1(\Omega)} \right),
$$
with~$C>0$ depending only on~$n$,~$s$ and~$g$.
By taking~$\Theta$ sufficiently large (in dependence of~$n$,~$s$ and~$g$ only), we can reabsorb the~$L^1$ norm of~$u$ on the left-hand side and obtain that
$$
\| u \|_{L^1(\Omega)} \le C \left(  d^s \left\| \Tail_s(u, \Omega_{\Theta d} \setminus \Omega;\,\cdot\,) \right\|_{L^1(\Omega)} + d |\Omega| \right).
$$
The conclusion follows by combining this estimate with~\eqref{CH:4:semiest}.
\end{proof}

As shown in the following Lemma, the integrability of the truncated tail is equivalent to~$L^1$ integrability
plus weighted integrability arbitrarily close to the boundary of the domain. 

\begin{lemma}\label{CH:4:tail_equiv_cond_Lemma}
Let~$n \ge 1$,~$s \in (0, 1)$,~$\Omega \Subset\Op\subseteq \R^n$ two bounded open sets, such that~$\Omega$ has Lipschitz boundary, and~$\varphi:\Co\Omega\to\R$.
Then,~$\Tail_s(\varphi, \Op \setminus \Omega;\,\cdot\,) \in L^1(\Omega)$
if and only if~$\varphi\in L^1(\Op\setminus\Omega)$
and~$\Tail_s(\varphi, \Omega_r \setminus \Omega;\,\cdot\,) \in L^1(\Omega\setminus\Omega_{-r})$, for some small~$r>0$.

Moreover, suppose that~$\varphi\in L^1(\Op\setminus\Omega)$ and let~$r>0$ be small. Then:
\begin{itemize}
\item[(i)] if~$\varphi\in W^{s,1}(\Omega_r\setminus\Omega)$,
then~$\Tail_s(\varphi, \Op \setminus \Omega;\,\cdot\,) \in L^1(\Omega)$;
\item[(ii)] if~$\varphi\in L^\infty(\Omega_r\setminus\Omega)$,
then~$\Tail_\sigma(\varphi, \Op \setminus \Omega;\,\cdot\,) \in L^1(\Omega)$,
for every~$\sigma\in(0,1)$.
\end{itemize}
\end{lemma}

\begin{proof}
To begin, let~$d:=\diam(\Op)$ and notice that
\bgs{
\frac{1}{|x-y|^{n+s}}\ge\frac{1}{d^{n+s}}\quad\mbox{for every }x\in\Omega\mbox{ and }y\in\Op\setminus\Omega.
}
Hence
\bgs{
\|\varphi\|_{L^1(\Op\setminus\Omega)}\le \frac{d^{n+s}}{|\Omega|}\left\|\Tail_s(\varphi, \Op \setminus \Omega;\,\cdot\,)\right\|_{L^1(\Omega)}.
}
Moreover, we clearly have
\bgs{
\left\|\Tail_s(\varphi, \Omega_r \setminus \Omega;\,\cdot\,)\right\|_{L^1(\Omega\setminus\Omega_{-r})}\le
\left\|\Tail_s(\varphi, \Op \setminus \Omega;\,\cdot\,)\right\|_{L^1(\Omega)},
}
for every small~$r>0$.

Now suppose that~$\varphi\in L^1(\Op\setminus\Omega)$ and let~$r>0$ be small.

If~$\Tail_s(\varphi, \Omega_r \setminus \Omega;\,\cdot\,) \in L^1(\Omega\setminus\Omega_{-r})$,
then~$\Tail_s(\varphi, \Op \setminus \Omega;\,\cdot\,) \in L^1(\Omega)$. Indeed, since
\bgs{
|x-y|\ge r\quad\mbox{for every }x\in\Omega\mbox{ and }y\in\Op\setminus\Omega_r,
}
we have
\eqlab{\label{CH:4:alatriste_e_le_maledette_code}
\left\|\Tail_s(\varphi, \Op \setminus \Omega_r;\,\cdot\,)\right\|_{L^1(\Omega)}
\le\frac{|\Omega|}{r^{n+s}}\|\varphi\|_{L^1(\Op\setminus\Omega_r)}.
}
Similarly,
\bgs{
\left\|\Tail_s(\varphi, \Omega_r \setminus \Omega;\,\cdot\,)\right\|_{L^1(\Omega_{-r})}
\le\frac{|\Omega_{-r}|}{r^{n+s}}\|\varphi\|_{L^1(\Omega_r\setminus\Omega)}
\le\frac{|\Omega|}{r^{n+s}}\|\varphi\|_{L^1(\Omega_r\setminus\Omega)}.
}
Therefore,
\bgs{
\left\|\Tail_s(\varphi, \Op \setminus \Omega;\,\cdot\,)\right\|_{L^1(\Omega)}&
=\left\|\Tail_s(\varphi, \Op \setminus \Omega_r;\,\cdot\,)\right\|_{L^1(\Omega)}
+\left\|\Tail_s(\varphi, \Omega_r \setminus \Omega;\,\cdot\,)\right\|_{L^1(\Omega_{-r})}\\
&
\qquad+\left\|\Tail_s(\varphi, \Omega_r\setminus\Omega;\,\cdot\,)\right\|_{L^1(\Omega\setminus\Omega_{-r})}\\
&
\le\frac{|\Omega|}{r^{n+s}}\|\varphi\|_{L^1(\Op\setminus\Omega)}
+\left\|\Tail_s(\varphi, \Omega_r\setminus\Omega;\,\cdot\,)\right\|_{L^1(\Omega\setminus\Omega_{-r})}.
}

If~$\varphi\in W^{s,1}(\Omega_r\setminus\Omega)$, then---since for small~$r>0$ the open set~$\Omega_r$ has Lipschitz boundary---by Corollary~\ref{CH:4:FHI_corollary} we obtain
\bgs{
\left\|\Tail_s(\varphi, \Omega_r \setminus \Omega;\,\cdot\,)\right\|_{L^1(\Omega)}
\le C(n,\Omega_r\setminus\Omega,s)\|\varphi\|_{W^{s,1}(\Omega_r\setminus\Omega)}.
}
Hence, recalling~\eqref{CH:4:alatriste_e_le_maledette_code}, we have
\bgs{
\left\|\Tail_s(\varphi, \Op \setminus \Omega;\,\cdot\,)\right\|_{L^1(\Omega)}
\le\frac{|\Omega|}{r^{n+s}}\|\varphi\|_{L^1(\Op\setminus\Omega_r)}+C\|\varphi\|_{W^{s,1}(\Omega_r\setminus\Omega)}.
}

If~$\varphi\in L^\infty(\Omega_r\setminus\Omega)$, then
\bgs{
\left\|\Tail_\sigma(\varphi, \Omega_r \setminus \Omega;\,\cdot\,)\right\|_{L^1(\Omega)}
\le\|\varphi\|_{L^\infty(\Omega_r\setminus\Omega)}\Per_\sigma(\Omega),
}
for every~$\sigma\in(0,1)$. Thus, we obtain point~$(ii)$ by exploiting~\eqref{CH:4:alatriste_e_le_maledette_code} again.
This concludes the proof of the Lemma.
\end{proof}

We conclude this Section by getting back to the geometric framework~$g=g_s$. We exploit Theorem~\ref{CH:4:Persdecreases}
in order to prove that the unique $s$-minimal set
in~$\Omega^M$ with respect to the exterior data~$\Sg(\varphi)\setminus\Omega^M$ is the subgraph of
the minimizer~$u_M$.

\begin{prop}\label{CH:4:rearrangement_truncated_Prop}
Let~$n \ge 1$,~$s \in (0, 1)$,~$M>0$, and~$\Omega \subseteq \R^n$ a bounded open set with Lipschitz boundary.
Let~$\varphi:\R^n\to\R$, such that~$\varphi=0$ a.e. in~$\Omega$ and let~$u_M$ be
the minimizer of~$\F^M_s(\,\cdot\,,\Omega)$ within~$\B_M\W^s_\varphi(\Omega)$.
Then,~$\Sg(u_M)$ is the unique set which is $s$-minimal in~$\Omega^M$
with respect to the exterior data~$\Sg(\varphi)\setminus\Omega^M$.
\end{prop}

\begin{proof}
Let~$E\subseteq\R^{n+1}$ be $s$-minimal in~$\Omega^M$,
with respect to the exterior data~$\Sg(\varphi)\setminus\Omega^M$---we
know that such a set exists by~\cite[Theorem~3.2]{CRS10}. Let~$w_E$
be the function defined in~\eqref{CH:4:rearr_func_def} and notice that the set~$E$
satisfies the hypothesis of Theorem~\ref{CH:4:Persdecreases}.
Hence,
\eqlab{\label{CH:4:saccodipulcimaledetto}
\Per_s(\Sg(w_E),\Omega^M)\le\Per_s(E,\Omega^M).
}
As a consequence, we conclude that $E=\Sg(w_E)$, since otherwise the inequality~\eqref{CH:4:saccodipulcimaledetto}
would be strict, thus contradicting the minimality of~$E$.
Recalling~\eqref{CH:4:u_iff_Sgu}, we have in particular that~$w_E\in\B_M\W^s_\varphi(\Omega)$.
Then, by identity~\eqref{CH:4:per_of_subgraph} and exploiting both the minimality of~$u_M$
and of~$E$,
we find that
\bgs{
0\ge\F^M_s(u_M,\Omega^M)-\F^M_s(w_E,\Omega^M)=
\Per_s(\Sg(u_M),\Omega^M)-\Per_s(\Sg(w_E),\Omega^M)\ge0.
}
Since~$u_M$ is the unique minimizer of~$\F^M_s(\,\cdot\,,\Omega)$ within~$\B_M\W^s_\varphi(\Omega)$,
this implies that~$u_M=w_E$, concluding the proof.
\end{proof}

\subsection{Proof of Theorem~\ref{CH:4:Dirichlet}}\label{CH:4:Proof_of_Dir_Sec}

%

Proposition~\ref{CH:4:ertyui} shows that, for each~$M > 0$, there exists a unique minimizer~$u_M$ of~$\F^M(\,\cdot\,, \Omega)$ within the space~$\B_M\W_\varphi^s(\Omega)$. To establish the existence of a minimizer of~$\F$, we now need~$u_M$ to converge as~$M\to\infty$. This is achieved through the uniform~$W^{s, 1}$ estimate of Proposition~\ref{CH:4:Ws1prop}, at the price of assuming some (weighted) integrability on the exterior datum in a sufficiently large neighborhood of~$\Omega$.
The minimality of the limit function~$u$ is then obtained as a consequence of the minimality of the
functions~$u_M$.

\begin{proof}[Proof of Theorem \ref{CH:4:Dirichlet}]
Let~$\Theta > 1$ be the constant given by Proposition~\ref{CH:4:Ws1prop}. For any~$M > 0$, the minimizer~$u_M$ satisfies the hypotheses of Proposition~\ref{CH:4:Ws1prop}. Therefore,
\begin{equation} \label{CH:4:Ws1estforuM}
\| u_M \|_{W^{s, 1}(\Omega)} \le C \left( \left\| \Tail_s(\varphi, \Omega_{\Theta \diam(\Omega)} \setminus \Omega;\,\cdot\,) \right\|_{L^1(\Omega)} + 1 \right),
\end{equation}
for some constant~$C > 0$ depending only on~$n$,~$s$,~$g$ and~$\Omega$, and, in particular, independent of~$M$.
By the compact fractional Sobolev embedding (see, e.g.,~\cite[Theorem~7.1]{HitGuide}),
we conclude that there exists a function~$u \in \W^s_\varphi(\Omega)$ to which~$\{ u_{M_j} \}$
converges in~$L^1(\Omega)$ and~a.e.~in~$\Omega$, for some diverging sequence~$\{ M_j \}_{j \in \N}$.
Letting~$M = M_j \rightarrow +\infty$ in~\eqref{CH:4:Ws1estforuM}, by Fatou's Lemma
we see that~$u$ satisfies~\eqref{CH:4:Ws1estformin}.
We are therefore left to show that~$u$ is a minimizer for~$\F$ in~$\W^s_\varphi(\Omega)$.

Take~$v \in \B\W^s_\varphi(\Omega)$. Then, for~$j$ large enough we have $M_j\ge\| v \|_{L^\infty(\Omega)}$,
and hence, by the minimality of~$u_{M_j}$ we get~$\F^{M_j}(u_{M_j}, \Omega) \le \F^{M_j}(v, \Omega)$. That is,
\begin{align*}
0 & \ge \A(u_{M_j}) + 2 \int_{\Omega} \left\{ \int_{\Omega_R \setminus \Omega} \G \left( \frac{u_{M_j}(x) - \varphi(y)}{|x-y|} \right) \frac{dy}{|x-y|^{n - 1 + s}} \right\} dx \\
& \quad - \A(v) - 2 \int_{\Omega} \left\{ \int_{\Omega_R \setminus \Omega} \G \left( \frac{v(x) - \varphi(y)}{|x-y|} \right) \frac{dy}{|x-y|^{n - 1+ s}} \right\} dx \\
& \quad + 2 \int_{\Omega} \left\{ \int_{\Co \Omega_R} \left[ \G \left( \frac{u_{M_j}(x) - \varphi(y)}{|x-y|} \right) - \G \left( \frac{v(x) - \varphi(y)}{|x-y|} \right) \right] \frac{dy}{|x-y|^{n - 1 + s}} \right\} dx,
\end{align*}
for any fixed~$R \in (0, \Theta \diam(\Omega)]$.
We now claim that letting~$j \rightarrow +\infty$ in the above formula,
we obtain the same inequality with~$u_{M_j}$ replaced by~$u$.

Indeed, the quantities on the first line can be dealt with by using Fatou's lemma.
Moreover, the Lipschitz character of~$\G$---see~\eqref{CH:4:Lip_Gcal}---and the fact that~$u_{M_j} \rightarrow u$ in~$L^1(\Omega)$ ensure that
\bgs{
 \int_{\Omega} &\left\{ \int_{\Co \Omega_R} \left| \G \left( \frac{u_{M_j}(x) - \varphi(y)}{|x-y|} \right) - \G \left( \frac{u(x) - \varphi(y)}{|x-y|} \right) \right| \frac{dy}{|x-y|^{n - 1 + s}} \right\} dx\\
&
\qquad\le\frac{\Lambda}{2}\int_\Omega|u_{M_j}(x)-u(x)|\left\{\int_{\Co\Omega_R}\frac{dy}{|x-y|^{n+s}}\right\}dx\\
&
\qquad\le CR^{-s}\|u_{M_j}-u\|_{L^1(\Omega)}\xrightarrow{j\to\infty}0.
}
Hence, the third line passes to the limit as well.
All in all, we have proved that~$u$ minimizes~$\F$ in~$\B\W_\varphi^s(\Omega)$. The minimality of~$u$
within the larger class~$\W^s_\varphi(\Omega)$ follows from the density of~$L^\infty(\Omega)\cap W^{s,1}(\Omega)$
in~$W^{s,1}(\Omega)$ and Lemma~\ref{CH:4:tartariccio}---see point~$(iv)$ of Remark~\ref{CH:4:tartapower_Remark}.
To conclude, the uniqueness of the minimizer follows by point~$(iii)$ of Remark~\ref{CH:4:tartapower_Remark}.
\end{proof}

\subsection{Boundedness results}\label{CH:4:Bded_proofs_Section}

The purpose of this section consists in proving that minimizers of~$\F$ are always locally bounded
and that they are globally bounded if the exterior data is bounded near the boundary of
the domain~$\Omega$.

More precisely, by exploiting a Stampacchia-type argument, we prove the following result:

\begin{prop} \label{CH:4:Linftylocprop}
Let~$n \ge 1$,~$s \in (0, 1)$,~$R > 0$, and~$u\in\W^s(B_{2R})$ be a minimizer of~$\F$ in~$B_{2 R}$. Then,
$$
\sup_{B_R} u \le C \left( R + \dashint_{B_{2 R}} u_+(x) \, dx \right),
$$
for some constant~$C > 0$ depending only on~$n$,~$s$ and~$g$.
\end{prop}

Clearly, Proposition~\ref{CH:4:Linftylocprop} implies that if~$u\in\W^s(\Omega)$ is a minimizer of~$\F$ in~$\Omega$,
then~$u\in L^\infty_{\loc}(\Omega)$.
Since the proof is rather lenghty and technical, we postpone it to Section~\ref{CH:4:Linftylocprop_proof}.

Moreover, we prove that a minimizer~$u$ of~$\F$ in~$\Omega$ belongs to~$L^\infty(\Omega)$, provided it is bounded,
outside~$\Omega$, in a sufficiently large neighborhood of~$\Omega$. Furthermore, we obtain an apriori
estimate on the~$L^\infty(\Omega)$ norm of~$u$ purely in terms of the exterior data.
That is, we show the validity of Theorem~\ref{CH:4:minareboundedthm} of the Introduction.

We establish this result by showing that, given any function~$u: \R^n \to \R$, bounded in~$B_R \setminus \Omega$ for some large~$R > 0$, the value~$\F^M(u, \Omega)$ decreases when~$u$ is truncated at a high enough level. This last statement can be made precise as follows.

For~$N \ge0$, we define
$$
u^{(N)} := \begin{cases}
\min \{ u, N \} & \quad \mbox{in } \Omega,\\
u & \quad \mbox{in } \Co \Omega.
\end{cases}
$$
Then, we have the following result.

\begin{prop} \label{CH:4:truncdecreaseprop}
Let~$n \ge 1$,~$s \in (0, 1)$, and~$M \ge 0$. Let~$\Omega \subseteq \R^n$ be a bounded open set with Lipschitz boundary and let~$R_0 > 0$ be such that~$\Omega \subseteq B_{R_0}$. Then, there exists a large constant~$\Theta > 1$, depending only on~$n$,~$s$ and~$g$, such that for every function~$u: \R^n \to \R$ bounded from above in~${B_{\Theta R_0} \setminus \Omega}$, it holds
\begin{equation} \label{CH:4:AandNdecrease}
\A(u^{(N)}, \Omega) \le \A(u, \Omega) \quad \mbox{and} \quad \Nl^M(u^{(N)}, \Omega) \le \Nl^M(u, \Omega)
\end{equation}
for every
\begin{equation} \label{CH:4:NgeR0+sup}
N \ge R_0 + \sup_{B_{\Theta R_0} \setminus \Omega} u.
\end{equation}
In particular,
$$
\F^M(u^{(N)}, \Omega) \le \F^M(u, \Omega)
$$
for every~$N$ satisfying~\eqref{CH:4:NgeR0+sup}.
\end{prop}

We observe that Proposition~\ref{CH:4:truncdecreaseprop}
directly implies Theorem~\ref{CH:4:minareboundedthm}, thanks to the uniqueness of the minimizer,
which is a consequence of the strict convexity of~$\F^M$---see point~$(iii)$ of Remark~\ref{CH:4:tartapower_Remark}.

By exploiting the interior local boundedness and by appropriately modifying the proof of Proposition~\ref{CH:4:truncdecreaseprop},
we are able to prove that, in order to ensure the global boundedness of
a minimizer of~$\F$, it is actually enough that~$u$ be bounded outside the domain~$\Omega$
in an arbitrarily small neighborhood of the boundary.
However, we remark that in this case, in general, we do not have a clean a priori bound on the~$L^\infty$ norm.

We first recall that if~$\Omega\subseteq\R^n$ is a bounded open set with~$C^2$ boundary,
then there exists~$r_0(\Omega)>0$ such that~$\Omega$ satisfies a uniform strict interior and strict
exterior ball condition of radius~$2r_0$. Then, if~$\bar{d}_\Omega$ denotes the signed distance function from~$\partial\Omega$, negative inside~$\Omega$,
we have that~$\bar{d}_\Omega\in C^2\big(N_{2r_0}(\partial\Omega)\big)$,
with
\bgs{
N_\varrho(\partial\Omega):=\left\{x\in\R^n\,|\,d(x,\partial\Omega)<\varrho\right\}=\left\{|\bar{d}_\Omega|<\varrho\right\}
\quad\forall\,\varrho>0.
}
For the details, we refer to Appendix \ref{CH:3:A2}---see in particular Remark \ref{CH:3:ext_unif_omega}.

The precise result is the following: 

\begin{theorem}\label{CH:4:Bdary_Bdedness_Thm}
Let~$n \ge 1$,~$s \in (0, 1)$ and~$\Omega\subseteq\R^n$ be a bounded open set with~$C^2$ boundary.
If~$u\in\W^s(\Omega)$ is a minimizer of~$\F$ in~$\Omega$ and~$u\in L^\infty(\Omega_d\setminus\Omega)$,
for some~$d\in(0,r_0)$,
then $u\in L^\infty(\Omega)$, with
\bgs{
\|u\|_{L^\infty(\Omega\setminus\Omega_{-\theta d})}
\le d+\max\left\{\|u\|_{L^\infty(\Omega_{-\theta d})},\|u\|_{L^\infty(\Omega_d\setminus\Omega)}\right\},
}
where~$\theta=\theta(n,s,g)\in(0,1)$ is a small positive constant.
\end{theorem}

We observe that if we further assume
that~$\Tail_s(u, \Omega_{\Theta \diam(\Omega)} \setminus \Omega;\,\cdot\,) \in L^1(\Omega)$,
then, by exploiting both the apriori~$L^1$ estimate of Proposition~\ref{CH:4:Ws1prop},
the estimate on~$\|u\|_{L^\infty(\Omega_{-\theta d})}$ given by Proposition~\ref{CH:4:Linftylocprop}---together with
a covering argument---and the estimate provided by Theorem~\ref{CH:4:Bdary_Bdedness_Thm},
we can obtain an apriori estimate on~$\|u\|_{L^\infty(\Omega)}$ purely in terms of the
exterior data and of the geometry of~$\Omega$.

\smallskip

We now proceed with the proofs of the aforementioned results.

To prove Proposition~\ref{CH:4:truncdecreaseprop}, we will make use of a couple of simple lemmas. First, we have the following elementary result on convex functions.

\begin{lemma} \label{CH:4:A+Blem}
Let~$\phi: \R \to \R$ be a convex function. Then, for every~$A, B, C, D \in \R$ satisfying~$\min \{ C, D\} \le A, B \le \max \{ C, D\}$ and~$A + B = C + D$, it holds
$$
\phi(A) + \phi(B) \le \phi(C) + \phi(D).
$$
\end{lemma}
\begin{proof}
Without loss of generality, we may suppose that~$A \le B$ and~$C \le D$. Since we have that~$C \le A \le B \le D$, there exist two values~$\lambda, \mu \in [0, 1]$ such that
$$
A = \lambda C + (1 - \lambda) D \quad \mbox{and} \quad B = \mu C + (1 - \mu) D.
$$
In view of the convexity of~$\phi$, it holds
\begin{equation} \label{CH:4:A+Blemtech}
\begin{aligned}
\phi(A) + \phi(B) & = \phi(\lambda C + (1 - \lambda) D) + \phi(\mu C + (1 - \mu) D) \\
& \le \lambda \phi(C) + (1 - \lambda) \phi(D) + \mu \phi(C) + (1 - \mu) \phi(D) \\
& = \left( \lambda + \mu \right) \phi(C) + \left( 2 - \lambda - \mu \right) \phi(D).
\end{aligned}
\end{equation}
By taking advantage of the fact that~$A + B = C + D$, we now observe that
$$
\lambda C + (1 - \lambda) D + \mu C + (1 - \mu) D = C + D,
$$
or, equivalently,
$$
(1 - \lambda - \mu) (C - D) = 0.
$$
Consequently, either~$C = D$ or~$\lambda + \mu = 1$ (or both). In any case, we conclude that the right-hand side of~\eqref{CH:4:A+Blemtech} is equal to~$\phi(C) + \phi(D)$, and from this the thesis follows.
\end{proof}

We use Lemma~\ref{CH:4:A+Blem} to obtain the following inequality for rather general convex functionals. In our later applications, we will simply take~$F(U; x, y) = \G(U/|x - y|)$.

\begin{lemma} \label{CH:4:Fminmaxlem}
Let~$F: \R \times \R^n \times \R^n \to \R$ be a measurable function, convex with respect to the first variable, i.e.~satisfying
\begin{equation} \label{CH:4:Fconvex}
F(\lambda u + (1 - \lambda) v; x, y) \le \lambda F(u; x, y) + (1 - \lambda) F(v; x, y)
\end{equation}
for every~$\lambda \in (0, 1)$,~$u, v \in \R$, and for~a.e.~$x, y \in \R^n$. Given a measurable set~$\U \subseteq \R^n \times \R^n$, consider the functional~$\F$ defined by
$$
\F(w) := \iint_{\U} F(u(x) - u(y); x, y) \, dx\,dy
$$
for every~$w: \R^n \to \R$. Then, for every~$u, v: \R^n \to \R$, it holds
\begin{equation} \label{CH:4:minmaxleuv}
\F(\min\{u, v\}) + \F(\max \{u, v\}) \le \F(u) + \F(v).
\end{equation}
\end{lemma}
\begin{proof}
For fixed~$(x, y) \in \U$, we write
\begin{align*}
A & := m(x) - m(y), \quad B := M(x) - M(y), \quad C := u(x) - u(y), \quad D := v(x) - v(y),
\end{align*}
and
$$
\phi(t) = \phi_{x, y}(t) := F(t; x, y) \quad \mbox{for every } t \in \R.
$$
Thanks to~\eqref{CH:4:Fconvex}, the function~$\varphi$ is convex. Also, we claim that
\begin{equation} \label{CH:4:CleA,BleD}
\min \{ C, D \} \le A, B \le \max \{ C, D \}
\end{equation}
and
\begin{equation} \label{CH:4:A+B=C+D}
A + B = C + D.
\end{equation}

Indeed, identity~\eqref{CH:4:A+B=C+D} is immediate since~$m + M \equiv u + v$. The inequalities in~\eqref{CH:4:CleA,BleD} are also obvious if~$u(x) \le v(x)$ and~$u(y) \le v(y)$ or if~$u(x) > v(x)$ and~$u(y) > v(y)$. On the other hand, when for example~$u(x) \le v(x)$ and~$u(y) > v(y)$, we have
$$
A = u(x) - v(y) \quad \mbox{and} \quad B = v(x) - u(y).
$$
Accordingly,
$$
C = u(x) - u(y) < u(x) - v(y) = A = u(x) - v(y) \le v(x) - v(y) = D
$$
and
$$
C = u(x) - u(y) \le v(x) - u(y) = B = v(x) - u(y) < v(x) - v(y) = D.
$$
Hence,~\eqref{CH:4:CleA,BleD} is proved in this case. Arguing analogously, one can check that~\eqref{CH:4:CleA,BleD} also holds when~$u(x) > v(x)$ and~$u(y) \le v(y)$.

Thanks to~\eqref{CH:4:CleA,BleD} and~\eqref{CH:4:A+B=C+D}, we may apply Lemma~\ref{CH:4:A+Blem} and deduce that
$$
\phi(A) + \phi(B) \le \phi(C) + \phi(D).
$$
That is,
\begin{align*}
F(m(x) - m(y); x, y) & + F(M(x) - M(y); x, y) \\
& \le F(u(x) - u(y); x, y) + F(v(x) - v(y); x, y).
\end{align*}
Inequality~\eqref{CH:4:minmaxleuv} then plainly follows by integrating the last formula in~$x$ and~$y$.
\end{proof}

With the aid of this last result, we can proceed to check the validity of Proposition~\ref{CH:4:truncdecreaseprop}.

\begin{proof}[Proof of Proposition~\ref{CH:4:truncdecreaseprop}]
Write~$v := u^{(N)}$ and~$R := \Theta R_0$, with~$\Theta \ge 2$ to be chosen later sufficiently large,
in dependence of~$n$,~$s$ and~$g$ only.
From Lemma~\ref{CH:4:Fminmaxlem}, it clearly follows that~$\A(v, \Omega) \le \A(u, \Omega)$.
Hence, we can focus on the inequality for the nonlocal part~$\Nl^M$.

Thanks to representation~\eqref{CH:4:nonlocal_explicit}, we have
$$
\Nl^M(v, \Omega) - \Nl^M(u, \Omega) = 2 \int_{\Omega} \int_{\Co \Omega} \left[ \G \left( \frac{v(x) - v(y)}{|x - y|} \right) - \G \left( \frac{u(x) - u(y)}{|x - y|} \right) \right] \frac{dx\,dy}{|x - y|^{n - 1 + s}}.
$$
Setting~$\Omega_+ := \{ x \in \Omega \,|\, u(x) > N \}$ and writing~$\Co \Omega = A_1 \cup A_2$, with~$A_1 := B_{R} \setminus \Omega$ and~$A_2 := \Co B_{R}$, we infer from the above identity that the second inequality in~\eqref{CH:4:AandNdecrease} is equivalent to
\begin{equation} \label{CH:4:Nlorderequiv}
\alpha_1 + \alpha_2 \le 0,
\end{equation}
where we set
\begin{align*}
\alpha_i := \int_{\Omega_+} \left\{ \int_{A_i} \left[ \G \left( \frac{N - u(y)}{|x - y|} \right) - \G \left( \frac{u(x) - u(y)}{|x - y|} \right) \right] \frac{dy}{|x - y|^{n - 1 + s}} \right\} dx,
\end{align*}
for~$i = 1, 2$.
 
First, we establish a (negative) upper bound for~$\alpha_1$. Let~$x \in \Omega_+$ and~$y \in A_1$. Since, by hypothesis~\eqref{CH:4:NgeR0+sup},~$u(y) \le N < u(x)$ and~$G$ is increasing, we have
$$
\G \left( \frac{N - u(y)}{|x - y|} \right) - \G \left( \frac{u(x) - u(y)}{|x - y|} \right) = \int_{\frac{u(x) - u(y)}{|x - y|}}^{\frac{N - u(y)}{|x - y|}} G(t) \, dt \le - G\left( \frac{N - u(y)}{|x - y|} \right) \frac{u(x) - N}{|x - y|},
$$
and consequently
$$
\alpha_1 \le - \int_{\Omega} \left( u(x) - N \right)_+ \left[ \int_{A_1} G \left( \frac{N - u(y)}{|x - y|} \right) \frac{dy}{|x - y|^{n + s}} \right] dx.
$$
In view of the fact that~$B_{2 R_0} \setminus B_{R_0} \subseteq A_1$ (as~$R \ge 2 R_0$) and, again,~\eqref{CH:4:NgeR0+sup} and the monotonicity of~$G$, we estimate
$$
\int_{A_1} G \left( \frac{N - u(y)}{|x - y|} \right) \frac{dy}{|x - y|^{n + s}} \ge G \left( \frac{R_0}{R_0 + R_0} \right) \frac{|B_{2 R_0} \setminus B_{R_0}|}{(R_0 + R_0)^{n + s}} \ge \frac{c_1}{R_0^s}
$$
for every~$x \in \Omega$ and for some constant~$c_1 > 0$ depending only on~$n$,~$s$ and~$g$. Accordingly,
\begin{equation} \label{CH:4:negcontr}
\alpha_1 \le - \frac{c_1}{R_0^s} \int_{\Omega} \left( u(x) - N \right)_+ dx.
\end{equation}

On the other hand, to control~$\alpha_2$ we simply use that~$\G$ is a globally Lipschitz function---see~\eqref{CH:4:Lip_Gcal}---and compute
\begin{align*}
\alpha_2 & \le\frac{\Lambda}{2} \int_{\Omega} \left( u(x) - N \right)_+ \left( \int_{\R^n \setminus B_R} \frac{dy}{|x - y|^{n + s}} \right) dx \\
& \le \frac{\Lambda}{2}\int_{\Omega} \left( u(x) - N \right)_+ \left( \int_{\R^n \setminus B_{R/2}} \frac{dz}{|z|^{n + s}} \right) dx \le \frac{C_2}{R^s} \int_{\Omega} \left( u(x) - N \right)_+ dx,
\end{align*}
for some constant~$C_2 > 0$ depending only on~$n$,~$s$ and~$g$. Notice that to get the second inequality we changed variables and took advantage of the inclusion~$B_{R/2}(x) \subseteq B_R$, which holds for all~$x \in \Omega \subseteq B_{R_0}$ since~$R \ge 2 R_0$. Combining this last estimate with~\eqref{CH:4:negcontr}, we obtain
$$
\alpha_1 + \alpha_2 \le - \left( \frac{c_1}{R_0^s} - \frac{C_2}{R^s} \right) \int_{\Omega} \left( u(x) - N \right)_+ dx,
$$
and~\eqref{CH:4:Nlorderequiv} follows provided we take~$R \ge (C_2 / c_1)^{1/s} R_0$.
\end{proof}

A suitable modification of the proof
of Proposition~\ref{CH:4:truncdecreaseprop} allows us to obtain
Theorem~\ref{CH:4:Bdary_Bdedness_Thm}.

\begin{proof}[Proof of Theorem~\ref{CH:4:Bdary_Bdedness_Thm}]
We recall that
\bgs{
u^{(N)}:=\chi_\Omega\min\{u,N\}+(1-\chi_\Omega)u.
}
We consider~$u^{(N)}$ with
\eqlab{\label{CH:4:sonny_rollins}
N\ge d+\max\left\{\sup_{\Omega_{-\theta d}} u,\sup_{\Omega_d\setminus\Omega} u\right\},
}
where~$\theta\le1/4$ will be chosen suitably small later. We fix~$M\ge0$ and we prove that
\bgs{
\F^M(u^{(N)},\Omega)\le\F^M(u,\Omega).
}
We remark that the analogous estimate holds true when we cut~$u$ from below, inside~$\Omega$.
Hence, by the minimality of~$u$ and the uniqueness of the minimizer---see points~$(ii)$ and~$(iii)$
of Remark~\ref{CH:4:tartapower_Remark}---this
implies the claim of the Theorem.

By arguing as in the proof of Proposition~\ref{CH:4:truncdecreaseprop},
we are left to prove that
\bgs{
I:=\int_{\Omega_+}\bigg\{\int_{\Co\Omega}\bigg[\G\left(\frac{N-u(y)}{|x-y|}\right)-\G\left(\frac{u(x)-u(y)}{|x-y|}\right)\bigg]
\frac{dy}{|x-y|^{n-1+s}}\bigg\}dx\le0,
}
where~$\Omega_+:=\{\xi\in\Omega\,|\,u(\xi)>N\}$. It is important to observe that by~\eqref{CH:4:sonny_rollins}
we have
\bgs{
\Omega_+\subseteq\Omega\setminus\Omega_{-\theta d}=\{\xi\in\Omega\,|\,d(\xi,\partial\Omega)\le\theta d\}.
}
As a consequence, since~$\theta\le1/4$, for every~$x\in\Omega_+$ we can find a
point~$z_x\in\Omega_d\setminus\Omega$ such that
\eqlab{\label{CH:4:sonny_rollins2}
B_\frac{\theta d}{2}(z_x)\subseteq B_{3\theta d}(x)\setminus\Omega
\subseteq B_d(x)\setminus\Omega\subseteq\Omega_d\setminus\Omega.
}
This is a consequence of the uniform interior and exterior ball conditions satisfied by~$\Omega$. More precisely, we
observe that
\bgs{
p:=x-\bar{d}_\Omega(x)\nabla\bar{d}_\Omega(x)\in\partial\Omega,
}
is the unique closest point to~$x$. That is,~$p$ is the unique point on~$\partial\Omega$ such that~$|x-p|=d(x,\partial\Omega)$.
Then,~$\Omega$ has an exterior tangent ball of radius~$r_0$ at~$p$. Notice that the center of the ball is
obtained by moving in direction~$\nu_\Omega(p)=\nabla\bar{d}_\Omega(x)$ of a distance~$r_0$. Hence, if we move only of a
distance~$\theta d/2$, we obtain the desired ball. All in all, we can write explicitely
\bgs{
z_x:=x+\left(\frac{\theta d}{2}-\bar{d}_\Omega(x)\right)\nabla\bar{d}_\Omega(x).
}
For the details about this kind of geometric considerations concerning the signed distance function, see Appendix \ref{CH:3:A2}.

Now we split $I=I_1+I_2$, with
\bgs{
I_1:=\int_{\Omega_+}\bigg\{\int_{B_d(x)\setminus\Omega}\bigg[\G\left(\frac{N-u(y)}{|x-y|}\right)-\G\left(\frac{u(x)-u(y)}{|x-y|}\right)\bigg]
\frac{dy}{|x-y|^{n-1+s}}\bigg\}dx,
}
and
\bgs{
I_2:=\int_{\Omega_+}\bigg\{\int_{\Co B_d(x)\setminus\Omega}\bigg[\G\left(\frac{N-u(y)}{|x-y|}\right)-\G\left(\frac{u(x)-u(y)}{|x-y|}\right)\bigg]
\frac{dy}{|x-y|^{n-1+s}}\bigg\}dx.
}
Since $\G$ is globally Lipschitz---see~\eqref{CH:4:Lip_Gcal}---we have
\eqlab{\label{CH:4:miles12}
I_2\le\frac{\Lambda}{2}\int_{\Omega_+}\big(u(x)-N\big)\bigg(\int_{\Co B_d(x)}\frac{dy}{|x-y|^{n+s}}\bigg)dx
=\frac{\Lambda\mathcal H^{n-1}(\mathbb S^{n-1})}{2s}d^{-s}.
}

As for $I_1$, let $x\in\Omega_+$ and $y\in B_d(x)\setminus\Omega$.
Since by \eqref{CH:4:sonny_rollins} we have $u(y)\le N<u(x)$ and $G$ is increasing, we obtain
\bgs{
I_1\le-\int_{\Omega_+}\big(u(x)-N\big)\bigg[\int_{B_d(x)\setminus\Omega}G\left(\frac{N-u(y)}{|x-y|}\right)
\frac{dy}{|x-y|^{n+s}}\bigg]dx.
}
Exploiting \eqref{CH:4:sonny_rollins} and the monotonicity of $G$, we see that
\bgs{
G\left(\frac{N-u(y)}{|x-y|}\right)\ge G\left(\frac{d}{d}\right)=G(1)>0,
}
for every $x\in\Omega_+$ and $y\in B_d(x)\setminus\Omega$.
Recalling \eqref{CH:4:sonny_rollins2} we thus obtain
\bgs{
\int_{\Omega_+}\big(u(x)&-N\big)\bigg[\int_{B_d(x)\setminus\Omega}G\left(\frac{N-u(y)}{|x-y|}\right)
\frac{dy}{|x-y|^{n+s}}\bigg]dx\\
&
\ge G(1)\int_{\Omega_+}\big(u(x)-N\big)\left[\int_{B_\frac{\theta d}{2}(z_x)}\frac{dy}{|x-y|^{n+s}}\right]dx\\
&
\ge G(1)\frac{\big|B_\frac{\theta d}{2}(z_x)\big|}{(3\theta d)^{n+s}}\int_{\Omega_+}\big(u(x)-N\big)\,dx\\
&
=\frac{G(1)|B_1|}{2^n 3^{n+s}}(\theta d)^{-s}\int_{\Omega_+}\big(u(x)-N\big)\,dx.
}
Therefore, using also \eqref{CH:4:miles12} we get
\bgs{
I\le-\left(\frac{G(1)|B_1|}{2^n 3^{n+s}}\theta^{-s}-\frac{\Lambda\mathcal H^{n-1}(\mathbb S^{n-1})}{2s}\right)
d^{-s}\int_{\Omega_+}\big(u(x)-N\big)\,dx,
}
which is negative, provided we take $\theta$ small enough. This concludes the proof.
\end{proof}

\subsubsection{Proof of the interior local boundedness}\label{CH:4:Linftylocprop_proof}

We get now to the proof of Proposition \ref{CH:4:Linftylocprop}.

\begin{proof}[Proof of Proposition~\ref{CH:4:Linftylocprop}]
Let~$0 < \varrho < \tau \le 2 R$ and~$\eta \in C^\infty(\R^n)$ be a
cutoff function satisfying~$0 \le \eta \le 1$ in~$\R^n$,~$\mbox{supp}(\eta) \Subset B_\tau$,~$\eta = 1$
in~$B_\varrho$ and~$| \nabla \eta | \le 2 / (\tau - \varrho)$ in~$\R^n$. For~$k \ge 0$, we consider the
functions~$w = w_k := (u - k)_+$ and~$v := u - \eta w$. Clearly,~$v = u$ in~$\Co B_\tau$ and therefore
\begin{equation} \label{CH:4:Iineq}
\iint_{Q(B_\tau)} \frac{\I(x, y)}{|x - y|^{n - 1 + s}} \, dx\, dy \ge 0,
\end{equation}
with
$$
\I(x, y) := \G \left( \frac{v(x) - v(y)}{|x - y|} \right) - \G \left( \frac{u(x) - u(y)}{|x - y|} \right).
$$

We consider the sets~$A(k) := \{ x \in \R^n \,|\, u > k \}$ and~$A(k, t) := B_t \cap A(k)$, for~$t > 0$. First of all, we claim that
\begin{equation} \label{CH:4:estinBr}
\I(x, y) \le - \frac{\Lambda}{2} \, \frac{|w(x) - w(y)|}{|x - y|} + \lambda \, \chi_{B_\varrho^2 \setminus (B_\varrho \setminus A(k, \varrho))^2}(x, y) \quad \mbox{for } x, y \in B_\varrho,
\end{equation}
with~$\lambda$ and~$\Lambda$ as defined in~\eqref{CH:4:intglelambda} and~\eqref{CH:4:gintegr} respectively.
Clearly,~\eqref{CH:4:estinBr} holds for every~$x, y \in \Co A(k)$, since~$\I(x, y) = 0$ for these points and~$w = 0$ in~$\Co A(k)$.
Furthermore, it is also valid for~$x, y \in A(k, \varrho)$, as indeed, by~\eqref{CH:4:GGbetterbound} we have
$$
\I(x, y) = - \G \left( \frac{u(x) - u(y)}{|x - y|} \right) = - \G \left( \frac{w(x) - w(y)}{|x - y|} \right) \le - \frac{\Lambda}{2} \, \frac{|w(x) - w(y)|}{|x - y|} + \lambda.
$$
By symmetry, we are left to check~\eqref{CH:4:estinBr} for~$x \in A(k, \varrho)$ and~$y \in B_\varrho \setminus A(k, \varrho)$.
In this case, using~$u(x) > k \ge u(y)$ along with~\eqref{CH:4:GGbetterbound}, we get
\begin{align*}
\I(x, y) & = \G \left( \frac{k - u(y)}{|x - y|} \right) - \G \left( \frac{u(x) - u(y)}{|x - y|} \right) \le \frac{\Lambda}{2} \left( \frac{k - u(y)}{|x - y|} - \frac{u(x) - u(y)}{|x - y|} \right) + \lambda \\
& = - \frac{\Lambda}{2} \, \frac{u(x) - k}{|x - y|} + \lambda = -\frac{\Lambda}{2} \, \frac{|w(x) - w(y)|}{|x - y|} + \lambda.
\end{align*}
Hence,~\eqref{CH:4:estinBr} is verified.

We now claim that
\begin{equation} \label{CH:4:estoutBr}
\I(x, y) \le \frac{\Lambda}{2} \left( \chi_{B_\tau^2} \frac{|w(x) - w(y)|}{|x - y|} + \frac{w(x)}{\max \{ \tau - \varrho, |x - y| \}} \right) \quad \mbox{for } x \in B_\tau, \, y \in \Co B_\varrho.
\end{equation}
We already observed that~$\I(x, y) = 0$ for every~$x, y \in \Co A(k)$. When~$x \in B_\tau \setminus A(k, \tau)$ and~$y \in A(k)$, then
\begin{align*}
u(y) - u(x) & \ge u(y) - u(x) - \eta(y) (u(y) - k) = (1 - \eta(y)) u(y) + k \eta(y) - u(x) \\
& \ge (1 - \eta(y)) u(y) + k \eta(y) - k = (1 - \eta(y)) (u(y) - k) \ge 0
\end{align*}
and therefore
$$
\I(x, y) = \G \left( \frac{u(y) - u(x) - \eta(y) (u(y) - k)}{|x - y|} \right) - \G \left( \frac{u(y) - u(x)}{|x - y|} \right) \le 0,
$$
by the monotonicity properties of~$\G$. We are thus left to deal with~$x \in A(k, \tau)$ and~$y \in \Co B_\varrho$. In this case, by the Lipschitz character of~$\G$ and the properties of~$\eta$,
\begin{align*}
\I(x, y) & \le \frac{\Lambda}{2} \, \frac{|\eta(x) w(x) - \eta(y) w(y)|}{|x - y|} \le \frac{\Lambda}{2} \, \frac{\eta(y) |w(x) - w(y)| + w(x) |\eta(x) - \eta(y)|}{|x - y|} \\
& \le \frac{\Lambda}{2} \left( \chi_{B_\tau^2}(y) \frac{|w(x) - w(y)|}{|x - y|} + \min \left\{ \frac{1}{\tau - \varrho}, \frac{1}{|x - y|} \right\} w(x) \right),
\end{align*}
and~\eqref{CH:4:estoutBr} follows.

By taking advantage of estimates~\eqref{CH:4:estinBr} and~\eqref{CH:4:estoutBr} in~\eqref{CH:4:Iineq}, by symmetry we deduce that
\begin{align*}
\iint_{B_\varrho^2} \frac{|w(x) - w(y)|}{|x - y|^{n + s}} \, dx\, dy & \le C \, \Bigg\{ \iint_{B_\tau^2 \setminus B_\varrho^2} \frac{|w(x) - w(y)|}{|x - y|^{n + s}} \, dx \,dy + \int_{A(k, \varrho)} \int_{B_\varrho} \frac{dx\, dy}{|x - y|^{n - 1 + s}} \\
& \quad + \int_{B_\tau} w(x) \left( \frac{1}{\tau - \varrho} \int_{B_{\tau - \varrho}} \frac{dz}{|z|^{n - 1 + s}} + \int_{\Co B_{\tau - \varrho}} \frac{dz}{|z|^{n + s}} \right) dx \Bigg\} \\
& \le  C \left\{\iint_{B_\tau^2 \setminus B_\varrho^2} \frac{|w(x) - w(y)|}{|x - y|^{n + s}} \, dx\, dy + |A(k, \varrho)| \varrho^{1 - s} + \frac{\| w \|_{L^1(B_\tau)}}{(\tau - \varrho)^s} \right\} ,
\end{align*}
where for the second inequality we also used Lemma~\ref{CH:4:dumb_kernel_lemma}. Adding to both sides~$C$ times the left-hand side and dividing by~$1 + C$, we get that
$$
[w]_{W^{s, 1}(B_\varrho)} \le \theta \left( [w]_{W^{s, 1}(B_\tau)} + |A(k, \tau)| \tau^{1 - s} + \frac{\| w \|_{L^1(B_\tau)}}{(\tau - \varrho)^s} \right)
$$
for every~$0 < \varrho < \tau \le 2 R$ and for some constant~$\theta \in (0, 1)$ depending only on~$n$,~$s$
and~$g$. Applying, e.g,~\cite[Lemma~1.1]{GG82}, we infer that
$$
[w]_{W^{s, 1}(B_{(\varrho + \tau) / 2})} \le C \left( |A(k, \tau)| \tau^{1 - s} + \frac{\| w \|_{L^1(B_\tau)}}{(\tau - \varrho)^s} \right).
$$
Let~$\eta$ be a cutoff acting between the balls~$B_\varrho$ and~$B_{(3 \varrho + \tau) / 4}$. Then, by the fractional Sobolev inequality (see,~e.g.,~\cite[Theorem~1]{bourgbrez} or~\cite[Theorem~6.5]{HitGuide}) and computations similar to other made previously, we have that
\begin{align*}
\| w \|_{L^{\frac{n}{n - s}}(B_\varrho)} & \le \| \eta w \|_{L^{\frac{n}{n - s}}(\R^n)} \le C \int_{\R^n} \int_{\R^n} \frac{|\eta(x) w(x) - \eta(y) w(y)|}{|x- y|^{n + s}} \, dx\,dy\\
& \le C \left( [w]_{W^{s, 1}(B_{(\varrho + \tau) / 2})} + \frac{\| w \|_{L^1(B_\tau)}}{(\tau - \varrho)^{s}} \right).
\end{align*}
Combining the last two inequalities and recalling that~$w = w_k$, we arrive at
\begin{equation} \label{CH:4:wstarlew}
\| w_k \|_{L^{\frac{n}{n - s}}(B_\varrho)} \le C \left( |A(k, \tau)| \tau^{1 - s} + \frac{\| w_k \|_{L^1(B_\tau)}}{(\tau - \varrho)^s} \right)
\end{equation}
for every~$0 < r < \tau \le 2 R$ and~$k \ge 0$.

Take now~$k > h \ge 0$. We have
$$
\| w_h \|_{L^1(B_\tau)} \ge \int_{A(k, \tau)} (u(x) - h) \, dx \ge (k - h) |A(k, \tau)|
$$
and
$$
\| w_h \|_{L^1(B_\tau)} \ge \int_{A(k, \tau)} (u(x) - h) \, dx \ge \int_{A(k, \tau)} (u(x) - k) \, dx = \| w_k \|_{L^1(B_\tau)}.
$$
Thanks to these relations,~\eqref{CH:4:wstarlew}, and H\"older's inequality, it is easy to see that
\begin{equation} \label{CH:4:preiter}
\varphi(k, \varrho) \le \frac{C}{(k - h)^{s/n}} \left( \frac{\tau^{1 - s}}{k - h} + \frac{1}{(\tau - \varrho)^s} \right) \varphi(h, \tau)^{1 + \frac{s}{n}},
\end{equation}
where we set~$\varphi(\ell, t) := \| w_\ell \|_{L^1(B_t)}$.

Consider two sequences~$\{ k_j \}$ and~$\{ r_j \}$ defined by~$k_j := M (1 - 2^{-j})$ and~$r_j := R (1 + 2^{-j})$ for every non-negative integer~$j$, where~$M > 0$ will be chosen later. By applying~\eqref{CH:4:preiter} with~$k = k_{j + 1}$,~$\varrho = r_{j + 1}$,~$h = k_j$, and~$\tau = r_j$, setting~$\varphi_j := \varphi(k_j, r_j)$, and taking~$M \ge R$, we find
$$
\varphi_{j + 1} \le \frac{C (2^j \varphi_j)^{1 + \frac{s}{n}}}{(R \sqrt[n]{M})^s}.
$$
Applying now,~e.g.,~\cite[Lemma~7.1]{G03}, we conclude that~$\varphi_j$ converges to~$0$---i.e.,~$u \le M$ in~$B_R$---, provided we choose~$M$ in such a way that
$$
\| u_+ \|_{L^1(B_{2 R})} = \varphi_0 \le c_\sharp R^n M,
$$
for some constant~$c_\sharp > 0$ depending only on~$n$,~$s$ and~$g$. This concludes the proof.
\end{proof}

\subsection{Geometric minimizers}\label{CH:4:Geom_min_proof_Section}

This section is concerned with the minimizers of the geometric situation, which corresponds to the
choice~$g=g_s$. Here, we provide the proofs of Theorems~\ref{CH:4:equiv_intro} and~\ref{CH:4:uniqueness_teo},
which we
have stated in the Introduction.

We begin by proving the equivalence between
funtional minimizers, geometric minimizers, and
the various notions of solutions
to the equation~$\h_s u=0$ in~$\Omega$.
This result is the consequence of the main theorems proved in this chapter, together with the
interior regularity ensured by~\cite{CaCo}.

\begin{proof}[Proof of Theorem~\ref{CH:4:equiv_intro}]
$(i)\,\implies\,(ii)$ follows by Lemma~\ref{CH:4:weak_implies_min_lemma}.

As for the implication~$(ii)\,\implies\,(iii)$, by Proposition~\ref{CH:4:Linftylocprop} we know that~$u\in L^\infty_{\loc}(\Omega)$.
Then, let~$\Omega_k\Subset\Omega$ be a sequence of bounded open sets with Lipschitz boundary, such that
\bgs{
\Omega_k\Subset\Omega_{k+1}\quad\mbox{and}\quad\bigcup_{k\in\mathbb N}\Omega_k=\Omega,
}
let~$M_k$ be a diverging sequence such that
\eqlab{\label{CH:4:cylind_approx}
M_k\ge\|u\|_{L^\infty(\Omega_k)},
}
and consider the cylinders~$\Op^k:=\Omega_k\times(-M_k,M_k)$. We prove that~$\Sg(u)$ is $s$-minimal
in every~$\Op^k$. Since~$\Op^k\nearrow\Omega^\infty$, this readily implies that~$\Sg(u)$
is locally $s$-minimal in~$\Omega^\infty$, as wanted.

Let~$E\subseteq\R^{n+1}$ be such that~$E\setminus\Op^k=\Sg(u)\setminus\Op^k$ and let~$w_E$ be the function
defined in~\eqref{CH:4:rearr_func_def}. We can suppose that~$\Per_s(E,\Op^k)<\infty$,
otherwise there is nothing to prove. Then,
by~\eqref{CH:4:cylind_approx}, we know that the set~$E$ satisfies~\eqref{CH:4:EboundedinOmega}
and hence Theorem~\ref{CH:4:Persdecreases} implies that
\eqlab{\label{CH:4:equiv_proof_eqn1}
\Per_s(\Sg(w_E),\Op^k)\le\Per_s(E,\Op^k).
}
Notice that, by Proposition~\ref{CH:4:per_of_subgraph_prop}, we have~$w_E\in\B_{M_k}\W^s_u(\Omega_k)$.
Thus, since also~$u\in\B_{M_k}\W^s(\Omega_k)$, by identity~\eqref{CH:4:per_of_subgraph}, by the minimality of~$u$
and recalling~\eqref{CH:4:equiv_proof_eqn1},
we obtain
\bgs{
\Per_s(\Sg(u),\Op^k)\le\Per_s(\Sg(w_E),\Op^k)\le\Per_s(E,\Op^k).
}
The arbitrariness of the set~$E$ implies that~$\Sg(u)$ is $s$-minimal in~$\Op^k$, as claimed.

Now we prove that~$(iii)\,\implies\,(iv)$. First of all, we observe that~\cite[Theorem~1.1]{CaCo}
guarantees that~$u\in C^\infty(\Omega)$. Therefore, given any~$x\in\Omega$,
we can find both an interior and an exterior tangent ball to~$\Sg(u)$ at the boundary
point~$(x,u(x))\in\partial\Sg(u)\cap\Omega^\infty$.
The Euler-Lagrange equation satisfied by $s$-minimal sets---see~\cite[Theorem~5.1]{CRS10}---and
identity~\eqref{CH:4:Geom_ELop_curvature}
then imply that
\bgs{
\h_s u(x)=\I_s[\Sg(u)](x,u(x))=0.
}

The implication~$(iv)\,\implies\,(i)$ follows from Proposition~\ref{CH:4:BHprop_curvature}.
Indeed, given~$v\in C^\infty_c(\Omega)$, we can find a bounded open set~$\Omega'$ such that
\bgs{
supp\, v\Subset\Omega'\Subset\Omega.
}
Then, since~$u$ is smooth in~$\Omega$, we have~$u\in BH(\Omega')$ and hence
Proposition~\ref{CH:4:BHprop_curvature} implies that
\bgs{
\langle \h_s u,v\rangle=\int_{\Omega'}\h_su(x)v(x)\,dx=0.
}

We observe that~$(iv)\,\implies\,(v)$ always holds true, thanks to Remark~\ref{CH:4:mah}.
Finally, if we assume that~$u\in L^1_{\loc}(\R^n)$, then
we have implication~$(v)\,\implies\,(i)$ by Theorem~\ref{CH:4:Gen_viscweak}.
This concludes the proof of the Theorem.
\end{proof}


We observe that Corollary~\ref{CH:4:equiv_intro_global_corollary} is a straightforward consequence of 
Theorem~\ref{CH:4:equiv_intro}.

We pass to the proof of the uniqueness of the locally
$s$-minimal set with exterior data given by the subgraph of a function that is bounded in a big
enough neighborhood of~$\Omega$.

\begin{proof}[Proof of Theorem~\ref{CH:4:uniqueness_teo}]
By~\cite[Lemma~3.3]{graph} we know that
\bgs{
\Omega \times (-\infty, -M_0) \subseteq E \cap \Omega^\infty \subseteq \Omega \times (-\infty, M_0).
}
We observe that, since~$E$ is locally $s$-minimal in~$\Omega^\infty$ and~$\Omega$ has regular boundary,
by Theorem~\ref{CH:2:confront_min_teo} and Remark~\ref{CH:2:rmk_from_compact_to_any_subset} we know that~$E$ is $s$-minimal in~$\Omega^{M_0}$. In particular,
we have~$\Per_s(E,\Omega^{M_0})<\infty$. Moreover,
since~$E$ satisfies the hypothesis of Theorem~\ref{CH:4:Persdecreases}, we get
\bgs{
\Per_s(\Sg(w_E),\Omega^{M_0})\le\Per_s(E,\Omega^{M_0}),
}
and~$u:=w_E\in\B_{M_0}\W^s_\varphi(\Omega)$.
The $s$-minimality of~$E$ implies that~$E=\Sg(u)$, since otherwise the inequality would be strict.
Thus,~$\Sg(u)$ is locally $s$-minimal in~$\Omega^\infty$ and, by Theorem~\ref{CH:4:equiv_intro},~$u$
minimizes~$\F_s$ in~$\W^s_\varphi(\Omega)$.
The conclusion then follows from the uniqueness of such minimizer.
\end{proof}

\section{Nonparametric Plateau problem with obstacles}\label{CH:4:NPPWO_SEC}

In this section we consider the
Plateau problem with (eventually discontinuous) obstacles.
Namely, besides imposing the exterior data condition
$$u=\varphi\quad\mbox{a.e. in }\Co\Omega,$$
we constrain the functions to lie above a fixed function which acts as an obstacle, that is
$$u\ge\psi\quad\mbox{a.e. in }A,$$
where~$A\subseteq\Omega$ is a fixed open set.

We stress that the purpose of the present section is only that of showing that the functional setting introduced in the previous sections
can be easily adapted to study the obstacle problem, so we do not aim at full generality in the statements
nor in the proofs of our results.
In particular, we limit ourselves to consider bounded obstacles and we
prove the existence of a solution only in the case where the exterior data is
bounded in a big enough neighborhood of the domain~$\Omega$.
Furthermore, we will not investigate the regularity properties
of such a solution and of the free boundary.

We begin by introducing appropriate functional spaces. Given a bounded open set~$\Omega\subseteq\R^n$
with Lipschitz boundary,~$s\in(0,1)$, an open set~$A\subseteq\Omega$, an obstacle function~$\psi\in L^\infty(A)$,
the exterior data~$\varphi:\Co\Omega\to\R$,  and~$M\ge\|\psi\|_{L^\infty(A)}$, we define the spaces
\bgs{
&\K^s(\Omega,\varphi,A,\psi):=
\{u\in\W^s_\varphi(\Omega)\,|\,u\ge\psi\mbox{ a.e. in }A\},\\
&
\B\K^s(\Omega,\varphi,A,\psi):=
\{u\in\B\W^s_\varphi(\Omega)\,|\,u\ge\psi\mbox{ a.e. in }A\},\\
&
\B_M\K^s(\Omega,\varphi,A,\psi):=\K^s(\Omega,\varphi,A,\psi)\cap\B_M\W^s(\Omega).
}

We say that a function~$u\in\K^s(\Omega,\varphi,A,\psi)$ solves the obstacle problem if~$u$ minimizes~$\F$
in~$\K^s(\Omega,\varphi,A,\psi)$, i.e. if
\bgs{
\iint_{Q(\Omega)} \left\{ \G \left( \frac{u(x) - u(y)}{|x - y|} \right) - \G \left( \frac{v(x) - v(y)}{|x - y|} \right) \right\} \frac{dx\,dy}{|x - y|^{n - 1 + s}} \le 0,
}
for every~$v\in\K^s(\Omega,\varphi,A,\psi)$. We remark that this definition is well posed, thanks
to Lemma~\ref{CH:4:tartariccio}.

The main result of this section is the following existence and uniqueness Theorem.

\begin{theorem}\label{CH:4:Dirichlet_obstacle}
Let~$n\ge1$,~$s\in(0,1)$,~$\Omega\subseteq\Rn$ a bounded open set with Lipschitz boundary,~$R_0>1$ be such that
$\Omega\subseteq B_{R_0}$ and let~$\Theta=\Theta(n,s,g)>1$ be as in Theorem~\ref{CH:4:minareboundedthm}.
Let~$A\subseteq\Omega$ be an open set and let~$\psi\in L^\infty(A)$.
For every~$\varphi:\Co\Omega\to\R$ such that~$\varphi\in L^\infty(B_{\Theta R_0}\setminus\Omega)$, there exists a
unique function $u\in\K^s(\Omega,\varphi,A,\psi)$ that solves the obstacle problem.
Moreover
\bgs{
\|u\|_{L^\infty(\Omega)}\le R_0+\max\big\{\|\varphi\|_{L^\infty(B_{\Theta R_0}\setminus\Omega)},\|\psi\|_{L^\infty(A)}\big\}.
}
\end{theorem}

The proof of this Theorem is the content of Section~\ref{CH:4:proof_dir_obst_sect}.
It is interesting to observe that
a solution exists without having to impose regularity assumptions on the domain~$A$ where the obstacle is defined,
nor on the obstacle function~$\psi$---besides boundedness.

Before going on, we mention that in Section~\ref{CH:4:obst_geom_min_sect} we consider the geometric case corresponding to the choice~$g=g_s$
and we show the connection between solutions of the functional obstacle problem
and of the geometric obstacle problem.

Now we point out that a solution of the obstacle problem is a supersolution of the equation~$\h u=0$
in the whole domain~$\Omega$ and a solution away from the contact set, that is, formally:


\bgs{
\h u\ge0\quad\mbox{ in }\Omega\quad\mbox{and}\quad
\h u=0\quad\mbox{ in }\Omega\setminus\{u=\psi\}.
}

More precisely, we have the following result:
\begin{prop}
Let~$n\ge1$,~$s\in(0,1)$,~$\Omega\subseteq\R^n$ a bounded open set with Lipschitz
boundary,~$\varphi:\Co\Omega\to\R$,~$A\subseteq\Omega$ an open set and~$\psi\in L^\infty(A)$.
Suppose that there exists a function~$u\in\K^s(\Omega,\varphi,A,\psi)$ that solves the obstacle problem. Then
\bgs{
\langle\h u,v\rangle\ge0\quad\forall\,v\in C_c^\infty(\Omega)\quad\mbox{s.t. }v\ge0.
}
Furthermore, if~$\Op\subseteq\Omega$ is an open set such that
\bgs{
\inf_{\Op\cap A} (u-\psi)\ge\delta,
}
for some~$\delta>0$, then
\bgs{
\langle\h u,v\rangle=0\quad\forall\,v\in C^\infty_c(\Op).
}
In particular, if~$\Op$ has Lipschitz boundary, then~$u$ minimizes~$\F$ in~$\W^s_u(\Op)$.

\end{prop}

\begin{proof}
First of all, notice that if~$v\in C^\infty_c(\Omega)$ is such that~$v\ge0$, then~$u+\eps v\in\K^s(\Omega,\varphi,A,\psi)$
for every~$\eps>0$. Thus, by the minimality of~$u$ and recalling identity~\eqref{CH:4:Lem_tartariccio_eqn1} in Lemma~\ref{CH:4:tartariccio}, we have
\bgs{
\F^0(u+\eps v,\Omega)-\F^0(u,\Omega)\ge 0.
}
Passing to the limit~$\eps\to0^+$ and recalling Lemma~\ref{CH:4:1varlem}, we find~$\langle\h u,v\rangle\ge0$, as claimed.

In order to prove that~$u$ is a solution away from the contact set, let~$v\in C^\infty_c(\Op)$ and
observe that for every~$|\eps|\le \delta/\|v\|_{L^\infty(\Op)}$ we have~$u+\eps v\in\K^s(\Omega,\varphi,A,\psi)$.
Roughly speaking, since we are away from the contact set, we are allowed to deform the function~$u$ both from above and from below. Hence, again by the minimality of~$u$ and exploiting Lemma~\ref{CH:4:1varlem}, we obtain~$\langle\h u,v\rangle=0$.

Finally, if~$\Op$ has Lipschitz boundary, then we conclude that~$u$ minimizes~$\F$
in~$\W^s_u(\Op)$ by Lemma~\ref{CH:4:weak_implies_min_lemma}.
\end{proof}

In particular, we observe that if~$A\Subset\Omega$ has Lipschitz boundary,
then~$u$ minimizes~$\F$ in~$\W^s_u(\Omega\setminus\overline{A})$.

\subsection{Proof of Theorem~\ref{CH:4:Dirichlet_obstacle}}\label{CH:4:proof_dir_obst_sect}

The argument is essentially the same one that we already employed to prove the existence of minimizers of~$\F$
in~$\W^s_\varphi(\Omega)$. We begin by considering the functions~$u_M$
that minimize~$\F^M(\,\cdot\,,\Omega)$ in~$\B_M\K^s(\Omega,\varphi,A,\psi)$, then we show that they stabilize,
by exploiting Proposition~\ref{CH:4:truncdecreaseprop}. 

\subsubsection*{Step 1}
First of all, we observe that
\bgs{
\K^s&(\Omega,\varphi,A,\psi)\subseteq\W^s_\varphi(\Omega),
\qquad\B\K^s(\Omega,\varphi,A,\psi)\subseteq\B\W^s_\varphi(\Omega)\\
&
\quad\mbox{and}\quad
\B_M\K^s(\Omega,\varphi,A,\psi)\subseteq\B_M\W^s_\varphi(\Omega)
}
are closed convex subsets.
As a consequence, by arguing as in the proof of Proposition~\ref{CH:4:ertyui} and exploiting the convexity of~$\F^M$ ensured by
Lemma~\ref{CH:4:conv_func}, we find that
for every~$M\ge\|\psi\|_{L^\infty(A)}$ there exists a unique~$u_M\in\B_M\K^s(\Omega,\varphi,A,\psi)$ such that
$$\F^M(u_M,\Omega)=\inf\left\{\F^M(v,\Omega)\,|\,v\in\B_M\K^s(\Omega,\varphi,A,\psi)\right\}.$$

\subsubsection*{Step 2}

Now we remark that, since the obstacle~$\psi$ is bounded, we can apply Proposition~\ref{CH:4:truncdecreaseprop}
to obtain an a priori bound on the~$L^\infty$ norm of the minimizers~$u_M$, provided~$M>0$ is big enough.
Let indeed
\bgs{
u^{(N)}:=\chi_\Omega\min\{u,N\}+(1-\chi_\Omega)u,
}
and notice that, if~$u\in\K^s(\Omega,\varphi,A,\psi)$ and~$N\ge\sup_A \psi$, then we clearly
have~$u^{(N)}\in\K^s(\Omega,\varphi,A,\psi)$. Therefore, if we consider
\bgs{
N:=R_0+\max\left\{\sup_{B_{\Theta R_0}\setminus\Omega}\varphi,\sup_A\psi\right\},
}
then by Proposition~\ref{CH:4:truncdecreaseprop} and by the uniqueness of the minimizer of the functional~$\F^M(\,\cdot\,,\Omega)$
in~$\B_M\K^s(\Omega,\varphi,A,\psi)$, we obtain
\bgs{
\sup_\Omega u_M\le R_0+\max\left\{\sup_{B_{\Theta R_0}\setminus\Omega}\varphi,\sup_A\psi\right\},
}
for every~$M\ge N$.
Since we can argue in the same way by truncating the functions from below, we find that
\eqlab{\label{CH:4:obst_pf_eqn1}
\|u_M\|_{L^\infty(\Omega)}\le R_0+
\max\left\{\|\varphi\|_{L^\infty(B_{\Theta R_0}\setminus\Omega)},\|\psi\|_{L^\infty(A)}\right\}=:N_0,
}
for every~$M\ge N_0$.

\subsubsection*{Step 3}
Fix~$M_0:=N_0+1$ and observe that~\eqref{CH:4:obst_pf_eqn1}
ensures that
\eqlab{\label{CH:4:obstacle_monster_power_mannaggia}
\|u_{M_0}\|_{L^\infty(\Omega)}\le N_0<M_0.
}
We claim that this implies that the function~$u:=u_{M_0}$
solves the obstacle problem. In order to prove this, let us consider~$v\in\B\K^s(\Omega,\varphi,A,\psi)$
and notice that by~\eqref{CH:4:obstacle_monster_power_mannaggia} we have
\bgs{
w:=t v+(1-t)u\in\B_{M_0}\K^s(\Omega,\varphi,A,\psi),
}
provided~$t\in(0,1)$ is small enough. Thus, by the minimality of~$u$ and exploiting the convexity of~$\F^{M_0}$,
we find
\bgs{
\F^{M_0}(u,\Omega)\le\F^{M_0}(w,\Omega)\le t\F^{M_0}(v,\Omega)+(1-t)\F^{M_0}(u,\Omega),
}
that is
\bgs{
\F^{M_0}(u,\Omega)\le\F^{M_0}(v,\Omega).
}
This shows that~$u$ minimizes~$\F^{M_0}(\,\cdot\,,\Omega)$ in~$\B\K^s(\Omega,\varphi,A,\psi)$.
Thanks to Lemma~\ref{CH:4:tartariccio}, this implies that~$u$ minimizes~$\F$ in the larger space~$\K^s(\Omega,\varphi,A,\psi)$
and hence solves the obstacle problem. Finally, the strict convexity of~$\F^M$
guarantees the uniqueness of such a solution---see also point~$(iii)$ of Remark~\ref{CH:4:tartapower_Remark}---concluding
the proof of the Theorem.

\subsection{Geometric obstacle problem}\label{CH:4:obst_geom_min_sect}

In this section we study the obstacle problem for the fractional perimeter in the unbounded domain~$\Omega^\infty$.
This problem has been recently considered---in the case of bounded domains---in~\cite{GeomObst}, where
the authors proved a regularity result for the solution.
Our aim consists in showing that, also in the presence of obstacles, minimizers of the functional problem are
minimizers for the geometric problem.

Again, we stress that we do not aim at full generality. In particular, we limit ourselves to give the definition of
a geometric minimizer in the setting that interests us, by considering only as domain
a bounded open set~$\Omega\subseteq\R^n$ with Lipschitz boundary and
bounded obstacles.

Let~$n\ge1$,~$s\in(0,1)$,~$\Omega\subseteq\R^n$ a bounded open set with Lipschitz
boundary,~$\varphi:\R^n\to\R$,~$A\subseteq\Omega$ an open set and~$\psi\in L^\infty(A)$. We define
\bgs{
\Op:=\left\{(x,x_{n+1})\in\R^{n+1}\,|\,x\in A\mbox{ and }x_{n+1}<\psi(x)\right\}.
}
We say that a set~$E\subseteq\R^{n+1}$ such that~$E\setminus\Omega^\infty=\Sg(\varphi)\setminus\Omega^\infty$
and~$\Op\subseteq E$ solves the geometric obstacle problem
if for every~$M\ge\|\psi\|_{L^\infty(A)}$ it holds~$\Per_s(E,\Omega^M)<\infty$, and
\bgs{
\Per_s(E,\Omega^M)\le\Per_s(F,\Omega^M),
}
for every~$F\subseteq\R^{n+1}$ such that~$F\setminus\Omega^M=E\setminus\Omega^M$
and~$\Op\subseteq F$.

\begin{remark}\label{CH:4:geom_obst_locmin_remark}
We observe that if~$E\subseteq\R^{n+1}$ solves the geometric obstacle problem, then it is locally $s$-minimal in the open
set~$\Omega^\infty\setminus\overline{\Op}$.

\end{remark}

By exploiting Theorem~\ref{CH:4:Persdecreases}
and Proposition~\ref{CH:4:per_of_subgraph_prop},
it is readily seen that if~$u$ solves the obstacle problem---with~$g=g_s$---then its subgraph
solves the geometric obstacle problem. 
\begin{prop}
Let~$n\ge1$,~$s\in(0,1)$,~$g=g_s$,~$\Omega\subseteq\Rn$ a bounded open set with Lipschitz boundary,~$R_0>1$ be such that $\Omega\subseteq B_{R_0}$ and let~$\Theta=\Theta(n,s,g_s)>1$ be as in Theorem~\ref{CH:4:minareboundedthm}.
Let~$A\subseteq\Omega$ be an open set,~$\psi\in L^\infty(A)$ and let~$\varphi:\Co\Omega\to\R$ such that~$\varphi\in L^\infty(B_{\Theta R_0}\setminus\Omega)$.
Let~$u\in\B\K^s(\Omega,\varphi,A,\psi)$ be the unique solution of the obstacle problem, as in
Theorem~\ref{CH:4:Dirichlet_obstacle}. Then,~$\Sg(u)$ solves the geometric obstacle problem.
\end{prop}

We conclude this section by proving that the subgraph of~$u$ is actually the unique solution to the geometric
obstacle problem.

In order to do this, we consider~$\Omega\subseteq\R^n$ to be a bounded open set with~$C^2$ boundary
and the domain of definition of the obstacle~$A\subseteq\Omega$
to be either~$A=\Omega$ or~$A\Subset\Omega$ with~$C^2$ boundary.
Since we are considering a bounded obstacle~$\psi\in L^\infty(A)$---and thanks to Remark~\ref{CH:4:geom_obst_locmin_remark}---it is easy to check that the argument of the proof of~\cite[Lemma~3.3]{graph} works also in this situation.

Therefore, there exists~$\tilde{R}(n,s,\Omega)>0$ such that, if~$\varphi\in L^\infty(B_{\tilde{R}}\setminus\Omega)$,
and~$E\subseteq\R^{n+1}$ solves the geometric obstacle problem, then
\bgs{
\Omega\times(-\infty,-M_0)\subseteq E\cap\Omega^\infty
\subseteq\Omega\times(-\infty,M_0),
}
for some~$M_0(n,s,\Omega,\varphi,A,\psi)>0$.
Let us now define~$R_s:=\max\big\{\Theta R_0,\tilde{R}\big\}$.
Then, we have the following uniqueness result:
\begin{prop}
Let~$n\ge1$,~$s\in(0,1)$,~$g=g_s$,~$\Omega\subseteq\R^n$ a bounded open set with~$C^2$ boundary.
Let~$A\subseteq\Omega$ be an open set such that either~$A=\Omega$ or~$A\Subset\Omega$ with~$C^2$ boundary,~$\psi\in L^\infty(A)$ and let~$\varphi:\Co\Omega\to\R$ such that~$\varphi\in L^\infty(B_{R_s}\setminus\Omega)$,
with~$R_s$ as defined above. Let~$u\in\B\K^s(\Omega,\varphi,A,\psi)$ be the unique solution of the obstacle problem, as in
Theorem~\ref{CH:4:Dirichlet_obstacle}. Then,~$\Sg(u)$ is the unique solution of the geometric obstacle problem.
\end{prop}

The proof follows by arguing as in the proof of Theorem~\ref{CH:4:uniqueness_teo} and exploiting the uniqueness
of the solution of the (functional) obstacle problem.

\section{Approximation results}\label{CH:4:Appro_Section}

In this section we collect some approximating results for the functionals~$\F^M(\,\cdot\,,\Omega)$.
These results are interesting for various reasons. First of all, they are meaningful in themselves and they somehow
complement the results proven in Section \ref{CH:2:GCAABSS_SEC}.
More precisely, in Proposition~\ref{CH:4:subgraph_appro_prop} we show that a subgraph having finite $s$-perimeter can be approximated with smooth subgraphs,
and not just with arbitrary smooth open sets as in Theorem \ref{CH:2:density_smooth_teo}.

Secondarily, we point out that the subgraphs of $\sigma$-harmonic functions are somehow less rigid than nonlocal minimal graphs.
Indeed, thanks to the surprising result proved in~\cite{DSV17}, it is always possible to approximate a nonlocal minimal graph with
$\sigma$-harmonic functions---see Theorem~\ref{CH:4:appro_nonlomin_fracharm_theorem}.
On the other hand---as observed in~\cite{DSV17}---the converse is not possible, because
by exploiting~\cite[Theorem~1.1]{DSV17} it is possible to construct $\sigma$-harmonic functions that oscillate wildly,
while nonlocal minimal graphs must satisfy uniform density estimates at the boundary points---see~\cite[Theorem~4.1]{CRS10}.

Finally, we prove that there is no gap phenomenon when we minimize~$\F$ with respect to regular exterior
data---see Proposition~\ref{CH:4:NoLavrentievgap_prop}.
Indeed, as we have remarked in the introduction, even when the exterior data is a smooth and compactly supported function,
the minimizer of~$\F$, in general, is not continuous across the boundary of the domain, because of stickiness effects which are
typically nonlocal.


Thus, it is natural to wonder whether the minimization of~$\F$ among functions which are smooth
in the whole of~$\R^n$ leads to a value which is strictly bigger than that obtained by
minimizing~$\F$ in the larger space~$\W^s_\varphi(\Omega)$.
Roughly speaking, given~$\varphi\in C^{0,1}(\R^n)$, we wonder whether the inequality
\bgs{
\inf\left\{\F^0(v,\Omega)\,|\,v\in C^{0,1}(\R^n)\mbox{ s.t. }v=\varphi\mbox{ in }\Co\Omega\right\}
\ge\inf_{v\in\W^s_\varphi(\Omega)}\F^0(v,\Omega)
}
can be strict. The answer is no. As shown by Proposition~\ref{CH:4:NoLavrentievgap_prop}, this inequality is actually always an equality.

\medskip

First of all, we remark that when we keep the exterior data fixed, then the approximation result follows
from Lemma~\ref{CH:4:tartariccio} and the density of~$C^\infty_c(\Omega)$ in~$W^{s,1}(\Omega)$.
We have alreday exploited this fact in the proof of the existence of minimizers---see also point~$(iv)$
of Remark~\ref{CH:4:tartapower_Remark}.

On the other hand, when we approximate also the exterior data we have the following useful result.

\begin{prop}\label{CH:4:conv_change_ext_data_prop}
Let~$n\ge1$,~$s\in(0,1)$,~$\Omega\subseteq\R^n$ a bounded open set with Lipschitz boundary.
Let~$u,\,u_k\in L^1_{\loc}(\R^n)\cap W^{s,1}(\Omega_d)$, for some~$d>0$, and suppose
that~$u_k\to u$ both in~$L^1_{\loc}(\R^n)$ and in~$W^{s,1}(\Omega_d)$. Then
\bgs{
\lim_{k\to\infty}\F^M(u_k,\Omega)=\F^M(u,\Omega),
}
for every~$M\ge0$.
\end{prop}

Before getting to the proof of Proposition~\ref{CH:4:conv_change_ext_data_prop},
we state some of its consequences.

If we consider a symmetric mollifier~$\eta\in C^\infty_c(\R^n)$ as in~\eqref{CH:4:symm_moll_def}
and we define the mollified functions~$u_\eps:=u\ast\eta_\eps$, then we obtain the following corollary
of Proposition~\ref{CH:4:conv_change_ext_data_prop}.

\begin{corollary}\label{CH:4:mollified_convergence_corollary}
Let~$n\ge1$,~$s\in(0,1)$,~$\Omega\subseteq\R^n$ a bounded open set with Lipschitz boundary
and let~$u\in L^1_{\loc}(\R^n)\cap W^{s,1}(\Omega_d)$, for some~$d>0$. Then
\bgs{
\lim_{\eps\to0}\F^M(u_\eps,\Omega)=\F^M(u,\Omega),
}
for every~$M\ge0$.
\end{corollary}

\begin{proof}
It is well known that~$u_\eps\to u$ in~$L^1_{\loc}(\R^n)$. On the other hand, since~$u\in W^{s,1}(\Omega_d)$,
we have also~$u_\eps\to u$ in~$W^{s,1}(\Omega_{d/2})$---see, e.g., Lemma~\ref{CH:4:convol_conv_frac}.
Hence, the conclusion follows from Proposition~\ref{CH:4:conv_change_ext_data_prop}.
\end{proof}

When we consider subgraphs of locally bounded functions, Proposition~\ref{CH:4:per_of_subgraph_prop}
and Corollary~\ref{CH:4:mollified_convergence_corollary}
straightforwardly imply the desired approximation result which complements Theorem \ref{CH:2:density_smooth_teo}.

\begin{prop}\label{CH:4:subgraph_appro_prop}
Let~$n\ge1$,~$s\in(0,1)$,~$\Omega\subseteq\R^n$ a bounded open set with Lipschitz boundary
and let~$u\in L^1_{\loc}(\R^n)\cap W^{s,1}(\Omega)\cap L^\infty_{\loc}(\Omega)$. Then,
for every open set~$\Op\Subset\Omega$ with Lipschitz boundary and every~$M\ge\|u\|_{L^\infty(\Op_r)}$,
with~$r:=d(\Op,\partial\Omega)/2$,
it holds
\eqlab{\label{CH:4:horkheimer_sucks}
\lim_{\eps\to0}\Per_s(\Sg(u_\eps),\Op^M)=\Per_s(\Sg(u),\Op^M).
}
Moreover, if~$u\in C(\mathcal U)$, for some open set~$\mathcal U\subseteq\R^n$, then
for every compact set~$K\Subset\mathcal U$ and every~$\delta>0$ we have
\bgs{
\partial\Sg(u_\eps)\cap K^\infty\subseteq N_\delta\big(\Sg(u)\big)\cap K^\infty,
}
for every~$\eps>0$ small enough.
\end{prop}

\begin{proof}
It is enough to notice that for every~$\eps>0$ small enough we have
\bgs{
\|u_\eps\|_{L^\infty(\Op)}\le\|u\|_{L^\infty(\Op_r)}\le M.
}
Then,~\eqref{CH:4:horkheimer_sucks} follows by making use of Corollary~\ref{CH:4:mollified_convergence_corollary}
and of identity~\eqref{CH:4:per_of_subgraph}.
To conclude, notice that~$u\in C(\mathcal U)$ implies that
\bgs{
\partial\Sg(u)\cap\mathcal U^\infty
=\{(x,u(x))\in\R^{n+1}\,|\,x\in\mathcal U\},
}
and similarly for~$u_\eps$. Thus, the uniform convergence of the boundaries follows
from the fact that~$u_\eps\to u$ locally uniformly in~$\mathcal U$.
\end{proof}

By exploiting~\cite[Theorem 1.1]{DSV17} to approximate the mollified functions~$u_\eps$,
it is immediate to see that we can find a sequence of $\sigma$-harmonic functions~$u_k$ such
that~$\F^M(u_k,\Omega)\to\F^M(u,\Omega)$. In particular, if the function~$u$ is bounded in~$\Omega$
and we take~$g=g_s$, then we can approximate the $s$-perimeter of the subgraph of~$u$
with the $s$-perimeter of the subgraphs of the functions~$u_k$.

We give a precise statement of this fact only in the case of nonlocal minimal graphs.

\begin{theorem}\label{CH:4:appro_nonlomin_fracharm_theorem}
Let~$n\ge1$,~$s\in(0,1)$,~$\Omega\subseteq\R^n$ a bounded open set with Lipschitz boundary,
let~$\varphi:\Co\Omega\to\R$ be such
that
\bgs{
\varphi\in L^1_{\loc}(\R^n\setminus\Omega)\cap W^{s,1}(\Omega_d\setminus\Omega)\cap L^\infty(\Omega_d\setminus\Omega),
}
for some~$d>0$ small, and let~$R_0>0$ such
that~$\Omega_d\Subset B_{R_0}$.
Let~$u\in\B\W^s_\varphi(\Omega)$ be the unique minimizer of~$\F_s$ in~$\W^s_\varphi(\Omega)$.
Then, for every fixed~$\sigma\in(0,1)$
and~$\ell\in\mathbb N$, there exists a sequence of compactly supported functions~$u_k\in H^\sigma(\R^n)\cap C^\sigma(\R^n)$ such that
\bgs{
&(i)\quad(-\Delta)^\sigma u_k=0\qquad\mbox{in }B_{k+R_0}\\
&
(ii)\quad u_k\to u\qquad\mbox{in }L^1_{\loc}(\R^n)\mbox{ and in }W^{s,1}(\Omega_{d/2})\\
&
(iii)\quad\lim_{k\to\infty}\|u_k-u\|_{C^\ell(\Omega')}=0\qquad\mbox{for every }\Omega'\Subset\Omega,\\
&
(iv)\quad\|u_k\|_{L^\infty(\Omega)}\le\|u\|_{L^\infty(\Omega_d)}+1,\\
&
(v)\quad\lim_{k\to\infty}\Per_s(\Sg(u_k),\Omega^M)=\Per_s(\Sg(u),\Omega^M),\qquad\mbox{for every }M\ge\|u\|_{L^\infty(\Omega_d)}+1.
}
Moreover, for every compact set~$K\Subset\Omega$ and every~$\delta>0$ it holds
\bgs{
\partial\Sg(u_k)\cap K^\infty\subseteq N_\delta\big(\Sg(u)\big)\cap K^\infty,
}
for every~$k$ big enough.
\end{theorem}

\begin{proof}
We begin by observing that, recalling Lemma~\ref{CH:4:tail_equiv_cond_Lemma}
and exploiting Theorem~\ref{CH:4:Dirichlet}, we know that there exists a unique function~$u\in\W^s_\varphi(\Omega)$
that minimizes~$\F$. Moreover, since~$\varphi$ is bounded near~$\partial\Omega$, by Theorem~\ref{CH:4:Bdary_Bdedness_Thm}
we know that~$u\in L^\infty(\Omega)$. Finally, by~\cite[Theorem~1.1]{CaCo} we have~$u\in C^\infty(\Omega)$.
We also remark that, since~$\varphi\in W^{s,1}(\Omega_d\setminus\Omega)$, it is readily seen that~$u\in W^{s,1}(\Omega_d)$.
As a consequence, we have that
\eqlab{\label{CH:4:horkheimer_reeeaaally_sucks}
&u_\eps\to u\qquad\mbox{in }L^1_{\loc}(\R^n)\mbox{ and in }W^{s,1}(\Omega_{d/2})\\
&
\lim_{\eps\to0}\|u_\eps-u\|_{C^\ell(\Omega')}=0\qquad\mbox{for every }\Omega'\Subset\Omega,\\
&
\|u_\eps\|_{L^\infty(\Omega)}\le\|u\|_{L^\infty(\Omega_d)}\qquad\mbox{for every }\eps>0\mbox{ small enough}.
}

Then, the claim follows by using~\cite[Theorem 1.1]{DSV17} to approximate the mollified functions~$u_\eps$
and a diagonal argument.
Indeed, fix~$\eps>0$ and notice that, since~$u_\eps$ is smooth in~$\R^n$, by~\cite[Theorem 1.1]{DSV17}
we can find for every~$k\in\mathbb N$ a compactly supported function~$u_k\in H^\sigma(\R^n)\cap C^\sigma(\R^n)$ such that
\bgs{
&(-\Delta)^\sigma u_k=0\qquad\mbox{in }B_{k+R_0}\\
&\|u_k-u_\eps\|_{C^\ell(B_{k+R_0})}<\frac{1}{e^k}.
}
In particular, this implies that
\bgs{
&u_k\to u_\eps\qquad\mbox{in }L^1_{\loc}(\R^n)\mbox{ and in }W^{s,1}(\Omega_{d}),\\
&\|u_k\|_{L^\infty(\Omega)}\le\|u_\eps\|_{L^\infty(\Omega)}+1\le\|u\|_{L^\infty(\Omega_d)}+1.
}
Therefore, after a diagonal argument and recalling~\eqref{CH:4:horkheimer_reeeaaally_sucks}, we obtain a sequence of compactly supported functions~$u_k\in H^\sigma(\R^n)\cap C^\sigma(\R^n)$ that
satisfies points~$(i),~(ii),~(iii)$ and~$(iv)$. Then, point~$(v)$ follows by points~$(ii)$ and~$(iv)$,
Proposition~\ref{CH:4:conv_change_ext_data_prop} and identity~\eqref{CH:4:per_of_subgraph}.

To conclude, notice that the locally uniform convergence of the boundaries follows from point~$(iii)$---used just for
the~$C^0$ norm, as in the proof of Proposition~\ref{CH:4:subgraph_appro_prop}..
\end{proof}

Now we provide the proof of Proposition~\ref{CH:4:conv_change_ext_data_prop}.

\begin{proof}[Proof of Proposition~\ref{CH:4:conv_change_ext_data_prop}]
First of all, we observe that by the Lipschitzianity of~$\G$---see~\eqref{CH:4:Lip_Gcal}---we have
\bgs{
\big|\A(u_k,\Omega)-\A(u,\Omega)|\le\frac{\Lambda}{2}\|u_k-u\|_{W^{s,1}(\Omega)}.
}
As for the nonlocal part, we will exploit identity~\eqref{CH:4:nonlocal_explicit} and, again, the Lipschitzianity of~$\G$.
We have
\bgs{
\big|\Nl^M(u_k,\Omega) & -\Nl^M(u,\Omega)\big| \le \int_{\Omega} \left\{ \int_{\Co \Omega} \left| 2 \, \G \left( \frac{u_k(x)-u_k(y)}{|x-y|} \right) -
\G \left( \frac{M+u_k(y)}{|x-y|} \right) \right. \right. \\
& \quad \left. \left. -\G \left( \frac{M-u_k(y)}{|x-y|} \right) -
2 \, \G \left( \frac{u(x)-u(y)}{|x-y|} \right) +
\G \left( \frac{M+u(y)}{|x-y|} \right) \right. \right. \\
& \quad \left. \left. +\G \left( \frac{M-u(y)}{|x-y|} \right)
\right| \frac{dy}{\kers} \right\} dx.
}
We split the domain~$\Co\Omega=\big(\Omega_r\setminus\Omega\big)\cup \big(B_R\setminus\Omega_r\big)
\cup\Co B_R$, with~$r\in(0,d)$ small enough such that~$\Omega_r$ has Lipschitz boundary, and~$R>0$ big---we will let~$R\to\infty$ in the end---and we treat the three cases differently.

We begin by observing that---by appropriately regrouping the terms, using the triangle inequality, the Lipschitzianity of~$\G$ and exploiting also Corollary~\ref{CH:4:FHI_corollary}---the double integral over~$\Omega\times\big(\Omega_r\setminus\Omega\big)$
can be estimated by
\bgs{
\Lambda\int_\Omega\left\{\int_{\Omega_r\setminus\Omega}\frac{|u_k(x)-u_k(y)-u(x)-u(y)|+|u_k(y)-u(y)|}{|x-y|^{n+s}}\,dy
\right\}dx\le C\|u_k-u\|_{W^{s,1}(\Omega_r)}.
}
Similarly, the double integral over~$\Omega\times\big(B_R\setminus\Omega_r\big)$ can be estimated by
\bgs{
\Lambda\int_\Omega&\left\{\int_{B_R\setminus\Omega_r}\frac{|u_k(x)-u_k(y)-u(x)-u(y)|+|u_k(y)-u(y)|}{|x-y|^{n+s}}\,dy
\right\}dx\\
&
\le\Lambda\left\{\int_\Omega|u_k(x)-u(x)|\left(\int_{B_R\setminus\Omega_r}\frac{dy}{|x-y|^{n+s}}\right)dx\right.\\
&\qquad
\left.+2\int_{B_R\setminus\Omega_r}|u_k(y)-u(y)|\left(\int_\Omega\frac{dx}{|x-y|^{n+s}}\right)dy\right\}\\
&
\le\frac{3\,\Lambda\,\Ha^{n-1}(\mathbb S^{n-1})}{s\,r^s}\|u_k-u\|_{L^1(B_R)}.
}
Now we observe that, since~$u_k\to u$ in~$L^1(\Omega)$, we have~$\|u_k\|_{L^1(\Omega)}\le 2\|u\|_{L^1(\Omega)}$
for all~$k$ big enough. Moreover, we take~$R_0>0$ such that~$\Omega\Subset B_{R_0}$ and~$R>R_0$. Then, by regrouping the terms in a different way, we estimate the double integral over~$\Omega\times\Co B_R$ with
\bgs{
\frac{\Lambda}{2}\int_\Omega&\left\{|u_k(x)-M|+|u_k(x)+M|+|u(x)-M|+|u(x)+M|\right\}\left(\int_{\Co B_R}\frac{dy}{|x-y|^{n+s}}\right)dx\\
&
\le\Lambda\int_\Omega\left\{|u_k(x)|+|u(x)|+2M\right\}\left(\int_{\Co B_{R-R_0}(x)}\frac{dy}{|x-y|^{n+s}}\right)dx\\
&
\le C\frac{\|u\|_{L^1(\Omega)}+M\,|\Omega|}{(R-R_0)^s}.
}
All in all, we have proved that
\bgs{
\big|\F^M(u_k,\Omega) -\F^M(u,\Omega)\big| \le C\left(\|u_k-u\|_{W^{s,1}(\Omega_r)}+\|u_k-u\|_{L^1(B_R)}
+\frac{\|u\|_{L^1(\Omega)}+M\,|\Omega|}{(R-R_0)^s}\right).
}
Passing first to the limit~$k\to\infty$, then to the limit~$R\to\infty$, concludes the proof of the Proposition.
\end{proof}

We conclude this section by proving that there is no gap phenomenon
in the minimization of~$\F$. This is a simple consequence of the density of~$C^\infty_c(\Omega)$
in~$W^{s,1}(\Omega)$, which sostantially means that functions in~$W^{s,1}(\Omega)$ do not have a well defined
trace. Roughly speaking, this implies that we can approximate any function~$u\in W^{s,1}(\Omega)$
with smooth functions that have a fixed boundary value.

\begin{prop}\label{CH:4:NoLavrentievgap_prop}
Let~$n\ge1$,~$s\in(0,1)$,~$\Omega\subseteq\R^n$ a bounded open set with Lipschitz boundary
and let~$\varphi:\R^n\to\R$ be such that~$\varphi\in C^{0,1}(\Omega_d)$, for some~$d>0$,
and~$\varphi\in L^1(\Omega_{\Theta\diam(\Omega)})$, with~$\Theta>1$ as in Theorem~\ref{CH:4:Dirichlet}.
Then,
\bgs{
\inf\left\{\F^0(v,\Omega)\,|\,v\in C^{0,1}(\Omega_d)\mbox{ s.t. }v=\varphi\mbox{ in }\Omega_d\setminus\Omega
\mbox{ and a.e. in }\Co\Omega_d\right\}=\min_{v\in\W^s_\varphi(\Omega)}\F^0(v,\Omega).
}
\end{prop}

\begin{proof}
Notice that~$\varphi\in C^{0,1}(\Omega_d)$ implies that~$\varphi\in W^{s,1}(\Omega_d)$, since
\bgs{
\int_{\Omega_d}\int_{\Omega_d}\frac{|\varphi(x)-\varphi(y)|}{|x-y|^{n+s}}\,dx\,dy
\le[\varphi]_{C^{0,1}(\Omega_d)}\int_{\Omega_d}\int_{\Omega_d}\frac{dx\,dy}{|x-y|^{n-1+s}},
}
which is finite, thanks to Lemma~\ref{CH:4:dumb_kernel_lemma}. Then, by recalling point~$(i)$ of Lemma~\ref{CH:4:tail_equiv_cond_Lemma}
and exploiting Theorem~\ref{CH:4:Dirichlet}, we know that there exists a unique function~$u\in\W^s_\varphi(\Omega)$
that minimizes~$\F$. By point~$(ii)$ of Remark~\ref{CH:4:tartapower_Remark}, this means that
\bgs{
\F^0(u,\Omega)=\inf_{v\in\W^s_\varphi(\Omega)}\F^0(v,\Omega).
}
Moreover, since~$w:=u-\varphi\in W^{s,1}(\Omega)$, by the density of~$C^\infty_c(\Omega)$ in~$W^{s,1}(\Omega)$,
we can find a sequence~$\{w_k\}\subseteq C^\infty_c(\Omega)$ such that
\bgs{
\lim_{k\to\infty}\|w_k-w\|_{W^{s,1}(\Omega)}=0.
}
If we extend the functions~$w_k$ by zero outside~$\Omega$,
this means that the functions~$v_k:=\varphi+w_k\in C^{0,1}(\Omega_d)$ converge to~$u$ in~$W^{s,1}(\Omega)$.
By Lemma~\ref{CH:4:tartariccio}, this implies
\bgs{
\lim_{k\to\infty}\F^0(v_k,\Omega)=\F^0(u,\Omega),
}
concluding the proof of the Proposition.
\end{proof}

\end{chapter}

\begin{chapter}[Bernstein-Moser-type results for nonlocal minimal graphs]{Bernstein-Moser-type results for nonlocal minimal graphs}\label{CH_Bern_Mos_result}
	
\minitoc


\section{Introduction and main results}


For simplicity, in this chapter sets that minimize~$\Per_s$ in all bounded open subsets of~$\R^{n + 1}$ will be simply called \emph{$s$-minimal} and their boundaries~\emph{$s$-minimal surfaces}.


In this brief chapter we are mostly interested in~$s$-minimal sets~$E \subseteq \R^{n + 1}$ that are subgraphs of a measurable function~$u: \R^n \to \R$, i.e., that satisfy
\begin{equation} \label{CH:5:Esubgraph}
E = \left\{ x = (x', x_{n + 1}) \in \R^n \times \R \,|\, x_{n + 1} < u(x') \right\}.
\end{equation}
We will call the boundaries of such extremal sets~\emph{$s$-minimal graphs}.

We observe that, differently from the previous chapters, we will use here the notation $x=(x',x_{n+1})\in\R^{n+1}$.

We recall that, if $u: \R^n \to \R$ is a function of class~$C^{1, 1}$ in a neighborhood of a point~$x' \in \R^n$,
and $E:=\Sg(u)$ as in \eqref{CH:5:Esubgraph}, then
\[
H_s[E](x',u(x')) = \cu_s u(x'),
\]
with
\[
\cu_s u(x') := 2 \, \PV \int_{\R^n} G_s \left( \frac{u(x')-u(y')}{|x' - y'|} \right) \frac{dy'}{|x' - y'|^{n+s}} 
\]
and
\begin{equation} \label{CH:5:Gdef}
G_s(t):=\int_0^t\frac{d\tau}{(1+\tau^2)^\frac{n+1+s}{2}} \quad \mbox{for } t\in\R.
\end{equation}

Taking advantage of the convexity of the energy functional associated to~$\cu_s$ and of a suitable rearrangement inequality, we have shown in Chapter \ref{CH_Nonparametric} that a set~$E$ given by~\eqref{CH:5:Esubgraph} for some function~$u: \R^n \to \R$ is~$s$-minimal if and only if~$u$ is a solution of
\begin{equation} \label{CH:5:Hu=0}
\cu_s u = 0 \quad \mbox{in } \R^n.
\end{equation}
There are several notions of solutions of~\eqref{CH:5:Hu=0}, such as smooth solutions, viscosity solutions, and weak solutions. However, all such definitions are equivalent under mild assumptions on~$u$---see Corollary \ref{CH:4:equiv_intro_global_corollary} for more details.
In what follows, a solution of~\eqref{CH:5:Hu=0} will always indicate a function~$u \in C^\infty(\R^n)$ that satisfies identity~\eqref{CH:5:Hu=0} pointwise. We stress that no growth assumptions at infinity are made on~$u$.

The main contribution of this chapter is the following result.

\begin{theorem} \label{CH:5:underPakmainthm}
Let~$n \ge \ell \ge 1$ be integers,~$s \in (0, 1)$, and suppose that
\begin{equation} \tag{$P_{s, \ell}$} \label{CH:5:Paellprop}
\mbox{there exist no singular~$s$-minimal cones in~$\R^\ell$.}
\end{equation}
Let~$u$ be a solution of~\eqref{CH:5:Hu=0} having~$n - \ell$ partial derivatives bounded on one side.

Then,~$u$ is an affine function.
\end{theorem}

We point out that throughout the chapter a \emph{cone} is any subset~$\C$ of the Euclidean space for which~$\lambda x \in \C$ for every~$x \in \C$ and~$\lambda > 0$. In addition, a \emph{singular} cone is a cone whose boundary is not smooth at the origin or, equivalently, any nontrivial cone that is not a half-space.

Characterizing the values of~$s$ and~$\ell$ for which~\eqref{CH:5:Paellprop} is satisfied represents a challenging open problem, whose solution would lead to fundamental advances in the understanding of the regularity properties enjoyed by nonlocal minimal surfaces. Currently, property~\eqref{CH:5:Paellprop} is known to hold in the following cases:
\begin{itemize}[leftmargin=*]
\item when~$\ell = 1$ or~$\ell = 2$, for every~$s \in (0, 1)$;
\item when~$3\le \ell \le 7$ and~$s \in (1 - \varepsilon_0, 1)$ for some small~$\varepsilon_0 \in (0,1]$ depending only on~$\ell$. 
\end{itemize}
Case~$\ell = 1$ holds by definition, while~$\ell = 2$ is the content of~\cite[Theorem~1]{SV13}. On the other hand, case~$3 \le \ell \le 7$ has been established in~\cite[Theorem~2]{regularity}---see also~\cite{CCS17} for a different approach yielding an explicit value for~$\varepsilon_0$ when~$\ell = 3$.

As a consequence of Theorem~\ref{CH:5:underPakmainthm} and the last remarks, we immediately obtain the following result.

\begin{corollary} \label{CH:5:mainthm}
Let~$n \ge \ell \ge 1$ be integers and~$s \in (0, 1)$. Assume that either
\begin{itemize}[leftmargin=*]
\item $\ell \in \{ 1, 2 \}$, or
\item $3\le \ell \le 7$ and~$s \in (1 - \varepsilon_0, 1)$, with~$\varepsilon_0=\varepsilon_0(\ell) > 0$ as in~\cite[Theorem~2]{regularity}.
\end{itemize}
Let~$u$ be a solution of~\eqref{CH:5:Hu=0} having~$n-\ell$ partial derivatives bounded on one side.

Then,~$u$ is an affine function.
\end{corollary}

We observe that
Theorem~\ref{CH:5:underPakmainthm} gives a new flatness result for~$s$-minimal graphs, under the assumption that~\eqref{CH:5:Paellprop} holds true. It can be seen as a generalization of the fractional~De~Giorgi-type lemma contained in~\cite[Theorem~1.2]{FV17}, which is recovered here taking~$\ell = n$. In this case, we indeed provide an alternative proof of said result.

On the other hand, the choice~$\ell = 2$ gives an improvement of~\cite[Theorem~4]{FarV17}, when specialized to~$s$-minimal graphs. In light of these observations, Theorem~\ref{CH:5:underPakmainthm} and Corollary~\ref{CH:5:mainthm} can be seen as a bridge between Bernstein-type theorems (flatness results in low dimensions) and Moser-type theorems (flatness results under global gradient bounds).

For classical minimal graphs---formally corresponding to the case~$s = 1$ here (see, e.g.,~\cite{Gamma,regularity})---the counterpart of Corollary~\ref{CH:5:mainthm} has been recently obtained by A. Farina in~\cite{F17}. In that case, the result is sharp and holds with~$\ell = \min \{ n, 7 \}$.
See also~\cite{F15} by the same author for a previous result established for~$\ell = 1$ and through a different argument.

Using the same ideas that lead to Theorem~\ref{CH:5:underPakmainthm}, we can prove the following rigidity result for entire~$s$-minimal graphs that lie above a cone.


\begin{theorem}\label{CH:5:growthTeo}
Let~$n\ge1$ be an integer and~$s\in(0,1)$. Let~$u$ be a solution of~\eqref{CH:5:Hu=0} and assume that there exists a constant~$C > 0$ for which
\begin{equation} \label{CH:5:ugecone}
u(x')\ge - C (1+|x'|)\quad\mbox{for every }x'\in\R^n.
\end{equation}
Then,~$u$ is an affine function.
\end{theorem}

Of course, the same conclusion can be drawn if~\eqref{CH:5:ugecone} is replaced by the specular
$$
u(x') \le C (1 + |x'|) \quad \mbox{for every } x' \in \R^n.
$$

For classical minimal graphs, the corresponding version of Theorem~\ref{CH:5:growthTeo} follows at once from the gradient estimate of Bombieri, De Giorgi \& Miranda~\cite{BDM69} and Moser's version of Bernstein's theorem~\cite{M61}. See for instance~\cite[Theorem~17.6]{Giusti} for a clean statement and the details of its proof.

In the nonlocal scenario, a gradient bound for~$s$-minimal graphs has been recently established in~\cite{CaCo}. However, this result is partly weaker than the one of~\cite{BDM69}, since it provides a bound for the gradient of a solution of~\eqref{CH:5:Hu=0} in terms of its oscillation, and not just of its supremum (or infimum) as in~\cite{BDM69}. Consequently, in~\cite{CaCo} a rigidity result analogous to Theorem~\ref{CH:5:growthTeo} is deduced, but with~\eqref{CH:5:ugecone} replaced by the stronger, two-sided assumption:~$|u(x')| \le C(1 + |x'|)$ for every~$x' \in \R^n$. Theorem~\ref{CH:5:growthTeo} thus improves~\cite[Theorem~1.5]{CaCo} directly. Moreover, our proof is different, as it relies on geometric considerations rather than uniform regularity estimates.

Theorem~\ref{CH:5:growthTeo} says in particular that there exist no non-flat~$s$-minimal subgraphs that contain a half-space. Actually, a more general result is true for~$s$-minimal sets that are not necessarily subgraphs, as shown by the following theorem.

\begin{theorem} \label{CH:5:nononflatthm}
Let~$n \ge 1$ be an integer and~$s \in (0, 1)$. If~$E$ is an~$s$-minimal set in~$\R^{n + 1}$ that contains a half-space, then~$E$ is a half-space.
\end{theorem}

Interestingly, Theorem~\ref{CH:5:nononflatthm} can be used to obtain a stronger version of Theorem~\ref{CH:5:growthTeo}, where the bound in~\eqref{CH:5:ugecone} is required to only hold at all points~$x'$ that lie in a half-space of~$\R^n$. See Remark~\ref{CH:5:stressedrmk} at the end of Section~\ref{CH:5:growthsec}.



The proof of Theorem~\ref{CH:5:underPakmainthm} is based on the extension to the fractional framework of a strategy devised by A. Farina for classical minimal graphs and previously unpublished. As a result, the ideas contained in the following sections can be used to obtain a different, easier proof of~\cite[Theorem~1.1]{F17}---since, by Simons' theorem (see, e.g.,~\cite[Theorem~28.10]{Maggi}), no singular classical minimal cones exist in dimension lower or equal to~$7$. Similarly, the same argument that we employ for Theorem~\ref{CH:5:growthTeo} can be successfully applied to classical minimal graphs, giving a different, more geometric, proof of~\cite[Theorem~17.6]{Giusti}.


The argument leading to Theorem~\ref{CH:5:underPakmainthm} relies on a general splitting result for blow-downs of~$s$-minimal graphs. Since it may have an interest on its own, we provide its statement here below.

\begin{theorem} \label{CH:5:singblowdownthm}
Let~$n \ge 1$ be an integer and~$s \in (0, 1)$. Let~$u$ be a solution of~\eqref{CH:5:Hu=0} and~$E$ as in~\eqref{CH:5:Esubgraph}. Assume that~$u$ is not affine and that, for some~$k \in \{ 1, \ldots, n - 1 \}$, the partial derivative~$\frac{\partial u}{\partial x_i}$ is bounded from below in~$\R^n$ for every~$i = 1, \ldots, k$.

Then, every blow-down limit~$\C \subseteq \R^{n + 1}$ of~$E$ is a cylinder of the form
$$
\C = \R^k \times P \times \R,
$$
for some singular~$s$-minimal cone~$P \subseteq \R^{n - k}$.
\end{theorem}

The notion of blow-down limit will be made precise in Section~\ref{CH:5:blowsec}.

\begin{remark}
As revealed by a simple inspection of its proof, Theorem~\ref{CH:5:singblowdownthm} still holds if we require any~$k$ directional derivatives~$\partial_{\nu_1}u,\ldots,\partial_{\nu_k}u$ (not necessarily the partial derivatives) to be bounded from below, provided that the directions~$\nu_1,\ldots,\nu_k$ are linearly independent. Consequently, one can similarly modify the statements of Theorem~\ref{CH:5:underPakmainthm} and Corollary~\ref{CH:5:mainthm} without affecting their validity.
\end{remark}

The remainder of the chapter is structured as follows. In Section~\ref{CH:5:blowsec} we gather some known facts about sets with finite perimeter, the regularity of~$s$-minimal surfaces, and their blow-downs. Section~\ref{CH:5:splittingsec} is devoted to the proof of Theorem~\ref{CH:5:singblowdownthm}, while in Section~\ref{CH:5:mainsec} we show how Theorem~\ref{CH:5:underPakmainthm} follows from it. Sections~\ref{CH:5:nononflatsec} and~\ref{CH:5:growthsec} contain the proofs of Theorems~\ref{CH:5:nononflatthm} and~\ref{CH:5:growthTeo}, respectively. The chapter is closed by Section~\ref{CH:5:constapp}, which includes the extension of a result due to~Chern~\cite{C65} to the framework of graphs having constant~$s$-mean curvature.

\section{Some remarks on nonlocal minimal surfaces and blow-down cones} \label{CH:5:blowsec}

As in the previous chapters, we implicitly assume that all the sets we consider
contain their measure theoretic interior, do not intersect their measure theoretic exterior, and are such that their topological boundary coincides with their measure theoretic boundary---see Remark \ref{CH:1:gmt_assumption} and Appendix \ref{CH:1:Appendix_meas_th_bdary} 
for the details.


We now recall some known results about the regularity of $s$-minimal surfaces, which will be often used without mention in the subsequent sections.

Let~$E\subseteq\R^{n+1}$ be an~$s$-minimal set. Then, its boundary~$\partial E$ is~$n$-rectifiable. Actually, by~\cite[Theorem~2.4]{CRS10},~\cite[Corollary~2]{SV13}, and~\cite[Theorem~1.1]{FV17},~$\partial E$ is locally of class~$C^\infty$, except possibly for a set of singular points~$\Sigma_E\subseteq\partial E$ satisfying
\[
\mathcal H^d(\Sigma_E)=0 \quad \mbox{for every }d>n-2.
\]
In particular, the set~$E$ has locally finite (classical) perimeter in~$\R^{n+1}$ and thus it makes sense to consider its reduced boundary~$\red \! E$.

Furthermore, thanks to the blow-up analysis developed in \cite{CRS10}---see in particular~\cite[Theorem~9.4]{CRS10}---and the tangential properties of the reduced boundary of a set of locally finite perimeter---see, e.g., \cite[Theorem~15.5]{Maggi}---we have that~$\red \! E$ is smooth and the singular set is given by
\[
\Sigma_E=\partial E\setminus\red \!E.
\]

Given a measurable set~$E \subseteq \R^{n + 1}$, a point~$x \in \R^{n + 1}$, and a real number~$r > 0$, we write
$$
E_{x, r} := \frac{E - x}{r}.
$$
We call any~$L^1_\loc$-limit~$E_{x, \infty}$ of~$E_{x, r_j}$ along a diverging sequence~$\{ r_j \}$ a \emph{blow-down limit} of~$E$ at~$x$.

Observe that doing a blow-down of a set~$E$ corresponds to the operation of looking at~$E$ from further and further away. As a result, in the limit one loses track of the point at which the blow-down was centered. That is, blow-down limits may depend on the chosen diverging sequence~$\{ r_j \}$ but not on the point of application~$x$. This fact is certainly well-known to the experts. Nevertheless, we include in the following Remark a brief justification of it for the convenience of the less experienced reader.

\begin{remark}\label{CH:5:tang_cone_infty}
Let~$x, y \in \R^{n + 1}$ and~$E \subseteq \R^{n + 1}$ be a measurable set. Assume that there exists a set~$F \subseteq \R^{n + 1}$ such that~$E_{x, r_j} \rightarrow F$ in~$L^1_\loc(\R^{n + 1})$ as~$j \rightarrow +\infty$, along a diverging sequence~$\{ r_j \}$. We claim that also
\begin{equation} \label{CH:5:E2conv}
E_{y, r_j} \rightarrow F \mbox{ in } L^1_\loc(\R^{n+1}) \mbox{ as } j \rightarrow +\infty.
\end{equation}
To verify this assertion, let~$R > 0$ be fixed and write~$f_j := \chi_{E_{x, r_j}}$ and~$f := \chi_{F}$. Notice that~$\chi_{E_{y, r_j}} = \tau_{v_j} f_j := f_j(\cdot - v_j)$, with~$v_j := (x - y) / r_j$. Since~$v_j \rightarrow 0$ as~$j \rightarrow 0$, we have
\begin{align*}
\left| (E_{y, r_j} \Delta F) \cap B_R \right| & = \| \chi_{E_{y, r_j}} - \chi_F \|_{L^1(B_R)} = \| \tau_{v_j} f_j - f \|_{L^1(B_R)} \\
& \le \| \tau_{v_j} f_j - \tau_{v_j} f \|_{L^1(B_R)} + \| \tau_{v_j} f - f \|_{L^1(B_R)}\\
& \le  \| f_j - f \|_{L^1(B_{R + 1})} + \| \tau_{v_j} f - f \|_{L^1(B_R)},
\end{align*}
provided~$j$ is sufficiently large. Claim~\eqref{CH:5:E2conv} follows since, by assumption,~$f_j \rightarrow f$ in~$L^1_\loc(\R^{n + 1})$ and~$R > 0$ is arbitrary.
\end{remark}

In light of this remark, we can assume blow-downs to be always centered at the origin. For simplicity of notation, we will write~$E_r := E_{0, r} = E/r$ and use~$E_\infty$ to indicate any blow-down limit.

The next lemma collects some known facts about blow-downs of~$s$-minimal sets.

\begin{lemma} \label{CH:5:blowdownlem}
Let~$E \subseteq \R^{n + 1}$ be a nontrivial~$s$-minimal set. Then, for every diverging sequence~$\{ r_j \}$, there exists a subsequence~$\{ r_{j_k} \}$ of~$\{ r_j \}$ and a set~$E_\infty \subseteq \R^{n + 1}$ such that~$E_{r_{j_k}} \rightarrow E_\infty$ in~$L^1_\loc(\R^{n + 1})$ as~$k \rightarrow +\infty$. The set~$E_\infty$ is a nontrivial~$s$-minimal cone. Furthermore,~$E_\infty$ is a half-space if and only if~$E$ is a half-space.
\end{lemma}
\begin{proof}
The existence of a limit of~$E_{r_j}$ (up to a subsequence) is a consequence of 
the fact that~$E_{r}$ is an~$s$-minimal set 
and of Proposition \ref{CH:2:compact_prop} and Remark \ref{CH:2:min_app_seq_rmk}.

The fact that~$E_\infty$ is~$s$-minimal is a consequence of the~$s$-minimality of the sets~$E_{r_{j_k}}$
and their~$L^1_\loc$ convergence to~$E_\infty$---see Corollary \ref{CH:2:local_minima_comp}.

Next we observe that, since~$E$ is nontrivial, we can find a point~$x\in\partial E$. Thanks to
Remark~\ref{CH:5:tang_cone_infty}, we then have that
$$
E_{x,r_{j_k}}\to E_\infty\mbox{ in }L^1_\loc(\R^{n+1})\mbox{ as }k\to\infty.
$$
Since~$0\in\partial E_{x,r_{j_k}}$ for every~$k\in\mathbb N$,
we can conclude that~$E_\infty$ is a cone by arguing as in~\cite[Theorem~9.2]{CRS10}.

The nontriviality of~$E_\infty$ can be established, for instance, by using the uniform density estimates of~\cite{CRS10}. Indeed,~$0 \in \partial E_{x,r_{j_k}}$ for every~$k \in \N$ and hence~\cite[Theorem~4.1]{CRS10} gives that~$\min \{ |E_{x, r_{j_k}} \cap B_1|, |B_1 \setminus E_{x, r_{j_k}}| \} \ge c$ for some constant~$c > 0$ independent of~$k$. As~$E_{x,r_{j_k}} \rightarrow E_\infty$ in~$L^1(B_1)$, it follows that both~$E_\infty$ and its complement have positive measure in~$B_1$. Consequently,~$E_\infty$ is neither the empty set nor the whole~$\R^{n + 1}$.

Finally, if~$E_\infty$ is a half-space, one can deduce the flatness of~$\partial E$ from the~$\varepsilon$-regularity theory of~\cite[Section~6]{CRS10} and the fact that~$\partial E_{r_{j_k}} \rightarrow \partial E_\infty$ in the Hausdorff sense, thanks to the uniform density estimates. See, e.g.,~\cite[Lemma~3.1]{FV17} for more details on this argument.
\end{proof}


%

\section{Proof of Theorem~\ref{CH:5:singblowdownthm}} \label{CH:5:splittingsec}

In this section we include a proof of the splitting result stated in the introduction, namely Theorem~\ref{CH:5:singblowdownthm}. The argument leading to it is based on the following classification result for nonlocal minimal cones that contain their translates. For classical minimal cones, it was proved in~\cite{GMM97}.

\begin{prop} \label{CH:5:GMMprop}
Let~$\C \subseteq \R^{n + 1}$ be an~$s$-minimal cone and assume that
\begin{equation} \label{CH:5:C+vinC}
\C + v \subseteq \C
\end{equation}
for some~$v \in \R^{n + 1} \setminus \{ 0 \}$. Then,~$\C$ is either a half-space or a cylinder in direction~$v$.
\end{prop}
\begin{proof}
First of all, we notice that, since~$\C$ is a cone and inclusion~\eqref{CH:5:C+vinC} holds true, the function~$w := - \nu_\C \cdot v$ satisfies
\begin{equation} \label{CH:5:wge0}
w \ge 0 \quad \mbox{in } \red \C.
\end{equation}
To see this, let~$x \in \red \C$ and observe that,~$\C$ being a cone, we have that~$\mu x \in \overline{\C}$ for every~$\mu > 0$. But then~$\mu x + v \in \overline{\C} + v$ and, using~\eqref{CH:5:C+vinC}, it follows that~$\mu x + v \in \overline{\C}$. Consequently,~$\mu \lambda x + \lambda v = \lambda (\mu x + v) \in \overline{\C}$ for every~$\lambda, \mu > 0$. Choosing~$\mu = 1/\lambda$ we get that~$x + \lambda v \in \overline{\C}$ for every~$\lambda > 0$, which gives that~$v$ points inside~$\overline{\C}$. Recalling that the normal~$\nu_\C$ points outside~$\C$, we are immediately led to~\eqref{CH:5:wge0}.

Now, by~\cite[Theorem~1.3(i)]{CaCo} we know that~$w$ solves
\begin{equation} \label{CH:5:WJharm}
\L w + c^2 w = 0 \quad \mbox{in } \red \C,
\end{equation}
where
\begin{align*}
\L w(x) & := \PV \int_{\red \C} \frac{w(y) - w(x)}{|x - y|^{n + 1 + s}} \, d\Haus^n_y,\\
c^2(x) & := \frac{1}{2} \int_{\red \C} \frac{|\nu_\C(x) - \nu_\C(y)|^2}{|x - y|^{n + 1 + s}} \, d\Haus^n_y,
\end{align*}
for every~$x \in \red \C$. As~$c^2 \ge 0$ in~$\red C$ and~\eqref{CH:5:wge0} holds true, we deduce from~\eqref{CH:5:WJharm} that~$w$ is~$\L$-superharmonic in~$\red \C$, i.e.,
$$
- \L w \ge 0 \quad \mbox{in } \red \C.
$$
By~\cite[Corollary~6.8]{CaCo} (and the lower perimeter bound reported in~\cite[Theorem~3.1]{CaCo}), we then infer that, for every point~$x \in \red \C$ and radius~$R > 0$, the function~$w$ satisfies
$$
\inf_{B_R(x) \cap \red \C} w \ge c_\star R^{1 + s} \int_{\red \C} \frac{w(y)}{(R + |y - x|)^{n + 1 + s}} \, d\Haus^n_y,
$$
for some constant~$c_\star \in (0, 1]$ depending only on~$n$ and~$s$.

Accordingly, either~$w = 0$ in the whole~$\red \C$ or~$\inf_{B_R(x) \cap \red \C} w \ge c_{x, R}$ for some constant~$c_{x, R} > 0$ and for every~$x \in \red \C$ and~$R > 0$. In the first case, it is easy to see that~$\C$ must be a cylinder in direction~$v$. If the second situation occurs, then~$\partial \C$ is a locally Lipschitz graph with respect to the direction~$v$ (see, e.g.,~\cite[Theorem~5.6]{M64}), and hence smooth, due to~\cite[Theorem~1.1]{FV17}. It being a cone, we conclude that~$\C$ must be a half-space.
\end{proof}

With this in hand, we may now proceed to prove the splitting result.

\begin{proof}[Proof of Theorem~\ref{CH:5:singblowdownthm}]
Let~$E$ denote the subgraph of~$u$, as defined by~\eqref{CH:5:Esubgraph}. We recall that, as observed right before the statement of Theorem~\ref{CH:5:underPakmainthm}, the set~$E$ is~$s$-minimal.

Let~$\C$ be a blow-down cone of~$E$. By definition, there exists a diverging sequence~${r_j}$ for which~$E_{r_j} = E/r_j \rightarrow \C$ in~$L^1_\loc(\R^{n + 1})$. As noticed in Lemma~\ref{CH:5:blowdownlem},~$\C$ is a nontrivial~$s$-minimal cone. Moreover,~$\C$ is not an half-space, since, otherwise,~$E$ would be a half-space too (again, by Lemma~\ref{CH:5:blowdownlem}),
contradicting the hypothesis that~$E$ is the subgraph of a non-affine function. We also recall that this is equivalent to the cone~$\C$ being singular.

As~$E$ is a subgraph, it follows that~$E - t e_{n + 1} \subseteq E$ for every~$t > 0$. This yields that~$E_{r_j} - e_{n + 1} \subseteq E_{r_j}$ for every~$j$. Hence, by~$L^1_\loc(\R^{n + 1})$ convergence,~$\C - e_{n + 1} \subseteq \C$. Since~$\C$ is not a half-space, by Proposition~\ref{CH:5:GMMprop} we conclude that~$\C$ is a cylinder in direction~$e_{n + 1}$, that is
\begin{equation} \label{CH:5:Cen+1cyl}
\C + \lambda e_{n + 1} = \C \quad \mbox{for every } \lambda \in \R,
\end{equation}
or, equivalently,~$\C = \C' \times \R$, for some singular~$s$-minimal cone~$\C' \subseteq \R^n$. Observe that the~$s$-minimality of~$\C'$ is a consequence of~\cite[Theorem~10.1]{CRS10}. Also note that to obtain~\eqref{CH:5:Cen+1cyl} we only took advantage of the fact that~$E$ is an~$s$-minimal subgraph and not the hypotheses on the partial derivatives of~$u$.

Let now~$i = 1, \ldots, k$ be fixed. By the bound from below on the partial derivative~$\frac{\partial u}{\partial x_i}$ and the fundamental theorem of calculus, there exists a constant~$\kappa > 0$ such that
$$
u(z' + t e_i) - u(z') = \int_{0}^{t} \frac{\partial u(z' + \tau e_i)}{\partial x_i} \, d\tau \ge - \kappa t
$$
for every~$z' \in \R^n$ and~$t > 0$. Let now~$u_j$ be the function defining the blown-down set~$E_{r_j}$. Clearly,~$u_j(z') = u(r_j z') / r_j$ and hence
$$
u_j(y' + e_i) - u_j(y') = \frac{u(r_j y' + r_j e_i) - u(r_j y')}{r_j} \ge - \kappa
$$
for every~$y' \in \R^n$ and~$j \in \N$. This means that~$E_j - \kappa e_{n + 1} + e_i \subseteq E_j$ for every~$j \ge 1$. Passing to the limit and using~\eqref{CH:5:Cen+1cyl}, we deduce that~$\C + e_i = \C - \kappa e_{n + 1} + e_i \subseteq \C$. Taking advantage once again of Proposition~\ref{CH:5:GMMprop} and of the fact that~$\C$ is not a half-space, we infer that~$\C$ is a cylinder in direction~$e_i$ for every~$i = 1, \ldots, k$. The conclusion of the theorem follows.
\end{proof}

\section{Proof of Theorem~\ref{CH:5:underPakmainthm}} \label{CH:5:mainsec}

First of all, we may assume that the partial derivatives of~$u$ bounded on one side are the first~$n - \ell$. Also, up to flipping the variable~$x_i$, for some~$i \in \{ 1, \ldots, n - \ell \}$, we may suppose that those partial derivatives are all bounded from below. All in all, we have that
$$
\frac{\partial u}{\partial x_i} \ge - \kappa \quad \mbox{for every } i = 1, \ldots, n - \ell,
$$
for some constant~$\kappa \ge 0$.

If~$u$ were not affine, then, by applying Theorem~\ref{CH:5:singblowdownthm} with~$k = n - \ell$, we would have that every blow-down cone~$\C$ of the set~$E$ defined by~\eqref{CH:5:Esubgraph} is given by
$$
\C = \R^k \times P \times \R,
$$
for some singular~$s$-minimal cone~$P \subseteq \R^{n - k} = \R^{\ell}$. As this contradicts assumption~\eqref{CH:5:Paellprop}, we conclude that~$u$ must be affine.

\section{Proof of Theorem~\ref{CH:5:nononflatthm}} \label{CH:5:nononflatsec}

Let~$\Pi$ be a half-space contained in~$E$. Without loss of generality, we may assume that~$\Pi = \{ x \in \R^n \,|\, x_{n + 1} < 0 \}$. Consider then a blow-down~$\C$ of~$E$, which is a nontrivial~$s$-minimal cone, by Lemma~\ref{CH:5:blowdownlem}. In particular, $\Pi \subseteq \C$ and~$0 \in \partial \Pi \cap \partial \C$. Using, e.g.,~\cite[Corollary~6.2]{CRS10}, we infer that~$\C = \Pi$ and therefore that~$E$ is half-space as well, thanks again to Lemma~\ref{CH:5:blowdownlem}.

\section{Proof of Theorem~\ref{CH:5:growthTeo}} \label{CH:5:growthsec}

Suppose by contradiction that the function~$u$ is not affine and denote with~$E$ its subgraph. Up to a translation of~$E$ in the vertical direction, hypothesis~\eqref{CH:5:ugecone} yields that~$E$ contains the cone
$$
\mathscr{D} := \left\{ x \in \R^{n + 1} \,|\, x_{n + 1} < - C |x'| \right\}.
$$
Consider now a blow-down~$\C$ of~$E$. On the one hand, we clearly have that~$\mathscr{D} \subseteq \C$. On the other hand, by arguing as in the beginning of the proof of Theorem~\ref{CH:5:singblowdownthm}, we have that~$\C$ must be a nontrivial vertical cylinder. More precisely,~$\C=\C'\times\R$, for some nontrivial singular~$s$-minimal cone~$\C'\subseteq\R^n$. These two facts imply that~$\C' = \R^n$, contradicting its nontriviality. This concludes the proof.

\begin{remark} \label{CH:5:stressedrmk}
By a refinement of this argument we can prove a stronger version of Theorem~\ref{CH:5:growthTeo}, where hypothesis~\eqref{CH:5:ugecone} is replaced by
\begin{equation} \label{CH:5:ugecone+}
u(x') \ge - C(1 + |x'|) \quad \mbox{for every } x' \in \R^n \mbox{ such that } x_1 < 0.
\end{equation}
Indeed, arguing by contradiction as before, we see that any blow-down of the subgraph of~$u$ is a cylinder of the form~$\C' \times \R$. In light of~\eqref{CH:5:ugecone+}, the cone~$\C'$ contains a half-space of~$\R^n$ and is thus flat, due to Theorem~\ref{CH:5:nononflatthm}. This leads to a contradiction.
\end{remark}

%

\section{Subgraphs of constant fractional mean curvature} \label{CH:5:constapp}


We recall---see Chapter \ref{CH_Nonparametric}---that, given a measurable function $u:\R^n\to\R$, we can understand~$\cu_s u$ as a linear form on the fractional Sobolev space~$W^{s, 1}(\R^n)$, setting
\[
\langle\cu_s u, v\rangle:=
\int_{\R^n}\int_{\R^n} G_s \! \left( \frac{u(x')-u(y')}{|x'-y'|} \right)
\! \left( v(x')-v(y') \right) \frac{dx' dy'}{|x'-y'|^{n+s}}
\]
for every~$v\in W^{s,1}(\R^n)$. This definition is indeed well-posed since~$G_s$ is bounded.

Let~$h$ be a real number. We say that a measurable function~$u:\R^n\to\R$ is a weak solution of~$\cu_s u = h$ in~$\R^n$ if it holds
\[
\langle\cu_s u,v\rangle=h\int_{\R^n}v(x')\,dx' \quad \mbox{for every } v\in W^{s,1}(\R^n).
\]
We remark that by the density of~$C^\infty_c(\R^n)$ in~$ W^{s,1}(\R^n)$,
it is equivalent to consider the test functions~$v$ to be smooth and compactly supported.

We now prove that if the~$s$-mean curvature of a global subgraph is constant, then this constant must be zero. More precisely, we have the following statement.

\begin{prop} \label{CH:5:c=0lem}
Let~$u:\R^n\to\R$ be a weak solution of~$\cu_s u=h$ in~$\R^n$, for some constant~$h \in \R$. Then~$h=0$.
\end{prop}

\begin{proof}
Recalling~\eqref{CH:5:Gdef}, we notice that
\[
|G_s(t)|\le \int_0^{+\infty}\frac{d\tau}{(1+\tau^2)^\frac{n+1+s}{2}}=\frac{\Lambda_{n,s}}{2}<+\infty \quad \mbox{for every } t\in\R.
\]
Suppose that~$h\ge0$---the case~$h\le0$ is analogous. Let~$R>0$ and consider the test function~$v=\chi_{B'_R}\in W^{s,1}(\R^n)$. We have
\[
|\langle\cu_s u,\chi_{B'_R}\rangle|\le\Lambda_{n,s}\int_{B'_R}\int_{\R^n\setminus B'_R}\frac{dx' dy'}{|x'-y'|^{n+s}}
= C R^{n-s},
\]
for some constant~$C>0$ depending only on~$n$ an~$s$. Since~$u$ weakly solves~$\cu_s = h$ in~$\R^n$, we deduce that
\[
h |B'_1| R^n=h\int_{\R^n}\chi_{B'_R}(x')\,dx'=\langle\cu_s u,\chi_{B'_R}\rangle\le C R^{n-s}
\]
for all~$R>0$, that is~$0 \le h R^s \le C / |B_1'|$. Letting~$R\to+\infty$ we conclude that~$h=0$.
\end{proof}

We point out that, as a consequence of Proposition~\ref{CH:5:c=0lem} and the results of Corollary \ref{CH:4:equiv_intro_global_corollary},
if a function~$u\in W^{s, 1}_\loc(\R^n)$ is a weak solution of~$\cu_s u = h$ in~$\R^n$, then the subgraph of~$u$ must be an~$s$-minimal set---thus extending to the nonlocal framework a celebrated result of Chern, namely the Corollary of Theorem~1 in~\cite{C65}.
 
We further remark that other definitions for solutions of the equation~$\cu_s u=h$ could have been considered, namely smooth
pointwise solutions and viscosity solutions (for a rigorous definition see Definition \ref{CH:4:VISCO_DEF}). However, it is readily seen that a smooth pointwise solution is also a viscosity solution. Moreover, Corollary \ref{CH:4:Visco_GLOB_CoRoLl} shows that
a viscosity solution is also a weak solution. Consequently, Proposition~\ref{CH:5:c=0lem} applies to these other two notions of solutions as well.

\end{chapter}

\begin{chapter}[A free boundary problem]{A free boundary problem: superposition
of nonlocal energy plus classical perimeter}\label{CH_FreeBdary_CHPT}

\minitoc

In this chapter we study the minimizers of the functional
\bgs{
	\Nl(u,\Omega)+\Per\big(\{u>0\},\Omega\big),
}
with $\Nl(u,\Omega)$ being, roughly speaking, the $\Omega$-contribution to the $H^s$ seminorm of a function $u:\R^n\to\R$

\smallskip

The main contributions of the present chapter
consist in establishing a monotonicity formula for the minimizers, in exploiting it to investigate the properties of blow-up limits and in proving a dimension reduction result. Moreover, we show that, when $s<1/2$, the perimeter dominates---in some sense---over the nonlocal energy. As a consequence, we obtain a regularity result for the free boundary $\{u=0\}$.

\section{Introduction: definitions and main results}

Let us begin by giving the rigorous definition of the functional that we are going to study.

Given $s\in(0,1)$ and a bounded open set $\Omega\subseteq\R^n$
with Lipschitz boundary,
we consider the functional
\begin{equation}\label{CH:6:functional}
\F_\Omega(u,E):=\iint_{\R^{2n}\setminus(\Co\Omega)^2}
\frac{|u(x)-u(y)|^2}{|x-y|^{n+2s}}\,dx\,dy+\Per(E,\Omega),
\end{equation}
where $E$ is the positivity set of the function $u:\R^n\to\R$, that is
\begin{equation*}
u\geq0\quad\textrm{a.e. in } E\qquad\textrm{and}\qquad 
u\leq0\quad\textrm{a.e. in }\Co E.
\end{equation*}
We call such a pair $(u,E)$ an {\emph{admissible pair}}.
Here above~$\Co E$ denotes the complement of~$E$ 
and $\Per(E,\Omega)$ denotes the (classical) 
perimeter of $E$ in $\Omega$.

Furthermore, we write
\begin{equation}\begin{split}\label{CH:6:def_N}
\Nl(u,\Omega)&:=\iint_{\R^{2n}\setminus(\Co\Omega)^2}
\frac{|u(x)-u(y)|^2}{|x-y|^{n+2s}}\,dx\,dy\\
&
=\iint_{\Omega\times\Omega}\frac{|u(x)-u(y)|^2}{|x-y|^{n+2s}}\,dx\,dy
+2\iint_{\Omega\times\Co\Omega}
\frac{|u(x)-u(y)|^2}{|x-y|^{n+2s}}\,dx\,dy,
\end{split}\end{equation}
for the nonlocal energy of $u$ appearing in the definition of $\F_\Omega$. Roughly speaking, this is the $\Omega$-contribution to the $H^s$ seminorm of $u$.

We will consider the following definition of minimizing pair.

\begin{defn}\label{CH:6:min_pair}
Given an admissible pair~$(u,E)$, we say that
a pair $(v,F)$ 
is an admissible competitor 
(for $\F_\Omega$ with respect to the pair $(u,E)$) if
\begin{equation}\begin{split}\label{CH:6:minim_def}
&\textrm{supp}(v-u)\Subset\Omega,\qquad F\Delta E\Subset\Omega,\\
&v-u\in H^s(\R^n)\qquad\textrm{and}\qquad \Per(F,\Omega)<+\infty.
\end{split}\end{equation}

We say that the admissible pair $(u,E)$
is minimizing in $\Omega$ if $\F_\Omega(u,E)<+\infty$ and
\begin{equation*}
\F_\Omega(u,E)\leq\F_\Omega(v,F),
\end{equation*}
for every admissible competitor~$(v,F)$.
\end{defn}


We observe that in Proposition \ref{CH:6:equiv_char} we will provide some equivalent characterizations of minimizing pairs.

\smallskip

In particular, we are interested in the following minimization problem, with respect to fixed ``exterior data''.
Given an admissible pair $(u_0,E_0)$ and a bounded open set $\Op\subseteq\R^n$ with Lipschitz boundary, such that
\begin{equation*}
\Omega\Subset\Op,\qquad \Nl(u_0,\Omega)<+\infty
\quad\textrm{ and }\quad \Per(E_0,\Op)<+\infty,
\end{equation*}
we want to find an admissible pair $(u,E)$ attaining the following infimum
\eqlab{\label{CH:6:Min_dir}
\inf\big\{\Nl(v,\Omega)+\Per(F,\Op)\,|\,(v,F)&\textrm{ admissible s.t. }v=u_0\textrm{ a.e. in }\Co\Omega\\
&\qquad\qquad
\textrm{ and }F\setminus\Omega=E_0\setminus\Omega\big\}.
}
Roughly speaking, as customary when dealing with minimization problems involving the classical perimeter, we are considering a (fixed) neighborhood $\Op$ of $\Omega$ (as small as we like)
in order to ``read'' the boundary data $\partial E_0\cap\partial\Omega$.

In Section~\ref{CH:6:sec:prelim} we prove the existence of pairs solving this Dirichlet problem. Moreover, we show
that a pair $(u,E)$ realizing the infimum in \eqref{CH:6:Min_dir}
is also a minimizing pair in the sense of Definition \ref{CH:6:min_pair}.


\smallskip

Concerning the minimizers of the functional $\F$, we also establish the following uniform energy estimates, which turn out to be important when proving the existence of blow-up limits.

\begin{theorem}\label{CH:6:TH:unif}
	Let~$(u,E)$ be a minimizing pair in~$B_2$. 
	Then
	\[ \iint_{\R^{2n}\setminus (\Co B_1)^2}\frac{|u(x)-u(y)|^2}{
		|x-y|^{n+2s}}\,dx\,dy + \Per(E,B_1)\le 
	C\left(1+\int_{\R^n}\frac{|u(y)|^2}{1+|y|^{n+2s}}\,dy\right),
	\]
	for some~$C=C(n,s)>0$.
\end{theorem}

\medskip

In order to study blow-up sequences, we will need a ``localized'' version of $\Nl(\,\cdot\,,\Omega)$ which is obtained through an extension technique studied
in \cite{CS07}. To be more precise, 
given a function~$u:\R^n\to\R$,
we consider the function~$\Ue:\R_+^{n+1}\to\R$,
where
\begin{equation*}
\R^{n+1}_+:=\{(x,z)\in\R^{n+1} \,|\, x\in\R^n,\,z>0\},
\end{equation*}
defined via the convolution with an appropriate Poisson kernel,
\begin{equation}\label{CH:6:extended_func_def}
\Ue(\,\cdot\,,z)=u\ast\mathcal K_s(\,\cdot\,,z),\quad\textrm{where}\quad\K_s(x,z):=c_{n,s}\frac{z^{2s}}{(|x|^2+z^2)^{(n+2s)/2}}.
\end{equation}
Here above,~$c_{n,s}>0$ is an appropriate normalizing constant. We observe that the extended function $\Ue$ is well defined, provided the function $u$ belongs to the weighted Lebesgue space
\[
L_s(\R^n):=\left\{u:\R^n\to\R\,\big|\,\int_{\R^n}\frac{|u(\xi)|}{1+|\xi|^{n+2s}}d\xi<+\infty\right\}.
\]
For a proof of this fact and for a detailed introduction to the extension operator, we refer the interested reader to \cite{Extension}.
In light of Remark \ref{CH:6:tail_energies_rmk}, we can thus consider the extended function of a minimizer.

We use capital letters, like $X=(x,z)$, to
denote points in $\R^{n+1}$.
Given a set $\Omega\subseteq\R^{n+1}$, we write
\begin{equation*}
\Omega_+:=\Omega\cap\{z>0\}\qquad\textrm{and}\qquad\Omega_0:=\Omega\cap\{z=0\}.
\end{equation*}
Moreover we identify the hyperplane $\{z=0\}\simeq\R^n$ via the projection function.

\smallskip

In particular, we exploit the energy naturally associated to the extension problem to define an extended functional,
which has a local behavior.

To be more precise, given a bounded open set $\Omega\subseteq\R^{n+1}$ 
with Lipschitz boundary, such that $\Omega_0\not=\emptyset$, we define
\begin{equation}\label{CH:6:extended_functional}
\lf_\Omega(\Vf,F):=c_{n,s}'\int_{\Omega_+}|\nabla\Vf|^2z^{1-2s}\,dX+\Per(F,\Omega_0),
\end{equation}
for $\Vf:\R^{n+1}_+\to\R$ and $F\subseteq\R^n\simeq\{z=0\}$
the positivity set of the trace of $\Vf$ on $\{z=0\}$, that is
\begin{equation*}
\Vf\big|_{\{z=0\}}\geq0\quad\textrm{a.e. in }F\quad\textrm{and}\quad
\Vf\big|_{\{z=0\}}\leq0\quad\textrm{a.e. in }\Co F.
\end{equation*}
We call such a pair $(\Vf,F)$ an {\emph{admissible pair}}
for the extended functional.

From now on, whenever considering the extended functional, unless otherwise stated we will implicitly assume that the open set $\Omega\subseteq\R^{n+1}$ is such that $\Omega_0\not=\emptyset$.

\begin{defn}\label{CH:6:adm-ext}
Given an admissible pair~$(\Uc,E)$ such that $\lf_\Omega(\Uc,E)<+\infty$, we say that a pair~$(\Vf,F)$
is an admissible competitor
(for $\lf_\Omega$ with respect to $(\Uc,E)$)
if~$\lf_\Omega(\Vf,F)<+\infty$ and
\begin{equation*}
\textrm{supp}\,(\Vf-\Uc)\Subset\Omega\qquad
\textrm{and}\qquad E\Delta F\Subset\Omega_0.
\end{equation*}

We say that an admssible pair $(\Uc,E)$ is minimal in~$\Omega$
if~$\lf_\Omega(\Uc,E)<+\infty$ and
\begin{equation*}
\lf_\Omega(\Uc,E)\leq\lf_\Omega(\Vf,F),
\end{equation*}
for every admissible competitor $(\Vf,F)$.
\end{defn}

We will study this extended functional in Section~\ref{CH:6:sec:ext}.
In particular, we relate minimizers of the extended functional $\lf$ with minimizers of the original functional $\F$, proving the following:

\begin{prop}\label{CH:6:Local_energy_prop}
	Let $(u,E)$ be an admissible pair for $\F$,
	according to Definition~\ref{CH:6:min_pair}, such
	that~$\F_{B_R}(u.E)<+\infty$.
	Then, the pair $(u,E)$ is minimizing in $B_R$ if and only 
	if the pair $(\Ue,E)$ is minimizing for $\lf_\Omega$,
	for every bounded open set $\Omega\subseteq\R^{n+1}$ 
	with Lipschitz boundary and such that~$\Omega_0\Subset B_R$.
\end{prop}

We now introduce the following notation
\begin{equation*}
\BaLL_r:=\{(x,z)\in\R^{n+1}\,|\,|x|^2+z^2<r^2\},\qquad\qquad
\BaLL_r^+:=\BaLL_r\cap\{z>0\},
\end{equation*}
and
\begin{equation*}
(\partial\BaLL_r)^+:=\partial\BaLL_r\cap\{z>0\}=
\{(x,z)\in\R^{n+1}_+\,|\,|x|^2+z^2=r^2\}.
\end{equation*}

The main reason for considering the extended functional consists in the fact that it allows us to obtain a Weiss-type monotonicity formula---by exploiting a scaled and ``corrected'' version of
the functional $\lf_{\BaLL_r}$.
More precisely:

\begin{theorem}[Weiss-type Monotonicity Formula]\label{CH:6:Monotonicity_teo}
	Let $(u,E)$ be a minimizing pair for $\F$ in $B_R$ and define the function $\Phi_u:(0,R)\to\R$ by
	\bgs{
		\Phi_u(r):=r^{1-n}&
		\left(c'_{n,s}\int_{\BaLL_r^+}|\nabla\Ue|^2z^{1-2s}\,dX+\Per(E,B_r)\right)\\
		&\qquad
		-c'_{n,s}\Big(s-\frac{1}{2}\Big)r^{-n}\int_{(\partial\BaLL_r)^+}\Ue^2z^{1-2s}\,d\Ha^n.
	}
	Then, the function $\Phi_u$ is increasing in $(0,R)$.	
	Moreover, $\Phi_u$ is constant in $(0,R)$ if and only if the
	extension $\Ue$ is homogeneous of degree $s-\frac{1}{2}$ in $\BaLL_R^+$
	and $E$ is a cone in $B_R$.
\end{theorem}

In order to prove the monotonicity formula, we will need
to construct appropriate competitors for the minimizing pair $(\Ue,E)$ of the extended functional.
For this, we need to consider the cone $E(r)$ spanned 
by the ``spherical slice'' $E\cap\partial B_r$, namely
\begin{equation}\label{CH:6:3.4bis}
E(r):=\{\lambda y\,|\,\lambda>0,\,y\in E\cap\partial B_r\}.
\end{equation}
In Section~\ref{CH:6:appA}, we show that this cone is indeed well defined for
a.e. $r>0$ and its perimeter in every ball $B_{\varrho}$ can be computed by means of a simple formula
(see Proposition \ref{CH:6:cones_from_slices}).
We mention that for the proof of Theorem \ref{CH:6:Monotonicity_teo}---which is in Section \ref{CH:6:MoNoToFoRSeC}---we will also need a result concerning the surface density of a Caccioppoli set, namely Corollary \ref{CH:6:guy_coroll}.

In order to study blow-up sequences, we
prove a general convergence result for minimizing pairs under appropriate conditions. More precisely:

\begin{theorem}[Proof in Section \ref{CH:6:CoNVoMinSeC}] \label{CH:6:TH:convmin}
	Let~$(\Ue_m,E_m)$ be a sequence of minimizing pairs
	in~$\BaLL_R^+$. Suppose that~$\Ue_m$ is the extension
	of~$u_m$, and
	\[ u_m\to u {\mbox{ in }}L^\infty(B_R), \quad 
	\Ue_m\to \Ue {\mbox{ in }}L^\infty(\BaLL_R^+),\quad
	{\mbox{ and }}\quad  \big|(E_m \Delta E)\cap B_R\big|\to0 
	\] 
	as~$m\to+\infty$, for some admissible pair~$(\Ue,E)$, 
	with~$\Ue$ continuous in~$\overline{\R^{n+1}_+}$, being
	$\Ue$ the extension function of $u$.
	Then, $(\Ue,E)$ is a minimizing pair in~$\BaLL_r^+$, for every $r\in(0,R)$.\\	
	Furthemore,
	\begin{equation}\label{CH:6:convmin1}
	\lim_{m\to+\infty} \int_{\BaLL_r^+}|\nabla \Ue_m|^2
	z^{1-2s}\,dX =  \int_{\BaLL_r^+}|\nabla \Ue|^2
	z^{1-2s}\,dX, \quad \forall\,r\in(0,R),
	\end{equation}
	and
	\begin{equation}\label{CH:6:convmin2}
	D\chi_{E_m}\stackrel{\ast}{\rightharpoonup}D\chi_E
	\quad\mbox{ and }\quad
	\big|D\chi_{E_m}\big|\stackrel{\ast}{\rightharpoonup}|D\chi_E|,\qquad\mbox{in }B_R.
	\end{equation}
	In particular,
	\begin{equation}\label{CH:6:convmin3}
	\lim_{m\to+\infty} \Per(E_m, B_r)= 
	\Per(E,B_r),
	\end{equation}
	for every $r\in(0,R)$ such that
	\[\Ha^{n-1}(\partial^*E\cap\partial B_r)=0.\]
\end{theorem}

Exploiting the results that we have mentioned so far, we are able to study blow-up limits. Let us first introduce some notation.

Given a function $u:\R^n\to\R$ and a set $E\subseteq\R^n$, we define
\begin{equation}\label{CH:6:rescaled_stuff}
u_\lambda(x):=\lambda^{\frac{1}{2}-s}u(\lambda x)\qquad\textrm{and}\qquad E_\lambda:=\frac{1}{\lambda}E,
\end{equation}
for every $\lambda>0$. We observe that the scaling introduced in \eqref{CH:6:rescaled_stuff} is consistent with the natural scaling of the functionals that we are considering---see Remark \ref{CH:6:scaling_minimality}.

Given a minimizing pair $(u,E)$, we are interested in the blow-up sequence, that is the sequence of pairs $(u_r,E_r)$ for $r\to0$. We observe that, as a consequence of the natural scaling of the functionals and of the monotonicity formula, blow-up limits possess homogeneity properties.

We thus introduce the following notion.
We say that the admissible pair~$(u,E)$ is a {\emph{minimizing cone}}
if it is a minimizing pair in $B_R$, for every $R>0$, and it is such that~$u$ is homogeneous 
of degree~$s-\frac12$ and~$E$ is a cone 
(that is, $\chi_E$ is homogeneous of degree~$0$).

With this, we can now state the following result:

\begin{theorem}[Proof in Section \ref{CH:6:BLOWupSeQSec}]\label{CH:6:TH:blow}
	Let $s>1/2$ and~$(u,E)$ be a minimizing pair in~$B_1$
	with~$0\in\partial E$. Let~$(u_r,E_r)$ be
	as in~\eqref{CH:6:rescaled_stuff}. Assume that $u\in C^{s-\frac12}(B_1)$.
	Then, there exist a minimizing cone~$(u_0,E_0)$
	and a sequence~$r_k\searrow 0$ such that~$u_{r_k}\to u_0$
	in~$L^\infty_{{\loc}}(\R^n)$ and~$E_{r_k}\xrightarrow{\loc} E_0$.
\end{theorem}

We point out that the assumption $u\in C^{s-\frac12}(B_1)$
in Theorem~\ref{CH:6:TH:blow} is clearly weaker than asking~$u$
to be~$C^{s-\frac12}$ in the whole of~$\R^n$, 
which is the requirement of~\cite[Theorem~1.3]{DSV}. In particular, in Theorem~\ref{CH:6:TH:blow} we are not even requiring~$u$
to be continuous outside~$B_1$. 

\smallskip

In Section \ref{CH:6:Regs<12free} we observe that in the case $s<1/2$ the perimeter is, in some sense, the leading term of the functional $\F_\Omega$. As a consequence, we are able to prove the following regularity result for the free boundary $\partial E$:

\begin{theorem}\label{CH:6:THMREGFREEBDARYS12}
	Let $s\in(0,1/2)$ and let $(u,E)$ be a minimizing pair in $\Omega$. Suppose that $u\in L^\infty_{\loc}(\Omega)$. 
	Then $E$ has almost minimal boundary in $\Omega$.
	
	More precisely, if $x_0\in\Omega$ and $d:=d(x_0,\Omega)/3$,
	then for every $r\in(0,d]$ it holds
	\begin{equation}\label{CH:6:almost_min}
	\Per(E,B_r(x_0))\le \Per(F,B_r(x_0))+ C\,r^{n-2s},\qquad\forall\,F\subseteq\R^n\mbox{ s.t. }E\Delta F\Subset B_r(x_0),
	\end{equation}
	where
	\[
	C=C\left(s,x_0,d,\|u\|_{L^\infty(B_{2d}(x_0))},\int_{\R^n}\frac{|u(y)|}{1+|y|^{n+2s}}\,dy\right)>0.
	\]
	Therefore
	\begin{itemize}
		\item[(i)] $\partial^*E$ is locally $C^{1,\frac{1-2s}{2}}$ in $\Omega$,
		
		\item[(ii)] the singular set $\partial E\setminus\partial^*E$ is such that
		\[
		\Ha^\sigma\big((\partial E\setminus\partial^*E)\cap\Omega\big)=0,\qquad\mbox{for every }\sigma>n-8.
		\]
	\end{itemize}
\end{theorem}

\smallskip

We conclude this Introduction by mentioning the following dimension reduction result for global minimizers.

Only in the following Theorem and in Section \ref{CH:6:DIMREDUCTIONSECTioN} we redefine
\[\F_\Omega(u,E):=(c_{n,s}')^{-1}\Nl(u,\Omega)+\Per(E,\Omega),\]
so that the corresponding extended functional is constant-free.

We say that an admissible pair $(u,E)$ is minimizing in $\R^n$ if
it minimizes $\F_\Omega$ in any bounded open subset $\Omega\subseteq\R^n$ (in
the sense of Definition \ref{CH:6:min_pair}).

\begin{theorem}\label{CH:6:DimREDtHM}
	Let $(u,E)$ be an admissible pair and define
	\begin{equation*}
	u^\star(x,x_{n+1}):=u(x)\qquad\mbox{and}\qquad
	E^\star:=E\times\R.
	\end{equation*}
	Then, the pair $(u,E)$ is minimizing in $\R^n$
	if and only if the pair $(u^\star,E^\star)$ is minimizing in $\R^{n+1}$.
\end{theorem}



%
%
%

\subsection{Notation and assumptions}

Throughout the chapter~$\Omega$ will be a bounded open set
with Lipschitz boundary, unless otherwise stated.

Like we did in the previous chapters, we will make the following assumption regarding the sets that we consider.

\subsubsection{Measure theoretic assumption}\label{CH:6:mta_subsec}

Let $F\subseteq\R^n$. Up to modifications in sets of measure zero, we can assume that $F$ coincides with
the set $F^{(1)}$ of points of density 1, which is a ``good representative'' for $F$ in its $L^1_{\loc}$ class.
In particular, we can thus assume that
$F$ contains its measure theoretic interior
\[F_{int}:=\{x\in\R^n\,|\,\exists\,r>0\textrm{ s.t. }|F\cap B_r(x)|=\omega_nr^n\}\subseteq F,\]
the complementary $\Co F$ contains its measure theoretic interior,
\[F_{ext}:=\{x\in\R^n\,|\,\exists\,r>0\textrm{ s.t. }|F\cap B_r(x)|=0\}\subseteq\Co F,\]
and the topological boundary of $F$ coincides with the measure theoretic boundary, $\partial F=\partial^-F$,
where
\begin{equation}
\partial^-F:=\R^n\setminus\big(F_{int}\cup F_{ext}\big)
=\{x\in\R^n\,|\,0<|F\cap B_r(x)|<\omega_nr^n\;\forall\,r>0\}.
\end{equation}
For the details, we refer to Appendix \ref{CH:1:Appendix_meas_th_bdary} and Section \ref{CH:6:appA}.

\section{Preliminary results}\label{CH:6:sec:prelim}

In this section we will prove some basic properties,
such as the existence of a minimizing pair $(u,E)$ for
the functional~$\F$
(using the direct method of Calculus of Variations) and
the $s$-harmonicity of the function $u$. We also establish
a comparison principle for minimizers.
Finally, we show that if~$(u,E)$ is minimizing 
in~$\Omega$, then it is minimizing in 
every~$\Omega'\Subset\Omega$.

\smallskip

We first point out the following useful remarks about the ``tail energies''.
Given $s\in(0,1)$ we define the weighted Lebesgue space
\[
L_s^2(\R^n):=\left\{u:\R^n\to\R\,\big|\,\int_{\R^n}\frac{|u(\xi)|^2}{1+|\xi|^{n+2s}}d\xi<+\infty\right\}.
\]

\begin{remark}\label{CH:6:tail_energies_rmk}
We observe that we have the continuous embedding
\[
L^2_s(\R^n)\subseteq L_s(\R^n).
\]
Indeed, if $u\in L^2_s(\R^n)$, then by Holder's inequality we have
\begin{equation*}\begin{split}
\int_{\R^n}\frac{|u(y)|}{1+|y|^{n+2s}}\,dy&=
\int_{\R^n}\frac{|u(y)|}{\big(1+|y|^{n+2s}\big)^\frac{1}{2}}\,\frac{dy}{\big(1+|y|^{n+2s}\big)^\frac{1}{2}}\\
&
\le\left(\int_{\R^n}\frac{|u(y)|^2}{1+|y|^{n+2s}}\,dy\right)^\frac{1}{2}
\left(\int_{\R^n}\frac{dy}{1+|y|^{n+2s}}\right)^\frac{1}{2}<+\infty.
\end{split}\end{equation*}
Moreover, it trivially holds true that
\[
L^2_s(\R^n)\subseteq L^2_{\loc}(\R^n).
\]
Finally, we point out that, if $\Omega\subseteq\R^n$ is a bounded open set and
if $u:\R^n\to\R$ is a measurable function, then
\begin{equation*}
\Nl(u,\Omega)<+\infty\quad\implies\quad u\in L^2_s(\R^n).
\end{equation*}
For the proof of this observation we refer, e.g., to Lemma \ref{CH:APP:usef_ineq_tail}.
\end{remark}

\subsection{Existence of a minimizing pair for the Dirichlet problem and $s$-harmonicity}

We begin by observing that, even if the choice of the 
neighborhood~$\Op\Supset\Omega$ for 
the Dirichlet problem is
arbitrary, it does not influence the minimization 
problem (provided that the positivity set of the exterior data is
regular enough).

\begin{remark}
Let $\Omega\Subset\Op'\Subset\Op$. 
Let~$E_0\subseteq\R^n$ be such that
\[\Per(E_0,\Op)<+\infty.\]
Then
\begin{equation}\label{CH:6:inter_per}
\Per(E,\Op)=\Per(E,\Op')+\Per(E_0,\Op\setminus\Op'),\qquad\forall\,E\subseteq\R^n\textrm{ s.t. }E\setminus\Omega=E_0\setminus\Omega.
\end{equation}
In particular, the minimization problem \eqref{CH:6:Min_dir} ``does not depend'' on the choice of $\Op\supset\supset\Omega$, in the sense that if
the exterior data $(u_0,E_0)$ is an admissible pair such that
\[\Nl(u_0,\Omega)<+\infty\quad\textrm{ and }
\quad \Per(E_0,\mathcal O)<+\infty,\]
then a pair $(u,E)$ realizes the infimum
\bgs{
\inf\big\{\Nl(v,\Omega)+\Per(F,\Op)\,|\,(v,F)&\textrm{ admissible s.t. }v=u_0\textrm{ a.e. in }\Co\Omega\\
&\qquad
\textrm{ and }F\setminus\Omega=E_0\setminus\Omega\big\}
}
if and only if it realizes the infimum
\bgs{
\inf\big\{\Nl(v,\Omega)+\Per(F,\Op')\,|\,(v,F)&\textrm{ admissible s.t. }v=u_0\textrm{ a.e. in }\Co\Omega\\
&\qquad
\textrm{ and }F\setminus\Omega=E_0\setminus\Omega\big\},
}
for every $\Omega\Subset\Op'\Subset\Op$.
\end{remark}

Given a fixed bounded open set $\Op\subseteq\R^n$ with Lipschitz boundary such that~$\Omega\Subset\Op$,
we denote
\begin{equation}\label{CH:6:functbar}
\overline\F_\Omega(u,E):=\Nl(u,\Omega)+\Per(E,\Op).\end{equation}
We notice that~$\overline\F_\Omega$ is the functional 
involved in the minimization of the Dirichlet problem~\eqref{CH:6:Min_dir}.

Now we show that Definition~\ref{CH:6:min_pair} is compatible with the
minimization of~$\overline\F_\Omega$, as given by~\eqref{CH:6:Min_dir}.

\begin{lemma}\label{CH:6:lemzero}
A pair $(u,E)$ realizing the infimum in \eqref{CH:6:Min_dir} 
is a minimizing pair in the sense of Definition~\ref{CH:6:min_pair}.
\end{lemma}

\begin{proof}
First of all, notice that
\[\Per(E,\Op)<+\infty.\]
Now let $(v,F)$ be an admissible competitor for $(u,E)$,
according to Definition~\ref{CH:6:min_pair}. Then
\[F\setminus\Omega'=E\setminus\Omega',\]
for some $\Omega'\Subset\Omega$ (with Lipschitz boundary),
thanks to~\eqref{CH:6:minim_def}.
So, by \eqref{CH:6:inter_per}, we have that
$$
\Per(F,\Op)=
\Per(F,\Omega)+\Per(F,\Op\setminus\Omega)\\
= \Per(F,\Omega)+\Per(E,\Op\setminus\Omega).$$
Therefore, recalling~\eqref{CH:6:functional} and~\eqref{CH:6:functbar},
we conclude that
\[\F_\Omega(v,F)-\F_\Omega(u,E)=\overline\F_\Omega(v,F)-
\overline\F_\Omega(u,E)\ge0,\]
which gives the desired result.
\end{proof}

\begin{defn}\label{CH:6:defmin2}
We will say that a pair~$(u,E)$ 
minimizing the Dirichlet problem in~\eqref{CH:6:Min_dir}
is a minimizing pair 
for~$\overline\F_\Omega$ (with respect to the exterior 
data~$(u_0,E_0)$).
\end{defn}

In particular, Lemma~\ref{CH:6:lemzero} says that a minimizing pair
according to Definition~\ref{CH:6:defmin2} is a minimizing
pair according to Definition~\ref{CH:6:min_pair}. 
Now we show that there exists a minimizer
for~$\overline\F_\Omega$, as given by
Definition~\ref{CH:6:defmin2}:

\begin{lemma}\label{CH:6:existence}
Let $\Op\subseteq\R^n$ be a bounded open set with Lipschitz 
boundary such that~$\Omega\Subset\Op$ and let~$(u_0,E_0)$ 
be an admissible pair for~\eqref{CH:6:Min_dir}
such that
\begin{equation}
\label{CH:6:uzero}
\Nl(u_0,\Omega)<+\infty\quad\textrm{ and }\quad 
\Per(E_0,\Op)<+\infty.\end{equation}
Then, there exists a minimizing pair~$(u,E)$ 
for~$\overline\F_\Omega$ with respect to the exterior 
data~$(u_0, E_0)$.
\end{lemma}

\begin{proof}
Since~$(u_0,E_0)$ is an admissible competitor, we have that
\bgs{
\inf\big\{\Nl(v,\Omega)&+\Per(F,\Op)\,|\,(v,F)\textrm{ admissible s.t. }v=u_0\textrm{ a.e. in }\Co\Omega\\
&\qquad\qquad\qquad\quad
\textrm{ and }F\setminus\Omega=E_0\setminus\Omega\big\}\\
&
\le\overline\F_\Omega(u_0,E_0)<+\infty,
}
thanks to~\eqref{CH:6:uzero}.

Now let $(u_k,E_k)$ be a minimizing sequence and notice that
\begin{equation*}
[u_k]_{H^s(\Omega)}^2+\Per(E_k,\Op)\leq\overline\F_\Omega(u_k,E_k)
\leq M\qquad\textrm{for every }k,
\end{equation*}
for some~$M>0$. 
Thus by compactness (see, e.g., \cite[Theorem 7.1]{HitGuide}
and \cite[Theorem 1.19]{Giusti}) we have that
\begin{equation*}\begin{split}
&u_k\to u\quad\textrm{in }L^2(\Omega)\quad\textrm{and a.e. in }\Omega,\\
&
\chi_{E_k}\to\chi_E\quad\textrm{in }L^1(\Op)\quad\textrm{and a.e. in }\Op\quad\textrm{ and }\quad
E_k\setminus\Omega=E_0\setminus\Omega,
\end{split}
\end{equation*}
as~$k\to+\infty$, up to subsequences. Since the functions~$u_k$
are fixed outside~$\Omega$, we actually have that~$u_k\to u$ 
a.e. in~$\R^n$.
Therefore, by Fatou's Lemma, we get
\begin{equation}\label{CH:6:efbhvkwndq}
\Nl(u,\Omega)\leq\liminf_{k\to+\infty}\Nl(u_k,\Omega).
\end{equation}
We remark that the perimeter functional $\Per(\,\cdot\,,\Op)$ is lower semicontinuous with respect to $L^1_{\loc}$ convergence
of sets (see, e.g., \cite[Theorem 1.9]{Giusti}).
This and~\eqref{CH:6:efbhvkwndq} imply
that~$\overline\F_\Omega(u,E)$ attains the desired minimum.

Hence, to complete the proof of Lemma~\ref{CH:6:existence},
we only need to check that
\begin{equation}\label{CH:6:ekfwe4836er}
u\geq0\quad\textrm{a.e. in }E\cap\Omega\quad\textrm{ and }
\quad u\leq0\quad\textrm{a.e. in }\Co E\cap\Omega,\end{equation}
to guarantee that $(u,E)$ is an admissible pair. 

To prove~\eqref{CH:6:ekfwe4836er}, we observe that,
for a.e.~$x\in E\cap\Omega$,
\begin{equation*}
u_k(x)\to u(x)\quad\textrm{and}\quad\chi_{E_k}(x)\to\chi_E(x)=1,
\end{equation*}
and hence $\chi_{E_k}(x)=1$ for every $k$ large enough.
Therefore, for a.e. such $x$, we have that~$u_k(x)\geq0$
for all $k$ large enough, and so also~$u(x)\geq0$, 
which proves the first part of~\eqref{CH:6:ekfwe4836er}. 
A similar argument holds for $\Co E\cap\Omega$, thus completing
the proof of~\eqref{CH:6:ekfwe4836er}.
\end{proof}

Thanks to Lemmata~\ref{CH:6:lemzero} and~\ref{CH:6:existence},
we obtain the existence of a minimizing pair
in the sense of Definition~\ref{CH:6:min_pair}.
In the next result we state the~$s$-harmonicity of the function~$u$
of a minimizing pair $(u,E)$:

\begin{lemma}
Let $(u,E)$ be a minimizing pair in $\Omega$,
according to Definition~\ref{CH:6:min_pair}.
If $\mathcal O\subseteq\Omega$ is an open set such that
\[
\inf_\Op |u|\ge\delta,
\]
for some $\delta>0$, then
\begin{equation*}
(-\Delta)^su(x)=0\qquad
\textrm{for any }x\in\mathcal O.
\end{equation*}
In particular, if $u\in C(\Omega)$, then $(-\Delta)^su=0$ in $\Omega\setminus\{u=0\}$.
\end{lemma}

The proof of the $s$-harmonicity of~$u$
is the same as in~\cite[Lemma~3.2]{DSV}, so we omit the proof here.
Roughly speaking, since the Euler-Lagrange 
functional associated to the functional $\Nl$ in~\eqref{CH:6:def_N} is the
fractional $s$-Laplacian,
the idea consists in considering small perturbating functions $u_\epsilon$ having as positivity set
the positivity set $E$ of $u$, so that when we look at the difference between the energies we get
\begin{equation*}
0\leq\F_\Omega(u_\epsilon,E)-\F_\Omega(u,E)=\Nl(u_\epsilon,\Omega)-\Nl(u,\Omega).
\end{equation*}

\begin{lemma}\label{CH:6:lem:comp}
Let $(u,E)$ be a minimizing pair for $\overline\F_\Omega$,
with respect to the exterior data $(u_0,E_0)$
and let $\alpha\in\R$.
If 
$$ u_0\geq\alpha \quad {\mbox{ a.e. in }}\Co\Omega \quad
{\mbox{(respectively $u_0\leq\alpha$ a.e. in~$\Co\Omega$}}),$$ then
$$u\geq\alpha \quad {\mbox{ a.e. in }}\R^n\quad
{\mbox{(respectively $u\leq\alpha$ a.e. in~$\R^n$).}}$$
\end{lemma}

The proof of the comparison principle in Lemma~\ref{CH:6:lem:comp}
is the same as in \cite[Lemma 3.3]{DSV}, so we omit it.

\subsection{Equivalent characterizations of a minimizing pair}

In this subsection, we give some equivalent definitions
of the notion of minimizing pair. 
 

First of all, notice that, if~$\Omega'\subseteq\Omega$, then
the functional in~\eqref{CH:6:def_N} can be written as
\begin{equation}\label{CH:6:different_sets}
\Nl(v,\Omega)=\Nl(v,\Omega')+[v]^2_{H^s(\Omega\setminus\Omega')}
+2\iint_{(\Omega\setminus\Omega')\times\Co\Omega}
\frac{|v(x)-v(y)|^2}{|x-y|^{n+2s}}\,dx\,dy.
\end{equation}
In particular, if $v=u$ a.e. in $\Co\Omega'$ and $\Nl(u,\Omega)<+\infty$, then
from \eqref{CH:6:different_sets} we see that
\begin{equation}\label{CH:6:useful_stupid}
\Nl(v,\Omega)<+\infty\quad\Longleftrightarrow\quad
\Nl(v,\Omega')<+\infty,
\end{equation}
and
\begin{equation}\label{CH:6:useful_koala}
\Nl(v,\Omega')-\Nl(u,\Omega')=\Nl(v,\Omega)-\Nl(u,\Omega).
\end{equation}

We also point out the following trivial but useful remark, which explains why in the definition of an admissible competitor
we ask $u-v\in H^s(\R^n)$. This is indeed equivalent to asking $\Nl(v,\Omega)<+\infty.$
\begin{remark}
Let $\Omega\subseteq\R^n$ be a bounded open set and let $u:\R^n\to\R$ be such that $\Nl(u,\Omega)<+\infty$.
Let $v:\R^n\to\R$ be such that $v=u$ a.e. in $\Co\Omega$. Then
\begin{equation}\label{CH:6:equiv_comp_assump}
\Nl(v,\Omega)<\infty\qquad\Longleftrightarrow\qquad u-v\in H^s(\R^n).
\end{equation}
First of all, we remark that
\[
[u-v]_{H^s(\R^n)}<+\infty\qquad\Longrightarrow\qquad
\|u-v\|_{L^2(\R^n)}=\|u-v\|_{L^2(\Omega)}<+\infty.
\]

This is a consequence of a fractional Poincar\'e type inequality---see, e.g., Proposition~\ref{CH:4:FPI}---which we can apply to the function
$w:=v-u$ thanks to the assumption
\[
v=u\quad\mbox{a.e. in }\Co\Omega,
\]
so it is enough to show that
\[
\Nl(v,\Omega)<\infty\qquad\Longleftrightarrow\qquad [u-v]_{H^s(\R^n)}<+\infty.
\]
This equivalence follows from the equality
\begin{equation*}
[u-v]_{H^s(\R^n)}^2=\iint_{\R^{2n}\setminus(\Co\Omega)^2}\frac{|u(x)-v(x)-u(y)+v(y)|^2}{|x-y|^{n+2s}}\,dx\,dy
=\Nl(u-v,\Omega)
\end{equation*}
and the ``triangle inequality''
\[
\Nl(u_1+u_2,\Omega)\le2\big(\Nl(u_1,\Omega)+\Nl(u_2,\Omega)\big).
\]
\end{remark}

As a consequence of formulas~\eqref{CH:6:inter_per}
and~\eqref{CH:6:different_sets}, we obtain
the following equivalent characterizations of minimizing pairs:

\begin{prop}\label{CH:6:equiv_char}
Let $(u,E)$ be an admissible pair
according to Definition~\ref{CH:6:min_pair}
such that $\F_\Omega(u,E)<+\infty$. 
Then, the following statements are equivalent:
\begin{itemize}
\item[(i)] the pair $(u,E)$ is minimizing in $\Omega$,
according to Definition~\ref{CH:6:min_pair},
\item[(ii)] for every open subset $\Omega'\Subset\Omega$ 
we have
\bgs{
\Nl(u,\Omega')+\Per(E,\Omega)& =\inf\big\{\Nl(v,\Omega')+
\Per(F,\Omega)\,|\,(v,F)\textrm{ admissible}\\
&\qquad\qquad\textrm{s.t. }
v=u\textrm{ a.e. in }\Co\Omega'\textrm{ and }
F\setminus\Omega'=E\setminus\Omega'\big\}.
}
\item[(iii)] the pair $(u,E)$ is minimizing in every open set $\Omega'\Subset\Omega$,
\item[(iv)] the pair $(u,E)$ is minimizing in every open set $\Omega'\subseteq\Omega$.
\end{itemize}

\end{prop}

\begin{proof}
We begin with the implication $(i)\,\Longrightarrow\,(ii)$.

Let $(v,F)$ be an admissible pair such that
\[
v=u\textrm{ a.e. in }\Co\Omega'\textrm{ and }
F\setminus\Omega'=E\setminus\Omega'.
\]
We can suppose that
\[
\Nl(v,\Omega')<+\infty\qquad\mbox{and}\qquad \Per(F,\Omega)<+\infty,
\]
otherwise there is nothing to prove. In particular, thanks to \eqref{CH:6:useful_stupid} we have
\[
\Nl(v,\Omega)<+\infty.
\]

Thus $(v,F)$ is an admissible competitor for $(u,E)$ in $\Omega$, according to Definition \ref{CH:6:min_pair}.\\
By minimality
of $(u,E)$ and equality \eqref{CH:6:useful_koala} we obtain
\[
\Nl(u,\Omega')+\Per(E,\Omega)-\Nl(v,\Omega')-\Per(F,\Omega)
=\F_\Omega(u,E)-\F_\Omega(v,F)\le0,
\]
as wanted.

As for the implication $(ii)\,\Longrightarrow\,(iii)$,
let $(v,F)$ be an admissible competitor for $(u,E)$ in $\Omega'$.\\
Then we can find an open set
$\Op\Subset\Omega'$ such that $v=u$ a.e. in $\Co\Op$ and $F\Delta E\subseteq\Op$.
Exploiting both \eqref{CH:6:inter_per}
and~\eqref{CH:6:different_sets}, we find
\[
\F_{\Omega'}(v,E)-\F_{\Omega'}(u,E)=\Nl(v,\Op)+\Per(F,\Omega)
-\Nl(u,\Op)-\Per(E,\Omega),
\]
which is nonnegative by $(ii)$.

The implication $(iii)\,\Longrightarrow\,(iv)$ is proved in the same way.
If $(v,F)$ is an admissible competitor for $(u,E)$ in $\Omega'$, then we can find
$\Op\Subset\Omega'$ such that supp$(v-u)\Subset\Op$ and $F\Delta E\Subset\Op$.\\
Then $(v,F)$ is an admissible competitor for $(u,E)$ in $\Op$. Exploiting the minimality assumed in $(iii)$
and using again both \eqref{CH:6:inter_per}
and~\eqref{CH:6:different_sets}, we thus obtain
\[
\F_{\Omega'}(v,E)-\F_{\Omega'}(u,E)=\F_\Op(v,F)-\F_\Op(u,E)\ge0.
\]

The last implication $(iv)\,\Longrightarrow\,(i)$ follows trivially by taking $\Omega'=\Omega$.
\end{proof}

\begin{remark}
Notice that point $(ii)$ of Proposition~\ref{CH:6:equiv_char} 
says that $(u,E)$ is a minimizing pair for $\overline\F_{\Omega'}$
for every open subset $\Omega'\Subset\Omega$ 
(with respect to the exterior data $(u,E)$).
\end{remark}

\section{The extended functional}\label{CH:6:sec:ext}

In this section we deal with the extended functional
defined in the Introduction. 
For an introduction to the extension operator, we refer the interested reader to \cite{Extension}.

We recall that, given a function $u:\R^n\to\R$,
we denote by $\Ue:\R^{n+1}_+\to\R$ the extended function
defined in~\eqref{CH:6:extended_func_def}, that is
\[
\Ue(x,z):=c_{n,s}z^{2s}\int_{\R^n}\frac{u(\xi)}{(|x-\xi|^2+z^2)^\frac{n+2s}{2}}d\xi,
\quad\mbox{for every }(x,z)\in\R^{n+1}_+.
\]
We observe that for the extended function $\Ue$ to be well defined, it is enough that $u\in L_s(\R^n)$.
Hence, in light of Remark \ref{CH:6:tail_energies_rmk}, if $u:\R^n\to\R$ is a measurable function such that
$\Nl(u,\Omega)<+\infty$, then the extended function $\Ue$ is well defined.


We start with some preliminary
observations:

\begin{remark}\label{CH:6:useful_rmk_energies}
If $\Nl(u,B_R)<+\infty$, then 
\begin{equation*} 
\int_{\Omega_+}|\nabla\Ue|^2z^{1-2s}\,dX<+\infty,
\end{equation*}
for every bounded open set $\Omega\subseteq\R^{n+1}$ with Lipschitz boundary 
and such that~$\Omega_0\Subset B_R$
(see~\cite[Proposition~7.1]{CRS10}).

In particular, if $(u,E)$ is an admissible pair for $\F$ s.t. $\F_{B_R}(u,E)<+\infty$, then $(\Ue,E)$
is an admissible pair for the extended functional $\lf$,
and $\lf_\Omega(\Ue,E)<+\infty$
for every bounded open set $\Omega\subseteq\R^{n+1}$ 
with Lipschitz boundary and such that~$\Omega_0\Subset B_R$.
\end{remark}

\begin{remark}\label{CH:6:useful_rmk_energies1}
Let $(u,E)$ be an admissible pair for $\F$ s.t. $\F_{B_R}(u,E)<+\infty$.
Let $\Omega\subseteq\R^{n+1}$ 
be a bounded open set with Lipschitz boundary ans such 
that~$\Omega_0\Subset B_R$, and let~$(\Vf,F)$
be an admissible competitor for $\lf_\Omega$, 
with respect to $(\Ue,E)$, according to Definition~\ref{CH:6:adm-ext}.
Define~$v:=\Vf\big|_{\{z=0\}}$.
Then $(v,F)$ is an admissible competitor for $\F_{B_R}$, 
with respect to $(u,E)$, according to Definition~\ref{CH:6:min_pair}.

Indeed, let
$\Omega'\subseteq\R^{n+1}$ be a bounded open set with 
Lipschitz boundary, such that~$\Omega\Subset\Omega'$ and
$\Omega'_0\Subset B_R$. From Remark~\ref{CH:6:useful_rmk_energies},
we know that
\begin{equation*}
\int_{\Omega'_+\setminus\Omega_+}|\nabla\Ue|^2z^{1-2s}\,dX<+\infty,
\end{equation*}
and hence, since $\lf(\Vf,\Omega)<+\infty$ and supp $(\Vf-\Ue)\Subset\Omega$, we get
\begin{equation*}
\int_{\Omega'_+}|\nabla\Vf|^2z^{1-2s}\,dX<+\infty.
\end{equation*}
It can be shown that this implies that $\Nl(v,\Omega)<+\infty$ (see e.g. the proof of \cite[Proposition 4.1]{DSV}).
Now, since $v=u$
in $\Co\Omega$ and $\Nl(u,B_R)<+\infty$, using \eqref{CH:6:useful_stupid} we get $\Nl(v,B_R)<+\infty$
and $u-v\in H^s(\R^n)$ as claimed.
\end{remark}

\subsection{An equivalent problem}

Now we show that we can use the extended functional $\lf$, 
defined in~\eqref{CH:6:extended_functional}, to obtain
an equivalent formulation of the minimization problem for $\F$.

We remark that, differently from the proof of \cite[Proposition~4.1]{DSV}, in our framework
we only ``localize'' the energy $\Nl$.

\begin{proof}[Proof of Proposition \ref{CH:6:Local_energy_prop}]
Let $r\in(0,R)$. From
\cite[Lemma 7.2]{CRS10} we know that
if $v:\R^n\to\R$ is such that
\begin{equation}\label{CH:6:hyp_Lemma_CRS10}
\Nl(v,B_r)<+\infty\qquad\textrm{and}\qquad\textrm{supp}(v-u)\Subset B_r,
\end{equation}
 then
\begin{equation}\label{CH:6:local_version}
\Nl(v,B_r)-\Nl(u,B_r)=c_{n,s}'\inf_{(\Omega,\Vf)\in\mathfrak J_v}\int_{\Omega_+}\big(|\nabla\Vf|^2-|\nabla\Ue|^2\big)z^{1-2s}\,dX,
\end{equation}
where the set $\mathfrak J_v$ consists of all the
couples~$(\Omega,\Vf)$, with $\Omega\subseteq\R^{n+1}$ a bounded open set with Lipschitz boundary such that $\Omega_0\subseteq B_r$ and $\Vf:\R^{n+1}\to\R$ such that
$\Vf-\Ue$ is compactly supported inside $\Omega$ and $\Vf\big|_{\{z=0\}}=v$.

Notice that for every such couple $(\Omega,\Vf)\in\mathfrak J_v$ we can prescribe without loss of generality
that $\Vf=\Ue$ outside $\Omega$.\\

$\Longrightarrow)\quad$ Let $(u,E)$ be a minimizing pair
for $\F$ in $B_r$, with $r\in(0,R)$.
We show that $(\Ue,E)$ is minimizing for $\lf_\Omega$ for every bounded open set $\Omega\subseteq\R^{n+1}$ with Lipschitz boundary and $\Omega_0\subseteq B_r$.\\
From Remark~\ref{CH:6:useful_rmk_energies}, we know that $(\Ue,E)$ is admissible for the
extended functional and $\lf_\Omega(\Ue,E)<+\infty$.

Now let
$(\Vf,F)$ be an admissible competitor
and define $v:=\Vf|_{\{z=0\}}$, so that $(\Omega,\Vf)\in\mathfrak J_v$.
Since $v-u=\Vf|_{\{z=0\}}-\Ue|_{\{z=0\}}$ is compactly supported in $\Omega_0\subseteq B_r$,
from Remark \ref{CH:6:useful_rmk_energies1}
 we see that $(v,F)$ is an admissible competitor for $\F$ in $B_r$.
Thus, using the minimality of $(u,E)$ and~\eqref{CH:6:local_version}, we obtain
\begin{eqnarray*}
0&\leq&\F_{B_r}(v,F)-\F_{B_r}(u,E)\\&=&\Nl(v,B_r)-\Nl(u,B_r)+\Per(F,B_r)-\Per(E,B_r)\\
&
=& c_{n,s}'\inf_{(\Omega,\Vf)\in\mathfrak J_v}\int_{\Omega_+}\big(|\nabla\Vf|^2-|\nabla\Ue|^2\big)z^{1-2s}\,dX
+\Per(F,B_r)-\Per(E,B_r)\\
&
\leq &\lf_\Omega(\Vf,F)-\lf_\Omega(\Ue,E).
\end{eqnarray*}
Since this holds for every admissible competitor, this shows that $(\Ue,E)$ is minimizing for $\lf_\Omega$.\\

$\Longleftarrow)\quad$ Suppose that $(\Ue,E)$ is minimizing for $\lf_\Omega$, for every $\Omega\subseteq\R^{n+1}$ as in
the statement of Proposition~\ref{CH:6:Local_energy_prop}.

Let $(v,F)$ be an admissible competitor for $\F$ in $B_R$. 
In particular, we have that
\begin{equation*}
\textrm{supp}\,(v-u)\Subset B_R\qquad\textrm{and}\qquad E\Delta F\Subset B_R,
\end{equation*}
hence we can suppose that
\begin{equation}\label{CH:6:stupid_incl}
\textrm{supp}\,(v-u)\Subset B_r\qquad\textrm{and}\qquad E\Delta F\Subset B_r,
\end{equation}
for some $r\in(0,R)$.

Notice that $v$ satisfies~\eqref{CH:6:hyp_Lemma_CRS10}
and that if $(\Omega,\Vf)\in\mathfrak J_v$, then $(\Vf,F)$ is an admissible
competitor for $\lf_\Omega$ with respect to $(\Ue,E)$ and $\Omega_0\subseteq B_r\Subset B_R$.

Thus, if $(\Omega,\Vf)\in\mathfrak J_v$, since $(\Ue,E)$ is minimizing for $\lf_\Omega$, we get
\begin{equation*}
\int_{\Omega_+}\big(|\nabla\Vf|^2-|\nabla\Ue|^2
\big)z^{1-2s}\,dX
+\Per(F,B_r)-\Per(E,B_r)
=\lf_\Omega(\Vf,F)-\lf_\Omega(\Ue,E)\geq0.
\end{equation*}
Since this holds true for every $(\Omega,\Vf)\in\mathfrak J_v$ and $(v,F)$
satisfies~\eqref{CH:6:stupid_incl}, we get from~\eqref{CH:6:local_version} that
\begin{equation*}\begin{split}
&\F_{B_R}(v,F)-\F_{B_R}(u,E)=\F_{B_r}(v,F)-\F_{B_r}(u,E)\\
&=c_{n,s}'\inf_{(\Omega,\Vf)\in\mathfrak J_v}\int_{\Omega_+}\big(|\nabla\Vf|^2-|\nabla\Ue|^2\big)z^{1-2s}\,dX
+\Per(F,B_r)-\Per(E,B_r)
\geq0.
\end{split}
\end{equation*}
This shows that $(u,E)$ is minimizing in $B_R$.
\end{proof}

%
%
%
%
%
%
%
%
%
%
%
%
%

\section{Monotonicity formula}\label{CH:6:MoNoToFoRSeC}

In this subsection, we obtain a monotonicity formula in the spirit
of~\cite{weiss}.
The main feature here is that we need to consider the associated
extension
problem to prove that some energy 
is monotone. As usual in this type of problems,
this will imply a homogeneity of the functions involved.
Other papers in which this approach has been 
exploited are~\cite{CSV, DSV}.
\smallskip

We introduce now some notation.
We say that a set $A\subseteq\R^n$ is a cone if $\lambda A=A$ 
for any~$\lambda>0$.
Notice that this is the same as asking $\chi_A$ to be 
homogeneous of degree~$0$, that is~$\chi_A(\lambda x)=
\chi_A(x)$ for any~$x\in\R^n$ and any~$\lambda>0$.

\smallskip

First of all we show that the functional $\lf$ 
possesses a natural scaling.
For this, recall the definition of the rescaled pairs $(u_\lambda,E_\lambda)$ given in \eqref{CH:6:rescaled_stuff}.


We recall also the notation
\begin{equation*}
\BaLL_r:=\{(x,z)\in\R^{n+1}\,|\,|x|^2+z^2<r^2\}\qquad\textrm{and}\qquad
\BaLL_r^+:=\BaLL_r\cap\{z>0\}.
\end{equation*}
We can now prove the following scaling result:

\begin{lemma}
Let $(u,E)$ be a minimizing pair for $\F$ in $B_R$. Define
\begin{equation}
\Gs_u(r):=r^{1-n}\lf_{\BaLL_r}(\Ue,E)=r^{1-n}
\Big(c'_{n,s}\int_{\BaLL_r^+}|\nabla\Ue|^2z^{1-2s}\,dX+\Per(E,B_r)\Big)
\end{equation}\label{CH:6:scal}
for any $r\in(0,R)$, where~$\lf$ has been introduced 
in~\eqref{CH:6:extended_functional}. Then, for any $\lambda>0$,
\begin{equation}\label{CH:6:scal2}
\Gs_u(\lambda r)=\Gs_{u_\lambda}(r).
\end{equation}
\end{lemma}

\begin{proof}
We know that the perimeter scales as
\begin{equation}\label{CH:6:uno}
\Per(E_\lambda,\Omega_\lambda)=\lambda^{1-n}\Per(E,\Omega).
\end{equation}
As for the energy of the extended functions,
it is enough to notice that if $\Ue_\lambda$
denotes the extension of $u_\lambda$ (as given by~\eqref{CH:6:extended_func_def}),
then 
\begin{equation}\label{CH:6:due}
\Ue_\lambda(X)=\lambda^{\frac{1}{2}-s}\Ue(\lambda X).
\end{equation}
Plugging~\eqref{CH:6:uno} and~\eqref{CH:6:due} into~\eqref{CH:6:scal},
we obtain the desired formula in~\eqref{CH:6:scal2}.
\end{proof}

Now we ``correct'' $\Gs_u$ by adding an appropriate term,
\begin{equation*}
\Phi_u(r):=\Gs_u(r)-\Cf_u(r),
\end{equation*}
where
\begin{equation*}
\mathfrak C_u:=c'_{n,s}\Big(s-\frac{1}{2}\Big)r^{-n}\int_{(\partial\BaLL_r)^+}\Ue^2z^{1-2s}\,d\Ha^n,
\end{equation*}
and we prove a monotonicity formula for $\Phi_u$. 
Here above, we used the notation
\begin{equation*}
(\partial\BaLL_r)^+:=\partial\BaLL_r\cap\{z>0\}=
\{(x,z)\in\R^{n+1}_+\,|\,|x|^2+z^2=r^2\}.
\end{equation*}

\begin{remark}
It is not difficult to see that~$\Phi_u$
has the same scale invariance property of~$\Gs_u$, i.e.
\begin{equation}\label{CH:6:scaling_wholephi}
\Phi_u(\lambda r)=\Phi_{u_\lambda}(r).
\end{equation}
\end{remark}

\begin{remark}\label{CH:6:scaling_minimality}
Before proving the Monotonicity Formula, we point out that $(u,E)$ is minimal in $\Omega$ if and only if $(u_\lambda,E_\lambda)$
is minimal in $\Omega_\lambda$, for every $\lambda>0$. This is a consequence of the homogeneous scaling
\[
\F_{\Omega_\lambda}(v_\lambda,F_\lambda)=\lambda^{1-n}\F_\Omega(v,F).
\]
Indeed, it is enough to notice that $(v,F)$ is an admissible competitor for $(u_\lambda,E_\lambda)$ in $\Omega_\lambda$
if and only if $(v_{1/\lambda},F_{1/\lambda})$ is an admissible competitor for $(u,E)$ in $\Omega$. Then, if $(u,E)$ is minimal
in $\Omega$, we find
\[
\F_{\Omega_\lambda}(v,F)=\lambda^{1-n}\F_\Omega(v_{1/\lambda},F_{1/\lambda})
\ge\lambda^{1-n}\F_\Omega(u,E)=\F_{\Omega_\lambda}(u_\lambda,E_\lambda).
\]
\end{remark}

%

\begin{proof}[Proof of Theorem \ref{CH:6:Monotonicity_teo}]
First of all, notice that $\Phi_u$ is differentiable a.e. in $(0,R)$.
We want to prove that there exists a subset $\mathcal G\subseteq(0,R)$ with $\mathcal L^1\big((0,R)\setminus\mathcal G\big)=0$
and such that
\begin{equation}\label{CH:6:claim_monoto}
\exists\,\frac{d}{dr}\Phi_u(r)\geq0\qquad\textrm{for every }r\in\mathcal G.
\end{equation}
We remark that, even if the function~$\Phi_u$
in general is not continuous, \eqref{CH:6:claim_monoto} is enough to
prove that~$\Phi_u$ is increasing in~$(0,R)$,
thanks to Lemma~\ref{CH:6:guy_lemma} and Corollary~\ref{CH:6:guy_coroll}.

Indeed, let
\[\theta_E(r):=\frac{\Per(E,B_r)}{r^{n-1}}\qquad\textrm{and}\qquad f(r):=\Phi_u(r)-\theta_E(r),\]
and notice that, since $f$ is continuous and differentiable a.e. in $(0,R)$, we can write
\begin{equation}\label{CH:6:crocco}
f(r_2)-f(r_1)=\int_{r_1}^{r_2}f'({\varrho})\,d{\varrho},\qquad\textrm{for every }0<r_1<r_2<R.
\end{equation}
Now suppose that \eqref{CH:6:claim_monoto} holds true
and notice that
\[\Phi_u'(r)=f'(r)+\theta_E'(r)\qquad\textrm{for a.e. }r\in(0,R).\]
Then, exploiting \eqref{CH:6:crocco} and formula \eqref{CH:6:formula_surf_dens_guy} we obtain
\begin{equation*}\begin{split}
\Phi_u(r_2)-\Phi_u(r_1)&=f(r_2)-f(r_1)+\theta_E(r_2)-\theta_E(r_1)\ge\int_{r_1}^{r_2}f'({\varrho})\,d{\varrho}
+\int_{r_1}^{r_2}\theta_E'({\varrho})\,d{\varrho}\\
&=\int_{r_1}^{r_2}\Phi'_u({\varrho})\,d{\varrho}
=\int_{(r_1,r_2)\cap\mathcal G}\Phi'_u({\varrho})\,d{\varrho}\ge0,
\end{split}\end{equation*}
for every $0<r_1<r_2<R$, thus proving the monotonicity of $\Phi_u$.

We also remark that, if we denote
\[\wp(r):=\Per(E,B_r),\]
then $\theta_E$ is differentiable at $r\in(0,R)$ if and only if $\wp$ is differentiable at $r$, and in this case
\[\theta_E'(r)=r^{1-n}\wp'(r)-(n-1)r^{-n}\wp(r).\]
\smallskip

Now we define the subset $\mathcal G\subseteq(0,R)$.\\
Notice that since $(u,E)$ is
an admissible pair, we have that~$|\{u<0\}\cap E|=0$.
Exploiting spherical
coordinates, we see that for a.e. $r>0$
\begin{equation*}
u(x)\geq0\qquad\textrm{for }\Ha^{n-1}\textrm{-a.e. }x\in E\cap\partial B_r.
\end{equation*}
In the same way, for a.e. $r>0$
\begin{equation*}
u(x)\leq0\qquad\textrm{for }\Ha^{n-1}\textrm{-a.e. }x
\in \Co E\cap\partial B_r.
\end{equation*}
All in all, we see that, for a.e. $r>0$,
\begin{equation}\label{CH:6:admissibility_on_slices}\begin{split}
& u(x)\geq0\qquad\textrm{for }
\Ha^{n-1}\textrm{-a.e. }x\in E\cap\partial B_r\\
{\mbox{and}}\quad &u(x)\leq0\qquad\textrm{for }
\Ha^{n-1}\textrm{-a.e. }x\in \Co E\cap\partial B_r.
\end{split}
\end{equation}
The set $\mathcal G$ is defined as the set of all those $r\in(0,R)$ which satisfy all the following properties:
\begin{itemize}
\item[(i)] \eqref{CH:6:admissibility_on_slices} holds true,

\item[(ii)] the functions $f$ and $\wp$ are differentiable at $r$

\item[(iii)] it holds
\[\Ha^{n-2}(\partial^*E\cap\partial B_r)<+\infty,\]
and $r$ is a Lebesgue point for the function
\[(0,R)\ni\varrho\longmapsto\Ha^{n-2}(\partial^*E\cap\partial B_\varrho),\]

\item[(iv)] the cone~$E(r)$ with vertex in~$0$ spanned by the spherical
slice~$E\cap\partial B_r$ (as defined in~\eqref{CH:6:3.4bis}) is a Caccioppoli set.
\end{itemize}

We remark that by Remark \ref{CH:6:layer_cake_surf}, Proposition \ref{CH:6:cones_from_slices} and Remark~\ref{CH:6:cone_competitor},
points (iii) and (iv) hold true for a.e. $r\in(0,R)$. Hence $\mathcal L^1\big((0,R)\setminus\mathcal G\big)=0$.
\bigskip

Now we prove claim \eqref{CH:6:claim_monoto}.\\
First of all, notice that thanks to the scaling property \eqref{CH:6:scaling_wholephi} we can assume without loss of generality
that $r=1$.
We have
\eqlab{\label{CH:6:monotone_piece1}
\Gs_u'(1)=c'_{n,s}\Big((1-n)&\int_{\BaLL_1^+}|\nabla\Ue|^2z^{1-2s}\,dX
+\int_{(\partial\BaLL_1)^+}|\nabla\Ue|^2z^{1-2s}\,d\Ha^n\Big)\\
&\qquad
+\frac{d}{dr}\frac{\Per(E,B_r)}{r^{n-1}}\Big|_{r=1}
}
and
\begin{equation}\label{CH:6:monotone_piece2}
\Cf_u'(1)=c'_{n,s}\Big(s-\frac{1}{2}\Big)\int_{(\partial\BaLL_1)^+}\big(2\Ue\,\Ue_\nu+(1-2s)\Ue^2\big)z^{1-2s}\,d\Ha^n,
\end{equation}
where $\Ue_\nu$ denotes the normal derivative of $\Ue$, so that the normal gradient is $\Ue_\nu(X)X$.\\
To prove~\eqref{CH:6:monotone_piece2} notice that changing variables $X=rY$, with $z=rw$ yields
\begin{equation*}
\Cf_u(r)=c'_{n,s}\Big(s-\frac{1}{2}\Big)r^{1-2s}\int_{(\partial\BaLL_1)^+}\Ue^2(rY)w^{1-2s}\,d\Ha^n(Y).
\end{equation*}
Then take the derivative in $r$ and set $r=1$.

To show that $\Phi'_u(1)\geq0$ we construct appropriate competitors for $(\Ue,E)$ and
compare the energies.\\
Given a small $\eps>0$, we consider the admissible competitor $(\Uc^\eps,E^\eps)$ for $(\Ue,E)$
defined as
\begin{equation*}
\Uc^\eps(X):=\left\{\begin{array}{cc}
(1-\eps)^{s-\frac{1}{2}}\Ue\big(\frac{1}{1-\eps}X\big)&\textrm{if }X\in\BaLL_{1-\eps}^+,\\
|X|^{s-\frac{1}{2}}\Ue\big(\frac{X}{|X|}\big)&\textrm{if }X\in\BaLL_1^+\setminus\BaLL_{1-\eps}^+,\\
\Ue(X)&\textrm{if }X\in\R^{n+1}_+\setminus\BaLL_1^+,
\end{array}\right.
\end{equation*}
and
\begin{equation*}
\chi_{E^\eps}(x):=\left\{\begin{array}{cc}
\chi_E\big(\frac{1}{1-\eps}x\big)&\textrm{if }x\in B_{1-\eps},\\
\chi_E\big(\frac{x}{|x|}\big)&\textrm{if }x\in B_1\setminus B_{1-\eps},\\
\chi_E(x)&\textrm{if }x\in\R^n\setminus B_1,
\end{array}\right.
\end{equation*}
that is
\begin{equation*}
E^\eps:=\big((1-\eps)E\cap B_{1-\eps}\big)
\cup \big(E(1)\cap(B_1\setminus B_{1-\eps})\big)
\cup(E\setminus B_1).
\end{equation*}
Let $u^\eps:=\Uc^\eps\big|_{\{z=0\}}$ be the trace of $\Uc^\eps$. It is clear that
\begin{equation*}
u^\eps\geq0\quad\textrm{a.e. in }
E^\eps\setminus(B_1\setminus B_{1-\eps})\quad\textrm{and}\quad
u^\eps\leq0\quad\textrm{a.e. in }\Co E^\eps\setminus(B_1\setminus B_{1-\eps}).
\end{equation*}
Moreover condition~\eqref{CH:6:admissibility_on_slices} (with $r=1$) guarantees that the same
holds also a.e. in $B_1\setminus B_{1-\eps}$,
so that $(\Uc^\eps,E^\eps)$ is an admissible pair for $\lf$.

By construction $(\Uc^\eps,E^\eps)$ is an admissible competitor for $(\Ue,E)$ in every bounded open set
$\Omega\subseteq\R^{n+1}$ with Lipschitz boundary and such that $\BaLL_1^+\Subset\Omega$ and $\Omega_0\Subset B_R$,
since
\begin{equation}\label{CH:6:over_cond_monotonicity}
\Uc^\eps=\Ue\quad\textrm{in}\quad\R^{n+1}_+\setminus\BaLL_1^+\qquad\textrm{and}\qquad E^\eps=E\quad\textrm{in} \quad\R^n\setminus B_1.
\end{equation}
In particular we can take $\Omega=\BaLL_{\varrho}$ for some ${\varrho}\in(1,R)$ (recall that we are assuming $1=r<R$).

Since
\[\Ha^{n-2}(\partial^*E\cap\partial B_1)<+\infty\quad\Longrightarrow\quad\Ha^{n-1}(\partial^*E\cap\partial B_1)=0,\]
by the definitions of $E(1)$ and $E^\eps$ and formulas \eqref{CH:6:slices_of_cone},
thanks to \cite[Theorem 16.16]{Maggi}
we have
\begin{equation}\label{CH:6:perimetercompetitormono}\begin{split}
\Per(E^\eps,B_{\varrho})&=\Per(E^\eps,B_1)+\Per(E,B_{\varrho}\setminus\overline{B_1})\\
&
=\Per(E_{1/(1-\eps)},B_{1-\eps})+\Per(E(1),B_1\setminus\overline{B_{1-\eps}})
+\Per(E,B_{\varrho}\setminus\overline{B_1}).
\end{split}\end{equation}
Thus, using also~\eqref{CH:6:over_cond_monotonicity}, we get from Proposition $\ref{CH:6:Local_energy_prop}$
\begin{equation}\label{CH:6:equation_mono1}
\lf_{\BaLL_1}(\Uc^\eps,E^\eps)-\lf_{\BaLL_1}(\Ue,E)
=\lf_{\BaLL_{\varrho}}(\Uc^\eps,E^\eps)-\lf_{\BaLL_{\varrho}}(\Ue,E)\geq0.
\end{equation}

We compute $\lf_{\BaLL_1}(\Uc^\eps,E^\eps)$ by splitting it in $\BaLL_{1-\eps}^+$ and $\BaLL_1^+\setminus\BaLL_{1-\eps}^+$.\\
Notice that in $\BaLL_{1-\eps}^+$ the pair $(\Uc^\eps,E^\eps)$ is just the rescaled pair $(\Ue_{1/(1-\eps)},E_{1/(1-\eps)})$.
Then
\begin{equation*}\begin{split}
\lf_{\BaLL_1}(\Uc^\eps,E^\eps)&=\lf_{\BaLL_{1-\eps}}(\Ue_{1/(1-\eps)},E_{1/(1-\eps)})
+\lf_{\BaLL_1\setminus\BaLL_{1-\eps}}(\Uc^\eps,E^\eps)\\
&
=(1-\eps)^{n-1}\Gs_{u_{1/(1-\eps)}}(1-\eps)+\lf_{\BaLL_1\setminus\BaLL_{1-\eps}}(\Uc^\eps,E^\eps)\\
&
=(1-\eps)^{n-1}\Gs_u(1)+\lf_{\BaLL_1\setminus\BaLL_{1-\eps}}(\Uc^\eps,E^\eps)
\end{split}\end{equation*}
Now we compute $\lf_{\BaLL_1\setminus\BaLL_{1-\eps}}(\Uc^\eps,E^\eps)$.\\
As for the perimeter, (recalling \eqref{CH:6:perimetercompetitormono}) by formula \eqref{CH:6:perimeter_cone} we have
\begin{equation*}\begin{split}
\Per(E^\eps,B_1\setminus \overline{B_{1-\eps}})&
=\Per(E(1),B_1)-\Per(E(1),B_{1-\eps})\\
&
=\frac{\Ha^{n-2}(\partial^*E\cap\partial B_1)}{n-1}
\big(1-(1-\eps)^{n-1}\big)\\
&
=\eps\,\Ha^{n-2}(\partial^*E\cap\partial B_1)+o(\eps).
\end{split}
\end{equation*}
Notice that in $\BaLL_1^+\setminus\BaLL_{1-\eps}^+$ we have
\bgs{
\nabla\Uc^\eps(X)&=\Big(s-\frac{1}{2}\Big)|X|^{s-\frac{3}{2}}\Ue\left(\frac{X}{|X|}\right)\frac{X}{|X|}
+|X|^{s-\frac{1}{2}}\frac{1}{|X|}\left[\nabla\Ue\left(\frac{X}{|X|}\right)\right.\\
&
\qquad\qquad\qquad
-\left.\left(\nabla\Ue\left(\frac{X}{|X|}\right)\cdot\frac{X}{|X|}\right)\frac{X}{|X|}\right]\\
&
=\Big(s-\frac{1}{2}\Big)|X|^{s-\frac{3}{2}}\Ue\left(\frac{X}{|X|}\right)\frac{X}{|X|}
+|X|^{s-\frac{1}{2}}\frac{1}{|X|}\Ue_\tau\left(\frac{X}{|X|}\right),
}
where $\Ue_\tau$ denotes the tangential gradient of $\Ue$ on $(\partial\BaLL_1)^+$.\\
Since $\Ue_\tau\cdot\frac{X}{|X|}=0$, this gives
\begin{equation}
|\nabla\Uc^\eps(X)|^2=
|X|^{2s-3}\left\{\Big(s-\frac{1}{2}\Big)^2\Ue^2\left(\frac{X}{|X|}\right)+\left|\Ue_\tau\left(\frac{X}{|X|}\right)\right|^2\right\}.
\end{equation}
Therefore
\begin{equation*}\begin{split}
\int_{\BaLL_1^+\setminus\BaLL_{1-\eps}^+}|\nabla\Uc^\eps|^2&z^{1-2s}\,dX
=\int_{1-\eps}^1dt\int_{(\partial\BaLL_t)^+}
|\nabla\Uc^\eps|^2z^{1-2s}\,d\Ha^n\\
&
=\int_{1-\eps}^1t^{2s-3}\,dt\int_{(\partial\BaLL_t)^+}
\left\{\Big(s-\frac{1}{2}\Big)^2\Ue^2\left(\frac{X}{|X|}\right)+\left|\Ue_\tau\left(\frac{X}{|X|}\right)\right|^2\right\}z^{1-2s}
\,d\Ha^n\\
&
=
\int_{1-\eps}^1t^{-2}\,dt\int_{(\partial\BaLL_1)^+}
\left\{\Big(s-\frac{1}{2}\Big)^2\Ue^2+|\Ue_\tau|^2\right\}z^{1-2s}\,d\Ha^n\\
&
=\Big(\frac{1}{1-\eps}-1\Big)
\int_{(\partial\BaLL_1)^+}
\left\{\Big(s-\frac{1}{2}\Big)^2\Ue^2+|\Ue_\tau|^2\right\}z^{1-2s}\,d\Ha^n\\
&
=\eps\int_{(\partial\BaLL_1)^+}
\left\{\Big(s-\frac{1}{2}\Big)^2\Ue^2+|\Ue_\tau|^2\right\}z^{1-2s}\,d\Ha^n+o(\eps).
\end{split}
\end{equation*}

Exploiting these computations, we get from~\eqref{CH:6:equation_mono1}
\bgs{
0&\leq\big((1-\eps)^{n-1}-1\big)\Gs_u(1)+\eps\Ha^{n-2}(\partial^*E\cap\partial B_1)\\
&\qquad\qquad
+\eps c'_{n,s}\int_{(\partial\BaLL_1)^+}
\left\{\Big(s-\frac{1}{2}\Big)^2\Ue^2+|\Ue_\tau|^2\right\}z^{1-2s}\,d\Ha^n+o(\eps).
}
Dividing by $\eps$ and passing to the limit $\eps\to0$ yields
\begin{equation}\label{CH:6:almost_monotone1}\begin{split}
c'_{n,s}\Big\{(1-n)\int_{\BaLL_1^+}&|\nabla\Ue|^2z^{1-2s}\,dX
+\int_{(\partial\BaLL_1)^+}
\Big\{\Big(s-\frac{1}{2}\Big)^2\Ue^2+|\Ue_\tau|^2\Big\}z^{1-2s}\,d\Ha^n
\Big\}\\
&
+(1-n)\Per(E,B_1)+\Ha^{n-2}(\partial^*E\cap\partial B_1)\geq0.
\end{split}
\end{equation}
Notice that
\begin{equation*}
\frac{d}{dr}\frac{\Per(E,B_r)}{r^{n-1}}\Big|_{r=1}-(1-n)\Per(E,B_1)=\frac{d}{dr}\Per(E,B_r)\Big|_{r=1}\geq0,
\end{equation*}
since $\Per(E,B_r)$ is increasing in $r$ and it is deifferentiable at $r=1$ by hypothesis. Actually, by Proposition \ref{CH:6:prop_cones_monot}
we have
\begin{equation}\label{CH:6:monotone_perimeter_claim}
\frac{d}{dr}\Per(E,B_r)\Big|_{r=1}\geq\Ha^{n-2}(\partial^*E\cap\partial B_1).
\end{equation}

Let $I$ denote the first line in~\eqref{CH:6:almost_monotone1}. Then we have
\begin{equation*}\begin{split}
0&
\leq I+(1-n)\Per(E,B_1)+\Ha^{n-2}(\partial^*E\cap\partial B_1)\\
&
=I+\frac{d}{dr}\frac{\Per(E,B_r)}{r^{n-1}}\Big|_{r=1}+\Big(\Ha^{n-2}(\partial^*E\cap\partial B_1)-\frac{d}{dr}\Per(E,B_r)\Big|_{r=1}\Big),
\end{split}
\end{equation*}
and hence, by \eqref{CH:6:monotone_perimeter_claim},
\[I+\frac{d}{dr}\frac{\Per(E,B_r)}{r^{n-1}}\Big|_{r=1}\ge
\frac{d}{dr}\Per(E,B_r)\Big|_{r=1}-\Ha^{n-2}(\partial^*E\cap\partial B_1)\ge0.\]
Therefore
\begin{equation*}\begin{split}
\Gs'_u(1)&=\frac{d}{dr}\frac{\Per(E,B_r)}{r^{n-1}}\Big|_{r=1}+I+c'_{n,s}\int_{(\partial\BaLL_1)^+}\Big\{|\Ue_\nu|^2-\Big(s-\frac{1}{2}\Big)^2\Ue^2\Big\}z^{1-2s}\,d\Ha^n\\
&
\geq c'_{n,s}\int_{(\partial\BaLL_1)^+}\Big\{|\Ue_\nu|^2-\Big(s-\frac{1}{2}\Big)^2\Ue^2\Big\}z^{1-2s}\,d\Ha^n\\
&\qquad\qquad
+\Big(\frac{d}{dr}\Per(E,B_r)\Big|_{r=1}-\Ha^{n-2}(\partial^*E\cap\partial B_1)\Big)
\end{split}\end{equation*}
and
\begin{equation*}\begin{split}
\Phi_u'&(1)=\Gs_u'(1)-\Cf_u'(1)\\
&
\geq c'_{n,s}\int_{(\partial\BaLL_1)^+}\Big\{|\Ue_\nu|^2-\Big(s-\frac{1}{2}\Big)^2\Ue^2
-2\Big(s-\frac{1}{2}\Big)\Ue\,\Ue_\nu
-\Big(s-\frac{1}{2}\Big)(1-2s)\Ue^2
\Big\}z^{1-2s}\,d\Ha^n\\
&
\qquad\qquad\qquad+\Big(\frac{d}{dr}\Per(E,B_r)\Big|_{r=1}-\Ha^{n-2}(\partial^*E\cap\partial B_1)\Big).
\end{split}\end{equation*}
Since
\begin{equation*}
|\Ue_\nu|^2-\Big(s-\frac{1}{2}\Big)^2\Ue^2
-2\Big(s-\frac{1}{2}\Big)\Ue\,\Ue_\nu-\Big(s-\frac{1}{2}\Big)(1-2s)\Ue^2
=\Big(\Ue_\nu-\Big(s-\frac{1}{2}\Big)\Ue\Big)^2,
\end{equation*}
we conclude
\eqlab{\label{CH:6:mono_fund_ineq}
\Phi'_u(1)&\geq
c'_{n,s}\int_{(\partial\BaLL_1)^+}
\Big(\Ue_\nu-\Big(s-\frac{1}{2}\Big)\Ue\Big)^2z^{1-2s}\,d\Ha^n\\
&\qquad\qquad\quad
+\Big(\frac{d}{dr}\Per(E,B_r)\Big|_{r=1}-\Ha^{n-2}(\partial^*E\cap\partial B_1)\Big)\geq0.
}
This proves \eqref{CH:6:claim_monoto}, concluding the proof of the monotonicity of $\Phi_u$.

\bigskip

We are left to prove that if $\Phi_u$ is constant, then $\Ue$ is homogeneous of degree $s-\frac{1}{2}$
in $\BaLL_R^+$ and $E$ is a cone in $B_R$
(the converse is a trivial consequence of the scaling invariance of $\Phi_u$).

First of all, notice that
\[\Phi_u\equiv c\quad\mbox{in }(0,R)\qquad\Longrightarrow\qquad\Phi_u'\equiv0\quad\mbox{in }(0,R),\]
hence from \eqref{CH:6:mono_fund_ineq} we find that
\begin{equation}\label{CH:6:mono_eq_end}
\nabla\Ue(X)\cdot X=\Big(s-\frac{1}{2}\Big)\Ue(X)\qquad\mbox{for a.e. }X\in\BaLL_R^+,
\end{equation}
and
\begin{equation}\label{CH:6:mono_fund_ineq2}
\frac{d}{dr}\Per(E,B_r)=\Ha^{n-2}(\partial^*E\cap\partial B_r)\qquad\mbox{for a.e. }r\in(0,R).
\end{equation}
Equality \eqref{CH:6:mono_eq_end} implies that $\Ue$ is homogeneous of degree $s-\frac{1}{2}$
in $\BaLL_R^+$ (see, e.g., \cite[Lemma 4.2]{DSV}).

Therefore, if we denote
\[f(r):=\Phi_u(r)-\theta_E(r),\]
thanks to the scaling invariance properties we have
\[f\equiv c'\quad\mbox{in }(0,R),\]
and hence
\[\Per(E,B_r)=r^{n-1}\big(\Phi_u(r)-f(r)\big)=r^{n-1}(c-c')\quad\forall\,r\in(0,R),\]
so that the function $\wp(r):=\Per(E,B_r)$ is continuous in $(0,R)$.
Thus, by \eqref{CH:6:mono_fund_ineq2} we obtain that $E$ is a cone in $B_R$, thanks to
Proposition \ref{CH:6:prop_cones_monot}.

This concludes the proof.
\end{proof}


\section{Blow-up sequence and homogeneous minimizers}

This section is devoted to the study of blow-up sequences. We begin by proving the uniform energy estimates of Theorem \ref{CH:6:TH:unif}.
Then we establish a convergence result for minimizing pairs---namely Theorem~\ref{CH:6:TH:convmin}---and finally we study blow-up limits, that is, we prove Theorem~\ref{CH:6:TH:blow} (exploiting also the monotonicity formula).

\subsection{Uniform energy estimates}

In this brief subsection, we provide the proof of the uniform energy estimates satisfied by minimizing pairs of the functional $\F$.

\begin{proof}[Proof of Theorem \ref{CH:6:TH:unif}]
The argument of the proof is the same as in the proof of \cite[Theorem 1.1]{DSV}, with a minor modification needed in order
to replace the fractional perimeter considered there with the classical perimeter.

More precisely, we need to replace formula $(7.7)$ of \cite{DSV} with the corresponding formula for the classical perimeter.
To this end, we consider the set
\[F:=B_1\cup(E\setminus B_1)\]
and notice that
\[\Per(F,B_{3/2})=\Ha^{n-1}(\partial^*F\cap\partial B_1)+\Per(E,B_{3/2}\setminus\overline{B_1})\]
and
\[\Per(E,B_{3/2})=\Per(E,B_1)+\Ha^{n-1}(\partial^*E\cap\partial B_1)+\Per(E,B_{3/2}\setminus\overline{B_1}).\]
Hence
\begin{equation}\label{CH:6:new77}\begin{split}
\Per(F,B_{3/2})&-\Per(E,B_{3/2})=\Ha^{n-1}(\partial^*F\cap\partial B_1) -\Per(E,B_1)-\Ha^{n-1}(\partial^*E\cap\partial B_1)\\
&
\le \Ha^{n-1}(\partial^*F\cap\partial B_1) -\Per(E,B_1)
\le\Ha^{n-1}(\partial B_1)-\Per(E,B_1).
\end{split}
\end{equation}
Then we can conclude the proof by arguing as in the proof of \cite[Theorem 1.1]{DSV}, substituting formula (7.7) there with
formula \eqref{CH:6:new77}, whenever needed.
\end{proof}

\subsection{Convergence of minimizers}\label{CH:6:CoNVoMinSeC}

In this subsection we establish some conditions that ensure the convergence of minimizing pairs---namely we prove Theorem \ref{CH:6:TH:convmin}. We will exploit this result in the particularly important case of blow-up sequences.

In order to prove Theorem~\ref{CH:6:TH:convmin}, we need the following glueing Lemma, which is a modification of
\cite[Lemma~6.2]{DSV} taking into account the classical perimeter
in place of (the extension of) the fractional perimeter.

\begin{lemma}\label{CH:6:superglue}
Let $(u_i,E_i)$ be admissible pairs such that
$\F_{B_R}(u_i,E_i)<+\infty$,
for $i=1,2$, and let $\Ue_i$ be the corresponding extension functions. Let $r\in(0,R)$ be such that
\[
\Ha^{n-1}(\partial^*E_i\cap\partial B_r)=0,\qquad\mbox{for }i=1,2,
\]
define
\[
F:=(E_1\cap B_r)\cup (E_2\setminus B_r),
\]
and fix $\varrho\in(r,R)$.

Then, for every small $\eps>0$, there exists a function $\Vf:\R^{n+1}_+\to\R$
such that $(\Vf,F)$ is an admissible pair with
\[
\Vf=\Ue_2\qquad\mbox{in a neighborhood of }(\partial\BaLL_\varrho)^+,
\]
and such that
\begin{equation}\label{CH:6:glueing_energy_ineq}\begin{split}
\int_{\BaLL_\varrho^+}\big(|\nabla\Vf|^2&-|\nabla\Ue_2|^2\big)z^{1-2s}\,dX\\
&
\le\int_{\BaLL_r^+}\big(|\nabla\Ue_1|^2-|\nabla\Ue_2|^2\big)z^{1-2s}\,dX
+C\eps^{-2}\int_{\BaLL_{r+\eps}^+\setminus\BaLL_{r-\eps}^+}|\Ue_1-\Ue_2|^2z^{1-2s}\,dX\\
&
\qquad\qquad+C\int_{\BaLL_{r+\eps}^+\setminus\BaLL_{r-\eps}^+}
\big(|\nabla\Ue_1|^2+|\nabla\Ue_2|^2\big)z^{1-2s}\,dX,
\end{split}\end{equation}
for some constant $C>0$, and
\begin{equation}\label{CH:6:glueing_perimeter_eq}
\Per(F,B_\varrho)-\Per(E_2,B_\varrho)=\Per(E_1,B_r)-\Per(E_2,B_r)+\Ha^{n-1}\big((E_1\Delta E_2)\cap\partial B_r\big).
\end{equation}
\end{lemma}

\begin{proof}
The construction of the function $\Vf$ and the proof of inequality \eqref{CH:6:glueing_energy_ineq}
are the same as in the proof of \cite[Lemma 6.2]{DSV}.

Equality \eqref{CH:6:glueing_perimeter_eq} follows from \cite[Theorem 16.16]{Maggi}.
\end{proof}

\begin{remark}\label{CH:6:local_min_ext_f}
We remark that, if $(\Ue,E)$ is a minimizing pair for the extended functional $\lf$ in $\BaLL_r^+$, then $(\Ue,E)$
is minimizing also in $\BaLL_\varrho^+$, for every $\varrho\in(0,r]$.
\end{remark}

\begin{proof}[Proof of Theorem~\ref{CH:6:TH:convmin}]
First of all, we remark that \eqref{CH:6:convmin1} follows by arguing
as in \cite[Lemma 8.3]{DSV}.

Now we prove that the pair $(\Ue,E)$ is minimizing in $\BaLL_r^+$ for every $r\in(0,R)$.\\
For this, we first point out that,
thanks to Remark \ref{CH:6:local_min_ext_f}, it is enough to prove that $(\Ue,E)$ is minimizing in $\BaLL_r^+$ for a.e.
$r\in(0,R)$.

Then we show that $(\Ue,E)$ is minimizing in $\BaLL_r^+$ for every $r\in\mathcal G$, where $\mathcal G$ is defined as
the set of all radii $r\in(0,R)$ such that
\begin{equation}\label{CH:6:goodradiuskoala}\begin{split}
&\Ha^{n-1}(\partial^*E_m\cap\partial B_r)=0=\Ha^{n-1}(\partial^*E\cap\partial B_r)\\
&\mbox{and}\quad
\lim_{m\to\infty}\Ha^{n-1}\big((E\Delta E_m)\cap\partial B_r\big)=0.
\end{split}\end{equation}
Notice that $\mathcal L^1\big((0,R)\setminus\mathcal G\big)=0$. Indeed, the first condition in
\eqref{CH:6:goodradiuskoala}
holds true for a.e. $r\in(0,R)$ (see e.g. Remark \ref{CH:6:layer_cake_surf}).
Moreover, notice that by hypothesis
\[
0=\lim_{m\to\infty}\big|(E\Delta E_m)\cap B_R\big|=
\int_0^R\Ha^{n-1}\big((E\Delta E_m)\cap\partial B_r\big)\,dr,
\]
so that also the second condition holds for a.e. $r\in(0,R)$.

Now fix a radius $r\in\mathcal G$ and let $(\Vf,F)$ be an admissible competitor for $(\Ue,E)$ in $\BaLL_r^+$.
In particular, notice that since $E\Delta F\Subset B_r$, by \eqref{CH:6:goodradiuskoala} we have
\begin{equation}\label{CH:6:flying_saucepan}
\Ha^{n-1}(\partial^*F\cap\partial B_r)=0\quad\mbox{ and }
\quad\lim_{m\to\infty}\Ha^{n-1}\big((F\Delta E_m)\cap\partial B_r\big)=0.
\end{equation}
Then we fix a radius $\varrho\in(r,R)$ and we consider the pairs $(\Vf_m,F_m)$ defined by using Lemma \ref{CH:6:superglue},
(with $(\Ue_1,E_1):=(\Vf,F)$ and $(\Ue_2,E_2):=(\Ue_m,E_m)$).\\
Notice that, by construction, the pair $(\Vf_m,F_m)$
is an admissible competitor for the pair $(\Ue_m,E_m)$ in $\BaLL_\varrho^+$. Hence the minimality
of $(\Ue_m,E_m)$ (see Remark \ref{CH:6:local_min_ext_f}) implies
\begin{equation}\label{CH:6:minimini}
c'_{n,s}\int_{\BaLL_\varrho^+}|\nabla\Ue_m|^2z^{1-2s}\,dX+\Per(E_m,B_\varrho)
\le c'_{n,s}\int_{\BaLL_\varrho^+}|\nabla\Vf_m|^2z^{1-2s}\,dX+\Per(F_m,B_\varrho).
\end{equation}
Moreover, by Lemma \ref{CH:6:superglue}, we have
\begin{equation}\label{CH:6:after_glue}\begin{split}
c'_{n,s}&\int_{\BaLL^+_\varrho}|\nabla\Vf_m|^2z^{1-2s}\,dX+\Per(F_m,B_\varrho)
\le c'_{n,s}\int_{\BaLL^+_r}|\nabla\Vf|^2z^{1-2s}\,dX+\Per(F,B_r)+c_m(\eps)\\
&
+c'_{n,s}\int_{\BaLL^+_\varrho}|\nabla\Ue_m|^2z^{1-2s}\,dX+\Per(E_m,B_\varrho)
-c'_{n,s}\int_{\BaLL^+_r}|\nabla\Ue_m|^2z^{1-2s}\,dX-\Per(E_m,B_r),
\end{split}\end{equation}
with
\begin{equation*}\begin{split}
c_m(\eps)&:=c'_{n,s}C\left(\eps^{-2}\int_{\BaLL_{r+\eps}^+\setminus\BaLL_{r-\eps}^+}|\Vf-\Ue_m|^2z^{1-2s}\,dX\right.\\
&\qquad
\left.+\int_{\BaLL_{r+\eps}^+\setminus\BaLL_{r-\eps}^+}\big(|\nabla\Vf|^2+|\nabla\Ue_m|^2\big)z^{1-2s}\,dX\right)+\Ha^{n-1}\big((F\Delta E_m)\cap\partial B_r\big)\\
&
=C'\left(\eps^{-2}\int_{\BaLL_{r+\eps}^+\setminus\BaLL_{r-\eps}^+}|\Ue-\Ue_m|^2z^{1-2s}\,dX
+\int_{\BaLL_{r+\eps}^+\setminus\BaLL_{r-\eps}^+}\big(|\nabla\Ue|^2+|\nabla\Ue_m|^2\big)z^{1-2s}\,dX\right)\\
&\qquad\qquad+\Ha^{n-1}\big((F\Delta E_m)\cap\partial B_r\big),
\end{split}
\end{equation*}
(where we have used that $\Vf=\Ue$ outside $\BaLL_{1-\eps}$, provided $\eps<\tilde{\eps}$, by definition of competitor).

Putting together \eqref{CH:6:minimini} and \eqref{CH:6:after_glue}, we find
\begin{equation}\label{CH:6:minimgluerel}\begin{split}
c'_{n,s}\int_{\BaLL^+_r}|\nabla\Ue_m|^2&z^{1-2s}\,dX+\Per(E_m,B_r)\\
&
\le c'_{n,s}\int_{\BaLL^+_r}|\nabla\Vf|^2z^{1-2s}\,dX+\Per(F,B_r)
+c_m(\eps).
\end{split}
\end{equation}

We remark that, arguing as in the proof of \cite[Theorem 1.2]{DSV} and recalling \eqref{CH:6:flying_saucepan},
we obtain
\[\lim_{\eps\to0}\,\lim_{m\to\infty}c_m(\eps)=0.\]
Thus, exploiting \eqref{CH:6:convmin1} and the lower semicontinuity of the perimeter, we obtain
\[
c'_{n,s}\int_{\BaLL^+_r}|\nabla\Ue|^2z^{1-2s}\,dX+\Per(E,B_r)
\le c'_{n,s}\int_{\BaLL^+_r}|\nabla\Vf|^2z^{1-2s}\,dX+\Per(F,B_r).
\]
The arbitrariness of the competitor $(\Vf,F)$ implies that $(\Ue,E)$ is minimizing in $\BaLL_r^+$.
\smallskip

We are left to prove \eqref{CH:6:convmin2}. Indeed, we point out that
\[
\big|D\chi_{E_m}\big|\stackrel{\ast}{\rightharpoonup}\big|D\chi_E\big|,\qquad\mbox{in }B_R
\]
implies
\eqref{CH:6:convmin3} (see, e.g., \cite[Remark 21.15]{Maggi}).

In order to prove \eqref{CH:6:convmin2}, we argue as in the proof of \cite[Theorem 21.14]{Maggi}.\\
A key observation is that, thanks to \eqref{CH:6:minimgluerel}, we have the (locally)
uniform boundedness
\begin{equation}\label{CH:6:locunifbd}
\sup_{m\in\mathbb N}\big|D\chi_{E_m}\big|(B_r)
=\sup_{m\in\mathbb N}\Per(E_m,B_r)
\le C(r)<+\infty,\qquad\forall\,r\in(0,R).
\end{equation}
From \eqref{CH:6:locunifbd} and the convergence $\big|(E_m\Delta E)\cap B_R\big|\to0$, we first get
\[
D\chi_{E_m}\stackrel{\ast}{\rightharpoonup}D\chi_E\qquad\mbox{in }B_R,
\]
by \cite[Theorem 12.15]{Maggi}.\\
Now notice that, in order to conclude the proof of \eqref{CH:6:convmin2},
it is enough to show that every subsequence of $\big|D\chi_{E_m}\big|$
has a subsequence which weakly-star converges to $\big|D\chi_E\big|$, in $B_R$.

We begin by remarking that every subsequence of $\big|D\chi_{E_m}\big|$
admits a weakly-star convergent subsequence, in $B_R$.
Indeed, given such a subsequence $\big|D\chi_{E_{m_h}}\big|$, thanks to \eqref{CH:6:locunifbd},
\cite[Theorem 4.33]{Maggi} implies that we can find a subsequence $m_{h_k}$ of $m_h$ such that
\[
\big|D\chi_{E_{m_{h_k}}}\big|\stackrel{\ast}{\rightharpoonup}\mu\qquad\mbox{in }B_R,
\]
for some Radon measure $\mu$ (which, a priori, might depend on the subsequence).

Finally, we claim that if for some subsequence we have
\begin{equation}\label{CH:6:marcondirondero}
\big|D\chi_{E_{m_h}}\big|\stackrel{\ast}{\rightharpoonup}\mu\qquad\mbox{in }B_R,
\end{equation}
for some Radon measure $\mu$,
then
\[
\mu=\big|D\chi_E\big|\qquad\mbox{in }B_R.
\]
In order to prove this claim, we first point out that
\begin{equation}\label{CH:6:pescatore}
\big|D\chi_E\big|\le\mu\qquad\mbox{in }B_R,
\end{equation}
by \cite[Proposition 4.30]{Maggi}.

Next we show that for every $x\in B_R$ we have
\begin{equation}\label{CH:6:criceto}
\mu\big(B_r(x)\big)\le\big|D\chi_E\big|\big(B_r(x)\big)=\Per(E,B_r(x)),\qquad\mbox{for a.e. }r\in(0,R-|x|).
\end{equation}
To prove \eqref{CH:6:criceto}, let $r<R-|x|$ be such that
\begin{equation*}\begin{split}
&\Ha^{n-1}\big(\partial^*E_m\cap\partial B_r(x)\big)=0=\Ha^{n-1}\big(\partial^*E\cap\partial B_r(x)\big)\\
&\mbox{and}\quad
\lim_{m\to\infty}\Ha^{n-1}\big((E\Delta E_m)\cap\partial B_r(x)\big)=0.
\end{split}\end{equation*}
As before, these conditions hold true for a.e. $r<R-|x|$.\\
Given such an $r$, we fix $\varrho\in(r,R-|x|)$ and we
consider the pair $(\Vf_h,F_h)$ defined by using Lemma~\ref{CH:6:superglue},
with $(\Ue_1,E_1):=(\Ue,E)$ and $(\Ue_2,E_2):=(\Ue_{m_h},E_{m_h})$ (up to a traslation).\\
In particular, by definition,
\[
F_h:=(E\cap B_r(x))\cup(E_{m_h}\setminus B_r(x)).
\]
Each pair $(\Vf_h,F_h)$ is an admissible competitor for the pair $(\Ue_{m_h},E_{m_h})$ in $\BaLL_\varrho^+(x,0)$.

Then, arguing as in the first part of the proof, we obtain
\begin{equation*}\begin{split}
c'_{n,s}\int_{\BaLL^+_r(x,0)}|\nabla\Ue_{m_h}|^2&z^{1-2s}\,dX+\Per(E_{m_h},B_r(x))\\
&
\le c'_{n,s}\int_{\BaLL^+_r(x,0)}|\nabla\Ue|^2z^{1-2s}\,dX+\Per(E,B_r(x))
+c_h(\eps),
\end{split}
\end{equation*}
(in place of \eqref{CH:6:minimgluerel}),
that is
\begin{equation}\label{CH:6:commedia_divertente}
\big|D\chi_{E_{m_h}}\big|\big(B_r(x)\big)=\Per(E_{m_h},B_r(x))
\le \Per(E,B_r(x))+\omega_h(\eps),
\end{equation}
where
\[\omega_h(\eps):=c_h(\eps)+c'_{n,s}
\int_{\BaLL^+_r(x,0)}\big(|\nabla\Ue|^2-|\nabla\Ue_{m_h}|^2\big)z^{1-2s}\,dX\]
and $c_h(\eps)$ is defined as before.

Arguing again as in Lemma 8.3 and the proof of Theorem 1.2 of \cite{DSV}, we obtain
\[
\lim_{\eps\to0}\,\lim_{h\to\infty}\omega_h(\eps)=0.
\]

Hence, by \eqref{CH:6:commedia_divertente} and \eqref{CH:6:marcondirondero} (see \cite[Proposition 4.26]{Maggi})
we obtain \eqref{CH:6:criceto}.\\
Therefore, if $x\in B_R$, then by \eqref{CH:6:pescatore} and
\eqref{CH:6:criceto} we find
\begin{equation}\label{CH:6:ognisera}
\big|D\chi_E\big|\big(B_r(x)\big)=\mu\big(B_r(x)\big),\qquad\mbox{for a.e. }r\in(0,R-|x|).
\end{equation}
By the Lebesgue-Besicovitch Theorem (see \cite[Theorem 5.8]{Maggi}), this implies that
\[
\big|D\chi_E\big|=\mu\qquad\mbox{in }B_R.
\]
Indeed, by \eqref{CH:6:ognisera} we have
\[
D_\mu\big|D\chi_E\big|(x)=\lim_{r\to0}\frac{\big|D\chi_E\big|\big(B_r(x)\big)}{\mu\big(B_r(x)\big)}
=1,\qquad\mbox{for }\mu\mbox{-a.e. }x\in\mbox{supp }\mu\cap B_R.
\]
Thus, since by \eqref{CH:6:pescatore} we have $\big|D\chi_E\big|\ll\mu$, we get
\[
\big|D\chi_E\big|=\big(D_\mu\big|D\chi_E\big|\big)\mu=\mu,\qquad\mbox{in }B_R.
\]

This concludes the proof of the claim and hence of \eqref{CH:6:convmin2}.
\end{proof}

\subsection{Blow-up sequence}\label{CH:6:BLOWupSeQSec}

This subsection is concerned with the existence of a blow-up limit and is dedicated to the proof of Theorem \ref{CH:6:TH:blow}.
We wll employ Theorem \ref{CH:6:TH:unif}, Theorem \ref{CH:6:TH:convmin} and also Theorem \ref{CH:6:Monotonicity_teo}.


In order to prove Theorem~\ref{CH:6:TH:blow} under the assumption that $u\in C^{s-\frac12}(B_1)$,
we need
the following estimate, wich improves the corresponding
estimate in~\cite{DSV} (see the first formula in display
on page~4595 there):

\begin{lemma} \label{CH:6:lem:ur}
Let $s\in(1/2,1)$, and let~$u:\R^n\to\R$ be such that $u\in L^2_s(\R^n)$, $u\in C^{s-\frac12}(B_1)$
and $u(0)=0$.
Let also~$u_r$ be as in~\eqref{CH:6:rescaled_stuff}.
Then, $u_r\in L^2_s(\R^n)$ for every $r\in(0,1)$, and
$$\int_{\R^n}\frac{|u_r(y)|^2}{1+|y|^{n+2s}}\,dy \le
C\left(\|u\|_{C^{s-\frac12}(B_r)}^2+(1-r)
\|u\|_{C^{s-\frac12}(B_1)}^2+
r\,\int_{\R^n}\frac{|u(y)|^2}{1+|y|^{n+2s}}\,dy\right),$$
for some~$C=C(n,s)>0$.
\end{lemma}

\begin{proof}
We write
\begin{equation}\label{CH:6:zeresima}
I:=\int_{\R^n}\frac{|u_r(y)|^2}{1+|y|^{n+2s}}\,dy
= I_1 + I_2 + I_3,\end{equation}
where
\begin{eqnarray*}
I_1&:=&\int_{B_1}\frac{|u_r(y)|^2}{1+|y|^{n+2s}}\,dy,\\
I_2&:=&\int_{B_{1/r}\setminus B_1}
\frac{|u_r(y)|^2}{1+|y|^{n+2s}}\,dy\\
{\mbox{and }} \qquad I_3&:=&
\int_{{\mathcal{C}}B_{1/r}}\frac{|u_r(y)|^2}{1+|y|^{n+2s}}\,dy.
\end{eqnarray*}
We start by estimating~$I_1$. For this, we notice that,
for any~$x$, $\tilde x\in B_{1}$, 
\begin{equation}\label{CH:6:9.1}
|u_r(x)-u_r(\tilde x)| = r^{\frac12-s}|u(rx)-u(r\tilde x)|
\le \|u\|_{C^{s-\frac12}(B_r)}|x-\tilde x|^{s-\frac12}.
\end{equation}
Moreover, since~$u(0)=0$,
we have that~$u_r(0)=0$, 
and so~\eqref{CH:6:9.1} implies that
\begin{equation}\label{CH:6:9.2}
|u_r(x)| \le \|u\|_{C^{s-\frac12}(B_r)}|x|^{s-\frac12},
\end{equation}
for any~$x\in B_{1}$.

As a consequence of~\eqref{CH:6:9.2},
\begin{equation}\begin{split}\label{CH:6:prima}
&I_1= \int_{B_1}\frac{|u_r(y)|^2}{1+|y|^{n+2s}}\,dy
\le \, \|u\|_{C^{s-\frac12}(B_r)}^2
\int_{B_1}\frac{|y|^{2s-1}}{1+|y|^{n+2s}}\,dy\\
&\qquad \le \|u\|_{C^{s-\frac12}(B_r)}^2
\int_{B_1}|y|^{2s-1}\,dy \le C\|u\|_{C^{s-\frac12}(B_r)}^2,
\end{split}\end{equation}
for some~$C>0$, possibly depending on~$n$ and~$s$.

To estimate~$I_2$, we exploit the change of variable~$x:=ry$
and we obtain that
\begin{equation}\label{CH:6:jgergerbgr}
I_2 = \int_{B_{1/r}\setminus B_1}
\frac{r^{1-2s}\,|u(ry)|^2}{1+|y|^{n+2s}}\,dy
= \int_{B_{1}\setminus B_r}
\frac{r\,|u(x)|^2}{r^{n+2s}+|x|^{n+2s}}\,dx.\end{equation}
Now, we use that $u\in C^{s-\frac12}(B_1)$ and the fact that~$0\in\partial E$
to see that
$$|u(x)|\le \|u\|_{C^{s-\frac12}(B_1)}|x|^{s-\frac12},$$
for any~$x\in B_1$. Plugging this information into~\eqref{CH:6:jgergerbgr},
we conclude that
\begin{equation}\begin{split}\label{CH:6:seconda}
& I_2\le r \|u\|_{C^{s-\frac12}(B_1)}^2
\int_{B_{1}\setminus B_r}
\frac{|x|^{2s-1}}{r^{n+2s}+|x|^{n+2s}}\,dx\\&\qquad\quad \le 
r \|u\|_{C^{s-\frac12}(B_1)}^2
\int_{B_{1}\setminus B_r}
\frac{|x|^{2s-1}}{|x|^{n+2s}}\,dx
=\omega_n(1-r)\|u\|_{C^{s-\frac12}(B_1)}^2.\end{split}
\end{equation}

It remains to estimate~$I_3$. To this end, we make the
change of variable~$x:=ry$ and we see that
\begin{equation}\begin{split}\label{CH:6:terza}
& I_3 =
\int_{{\mathcal{C}}B_{1/r}}\frac{r^{1-2s}\,|u(ry)|^2
}{1+|y|^{n+2s}}\,dy =
\int_{{\mathcal{C}}B_{1}}\frac{r\,|u(x)|^2
}{r^{n+2s}+|x|^{n+2s}}\,dx\\
&\qquad \le r\,\int_{{\mathcal{C}}B_{1}}\frac{|u(x)|^2
}{|x|^{n+2s}}\,dx \le Cr\,
\int_{{\mathcal{C}}B_{1}}\frac{|u(x)|^2
}{1+|x|^{n+2s}}\,dx\le Cr\int_{\R^n}\frac{|u(x)|^2
}{1+|x|^{n+2s}}\,dx,
\end{split}\end{equation}
for some~$C>0$, possibly depending on~$n$ and~$s$.

Putting together~\eqref{CH:6:prima}, \eqref{CH:6:seconda}
and~\eqref{CH:6:terza}, and recalling~\eqref{CH:6:zeresima},
we obtain the desired estimate.
\end{proof}

We can now complete the proof of Theorem~\ref{CH:6:TH:blow}.

\begin{proof}[Proof of Theorem~\ref{CH:6:TH:blow}]
As a first step we claim that there exist a function $u_0\in C_{\loc}^{s-\frac{1}{2}}(\R^n)$
and a sequence $r_k\searrow0$ such that $u_{r_k}$ converges to $u_0$ locally uniformly in $\R^n$, that is
\begin{equation}\label{CH:6:trace_conv_eq}
\lim_{k\to\infty}\|u_{r_k}-u_0\|_{C^0(B_R)}=0,\qquad\forall\,R>0.
\end{equation}
Indeed, arguing as in \eqref{CH:6:9.1} and \eqref{CH:6:9.2}, we see that if $r<1/R$, then $u_r\in C^{s-\frac{1}{2}}(B_R)$,
with
\begin{equation}\label{CH:6:unif_bded_holder}
\sup_{r<1/R}\|u_r\|_{C^{s-\frac{1}{2}}(B_R)}\le C_R<+\infty.
\end{equation}
More precisely, if $x,\tilde{x}\in B_R$ and $r<1/R$, then
\begin{equation*}
|u_r(x)-u_r(\tilde{x})|\le\|u\|_{C^{s-\frac{1}{2}}(B_1)}|x-\tilde{x}|^{s-\frac{1}{2}}.
\end{equation*}
Hence, since $u_r(0)=r^{\frac{1}{2}-s}u(0)=0$,
\begin{equation}\label{CH:6:growth}
|u_r(x)|\le\|u\|_{C^{s-\frac{1}{2}}(B_1)}|x|^{s-\frac{1}{2}},\qquad\forall\,x\in B_R,
\end{equation}
if $r<1/R$. In particular,
\[
\sup_{B_R}|u_r|\le\|u\|_{C^{s-\frac{1}{2}}(B_1)}R^{s-\frac{1}{2}},
\]
for every $r<1/R$, concluding the proof of \eqref{CH:6:unif_bded_holder}.\\
Then we get the claim by Ascoli-Arzel\'a Theorem, via a diagonal argument.

We also point out that from \eqref{CH:6:growth} we obtain
\begin{equation}\label{CH:6:global_growth}
|u_0(x)|\le\|u\|_{C^{s-\frac{1}{2}}(B_1)}|x|^{s-\frac{1}{2}},\qquad\forall\,x\in\R^n.
\end{equation}

\smallskip

The second step consists in showing the convergence of the positivity sets.\\
We begin by recalling that, as observed in Remark \ref{CH:6:tail_energies_rmk}, $u\in L^2_s(\R^n)$. Thus, by Lemma \ref{CH:6:lem:ur} we have
\begin{equation*}
\sup_{r\in(0,1)}\int_{\R^n}\frac{|u_r(y)|^2}{1+|y|^{n+2s}}\,dy\le\Lambda<+\infty,
\end{equation*}
with $\Lambda=\Lambda(n,s,u)>0$.
Next we recall that, thanks to Remark \ref{CH:6:scaling_minimality}, the pair $(u_r,E_r)$ is minimal in $B_{1/r}$
and hence also in $B_{2R}$, if $r<1/2R$.
Therefore, by Theorem \ref{CH:6:TH:unif} we obtain
\begin{equation*}\begin{split}
\iint_{\R^{2n}\setminus(\Co B_R)^2}&\frac{|u_r(x)-u_r(y)|^2}{|x-y|^{n+2s}}\,dx\,dy+\Per(E_r,B_R)\\
&
\le C_R\left(1+\int_{\R^n}\frac{|u_r(y)|^2}{1+|y|^{n+2s}}\,dy\right)\le C_R(1+\Lambda)<+\infty.
\end{split}\end{equation*}
In particular, we have
\[
\sup_{k\in\mathbb N}\Per(E_{r_k},B_R)\le C_R(1+\Lambda)<+\infty,\qquad\forall\,R>0.
\]
Thus by compactness (see, e.g., \cite[Corollary 12.27]{Maggi}), up to a subsequence, we get
\[
\chi_{E_{r_k}}\to\chi_{E_0},\qquad\mbox{ both in }L^1_{\loc}(\R^n)\mbox{ and a.e. in }\R^n,
\]
for some set $E_0\subseteq\R^n$ of locally finite perimeter.
Arguing as in the end of the proof of Lemma \ref{CH:6:existence}, we see that $(u_0,E_0)$ is an admissible pair.

\smallskip

As a third step, let $\Ue_r$ and $\Ue_0$ be the extension functions of $u_r$ and $u_0$ respectively.\\
We claim that
\begin{equation}\label{CH:6:unif_conv_blow}
\lim_{k\to\infty}\|\Ue_{r_k}-\Ue_0\|_{L^\infty(\Qc_R)}=0,\qquad\forall\,R>0.
\end{equation}
We first remark that if $w_k:=u_{r_k}-u_0$ and $\We_k$ is the extension function of $w_k$,
then
\[
\We_k=\Ue_{r_k}-\Ue_0.
\]
Hence, by \cite[Lemma 3.1]{DV-onephase}
we find
\begin{equation*}
\|\Ue_{r_k}-\Ue_0\|_{L^\infty(\Qc_R)}=
\|\We_k\|_{L^\infty(\Qc_R)}\le
C_R\left(\|w_k\|_{L^\infty(B_{2R})}+\int_{\R^n\setminus B_{2R}}\frac{|w_k(y)|}{|y|^{n+2s}}\,dy\right)
\end{equation*}
By \eqref{CH:6:trace_conv_eq} we know that
\[
\lim_{k\to\infty}\|w_k\|_{L^\infty(B_{2R})}=0.
\]
Hence, in order to prove \eqref{CH:6:unif_conv_blow} we only need to show that
\begin{equation}\label{CH:6:rtyupoi}
\lim_{k\to\infty}\int_{\R^n\setminus B_{2R}}\frac{|w_k(y)|}{|y|^{n+2s}}\,dy=0.
\end{equation}
First of all, we remark that by Lemma \ref{CH:6:lem:ur} and Fatou's Lemma we obtain
\begin{equation*}
\int_{\R^n}\frac{|u_0(y)|^2}{1+|y|^{n+2s}}\,dy\le\liminf_{k\to\infty}
\int_{\R^n}\frac{|u_{r_k}(y)|^2}{1+|y|^{n+2s}}\,dy\le\Lambda,
\end{equation*}
and hence
\begin{equation}\label{CH:6:spaghettiincident}
\int_{\R^n}\frac{|w_k(y)|^2}{1+|y|^{n+2s}}\,dy\le
2\left(\int_{\R^n}\frac{|u_0(y)|^2}{1+|y|^{n+2s}}\,dy
+\int_{\R^n}\frac{|u_{r_k}(y)|^2}{1+|y|^{n+2s}}\,dy\right)
\le4\Lambda,
\end{equation}
for every $k\in\mathbb N$.
We also remark that
\begin{equation}\label{CH:6:cider}
\frac{1}{|y|^{n+2s}}
\le C_R
\frac{1}{1+|y|^{n+2s}},\qquad\forall\,y\in\Co B_{2R}.
\end{equation}
Now let $\varrho>2R$. Then, by Holder's inequality, \eqref{CH:6:cider} and \eqref{CH:6:spaghettiincident},
we obtain
\begin{equation*}\begin{split}
\int_{\R^n\setminus B_{2R}}&\frac{|w_k(y)|}{|y|^{n+2s}}\,dy=
\int_{B_\varrho\setminus B_{2R}}\frac{|w_k(y)|}{|y|^{n+2s}}\,dy+
\int_{\Co B_\varrho}\frac{|w_k(y)|}{|y|^{n+2s}}\,dy\\
&
\le C_{R,\varrho}\|w_k\|_{L^\infty(B_\varrho)}+
\left(\int_{\Co B_\varrho}\frac{|w_k(y)|^2}{|y|^{n+2s}}\,dy\right)^\frac{1}{2}
\left(\int_{\Co B_\varrho}\frac{dy}{|y|^{n+2s}}\right)^\frac{1}{2}\\
&
\le C_{R,\varrho}\|w_k\|_{L^\infty(B_\varrho)}+
\left(C_R\int_{\Co B_\varrho}\frac{|w_k(y)|^2}{1+|y|^{n+2s}}\,dy\right)^\frac{1}{2}
\left(\frac{\omega_n}{2s}\varrho^{-2s}\right)^\frac{1}{2}\\
&
\le C_{R,\varrho}\|w_k\|_{L^\infty(B_\varrho)}+2\left(\frac{\omega_nC_R\Lambda}{2s}\right)^\frac{1}{2}\varrho^{-s}.
\end{split}
\end{equation*}
By \eqref{CH:6:trace_conv_eq}, passing to the limit $k\to\infty$ yields
\[
\limsup_{k\to\infty}
\int_{\R^n\setminus B_{2R}}\frac{|w_k(y)|}{|y|^{n+2s}}\,dy
\le2\left(\frac{\omega_nC_R\Lambda}{2s}\right)^\frac{1}{2}\varrho^{-s},
\]
for every $\varrho>2R$. Then, passing to the limit $\varrho\to\infty$ proves \eqref{CH:6:rtyupoi}
and hence also \eqref{CH:6:unif_conv_blow}.

\smallskip

The final step consists in showing that $(u_0,E_0)$ is a minimizing cone.\\
We first remark that
 $\Ue_0$ is continuous in $\overline{\R^{n+1}_+}$.
 This can be proved by arguing as in the proof of \cite[Theorem 1.3]{DSV}, by exploiting
 \eqref{CH:6:trace_conv_eq} and \eqref{CH:6:global_growth}.

Now we can apply Theorem \ref{CH:6:TH:convmin} to conclude that the pair $(\Ue_0,E_0)$ is
minimizing in $\BaLL^+_R$, for every $R>0$ and hence,
by Proposition \ref{CH:6:Local_energy_prop}, the pair $(u_0,E_0)$ is minimizing in $B_R$,
for every $R>0$.

We are left to show that $(u_0,E_0)$ is a cone. For this, we are going to use
Theorem \ref{CH:6:Monotonicity_teo}.\\
Since $\Phi_u$ is monotone in $(0,1)$, there exists the limit
\begin{equation}\label{CH:6:gianna}
\lim_{r\to0}\Phi_u(r)=:\Phi\in\R.
\end{equation}
Now, if $\varrho>0$ is such that
\[
\Ha^{n-1}(\partial^*E_0\cap\partial B_\varrho)=0,
\]
then by \eqref{CH:6:unif_conv_blow}, \eqref{CH:6:convmin1} and \eqref{CH:6:convmin3} we obtain
\[
\lim_{k\to\infty}\Phi_{u_{r_k}}(\varrho)=\Phi_{u_0}(\varrho).
\]
Hence, by \eqref{CH:6:gianna} and the scaling invariance \eqref{CH:6:scaling_wholephi}, we get
\[
\Phi_{u_0}(\varrho)=\lim_{k\to\infty}\Phi_{u_{r_k}}(\varrho)
=\lim_{k\to\infty}\Phi_{u}(r_k\varrho)=\Phi,
\]
that is
\begin{equation}\label{CH:6:hopefully_the_last_eq}
\Phi_{u_0}(\varrho)=\Phi,\qquad\mbox{for a.e. }\varrho>0.
\end{equation}
Since $\Phi_{u_0}$ is increasing in $(0,+\infty)$, \eqref{CH:6:hopefully_the_last_eq}
actually implies that
\[
\Phi_{u_0}\equiv\Phi,\qquad\mbox{in }(0,+\infty).
\]
Therefore, by Theorem \ref{CH:6:Monotonicity_teo} we have that $u_0$ is homogeneous of degree
$s-\frac{1}{2}$ in $\R^n$ and $E_0$ is a cone.

This concludes the proof of Theorem~\ref{CH:6:TH:blow}.
\end{proof}

\section{Regularity of the free boundary when $s<1/2$}\label{CH:6:Regs<12free}

We observe that in the case $s<1/2$ the perimeter is, in some sense, the leading term of the functional $\F$. More precisely, by comparing the energy of a minimizing pair with the energy of a simple competitor, we obtain the following estimate.

\begin{theorem}\label{CH:6:sadestimateteo}
Let $(u,E)$ be a minimizing pair in $\Omega$, with $s\in(0,1/2)$, and suppose that $u\in L^\infty_{\loc}(\Omega)$
Let $x_0\in\Omega$ and let $d:=d(x_0,\partial\Omega)/3$. Let $r\in(0,d]$ and define
\begin{equation}\label{CH:6:triste_competitor}
u_*:=\left\{\begin{array}{cc}
0 & \mbox{in }B_r(x_0),\\
u & \mbox{in }\R^n\setminus B_r(x_0).
\end{array}\right.
\end{equation}
Then
\begin{equation}\label{CH:6:local_en_est_2p}\begin{split}
\Nl(u,B_r(x_0))&\le\Nl(u_*,B_r(x_0))\\
&
\le
2\left(\Per_{2s}(B_1)\|u\|_{L^\infty(B_{2d}(x_0))}^2
+r^{2s}|B_1|C_0\int_{\R^n}\frac{|u(y)|^2}{1+|y|^{n+2s}}\,dy
\right)r^{n-2s}
\end{split}
\end{equation}
where
\begin{equation}\label{CH:6:constant_estimate}
C_0=C_0(s,x_0,d):=\sup\left\{\frac{1+|y|^{n+2s}}{|x-y|^{n+2s}}\;
:\;x\in\overline{B_d(x_0)},y\in\R^n\setminus B_{2d}(x_0)\right\}.
\end{equation}
\end{theorem}

\begin{proof}
First of all, notice that the function $u_*$ defined in \eqref{CH:6:triste_competitor}
is such that
\[
u_*\ge0\quad\mbox{a.e. in }E\qquad\mbox{and}\qquad u_*\le0\quad\mbox{a.e. in }\Co E,
\]
hence $(u_*,E)$ is an admissible pair. Moreover
\[
\mbox{supp}(u_*-u)\Subset\Omega
\]
by definition of $u_*$, so that $(u_*,E)$ is an admissible competitor for $(u,E)$.
Thus, since
\[u=u_*\quad\mbox{in }\Co B_r(x_0),\]
by minimality of $(u,E)$ we get
\[
\Nl(u,B_r(x_0))-\Nl(u_*,B_r(x_0))=
\Nl(u,\Omega)-\Nl(u_*,\Omega)=
\F_\Omega(u,E)-\F_\Omega(u_*,E)
\le0.\]

We recall that
\[
\int_{B_r(x_0)}\int_{\Co B_r(x_0)}\frac{dx\,dy}{|x-y|^{n+2s}}
=\Per_{2s}(B_r(x_0))=r^{n-2s}\Per_{2s}(B_1).
\]
Now we can estimate the energy of $u_*$ as follows:
\begin{equation*}\begin{split}
\Nl&(u_*,B_r(x_0))=2\int_{B_r(x_0)}dx\int_{\Co B_r(x_0)}\frac{|u(y)|^2}{|x-y|^{n+2s}}\,dy\\
&
=2\int_{B_r(x_0)}\Big(
\int_{B_{2d}(x_0)\setminus B_r(x_0)}\frac{|u(y)|^2}{|x-y|^{n+2s}}\,dy+
\int_{\Co B_{2d}(x_0)}\frac{|u(y)|^2}{|x-y|^{n+2s}}\,dy
\Big)dx\\
&
\le2\int_{B_r(x_0)}\Big(\|u\|_{L^\infty(B_{2d}(x_0))}^2\int_{B_{2d}(x_0)\setminus B_r(x_0)}\frac{dy}{|x-y|^{n+2s}}
+C_0\int_{\Co B_{2d}(x_0)}\frac{|u(y)|^2}{1+|y|^{n+2s}}\,dy\Big)dx\\
&
\le2\|u\|_{L^\infty(B_{2d}(x_0))}^2\Per_{2s}(B_r(x_0))
+2C_0|B_r(x_0)|\int_{\R^n}\frac{|u(y)|^2}{1+|y|^{n+2s}}\,dy,
\end{split}\end{equation*}
proving \eqref{CH:6:local_en_est_2p} and concluding the proof of the Theorem.
\end{proof}


Since the nonlocal energy $\Nl(u,B_r)$ of a minimizing pair $(u,E)$ goes to zero at least as a power $r^{n-2s}$, we can prove that the set $E$ is almost minimal---in the sense of \cite{Tamanini}---and hence the free boundary $\partial E$ enjoys some regularity properties.

\begin{proof}[Proof of Theorem \ref{CH:6:THMREGFREEBDARYS12}]
First of all, we can assume that $F$ has finite perimeter in $B_r(x_0)$, otherwise there is nothing to prove.
Now
let $u_*$ be the function defined in \eqref{CH:6:triste_competitor}. Notice that, since
\[
E\Delta F\Subset B_r(x_0),
\]
then, by definition of $u_*$,
\[
u_*\ge0\quad\mbox{a.e. in }F\qquad\mbox{and}\qquad u_*\le0\quad\mbox{a.e. in }\Co F,
\]
so that $(u_*,F)$ is an admissible pair and is actually an admissible competitor for $(u,E)$.

Therefore, the minimality of $(u,E)$ implies
\[
0\ge\F_\Omega(u,E)-\F_\Omega(u_*,F)=\F_{B_r(x_0)}(u,E)-\F_{B_r(x_0)}(u_*,F).
\]
Hence
\begin{equation*}\begin{split}
\Per(E,B_r(x_0))&\le \Per(F,B_r(x_0))+\Nl(u_*,B_r(x_0))-\Nl(u,B_r(x_0))\\
&
\le \Per(F,B_r(x_0))+2
\int_{B_r(x_0)}dx\int_{\Co B_r(x_0)}\frac{|u(y)|^2-|u(x)-u(y)|^2}
{|x-y|^{n+2s}}\,dy\\
&
\le \Per(F,B_r(x_0))+2
\int_{B_r(x_0)}dx\int_{\Co B_r(x_0)}\frac{2|u(x)|\,|u(y)|}
{|x-y|^{n+2s}}\,dy\\
&
\le \Per(F,B_r(x_0))+4\|u\|_{L^\infty(B_d(x_0))}
\int_{B_r(x_0)}dx\int_{\Co B_r(x_0)}\frac{|u(y)|}
{|x-y|^{n+2s}}\,dy.
\end{split}\end{equation*}
Estimating the last double integral as in the proof of Theorem \ref{CH:6:sadestimateteo}, we find
\[
\int_{B_r(x_0)}dx\int_{\Co B_r(x_0)}\frac{|u(y)|}
{|x-y|^{n+2s}}\,dy\le C r^{n-2s},
\]
concluding the proof of \eqref{CH:6:almost_min}.

The claims about the regularity of $\partial E$ follow from classical properties of almost minimal sets---see, e.g., \cite{Tamanini}.
\end{proof}

\section{Dimensional reduction}\label{CH:6:DIMREDUCTIONSECTioN}

In this Section we prove a dimensional reduction result in the style of Federer---namely Theorem \ref{CH:6:DimREDtHM}.
In order to do this, we need
to slightly modify the functonal $\F$ by multiplying $\Nl$ with the dimensional constant
$(c_{n,s}')^{-1}$, so that the corresponding extended functional is ``constant-free''.

More precisely, only in this Section we will redefine
\[\F_\Omega(u,E):=(c_{n,s}')^{-1}\Nl(u,\Omega)+\Per(E,\Omega).\]
We say that an admissible pair $(u,E)$ is minimizing in $\R^n$ if
it minimizes $\F_\Omega$ in any bounded open subset $\Omega\subseteq\R^n$ (in
the sense of Definition \ref{CH:6:min_pair}).

The corresponding extended functional then becomes
\[
\lf_\Omega(\Vf,F):=\int_{\Omega_+}|\nabla\Vf|^2z^{1-2s}\,dX+\Per(F,\Omega_0),
\]
for $\Omega\subseteq\R^{n+1}$.

\begin{proof}[Proof of Theorem \ref{CH:6:DimREDtHM}]
The proof is basically a combination of the proof of \cite[Theorem 2.2]{DSV}
and \cite[Lemma 28.13]{Maggi}.
Before going into the details of the proof, we point out some notation which we use only here.
We denote by
$\Pe(F,\Op)$
the perimeter of a set $F\subseteq\R^{n+1}$ in an open set $\Op\subseteq\R^{n+1}$.\\
We write $\Xe:=(x,x_{n+1},z)$ and, with a slight abuse of notation,
\[
\BaLL_R^+\times(-a,a):=\{\Xe=(x,x_{n+1},z)\in\R^{n+2}\,|\,X=(x,z)\in\BaLL_R^+,\,|x_{n+1}|<a\},
\]
``reversing'' for notational simplicity the order of $x_{n+1}$ and $z$ in the domains.\\
If $\Vf:\R^{n+2}\to\R$, $\Xe\longmapsto\Vf(\Xe)$, we write
$\nabla_\Xe\Vf$ for the ``full'' gradient of $\Vf$ and
\[\nabla_X\Vf:=(\partial_1\Vf,\dots,\partial_n\Vf,\partial_z\Vf).\]
In particular, notice that for every fixed $x_{n+1}$ we have
\begin{equation}\label{CH:6:dimred1}
|\nabla_\Xe\Vf|^2=\sum_{i=1}^{n+1}|\partial_i\Vf|^2+|\partial_z\Vf|^2
\ge\sum_{i=1}^n|\partial_i\Vf|^2+|\partial_z\Vf|^2=|\nabla_X\Vf|^2.
\end{equation}
We also remark that if $\Ue$ and $\Ue^\star$ denote the extension functions of $u$ and $u^\star$ respectively,
then
\begin{equation}\label{CH:6:dimred0}
\Ue^\star(x,x_{n+1},z)=\Ue(x,z).
\end{equation}

\smallskip

We first prove, by slicing, that if $(u,E)$ is minimizing in $\R^n$,
then $(u^\star,E^\star)$ is minimizing in $\R^{n+1}$.\\
Fix $a,R>0$ and let $(\Vf,F)$ be a competitor for $(\Ue^\star,E^\star)$ in $\BaLL_R^+\times(-a,a)$.
For every $|t|<a$ we define the hyperplane slices
\[
\Vf_t(x,z):=\Vf(x,t,z)\qquad\mbox{and}\qquad F_t:=\{x\in\R^n\,|\,(x,t)\in F\}.
\]
By \cite[Theorem 18.11]{Maggi}, the slice $F_t$ has locally finite perimeter in $\R^n$ for
a.e. $t\in(-a,a)$. Moreover, since $F$ is the positivity set of $\Vf$, we have
\[
\Vf_t\big|_{\{z=0\}}\ge0\quad\mbox{a.e. in }F_t\qquad
\mbox{and}\qquad
\Vf_t\big|_{\{z=0\}}\le0\quad\mbox{a.e. in }\R^n\setminus F_t,
\]
for a.e. $t\in(-a,a)$. Furthermore
\[
\mbox{supp}(\Vf_t-\Ue)\Subset \BaLL_R\qquad
\mbox{and}\qquad F_t\Delta E\Subset B_R,
\]
for a.e. $t\in(-a,a)$. Hence $(\Vf_t,F_t)$ is an admissible competitor for $(\Ue,E)$ in $\BaLL_R$,
for a.e. $t\in(-a,a)$ and so the minimality of $(u,E)$ implies that
\begin{equation}\label{CH:6:dimred2}
\int_{\BaLL_R^+}|\nabla_X\Vf_t|^2z^{1-2s}\,dX+\Ha^{n-1}(\partial^*F_t\cap B_R)
\ge \int_{\BaLL_R^+}|\nabla_X\Ue|^2z^{1-2s}\,dX+\Ha^{n-1}(\partial^*E\cap B_R),
\end{equation}
for a.e. $t\in(-a,a)$. By formula $(18.25)$ of \cite{Maggi}, we have
\begin{equation}\label{CH:6:dimred3}\begin{split}
\int_{-a}^a\Ha^{n-1}(\partial^*F_t\cap B_R)\,dt&=\int_{\partial^*F\cap\big(B_R\times(-a,a)\big)}|\pi\nu_E|\,d\Ha^n
\le\Ha^n\big(\partial^*F\cap B_R\times(-a,a)\big)\\
&
=\Pe\big(F,B_R\times(-a,a)\big),
\end{split}\end{equation}
where $\pi:\R^{n+1}=\R^n\times\R\to\R^n$, $\pi(x,x_{n+1}):=x$.

By \eqref{CH:6:dimred1} and \eqref{CH:6:dimred3} we obtain
\begin{equation}\label{CH:6:dimred4}\begin{split}
\int_{\BaLL_R^+\times(-a,a)}&|\nabla_\Xe\Vf|^2z^{1-2s}\,d\Xe+\Pe\big(F,B_R\times(-a,a)\big)\\
&
\ge
\int_{-a}^a\left(\int_{\BaLL_R^+}|\nabla_X\Vf_t|^2z^{1-2s}\,dX+\Ha^{n-1}(\partial^*F_t\cap B_R)\right)dt.
\end{split}\end{equation}
On the other hand, by \eqref{CH:6:dimred0} and formula (28.38) of \cite{Maggi}, we have
\begin{equation}\label{CH:6:dimred5}\begin{split}
\int_{\BaLL_R^+\times(-a,a)}&|\nabla_\Xe\Ue^\star|^2z^{1-2s}\,d\Xe+\Pe\big(E^\star,B_R\times(-a,a)\big)\\
&
=\int_{-a}^a\left(\int_{\BaLL_R^+}|\nabla_X\Ue|^2z^{1-2s}\,dX+\Ha^{n-1}(\partial^*E\cap B_R)\right)dt\\
&
=2a\left(\int_{\BaLL_R^+}|\nabla_X\Ue|^2z^{1-2s}\,dX+\Ha^{n-1}(\partial^*E\cap B_R)\right).
\end{split}\end{equation}
Exploiting \eqref{CH:6:dimred2}, \eqref{CH:6:dimred4} and \eqref{CH:6:dimred5} we finally get
\begin{equation*}\begin{split}
\int_{\BaLL_R^+\times(-a,a)}&|\nabla_\Xe\Vf|^2z^{1-2s}\,d\Xe+\Pe\big(F,B_R\times(-a,a)\big)\\
&
\ge\int_{\BaLL_R^+\times(-a,a)}|\nabla_\Xe\Ue^\star|^2z^{1-2s}\,d\Xe+\Pe\big(E^\star,B_R\times(-a,a)\big).
\end{split}\end{equation*}
This proves that the pair $(\Ue^\star,E^\star)$ is minimizing in $\R^{n+2}_+$ and hence that $(u^\star,E^\star)$
is minimizing in $\R^{n+1}$, as claimed.

\smallskip

Now let $(u^\star,E^\star)$ be minimizing in $\R^{n+1}$ and suppose that $(u,E)$ is not minimizing in $\R^n$.\\
Then we can find $R>0$ and an admissible competitor $(\Vf,F)$ for $(\Ue,E)$ in $\BaLL_R^+$, such that
\begin{equation}\label{CH:6:dimred8}
\int_{\BaLL_R^+}|\nabla\Vf|^2z^{1-2s}\,dX+\Per(F,B_R)+\eps\le
\int_{\BaLL_R^+}|\nabla\Ue|^2z^{1-2s}\,dX+\Per(E,B_R),
\end{equation}
for some $\eps>0$. Now we exploit \cite[Corollary 5.2]{DSV} in order to construct a competitor for
$(\Ue^\star,E^\star)$.

More precisely, fix $a>0$ (which in the end of the argument will be taken arbitrarily large) and let $\mathcal Z:\BaLL_R^+\times(-a,a)\to\R$ be the function
constructed in \cite[Corollary 5.2]{DSV}, from $\mathcal U:=\Ue$ and $\Vf$.
Then define the set
\[
G:=\big(F\times(-a,a)\big)\cup\big(E\times(\R\setminus(-a,a))\big)\subseteq\R^{n+1}
\]
and notice that thanks to (5.8) in \cite{DSV}, the pair $(\mathcal Z,G)$ is an admissible competitor for $(\Ue^\star,E^\star)$
in $\BaLL_R^+\times(-a-1,a+1)$.

Arguing as in Step three of the proof of \cite[Lemma 28.13]{Maggi}, we find
\begin{equation}\label{CH:6:dimred6}\begin{split}
\Pe\big(G,B_R\times(-a-1&,a+1)\big)-\Pe\big(E^\star,B_R\times(-a-1,a+1)\big)\\
&
\le2a\left(\Per(F,B_R)-\Per(E,B_R)\right)+2\Ha^n(B_R).
\end{split}\end{equation}
Moreover, by (5.10) in \cite{DSV} the energy of $\mathcal Z$ is
\begin{equation}\label{CH:6:dimred7}\begin{split}
\int_{\BaLL_R^+\times(-a-1,a+1)}&|\nabla_\Xe\mathcal Z|^2z^{1-2s}\,d\Xe\\
&
=2\int_{\BaLL_R^+\times(a-1,a+1)}|\nabla_\Xe\mathcal Z|^2z^{1-2s}\,d\Xe
+2(a-1)\int_{\BaLL_R^+}|\nabla\Vf|^2z^{1-2s}\,dX,
\end{split}\end{equation}
with
\[
2\int_{\BaLL_R^+\times(a-1,a+1)}|\nabla_\Xe\mathcal Z|^2z^{1-2s}\,d\Xe=:C(\mathcal Z)
\]
independent of $a$ by (5.9) in \cite{DSV}.

Therefore, from \eqref{CH:6:dimred5}, \eqref{CH:6:dimred6} and \eqref{CH:6:dimred7} we obtain
\begin{equation}\label{CH:6:dimred9}\begin{split}
&\int_{\BaLL_R^+\times(-a-1,a+1)}|\nabla_\Xe\mathcal Z|^2z^{1-2s}\,d\Xe+
\Pe\big(G,B_R\times(-a-1,a+1)\big)\\
&\qquad\qquad
-\int_{\BaLL_R^+\times(-a-1,a+1)}|\nabla_\Xe\Ue^\star|^2z^{1-2s}\,d\Xe
-\Pe\big(E^\star,B_R\times(-a-1,a+1)\big)\\
&
\le
2(a-1)\Big(
\int_{\BaLL_R^+}|\nabla\Vf|^2z^{1-2s}\,dX+\Per(F,B_R)-\int_{\BaLL_R^+}|\nabla\Ue|^2z^{1-2s}\,dX-\Per(E,B_R)
\Big)+C,
\end{split}\end{equation}
where
\[
C:=C(\mathcal Z)+2\big(\Per(F,B_R)-\Per(E,B_R)+|B_R|\big)+4\int_{\BaLL_R^+}|\nabla\Ue|^2z^{1-2s}\,dX,
\]
which is independent of $a$.

Finally, by \eqref{CH:6:dimred8} and \eqref{CH:6:dimred9} we get
\begin{equation*}\begin{split}
&\int_{\BaLL_R^+\times(-a-1,a+1)}|\nabla_\Xe\mathcal Z|^2z^{1-2s}\,d\Xe+
\Pe\big(G,B_R\times(-a-1,a+1)\big)\\
&\qquad\qquad
-\int_{\BaLL_R^+\times(-a-1,a+1)}|\nabla_\Xe\Ue^\star|^2z^{1-2s}\,d\Xe
-\Pe\big(E^\star,B_R\times(-a-1,a+1)\big)\\
&
\qquad\le-2(a-1)\eps+C<0,
\end{split}\end{equation*}
provided we take $a$ big enough.

This contradicts the minimality of $(\Ue^\star,E^\star)$, concluding the proof.
\end{proof}

\section{Slicing the perimeter and cones}\label{CH:6:appA}

In this section we collect some (more or less known) results about Caccioppoli sets which we used throughout the chapter.
In particular, we recall the coarea formula (see \cite[Theorem 18.8]{Maggi}), which we then
exploit to  construct a cone starting from a ``spherical slice'' of a Caccioppoli set
and to prove a useful formula to compute the perimeter of such a cone.\\
This construction is used 
in the proof of the monotonicity formula in
Theorem~\ref{CH:6:Monotonicity_teo}.

\begin{theorem}[Coarea formula]
If $M$ is a locally $\Ha^{n-1}$-rectifiable set in $\R^n$ and $u:\R^n\to\R$ is a Lipschitz function, then
\begin{equation}
\int_{\R}\Ha^{n-2}\big(M\cap\{u=t\}\big)\,dt
=\int_M|\nabla^Mu|\,d\Ha^{n-1},
\end{equation} 
where
\begin{equation*}
\nabla^Mu(x)=\nabla u(x)-(\nabla u(x)\cdot\nu_M(x))\nu_M(x)
\end{equation*}
is the tangential gradient of $u$.
In particular, if $g:M\to[-\infty,+\infty]$ is a Borel function such that $g\geq0$, then
\begin{equation}\label{CH:6:coarea2}
\int_{\R}dt\int_{M\cap\{u=t\}}g\,d\Ha^{n-2}=\int_Mg|\nabla^Mu|\,d\Ha^{n-1}.
\end{equation}
\end{theorem}

Now we recall that, as remarked in Section \ref{CH:6:mta_subsec}, given a set $E\subseteq\R^n$
we can always find a set $\tilde{E}$ such that
\[|\tilde{E}\Delta E|=0\]
and
\begin{equation}
E_{int}\subseteq\tilde{E},\qquad E_{ext}\subseteq\Co\tilde{E}\qquad\textrm{and}\qquad\partial\tilde{E}=\partial^-E.
\end{equation}
Such a set $\tilde{E}$ is given e.g. by the set of points of density 1 of $E$, that is
\[E^{(1)}:=\Big\{x\in\R^n\,|\,\exists\,\lim_{r\to0^+}\frac{|E\cap B_r(x)|}{\omega_nr^n}=1\Big\}\]
(see, e.g., Appendix \ref{CH:1:Appendix_meas_th_bdary}). In \cite{Visintin} it is also shown that the measure theoretic boundary $\partial^-E$ has
a nice characterization as
the smallest topological boundary among the topological boundaries in the equivalence class of $E$ in $L^1_{\loc}$,
that is
\begin{equation}\label{CH:6:BORDO}\partial^-E=\bigcap_{|F\Delta E|=0}\partial F=\partial E^{(1)}.\end{equation}
If, furthermore, $E$ is a Caccioppoli set, then $\partial^-E$ is the support of the Radon measure $D\chi_E$,
\begin{equation*}
\partial^-E=\textrm{supp }D\chi_E\end{equation*}
(see, e.g., \cite[Proposition 12.19]{Maggi}).\\
In this sense, the set $E^{(1)}$ is a ``good representative'' for $E$ in its $L_{\loc}^1$ equivalence class.

Recall also that the reduced boundary of a Caccioppoli set $E\subseteq\R^n$
\begin{equation*}
\partial^*E:=\left\{x\in\textrm{supp }D\chi_E {\mbox{ s.t. }}
\exists\lim_{{\varrho}\to0^+}
\frac{D\chi_E(B_{\varrho}(x))}{|D\chi_E|(B_{\varrho}(x))}=:\nu_E(x)\in\mathbb S^{n-1}\right\}
\end{equation*}
is locally $\Ha^{n-1}$-rectifiable by De Giorgi's structure Theorem (see, e.g., \cite[Theorem 15.9]{Maggi}).
The Borel function $\nu_E:\partial^*E\to\mathbb S^{n-1}$ is the (measure theoretic)
outer unit normal to $E$. Also notice that by Lebesgue-Besicovitch differentiation Theorem we have
\begin{equation*}
D\chi_E=\nu_E|D\chi_E|\llcorner\partial^*E.
\end{equation*}
Moreover De Giorgi's structure Theorem also says that
\begin{equation*}
|D\chi_E|=\Ha^{n-1}\llcorner\partial^*E,\qquad D\chi_E=\nu_E\Ha^{n-1}\llcorner\partial^*E,
\end{equation*}
so that, in particular,
\begin{equation}
\Per(E,B)=|D\chi_E|(B)=\Ha^{n-1}(\partial^*E\cap B),
\end{equation}
for any Borel set $B\subseteq\R^n$.


\begin{remark}\label{CH:6:layer_cake_surf}
Let $E\subseteq\R^n$ be a set having finite perimeter in $B_{\tilde R}$ and let $R<\tilde R$.
Using formula~\eqref{CH:6:coarea2} for~$M=\partial^*E$, with $u(x)=|x|$ and
$g=\chi_{B_R}$,
we obtain
\begin{equation*}
\int_0^R\Ha^{n-2}(\partial^*E\cap\partial B_t)\,dt=\int_{\partial^*E\cap B_R}|\nabla^{\partial^*E}u|\,d\Ha^{n-1}
\leq \Per(E,B_R)<+\infty.
\end{equation*}
As a consequence
the function
\[h:r\longmapsto\Ha^{n-2}(\partial^*E\cap\partial B_r)\]
is such that $h\in L^1(0,R)$ and
\begin{equation}\label{CH:6:hausdorff_slice}
\Ha^{n-2}(\partial^*E\cap \partial B_r)<+\infty,
\end{equation}
for a.e. $r>0$. Notice that for any $r$ such that \eqref{CH:6:hausdorff_slice} holds true, we have
\begin{equation}\label{CH:6:haus_slice}
\Ha^{n-1}(\partial^*E\cap\partial B_r)=0.
\end{equation}
Hence \eqref{CH:6:haus_slice} also holds true for a.e. $r>0$.

Furthermore, we remark that since $h\in L^1(0,R)$, a.e. $r\in(0,R)$ is a Lebesgue point for $h$.
\end{remark}


We now recall the following result (see e.g. Lemma~4.2.1 on page~102 of~\cite{AMB}):

\begin{lemma}\label{CH:6:90}
Let~$x\in\Omega$ and let~$E$ be a set of finite perimeter in~$\Omega$.
For a.e.~$\varrho\in\big(0,\,d(x,\partial\Omega)\big)$ there exists a set~$E_\varrho$
which has finite perimeter in~$\Omega$, such that~$E\Delta E_\varrho$ is contained in~$B_\varrho(x)$
and
\begin{equation}\label{CH:6:Ambro_formula}
P\big(E_\varrho,\,\overline{B_\varrho(x)}\big)\le\frac{\varrho}{n-1}\,\frac{d}{d\varrho}
\Per(E,\,{B_\varrho(x)}).
\end{equation}
\end{lemma}

As a matter of fact, taking~$x:=0$
up to a translation,
the set~$E_\varrho$ 
given in Lemma~\ref{CH:6:90} is exactly the cone
defined in~\eqref{CH:6:3.4bis} 
(inside~$\overline{B_\varrho}$), see
the formula in display after~(2.8) on page~104 of~\cite{AMB}.\\
More precisely, we recall that we always suppose that the ``good representative''
of a set is chosen, by taking the points of Lebesgue density~$1$. In this sense,
formula~\eqref{CH:6:3.4bis} has to be interpreted as
$$ E(r):=\{\lambda y\,|\,\lambda>0,\,y\in E^{(1)}\cap\partial B_r\}.$$
Lemma \ref{CH:6:90} then guarantees that for a.e. $r\in\big(0,d(0,\Omega)\big)$ the cone $E(r)$ is a Caccioppoli set.

We also observe that the cone structure of~$E(r)$,
together with~\eqref{CH:6:BORDO},
implies that
\begin{equation}\label{CH:6:slices_of_cone}
\begin{split}
&\partial E(r)\cap\partial B_t=\frac{t}{r}\big(\partial E^{(1)}\cap\partial B_r\big)
=\frac{t}{r}\big(\partial^- E\cap\partial B_r\big)
\\{\mbox{and }}\;
&\partial^*E(r)\cap\partial B_t=\frac{t}{r}\big(\partial^*E\cap\partial B_r\big).
\end{split}\end{equation}
The cone structure of~$E(r)$ also implies that
\begin{equation}\label{CH:6:normal_cone}
x\cdot\nu_{E(r)}(x)=0\qquad\textrm{for }\Ha^{n-1}\textrm{-a.e. }x\in\partial^*E(r),
\end{equation}
see, e.g., \cite[Proposition 28.8]{Maggi}.

With these pieces of information we obtain that:

\begin{prop}\label{CH:6:cones_from_slices}
Let $E\subseteq\R^n$ be a Caccioppoli set. Then for a.e. $r>0$ the cone $E(r)$ is a 
Caccioppoli set and
\begin{equation}\label{CH:6:perimeter_cone}
\Per(E(r),B_\varrho)=\frac{\Ha^{n-2}(\partial^*E\cap\partial B_r)}{(n-1)r^{n-2}}\varrho^{n-1},
\end{equation}
for every $\varrho>0$.
\end{prop}

\begin{proof}
The computation relies on~\eqref{CH:6:slices_of_cone} and~\eqref{CH:6:normal_cone} and uses
the coarea formula with $M=\partial^*E(r)$
and~$u(x)=|x|$, so that $\nabla u(x)=\frac{x}{|x|}$.
Indeed,
\begin{eqnarray*}
&& \Per(E(r),B_{\varrho})=\int_{\partial^*E(r)}\chi_{B_{\varrho}}\,d\Ha^{n-1}
=\int_{\partial^*E(r)}\chi_{B_{\varrho}}|\nabla^{\partial^*E(r)}u|\,d\Ha^{n-1}\\
&&\qquad\qquad
=\int_0^{\varrho}\Ha^{n-2}(\partial^*E(r)\cap\partial B_t)\,dt
=\int_0^{\varrho}\Big(\frac{t}{r}\Big)^{n-2}\Ha^{n-2}(\partial^*E\cap\partial B_r)\,dt\\
&&\qquad\qquad
=\frac{\Ha^{n-2}(\partial^*E\cap\partial B_r)}{(n-1)r^{n-2}}\int_0^{\varrho}\frac{d}{dt}t^{n-1}\,dt,
\end{eqnarray*}
proving~\eqref{CH:6:perimeter_cone}.\end{proof}

\begin{remark}\label{CH:6:cone_competitor}
The same argument shows that if $E\subseteq\R^n$ has
finite perimeter in $B_R$, then for a.e. $r\in(0,R)$ the cone
$E(r)$ is a Caccioppoli set and satisfies formula~\eqref{CH:6:perimeter_cone}.
\end{remark}

We remark that, as a consequence of formulas \eqref{CH:6:Ambro_formula} and \eqref{CH:6:perimeter_cone},
we obtain that
\[\Ha^{n-2}(\partial^*E\cap\partial B_r)\le\frac{d}{dr}\Per(E,B_r),\qquad\mbox{for a.e. }r>0.\]
We now prove this inequality by exploiting the coarea formula.

\begin{prop}\label{CH:6:prop_cones_monot}
Let $E\subseteq\R^n$ be a set having finite perimeter in $B_R$. Then
\begin{equation}\label{CH:6:rad_der_per}
\Ha^{n-2}(\partial^*E\cap\partial B_r)\le\frac{d}{dr}\Per(E,B_r),
\end{equation}
for a.e. $r\in(0,R)$. Moreover, 
the following are equivalent:
\begin{itemize}
\item[(i)] the set $E$ is a cone in $B_R$, i.e. there exists a cone $C\subseteq\R^n$ such that
\[\big|(E\Delta C)\cap B_R\big|=0,\]

\item[(ii)] the function
\[(0,R)\ni r\longmapsto \Per(E,B_r)\]
is continuous and
\begin{equation}\label{CH:6:cone_iff}
\Ha^{n-2}(\partial^*E\cap\partial B_r)=\frac{d}{dr}\Per(E,B_r),\qquad\mbox{for a.e. }r\in(0,R).
\end{equation}
\end{itemize}
\end{prop}

\begin{proof}
We define the functions
\[h(r):=\Ha^{n-2}(\partial^*E\cap\partial B_r)\qquad\mbox{and}\qquad\wp(r):=\Per(E,B_r).\]
Then $h\in L^1(0,R)$ (see Remark \ref{CH:6:layer_cake_surf}) and $\wp$ is differentiable almost everywhere
in $(0,R)$, since it is monotone non-decreasing.
Let
\[\mathcal G:=\left\{ r\in(0,R)\,|\,r\mbox{ is a Lebesgue point of }h\mbox{ and }\exists\,\wp'(r)\right\},\]
and notice that $\mathcal L^1\big((0,R)\setminus\mathcal G\big)=0$.
We also remark that
\[r\in\mathcal G\quad\Longrightarrow\quad\Ha^{n-1}(\partial^*E\cap\partial B_r)=0.\]

We prove that the inequality \eqref{CH:6:rad_der_per} holds true for every $r\in\mathcal G$.
To this end, we use the coarea formula for $\partial^*E$, with $u(x):=|x|$.
Notice that
\[
|\nabla^{\partial^*E}u(x)|=\sqrt{1-\Big(\frac{x}{|x|}\cdot\nu_E(x)\Big)^2}\leq1.
\]
Thus
\begin{equation*}
\begin{split}
\Per(E,B_{r+\eps})&-\Per(E,B_r)=\Per(E,B_{r+\eps}\setminus \overline{B_r})=
\int_{\partial^*E\cap(B_{r+\eps}\setminus \overline{B_r})}d\Ha^{n-1}\\
&
\geq\int_{\partial^*E\cap(B_{r+\eps}\setminus\overline{B_r})}
\sqrt{1-\Big(\frac{x}{|x|}\cdot\nu_E(x)\Big)^2}d\Ha^{n-1}
=\int_r^{r+\eps}\Ha^{n-2}(\partial^*E\cap\partial B_t)\,dt,
\end{split}
\end{equation*}
for every $\eps>0$ small enough.
Since $r\in\mathcal G$, dividing by $\eps$ and passing to the limit $\eps\to0^+$ 
yields~\eqref{CH:6:rad_der_per}.

Now we prove that $(i)$ implies $(ii)$. First of all, notice that since $\lambda C=C$ for every $\lambda>0$,
we have
\[\Per(E,B_\varrho)=\Per(C,B_\varrho)=P\Big(\frac{\varrho}{r}C,\frac{\varrho}{r}B_r\Big)
=\Big(\frac{\varrho}{r}\Big)^{n-1}\Per(C,B_r)=\Big(\frac{\varrho}{r}\Big)^{n-1}\Per(E,B_r),\]
for every $r,\varrho\in(0,R)$. Hence
\[\lim_{\varrho\to r}\Per(E,B_\varrho)=\lim_{\varrho\to r}\Big(\frac{\varrho}{r}\Big)^{n-1}\Per(E,B_r)=\Per(E,B_r),\]
proving that $\wp$ is continuous in $(0,R)$.

Since $E$ is a cone in $B_R$, we have by \cite[Proposition 28.8]{Maggi} that
\[x\cdot\nu_E(x)=0\qquad\mbox{for }\Ha^{n-1}\mbox{-a.e. }\partial^*E\cap B_R.\]
Hence, if $u(x):=|x|$, then we find
\[|\nabla^{\partial^*E}u(x)|=1\qquad\mbox{for }\Ha^{n-1}\mbox{-a.e. }\partial^*E\cap B_R.\]
Therefore, the coarea formula implies that
\[
\Per(E,B_{r+\eps})-\Per(E,B_r)=\int_r^{r+\eps}\Ha^{n-2}(\partial^*E\cap\partial B_t)\,dt,
\]
for every $r\in\mathcal G$ and $\eps>0$ small enough. Dividing by $\eps$ and passing to the limit $\eps\to0^+$ thus proves
\eqref{CH:6:cone_iff}.

We are left to show that $(ii)$ implies $(i)$. To this end, first notice that since $\wp$ is continuous and differentiable a.e. in $(0,R)$,
by the Fundamental Theorem of Calculus we have
\begin{equation}\label{CH:6:batman}
\Per(E,B_r)-\Per(E,B_\varrho)=\int^r_\varrho\frac{d}{dt}\Per(E,B_t)\,dt,
\end{equation}
for every $0<\varrho<r<R$. Then, from \eqref{CH:6:cone_iff} and \eqref{CH:6:batman} we get
\begin{equation}\label{CH:6:ErikLensherr1}
\int_{\partial^*E\cap(B_r\setminus B_\varrho)}\,d\Ha^{n-1}=\Per(E,B_r)-\Per(E,B_\varrho)=\int^r_\varrho\Ha^{n-2}(\partial^*E\cap \partial B_t)\,dt.
\end{equation}
Therefore, by exploiting the coarea formula,  
from \eqref{CH:6:ErikLensherr1} we obtain
\begin{equation*}
\int_{\partial^*E\cap(B_r\setminus B_\varrho)}
\sqrt{1-\Big(\frac{x}{|x|}\cdot\nu_E(x)\Big)^2}\,d\Ha^{n-1}
=\int_{\partial^*E\cap(B_r\setminus B_\varrho)}\,d\Ha^{n-1},
\end{equation*}
for every $0<\varrho<r<R$.
Thus
\[x\cdot\nu_E(x)=0\qquad\mbox{for }\Ha^{n-1}\mbox{-a.e. }x\in\partial^*E\cap B_R.\]
By \cite[Proposition 28.8]{Maggi}, this implies that $E^{(1)}$ is a cone in $B_R$, concluding the proof.
\end{proof}

\section{The surface density of a Caccioppoli set}

The following Lemma is a variation of \cite[Lemma 5.1]{GuyDavid} and~\cite[Exercises~3.2.4 and~1.3.6]{AMB}.

\begin{lemma}\label{CH:6:guy_lemma}
Let~$\wp:(0,R)\to\R$ be a monotone non-decreasing
function and let~$
\beta\in C^1\big((0,R),\,(0,+\infty)\big)$.
Then
\begin{equation}\label{CH:6:LA1}
\beta(t_2)\wp(t_2)-\beta(t_1)\wp(t_1)=\int_{[t_1,t_2)}\beta(r)\,dD\wp(r)+
\int_{[t_1,t_2)}\beta'(r)\wp(r)\,dr,
\end{equation}
for every $0<t_1<t_2<R$.
Moreover $\wp$ is differentiable a.e. in $(0,R)$ and
\begin{equation}\label{CH:6:LA2}
\beta(t_2)\wp(t_2)-\beta(t_1)\wp(t_1)\ge\int_{t_1}^{t_2}\big[\beta(r)\wp'(r)+\beta'(r)\wp(r)\big]\,dr,
\quad\textrm{for every }0<t_1<t_2<R.
\end{equation}
\end{lemma}

\begin{proof} We start by proving~\eqref{CH:6:LA1}. For this, we define
\begin{equation}\label{CH:6:ALFA}
\alpha:=\beta\wp.\end{equation} By construction, $\alpha\in BV(0,R)$. We also set~$\alpha_\star(t):=D\alpha\big( [0,t)\big)$
and we claim that
\begin{equation}\label{CH:6:DISTR}
{\mbox{the distributional derivative of~$\alpha_\star$ is equal to~$D\alpha$.}}\end{equation}
To check this, we observe that, by Fubini's Theorem, for any~$\phi\in C^\infty_c(0,R)$,
\begin{eqnarray*} &&-\int_{[0,R)} \phi(\tau)\,dD\alpha(\tau)=
\int_{[0,R)} \big(\phi(R)-\phi(\tau)\big)\,dD\alpha(\tau)
=\int_{[0,R)} \left(\int_{[\tau,R)} \phi'(t)\,dt\right)\,dD\alpha(\tau)
\\ &&\qquad=\int_{[0,R)}\left(\int_{[0,t)} \phi'(t)\,dD\alpha(\tau)\right)\,dt=
\int_{[0,R)} \phi'(t)\,D\alpha\big( [0,t)\big)\,dt=
\int_{[0,R)}  \phi'(t)\,\alpha_\star(t)\,dt.
\end{eqnarray*}
This proves~\eqref{CH:6:DISTR}.

Now we claim that there exists~$c\in\R$ such that, a.e.~$t\in(0,R)$,
\begin{equation}\label{CH:6:OIA}
\alpha(t)=c+D\alpha\big( [0,t)\big).
\end{equation}
To this end, we set~$\gamma(t):=\alpha(t)-\alpha_\star(t)$.
Since~$\alpha_\star$ is monotone non-decreasing, we see that~$\gamma\in BV(0,R)$.
Also, by~\eqref{CH:6:DISTR}, we have that the
distributional derivative of~$\gamma$ vanishes identically, hence~$D\gamma=0$ and
therefore~$\gamma$ is constant. This implies~\eqref{CH:6:OIA}, as desired.

Now, from~\eqref{CH:6:OIA}, it follows that
\begin{equation}\label{CH:6:LAPWO} \alpha(t_2)-\alpha(t_1)=
D\alpha\big( [0,t_2)\big)-D\alpha\big( [0,t_1)\big)=
D\alpha\big( [t_1,t_2)\big)
=
\int_{[t_1,t_2)} dD\alpha(t).
\end{equation}
{F}rom this and~\eqref{CH:6:ALFA}, we obtain~\eqref{CH:6:LA1}. Now we prove~\eqref{CH:6:LA2}.
For this, we use the Lebesgue-Besicovitch
Theorem (see, e.g.,~\cite[Theorem~5.8]{Maggi})
to write
\begin{equation}\label{CH:6:QUI} D\wp= \Psi {\mathcal{L}}^1+D^s\wp,\end{equation}
with~$D^s\wp$ is the singular part of~$D\wp$,
that is a measure supported in a set of zero Lebesgue measure, and
(see, e.g.,~\cite[Corollary~5.11]{Maggi})
$$ \Psi(t):=\lim_{{\varrho}\to0^+} \frac{D\wp\big((t-{\varrho},t+{\varrho})\big)}{2{\varrho}}.$$
We define
\bgs{
\mathcal{G}:=\left\{ t\in(0,R) {\mbox{ s.t. }} t {\mbox{ is a Lebesgue point of }}\Psi
{\mbox{ and }}\lim_{{\varrho}\to0^+} \frac{D^s\wp\big(
(t-{\varrho},t+{\varrho})\big)}{{\varrho}}=0\right\}
}
and
\bgs{
{\mathcal{B}}:=\left\{ t\in(0,R) {\mbox{ s.t. }} \lim_{{\varrho}\to0^+} \frac{D^s\wp\big(
(t-{\varrho},t+{\varrho})\big)}{{\varrho}}\ne 0\right\}.
}
Since~$\wp$ is non-decreasing, we have that
$$ {\mathcal{B}}=\left\{ t\in(0,R) {\mbox{ s.t. }} \lim_{{\varrho}\to0^+} \frac{D^s\wp\big(
(t-{\varrho},t+{\varrho})\big)}{{\varrho}}> 0\right\},$$
hence~${\mathcal{B}}$ is a subset of the support of~$D^s\wp$, and so it has zero Lebesgue measure.
Consequently,
\begin{equation}\label{CH:6:THS0}
{\mbox{$\mathcal{G}$ has full Lebesgue measure in~$(0,R)$.}}
\end{equation}
Now we claim that
\begin{equation}\label{CH:6:THS}
{\mbox{for any~$t\in\mathcal{G}$, the function~$\wp$ is differentiable at~$t$ and~$\wp'(t)=\Psi(t)$.}}
\end{equation}
To check this, we exploit~\eqref{CH:6:OIA} (here, by choosing~$\beta:=1$) and we write that
\begin{equation}\label{CH:6:LI} \wp(t)=c+D\wp\big( [0,t)\big),\end{equation}
for some~$c\in\R$. Then, by~\eqref{CH:6:QUI}, we infer that
\begin{equation*} D\wp\big( [0,t)\big) = 
\int_{[0,t)}\Psi(\tau)\,d\tau+D^s\wp\big( [0,t)\big),\end{equation*}
and therefore, by~\eqref{CH:6:LI},
$$ \wp(t)=c+\int_{[0,t)}\Psi(\tau)\,d\tau+D^s\wp\big( [0,t)\big).$$
As a consequence,
if~$t\in\mathcal{G}$ we have that
\[ \lim_{\varrho\to0^+} \frac{\wp(t+\varrho)-\wp(t)}{\varrho}=
\lim_{\varrho\to0^+} \frac{1}{\varrho}\left[
\int_{[t,t+\varrho)}\Psi(\tau)\,d\tau+D^s\wp\big( [t,t+\varrho)\big)
\right]=\Psi(t)+0,
\]
and this proves~\eqref{CH:6:THS}.

In view of~\eqref{CH:6:THS0} and~\eqref{CH:6:THS}, we obtain that
\begin{equation*}
{\mbox{the function~$\wp$ is differentiable a.e. in~$(0,R)$, with~$\wp'=\Psi$.}}
\end{equation*}
This and~\eqref{CH:6:QUI} give that
$$ D\wp= \wp' {\mathcal{L}}^1+D^s\wp.$$
Hence, by~\eqref{CH:6:ALFA},
$$ D\alpha=D\beta\wp+\beta D\wp= \beta'\wp {\mathcal{L}}^1
+\beta\big(\wp' {\mathcal{L}}^1+D^s\wp\big)=
(\beta'\wp+\beta\wp') {\mathcal{L}}^1+
\beta D^s\wp
.$$
Accordingly, in view of~\eqref{CH:6:LAPWO},
and using that~$\beta\ge0$,
\begin{eqnarray*}
\alpha(t_2)-\alpha(t_1)&=&
\int_{[t_1,t_2)} 
(\beta'\wp+\beta\wp')(t)\,dt+
\int_{[t_1,t_2)}
\beta (t)\,dD^s\wp(t)\\
&\ge&\int_{[t_1,t_2)} 
(\beta'\wp+\beta\wp')(t)\,dt.\end{eqnarray*}
This completes the proof of~\eqref{CH:6:LA2}.
\end{proof}

In particular, by applying Lemma \ref{CH:6:guy_lemma}
to the ``surface density'' of $F\subseteq\R^n$ in 0,
\[\theta_F(r):=\frac{\Per(F,B_r)}{r^{n-1}},\]
we obtain the following result:

\begin{corollary}\label{CH:6:guy_coroll}
Let $F\subseteq\R^n$ be a set having finite perimeter in $B_R$ and let
\[\wp(r):=\Per(F,B_r),\qquad \theta_F(r):=\frac{\Per(F,B_r)}{r^{n-1}}.\]
Then the function $\theta_F$
is differentiable a.e.
in $(0,R)$, with
\[\theta_F'(r)=r^{1-n}\wp'(r)-(n-1)r^{-1}\theta_F(r)\qquad\textrm{for a.e. }r\in(0,R).\]
Moreover
\begin{equation}\label{CH:6:formula_surf_dens_guy}
\theta_F(t_2)-\theta_F(t_1)\ge\int_{t_1}^{t_2}\theta_F'(r)\,dr,
\quad\textrm{for every }0<t_1<t_2<R.
\end{equation}
\end{corollary}

\end{chapter}

\begin{chapter}[The Phillip Island penguin parade]{The Phillip Island penguin parade
(a mathematical treatment)}\label{CH_PEngUInS}

\minitoc

%

\section{Introduction}\label{CH:7:INTRO}

The goal of this chapter is to provide a simple, but rigorous, mathematical
model which describes the formation of groups of penguins
on the shore at sunset. 


The results that we obtain are the following. First of all, we provide
the construction of a mathematical model
to describe the formation of groups of penguins on the shore
and their march towards their burrows; this model is based
on systems of ordinary differential equations, with a number of degree
of freedom which is variable in time (we show that the model
admits a unique solution, which needs to be appropriately defined).
Then, we give some rigorous mathematical results
which provide sufficient conditions for a group of penguins to reach the burrows.
In addition, we provide some numerical simulations which
show that the mathematical model well predicts, at least
at a qualitative level,
the formation of clusters of penguins and their march towards the burrows;
these simulations are easily implemented by images and videos.

It would be
desirable to have empirical data about the formation of penguins clusters
on the shore and their movements, in order to compare and
adapt the model to experimental
data and possibly give a quantitative description of concrete
scenarios.
\medskip

The methodology used is based on direct observations on site, strict
interactions with experts in biology
and penguin ecology, 
mathematical formulation
of the problem and rigorous deductive arguments, and 
numerical simulations.
\medskip

In this introduction, we will describe the elements which lead
to the construction of the model, presenting its basic features and also its limitations.
Given the interdisciplinary flavor of the subject, it is
not possible to completely split the biological discussion
from the mathematical formulation, but we can mention that
the main mathematical equation is given in formula~\eqref{CH:7:EQ}.
Before~\eqref{CH:7:EQ}, the main information coming from live observations
are presented. After~\eqref{CH:7:EQ}, the mathematical quantities
involved in the equation are discussed and elucidated.
The existence and uniqueness theory for equation~\eqref{CH:7:EQ}
is presented in Section~\ref{CH:7:EXUT}.
Some rigorous mathematical results about equation~\eqref{CH:7:EQ}
are given in Section~\ref{CH:7:HOME}. Roughly speaking,
these are results which give sufficient conditions
on the initial data of the system and on the external environment
for the successful homecoming of the penguins, and their precise
formulation requires the development of the mathematical framework
in~\eqref{CH:7:EQ}.

In Section~\ref{CH:7:VIDEO} we present numerics, images and videos which support our intuition and
set the mathematical model of~\eqref{CH:7:EQ} into
a concrete framework that is easily comparable with the
real-world phenomenon.
\medskip

Prior to this, we think
it is important to describe our experience of
the penguins parade in Phillip Island, both to allow the
reader who is not familiar with the event to concretely take part in it,
and to describe some peculiar environmental aspects which are
crucial to understand our description (for instance,
the weather in Phillip Island is completely different from the Antarctic one,
so many of our considerations are meant to be limited to this particular 
habitat) -- also, our personal experience
in this bio-mathematical adventure is a crucial point, in our opinion,
to describe how
scientific curiosity can trigger academic
activities. 

\subsection{Description of the penguins parade}\label{CH:7:QUIQUI}
An  extraordinary event in the state of Victoria, Australia,
consists in the march of the little penguins
(whose scientific name is {\it Eudyptula minor}) who live
in Phillip Island. At sunset, when it gets too dark for the little penguins to
hunt their food in the sea, they come out
to return to their homes
(which are small cavities in the terrain, that are located
at some dozens of meters from the water edge).
What follows is
the mathematical description that came out
of the observations on site at Phillip Island, enriched by the scientific discussions
we later had with 
penguin ecologists.
\medskip

By watching the penguins parade in
Phillip Island, it
seemed to us that some simple
features appeared in the very unusual pattern followed by the little penguins. First of all, they have the strong tendency to
gather together in a sufficiently large number before starting their
march home. They have the tendency to march on a straight line,
compactly arranged
in a cluster, or group. To make this group, they move back
and forth, waiting for other fellows or even going back to
the sea if no other mate is around.

If a little penguin remains isolated,
some parameters in the model proposed may lead
to a complete stop of the individual. More precisely,
in the model that we propose, there is a term which makes
the velocity vanish. In practice,
this interruption in the penguin's movement is not
due to physical impediments, but rather to the fact that there is no other
penguin in a sufficiently small neighbourhood: in this sense, at a mathematical level,
a quantified version of the notion of ``isolation'' leads the penguin to stop.
\medskip

Of course, from the point of view of ethology,
it would be desirable to have further non-invasive tests to measure
how the situation that we describe
is felt by the penguin
at an emotional level (at the moment, we are not aware of experiments like
this in the literature). Also, it would be highly desirable to have some
precise experiments to determine how many penguins do not manage
to return to their burrows within
a certain time
after dusk and stay either in the water or in the vicinity of the shore.

On one hand, in our opinion, it is likely that rigorous experiments on site will
demonstrate that the phenomenon for which an isolated penguin stops
is rather uncommon, but not completely exceptional, in nature. 
On the other hand, our model is general
enough to take into account
the possibility that a penguin stops its march, and, at a quantitative level,
we emphasized this feature in the pictures of
Section \ref{CH:7:VIDEO} to make the situation visible.

The reader who does not want to take into account the stopping function in the model
can just set this function to be identically equal to~$1$
(the mathematical formulation of this remark will be given
after formula~\eqref{CH:7:PAN}). In this particular case,
our model will still exhibit the formation of groups of penguins moving together.
\medskip

Though no experimental test has been run on the emotive feelings of penguins during their homecoming,
in the parade that we have seen live it indeed happened
that one little penguin remained isolated from the others: even though (s)he was absolutely fit
and no concrete obstacle
was obstructing the motion, (s)he got completely stuck for
half an hour and the staff of the Nature Park had
to go and provide assistance. We stress again that
the fact that the penguin stopped did not seem to
be caused by any physical impediment
(as confirmed to us by the Ranger on site), since no extreme environmental condition
was occurring,
the animal
was not underweight, and was able to come out of the water
and move effortlessly on the shore autonomously
for about~$15$ meters, before suddenly stopping.
\medskip

For a short video (courtesy of Phillip Island Nature Parks)
of the little penguins parade, in which the formation
of groups is rather evident, see e.g. the file {\tt Penguins1.MOV},
available at the webpage\\
\small
\url{https://youtu.be/x488k4n3ip8}
\normalsize
\medskip

The simple features listed above are likely to be a consequence of
the morphological structure of the little penguins and of
the natural environment. As a matter of fact,
little penguins are a marine-terrestrial species.
They are highly efficient swimmers
but possess a rather inefficient form of locomotion on land
(indeed, flightless penguins, as the ones
in Phillip Island, 
waddle, more than walk). At dusk,
about 80 minutes after sunset according to the data collected
in~\cite{CH5},
little penguins return ashore
after their fishing activity in the sea.
Since their bipedal locomotion is slow and rather goofy
(at least from the human subjective perception, but also in comparison with
the velocity or agility that is well known to be
typical of predators in nature), and the
easily
recognizable countershading of the penguins
is likely to make them visible to predators,
the transition
between the marine and terrestrial environment
may be particularly stressful
for the penguins 
(see~\cite{CH6}) 
and this fact is probably related to the
formation of penguins groups 
(see e.g.~\cite{CH1}). 
Thus, in our opinion, the rules
that we have listed may be seen as the outcome
of the difficulty of the little penguins to perform their transition
from a more favourable environment to an habitat in which
their morphology turns out to be suboptimal.
\medskip

At the moment,
there seems to be no complete experimental evidence measuring the subjective
perceptions of the penguins with respect to the surrounding environments. Nevertheless,
given the
swimming
ability of the penguins and the environmental conditions,
one may well conjecture that
an area of high potential danger for a penguin is the one adjacent to the
shore-line,
since this is a habitat which provides little or no shelter,
and it is also in an area
of reduced visibility.
As a matter of fact, to protect the penguins in this critical area next to the water edge,
the Rangers in Phillip Island implemented a control on the
presence of the foxes in the proximity of the shore, with the aim of
limiting the number of possible predators.

\subsection{Comparison with the existing literature}
We observe that, to the best of our knowledge,
there is still no specific
mathematical attempt to describe in a concise
way the penguins parade. The mathematical
literature of penguins has mostly focused on the
description of the heat flow in the penguins feathers
(see~\cite{MR2899094}),
on the numerical analysis
to mark animals for later identification
(see~\cite{MR2926716}),
on the statistics of the 
Magellanic penguins at sea
(see~\cite{MR2718048}),
on the hunting strategies of fishing penguins (see~\cite{fish}),
and
on the isoperimetric arrangement of the 
Antarctic penguins to prevent the heat dispersion
caused by the polar wind and on the crystal
structures and solitary waves produced by such arrangements
(see~\cite{1367-2630-15-12-125022}
and~\cite{MR3235841}). We remark that
the climatic situation in Phillip Island
is rather different from the Antarctic one and, given
the very mild temperatures of the area, we do not think that
heat considerations should affect too much the behaviour and
the moving
strategies of the Victorian little penguins and their tendency
to cluster seems more likely to be a defensive
strategy against possible predators.
\medskip

Though no mathematical formulation of the little penguins parade
has been given till now, a series
of experimental analysis has been recently performed
on the specific environment of Phillip Island.
We recall, in particular,~\cite{CH1}, in which
the association of the little penguins in groups is described, by collecting data
spanning over several years,
\cite{CH2}, in which there is a description of
the effect of fog on the orientation
of the little penguins
(which may actually not come back home
in conditions of poor visibility), \cite{CH3} and~\cite{CH4},
which presents a data analysis
to show the fractal structure
in space and time for the foraging of the little penguins, also
in relation to L\'evy flights and fractional Brownian motions.

For an exhaustive list of publications focused on the behaviour
of the little penguins of Phillip Island, we refer to 
the web page\\
\small
\url{https://www.penguins.org.au/conservation/research/publications/}\\
\normalsize
This pages contains more than~$160$ publications related
to the environment of Phillip Island, with special emphasis on the biology of little penguins.

We recall that there is also
a wide literature from the point of view
of biology and ethology focused on collective
mathematical behaviours, also in terms of
formation of groups and hierarchies (see e.g.~\cite{MX13a}
\cite{MX13b} and~\cite{MX15}).   

The mathematical literature studying the collective behaviour
of animal groups is also rather broad: we mention in
particular~\cite{BAL08},
which studied the local rules of interaction
of individual birds in airborne flocks,~\cite{MR2744705},
which analyzed
the self-organization from a microscopic to a macroscopic scale,~\cite{MR3304346}, which took into account
movements with a speed depending on an additional variable,
and~\cite{BV412}
for different models on
opinion formation within an interacting group.\medskip

We remark that our model is 
specifically
tailored on the
Phillip Island penguins :
for instance, other colonies of penguins, such as those in St Kilda,
exhibit behaviours different from those in Phillip Island, due to the different
environmental conditions, see e.g. the scientific report by~\cite{Giling}
for additional information on
the penguins colony on the St Kilda breakwater.

\subsection{Mathematical formulation}\label{CH:7:4365664565674}

In this section we provide a mathematical description of the
penguins parade, which was described
in Section~\ref{CH:7:QUIQUI}. 
The idea for providing an equation for this parade is to prescribe
that the velocity of a group of penguins which travels in line
is influenced by the natural environment and by the position
of the other visible groups. Anytime a group is formed,
the equation needs to be modified to encode the formation
of this new structure.
The main mathematical notation is described in
Table~\ref{CH:7:t24t}.

\begin{table}\begin{center}
\begin{tabular}{|r|l|}
  \hline
  $p_i(t)$ & one-dimensional position of the $i$th group of penguins
at time~$t$ \\
  \hline
  $w_i(t)$ & number of penguins belonging to the $i$th group of penguins at time~$t$ \\
  \hline
$f$ & function describing the environment (sea, shore,
presence of predators, etc.)\\
  \hline
$\pan_i$ & stopping function\\
  \hline
$\e$ & speed of a solitary penguin in a neutral condition
(may be zero)\\
  \hline
$\vel_i$ & strategic speed
of the $i$th group of penguins \\
$\,$& (depending on the position
of the penguins, on the size of the group and on time)
\\
  \hline
$v$ & speed of ``large'' penguins groups\\
  \hline
$m_i$ & influence of the ``visible'' penguins ahead and behind
on the speed of the $i$th group\\
  \hline
$\sgotica$ & eye-sight of the penguins\\
  \hline
\end{tabular}
\caption{Notation.\label{CH:7:t24t}}\end{center}\end{table}
\bigskip

In further details, to
translate  into a mathematical framework the
simple observations on the penguins behaviour that we
listed in Subsection~\ref{CH:7:QUIQUI},
we propose the following equation:
\begin{equation} \label{CH:7:EQ}
\dot{p}_i(t) = \pan_i\big(p(t),w(t);t\big) \,\Big( \eps+\vel_i\big(
p(t),w(t);t\big)
\Big)+f\big(p_i(t),t\big).\end{equation}
The variable~$t\ge0$ represents time
and~$p(t)$ is a vector valued function of time,
that takes into account the positions of the different groups of penguins.
Roughly speaking, at time~$t$, there are~$n(t)$ groups of penguins,
therefore~$p(t)$ is an array with~$n(t)$ components,
and so we will write
\begin{equation}\label{CH:7:pitt}
p(t)=\big(p_1(t),\dots,p_{n(t)}(t)\big).\end{equation}
We stress that~$n(t)$ may vary in time
(in fact, it will be taken to be piecewise constant),
hence the spatial dimension of the image of~$p$ is also a function of time.
For any~$i\in \{1,\dots,n(t)\}$, the $i$th group of penguins
contains a number of penguins denoted by~$w_i(t)$
(thus, the number of penguins belonging to each group is also
a function of time).
\medskip

In further detail, the following notation is used.
The function~$n:[0,+\infty)\to\N_0$, where $\N_0:=\N\setminus\{0\}$,
is piecewise constant and nonincreasing, namely
there exist a (possibly finite)
sequence~$0=t_0<t_1<\dots<t_j<\dots$ 
and integers~$n_1>\dots>n_j>\dots$
such that
\begin{equation}\label{CH:7:ndit}
{\mbox{$n(t)=n_j\in\N_0$ for any~$t\in(t_{j-1},t_j)$.}}
\end{equation}
%

In this model, for simplicity, the spatial occupancy
of a cluster of penguins coincide with that of a single penguin:
of course, in reality,
there is a small
repulsion playing among the penguins, which cannot stay too close to one another.
This additional complication may also be taken into account in our model,
by enlarging the spatial size of the cluster in dependence of the
numerousness of the penguins in the group. In any case, for practical
purposes, we think it is not too inaccurate to identify a group of penguins with just
a single element, since the scale at which the parade occurs (several dozens of meters)
is much larger than the size of a single penguin (little penguins are only about 30 cm. tall).

We also consider the array~$w(t)=
\big(w_1(t),\dots,w_{n(t)}(t)\big)$.
We assume that~$w_i$ is piecewise constant, namely, $w_i(t)=\bar w_{i,j}$
for any~$t\in(t_{j-1},t_j)$, for some~$\bar w_{i,j}\in\N_0$,
namely the number of little penguins in each group remains
constant, till the next penguins join the group at time~$t_j$ 
(if, for the sake of simplicity,
one wishes to think that initially all the little penguins are
separated one from the other, one may also suppose that~$w_i(t)=1$
for all~$i\in\{1,\dots,n_1\}$ and~$t\in[0,t_1)$).

By possibly renaming the variables, we suppose that the initial position of the groups
is increasing with respect to the index, namely
\begin{equation}\label{CH:7:INI P}
p_1(0)<\dots<p_{n_1}(0).
\end{equation}
The parameter~$\eps\ge0$ represents
a drift velocity of the penguins towards their house,
which is located
at the point~$H\in(0,+\infty)$. 
The parameter~$\e$, from the biological point of view,
represents the fact that each penguin, in a neutral
situation, has a natural tendency
to move towards its burrow. We can also
allow~$\e=0$ in our treatment
(namely, the existence and uniqueness theory
in Section~\ref{CH:7:EXUT} remains unchanged if~$\e=0$
and the rigorous results in Section~\ref{CH:7:HOME}
present cases in which they still hold true when~$\e=0$,
compare in particular with assumptions \eqref{CH:7:jk:11} and~\eqref{CH:7:03pr273r047547358935896}).

For concreteness, if~$p_i(T)=H$ for some~$T\ge0$,
we can set~$p_i(t):=H$ for all~$t\ge T$ and remove~$p_i$ from the equation of motion -- that is, the penguin has safely come back
home.

For any~$i\in\{1,\dots,n(t)\}$,
the quantity~$\vel_i\big(p(t),w(t);t\big)$
represents the strategic velocity of the~$i$th
group of penguins and 
it can be considered as a function with domain varying in time
$$ \vel_i(\cdot,\cdot;t):\R^{n(t)}\times\N^{n(t)}\to \R,$$
i.e.
$$ \vel_i(\cdot,\cdot;t):\R^{n_j}\times\N^{n_j}\to \R
\quad{\mbox{ for any }}t\in(t_{j-1},t_j),$$
and,
for any~$(p,w)=
(p_1,\dots,p_{n(t)},w_1,\dots,w_{n(t)})
\in\R^{n(t)}\times\N^{n(t)}$,
it is of the form
\begin{equation}\label{CH:7:VELOCITY}
\vel_i\big(p,w;t\big):=\Big(1-\mu\big(w_i\big)\Big)
\,m_i\big(p,w;t\big)
+v\mu\big(w_i\big).\end{equation}
In this setting, for any~$(p,w)=
(p_1,\dots,p_{n(t)},w_1,\dots,w_{n(t)})
\in\R^{n(t)}\times\N^{n(t)}$,
we have that
\begin{equation}\label{CH:7:emme i} 
m_i\big(p,w;t\big):=
\sum_{j\in\{1,\dots,n(t)\}} {\rm sign}\,(p_j-p_i)\;w_j\;
\sgotica(|p_i-p_j|),
\end{equation}
where~$\sgotica\in {\rm Lip}([0,+\infty))$ is
nonnegative and nonincreasing and, as usual,
we denoted the ``sign function'' as
$$ \R\ni r\mapsto {\rm sign}\,(r):=\left\{
\begin{matrix}
1 & {\mbox{ if }} r>0,\\
0 & {\mbox{ if }} r=0,\\
-1& {\mbox{ if }} r<0.
\end{matrix}
\right.$$
Also, for any~$\ell\in\N$, we set
\begin{equation}\label{CH:7:MU} \mu(\ell):=\left\{
\begin{matrix}
1 & {\mbox{ if }} \ell\ge\kappa,\\
0 & {\mbox{ if }} \ell\le\kappa-1,
\end{matrix}
\right.
\end{equation}
for a fixed~$\kappa \in\N$, with~$\kappa\ge2$,
and~$v>\eps$.

In our framework, the meaning of the
strategic velocity of the~$i$th
group of penguins is the following.
When the group of penguins is too small
(i.e. it contains less than~$\kappa$ little penguins),
then the term involving~$\mu$ vanishes, thus the
strategic velocity reduces to the term given by~$m_i$;
this term, in turn, takes into account the position of the other
groups of penguins.
That is, each penguin is endowed with a ``eye-sight''
(i.e., the capacity of seeing the other penguins that are ``sufficiently close''
to them), which is
modelled by the function~$\sgotica$
(for instance, if~$\sgotica$ is identically equal to~$1$, then
the penguin has a ``perfect eye-sight'';
if~$\sgotica(r)=e^{-r^2}$, then the penguin sees close
objects much better than distant ones; if~$\sgotica$
is compactly supported, then the penguin does not see too far
objects, etc.). Based on the position of the other mates that (s)he
sees, the penguin has the tendency to move either forward
or backward (the more penguins (s)he sees ahead, the more
(s)he is inclined to move forward, 
the more penguins (s)he sees behind, the more
(s)he is inclined to move backward, and nearby penguins
weight more than distant ones, due to the monotonicity
of~$\sgotica$). This strategic tension coming from the position
of the other penguins is encoded by the function~$m_i$.
\medskip

The eye-sight function can be also considered as a modification of the
interaction model based simply on metric distance.
Another interesting feature
which has been observed in several animal groups (see
e.g.~\cite{BAL08}), 
is the so-called
``topological interaction'' model, in which every agent interacts only with a fixed number
of agents, among the ones which are closer. A modification of the function~$\sgotica$
can also take into account this possibility. It is of course very interesting to investigate
by direct observations how much topological, quantitative and metric considerations
influence the formation and the movement of little penguin clusters.
\medskip

When the group of penguins is sufficiently large
(i.e. it contains at least~$\kappa$ little penguins),
then the term involving~$\mu$ is equal to~$1$;
in this case,
the
strategic velocity is~$v$
(that is, when the group of penguins is 
sufficiently rich in population, its strategy
is to move forward with cruising speed equal to~$v$). 

The function~$\pan_i\big(p(t),w(t);t\big)$
describes the case of extreme isolation of the $i$th individual 
from the rest of the herd. Here,
we take~$\overline d>\underline d>0$,
a nonincreasing
function~$\varphi\in {\rm Lip}(\R,[0,1])$,
with~$\varphi(r)=1$ if~$r\le \underline d$
and~$\varphi(r)=0$ if~$r\ge \overline d$,
and, for any~$\ell\in\N_0$,
\begin{equation}\label{CH:7:OMEGA} \wgotica(\ell) :=
\left\{
\begin{matrix}
1 & {\mbox{ if }} \ell\ge2,\\
0 & {\mbox{ if }} \ell=1,
\end{matrix}
\right.\end{equation}
and we take as stopping function
the function with variable domain
$$ \pan_i(\cdot,\cdot;t):\R^{n(t)}\times\N^{n(t)}\to[0,1],$$
i.e.
$$ \pan_i(\cdot,\cdot;t):\R^{n_j}\times\N^{n_j}\to[0,1]
\quad{\mbox{ for any }}t\in(t_{j-1},t_j),$$
given,
for any~$(p,w)=
(p_1,\dots,p_{n(t)},w_1,\dots,w_{n(t)})
\in\R^{n(t)}\times\N^{n(t)}$, by
\begin{equation} \label{CH:7:PAN}
\pan_i\big(p,w;t\big) :=
\max\Big\{ \wgotica(w_i),\;
\max_{{j\in\{1,\dots,n(t)\}}\atop{j\ne i}}
\varphi\big(|p_i-p_j|\big)
\Big\}.\end{equation}
Here the notation ``${\rm Lip}$''
stands for bounded and Lipschitz continuous functions.
\medskip

The case of~$\varphi$
identically equal to~$1$ can be also comprised 
in our setting. In this case, also $\pan_i$ is identically one
(which corresponds to the case in which the stopping function
has no effect).
\medskip

The stopping function describes the fact that the group may present the tendency
to suddenly
stop. This happens
when the group contains only one element (i.e., $\wgotica_i=0$)
and the other groups are far apart (at distance
larger than~$\overline d$).

Conversely, if the group contains at least two little penguins,
or if there is at least another group sufficiently close
(say at distance smaller than~$\underline d$), then
the group is self-confident,
namely the function~$\pan_i\big(p(t),w(t);t\big)$
is equal to~$1$ and the total 
intentional velocity of the group coincides
with the strategic velocity.

Interestingly, the stopping function~$\pan_i$ may
be independent of the eye-sight function~$\sgotica$: namely
a little penguin can stop if (s)he 
feels too much exposed, even if (s)he can see
other little penguins (for instance, if~$\sgotica$
is identically equal to~$1$, the little penguin always sees the other
members of the herd, still (s)he can stop if they are too far apart).

The function~$f\in{\rm Lip}(\R\times[0,+\infty))$
takes into account the environment. For a neutral environment,
one has that this term vanishes (where neutral means here
that the environment does not favour or penalize the homecoming of the penguins). In practice,
it may take into account the ebb and flow
of the sea on the foreshore (where the little penguins parade
starts), the possible
ruggedness of the terrain, the
presence of predators, etc. (as a variation, one
can consider also a stochastic version of this term).
This environment
function can take into account several characteristics
at the same time.
For example, a possible situation that we wish to model is that
in which
the sea occupies the spatial region~$(-\infty,0)$,
producing waves that are periodic in time,
with frequency~$\varpi$
and amplitude~$\delta$; suppose also that the 
shore is located in the spatial region~$(-\infty,0)$, presenting a steep hill in the region~$(1,2)$
which can slow down the motion of the penguins,
whose
burrows are located at the point~$4$. In this setting, a possible choice of the environment function~$f$
is
$$ \R\times[0,+\infty)\ni(p,t)\longmapsto f(p,t)=\delta\sin(\varpi t+\phi)\,\chi_{(-\infty,0)}(p)- h\,\chi_{(1,2)}(p).$$
In this notation~$h>0$ is a constant that takes into account ``how steep'' the hill
located in
the region~$(1,2)$ is, 
$\phi\in\R$ is an initial phase of the wave in the sea,
and $\chi_E$ is the characteristic function of a set~$E$, namely
$$ \chi_E(x):=\left\{
\begin{matrix}
1 & {\mbox{ if }} x\in E,\\
0& {\mbox{ if }} x\not\in E.
\end{matrix}
\right. $$
Given the interpretations above, equation~\eqref{CH:7:EQ}
tries to comprise the pattern that we described in words and
to set the scheme of motion of the little penguins
into a mathematical framework.

\subsection{Preliminary presentation of the mathematical results}

In this chapter, three main mathematical results will be presented.
First of all, in Section~\ref{CH:7:EXUT}, we provide an existence
and uniqueness theory for the solutions of equation~\eqref{CH:7:EQ}.

{F}rom the mathematical viewpoint, we remark that~\eqref{CH:7:EQ}
does not fall into the classical framework of 
the standard Cauchy initial value problem for
ordinary differential
equations (compare e.g. with formula~(2.3)
and Theorem~2.1 in~\cite{BARBU}), since the right hand side of the equation is not
Lipschitz continuous (and, in fact, it is not even continuous).
This mathematical complication is indeed
the counterpart of the real motion of the little penguins 
in the parade, which have the tendency to change their speed
rather abruptly to maintain contact with the other elements of 
the herd. That is, on view, it does not
seem unreasonable to model, as a simplification, 
the speed of the penguin as
a discontinuous function, to take into account the sudden
modifications of the waddling
according to the position of the other penguins,
with the conclusive aim of gathering together a sufficient
number of penguins in a group which eventually will
march concurrently in the direction of their 
burrows.

Then, in Section~\ref{CH:7:HOME} we provide two
rigorous results which guarantee suitable conditions under which
all the penguins, or some of them, safely return to their burrows.
In Theorem~\ref{CH:7:THM:1} we establish that
if the sum of the drift velocity
and the environmental function is
strictly positive and if there is a time (which can be the initial
time or a subsequent one) for which the group at the end of the line
consists of at least two penguins, then all the penguins
reach their burrows in a finite time, which can be explicitly estimated.

Also, in Theorem~\ref{CH:7:NONAMETH} we prove that
if the sum of the drift and cruise velocities
and of the environmental function is
strictly positive and if there is a time for which
one of the penguins group is sufficiently
numerous, then all the penguins of this group and of the groups
ahead safely return home
in a finite time, which can be explicitly estimated.

Rigorous statements and proofs will be given in Sections~\ref{CH:7:EXUT} and~\ref{CH:7:HOME}.

\subsection{Detailed organization of the chapter}

The mathematical treatment of equation~\eqref{CH:7:EQ} that
we provide in this chapter is the following.

In Section~\ref{CH:7:EXUT},
we provide a notion of solution for which~\eqref{CH:7:EQ} is uniquely
solvable in the appropriate setting. This notion of solution will
be obtained by a ``stop-and-go'' procedure, which
is compatible with the idea that when two (or more)
groups of penguins meet, they form a new, bigger group
which will move coherently in the sequel of the march.

In Section~\ref{CH:7:HOME}, we discuss a couple of
concrete examples in which the penguins are able to safely return home: namely, 
we show that there are
``nice'' conditions in which the strategy of the penguins
allows a successful homecoming.

In Section~\ref{CH:7:VIDEO}, we present a series
of numerical simulations
to compare our mathematical model with the real-world experience.
This part also contains some figures
produced by the numerics.

Several possible structural
generalizations of the model proposed are presented
in Section~\ref{CH:7:S-03385787757}. Furthermore,
the model that we propose can be easily generalized
to a multi-dimensional setting, as discussed in
Section~\ref{CH:7:CA823oa-M}.

The conclusions of our work will be summarized in Section~\ref{CH:7:CA823oa}.

\section{Existence and uniqueness theory for equation~\eqref{CH:7:EQ}}\label{CH:7:EXUT}

We stress that equation~\eqref{CH:7:EQ}
does not lie within the setting of ordinary differential
equations, since the right hand side is not Lipschitz
continuous (due to the discontinuity of the functions~$w$
and~$m_i$, and in fact the right hand side also involves
functions with domain varying in time).
As far as we know, the weak formulations of 
ordinary differential
equations as the ones treated by~\cite{MR1022305} do not take
into consideration the setting of equation~\eqref{CH:7:EQ},
so we briefly discuss here a direct
approach to the existence and uniqueness theory for
such equation. To this end, and to clarify our
direct approach, we present two illustrative examples (see
e.g.~\cite{MR1028776}).

\begin{example}{\rm Setting~$x:[0,+\infty)\to\R$, the ordinary differential
equation
\begin{equation}\label{CH:7:L1} \dot x(t) =\left\{
\begin{matrix}
-1 & {\mbox{ if }} x(t)\ge0,\\
1 & {\mbox{ if }} x(t)<0
\end{matrix}
\right. \end{equation}
is not well posed. Indeed, taking an initial datum~$x(0)<0$,
it will evolve with the formula~$x(t)=t+x(0)$ for any~$t\in[0,-x(0)]$
till it hits the zero value. At that point, equation~\eqref{CH:7:L1}
would prescribe a negative velocity, which becomes contradictory
with the positive velocity prescribed to the negative coordinates.
}\end{example}

\begin{example}\label{CH:7:EXMPL2}{\rm The ordinary differential
equation
\begin{equation}\label{CH:7:L2} \dot x(t) =\left\{
\begin{matrix}
-1 & {\mbox{ if }} x(t)>0,\\
0 & {\mbox{ if }} x(t)=0,\\
1 & {\mbox{ if }} x(t)<0
\end{matrix}
\right. \end{equation}
is similar to the one in~\eqref{CH:7:L1},
in the sense that it does not fit into the standard
theory of ordinary differential equations, due to the lack of
continuity of the right hand side. But, differently from
the one in~\eqref{CH:7:L1}, it can be set into an existence and
uniqueness theory by a simple ``reset'' algorithm.

Namely, taking an initial datum~$x(0)<0$,
the solution
evolves with the formula~$x(t)=t+x(0)$ for any~$t\in[0,-x(0)]$
till it hits the zero value. At that point, equation~\eqref{CH:7:L2}
would prescribe a zero velocity, thus a natural way to continue the
solution is to take~$x(t)=0$ for any~$t\in[-x(0),+\infty)$
(similarly,
in the case of positive 
initial datum~$x(0)>0$, a natural way to continue the solution
is~$x(t)=-t+x(0)$ for any~$t\in[0,x(0)]$
and~$x(t)=0$ for any~$t\in[x(0),+\infty)$).
The basic idea for this continuation method
is to flow the equation according
to the standard Cauchy theory of
ordinary differential equations
for as long as possible, and then, when the classical theory
breaks, ``reset'' the equation with respect of the datum
at the break time (this method is not universal
and indeed it does not work for~\eqref{CH:7:L1},
but it produces a natural global solution for~\eqref{CH:7:L2}).
}\end{example}

In the light of Example~\ref{CH:7:EXMPL2},
we now present a framework in which equation~\eqref{CH:7:EQ}
possesses a unique solution (in a suitable ``reset'' setting).
To this aim, we first notice that 
the initial number of groups of penguins is fixed to be equal to~$n_1$
and each group is given by a fixed number of little penguins packed
together (that is,
the number of little penguins in the $i$th initial group being equal to~$\bar w_{i,1}$
and~$i$ ranges from~$1$ to~$n_1$). So, we set~$\bar w_1:=(\bar w_{1,1},\dots,\bar w_{n_1,1})$
and~$\bar\wgotica_{i,1}=\wgotica(\bar w_{i,1})$, where~$\wgotica$ was defined in~\eqref{CH:7:OMEGA}.
For any~$p=(p_1,\dots,p_{n_1})\in\R^{n_1}$, let also
\begin{equation}\label{CH:7:PANNUOVA}  \pan_{i,1}(p)
:=
\max\Big\{ \bar\wgotica_{i,1},\;
\max_{ {j\in\{1,\dots,n_1\}}\atop{j\ne i}}
\varphi\big(|p_i-p_j|\big)
\Big\}. \end{equation} 
The reader may compare this definition with the one in~\eqref{CH:7:PAN}.
For any~$i\in\{1,\dots,n_1\}$ we also
set
$$\bar\mu_{i,1}:=\mu(\bar w_{i,1}),$$ where~$\mu$ is the function defined in~\eqref{CH:7:MU},
and, for any~$p=(p_1,\dots,p_{n_1})\in\R^{n_1}$,
\begin{equation*} 
\bar m_{i,1}(p):=
\sum_{j\in\{1,\dots,n_1\}} {\rm sign}\,(p_j-p_i)\;\bar w_{j,1}\;
\sgotica(|p_i-p_j|)
.\end{equation*}
This definition has to be compared with~\eqref{CH:7:emme i}. Recalling~\eqref{CH:7:INI P}
we also set
$$ {\mathcal{D}}_1:=\{p=(p_1,\dots,p_{n_1})\in\R^{n_1} {\mbox{ s.t. }}
p_1<\dots<p_{n_1}\}.$$
We remark that if~$p\in{\mathcal{D}}_1$ then
\begin{equation*} 
\bar m_{i,1}(p)=
\sum_{j\in\{i+1,\dots,n_1\}} \bar w_{j,1}\;
\sgotica(|p_i-p_j|)
-\sum_{j\in\{1,\dots,i-1\}} \bar w_{j,1}\;
\sgotica(|p_i-p_j|)
\end{equation*}
and therefore
\begin{equation}\label{CH:7:VELOCITY2}
{\mbox{$\bar m_{i,1}(p)$ is bounded
and Lipschitz for any~$p\in{\mathcal{D}}_1$.}}\end{equation}
Then, we set
\begin{equation*}
\vel_{i,1}(p):=(1-\bar\mu_{i,1})
\,\bar m_{i,1}(p)
+v\bar\mu_{i,1}.\end{equation*}
This definition has to be compared with the one in~\eqref{CH:7:VELOCITY}.
Notice that, in view of~\eqref{CH:7:VELOCITY2}, we have that
\begin{equation}\label{CH:7:VELOCITY3}
{\mbox{$\vel_{i,1}(p)$ is bounded and Lipschitz
for any~$p\in{\mathcal{D}}_1$.}}\end{equation}
So, we set
$$ G_{i,1}(p,t):=\pan_{i,1}(p)\,\big(\eps+\vel_{i,1}(p)\big)+
f(p_i,t).$$
{F}rom~\eqref{CH:7:PANNUOVA} and~\eqref{CH:7:VELOCITY3}, we have that~$G_{i,1}$
is bounded and Lipschitz in~${\mathcal{D}}_1\times[0,+\infty)$. Consequently, from
the global existence and uniqueness of solutions of ordinary differential equations,
we have that there exist~$t_1\in(0,+\infty]$ and a solution~$p^{(1)}(t)=(p^{(1)}_1(t),\dots,p^{(1)}_{n_1}(t))\in{\mathcal{D}}_1$
of the Cauchy problem
\begin{align*}
\begin{cases}
\,\dot p^{(1)}_i (t)= G_{i,1}\big( p^{(1)}(t),\,t\big) \qquad {\mbox{ for }}t\in(0,t_1), \\
\,p^{(1)}(0) \qquad {\mbox{ given in }}{\mathcal{D}}_1 \\
\end{cases}
\end{align*}
and
\begin{equation}\label{CH:7:STOP}
p^{(1)}(t_1)\in\partial {\mathcal{D}}_1,\end{equation} see e.g. 
Theorem 1.4.1 in the book~\cite{MR3244289}.

Notice that, as customary in the mathematical literature, we denoted by~$\partial$ the ``topological
boundary'' of a set. In particular,
\begin{eqnarray*} \partial {\mathcal{D}}_1&=&\{p=(p_1,\dots,p_{n_1})\in\R^{n_1} {\mbox{ s.t. }}
p_1\le\dots\le p_{n_1}\\&&\quad
{\mbox{ and  there exists~$i\in\{1,\dots,n_1-1\}$ s.t. }}p_i=p_{i+1}
\}.\end{eqnarray*}
The idea for studying the Cauchy problem in our framework is thus that,
as long as the trajectory of the system stays in the interior
of the domain~${\mathcal{D}}_1$,
the forcing term remains
uniformly Lipschitz, thus the flow does not develop any singularity.
Hence the trajectory exists
and it is defined up to the time (if any) in which it meets the boundary of
the domain~${\mathcal{D}}_1$, that, in the biological framework, corresponds to
the situation in which two (or more) penguins meet (i.e., they occupy the same position
at the same time). In this case, the standard flow procedure of the ordinary differential
equation is stopped, we will merge the joint penguins
into a common cluster, and then repeat the argument.\medskip

In further detail, the solution of~\eqref{CH:7:EQ} will be taken to be~$p^{(1)}$ in~$[0,t_1)$,
that is, we set~$p(t):=p^{(1)}(t)$ for any~$t\in[0,t_1)$.
We also set that~$n(t):=n_1$ and~$w(t):=(\bar w_{1,1},\dots,\bar w_{n_1,1})$.
With this setting, we have that~$p$ is a solution of equation~\eqref{CH:7:EQ}
in the time range~$t\in(0,t_1)$
with prescribed initial datum~$p(0)$. Condition~\eqref{CH:7:STOP} allows us to perform
our ``stop-and-go'' reset procedure as follows: we denote by~$n_2$
the number of distinct points in the set~$\{p^{(1)}_1(t_1),\dots,p^{(1)}_{n_1}(t_1) \}$.
Notice that~\eqref{CH:7:STOP} says that if~$t_1$ is finite then~$n_2\le n_1-1$
(namely, at least two penguins have reached the same position). In this way,
the set of points~$\{p^{(1)}_1(t_1),\dots,p^{(1)}_{n_1}(t_1) \}$
can be identified by the set of~$n_2$ distinct points, that we denote by~$\{p^{(2)}_1(t_1),\dots,p^{(2)}_{n_2}(t_1) \}$, with the convention that
$$ p^{(2)}_1(t_1)<\dots<p^{(2)}_{n_2}(t_1).$$
For any~$i\in\{1,\dots,n_2\}$, we also set
$$ \bar w_{i,2} := \sum_{{j\in\{1,\dots,n_1\}}\atop{
p^{(1)}_j(t_1)=p^{(2)}_i(t_1)
}} \bar w_{j,1}.$$
This says that the new group of penguins indexed by~$i$
contains all the penguins that have reached that position at time~$t_1$.

Thus, having the ``new number of groups'', that is~$n_2$,
the ``new number of little penguins in each group'', that is~$\bar w_2=(\bar w_{1,2},\dots,\bar w_{n_2,2})$,
and the ``new initial datum'', that is~$p^{(2)}(t_1)=\big(
p^{(2)}_1(t_1),\dots,p^{(2)}_{n_2}(t_1)\big)$, we can solve a new differential equation
with these new parameters, exactly in the same way as before, and keep iterating this process.

Indeed, recursively, we suppose that we have 
found~$t_1<t_2<\dots<t_k$,
$p^{(1)}:[0,t_1]\to\R^{n_1}$, $\dots$, $p^{(k)}:[0,t_k]\to\R^{n_k}$
and~$\bar w_1\in\N_0^{n_1}$, $\dots$,
$\bar w_k\in\N_0^{n_k}$ such that, setting
\begin{eqnarray*}&& p(t):=p^{(j)}(t)\in
{\mathcal{D}}_j, \qquad 
n(t):=n_j\\{\mbox{and }}\quad&&
w(t):=\bar w_{j}\qquad{\mbox{for $t\in [t_{j-1},t_j)$ and~$j\in\{1,\dots,k\}$,}} \end{eqnarray*}
one has that~$p$ solves~\eqref{CH:7:EQ} in each interval~$(t_{j-1},t_j)$
for~$j\in\{1,\dots,k\}$, with the ``stop condition''
$$ p^{(j)}(t_j)\in\partial {\mathcal{D}}_j,$$
where
$$ {\mathcal{D}}_j:=\{p=(p_1,\dots,p_{n_j})\in\R^{n_j} {\mbox{ s.t. }}
p_1<\dots<p_{n_j}\}.$$
Then, since~$
p^{(k)}(t_k)\in\partial {\mathcal{D}}_k$,
if~$t_k$ is finite, we find~$n_{k+1}\le n_k-1$
such that
the set of points~$\{p^{(k)}_1(t_k),\dots,p^{(k)}_{n_k}(t_k) \}$
coincides
with
a set of~$n_{k+1}$ distinct points, that we denote by~$\{p^{(k+1)}_1(t_k),\dots,p^{(k+1)}_{n_k}(t_k) \}$, 
with the convention that
$$ p^{(k+1)}_1(t_k)<\dots<p^{(k+1)}_{n_k}(t_k).$$
For any~$i\in\{1,\dots,n_{k+1}\}$, we set
\begin{equation}\label{CH:7:numero} \bar w_{i,k+1} := \sum_{{j\in\{1,\dots,n_k\}}\atop{
p^{(k)}_j(t_k)=p^{(k+1)}_i(t_k)
}} \bar w_{j,k}.\end{equation}
It is useful
to observe that, in light of~\eqref{CH:7:numero},
$$ \sum_{i\in \{1,\dots,n_{k+1}\} } \bar w_{i,k+1}=
\sum_{i\in \{1,\dots,n_{k}\} } \bar w_{i,k},$$ \label{CH:7:OGG}
which says that the total number of little penguins remains always the same
(more precisely, the sum of all the little penguins in all groups
is constant in time).

Let also~$ \bar\wgotica_{i,k+1}=\wgotica(\bar w_{i,k+1})$.
Then, for any~$i\in\{1,\dots,n_{k+1}\}$ and any~$p=
(p_1,\dots,p_{n_{k+1}})\in\R^{n_{k+1}}$, we set
\begin{equation*}  \pan_{i,k+1}(p)
:=
\max\Big\{ \bar\wgotica_{i,k+1},\;
\max_{ {j\in\{1,\dots,n_{k+1}\}}\atop{j\ne i}}
\varphi\big(|p_i-p_j|\big)
\Big\}. \end{equation*} 
For any~$i\in\{1,\dots,n_{k+1}\}$ we also
define
$$\bar\mu_{i,k+1}:=\mu(\bar w_{i,k+1}),$$ where~$\mu$ is the function defined in~\eqref{CH:7:MU}
and, for any~$p\in\R^{n_{k+1}}$,
\begin{equation*} 
\bar m_{i,k+1}(p):=
\sum_{j\in\{1,\dots,n_{k+1}\}} {\rm sign}\,(p_j-p_i)\;\bar w_{j,k+1}\;
\sgotica(|p_i-p_j|)
.\end{equation*}
We notice that~$
\bar m_{i,k+1}(p)$ is bounded and Lipschitz for any~$p\in
{\mathcal{D}}_{k+1}:=\{p=(p_1,\dots,p_{n_{k+1}})\in\R^{n_{k+1}} {\mbox{ s.t. }}
p_1<\dots<p_{n_{k+1}}\}$.

We also define
\begin{equation*}
\vel_{i,k+1}(p):=(1-\bar\mu_{i,k+1})
\,\bar m_{i,k+1}(p)
+v\bar\mu_{i,k+1}\end{equation*}
and
$$ G_{i,k+1}(p,t):=\pan_{i,k+1}(p)\,\big(\eps+\vel_{i,k+1}(p)\big)+
f(p_i,t).$$
In this way,
we have that~$G_{i,k+1}$
is bounded and Lipschitz in~${\mathcal{D}}_{k+1}\times[0,+\infty)$
and so we find the next solution~$p^{(k+1)}(t)=(p^{(k+1)}_1(t),\dots,p^{(k+1)}_{n_{k+1}}(t))\in {\mathcal{D}}_{k+1}$
in the interval~$(t_k,t_{k+1})$, with~$p^{(k+1)}(t_{k+1})\in\partial
{\mathcal{D}}_{k+1}$, by solving the ordinary differential equation
$$ \dot p^{(k+1)}_i(t)= G_{i,k+1}\big(p^{(k+1)}(t),t\big). $$
This completes the iteration argument
and provides the desired notion of solution for equation~\eqref{CH:7:EQ}.

\section{Examples of safe return home}\label{CH:7:HOME}

Here, we provide some sufficient conditions for
the penguins to reach their home, located at the point~$H$, which is taken to be ``far away with respect
to the initial position of the penguins'', namely we suppose that
$$ H>\max_{i\in\{1,\dots,n(0)\}}p_i(0),$$
and~$\e$ has to be thought sufficiently small.
Let us mention that,
in the parade that we saw live,
one little penguin remained stuck and did not manage to return home -- so,
giving a mathematical treatment of the case in which
the strategy of the penguins turns out to be successful
somehow reassured
us on the fate of the species.

To give a mathematical framework of the notion of homecoming,
we introduce the function
$$ [0,+\infty)\ni t\mapsto {\mathcal{N}}(t):=
\sum_{{j\in\{1,\dots,n(t)\}}\atop{p_j(t)=H}} w_j(t).$$
In the setting of Subsection~\ref{CH:7:4365664565674},
the function~${\mathcal{N}}(t)$ represents the number of penguins
that have safely returned home at
time~$t$.

For counting reasons, we also point out that
the total number of penguins is constant and given by
$$ {\mathcal{M}}:=\sum_{ {j\in\{1,\dots,n(0)\}} } w_j(0)
=\sum_{ {j\in\{1,\dots,n(t)\}} } w_j(t),$$
for any~$t\ge0$.

The first result that we present says that if at some time
the group of penguins that stay
further behind
gathers into a group of at least two elements, then
all the penguins will manage to eventually return home.
The mathematical setting goes as follows:

\begin{theorem}\label{CH:7:THM:1}
Let~$t_o\ge0$ and assume that
\begin{equation}\label{CH:7:jk:11}
\eps+\inf_{(r,t)\in \R\times[t_o,+\infty)}
f(r,t)\ge \iota 
\end{equation}
for some~$\iota>0$, and
\begin{equation}\label{CH:7:jk:12}
w_1(t_o)\ge2.
\end{equation}
Then, there exists~$T\in \left[ t_o, \,t_o+\frac{H-p_1(t_o)}{\iota}\right]$
such that
$$ {\mathcal{N}}(T)={\mathcal{M}}.$$
\end{theorem}

\begin{proof} We observe that~$w_1(t)$ is nondecreasing in~$t$,
by~\eqref{CH:7:numero}, and therefore~\eqref{CH:7:jk:12}
implies that~$w_1(t)\ge2$ for any~$t
\ge t_o$. Consequently, from~\eqref{CH:7:OMEGA},
we obtain that~$ \wgotica(w_1(t))=1$
for any~$t
\ge t_o$. This and~\eqref{CH:7:PAN} give that~$
\pan_1\big(p,w(t);t\big) =1$ for any~$t\ge t_o$ and any~$p\in\R^{n(t)}$.
Accordingly, the equation of motions in~\eqref{CH:7:EQ} 
gives that, for any~$t\ge t_o$,
\begin{equation*}
\dot{p}_1(t) =\eps+\vel_1\big(
p(t),w(t);t\big)+f\big(p_1(t),t\big)\ge \eps+f\big(p_1(t),t\big)\ge
\iota,
\end{equation*}
by~\eqref{CH:7:jk:11}. That is, for any~$j\in\{1,\dots,n(t)\}$,
$$ p_j(t)\ge p_1(t)\ge \min\{ H, \; p_1(t_o)+\iota \,(t-t_o) \},$$
which gives the desired result.
\end{proof}

A simple variation of Theorem~\ref{CH:7:THM:1}
says that if, at some time, a group of little penguins
reaches a sufficiently large size, then all the penguins in this group
(as well as the ones ahead) safely reach their home.
The precise statement (whose proof is similar to the one of
Theorem~\ref{CH:7:THM:1}, up to technical modifications, and is therefore
omitted) goes as follows:

\begin{theorem}\label{CH:7:NONAMETH}
Let~$t_o\ge0$ and assume that
\begin{equation}\label{CH:7:03pr273r047547358935896}
\eps+v+\inf_{(r,t)\in \R\times[t_o,+\infty)}
f(r,t)\ge \iota 
\end{equation}
for some~$\iota>0$, and
\begin{equation*}
w_{j_o}(t_o)\ge\kappa,
\end{equation*}
for some~$j_o\in\{1,\dots,n(t_o)\}$, where~$\kappa$ is defined in~\eqref{CH:7:MU}.

Then, there exists~$T\in \left[ t_o, \,t_o+\frac{H-p_{j_o}(t_o)}{\iota}\right]$
such that
$$ {\mathcal{N}}(T)\ge 
\sum_{ {j\in\{j_o,\dots,n(t_o)\}} } w_j(t_o)
.$$
\end{theorem}

\section{Pictures, videos and numerics}\label{CH:7:VIDEO}

In this section, we present some simple numerical experiments
to facilitate the intuition at the base of the model presented in~\eqref{CH:7:EQ}.
These simulations may actually
show some of the typical treats of the little penguins parade,
such as the
oscillations and sudden change of direction, the gathering of the penguins into clusters and the possibility that some elements of the herd remain isolated,
either
on the land or in the sea.

The possibility that a penguin remains isolated
also in the sea may actually occur in the real-world experience, as demonstrated
by the last penguin in the herd on the video (courtesy of Phillip Island Nature Parks)
named {\tt Penguins2.MOV} available 
online at the webpage\\
\small
\url{https://youtu.be/dVk1uYbH_Xc}
\normalsize

In our simulations, for the sake of simplicity, we considered 20 penguins
returning to their burrows from the shore -- some of the penguins may start
their trip from the sea (that occupies the region below level~$0$
in the simulations) in which waves and currents may affect the
movements of the animals. The pictures that we produce (see Section~\ref{CH:7:APPB}) have the time variable
on the horizontal axis and the space variable on the vertical axis
(with the burrow of the penguins community set at level~$4$
for definiteness). The pictures are,
somehow, self-explanatory. For instance, in Figure~\ref{CH:7:MP5},
we present a case in which, fortunately, all the little penguins
manage to safely return home,
after having gathered into groups:
as a matter of fact, in the first of these pictures
all the penguins safely reach home together at the same time
(after having rescued the first penguin, who stayed still for a long
period due to isolation); on the other hand,
the second of these pictures
shows
that a first group of penguins, which was originated
by the animals that were on the land at the initial time,
reaches home
slightly before the second group of penguins, which was originated
by the animals that were in the sea
at the initial time (notice also that
the motion of the penguins in the sea
appears to be affected by waves and currents).

We also observe a different scenario depicted in Figure~\ref{CH:7:MK1}
(with two different functions to represent the currents in the sea):
in this situation, a big group of $18$ penguins
gathers together (collecting also penguins who were initially in the water)
and safely returns home. Two penguins remain isolated in the water,
and they keep slowly moving towards their final destination (that they
eventually reach after a longer time).

Similarly, in Figure~\ref{CH:7:MK12}, almost all the penguins gather into a single group
and reach home, while two penguins get together in the sea, they come to
the shore and slowly waddle towards their final destination,
and one single penguin remains isolated in the water, moved by the currents.

The situation in Figure~\ref{CH:7:MK4} is slightly different, since
the last penguin at the beginning moves towards the others, but (s)he
does not manage to join the forming group by the time the other penguins
decide to move consistently towards their burrows -- so, unfortunately
this last penguin, in spite of the initial effort, finally remains in the water.
  

With simple modifications of the function $f$, one can also consider the
case in which the waves of the sea change with time and their influence
may become more (or less) relevant for the swimming of the little penguins:
as an example of this feature, see Figures \ref{CH:7:ON-1}
and~\ref{CH:7:ON-3}.

In Figures \ref{CH:7:Newpic1} and \ref{CH:7:Newpic2} we give some examples of what happens when varying the parameters
that we used in the numerics of the other figures.
For example we consider different values of $\kappa$, the parameter which encodes when a group of penguins is
big enough to be self confident and waddle home without being influenced by the other groups of penguins in sight.

By considering small values of $\kappa$ we can represent a strong preference of the penguins to go straight towards their homes, instead of first trying to form a large group.
This situation is depicted in the second picture of Figure \ref{CH:7:Newpic1} where we see that after a few time the penguins form two
distinct small groups and go towards home without trying to form a unique large group together.

On the contrary, considering a large value of $\kappa$ represents the preference of the penguins to gather in a very large group before starting their march towards home, like in the first picture of Figure \ref{CH:7:Newpic1}. This situation could represent for example the penguins being
timorous because of the presence of predators.

We think that the case in which one penguin, or a small number of penguins,
remain(s) in water even after 
the return of the main group is worth of further investigation also by means of concrete experiments.
One possible scenario is that the penguins in the water will just wait long enough
for other penguins to get close to the shore and join them to form a new group; on the other
hand, if all the other penguins have already returned, the few ones remained in the water
will have to accept the risk of returning home isolated from the other conspecifics
and in an unprotected situation, 
and we think that interesting biological features could be detected in this case.

Finally, we recall that
once a group of little penguins is created, then it moves consistently
altogether. This is of course a simplifying assumption,
and it might happen in reality that one or a few penguins
leave a large group after its formation
-- perhaps because one penguin is 
slower than the other penguins
of the group, perhaps because (s)he gets distracted
by other events on the beach, or simply because
(s)he feels too exposed being at the side of the group
and may prefer to form a new group in which (s)he
finds a more central and protected position. 
We plan to describe
this case in detail in a forthcoming project (also possibly 
in light of morphological and social considerations
and taking into account a possible randomness
in the system).

The situation in which one little penguin seems to think about
leaving an already formed group can be observed in the video (courtesy of Phillip Island Nature Parks)
named {\tt Penguins2.MOV} and available online at
\\
\small
\url{https://youtu.be/dVk1uYbH_Xc}\\
\normalsize
(see in particular the behaviour of the
second penguin from the bottom, i.e.
the last penguin of 
the already formed large cluster).
\medskip

We point out that all these pictures have been easily
obtained by short programs
in MATLAB.
\medskip

We describe here the algorithm of
the basic program, with waves of constant size and standard behaviour of all the little penguins. The modified versions (periodic strong waves, tired little penguins and so on) can be easily inferred from it.
		
We take into account~$N$ little penguins,
we set their house at $H=4$ and the sea below the location~$0$. Strong waves can go beyond the location~$0$ in some cases, but in the standard program we just consider normal ones. We take a small $\varepsilon$ to represent the natural predisposition to go home of the little penguins, and we define a constant $\delta=(N+1)\varepsilon$ that we need to define the velocity of the little penguins. We define the waves as $\text{WAVE}=\delta\sin(T)$, where $T$ is the array of times. The speed of the animals is related to the one of waves in such a way that it becomes the strongest just when the little penguins form a group that is big enough.

The program starts with a ``for" loop that counts all the animals in a range near the chosen little penguin. This ``for" loop gives us two values: 
the indicator of
the parameter PAN (short for ``panic'')
and the function W, that represents the number of animals in the same position of the one we are considering. We needed this function since we have seen that when the little penguins form a group that is big enough, they proceed towards their home with a cruise speed that is bigger than it was before. We define this cruise speed as vc (short for ``velocity'')
in the program.

Then we start computing the speed $V$ of the little penguin. If PAN is equal to zero, the little penguin freezes. His velocity is zero if he is on the shore (namely his position is greater or equal than zero), or it is given by the waves if he is in the water. It is worth noting that at each value of time the ``for" loop counts the value of PAN, hence a little penguin can leave the 
stopping
condition 
if he sees some mates and start moving again.
		
If PAN is not zero we have mainly two cases, according to the fact that a big group is formed or not. If this has happened, namely $W>\frac{N}{2}$, then the little penguin we are considering is in the group, so he goes towards home with a cruise speed vc, possibly modified by the presence of waves. If the group is not formed yet, the animal we are considering is surrounded by some mates, but they are not enough to proceed straight home. His speed is positive or negative, namely he moves forward or backward, in dependence of the amount of little penguins that he has ahead of him or behind him. Its speed is given by:
\[
V=\varepsilon+M
\]
where $M$ is the number of penguins ahead of him minus the number of animals behind him multiplied by $\frac{\delta}{N}$, and $\epsilon$ has been defined before. As in the other cases, the speed can be modified by the presence of waves if the position is less than zero.
		
Now that we have computed the speed of the animal, we can obtain his position $P$ after a discrete time interval $t$ by considering $P(k+1)=P(k)+Vt$.
		
The last ``for" loop is done in order to put in the same position
two animals that are closed enough. 
Then we reset the counting variables PAN, W and M and we restart the loop.

For completeness,
we made the source codes of
all the programs
available
on the webpage\\
\small
\url{https://www.dropbox.com/sh/odgic3a0ke5qp0q/AABIMaasAcTwZQ3qKRoB--xra?dl=0}
\normalsize

An example of the code is given in Section~\ref{CH:7:APPA}.
The simplicity of these programs shows that the model in~\eqref{CH:7:EQ}
is indeed very simple to implement numerically, still producing 
sufficiently ``realistic'' results
in terms of cluster formation and cruising speed of the groups.
The parameters in the code
are chosen as examples, producing simulations that show some features similar to those observed on site and in the videos. From one picture to another, what is varying is the initial conditions and the environment function (minor modifications in the code would allow also to change the number of penguins, their eye-sight, the drift and cruise velocities, the stopping function, and also to take into account multi-dimensional cases). 

Also, these pictures can be easily translated into animations. Simple videos
that we have
obtained by these numerics are available from 
the webpage\\
\small
\url{https://www.youtube.com/playlist?list=PLASZVs0A5ReZgEinpnJFat66lo2kIkWTS }
\normalsize

The source codes of the animations are available 
online at\\
\small
\url{https://www.dropbox.com/s/l1z5riqtc8jzxbs/scatter.txt?dl=0 }
\normalsize

\section{Discussion on the model proposed: simplifications,
generalizations and further directions of investigation}\label{CH:7:S-03385787757}

We stress that the model proposed in~\eqref{CH:7:EQ} is of course a
dramatic
simplification of ``reality''.
As often happens in science indeed, several simplifications
have been adopted in order
to allow a rigorous mathematical treatment
and handy numerical computations: nevertheless the model is
already rich enough to detect some specific features
of the little penguins parade, such as the formation of groups, the oscillatory 
waddling of the penguins and the possibility
{of isolated and exposed
individuals}.
Moreover, our model is flexible enough to allow
specific distinctions between the single penguins (for instance,
with minor modifications, one can take into account the possibility
that different penguins have a different eye-sight, that they have a different reaction
to isolation, or that they exhibit some specific social behaviour
that favours the formation of clusters selected by specific characteristics); 
similarly, the modeling of the habitat may also encode different possibilities
(such as the burrows of the penguins being located in different places),
and multi-dimensional models can be also constructed using similar ideas (see
Section~\ref{CH:7:CA823oa-M} for details).

We observe that one can replace the quantities $v,\sgotica,\mu,\kappa,\varphi$ 
with $v_i,\sgotica_i,\mu_i,\kappa_i,\varphi_i$ if one wants to customize these features for every group.
\medskip

Furthermore, natural modifications lead to the possibility that
one or a few penguins may leave an already formed group:
for
instance, rather than forming one single group, 
the model can still consider the penguins of the
cluster as separate elements, each one with its
own peculiar behaviour. At the moment, for simplicity,
we considered here the basic model in which, 
once a cluster is made up,
it keeps moving without losing any of its elements -- we plan to
address in a future project
in detail the case of groups which may also
decrease the number of components, possibly in dependence
of random fluctuations or social considerations among
the members of the group.\medskip


In addition, for simplicity,
in this chapter
we modelled each group to be located at a precise point:
though this is not a
completely unrealistic assumption (given that the scale of
the individual
penguin is much smaller than that of the beach), one can also easily modify
this feature by locating a cluster in a region comparable to its size.\medskip

In future projects, we plan to introduce
other more sophisticated models, 
also taking into account
stochastic oscillations and optimization methods, 
and, in the long run,
to use these models
in a detailed experimental confrontation taking advantage of
the automated monitoring systems which is
under development in Phillip Island.\medskip

The model that we propose here is also flexible enough
to allow quantitative modifications of all the parameters involved.
This is quite important, since these parameters may vary due
to different conditions of the environment. For instance,
the eye-sight of the penguins can be reduced by the
fog (see~\cite{CH2}),
and by the effect of moonlight and artificial light (see~\cite{CH5}).

Similarly, the number of penguins in each group
and the velocity of the herd may vary due to 
structural changes of the beach:
roughly speaking,
from the empirical data, penguins typically gather into groups of~$ 5$--$10$
individuals (but we have also observed
much larger groups forming on the beach)
within 40 second intervals, see~\cite{CH1},
but the way these groups are built varies year by year
and, for instance, the number of individuals which
always gather into the same group changes year by year
in strong dependence with the breeding success
of the season, see again~\cite{CH1}. Also,
tidal phenomena may change the number of little
penguins in each group and the velocity of the group,
since the change of the beach width alters the perception
of the risk of the penguins. For instance, a low tide produces a larger
beach, with higher potential risk of predators, thus making
the penguins gather in groups of larger size, see~\cite{CH6}.

\section{Multi-dimensional models}\label{CH:7:CA823oa-M}

It is interesting to remark that the model in~\eqref{CH:7:EQ}
can be easily generalized to the multi-dimensional case.
That is, for any~$i\in\{1,\dots,n(t)\}$
the $i$th coordinate~$p_i$ can be taken to have image in some~$\R^{d}$.
More generally, the dimension of the target space can also vary
in time, by allowing 
for any~$i\in\{1,\dots,n(t)\}$
the $i$th coordinate~$p_i$ to range in some~$\R^{d_i(t)}$,
with~$d_i(t)$ piecewise constant, namely~$d_i(t)=d_{i,j}\in\N_0$ for
any~$t\in(t_{j-1},t_j)$ (compare with~\eqref{CH:7:ndit}).

This modification just causes a small notational complication
in~\eqref{CH:7:pitt}, since each~$p_i(t)$ would now be a vector
in~$\R^{d_i(t)}$ and the array~$p(t)$ would now be
of dimension~$d_1(t)+\dots+d_{n(t)}$.
While we do not indulge here in this generalization,
we observe that such mathematical
extension can be useful, in practice, to consider
the specific location of the
burrows and describe for instance
the movements of the penguins
on the beach (say, a two-dimensional surface) which, as time flows,
gather together in a single queue and move in the end
on a one-dimensional line.

Of course, the rigorous results in Section~\ref{CH:7:HOME}
need to be structurally modified in higher dimension, since several
notions of ``proximity'' of groups, ``direction of march''
and ``orientation of the eye-sight'' can be considered.

\section{Conclusions}\label{CH:7:CA823oa}

As a result of our direct observation at Phillip Island
and a series of scientific discussions with penguin ecologists,
we provide a simple, 
but rigorous, mathematical
model which aims to describe the 
formation of groups of penguins
on the shore at sunset and the return to
their burrows. 

The model is proved to possess existence
and uniqueness of solutions and
quantitative results on the homecoming
of the penguins are given.

The framework is general enough to show the formation
of groups of penguins marching together -- as well as
the possibility that some penguins remain isolated from the rest
of the herd.

The model is also numerically implemented
in simple and explicit simulations.

We believe that the method proposed can
be suitably compared with the real penguins parade,
thus triggering a specific field work on this rather peculiar
topic. Indeed, at the moment, a precise collection of data
focused on the penguins parade
seems to be still missing in the literature, and we think that
a mathematical formulation provides the necessary setting
for describing specific behaviours in ethology,
such as the formation of groups and the possible
isolation of penguins, in a rigorous and quantitative way.

Given the simple and quantitative mathematical
setting that we introduced here, we also believe that
our formulation can be easily modified and improved to
capture possible additional details of the penguins march provided by
the biological data
which may be collected in future specialized field work.

We hope that this problem will also take
advantage of
statistically sound observations by ecologists,
possibly taking into account the speed of the penguins
in different environments, the formation of groups
of different size, the velocity of each group depending
on its size and the links between group formations motivated
by homecoming
and the social structures of the penguin population.

Due to the lack of available biological theories and precise experimental data, the form of some of the functions considered in this chapter should just be considered as an example. This applies in particular to the strategic velocity function, to the eye-sight function and to the stopping function, and it would be desirable to run experiments to provide a better quantification of these notions. 

It would be also interesting to detect how changes in the environment,
such as modified visibility or presence of predators,
influence the formation of groups, their size and their speed.
In general, we think that it would
be very important to provide precise conditions for clustering and to explore
these conditions systematically.

In addition, it would be interesting to
adapt models of this type to social studies, politics and evolutionary biology, in order to describe and quantify the phenomenon
of ``front runners'' which ``wait for the formation of groups of considerable size'' in order to ``more safely
proceed towards their goal''.


\section{Example of a program list}\label{CH:7:APPA}

\small

\begin{verbatim}
H=4;	% Position of the burrow of the penguins community
S=-2;	% The sea lies in the region (-\infty,0]. For simplicity we assume
      that penguins start near the shore, that is, the initial position
      of each penguin is at least S
eps=0.005;	% Drift velocity of the penguins
vc=0.05;		% Cruising speed of a big enough raft of penguins
N=20;		% Number of penguins
delta=(N+1)*eps;		% This parameter is used to compute the strategic
                                velocity of a penguin.
% These parameters define the time interval
TMAX=(H-S)/(2*eps);
t=0.01;
T=(0:t:TMAX);
TG=T(1:1,1:12000);
P=zeros(N,length(T));
% The following is the array of the initial positions of the N penguins
P(:,1)=[-1.95 -1.5 -1.05 -0.6 -0.55 -0.4 -0.2 0.1 0.2 0.4 0.8 0.85 0.9
                               1 1.1 1.15 1.2 1.65 3 3.4];
s=(H-S)/3;	% The parameter encoding the eye-sight of the penguins
pgot=(H-S)/12;	% The parameter representing the stopping function
M=zeros(1,N);
V=M;
PAN=-1;
W=0;

WAVE=sin(T)*delta;  % The "environment function". In this case only
                    waves are taken into account

for k=1 : length(T)-1
  for i=1 : N
    if P(i,k)<H
      for j=1: N  % This cycle checks if the ith penguin is in panic
        if -pgot<P(i,k)-P(j,k) & P(i,k)-P(j,k)<pgot
          PAN=PAN+1;
          if P(i,k)==P(j,k)
            W=W+1;  % This counts the number of penguins in the same
                    position of the ith penguin, that is the dimension
                    of the raft
          end
        end
      end
      if PAN==0  % The ith penguin is stuck because of panic
        if -3.5<P(i,k) & P(i,k)<0
          V(i)=-WAVE(k);
          P(i,k+1)=P(i,k)+V(i)*t;
        else			
          P(i,k+1)=P(i,k);
        end
      else 
        if W>N/2  % The ith penguin is a member of a big enough raft,
               so it tends to go home, forgetful of the other penguins
          if -3.5<P(i,k) & P(i,k)<0  % The environment can still affect
                                     the movement of the raft
            V(i)=vc-WAVE(k);
          else			
            V(i)=vc;  % If the environment does not affect the movement,
                      the penguin moves at cruise velocity
          end
        else  % The raft is not big enough, so the strategic velocity
        of the ith penguin is influenced by the other penguins in sight
          for j=1 : N			
            if -s<P(i,k)-P(j,k) & P(i,k)-P(j,k)<0	
              M(i)=M(i)+delta/N;  % Each penguin in sight ahead adds a
                              delta/N to the strategic velocity of the
                              ith penguin
            else
              if 0<P(i,k)-P(j,k) & P(i,k)-P(j,k)<s		
                M(i)=M(i)-delta/N;  % Each penguin in sight behind
                              subtracts a delta/N from the strategic
                              velocity of the ith penguin
              end
            end
          end
          if -3.5<P(i,k) & P(i,k)<0
            V(i)=eps+M(i)-WAVE(k);
          else
            V(i)=eps+M(i);
          end
        end
        P(i,k+1)=P(i,k)+V(i)*t;
      end
    else
      P(i,k+1)=H;
    end
    PAN=-1;
    W=0;	
  end
  M=zeros(1,N);
  for i=2 : N
    for j=1 : i-1
      if -0.011<P(j,k+1)-P(i,k+1) & P(j,k+1)-P(i,k+1)<0.011
        P(j,k+1)=P(i,k+1);  % For simplicity, we assume that penguins
                  close enough occupy the same position, forming a raft
                  and moving together
      end
    end
  end
end
Q=P(1:N,1:length(TG));
plot(TG,Q)
\end{verbatim}
\normalsize


\section{Figures}\label{CH:7:APPB}

\begin{figure}[!ht]
\centering
\begin{tabular}{c}
\includegraphics[width=8.2cm]{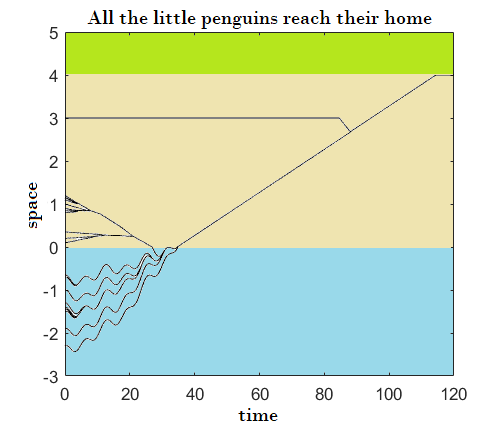}\\ \\ \\ 

\includegraphics[width=8.6cm]{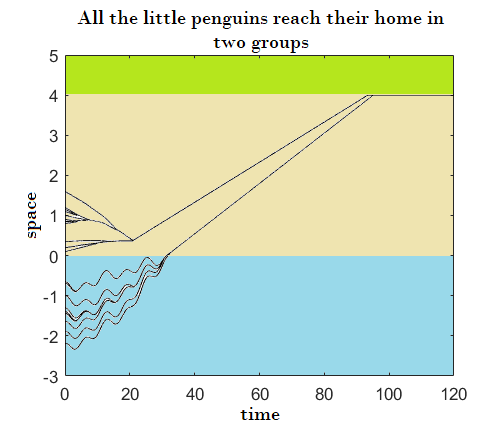}\end{tabular}
\caption{{{All the little penguins safely return home.}}}
\label{CH:7:MP5}
\end{figure}

\newpage

\begin{figure}[!hb]
    \centering
    \includegraphics[width=8.7cm]{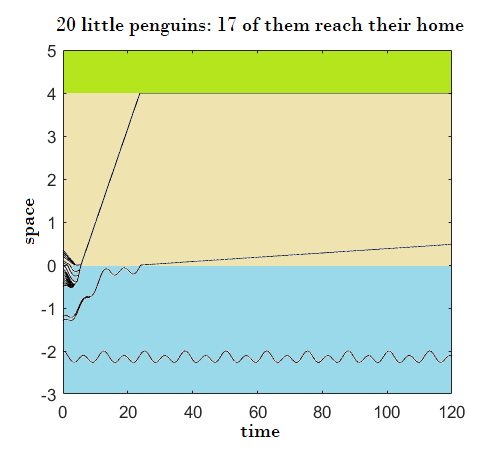}
    \caption{{{One penguin remains in the water.}}}
    \label{CH:7:MK12}
\end{figure}

\begin{figure}[!hb]
    \centering
    \includegraphics[width=8.6cm]{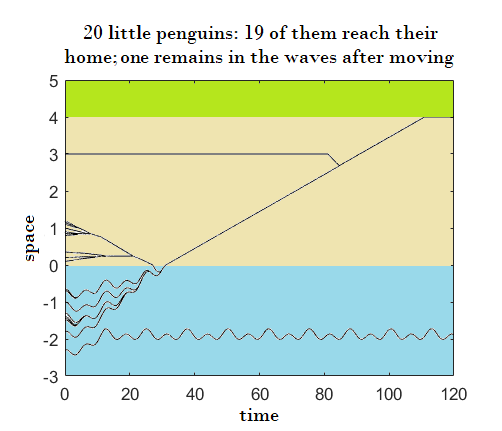}
    \caption{{{One penguin moves towards the others but
    remains in the water.}}}
    \label{CH:7:MK4}
\end{figure}

\newpage

\begin{figure}[!h]
\centering
\begin{tabular}{c}
\includegraphics[width=8.5cm]{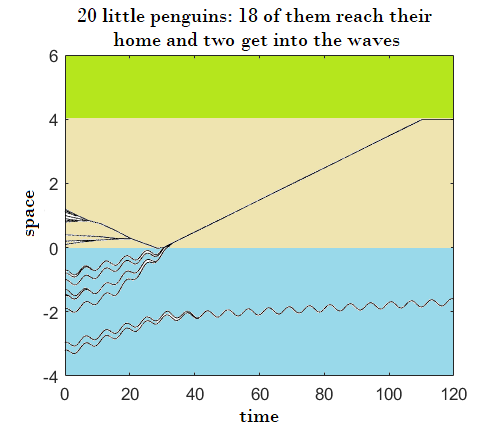} \\ \\ \\ 
\includegraphics[width=8.5cm]{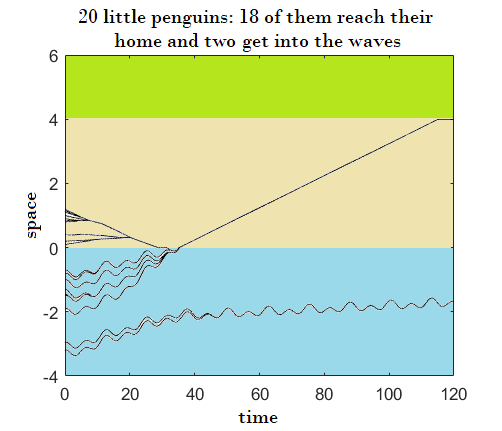}\end{tabular}
\caption{{{Two penguins are still in the water after a long time.}}}
\label{CH:7:MK1}
\end{figure}

\newpage

\begin{figure}[!hb]
    \centering
    \begin{tabular}{c}
    \includegraphics[width=8.5cm]{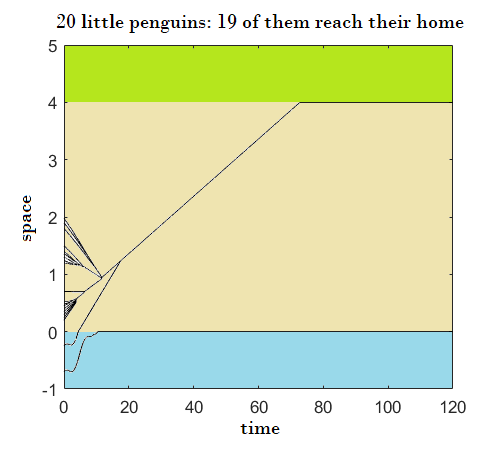}\\ \\ \\ 
    \includegraphics[width=8.5cm]{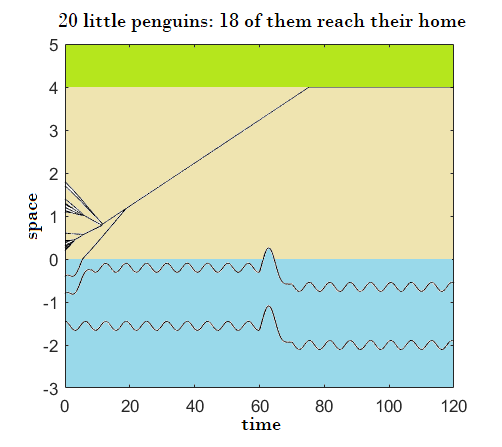}\end{tabular}
    \caption{{{Effect of the waves on the movement of the
    penguins in the sea.}}}
    \label{CH:7:ON-1}
\end{figure}

\newpage

\begin{figure}[!hb]
    \centering
    \includegraphics[width=8.6cm]{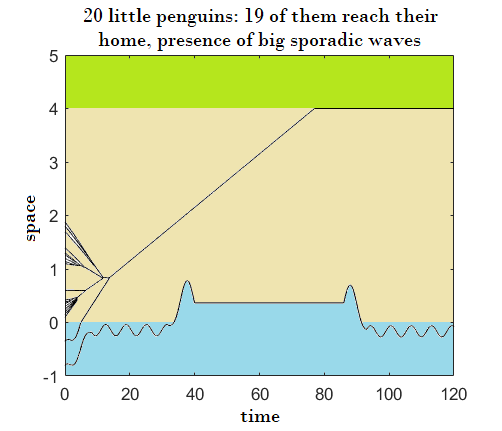}
    \caption{{{Effect of the waves on the movement of
    the penguins in the sea.}}}
    \label{CH:7:ON-3}
\end{figure}

\begin{figure}[!hb]
    \centering
    \includegraphics[width=8.6cm]{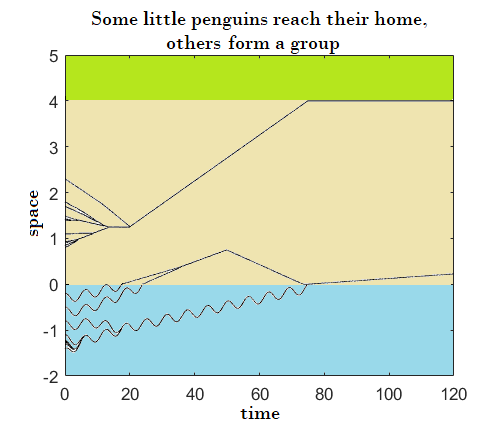}
    \caption{{{The penguins form smaller groups and move towards their home.}}}
    \label{CH:7:Newpic2}
\end{figure}

\newpage

\begin{figure}[!h]
\centering
\begin{tabular}{c}
\includegraphics[width=8.1cm]{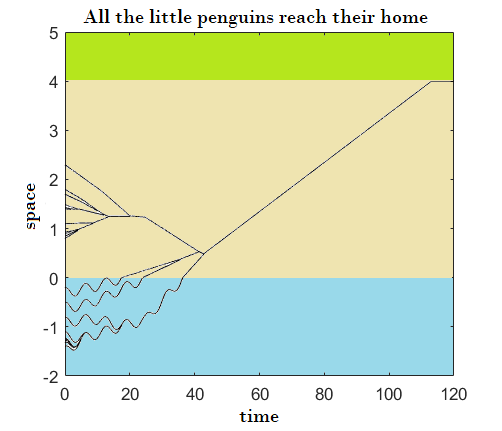} \\ \\ \\ 
\includegraphics[width=8.6cm]{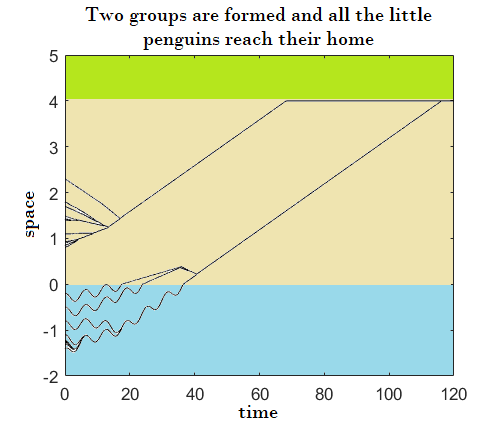}\end{tabular}
\caption{{{The penguins form groups of different sizes and reach their home.}}}
\label{CH:7:Newpic1}
\end{figure}

\end{chapter}

\backmatter

\begin{appendix}

\begin{chapter}{Measure theoretic boundary}\label{CH:1:Appendix_meas_th_bdary}

Since
\begin{equation}\label{CH:1:fin_spazz_basta1}
|E\Delta F|=0\quad\Longrightarrow\quad \Per(E,\Omega)=\Per(F,\Omega)\quad\textrm{and}\quad \Per_s(E,\Omega)=\Per_s(F,\Omega),
\end{equation}
we can modify a set making its topological boundary as big as we want, without changing its (fractional) perimeter.\\
For example, let $E\subseteq\R^n$ be a bounded open set with Lipschitz boundary. Then, if we set
\[
F:=(E\setminus\mathbb Q^n)\cup(\mathbb Q^n\setminus E),
\]
we have $|E\Delta F|=0$ and hence we get \eqref{CH:1:fin_spazz_basta1}.
However $\partial F=\R^n$.

For this reason one considers measure theoretic notions of interior, exterior and boundary, which solely depend on the class
of $\chi_E$ in $L^1_{loc}(\R^n)$.\\
In some sense, by considering the measure theoretic boundary $\partial^-E$ defined below
we can also minimize the size of the topological boundary (see \eqref{CH:1:ess_bdry_intersect}). Moreover, this measure theoretic boundary is actually the
topological boundary of a set which is equivalent to $E$. Thus we obtain a ``good'' representative for the class of $E$.

We refer to \cite[Section 3.2]{Visintin} (see also step two in the proof of~\cite[Proposition~12.19]{Maggi} and \cite[Proposition 3.1]{Giusti}). For some details about the good representative of an $s$-minimal set, see the Appendix of \cite{graph}.
\begin{defn}
Let $E\subseteq\R^n$. For every $t\in[0,1]$ we define the set
\begin{equation}\label{CH:1:density_t}
E^{(t)}:=\left\{x\in\R^n\,\big|\,\exists\,\lim_{r\to0}\frac{|E\cap B_r(x)|}{\omega_nr^n}=t\right\},
\end{equation}
of points density $t$ of $E$.
We also define the essential boundary of $E$ as
\begin{equation*}
\partial_eE:=\R^n\setminus\big(E^{(0)}\cup E^{(1)}\big).
\end{equation*}
\end{defn}

Using the Lebesgue's points Theorem for the characteristic function $\chi_E$, we see that the limit in \eqref{CH:1:density_t} exists for a.e. $x\in\R^n$ and
\begin{equation*}
\lim_{r\to0}\frac{|E\cap B_r(x)|}{\omega_nr^n}=\left\{\begin{array}{cc}1&\textrm{for a.e. }x\in E,\\
0&\textrm{for a.e. }x\in\Co E.
\end{array}
\right.
\end{equation*}
So
\begin{equation*}
|E\Delta E^{(1)}|=0,\qquad|\Co E\Delta E^{(0)}|=0\qquad\textrm{and }|\partial_eE|=0.
\end{equation*}
In particular a set $E$ is equivalent to the set $E^{(1)}$ of its points of density 1.\\
Roughly speaking, the sets $E^{(0)}$ and $E^{(1)}$ can be thought of as a measure theoretic version
of, respectively, the exterior and the interior of the set $E$.
However, notice that both $E^{(1)}$ and $E^{(0)}$ in general are not open.

\smallskip

We have another natural way to define measure theoretic versions of interior, exterior and boundary.
\begin{defn}
Given a set $E\subseteq\R^n$, we define the measure theoretic interior and exterior of $E$ by
\begin{equation*}
E_{int}:=\{x\in\R^n\,|\,\exists\, r>0,\,|E\cap B_r(x)|=\omega_nr^n\}
\end{equation*}
and
\begin{equation*}
E_{ext}:=\{x\in\R^n\,|\,\exists\, r>0,\,|E\cap B_r(x)|=0\},
\end{equation*}
respectively.
Then we define the measure theoretic boundary of $E$ as
\begin{equation*}\begin{split}
\partial^-E&:=\R^n\setminus(E_{ext}\cup E_{int})\\
&
=\{x\in\R^n\,|\,0<|E\cap B_r(x)|<\omega_nr^n\textrm{ for every }r>0\}.
\end{split}
\end{equation*}
\end{defn}
Notice that $E_{ext}$ and $E_{int}$ are open sets and hence $\partial^-E$ is closed. Moreover, since
\begin{equation}\label{CH:1:density_subsets}
E_{ext}\subseteq E^{(0)}\qquad\textrm{and}\qquad E_{int}\subseteq E^{(1)},
\end{equation}
we get
\begin{equation*}
\partial_eE\subseteq\partial^-E.
\end{equation*}
We observe that
\begin{equation}\label{CH:1:ess_bdry_top1}
F\subseteq\R^n\textrm{ s.t. }|E\Delta F|=0\quad\Longrightarrow\quad\partial^-E\subseteq\partial F.
\end{equation}
Indeed, if $|E\Delta F|=0$, then $|F\cap B_r(x)|=|E\cap B_r(x)|$ for every $r>0$. Thus for any $x\in\partial^-E$ we have
\begin{equation*}
0<|F\cap B_r(x)|<\omega_nr^n,
\end{equation*}
which implies
\begin{equation*}
F\cap B_r(x)\not=\emptyset\quad\textrm{and}\quad\Co F\cap B_r(x)\not=\emptyset\quad\textrm{for every }r>0,
\end{equation*}
and hence $x\in\partial F$.

In particular, $\partial^-E\subseteq\partial E$.

Moreover
\begin{equation}\label{CH:1:ess_bdry_top2}
\partial^-E=\partial E^{(1)}.
\end{equation}
Indeed, since $|E\Delta E^{(1)}|=0$, we already know that $\partial^-E\subseteq\partial E^{(1)}$.
The converse inclusion follows from \eqref{CH:1:density_subsets} and the fact that both $E_{ext}$ and $E_{int}$ are open.\\
From \eqref{CH:1:ess_bdry_top1} and \eqref{CH:1:ess_bdry_top2} we obtain
\begin{equation}\label{CH:1:ess_bdry_intersect}
\partial^-E=\bigcap_{F\sim E}\partial F,
\end{equation}
where the intersection is taken over all sets $F\subseteq\R^n$ such that $|E\Delta F|=0$,
so we can think of 
$\partial^-E$ as a way to minimize the size of the topological boundary of $E$.
In particular
\begin{equation*}
F\subseteq\R^n\textrm{ s.t. }|E\Delta F|=0\quad\Longrightarrow\quad\partial^-F=\partial^-E.
\end{equation*}

From \eqref{CH:1:density_subsets} and \eqref{CH:1:ess_bdry_top2} we see that we can take $E^{(1)}$ as ``good'' representative
for $E$, obtaining Remark \ref{CH:1:gmt_assumption}.

\smallskip

Recall that the support of a Radon measure $\mu$ on $\R^n$ is defined as the set
\begin{equation*}
\textrm{supp }\mu:=\{x\in\R^n\,|\,\mu(B_r(x))>0\textrm{ for every }r>0\}.
\end{equation*}
Notice that, being the complementary of the union of all open sets of measure zero, it is a closed set.
In particular, if $E$ is a Caccioppoli set, we have
\begin{equation}\label{CH:1:support_perimeter}
\textrm{supp }|D\chi_E|=\{x\in\R^n\,|\,\Per(E,B_r(x))>0\textrm{ for every }r>0\},
\end{equation}
and
it is easy to verify that
\begin{equation*}
\partial^-E=\textrm{supp }|D\chi_E|=\overline{\partial^*E},
\end{equation*}
where $\partial^*E$ denotes the reduced boundary (see, e.g., \cite[Chapter 15]{Maggi}).
Moreover, $\partial^*E\subseteq\partial_eE$ and by Federer's Theorem (see, e.g., \cite[Theorem 16.2]{Maggi})
we have
\[
\Ha^{n-1}(\partial_eE\setminus\partial^*E)=0.
\]

\begin{figure}[htbp]\label{CH:1:Bdaries_ex}
\begin{center}
\includegraphics[width=125mm]{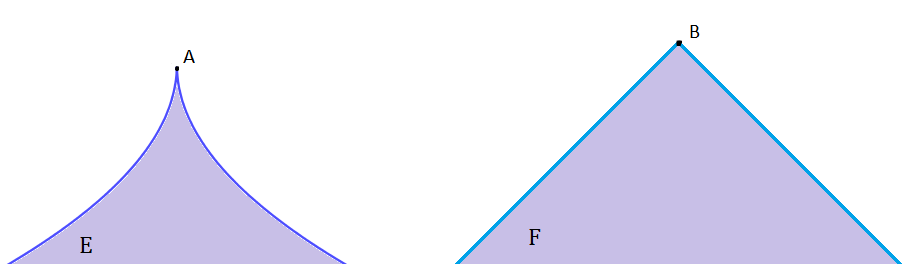}
\caption{{\it The point $A$ belongs to $\partial^-E$ but $A\not\in\partial_e E$.
The point $B$ belongs to $\partial_e F$ but $B\not\in\partial^*F$.}}
\end{center}
\end{figure}

We remark that in general the inclusions
\begin{equation*}
\partial^*E\subseteq\partial_eE\subseteq\partial^-E\subseteq\partial E
\end{equation*}
are all strict.
Indeed, we have already observed in the previous discussion that in general $\partial^-E$
is much smaller than the topological boundary $\partial E$.
In order to have an example of a point $p\in\partial^-E\setminus\partial_eE$ it is enough to consider sublinear cusps. For example,
if $E:=\{(x,y)\in\R^2\,|\,y<-|x|^\frac{1}{2}\}$ and $p:=(0,0)$, then it is easy to verify that $p\in E^{(0)}$ and hence
$p\not\in\partial_e E$. On the other hand, $p\in\partial^-E$.
Finally, the vertex of an angle is an example of a point $p\in\partial_eE\setminus\partial^*E$
(see, e.g., \cite[Example 15.4]{Maggi}).

\end{chapter}

\begin{chapter}{Some geometric observations}\label{CH:3:appendicite}

We collect here some useful results and observations of a geometric nature, concerning in particular the signed distance function.

\section{Signed distance function}\label{CH:1:Appendix_distance_function}

Given $\emptyset\not=E\subseteq\R^n$, the distance function from $E$ is defined as
\begin{equation*}
d_E(x)=d(x,E):=\inf_{y\in E}|x-y|,\qquad\textrm{for }x\in\R^n.
\end{equation*}
The signed distance function from $\partial E$, negative inside $E$, is then defined as
\begin{equation*}
\bar{d}_E(x)=\bar{d}(x,E):=d(x,E)-d(x,\Co E).
\end{equation*}
For the details about the main properties of the signed distance function we refer, e.g., to \cite{GilTru,Ambrosio} and \cite{Bellettini}.

We also define the sets
\begin{equation*}
E_r:=\{x\in\R^n\,|\,\bar{d}_E(x)<r\},
\end{equation*}
for every $r\in\R$, and
\begin{equation*}
N_\varrho(\partial E):=\{|\bar{d}_E|<\varrho\}=\{x\in\R^n\,|\,d(x,\partial E)<\varrho\},
\end{equation*}
for every $\varrho>0$, which is usually called the tubular $\varrho$-neighborhood of $\partial E$.

Let $\Omega\subseteq\R^n$ be a bounded open set with Lipschitz boundary. By definition we can locally describe $\Omega$ near its boundary as the subgraph of appropriate Lipschitz functions.
To be more precise, we can find a finite open covering $\{C_{\varrho_i}\}_{i=1}^m$ of $\partial\Omega$ made of cylinders,
and Lipschitz functions $\varphi_i:B'_{\varrho_i}\longrightarrow\R$ such that $\Omega\cap C_{\varrho_i}$ is the subgraph of
$\varphi_i$.
That is, up to rotations and translations,
\begin{equation*}
C_{\varrho_i}=\{(x',x_n)\in\R^n\,|\,|x'|<\varrho_i,\,|x_n|<\varrho_i\},
\end{equation*}
and
\begin{equation*}\begin{split}
\Omega\cap C_{\varrho_i}&=\{(x',x_n)\in\R^n\,|\,x'\in B'_{\varrho_i},\,-\varrho_i<x_n<\varphi_i(x')\},\\
&
\partial\Omega\cap C_{\varrho_i}=\{(x',\varphi_i(x'))\in\R^n\,|\,x'\in B_{\varrho_i}'\}.
\end{split}
\end{equation*}
Let $L$ be the sup of the Lipschitz constants of the functions $\varphi_i$.

We observe that \cite[Theorem 4.1]{Doktor} guarantees that also the bounded open sets $\Omega_r$ have
Lipschitz boundary, when $r$ is small enough, say $|r|<r_0$.\\
Moreover these sets $\Omega_r$
can locally be described, in the same cylinders $C_{\varrho_i}$ used for $\Omega$, as subgraphs of
Lipschitz functions $\varphi_i^r$ which approximate $\varphi_i$ (see \cite{Doktor} for the precise statement) and whose Lipschitz constants are less than or equal to $L$.\\
Notice that
\begin{equation*}
\partial\Omega_r=\{\bar{d}_\Omega=r\}.
\end{equation*}
Now, since in $C_{\varrho_i}$ the set $\Omega_r$ coincides with the subgraph of $\varphi_i^r$, we have
\begin{equation*}
\Ha^{n-1}(\partial\Omega_r\cap C_{\varrho_i})=\int_{B_{\varrho_i}'}\sqrt{1+|\nabla\varphi_i^r|^2}\,dx'\leq M_i,
\end{equation*}
with $M_i$ depending on $\varrho_i$ and $L$ but not on $r$.\\
Therefore
\begin{equation*}
\Ha^{n-1}(\{\bar{d}_\Omega=r\})\leq\sum_{i=1}^m\Ha^{n-1}(\partial\Omega_r\cap C_{\varrho_i})\leq\sum_{i=1}^mM_i
\end{equation*}
independently on $r$,
proving the following
\begin{prop}\label{CH:1:bound_perimeter_unif}
Let $\Omega\subseteq\R^n$ be a bounded open set with Lipschitz boundary. Then there exists $r_0=r_0(\Omega)>0$ such that
$\Omega_r$ is a bounded open set with Lipschitz boundary for every $r\in(-r_0,r_0)$ and
\begin{equation*}
\sup_{|r|<r_0}\Ha^{n-1}(\{\bar{d}_\Omega=r\})<\infty.
\end{equation*}
\end{prop}

\subsection{Smooth domains}\label{CH:3:A2}
In this section we collect some properties of the signed distance function from the boundary of a regular open set.


We begin by recalling the notion of (uniform) interior ball condition.
\begin{defn}
We say that an open set $\mathcal O$ satisfies an interior ball condition at $x\in\partial\mathcal O$ if
there exists a ball $B_r(y)$ s.t.
\begin{equation*}
B_r(y)\subseteq\mathcal O\qquad\textrm{and}\qquad x\in\partial B_r(y).
\end{equation*}
We say that the condition is ``strict'' if $x$ is the only tangency point, i.e.
\[\partial B_r(y)\cap\partial\mathcal O=\{x\}.\]
The open set $\mathcal O$ satisfies a uniform (strict) interior ball condition of radius $r$ if
it satisfies the (strict) interior ball condition at every point of $\partial\mathcal O$,
with an interior tangent ball of radius at least $r$.\\
In a similar way one defines exterior ball conditions.
\end{defn}
We remark that
if $\mathcal O$ satisfies an interior ball condition of radius $r$ at $x\in\partial\mathcal O$,
then the condition is strict for every radius $r'<r$.

\begin{remark}\label{CH:3:ext_unif_omega}
Let $\Omega\subseteq\R^n$ be a bounded open set with $C^2$ boundary. It is well
known that $\Omega$ satisfies a uniform interior and exterior ball condition. We fix $r_0=r_0(\Omega)>0$
such that $\Omega$ satisfies a strict interior and a strict exterior ball contition of radius $2r_0$
at every point $x\in\partial\Omega$.
Then
\begin{equation}\label{CH:3:r01}
\bar{d}_\Omega\in C^2(N_{2r_0}(\partial\Omega)),
\end{equation}
(see, e.g., \cite[Lemma 14.16]{GilTru}).
\end{remark}

We remark that the distance function $d(\,\cdot\,,E)$ is differentiable at $x\in\R^n\setminus\overline E$ if
and only if there is a unique point $y\in\partial E$ of minimum distance, i.e.
\[d(x,E)=|x-y|.\]
In this case, the two points $x$ and $y$ are related by the formula
\[y=x-d(x,E)\nabla d(x,E).\]

This generalizes to the signed distance function. In particular, if $\Omega$ is bounded and has $C^2$ boundary, then we can
define a $C^1$ projection function from the tubular $2r_0$-neighborhood $N_{2r_0}(\partial\Omega)$ onto $\partial\Omega$ by
assigning to a point $x$ its unique nearest point $\pi(x)$, that is
\[\pi:N_{2r_0}(\partial\Omega)\longrightarrow\partial\Omega,\qquad\pi(x):=x-\bar{d}_\Omega(x)\nabla\bar{d}_\Omega(x).\]
We also remark that
on $\partial\Omega$ we have that $\nabla\bar{d}_\Omega=\nu_\Omega$ and that
\[\nabla\bar{d}_\Omega(x)=\nabla\bar{d}_\Omega(\pi(x))=\nu_\Omega(\pi(x)),\qquad\forall x\in N_{2r_0}(\partial\Omega).\]
Thus $\nabla\bar{d}_\Omega$ is a vector field which extends
the outer unit normal to a tubular neighborhood of $\partial\Omega$, in a $C^1$ way.

Notice that given a point $y\in\partial\Omega$, for every $|\delta|<2r_0$ the point $x:=y+\delta\nu_\Omega(y)$ is such that $\bar{d}_\Omega(x)=\delta$ (and $y$ is its unique nearest point).
Indeed, we consider for example $\delta\in(0,2r_0)$. Then we can find an exterior tangent ball
\[B_{2r_0}(z)\subseteq\Co\Omega,\qquad\partial B_{2r_0}(z)\cap\partial\Omega=\{y\}.\]
Notice that the center of the ball must be
\[z=y+2r_0\nu_\Omega(y).\]
Then, for every $\delta\in(0,2r_0)$ we have
\[B_\delta(y+\delta\nu_\Omega(y))\subseteq B_{2r_0}(y+2r_0\nu_\Omega(y))\subseteq
\Co\Omega,\qquad\partial B_\delta(y+\delta\nu_\Omega(y))\cap\partial\Omega=\{y\}.\]
This proves that
\[|\bar{d}_\Omega(y+\delta\nu_\Omega(y))|=d(x,\partial\Omega)=\delta.\]
Finally, since the point $x$ lies outside $\Omega$, its signed distance function is positive.

\begin{remark}\label{CH:3:c21}
Since $|\nabla\bar{d}_\Omega|=1$, the bounded open sets
\begin{equation*}
\Omega_\delta:=\{\bar{d}_\Omega<\delta\}
\end{equation*}
have $C^2$ boundary
\begin{equation*}
\partial\Omega_\delta=\{\bar{d}_\Omega=\delta\},
\end{equation*}
for every $\delta\in(-2r_0,2r_0)$.
\end{remark}

As a consequence, we know that for every $|\delta|<2r_0$ the set $\Omega_\delta$
satisfies a uniform interior and exterior ball condition of radius $r(\delta)>0$.
Moreover, we have that $r(\delta)\geq r_0$ for every $|\delta|\leq r_0$
(see also \cite[Appendix A]{MR3436398}
for related results).

\begin{lemma}\label{CH:3:geomlem}
Let $\Omega\subseteq\R^n$ be a bounded open set with $C^2$ boundary.
Then for every $\delta\in[-r_0,r_0]$ the set $\Omega_\delta$ 
satisfies a uniform interior and exterior ball condition of radius at least $r_0$, i.e.
\begin{equation*}
r(\delta)\geq r_0\qquad\textrm{for every }|\delta|\leq r_0.
\end{equation*}
\end{lemma}
\begin{proof}
Take for example $\delta\in[-r_0,0)$ and let $x\in\partial\Omega_\delta=\{\bar{d}_\Omega=\delta\}$.
We show that $\Omega_\delta$ has an interior tangent ball of radius $r_0$ at $x$. The other cases are proven in a similar way.

Consider the projection $\pi(x)\in\partial\Omega$ and the point
\[x_0:=x-r_0\nabla\bar{d}_\Omega(x)=\pi(x)-(r_0+|\delta|)\nu_\Omega(\pi(x)).\]
Then
\[B_{r_0}(x_0)\subseteq\Omega_\delta\quad\textrm{ and }\quad x\in\partial B_{r_0}(x_0)\cap\partial\Omega_\delta.\]
Indeed, notice that, as remarked above,
\[d(x_0,\partial\Omega)=|x_0-\pi(x)|=r_0+|\delta|.\]
Thus, by the triangle inequality we have that
\[d(z,\partial\Omega)\ge d(x_0,\partial\Omega)-|z-x_0|>|\delta|,\qquad\textrm{ for every }z\in B_{r_0}(x_0),\]
so $B_{r_0}\subseteq\Omega_\delta$. Moreover, by definition of $x_0$ we have
\[x\in\partial B_{r_0}(x_0)\cap\partial\Omega_\delta\]
and the desired result follows.
\end{proof}

To conclude, we remark that the sets $\overline{\Omega_{-\delta}}$ are retracts of $\Omega$, for every $\delta\in(0,r_0]$.
Indeed, roughly speaking, each set $\overline{\Omega_{-\delta}}$ is obtained by deforming $\Omega$ in normal direction,
towards the interior.
An important consequence is that if $\Omega$ is connected then $\overline{\Omega_{-\delta}}$ is path connected.

To be more precise, we have the following:

\begin{prop}\label{CH:3:retract}
Let $\Omega\subseteq\Rn$ be a bounded open set with $C^2$ boundary.
Let $\delta\in(0,r_0]$ and define
\[\mathcal D:\Omega\longrightarrow\overline{\Omega_{-\delta}},\qquad\mathcal D(x):=
\left\{\begin{split}&x,& &x\in\Omega_{-\delta},\\
&x-\big(\delta+\bar{d}_\Omega(x)\big)\nabla\bar{d}_\Omega(x),& & x\in\Omega\setminus\Omega_{-\delta}.\end{split}\right.\]
Then $\mathcal D$ is a retraction of $\Omega$ onto $\overline{\Omega_{-\delta}}$, i.e. it is continuous and
$\mathcal D(x)=x$ for every $x\in\overline{\Omega_{-\delta}}$.
In particular, if $\Omega$ is connected, then $\overline{\Omega_{-\delta}}$ is path connected.
\end{prop}

\begin{proof}
Notice that the function
\[\Phi(x):=x-\big(\delta+\bar{d}_\Omega(x)\big)\nabla\bar{d}_\Omega(x)\]
is continuous in $\Omega\setminus\Omega_{-\delta}$ and $\Phi(x)=x$ for every $x\in\partial\Omega_{-\delta}$.
Therefore the function $\mathcal D$ is continuous.

We are left to show that
\[\mathcal D(\Omega\setminus\Omega_{-\delta})\subseteq\partial\Omega_{-\delta}.\]
For this, it is enough to notice that
\[\mathcal D(x)=\pi(x)-\delta\nu_\Omega(\pi(x))\qquad\textrm{for every }x\in\Omega\setminus\Omega_{-\delta}.\]
To conclude, suppose that $\Omega$ is connected and recall that if an open set $\Omega\subseteq\R^n$ is connected, then it is also path connected.
Thus $\overline{\Omega_{-\delta}}$,
being the continuous image of a path connected space, is itself
path connected.
\end{proof}

\section{Sliding the balls}

We now point out the following useful geometric result, which has been exploited in Chapter \ref{Asympto0_CH_label}.

\begin{lemma}\label{CH:3:slidetheballs}
Let $F\subseteq\Rn$ be such that\footnote{Concerning the statement of Lemma~\ref{CH:3:slidetheballs},
we recall that the notation~$\overline{F}$ denotes the closure of the set~$F$,
when~$F$ is modified, up to sets of measure zero,
in such a way that~$F$ is assumed to contain its measure theoretic interior~$F_{int}$
and to have empty intersection with the exterior~$F_{ext}$, 
according to the setting described in Remark~\ref{CH:1:gmt_assumption}.
For instance, if~$F$ is a segment in~$\R^2$, this convention implies that~$F_{int}=\emptyset$, $F_{ext}=\R^2$ and so~$F$ and~$\overline{F}$ in this case also reduce to the empty set.}
\[B_\delta(p)\subseteq F_{ext}\quad\textrm{for some }\delta>0\qquad\textrm{and}\qquad q\in\overline{F},\]
and let $c:[0,1]\longrightarrow\Rn$ be a continuous curve connecting $p$ to $q$, that is
\[c(0)=p\qquad\textrm{and}\qquad c(1)=q.\]
Then there exists $t_0\in[0,1)$ such that $B_\delta\big(c(t_0)\big)$ is an exterior tangent ball to $F$,
that is
\eqlab{\label{CH:3:slide_ext_tg}
B_\delta\big(c(t_0)\big)\subseteq F_{ext}\qquad\textrm{and}\qquad\partial B_\delta\big(c(t_0)\big)\cap\partial F\not=\emptyset.}
\end{lemma}

\begin{proof}
Define
\eqlab{\label{CH:3:slide_t}
t_0:=\sup\Big\{\tau\in[0,1]\,\big|\,\bigcup_{t\in[0,\tau]}B_\delta\big(c(t)\big)\subseteq F_{ext}\Big\}.
}
We begin by proving that
\eqlab{\label{CH:3:jardine}
B_\delta\big(c(t_0)\big)\subseteq F_{ext}.
}
If $t_0=0$, this is trivially true by hypothesis.
Thus, suppose that $t_0>0$ and assume by contradiction that
\bgs{
B_\delta\big(c(t_0)\big)\cap \overline{F}\not=\emptyset.
}
Then there exists a point
\[
y\in\overline{F}=F_{int}\cup\partial F\quad\mbox{s.t.}\quad d:=|y-c(t_0)|<\delta.
\]
By exploiting the continuity of $c$, we can find $t\in[0,t_0)$
such that
\[
|y-c(t)|\le|y-c(t_0)|+|c(t_0)-c(t)|\le d+\frac{\delta-d}{2}<\delta,
\]
and hence $y\in B_\delta\big(c(t)\big)$.
However, this is in contradiction with the fact that, by definition of $t_0$, we have $B_\delta\big(c(t)\big)\subseteq F_{ext}$.
This concludes the proof of \eqref{CH:3:jardine}.

We point out that, since $q\in\overline F$, by \eqref{CH:3:jardine} we have that $t_0<1$.

Now we prove that $t_0$ as defined in \eqref{CH:3:slide_t} satisfies \eqref{CH:3:slide_ext_tg}.

Notice that by \eqref{CH:3:jardine} we have
\eqlab{\label{CH:3:slide_pf}
\overline{B_\delta\big(c(t_0)\big)}\subseteq \overline{F_{ext}}=F_{ext}\cup\partial F.
}
Suppose that
\[\partial B_\delta\big(c(t_0)\big)\cap\partial F=\emptyset.\]
Then \eqref{CH:3:slide_pf} implies that
\[\overline{B_\delta\big(c(t_0)\big)}\subseteq F_{ext},\]
and, since $F_{ext}$ is an open set, we can find $\tilde\delta>\delta$ such that
\[B_{\tilde\delta}\big(c(t_0)\big)\subseteq F_{ext}.\]
By continuity of $c$ we can find $\eps\in(0,1-t_0)$ small enough such that
\[|c(t)-c(t_0)|<\tilde\delta-\delta,\qquad\forall\,t\in[t_0,t_0+\eps].\]
Therefore
\[B_\delta\big(c(t)\big)\subseteq B_{\tilde\delta}\big(c(t_0)\big)\subseteq F_{ext},\qquad\forall\,t\in[t_0,t_0+\eps],\]
and hence
\[\bigcup_{t\in[0,t_0+\eps]}B_\delta\big(c(t)\big)\subseteq F_{ext},\]
which is in contradiction with
the definition of $t_0$. Thus
\[\partial B_\delta\big(c(t_0)\big)\cap\partial F\not=\emptyset,\]
which concludes the proof.
\end{proof}

\end{chapter}

\begin{chapter}{Collection of useful results on nonlocal minimal surfaces}\label{CH:3:appendicite2}

Here, we collect some auxiliary results on nonlocal minimal surfaces.
In particular, we recall the representation of
the fractional mean curvature when the set is a graph and
a useful and general version of the maximum principle.

\section{Explicit formulas for the fractional mean curvature of a graph}
We denote 
\[Q_{r,h}(x):=B'_r(x')\times(x_n-h,x_n+h),\] for $x\in\R^n,$ $r,h>0$. If $x=0$, we write $Q_{r,h}:=Q_{r,h}(0)$. Let also 
\[g_s(t):=\frac{1}{(1+t^2)^\frac{n+s}{2}}\qquad\textrm{and}\qquad G_s(t):=\int_0^tg_s(\tau)\,d\tau.\]
Notice that
\[0<g_s(t)\leq1,\quad\forall\,t\in\R\qquad\textrm{and}\qquad\int_{-\infty}^{+\infty}g_s(t)\,dt<\infty,\]
for every $s\in(0,1)$.

In this notation, we can write the fractional mean curvature of a supergraph as follows:

\begin{prop}
Let $F\subseteq\R^n$ and $p\in\partial F$ such that
\[F\cap Q _{r,h}(p)=\{(x',x_n)\in\R^n\,|\,x'\in B'_r(p'),\,v(x')<x_n<p_n+h\},\]
for some $v\in C^{1,\alpha}(B'_r(p'))$. Then for every $s\in(0,\alpha)$
\eqlab{\label{CH:3:complete_curv_formula}\I_s[F](p)&
=2\int_{B'_r(p')}\Big\{G_s\Big(\frac{v(y')-v(p')}{|y'-p'|}\Big)
-G_s\Big(\nabla v(p')\cdot\frac{y'-p'}{|y'-p'|}\Big)\Big\}\frac{dy'}{|y'-p'|^{n-1+s}}\\
&
\qquad\qquad+\int_{\R^n\setminus Q_{r,h}(p)}\frac{\chi_{\Co F}(y)-\chi_F(y)}{|y-p|^{n+s}}\,dy.}
\end{prop}

This explicit formula was introduced in \cite{regularity} (see also \cite{Abaty}) when $\nabla v(p)=0$. In \cite{bootstrap}, the reader can find the formula for the case of non-zero gradient. 

\begin{remark}
In the right hand side of \eqref{CH:3:complete_curv_formula}
there is no need to consider the principal value, since the integrals are summable.
Indeed,
\bgs{
\Big|G_s\Big(&\frac{v(y')-v(p')}{|y'-p'|}\Big)
-G_s\Big(\nabla v(p')\cdot\frac{y'-p'}{|y'-p'|}\Big)\Big|
=\Big|\int_{\nabla v(p')\cdot\frac{y'-p'}{|y'-p'|}}^{\frac{v(y')-v(p')}{|y'-p'|}}g_s(t)\,dt\Big|\\
&
\leq\Big|\frac{v(y')-v(p')-\nabla v(p')\cdot(y'-p')}{|y'-p'|}\Big|\leq \|v\|_{C^{1,\alpha}(B'_r(p'))}|y'-p'|^\alpha,
}
for every $y'\in B'_r(p')$.
As for the last inequality, notice that by the Mean value Theorem we have
\[v(y')-v(p')=\nabla v(\xi)\cdot(y'-p'),\]
for some $\xi\in B'_r(p')$ on the segment with end points $y'$ and $p'$. Thus
\bgs{|v(y')-v(p')&-\nabla v(p')\cdot(y'-p')|=|(\nabla v(\xi)-\nabla v(p'))\cdot(y'-p')|\\
&
\leq|\nabla v(\xi)-\nabla v(p')||y'-p'|\leq\|\nabla v\|_{C^{0,\alpha}(B'_r(p'))}|\xi-p'|^\alpha|y'-p'|\\
&
\leq\|v\|_{C^{1,\alpha}(B'_r(p'))}|y'-p'|^{1+\alpha}.
}
We denote for simplicity
\eqlab{ \label{CH:3:mathcalg} \mathcal G(s,v,y',p'):= G_s\Big(&\frac{v(y')-v(p')}{|y'-p'|}\Big)
-G_s\Big(\nabla v(p')\cdot\frac{y'-p'}{|y'-p'|}\Big).}
With this notation, we have
\eqlab{\label{CH:3:Holder_useful} |\mathcal G(s,v,y',p')|  \leq \|v\|_{C^{1,\alpha}(B'_r(p'))}|y'-p'|^\alpha.}
\end{remark}


\section[Interior regularity theory]{Interior regularity theory and
its influence on the Euler-Lagrange equation inside the domain}\label{CH:3:brr2}

In this Appendix we give a short review of the the Euler-Lagrange equation 
in the interior of the domain. In particular, by exploiting results which give an improvement of
the regularity of $\partial E$, we show that an $s$-minimal set is a classical solution
of the Euler-Lagrange equation almost everywhere.

First of all, we recall the definition of supersolution.
\begin{defn}
Let $\Omega\subseteq\R^n$ be an open set and let $s\in(0,1)$. A set $E$
is an $s$-supersolution in $\Omega$ if $\Per_s(E,\Omega)<\infty$ and
\begin{equation}\label{CH:3:supersolution}
\Per_s(E,\Omega)\leq \Per_s(F,\Omega)\quad\textrm{for every set }E\textrm{ s.t. }E\subseteq F\textrm{ and }F\setminus\Omega=E\setminus\Omega.
\end{equation}
\end{defn}
We remark that \eqref{CH:3:supersolution} is equivalent to
\begin{equation*}
A\subseteq\Co E\cap\Omega\qquad\Longrightarrow\qquad \Ll_s(A,E)-\Ll_s(A,\Co(E\cup A))\leq0.
\end{equation*}
In a similar way one defines $s$-subsolutions.

In \cite{CRS10} it is shown that a set $E$ which is an $s$-supersolution in $\Omega$
is also a viscosity supersolution of the equation $\I_s[E]=0$ on $\partial E\cap\Omega$.
To be more precise

\begin{theorem}[Theorem 5.1 of \cite{CRS10}]\label{CH:3:viscsol}
Let $E$ be an $s$-supersolution in the open set $\Omega$. If $x_0\in\partial E\cap \Omega$ and $E$ has an interior tangent ball
at $x_0$, contained in $\Omega$,
i.e.
\begin{equation*}
B_r(y)\subseteq E\cap\Omega\quad\textrm{s.t.}\quad x_0\in\partial E\cap\partial B_r(y),
\end{equation*}
then
\begin{equation}\label{CH:3:supersolution_ineq}
\liminf_{\varrho\to0^+}\I_s^\varrho[E](x_0)\geq0.
\end{equation}

\end{theorem}

In particular, $E$ is a viscosity supersolution in the following sense.

\begin{corollary}
Let $E$ be an $s$-supersolution in the open set $\Omega$ and let $F$ be an open set such that $F\subseteq E$.
If $x\in(\partial E\cap\partial F)\cap\Omega$ and $\partial F$ is $C^{1,1}$ near $x$,
then $\I_s[F](x)\geq0$.
\end{corollary}
\begin{proof}
Since $\partial F$ is $C^{1,1}$ near $x$, $F$ has an interior tangent ball at $x$. In particular, notice that this ball
is tangent also to $E$ at $x$ (from the inside).
Thus by Theorem \ref{CH:3:viscsol}
\begin{equation*}
\liminf_{\varrho\to0^+}\I_s^\varrho[E](x)\geq0.
\end{equation*}
Now notice that
\begin{equation*}
F\subseteq E\qquad\Longrightarrow\qquad\chi_{\Co F}-\chi_F\geq\chi_{\Co E}-\chi_E,
\end{equation*}
so
\begin{equation*}
\I_s^\delta[F](x)\geq\I_s^\delta[E](x)\qquad\forall\,\delta>0.
\end{equation*}
Since $\I_s[F](x)$ is well defined, it is then enough to pass to the limit $\delta\to0$.
\end{proof}

\begin{remark}
Similarly, for an $s$-subsolution $E$ which has an exterior tangent ball at $x_0$ we obtain
\begin{equation}\label{CH:3:subsolution_ineq}
\limsup_{\varrho\to0^+}\I_s^\varrho[E](x_0)\leq0.
\end{equation}
\end{remark}

Now we recall the following two regularity results.
If $E$ is $s$-minimal, having a tangent ball (either interior or exterior) at some point $x_0\in\partial E\cap\Omega$ is enough (via an improvement of flatness result) to have $C^{1,\alpha}$ regularity in a neighborhood of $x_0$ (see \cite[Corollary 6.2]{CRS10}).
%
Moreover, bootstrapping arguments prove that $C^{0,1}$ regularity
guarantees $C^\infty$ regularity (according to \cite[Theorem 1.1]{FV17}).

%

It is also convenient to recall the notion of locally $s$-minimal set, which is useful when considering an unbounded domain $\Omega$.\\
We say that a set $E\subseteq\R^n$ is locally $s$-minimal in an open set $\Omega\subseteq\R^n$ if $E$ is $s$-minimal in every bounded
open set $\Omega'\Subset\R^n$.

Exploiting the regularity results that we recalled above, we obtain the following:
\begin{theorem}\label{CH:3:EL_inside}
Let $\Omega\subseteq\R^n$ be an open set and let $E$ be locally $s$-minimal in $\Omega$. If $x_0\in\partial E\cap\Omega$
and $E$ has either an interior or exterior tangent ball at $x_0$,
then there exists $r>0$ such that $\partial E\cap B_r(x_0)$
is $C^\infty$ and
\begin{equation}\label{CH:3:EL_eq_inside}
\I_s[E](x)=0\qquad\textrm{for every}\quad x\in\partial E\cap B_r(x_0).
\end{equation}
\end{theorem}
\begin{proof}
Since $x_0\in\partial E\cap \Omega$ and $\Omega$ is open, we can find $r>0$ such that $B_r(x_0)\Subset\Omega$.\\
The set $E$ is then $s$-minimal in $B_r(x_0)$. Moreover, by hypothesis we have a tangent ball (either interior or exterior)
to $E$ at $x_0$. Also notice that we can suppose that the tangent ball is contained in $B_r(x_0)$.\\
Thus, by \cite[Corollary 6.2]{CRS10} and \cite[Theorem 1.1]{FV17}, we know that $\partial E$ is $C^\infty$ in
$B_r(x_0)$ (up to taking another $r>0$ small enough).

In particular, $\I_s[E](x)$ is well defined for every $x\in\partial E\cap B_r(x_0)$ and $E$ has both an interior and an exterior tangent
ball at every $x\in\partial E\cap B_r(x_0)$ (both contained in $B_r(x_0)$).\\
Therefore, since an $s$-minimal set is both an $s$-supersolution and an $s$-subsolution, by
\eqref{CH:3:supersolution_ineq} and \eqref{CH:3:subsolution_ineq}, we obtain
\begin{equation*}
0\leq\liminf_{\varrho\to0^+}\I_s^\varrho[E](x)=\I_s[E](x)=\limsup_{\varrho\to0^+}\I_s^\varrho[E](x)\leq0,
\end{equation*}
for every $x\in\partial E\cap B_r(x_0)$, proving \eqref{CH:3:EL_eq_inside}.
\end{proof}

Furthermore, we recall that if $E\subseteq\R^n$ is $s$-minimal in $\Omega$, then the singular set
$\Sigma(E;\Omega)\subseteq\partial E\cap\Omega$
has Hausdorff dimension at most $n-3$ (by the dimension reduction argument developed in \cite[Section 10]{CRS10}
and \cite[Corollary 2]{SV13}).

Now suppose that $E$ is locally $s$-minimal in an open set $\Omega$.
We observe that
we can find a sequence of bounded open sets with Lipschitz
boundaries $\Omega_k\Subset\Omega$ such that $\bigcup\Omega_k=\Omega$
(see, e.g., Corollary \ref{CH:2:regular_approx_open_sets_coroll}). Since $E$ is $s$-minimal in each $\Omega_k$ and
$\Sigma(E;\Omega)=\bigcup\Sigma(E;\Omega_k)$, we get in particular
\begin{equation}\label{CH:3:singset}
\Ha^{n-2}(\Sigma(E;\Omega))\leq\sum_{k=1}^\infty\Ha^{n-2}(\Sigma(E;\Omega_k))=0
\end{equation}
(and indeed $\Sigma(E;\Omega)$ has Hausdorff dimension at most $n-3$, since we have inequality \eqref{CH:3:singset}
with $n-d$ in place of $n-2$, for every $d\in[0,3)$).

As a consequence, a (locally) $s$-minimal set is a classical solution of the Euler-Lagrange equation,
in the following sense
\begin{theorem}\label{CH:3:classicalsense}
Let $\Omega\subseteq\R^n$ be an open set and let $E$ be locally $s$-minimal in $\Omega$.
Then
\begin{equation*}
\I_s[E](x)=0\qquad\textrm{for every }x\in(\partial E\cap\Omega)\setminus\Sigma(E;\Omega),
\end{equation*}
and hence in particular for $\Ha^{n-1}$-a.e. $x\in\partial E\cap\Omega$.
\end{theorem}


\section[Boundary Euler-Lagrange inequalities]{Boundary Euler-Lagrange inequalities for the fractional perimeter}\label{CH:3:appendicite3}
We recall that a set $E$ is locally $s$-minimal in an open set $\Omega$ if it is $s$-minimal in every bounded open set compactly contained in $\Omega$.
In this section we show that the Euler-Lagrange equation of a locally $s$-minimal set $E$ holds (at least as an inequality)
also at a point $p\in\partial E\cap\partial\Omega$,
provided that the boundary $\partial E$ and the boundary $\partial\Omega$ do not intersect
``transversally'' in $p$.

To be more precise, we prove the following

\begin{theorem}\label{CH:3:EL_boundary_coroll}
Let $s\in(0,1)$. Let $\Omega\subseteq\R^n$ be an open set and let $E\subseteq\R^n$
be locally $s$-minimal in $\Omega$. Suppose that $p\in\partial E\cap\partial\Omega$ is such that
$\partial\Omega$ is $C^{1,1}$ in $B_{R_0}(p)$, for some $R_0>0$. Assume also that
\eqlab{\label{CH:3:obstaclehp}
B_{R_0}(p)\setminus\Omega\subseteq\Co E.
}
Then
\begin{equation*}
\I_s[E](p)\leq0.
\end{equation*}
Moreover, if there exists $R\in(0,R_0)$ such that
\begin{equation}\label{CH:3:detach_hp}
\partial E\cap\big(\Omega\cap B_r(p)\big)\not=\emptyset\qquad\textrm{for every }r\in(0,R),
\end{equation}
then $$\I_s[E](p)=0.$$
\end{theorem}

We remark that
by hypothesis the open set $B_{R_0}(p)\setminus\overline{\Omega}$ is tangent to $E$ at $p$, from the outside.
Therefore, either \eqref{CH:3:detach_hp} holds true, meaning roughly speaking that the boundary of $E$ detaches
from the boundary of $\Omega$ at $p$ (towards the interior of $\Omega$),
or $\partial E$ coincides with $\partial\Omega$ near $p$.

\begin{figure}[htbp]
\begin{center}
\includegraphics[width=120mm]{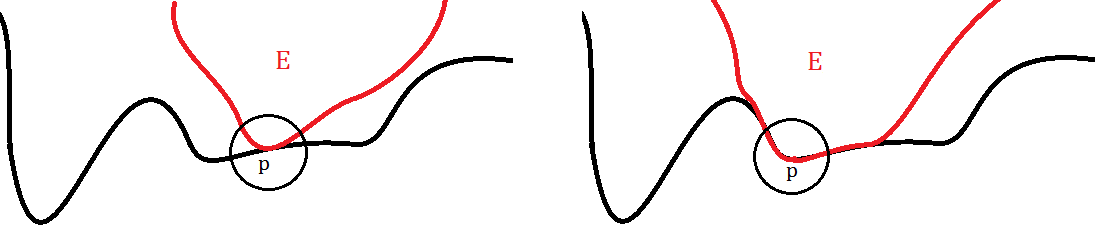}
\caption{\it Examples of a set which satisfies \eqref{CH:3:detach_hp} (on the left)
and of a set whose boundary sticks to that of $\Omega$ near $p$ (on the right)}
\end{center}
\end{figure}

Roughly speaking, the idea of the proof of Theorem \ref{CH:3:EL_boundary_coroll} is the following.
The set $\mathcal O:=B_{R_0}(p)\setminus\overline{\Omega}$
plays the role of an obstacle in the minimization of the $s$-perimeter in $B_{R_0}(p)$.
The (local) minimality of $E$ in $\Omega$, together with hypothesis
\eqref{CH:3:obstaclehp}, implies that $E$ solves this geometric obstacle type problem, which has
been investigated in \cite{GeomObst}. As a consequence, the set $E$ is a viscosity subsolution in $B_{R_0}(p)$ and we obtain that
$\I_s[E](p)\le0$.
Furthermore, the regularity result proved in \cite{GeomObst} guarantees that $\partial E$ is $C^{1,\sigma}$, with $\sigma>s$, near $p$. Thus, if $\partial E$ satisfies \eqref{CH:3:detach_hp}, then we can exploit the Euler-Lagrange equation inside $\Omega$
and the continuity of $\I_s[E]$ to prove that $\I_s[E](p)=0$.

We now proceed to give a rigorous proof of Theorem \ref{CH:3:EL_boundary_coroll}.

\begin{proof}[Proof of Theorem \ref{CH:3:EL_boundary_coroll}]
We begin by observing that we can find a bounded and connected open set $\Omega'\subseteq\Omega$ such that
\begin{equation*}
\partial\Omega'\textrm{ is }C^{1,1}\qquad\textrm{and}\qquad
\Omega'\cap B_\frac{R_0}{2}(p)=\Omega\cap B_\frac{R_0}{2}(p).
\end{equation*}
Then, since $E$ is locally $s$-minimal in $\Omega$, we know that it is locally $s$-minimal also in $\Omega'$.
Hence, since $\Omega'$ is bounded and has regular boundary, by Theorem \ref{CH:2:confront_min_teo} we find that
$E$ is actually $s$-minimal in $\Omega'$.
Moreover $p\in\partial E\cap\partial\Omega'$
and
\begin{equation*}
B_\frac{R_0}{2}(p)\setminus\Omega'=B_\frac{R_0}{2}(p)\setminus\Omega
\subseteq B_{R_0}(p)\setminus\Omega\subseteq\Co E.
\end{equation*}
Therefore, we can suppose without loss of generality that $\Omega$ is a bounded and connected open set with $C^{1,1}$ boundary
$\partial\Omega$ and that $E$ is $s$-minimal in $\Omega$.

\smallskip

As observed in the proof of \cite[Theorem 5.1]{graph}, the minimality of $E$ and hypothesis
\eqref{CH:3:obstaclehp} imply that the set $\Co E$ is a solution, in $B_\frac{R_0}{4}(p)$, of the geometric obstacle type problem considered in \cite{GeomObst}.

More precisely, we remark that we can find a bounded and connected open set $\mathcal O$ with $C^{1,1}$ boundary, such that
\bgs{
\mathcal O\cap B_\frac{R_0}{4}(p)=B_\frac{R_0}{4}(p)\setminus\overline{\Omega}.
}
Then hypothesis \eqref{CH:3:obstaclehp} guarantees that
\bgs{
\mathcal O\cap B_\frac{R_0}{4}(p)\subseteq\Co E.
}
Now, by arguing as in the proof of \cite[Theorem 5.1]{graph},
we find that the minimality of $E$ (hence also of $\Co E$) in $\Omega$
implies that
\bgs{
\Per_s\Big(\Co E,B_\frac{R_0}{4}(p)\Big)
\le \Per_s\Big(F,B_\frac{R_0}{4}(p)\Big),
}
for every $F\subseteq\R^n$ such that
\bgs{
F\setminus B_\frac{R_0}{4}(p)=\Co E\setminus B_\frac{R_0}{4}(p)
\quad\mbox{and}\quad\mathcal O\cap B_\frac{R_0}{4}(p)\subseteq F.
}
In particular, as observed in \cite{GeomObst} (see the comment (2.2) there), the set $\Co E$ is a viscosity supersolution in $B_\frac{R_0}{4}(p)$, meaning that
the set $E$ is a viscosity subsolution in $B_\frac{R_0}{4}(p)$. Now, since the set $\Omega$ has $C^{1,1}$ boundary,
we can find an exterior tangent ball at $p\in\partial\Omega$.
By hypothesis \eqref{CH:3:obstaclehp}, this means that we can find an exterior tangent ball at $p\in\partial E$ and hence
we have
\eqlab{\label{CH:3:obst1}
\limsup_{\varrho\to0^+}\I_s^\varrho[E](p)\le0.
}
Furthermore, \cite[Theorem 1.1]{GeomObst} guarantees that $\partial E$ is $C^{1,\sigma}$ in $B_{R_0'}(p)$
for some $R'_0\in(0,R_0)$, and $\sigma:=\frac{1+s}{2}$ (see also \cite[Theorem 5.1]{graph}).
In particular, since $\sigma>s$, we know that the $s$-fractional mean curvature of $E$ is well defined at $p$.
Therefore \eqref{CH:3:obst1} actually implies that~$\I_s[E](p)\le0$, as claimed.

Now we suppose in addition
that \eqref{CH:3:detach_hp} holds true, i.e. that
\begin{equation*}
\partial E\cap\big(\Omega\cap B_r(p)\big)\not=\emptyset\qquad\textrm{for every }r\in(0,R),
\end{equation*}
with $R<R'_0$.
By \cite[Theorem 1.1]{FV17} we know that $\partial E\cap\big(B_R(p)\cap\Omega\big)$ is $C^\infty$.
In particular, as observed in Theorem \ref{CH:3:EL_inside}, we know that every point $x\in\partial E\cap\big(B_R(p)\cap\Omega\big)$ satisfies the Euler-Lagrange
equation in the classical sense, i.e.
\begin{equation}\label{CH:3:EL_proof_eq1}
\I_s[E](x)=0\qquad\textrm{for every }x\in\partial E\cap\big(B_R(p)\cap\Omega\big).
\end{equation}
Since $\partial E\cap B_R(p)$ is $C^{1,\sigma}$, with $\sigma>s$,
we also know that $\I_s[E]\in C(\partial E\cap B_R(p))$ (by, e.g., Proposition \ref{CH:3:rsdfyish} or \cite[Lemma 3.4]{graph}).
Finally, we observe that by \eqref{CH:3:detach_hp} we can find a sequence of points $x_k\in\partial E\cap\big(B_R(p)\cap\Omega\big)$
such that $x_k\longrightarrow p$.
Then, by the continuity of $\I_s[E]$ and \eqref{CH:3:EL_proof_eq1} we get
\begin{equation*}
\I_s[E](p)=\lim_{k\to\infty}\I_s[E](x_k)=0,
\end{equation*}
concluding the proof.
\end{proof}


\section{A maximum principle}
By exploiting the Euler-Lagrange equation, we can compare an $s$-minimal set with half spaces.
We show that if $E$ is $s$-minimal in $\Omega$ and the exterior data $E_0:=E\setminus\Omega$ lies above a half-space, then
also $E\cap\Omega$ must lie above that same half-space. This is indeed
a very general principle, that we now discuss in full detail.
To this aim, it is convenient to point out that
if $E\subseteq F$ and the boundaries of the two sets touch at a common point $x_0$ where the $s$-fractional mean curvatures coincide, then the two sets must be equal.
The precise result goes as follows:

\begin{lemma}\label{CH:3:curv_rigidity}
Let $E,F\subseteq\R^n$ be such that $E\subseteq F$ and $x_0\in\partial E\cap\partial F$. Then
\begin{equation}\label{CH:3:confront_curv_ineq}
\I_s^\varrho[E](x_0)\geq\I_s^\varrho[F](x_0)\qquad\textrm{for every }\varrho>0.
\end{equation}
Furthermore, if
\begin{equation}\label{CH:3:ineq_for_curvs}
\liminf_{\varrho\to0^+}\I_s^\varrho[F](x_0)\geq a\quad\textrm{and}\quad\limsup_{\varrho\to0^+}\I_s^\varrho[E](x_0)\leq a,
\end{equation}
then $E=F$, the fractional mean curvature is well defined in $x_0$ and $\I_s[E](x_0)=a$.

\begin{proof}
To get $(\ref{CH:3:confront_curv_ineq})$ it is enough to notice that
\[
E\subseteq F\quad\Longrightarrow\quad\big(\chi_{\Co E}(y)-\chi_E(y)\big)
\geq\big(\chi_{\Co F}(y)-\chi_F(y)\big)\qquad\forall\,y\in\R^n.
\]
Now suppose that $(\ref{CH:3:ineq_for_curvs})$ holds true. Then by $(\ref{CH:3:confront_curv_ineq})$ we find that
\[
\exists\,\lim_{\varrho\to0^+}\I_s[E](x_0)=\lim_{\varrho\to0^+}\I_s[F](x_0)=a.
\]

To conclude, notice that if the two curvatures are well defined (in the principal value sense) in $x_0$ and are equal,
then
\begin{equation*}\begin{split}
0\leq\int_{\Co B_\varrho(x_0)}&\frac{\big(\chi_{\Co E}(y)-\chi_E(y)\big)-\big(\chi_{\Co F}(y)-\chi_F(y)\big)}{|x_0-y|^{n+s}}dy\\
&
=\I_s^\varrho[E](x_0)-\I_s^\varrho[F](x_0)\xrightarrow{\varrho\to0^+}0,
\end{split}\end{equation*}
which implies that $\chi_E(y)=\chi_F(y)$ for a.e. $y\in\R^n$, i.e. $E=F$.
\end{proof}\end{lemma}

\begin{prop}\label{CH:3:maximum_principle}[Maximum Principle]
Let $\Omega\subseteq\R^n$ be a bounded open set with $C^{1,1}$ boundary. Let $s\in(0,1)$ and let $E$ be 
$s$-minimal in $\Omega$. If
\begin{equation}\label{CH:3:ext_data_incl}
\{x\cdot\nu\leq a\}\setminus\Omega\subseteq\Co E,\end{equation}
for some $\nu\in\mathbb S^{n-1}$ and $a\in\R$, then
\[\{x\cdot\nu\leq a\}\subseteq \Co E.\]

\begin{proof}
First of all, we remark that up to a rotation and translation, we can suppose that $\nu=e_n$ and $a=0$.
Furthermore we can assume that
\[
\inf_{x\in\overline{\Omega}}x_n<0,
\]
otherwise there is nothing to prove.

If $E\cap\Omega=\emptyset$, i.e. $\Omega\subseteq\Co E$, we are done.
Thus we can suppose that $E\cap\Omega\not =\emptyset$.\\
Since $\overline{E}\cap\overline{\Omega}$ is compact, we 
have
\[
b:=\min_{x\in\overline{E}\cap\overline{\Omega}}x_n\in\R.
\]
Now we consider the set of points which realize the minimum
above, namely we set
$$\mathcal P:=\{p\in\overline{E}\cap\overline{\Omega}\,|\,p_n=b\}.$$
Notice that
\begin{equation}\label{CH:3:confronto_first_incl}
\big\{x_n\leq\min\{b,0\}\big\}\subseteq\Co E,
\end{equation}
so we are reduced to prove that $b\geq0$.

We argue by contradiction and suppose that $b<0$. We will prove that $\mathcal P=\emptyset$.
We remark that $\mathcal P\subseteq\partial E\cap\overline{\Omega}$.

Indeed, if $p\in\mathcal P$, then by $(\ref{CH:3:confronto_first_incl})$ we have that
$B_\delta(p)\cap\{x_n\leq b\}\subseteq\Co E$
for every $\delta >0$, so $|B_\delta(p)\cap\Co E|\geq\frac{\omega_n}{2}\delta^n$
and $p\not\in E_{int}$. Therefore, since $\overline{E}=E_{int}\cup\partial E$, we 
find that $p\in\partial E$.

Roughly speaking, we are sliding upwards the half-space $\{x_n\leq t\}$ until we first touch the set $\overline{E}$. Then the contact points must belong to the boundary of $E$.

Notice that the points of $\mathcal P$ can be either inside $\Omega$ or on $\partial\Omega$.
In both cases we can use the Euler-Lagrange equation to get a contradiction. The precise argument goes as follows.

First, if $p=(p',b)\in\partial E\cap\Omega$, then since $H:=\{x_n\leq b\}\subseteq\Co E$,
we can find an exterior tangent ball to $E$ at $p$ (contained in $\Omega$), so $\I_s[E](p)=0$.

On the other hand, if $p\in\partial E\cap\partial\Omega$, then
$B_{|b|}(p)\setminus\Omega\subseteq\Co E$ and hence (by \cite[Theorem 5.1]{graph})
$\partial E\cap B_r(p)$ is $C^{1,\frac{s+1}{2}}$ for some $r\in(0,|b|)$, and $\I_s[E](p)\leq0$ by Theorem \eqref{CH:3:EL_boundary_coroll} .

In both cases, we have that
\[
p\in\partial H\cap\partial E,\quad H\subseteq \Co E\quad\textrm{and}\quad\I_s[\Co E](p)=-\I_s[E](p)\geq0=\I_s[H](p),
\]
and hence Lemma \ref{CH:3:curv_rigidity} implies $\Co E=H$. However, since $b<0$, this contradicts
$(\ref{CH:3:ext_data_incl})$.

This proves that $b\geq0$, thus concluding the proof.
\end{proof}
\end{prop}

{F}rom this, we obtain a strong comparison principle with planes,
as follows:

\begin{corollary}
Let $\Omega\subseteq\R^n$ be a bounded open set with $C^{1,1}$ boundary.
Let $E\subseteq\R^n$ be $s$-minimal in $\Omega$, with
$\{x_n\leq0\}\setminus\Omega\subseteq\Co E$.
Then

$(i)\quad$ if $|(\Co E\setminus\Omega)\cap\{x_n>0\})|=0$, then $E=\{x_n>0\}$;

$(ii)\quad$ if $|(\Co E\setminus\Omega)\cap\{x_n>0\}|>0$,
then for every $x=(x',0)\in\Omega\cap\{x_n=0\}$ there exists $\delta_x\in(0,d(x,\partial\Omega))$ s.t.
$B_{\delta_x}(x)\subseteq\Co E$. Thus
\begin{equation}
\{x_n\leq0\}\cup\bigcup_{(x',0)\in\Omega}B_{\delta_x}(x)\subseteq\Co E.
\end{equation}

\begin{proof}
First of all, Proposition \ref{CH:3:maximum_principle} guarantees that
\begin{equation*}
\{x_n\leq0\}\subseteq\Co E.
\end{equation*}

$(i)\quad$ Notice that since $E$ is $s$-minimal in $\Omega$, also $\Co E$ is $s$-minimal in $\Omega$.\\
Thus, since $\{x_n>0\}\setminus\Omega\subseteq E=\Co(\Co E)$, we can use again Proposition \ref{CH:3:maximum_principle}
(notice that $\{x_n=0\}$ is a set of measure zero) to
get $\{x_n>0\}\subseteq E$, proving the claim.

$(ii)\quad$ Let $x\in\{x_n=0\}\cap \Omega$.

We argue by contradiction. Suppose that $|B_\delta(x)\cap E|>0$ for every $\delta>0$.\\
Notice that,
since $B_\delta(x)\cap\{x_n\leq0\}\subseteq\Co E$ for every $\delta>0$, this implies that $x\in\partial E\cap \Omega$.
Moreover, we can find an exterior tangent ball to $E$ in $x$, namely
\begin{equation*}
B_\eps(x-\eps\,e_n)\subseteq\{x_n\leq0\}\cap\Omega\subseteq\Co E\cap\Omega.
\end{equation*}
Thus the Euler-Lagrange equation gives $\I_s[E](x)=0$.

Let $H:=\{x_n\leq0\}$.
Since $x\in\partial H$, $H\subseteq\Co E$ and also $\I_s[H](x)=0$,
Lemma \ref{CH:3:curv_rigidity} implies $\Co E=H$.
However this contradicts the hypothesis
\begin{equation*}
|(\Co E\setminus\Omega)\cap\{x_n>0\}|>0,
\end{equation*}
which completes the proof.
\end{proof}\end{corollary}

%
%
\end{chapter}

\begin{chapter}{Some auxiliary results}

\section{Useful integral inequalities} \label{CH:4:app}

We collect here some useful inequalities which we have exploited at various places within the thesis.

We begin with the following simple integral inequality.

%


\begin{lemma}\label{CH:4:dumb_kernel_lemma}
	Let~$n \ge1$,~$s \in (0, 1)$ and~$A, B\subseteq\Rn$ be bounded sets. Then
	$$
	\int_A\int_B\dkers \le \frac{\Ha^{n-1}(\mathbb S^{n-1})}{1-s} \min \big\{ |A|, |B| \big\} \diam (A\cup B)^{1-s}.
	$$
\end{lemma}
\begin{proof}
	Suppose without loss of generality that~$|A| \le |B|$ and set~$D := \diam(A \cup B)$. Then, by changing variables conveniently we estimate
	$$
	\int_A\int_B\dkers \le  \int_A \Big( \int_{B_D} \frac{dz}{|z|^{n - 1 + s}} \Big) \, dx = \Ha^{n-1}(\mathbb S^{n-1}) |A| \int_0^D \frac{d\varrho}{\varrho^{s}},
	$$
	which directly leads to the conclusion.
\end{proof}

Now we prove that a measurable function with finite $W^{s,p}$-seminorm is actually $L^p$-summable and hence belongs to
the fractional Sobolev space $W^{s,p}$.
The proof follows by arguing as in the proof \cite[Theorem 8.2]{HitGuide} 
(see in particular the formula (8.3) there).

\begin{lemma}\label{CH:A:usef_ineq_hit}
	Let $p\in[1,\infty),\,s\in(0,1)$ and let $\Omega\subseteq\R^n$ be a bounded open set. Let $u:\Omega\to\R$ be a measurable
	function such that
	\[
	[u]_{W^{s,p}(\Omega)}^p=\int_\Omega\int_\Omega\frac{|u(x)-u(\xi)|^p}{|x-\xi|^{n+sp}}\,dx\,d\xi<+\infty.
	\]
	Then $u\in W^{s,p}(\Omega)$. More precisely, if $E\subseteq\Omega$ is any measurable set such that
	\eqlab{\label{CH:A:weak_wsp_hp}
		|E|>0\quad\mbox{ and }\quad \int_{E}|u(\xi)|\,d\xi<+\infty,
	}
	then, if we denote
	\bgs{
		M_E:=\int_{E}u(\xi)\,d\xi,
	}
	we have
	\eqlab{\label{CH:A:eqns6}
		\|u\|_{L^p(\Omega)}^p\le\frac{2^{p-1}}{|E|}\left\{(\diam \Omega)^{n+sp}\,[u]^p_{W^{s,p}(\Omega)}
		+|\Omega|\,\frac{|M_E|^p}{|E|^{p-1}}\right\}.
	}
\end{lemma}

\begin{proof}
	First of all, we remark that since $u$ is measurable there exists at least one set $E$ satisfying \eqref{CH:A:weak_wsp_hp}.
	Indeed, for every $k\in\N$ we can consider the set
	\bgs{
		E_k:=\left\{x\in\Omega\,|\,|u(x)|\le k\right\},
	}
	which is measurable. Since $u$ is finite almost everywhere in $\Omega$, there exists $h\in\N$ such that $|E_h|>0$. Then,
	notice that
	\bgs{
		\int_{E_h}|u(\xi)|\,d\xi\le|E_h|h\le|\Omega|h<+\infty,
	}
	so that $E_h$ satisfies \eqref{CH:A:weak_wsp_hp}.
	
	Now let $E$ be any set satisfying \eqref{CH:A:weak_wsp_hp} and define the constant
	\bgs{
		c:=\frac{1}{|E|}\int_{E} u(\xi)\,d\xi=\frac{M_E}{|E|},
	}
	which is finite by hypothesis.
	
	By exploiting Holder's inequality we find
	\bgs{
		|u(x)-c|^p=\frac{1}{|E|^p}\Big|\int_{E}\big(u(x)-u(\xi)\big)d\xi\Big|^p
		\le\frac{1}{|E|}\int_{E}|u(x)-u(\xi)|^p\,d\xi,
	}
	for every $x\in\Omega$.
	Integrating in $x$ over $\Omega$ we obtain
	\bgs{
		\int_\Omega|u(x)-c|^p\,dx\le
		\frac{1}{|E|}\int_\Omega\int_{E}|u(x)-u(\xi)|^p\,dx\,d\xi.
	}
	Since $|x-\xi|\le \mbox{diam }\Omega$ for every $x\in\Omega$ and $\xi\in E\subseteq\Omega$, we conclude that
	\bgs{
		\int_\Omega|u(x)-c|^p\,dx\le
		\frac{1}{|E|}\int_\Omega\int_{E}|u(x)-u(\xi)|^p\,dx\,d\xi
		\le\frac{(\mbox{diam }\Omega)^{n+sp}}{|E|}[u]^p_{W^{s,p}(\Omega)}.
	}
	Finally, we observe that
	\bgs{
		\int_\Omega|c|^p\,dx=|\Omega|\left(\frac{|M_E|}{|E|}\right)^p.
	}
	Therefore
	\bgs{
		\|u\|^p_{L^p(\Omega)}&\le2^{p-1}\int_\Omega|u(x)-c|^p\,dx
		+2^{p-1}\int_\Omega|c|^p\,dx\\
		&
		\le2^{p-1}\frac{(\mbox{diam }\Omega)^{n+sp}}{|E|}[u]^p_{W^{s,p}(\Omega)}
		+2^{p-1}|\Omega|\,\frac{|M_E|^p}{|E|^p},
	}
	proving \eqref{CH:A:eqns6} and concluding the proof of the Lemma.
\end{proof}

Now we prove a ``global version'' of Lemma \ref{CH:A:usef_ineq_hit} in which we use the nonlocal functional
\[
\Nl_s(u,\Omega):=\iint_{\R^{2n}\setminus(\Co\Omega)^2}\frac{|u(x)-u(y)|^2}{|x-y|^{n+2s}}\,dx\,dy,
\]
with $s\in(0,1)$,
in place of the Gagliardo seminorm.
We recall the following definition,
\[
L_s^2(\R^n):=\left\{u:\R^n\to\R\,\big|\,\|u\|_{L^2_s(\R^n)}^2:=\int_{\R^n}\frac{|u(\xi)|^2}{1+|\xi|^{n+2s}}d\xi<\infty\right\}.
\]

\begin{lemma}\label{CH:APP:usef_ineq_tail}
	Let $\Omega\subseteq\R^n$ be a bounded open set and let $s\in(0,1)$. If $u:\R^n\to\R$ is a measurable function such that $\Nl_s(u,\Omega)<\infty$,
	then $u\in L^2_s(\R^n)$. More precisely, if $E\subseteq\Omega$ is any measurable set such that
	\eqlab{\label{CH:ApP:weak_wsp_hp2}
		|E|>0\quad\mbox{ and }\quad \int_{E}|u(\xi)|\,d\xi<\infty,
	}
	then, if we denote
	\bgs{
		M_E:=\int_{E}u(\xi)\,d\xi,
	}
	we have
	\bgs{
		\|u\|_{L^2_s(\R^n)}^2\le \frac{C}{|E|} \left\{\Nl_s(u,\Omega)+\frac{M_E^2}{|E|}\right\},
	}
	for some $C=C(n,s,\Omega)>0$.
\end{lemma}

\begin{proof}
	The proof is similar to that of Lemma \ref{CH:A:usef_ineq_hit}. Again, since $u$ is measurable we know that there exists
	at least one set $E\subseteq\Omega$ satisfying \eqref{CH:ApP:weak_wsp_hp2}.
	
	Now we take a set $E\subseteq\Omega$ which satisfies \eqref{CH:ApP:weak_wsp_hp2}, we define the constant
	\bgs{
		c:=\frac{1}{|E|}\int_{E}u(\xi)\,d\xi=\frac{M_E}{|E|},
	}
	and we remark that
	\bgs{
		|u(x)-c|^2\le\frac{1}{|E|}\int_{E}|u(x)-u(\xi)|^2\,d\xi,
	}
	for every $x\in\R^n$. Integrating in $x$ over $\R^n$, against the weight $1/(1+|x|^{n+2s})$, we find
	\eqlab{\label{CH:APp:eqns37}
		\int_{\R^n}\frac{|u(x)-c|^2}{1+|x|^{n+2s}}\,dx
		\le
		\frac{1}{|E|}\int_{E}\int_{\R^n}\frac{|u(x)-u(\xi)|^2}{1+|x|^{n+2s}}\,d\xi\,dx.
	}
	Now notice that, since $\Omega$ is bounded, there exists a constant $C=C(n,s,\Omega)>0$ such that for every $\xi\in \Omega$
	and every $x\in\R^n$, it holds
	\bgs{
		\frac{1}{1+|x|^{n+2s}}\le C\frac{1}{|x-\xi|^{n+2s}}.
	}
	Thus, from \eqref{CH:APp:eqns37} we obtain
	\bgs{
		\int_{\R^n}\frac{|u(x)-c|^2}{1+|x|^{n+2s}}\,dx
		\le \frac{C}{|E|}\int_{E}\int_{\R^n}\frac{|u(x)-u(\xi)|^2}{|x-\xi|^{n+2s}}\,d\xi\,dx
		\le \frac{C}{|E|}\Nl_s(u,\Omega).
	}
	Finally, notice that
	\bgs{
		\int_{\R^n}\frac{|u(x)|^2}{1+|x|^{n+2s}}\,dx\le
		2\int_{\R^n}\frac{|u(x)-c|^2}{1+|x|^{n+2s}}\,dx
		+2\int_{\R^n}\frac{|c|^2}{1+|x|^{n+2s}}\,dx,
	}
	and
	\bgs{
		\int_{\R^n}\frac{|c|^2}{1+|x|^{n+2s}}\,dx
		= \frac{M_E^2}{|E|^2}\int_{\R^n}\frac{1}{1+|x|^{n+2s}}\,dx.
	}
	This concludes the proof of the Lemma.
\end{proof}

\subsection{Fractional Hardy-type inequality}

We point out the following fractional Hardy-type inequality, which is stated, e.g., in \cite{Dy04}---see formula (17) there. Since the proof for the case $p=1$ is hard to find in the literature, we provide a simple argument based on the fractional Hardy inequality on half-spaces ensured by~\cite[Theorem~1.1]{FS10}.

We recall that $\bar{d}_\Omega$ denotes the signed distance function from $\partial\Omega$, negative inside $\Omega$---see Appendix \ref{CH:1:Appendix_distance_function}. Let us also observe that
\[
|\bar{d}_\Omega(x)|=\mbox{dist}(x,\partial\Omega).
\]

\begin{theorem}\label{CH:4:FHI}
	Let~$n \ge 1$, $p \ge 1$ and let $s \in (0, 1)$ be such that~$s p < 1$. Let $\Omega \subseteq \R^n$ be a bounded open set with Lipschitz boundary. Then, there exists a constant $C=C(n,s,p,\Omega)\geq1$ such that
	\begin{equation} \label{CH:ApP:hardyine}
	\int_\Omega \frac{|u(x)|^p}{|\bar{d}_\Omega(x)|^{s p}} \, dx \le C \| u \|_{W^{s, p}(\Omega)}^p
	\end{equation}
	for every~$u \in W^{s, p}(\Omega)$.
\end{theorem}
\begin{proof}
	We first prove \eqref{CH:ApP:hardyine} for a function $u\in C^\infty_c(\Omega)$, then we extend it to the whole space $W^{s,p}(\Omega)$ by density.
	
	Let~$\{ B^{(j)} \}_{j = 1}^N$ be a sequence of balls of the form~$B^{(j)} = B_r(x^{(j)})$, with $N \in \N$, $r > 0$, and~$x^{(j)} \in \partial \Omega$, for which there exist Lipschitz isomorphisms
	$$
	T_j : B'_2 \times (-2, 2) \longrightarrow 2 B^{(j)} := B_{2 r}(x^{(j)})
	$$
	satisfying
	\begin{align*}
	T_j(U_2) & = 2 B^{(j)}, && \hspace{-50pt} \mbox{with } U_2 := B'_2 \times (-2, 2),\\
	T_j(U^+_2) & = \Omega \cap 2 B^{(j)}, && \hspace{-50pt} \mbox{with } U^+_2 := B_2' \times (0, 2),\\
	T_j(U^0_2) & = \partial \Omega \cap 2 B^{(j)}, && \hspace{-50pt} \mbox{with } U^0_2 := B_2' \times \{ 0 \},
	\end{align*}
	and such that~$\partial \Omega \subseteq \cup_{j = 1}^N B^{(j)}$.
	
	Let~$\varepsilon > 0$ be such that~$\Omega \setminus \cup_{j = 1}^N B^{(j)} \Subset \Omega_{- \varepsilon}$ and set~$B^{(0)} := \Omega_{- \varepsilon}$. Clearly,
	\begin{equation} \label{CH:ApP:B0covered}
	\int_{B^{(0)}} \frac{|u(x)|^p}{|\bar{d}_\Omega(x)|^{s p}} \, dx \le \varepsilon^{- s p} \int_{B^{(0)}} |u(x)|^p \, dx \le C \| u \|_{L^p(\Omega)}^p,
	\end{equation}
	where, from now on,~$C$ denotes any constant larger than~$1$, whose value depend at most on~$n$,~$s$,~$p$, and~$\Omega$.
	
	Notice that~$\{ B^{(j)} \}_{j = 0}^N$ is an oper cover of~$\Omega$ and let~$\{ \eta_j \}_{j = 0}^N$ be a smooth partition of unity on~$\Omega$ subordinate to~$\{ B^{(j)} \}_{j = 0}^N$.
	
	For~$j =1, \ldots, N$, we define~$v_j := \eta_j u \in C^\infty_c(\Omega\cap B^{(j)})$. Changing variables through~$T_j$, we have
	$$
	\int_{\Omega \cap B^{(j)}} \frac{|v_j(x)|^p}{|\bar{d}_\Omega(x)|^{s p}} \, dx = \int_{T_j^{-1}(\Omega \cap B^{(j)})} \frac{|v_j(T_j(\bar{x}))|^p}{|\bar{d}_\Omega(T_j(\bar{x}))|^{s p}} |\det D T_j(\bar{x})| \, d\bar{x}.
	$$
	Notice that for every~$x \in \Omega \cap B^{(j)}$ there exists~$D_j(x) \in \partial \Omega \cap 2 B^{(j)}$ such that~$|\bar{d}_\Omega(x)| = |x - D_j(x)|$. Since~$T_j$ is bi-Lipschitz and~$T_j^{-1}(D_j(x)) \in B_2' \times \{ 0 \}$, we have
	$$
	\begin{aligned}
	|\bar{d}_\Omega(T_j(\bar{x}))| & = |T_j(\bar{x}) - D_j(T_j(\bar{x}))| = |T_j(\bar{x}) - T_j(T_j^{-1}(D_j(T_j(\bar{x}))))| \\
	& \ge C^{-1} |\bar{x} - T_j^{-1}(D_j(T_j(\bar{x})))| \ge C^{-1} \bar{x}_n
	\end{aligned}
	$$
	for every~$\bar{x} \in T_j^{-1}(\Omega \cap B^{(j)})$. Accordingly, writing~$w_j := v_j \circ T_j$ we get
	$$
	\int_{\Omega \cap B^{(j)}} \frac{|v_j(x)|^p}{|\bar{d}_\Omega(x)|^{s p}} \, dx \le C \int_{U^+_2} \frac{|w_j(\bar{x}))|^p}{|\bar{x}_n|^{s p}} \, d\bar{x}.
	$$
	
	Let us observe that $w_j$ is supported inside~$T_j^{-1}(\Omega \cap B^{(j)})$.
	We now employ the fractional Hardy inequality on half-spaces---e.g.,~\cite[Theorem~1.1]{FS10}---and deduce that
	\begin{equation} \label{CH:ApP:hardytech1}
	\int_{\Omega \cap B^{(j)}} \frac{|v_j(x)|^p}{|\bar{d}_\Omega(x)|^{s p}} \, dx \le C \int_{\R^n_+} \int_{\R^n_+} \frac{|w_j(\bar{x}) - w_j(\bar{y})|^p}{|\bar{x} - \bar{y}|^{n + s p}} \, d\bar{x}\, d\bar{y},
	\end{equation}
	where~$\R^n_+ = \{ z \in \R^n \,|\, z_n > 0 \}$ and it is understood that $w_j$
	is extended by 0 in $\R^n_+\setminus U_2^+$. We point out that---since $T_j^{-1}(B^{(j)}) \Subset U_2$ and $T_j^{-1}(\Omega \cap B^{(j)})\subseteq U_2^+$---we have
	\[
	\mbox{dist}\big(T_j^{-1}(\Omega \cap B^{(j)}),\R^n_+ \setminus U_2^+\big)>0.
	\]
	Thus, using that~$w_j$ is supported inside~$T_j^{-1}(\Omega \cap B^{(j)})$,
	we estimate
	\begin{equation} \label{CH:ApP:hardytech2}
	\begin{aligned}
	\int_{\R^n_+} \int_{\R^n_+} \frac{|w_j(\bar{x}) - w_j(\bar{y})|^p}{|\bar{x} - \bar{y}|^{n + s p}} \, d\bar{x}\, d\bar{y} & \le \int_{U_2^+} \int_{U_2^+} \frac{|w_j(\bar{x}) - w_j(\bar{y})|^p}{|\bar{x} - \bar{y}|^{n + s p}} \, d\bar{x} \,d\bar{y} \\
	& \quad + 2 \int_{T_j^{-1}(\Omega \cap B^{(j)})} \left( \int_{\R^n_+ \setminus U_2^+} \frac{ |w_j(\bar{x})|^p}{|\bar{x} - \bar{y}|^{n + s p}} \, d\bar{y} \right) d\bar{x} \\
	& \le \int_{U_2^+} \int_{U_2^+} \frac{|w_j(\bar{x}) - w_j(\bar{y})|^p}{|\bar{x} - \bar{y}|^{n + s p}} \, d\bar{x}\, d\bar{y} + C \| w_j \|_{L^p(U_2^+)}^p.
	\end{aligned}
	\end{equation}
	
	By combining~\eqref{CH:ApP:hardytech1} with~\eqref{CH:ApP:hardytech2} and switching back to the variables in~$\Omega$, we easily find that
	$$
	\int_{\Omega \cap B^{(j)}} \frac{|v_j(x)|^p}{|\bar{d}_\Omega(x)|^{s p}} \, dx \le C \left( \int_{\Omega \cap 2 B^{(j)}} \int_{\Omega \cap 2 B^{(j)}} \frac{|v_j(x) - v_j(y)|^p}{|x - y|^{n + s p}} \, dx \,dy + \| v_j \|_{L^p(\Omega \cap 2 B^{(j)})}^p \right).
	$$
	Recalling that~$v_j = \eta_j u$ and~supp$(\eta_j) \Subset B^{(j)}$, a simple computation then leads us to
	$$
	\int_{\Omega \cap B^{(j)}} \frac{|v_j(x)|^p}{|\bar{d}_\Omega(x)|^{s p}} \, dx \le C \| u \|_{W^{s, p}(\Omega)}^p \quad \mbox{for all } j = 1, \ldots, N.
	$$
	Then, estimate~\eqref{CH:ApP:hardyine} for $u\in C^\infty_c(\Omega)$ follows by putting together this with~\eqref{CH:ApP:B0covered} and using that~$\{ \eta_j \}$ is a partition of unity, whereas the general case of $u\in W^{s,p}(\Omega)$ is obtained by density.
	More precisely, let $u\in W^{s,p}(\Omega)$ and notice that by the density of $C^\infty_c(\Omega)$ in $W^{s,p}(\Omega)$---see, e.g., Theorem \ref{CH:4:smooth_cpt_dense}---we can find $\{u_k\}\subseteq C_c^\infty(\Omega)$ such
	that
	\bgs{
		\lim_{k\to\infty}\|u-u_k\|_{W^{s,p}(\Omega)}=0.
	}
	Up to passing to a subsequence, we can further suppose that
	$u_k\to u$ a.e. in $\Omega$.
	Then, by Fatou's Lemma we find
	\bgs{
		\int_\Omega\frac{|u(x)|^p}{|\bar{d}_\Omega(x)|^{sp}}dx\le
		\liminf_{k\to\infty}\int_\Omega\frac{|u_k(x)|^p}{|\bar{d}_\Omega(x)|^{sp}}dx\le
		\lim_{k\to\infty} C\|u_k\|_{W^{s,p}(\Omega)}^p
		= C\|u\|_{W^{s,p}(\Omega)}^p,
	}
	concluding the proof of the Theorem.
\end{proof}

\begin{corollary}\label{CH:4:FHI_corollary}
	Let $n\ge1$ and let $\Omega\subseteq\Rn$ be a bounded open set with Lipschitz boundary. Let $p\geq 1$ and $s\in(0,1)$ be such that
	$sp<1$. Then
	\eqlab{
		\int_\Omega\left(\int_{\Co\Omega}\frac{|u(x)|^p}{|x-y|^{n+sp}}dy\right)dx\le
		C(n,s,p,\Omega)\|u\|^p_{W^{s,p}(\Omega)},
	}
	for every $u\in W^{s,p}(\Omega)$.
\end{corollary}
\begin{proof}
	It is enough to notice that
	\bgs{
		\int_\Omega\left(\int_{\Co\Omega}\frac{|u(x)|^p}{|x-y|^{n+sp}}dy\right)dx&\le
		\int_\Omega\left(\int_{\Co B_{|\bar{d}_\Omega(x)|}(x)}\frac{dy}{|x-y|^{n+sp}}\right)|u(x)|^p\,dx\\
		&
		=\frac{\Ha^{n-1}(\mathbb S^{n-1})}{sp}\int_\Omega
		\frac{|u(x)|^p}{|\bar{d}_\Omega(x)|^{sp}}dx.
	}
	Then the conclusion follows from Theorem \ref{CH:4:FHI}.
\end{proof}

\subsection{Fractional Poincar\'e-type inequality}

For the convenience of the reader, we provide a proof of the following well known fractional Poincar\'e-type inequality.

\begin{prop}\label{CH:4:FPI}
	Let $\Omega\subseteq\Op\subseteq\Rn$ be bounded open sets such that $|\Op\setminus\Omega|>0$ and let $p\in[1,\infty)$ and $s\in(0,1)$.
	Let $u:\Op\lra\R$ be such that $u=0$ almost everywhere in $\Op\setminus\Omega$. Then
	\eqlab{\label{CH:4:FPI_eq_ineq}
		\|u\|_{L^p(\Omega)}^p
		\le\frac{(\diam  \Op)^{n+sp}}{|\Op\setminus\Omega|}\int_\Omega\int_{\Op\setminus\Omega}\frac{|u(x)|^p}{|x-y|^{n+sp}}\,dx\,dy
		\le\frac{(\diam  \Op)^{n+sp}}{|\Op\setminus\Omega|}\,[u]^p_{W^{s,p}(\Op)}.
	}
\end{prop}
\begin{proof}
	Notice that
	\bgs{
		|u(x)|=|u(x)-u(y)|\qquad\textrm{for almost every }(x,y)\in\Omega\times(\Op\setminus\Omega).
	}
	Hence
	\bgs{
		|u(x)|^p=\frac{1}{|\Op\setminus\Omega|}\int_{\Op\setminus\Omega}|u(x)-u(y)|^p\,dy
		=\frac{1}{|\Op\setminus\Omega|}\int_{\Op\setminus\Omega}\frac{|u(x)-u(y)|^p}{|x-y|^{n+sp}}|x-y|^{n+sp}\,dy.
	}
	Since
	\bgs{
		|x-y|\le\textrm{diam }\Op\qquad\forall\,(x,y)\in\Omega\times(\Op\setminus\Omega),
	}
	we obtain
	\bgs{
		|u(x)|^p\le\frac{(\textrm{diam }\Op)^{n+sp}}{|\Op\setminus\Omega|}\int_{\Op\setminus\Omega}
		\frac{|u(x)-u(y)|^p}{|x-y|^{n+sp}}\,dy.
	}
	Integrating over $\Omega$ gives
	\bgs{
		\|u\|^p_{L^p(\Omega)}\le\frac{(\textrm{diam }\Op)^{n+sp}}{|\Op\setminus\Omega|}
		\int_\Omega\int_{\Op\setminus\Omega}\frac{|u(x)-u(y)|^p}{|x-y|^{n+sp}}\,dx\,dy,
	}
	hence the claim.
\end{proof}

\section{Density of compactly supported smooth functions}\label{CH:4:Density_appendix}

As customary, we denote by $W^{s,p}_0(\Omega)$ the closure of $C^\infty_c(\Omega)$
in $W^{s,p}(\Omega)$ with respect to the usual $W^{s,p}$-norm.

The aim of this section consists in providing a proof of the well known fact that, when $sp<1$, the space $C^\infty_c(\Omega)$ is dense in $W^{s,p}(\Omega)$. Roughly speaking, this means that, in this case, the space $W^{s,p}(\Omega)$ has no well defined trace on $\partial\Omega$.

\begin{theorem}\label{CH:4:smooth_cpt_dense}
	Let $\Omega\subseteq\Rn$ be a bounded open set with Lipschitz boundary and let $p\in[1,\infty)$, $s\in(0,1)$. Then
	\bgs{
		sp<1\quad\Longrightarrow\quad W^{s,p}_0(\Omega)=W^{s,p}(\Omega),
	}
	i.e. $C_c^\infty(\Omega)$ is dense in $W^{s,p}(\Omega)$.
\end{theorem}

The proof of this well known theorem is the consequence of the following results.

\begin{lemma}
	Let $\Omega\subseteq\Rn$ be an open set and let $p\in[1,\infty)$, $s\in(0,1)$. Then
	\bgs{
		W^{s,p}(\Omega)\cap L^\infty(\Omega)\quad\textrm{is dense in}\quad W^{s,p}(\Omega).
	}
\end{lemma}

\begin{proof}
	Given $u\in W^{s,p}(\Omega)$, consider the functions
	\begin{equation*}
	u_k:=\left\{\begin{array}{cc}
	u & \textrm{in }\{|u|\le k\},\\
	k & \textrm{in }\{u\ge k\},\\
	-k & \textrm{in }\{u\le-k\}.
	\end{array}\right.
	\end{equation*}
	Then
	\bgs{
		|u_k|^p\le|u|^p\quad\textrm{a.e. in }\Omega\quad\textrm{and}\quad u_k\lra u\quad\textrm{a.e. in }\Omega,
	}
	hence
	\bgs{
		\lim_{k\to\infty}\|u-u_k\|_{L^p(\Omega)}=0,
	}
	by the dominated convergence Theorem.
	Similarly, since
	\bgs{
		\frac{|u_k(x)-u_k(y)|^p}{|x-y|^{n+sp}}\le\frac{|u(x)-u(y)|^p}{|x-y|^{n+sp}}\quad\textrm{for a.e. }(x,y)\in\Omega\times\Omega,
	}
	by using again the dominated convergence Theorem, we find
	\bgs{
		\lim_{k\to\infty}[u-u_k]_{W^{s,p}(\Omega)}=0,
	}
	concluding the proof.
\end{proof}

Now we consider a symmetric mollifier, that is $\eta\in C_c^\infty(\Rn)$ such that
\eqlab{\label{CH:4:symm_moll_def}
	\eta\ge0,\quad\int_{\Rn}\eta\,dx=1,\quad\eta(x)=\eta(-x)\quad\textrm{and supp }\eta\subseteq B_1.
}
We set
\bgs{
	\eta_\eps(x):=\frac{1}{\eps^n}\eta\Big(\frac{x}{\eps}\Big),
}
for every $\eps\in(0,1)$.

We recall the following well known result:

\begin{lemma}\label{CH:4:convol_conv_frac}
	Let $\Omega\subseteq\Rn$ be an open set and let $p\in[1,\infty)$, $s\in(0,1)$.
	Then for every $u\in W^{s,p}(\Omega)$ it holds
	\bgs{
		\lim_{\eps\to0^+}\|u-u\ast\eta_\eps\|_{W^{s,p}(\Omega')}=0\qquad\forall\,\Omega'\Subset\Omega.
	}
\end{lemma}

We only observe that the proof of Lemma \ref{CH:4:convol_conv_frac} can be obtained by arguing as in the proof of point $(i)$ of Lemma \ref{CH:2:dens_lemma}.

\begin{prop}
	Let $\Omega\subseteq\Rn$ be a bounded open set with Lipschitz boundary and let $p\in[1,\infty)$, $s\in(0,1)$. Then
	\bgs{
		sp<1\quad\Longrightarrow\quad C_c^\infty(\Omega)\quad\textrm{is dense in}\quad W^{s,p}(\Omega)\cap L^\infty(\Omega).
	}
\end{prop}

\begin{proof}
	Let $\sigma:=sp\in(0,1)$ and let $u\in W^{s,p}(\Omega)\cap L^\infty(\Omega)$.
	For $\delta>0$ small enough, let
	\bgs{
		u_\delta:=u\chi_{\Omega_{-\delta}}.
	}
	Then
	\eqlab{\label{CH:4:abacab}
		\lim_{\delta\to0^+}\|u-u_\delta\|_{W^{s,p}(\Omega)}=0.
	}
	Indeed
	\bgs{
		\|u-u_\delta\|_{L^p(\Omega)}^p\le\|u\|_{L^\infty(\Omega)}^p|\Omega\setminus\Omega_{-\delta}|\xrightarrow{\delta\to0^+}0.
	}
	We remark that, since $\Omega$ is bounded and has Lipschitz boundary, and since $\sigma\in(0,1)$, by Lemma \ref{CH:2:UniFNeiGHboEstiLeMmA}
	we have
	\eqlab{\label{CH:4:artropodi}
		\int_{\Omega_{-\delta}}\int_{\Omega\setminus\Omega_{-\delta}}\frac{dx\,dy}{|x-y|^{n+\sigma}}\le C(n,\Omega,\sigma)\delta^{1-\sigma}.
	}
	Then
	\bgs{
		\int_\Omega\int_\Omega&\frac{|u(x)(1-\chi_{\Omega_{-\delta}}(x))-u(y)(1-\chi_{\Omega_{-\delta}}(y))|^p}{|x-y|^{n+\sigma}}dx\,dy\\
		&\qquad\qquad
		=2\int_{\Omega_{-\delta}}\int_{\Omega\setminus\Omega_{-\delta}}\frac{|u(y)|^p}{|x-y|^{n+\sigma}}dx\,dy
		+[u]^p_{W^{s,p}(\Omega\setminus\Omega_{-\delta})}\\
		&\qquad\qquad
		\le2\|u\|_{L^\infty(\Omega)}^pC(n,\Omega,\sigma)\delta^{1-\sigma}+[u]^p_{W^{s,p}(\Omega\setminus\Omega_{-\delta})}.
	}
	Notice that, since $|\Omega\setminus\Omega_{-\delta}|\xrightarrow{\delta\to0^+}0$, we get by the dominated convergence Theorem
	\bgs{
		\lim_{\delta\to0^+}[u]^p_{W^{s,p}(\Omega\setminus\Omega_{-\delta})}=0.
	}
	Therefore
	\bgs{
		\lim_{\delta\to0^+}[u-u_\delta]_{W^{s,p}(\Omega)}=0,
	}
	proving \eqref{CH:4:abacab}.
	
	Now we consider the $\eps$-regularization of the function $u_\delta$.\\
Notice that for every $\eps\in(0,\delta/4)$
	\bgs{
		\textrm{supp}(u_\delta\ast\eta_\eps)\Subset\Omega_{-\frac{\delta}{2}},
	}
	since the $\eps$-convolution enlarges the support at most to an $\eps$-neighborhood of the original function.
	It is well known that---since $u_\delta$ is compactly supported inside $\Omega$---we have
	\eqlab{\label{CH:4:cefalopodi}
		\lim_{\eps\to0^+}\|u_\delta-u_\delta\ast\eta_\eps\|_{L^p(\Omega)}=0.
	}
	Moreover
	\bgs{
		\|u_\delta\ast\eta_\eps\|_{L^\infty(\Omega)}\le\|u_\delta\|_{L^\infty(\Omega)}\le\|u\|_{L^\infty(\Omega)}.
	}
	Thus, by \eqref{CH:4:artropodi}
	\begin{equation*}
	[u_\delta-u_\delta\ast\eta_\eps]^p_{W^{s,p}(\Omega)}
	\le[u_\delta-u_\delta\ast\eta_\eps]^p_{W^{s,p}(\Omega_{-\delta/2})}
	+2\|u\|_{L^\infty(\Omega)}^pC(n,\Omega,\sigma)\Big(\frac{\delta}{2}\Big)^{1-\sigma}.
	\end{equation*}
	
	By Lemma \ref{CH:4:convol_conv_frac} we have
	\bgs{
		\lim_{\eps\to0^+}[u_\delta-u_\delta\ast\eta_\eps]_{W^{s,p}(\Omega_{-\delta/2})}=0.
	}
	Hence, recalling \eqref{CH:4:cefalopodi},
	we can find $\eps_\delta\in(0,\delta/4)$ small enough such that, if we set
	\bgs{
		\tilde u_\delta:=u_\delta\ast\eta_{\eps_\delta}\in C_c^\infty(\Omega),
	}
	then
	\eqlab{\label{CH:4:tartaruga}
		\|u_\delta-\tilde u_\delta\|_{L^p(\Omega)}\le\delta\quad\textrm{and}
		\quad[u_\delta-\tilde u_\delta]^p_{W^{s,p}(\Omega)}\le\delta+ C\delta^{1-\sigma}.
	}
	Then, by \eqref{CH:4:tartaruga} and \eqref{CH:4:abacab} we obtain
	\bgs{
		\lim_{\delta\to0^+}\|u-\tilde u_\delta\|_{W^{s,p}(\Omega)}=0,
	}
	concluding the proof.
\end{proof}

\end{chapter}

\end{appendix}


\begin{thebibliography}{100}

\bibitem{MR3393247}
Nicola Abatangelo.
\newblock Large {$S$}-harmonic functions and boundary blow-up solutions for the
  fractional {L}aplacian.
\newblock {\em Discrete Contin. Dyn. Syst.}, 35(12):5555--5607, 2015.

\bibitem{Abaty}
Nicola Abatangelo and Enrico Valdinoci.
\newblock A notion of nonlocal curvature.
\newblock {\em Numer. Funct. Anal. Optim.}, 35(7-9):793--815, 2014.

\bibitem{AMB}
Luigi Ambrosio.
\newblock {\em Corso introduttivo alla teoria geometrica della misura ed alle
  superfici minime}.
\newblock Appunti dei Corsi Tenuti da Docenti della Scuola. [Notes of Courses
  Given by Teachers at the School]. Scuola Normale Superiore, Pisa, 1997.

\bibitem{Ambrosio}
Luigi Ambrosio and Norman Dancer.
\newblock {\em Calculus of variations and partial differential equations}.
\newblock Springer-Verlag, Berlin, 2000.
\newblock Topics on geometrical evolution problems and degree theory, Papers
  from the Summer School held in Pisa, September 1996, Edited by G. Buttazzo,
  A. Marino and M. K. V. Murthy.

\bibitem{Gamma}
Luigi Ambrosio, Guido De~Philippis, and Luca Martinazzi.
\newblock Gamma-convergence of nonlocal perimeter functionals.
\newblock {\em Manuscripta Math.}, 134(3-4):377--403, 2011.

\bibitem{ACKS}
I.~Athanasopoulos, L.~A. Caffarelli, C.~Kenig, and S.~Salsa.
\newblock An area-{D}irichlet integral minimization problem.
\newblock {\em Comm. Pure Appl. Math.}, 54(4):479--499, 2001.

\bibitem{BAL08}
M.~Ballerini, N.~Cabibbo, R.~Candelier, A.~Cavagna, E.~Cisbani, I.~Giardina,
  V.~Lecomte, A.~Orlandi, G.~Parisi, A.~Procaccini, M.~Viale, and
  V.~Zdravkovic.
\newblock Interaction ruling animal collective behavior depends on topological
  rather than metric distance: Evidence from a field study.
\newblock {\em Proceedings of the National Academy of Sciences},
  105(4):1232--1237, 2008.

\bibitem{BARBU}
Viorel Barbu.
\newblock {\em Differential equations}.
\newblock Springer Undergraduate Mathematics Series. Springer, Cham, 2016.
\newblock Translated from the 1985 Romanian original by Liviu Nicolaescu.

\bibitem{bootstrap}
Bego{\~{n}}a {Barrios}, Alessio {Figalli}, and Enrico {Valdinoci}.
\newblock {Bootstrap regularity for integro-differential operators and its
  application to nonlocal minimal surfaces.}
\newblock {\em {Ann. Sc. Norm. Super. Pisa, Cl. Sci. (5)}}, 13(3):609--639,
  2014.

\bibitem{Bellettini}
Giovanni Bellettini.
\newblock {\em Lecture notes on mean curvature flow, barriers and singular
  perturbations}, volume~12 of {\em Appunti. Scuola Normale Superiore di Pisa
  (Nuova Serie) [Lecture Notes. Scuola Normale Superiore di Pisa (New
  Series)]}.
\newblock Edizioni della Normale, Pisa, 2013.

\bibitem{MX13a}
Andrew Berdahl, Colin~J. Torney, Christos~C. Ioannou, Jolyon~J. Faria, and
  Iain~D. Couzin.
\newblock Emergent sensing of complex environments by mobile animal groups.
\newblock {\em Science}, 339(6119):574--576, 2013.

\bibitem{MR3304346}
Andrea~L. Bertozzi, Jesus Rosado, Martin~B. Short, and Li~Wang.
\newblock Contagion shocks in one dimension.
\newblock {\em J. Stat. Phys.}, 158(3):647--664, 2015.

\bibitem{BDM69}
E.~Bombieri, E.~De~Giorgi, and M.~Miranda.
\newblock Una maggiorazione a priori relativa alle ipersuperfici minimali non
  parametriche.
\newblock {\em Arch. Rational Mech. Anal.}, 32:255--267, 1969.

\bibitem{BBM}
Jean Bourgain, Ha\"im Brezis, and Petru Mironescu.
\newblock Limiting embedding theorems for {$W^{s,p}$} when {$s\uparrow1$} and
  applications.
\newblock {\em J. Anal. Math.}, 87:77--101, 2002.
\newblock Dedicated to the memory of Thomas H. Wolff.

\bibitem{LuCla}
Claudia Bucur and Luca Lombardini.
\newblock Asymptotics as $s\searrow0$ of the nonlocal nonparametric {P}lateau
  problem with obstacles.
\newblock {\em In preparation}.

\bibitem{BLV16}
Claudia Bucur, Luca Lombardini, and Enrico Valdinoci.
\newblock Complete stickiness of nonlocal minimal surfaces for small values of
  the fractional parameter.
\newblock {\em arXiv preprint arXiv:1612.08295}, 2016.
\newblock Accepted for publication by {A}nnales de l'{I}nstitut {H}enri
  {P}oincar\'e {A}nalyse {N}on {L}in\'eaire.

\bibitem{bucval}
Claudia Bucur and Enrico Valdinoci.
\newblock {\em Nonlocal diffusion and applications}, volume~20 of {\em Lecture
  Notes of the Unione Matematica Italiana}.
\newblock Springer, [Cham]; Unione Matematica Italiana, Bologna, 2016.

\bibitem{CCS17}
Xavier Cabr\'e, Eleonora Cinti, and Joaquim Serra.
\newblock Stable $s$-minimal cones in $\mathbb{R}^3$ are flat for $s \sim 1$.
\newblock {\em arXiv preprint arXiv:1710.08722}, 2017.

\bibitem{CaCo}
Xavier Cabr\'e and Matteo Cozzi.
\newblock A gradient estimate for nonlocal minimal graphs.
\newblock {\em arXiv preprint arXiv:1711.08232}, 2017.

\bibitem{GeomObst}
Luis Caffarelli, Daniela De~Silva, and Ovidiu Savin.
\newblock Obstacle-type problems for minimal surfaces.
\newblock {\em Communications in Partial Differential Equations},
  41(8):1303--1323, 2016.

\bibitem{CRS10}
Luis Caffarelli, Jean-Michel Roquejoffre, and Ovidiu Savin.
\newblock Nonlocal minimal surfaces.
\newblock {\em Comm. Pure Appl. Math.}, 63(9):1111--1144, 2010.

\bibitem{CSV}
Luis Caffarelli, Ovidiu Savin, and Enrico Valdinoci.
\newblock Minimization of a fractional perimeter-{D}irichlet integral
  functional.
\newblock {\em Ann. Inst. H. Poincar\'{e} Anal. Non Lin\'{e}aire},
  32(4):901--924, 2015.

\bibitem{CS07}
Luis Caffarelli and Luis Silvestre.
\newblock An extension problem related to the fractional {L}aplacian.
\newblock {\em Comm. Partial Differential Equations}, 32(7-9):1245--1260, 2007.

\bibitem{uniform}
Luis Caffarelli and Enrico Valdinoci.
\newblock Uniform estimates and limiting arguments for nonlocal minimal
  surfaces.
\newblock {\em Calc. Var. Partial Differential Equations}, 41(1-2):203--240,
  2011.

\bibitem{regularity}
Luis Caffarelli and Enrico Valdinoci.
\newblock Regularity properties of nonlocal minimal surfaces via limiting
  arguments.
\newblock {\em Adv. Math.}, 248:843--871, 2013.

\bibitem{C65}
Shiing-shen Chern.
\newblock On the curvatures of a piece of hypersurface in euclidean space.
\newblock {\em Abh. Math. Sem. Univ. Hamburg}, 29:77--91, 1965.

\bibitem{CH2}
Andre Chiaradia, John McBride, Tanya Murray, and Peter Dann.
\newblock Effect of fog on the arrival time of little penguins {E}udyptula
  minor: a clue for visual orientation?
\newblock {\em J. of Ornithol.}, 148(2):229--233, 2007.

\bibitem{CSV16}
Eleonora Cinti, Joaquim Serra, and Enrico Valdinoci.
\newblock Quantitative flatness results and ${BV}$-estimates for stable
  nonlocal minimal surfaces.
\newblock {\em arXiv preprint arXiv:1602.00540}, 2016.
\newblock To appear in Journal of Differential Geometry.

\bibitem{mattheorem}
Matteo Cozzi.
\newblock On the variation of the fractional mean curvature under the effect of
  {$C^{1,\alpha}$} perturbations.
\newblock {\em Discrete Contin. Dyn. Syst.}, 35(12):5769--5786, 2015.

\bibitem{CDV17}
Matteo Cozzi, Serena Dipierro, and Enrico Valdinoci.
\newblock Planelike interfaces in long-range {I}sing models and connections
  with nonlocal minimal surfaces.
\newblock {\em J. Stat. Phys.}, 167(6):1401--1451, 2017.

\bibitem{CFL18}
Matteo Cozzi, Alberto Farina, and Luca Lombardini.
\newblock Bernstein-{M}oser-type results for nonlocaal minimal graphs.
\newblock {\em arXiv preprint arXiv:1807.05774}, 2018.

\bibitem{MR2744705}
Emiliano Cristiani, Benedetto Piccoli, and Andrea Tosin.
\newblock Modeling self-organization in pedestrians and animal groups from
  macroscopic and microscopic viewpoints.
\newblock In {\em Mathematical modeling of collective behavior in
  socio-economic and life sciences}, Model. Simul. Sci. Eng. Technol., pages
  337--364. Birkh\H{a}user Boston, Inc., Boston, MA, 2010.

\bibitem{CH1}
T.~Daniel, A.~Chiaradia, M.~Logan, G.~Quinn, and R.~Reina.
\newblock Synchronized group association in little penguins, {E}udyptula
  {M}inor.
\newblock {\em Anim. behav.}, 74(5):1241--1248, 2007.

\bibitem{GuyDavid}
Guy David.
\newblock {$C^{1+\alpha}$}-regularity for two-dimensional almost-minimal sets
  in {$\Bbb R^n$}.
\newblock {\em J. Geom. Anal.}, 20(4):837--954, 2010.

\bibitem{Davila}
Juan D\'avila.
\newblock On an open question about functions of bounded variation.
\newblock {\em Calc. Var. Partial Differential Equations}, 15(4):519--527,
  2002.

\bibitem{DdPW18}
Juan D\'avila, Manuel del Pino, and Juncheng Wei.
\newblock Nonlocal {$s$}-minimal surfaces and {L}awson cones.
\newblock {\em J. Differential Geom.}, 109(1):111--175, 2018.

\bibitem{Demengel}
Fran\c{c}oise Demengel.
\newblock Fonctions \`a hessien born\'e.
\newblock {\em Ann. Inst. Fourier (Grenoble)}, 34(2):155--190, 1984.

\bibitem{HitGuide}
Eleonora Di~Nezza, Giampiero Palatucci, and Enrico Valdinoci.
\newblock Hitchhiker's guide to the fractional {S}obolev spaces.
\newblock {\em Bull. Sci. Math.}, 136(5):521--573, 2012.

\bibitem{MR1022305}
R.~J. DiPerna and P.-L. Lions.
\newblock Ordinary differential equations, transport theory and {S}obolev
  spaces.
\newblock {\em Invent. Math.}, 98(3):511--547, 1989.

\bibitem{DFPV13}
Serena {Dipierro}, Alessio {Figalli}, Giampiero {Palatucci}, and Enrico
  {Valdinoci}.
\newblock {Asymptotics of the $s$-perimeter as $s \searrow 0$.}
\newblock {\em {Discrete Contin. Dyn. Syst.}}, 33(7):2777--2790, 2013.

\bibitem{DMLV16}
Serena Dipierro, Pietro Miraglio, Luca Lombardini, and Enrico Valdinoci.
\newblock The {P}hillip {I}sland penguin parade (a mathematical treatment).
\newblock {\em arXiv preprint arXiv:1611.08715}, 2018.
\newblock To appear in The ANZIAM Journal.

\bibitem{DSV}
Serena Dipierro, Ovidiu Savin, and Enrico Valdinoci.
\newblock A nonlocal free boundary problem.
\newblock {\em SIAM J. Math. Anal.}, 47(6):4559--4605, 2015.

\bibitem{graph}
Serena Dipierro, Ovidiu Savin, and Enrico Valdinoci.
\newblock Graph properties for nonlocal minimal surfaces.
\newblock {\em Calc. Var. Partial Differential Equations}, 55(4):Art. 86, 25,
  2016.

\bibitem{DSV17}
Serena Dipierro, Ovidiu Savin, and Enrico Valdinoci.
\newblock All functions are locally {$s$}-harmonic up to a small error.
\newblock {\em J. Eur. Math. Soc. (JEMS)}, 19(4):957--966, 2017.

\bibitem{boundary}
Serena Dipierro, Ovidiu Savin, and Enrico Valdinoci.
\newblock Boundary behavior of nonlocal minimal surfaces.
\newblock {\em J. Funct. Anal.}, 272(5):1791--1851, 2017.

\bibitem{DV-onephase}
Serena Dipierro and Enrico Valdinoci.
\newblock Continuity and density results for a one-phase nonlocal free boundary
  problem.
\newblock {\em Ann. Inst. H. Poincar\'{e} Anal. Non Lin\'{e}aire},
  34(6):1387--1428, 2017.

\bibitem{DV18}
Serena Dipierro and Enrico Valdinoci.
\newblock Nonlocal minimal surfaces: interior regularity, quantitative
  estimates and boundary stickiness.
\newblock In Giampiero Palatucci and Tuomo Kuusi, editors, {\em Recent
  Developments in the Nonlocal Theory}, chapter~4, pages 165--209. De Gruyter
  Open, 2018.

\bibitem{Doktor}
Pavel {Doktor}.
\newblock {Approximation of domains with {L}ipschitzian boundary.}
\newblock {\em {\v{C}as. P\v{e}stov\'an\'{\i} Mat.}}, 101:237--255, 1976.

\bibitem{MR2899094}
Ning Du, Jintu Fan, Huijun Wu, Shuo Chen, and Yang Liu.
\newblock An improved model of heat transfer through penguin feathers and down.
\newblock {\em J. Theoret. Biol.}, 248(4):727--735, 2007.

\bibitem{Dy04}
Bart\l~omiej Dyda.
\newblock A fractional order {H}ardy inequality.
\newblock {\em Illinois J. Math.}, 48(2):575--588, 2004.

\bibitem{Falconer}
Kenneth Falconer.
\newblock {\em Fractal geometry}.
\newblock John Wiley \& Sons, Ltd., Chichester, third edition, 2014.
\newblock Mathematical foundations and applications.

\bibitem{FarRog}
Daniel Faraco and Keith~M. Rogers.
\newblock The {S}obolev norm of characteristic functions with applications to
  the {C}alder\'on inverse problem.
\newblock {\em Q. J. Math.}, 64(1):133--147, 2013.

\bibitem{F15}
Alberto Farina.
\newblock A {B}ernstein-type result for the minimal surface equation.
\newblock {\em Ann. Sc. Norm. Super. Pisa Cl. Sci. (5)}, 14(4):1231--1237,
  2015.

\bibitem{F17}
Alberto Farina.
\newblock A sharp bernstein-type theorem for entire minimal graphs.
\newblock {\em arXiv preprint arXiv:1707.0411}, 2017.
\newblock To appear in Calc. Var. Partial Differential Equations.

\bibitem{FarV17}
Alberto Farina and Enrico Valdinoci.
\newblock Flatness results for nonlocal minimal cones and subgraphs.
\newblock {\em arXiv preprint arXiv:1706.05701}, 2017.

\bibitem{MX15}
D.~R. Farine, P.~O. Montiglio, and O.~Spiegel.
\newblock From individuals to groups and back: the evolutionary implications of
  group phenotypic composition.
\newblock {\em Trends in Ecol. and Evolut.}, 30(10):609--621, 2015.

\bibitem{MR2985500}
Patricio Felmer and Alexander Quaas.
\newblock Boundary blow up solutions for fractional elliptic equations.
\newblock {\em Asymptot. Anal.}, 78(3):123--144, 2012.

\bibitem{FV17}
Alessio Figalli and Enrico Valdinoci.
\newblock Regularity and {B}ernstein-type results for nonlocal minimal
  surfaces.
\newblock {\em J. Reine Angew. Math.}, 729:263--273, 2017.

\bibitem{MR1028776}
A.~F. Filippov.
\newblock {\em Differential equations with discontinuous righthand sides},
  volume~18 of {\em Mathematics and its Applications (Soviet Series)}.
\newblock Kluwer Academic Publishers Group, Dordrecht, 1988.
\newblock Translated from the Russian.

\bibitem{FSV15}
Alessio Fiscella, Raffaella Servadei, and Enrico Valdinoci.
\newblock Density properties for fractional {S}obolev spaces.
\newblock {\em Ann. Acad. Sci. Fenn. Math.}, 40(1):235--253, 2015.

\bibitem{FS10}
Rupert~L. Frank and Robert Seiringer.
\newblock Sharp fractional {H}ardy inequalities in half-spaces.
\newblock In {\em Around the research of {V}ladimir {M}az'ya. {I}}, volume~11
  of {\em Int. Math. Ser. (N. Y.)}, pages 161--167. Springer, New York, 2010.

\bibitem{1367-2630-15-12-125022}
R.~C. Gerum, B.~Fabry, C.~Metzner, M.~Beaulieu, A.~Ancel, and D.~P. Zitterbart.
\newblock The origin of traveling waves in an emperor penguin huddle.
\newblock {\em New J. of Phys.}, 15(12):125022, 2013.

\bibitem{fish}
Y.~Gheraibia and A.~Moussaoui.
\newblock {\em Penguins Search Optimization Algorithm (PeSOA)}.
\newblock in Recent Trends in Applied Artificial Intelligence: 26th
  International Conference on Industrial, Engineering and Other Applications of
  Applied Intelligent Systems, IEA/AIE 2013, Amsterdam, The Netherlands, June
  17-21, 2013. Springer Berlin Heidelberg, 2013.

\bibitem{GG82}
Mariano Giaquinta and Enrico Giusti.
\newblock On the regularity of the minima of variational integrals.
\newblock {\em Acta Math.}, 148:31--46, 1982.

\bibitem{GiaMart12}
Mariano Giaquinta and Luca Martinazzi.
\newblock {\em An introduction to the regularity theory for elliptic systems,
  harmonic maps and minimal graphs}, volume~11 of {\em Appunti. Scuola Normale
  Superiore di Pisa (Nuova Serie) [Lecture Notes. Scuola Normale Superiore di
  Pisa (New Series)]}.
\newblock Edizioni della Normale, Pisa, second edition, 2012.

\bibitem{GilTru}
David Gilbarg and Neil~S. Trudinger.
\newblock {\em Elliptic partial differential equations of second order}.
\newblock Classics in Mathematics. Springer-Verlag, Berlin, 2001.
\newblock Reprint of the 1998 edition.

\bibitem{Giling}
D.~Giling, R.~D. Reina, and Z.~Hogg.
\newblock Anthropogenic influence on an urban colony of the little penguin {\it
  eudyptula minor}.
\newblock {\em Marine Freshwater Res.}, 59:647--651, 2008.

\bibitem{Giusti}
Enrico Giusti.
\newblock {\em Minimal surfaces and functions of bounded variation}, volume~80
  of {\em Monographs in Mathematics}.
\newblock Birkh\"auser Verlag, Basel, 1984.

\bibitem{G03}
Enrico Giusti.
\newblock {\em Direct methods in the calculus of variations}.
\newblock World Scientific Publishing Co., Inc., River Edge, NJ, 2003.

\bibitem{GMM97}
E.~Gonzalez, U.~Massari, and M.~Miranda.
\newblock On minimal cones.
\newblock {\em Appl. Anal.}, 65(1-2):135--143, 1997.

\bibitem{BV412}
R.~Hegselmann and U.~Krause.
\newblock Opinion dynamics and bounded confidence: models, analysis and
  simulation.
\newblock {\em J. Artif. Soc. and Soc. Simul.}, 5(3):1--33, 2002.

\bibitem{MR3244289}
Qingkai Kong.
\newblock {\em A short course in ordinary differential equations}.
\newblock Universitext. Springer, Cham, 2014.

\bibitem{CH6}
S.~M. Laaksonen, A.~Chiaradia, and R.~D. Reina.
\newblock Behavioural plasticity of a multihabitat animal, the little penguin,
  eudyptula minor, in response to tidal oscillations on its interhabitat
  transitions.
\newblock {\em Preprint}, 2016.

\bibitem{Lin16}
Erik Lindgren.
\newblock H\"older estimates for viscosity solutions of equations of fractional
  {$p$}-{L}aplace type.
\newblock {\em NoDEA Nonlinear Differential Equations Appl.}, 23(5):Art. 55,
  18, 2016.

\bibitem{Extension}
Luca Lombardini.
\newblock An introduction to the extension operator of the fractional
  {L}aplacian.
\newblock {\em in preparation}.

\bibitem{mine_cyl_stuff}
Luca Lombardini.
\newblock Approximation of sets of finite fractional perimeter by smooth sets
  and comparison of local and global $s$-minimal surfaces.
\newblock {\em Interfaces Free Bound.}, 20(2):261--296, 2018.

\bibitem{Myfractal}
Luca Lombardini.
\newblock Fractional perimeters from a fractal perspective.
\newblock {\em arXiv preprint arXiv:1508.06241}, 2018.
\newblock To appear in Advanced Nonlinear Studies.

\bibitem{CH3}
A.~J.~J. Macintosh, L.~Pelletier, A.~Chiaradia, A.~Kato, and Y.~Ropert-Coudert.
\newblock Temporal fractals in seabird foraging behaviour: diving through the
  scales of time.
\newblock {\em Sci. Rep.}, 3(1884):1--10, 2013.

\bibitem{Maggi}
Francesco Maggi.
\newblock {\em Sets of finite perimeter and geometric variational problems},
  volume 135 of {\em Cambridge Studies in Advanced Mathematics}.
\newblock Cambridge University Press, Cambridge, 2012.
\newblock An introduction to geometric measure theory.

\bibitem{Mattila}
Pertti Mattila.
\newblock {\em Geometry of sets and measures in {E}uclidean spaces}, volume~44
  of {\em Cambridge Studies in Advanced Mathematics}.
\newblock Cambridge University Press, Cambridge, 1995.
\newblock Fractals and rectifiability.

\bibitem{bourgbrez}
V.~{Maz'ya} and T.~{Shaposhnikova}.
\newblock {On the Bourgain, Brezis, and Mironescu theorem concerning limiting
  embeddings of fractional Sobolev spaces.}
\newblock {\em {J. Funct. Anal.}}, 195(2):230--238, 2002.

\bibitem{MX13b}
N.~Miller, S.~Garnier, and I.~D. Couzin.
\newblock Both information and social cohesion determine collective decisions
  in animal groups.
\newblock {\em Proc. Natl. Acad. Sci. USA}, 110(13):5263--5268.

\bibitem{M64}
Mario Miranda.
\newblock Distribuzioni aventi derivate misure insiemi di perimetro localmente
  finito.
\newblock {\em Ann. Scuola Norm. Sup. Pisa (3)}, 18:27--56, 1964.

\bibitem{M61}
J\"urgen Moser.
\newblock On {H}arnack's theorem for elliptic differential equations.
\newblock {\em Comm. Pure Appl. Math.}, 14:577--591, 1961.

\bibitem{Ponce}
Augusto~C. Ponce.
\newblock A new approach to {S}obolev spaces and connections to
  {$\Gamma$}-convergence.
\newblock {\em Calc. Var. Partial Differential Equations}, 19(3):229--255,
  2004.

\bibitem{MR3235841}
Romain Ragonnet, Romain Jumentier, and Bernhard Beckermann.
\newblock La marche de l'empereur.
\newblock {\em Matapli}, (102):71--82, 2013.

\bibitem{CH4}
A.~M. Reynolds, Y.~Ropert-Coudert, A.~Kato, A.~Chiaradia, and A.~J.~J.
  MacIntosh.
\newblock A priority-based queuing process explanation for scale-free foraging
  behaviours.
\newblock {\em Anim. Behav.}, (108):67--71, 2015.

\bibitem{CH5}
A.~Rodr\'{\i}guez, A.~Chiaradia, P.~Wasiak, L.~Renwick, and P.~Dann.
\newblock Waddling on the dark side: Ambient light affects attendance behavior
  of little penguins.
\newblock {\em J. Biol. Rhythms}, 31(2):194--204, 2016.

\bibitem{MR3168912}
Xavier Ros-Oton and Joaquim Serra.
\newblock The {D}irichlet problem for the fractional {L}aplacian: regularity up
  to the boundary.
\newblock {\em J. Math. Pures Appl. (9)}, 101(3):275--302, 2014.

\bibitem{MR3436398}
Xavier Ros-Oton and Enrico Valdinoci.
\newblock The {D}irichlet problem for nonlocal operators with singular kernels:
  convex and nonconvex domains.
\newblock {\em Adv. Math.}, 288:732--790, 2016.

\bibitem{SV12}
Ovidiu Savin and Enrico Valdinoci.
\newblock {$\Gamma$}-convergence for nonlocal phase transitions.
\newblock {\em Ann. Inst. H. Poincar\'e Anal. Non Lin\'eaire}, 29(4):479--500,
  2012.

\bibitem{SV13}
Ovidiu Savin and Enrico Valdinoci.
\newblock Regularity of nonlocal minimal cones in dimension 2.
\newblock {\em Calc. Var. Partial Differential Equations}, 48(1-2):33--39,
  2013.

\bibitem{SV14}
Ovidiu Savin and Enrico Valdinoci.
\newblock Density estimates for a variational model driven by the {G}agliardo
  norm.
\newblock {\em J. Math. Pures Appl. (9)}, 101(1):1--26, 2014.

\bibitem{Sickel}
Winfried Sickel.
\newblock Pointwise multipliers of {L}izorkin-{T}riebel spaces.
\newblock In {\em The {M}az'ya anniversary collection, {V}ol. 2 ({R}ostock,
  1998)}, volume 110 of {\em Oper. Theory Adv. Appl.}, pages 295--321.
  Birkh\"auser, Basel, 1999.

\bibitem{MR2926716}
L.~A. Sidhu, E.~A. Catchpole, and P.~Dann.
\newblock Modelling banding effect and tag loss for {L}ittle {P}enguins {\it
  {e}udyptula minor}.
\newblock {\em ANZIAM J. Electron. Suppl.}, 52((C)):C206--C221, 2010.

\bibitem{MR2718048}
Elizabeth~A. Skewgar.
\newblock {\em Behavior of {M}agellanic penguins at sea}.
\newblock ProQuest LLC, Ann Arbor, MI, 2009.
\newblock Thesis (Ph.D.)--University of Washington.

\bibitem{Tamanini}
Italo Tamanini.
\newblock Regularity results for almost minimal oriented hypersurfaces in
  $\mathbb{R}^n$.
\newblock {\em Quaderni del {D}ipartimento di {M}atematica, {U}niversit\'a di
  {L}ecce}, 1984.

\bibitem{V13}
Enrico Valdinoci.
\newblock A fractional framework for perimeters and phase transitions.
\newblock {\em Milan J. Math.}, 81(1):1--23, 2013.

\bibitem{Visintin}
Augusto Visintin.
\newblock Generalized coarea formula and fractal sets.
\newblock {\em Japan J. Indust. Appl. Math.}, 8(2):175--201, 1991.

\bibitem{weiss}
Georg~Sebastian Weiss.
\newblock Partial regularity for a minimum problem with free boundary.
\newblock {\em J. Geom. Anal.}, 9(2):317--326, 1999.

\end{thebibliography}

\end{document}